\documentclass[12pt]{article}

\usepackage[margin=1in]{geometry}  
\usepackage{graphicx}              
\usepackage{amsmath}               
\usepackage{amsfonts}              
\usepackage{amsthm}                
\usepackage{amssymb}

\usepackage{tikz}
\usetikzlibrary{arrows}
\newcommand{\midarrow}{\tikz \draw[-tt] (0,0) -- +(.1,0);}

\newcommand{\arr}[2]{ \draw (#1)-- node[sloped,allow upside down] {\midarrow} (#2); }

\newlength{\arrowsizez}  
\pgfarrowsdeclare{tt}{tt}{  
  \setlength{\arrowsizez}{0.4pt}  
  \addtolength{\arrowsizez}{.5\pgflinewidth}  
  \pgfarrowsrightextend{0}  
  \pgfarrowsleftextend{-5\arrowsizez}  
}{  
  \setlength{\arrowsizez}{0.4pt}  
  \addtolength{\arrowsizez}{.5\pgflinewidth}  
  \pgfpathmoveto{\pgfpoint{-5\arrowsizez}{2\arrowsizez}}  
  \pgfpathlineto{\pgfpointorigin}  
  \pgfpathlineto{\pgfpoint{-5\arrowsizez}{-2\arrowsizez}}  
  \pgfusepathqstroke  
}

\renewcommand{\AA}{\mathcal{A}}
\newcommand{\GG}{G_{\mathcal{A}}}
\newcommand{\AAN}[1]{\mathcal{A}_{#1}}
\newcommand{\ZZ}[1]{\mathbb Z / {#1} \mathbb Z}

\newcommand{\om}{{\omega}}

\newcommand{\BB}{\mathcal{B}}

\newcommand{\AAbar}{\overline{\mathcal{A}}}
\newcommand{\sqr}[1]{\sqrt[{#1}]{1}}

\begin{document}

\nocite{*}

\title{Abelian automorphism groups of cubic fourfolds.}

\author{Evgeny Mayanskiy}

\maketitle

\begin{abstract}
  We list all finite abelian groups which act effectively on smooth cubic fourfolds. 
\end{abstract}

\section{Introduction.}

Automorphism groups are important geometric invariants of algebraic varieties. Let us mention that for cubic surfaces they were studied and computed by Segre \cite{Segre} and Dolgachev and Verra \cite{DolgachevVerra} (see also \cite{CAG} and \cite{DolgachevIskovskih}). The case of $K3$ surfaces was considered by Nikulin \cite{Nikulin}.\\

Cubic fourfolds form another interesting class of algebraic varieties. Their {\it symplectic} automorphisms (i.e. cyclic automorphism groups) were described in a recent work of Lie Fu \cite{LieFu}. The prime order automorphisms of cubic hypersurfaces in general is a subject of an earlier paper by Gonz\'{a}lez-Aguilera and Liendo \cite{Liendo}, where the authors in particular listed all prime order automorphisms of cubic threefolds and cubic fourfolds.\\ 

In this note we address the problem of computing all {\it abelian} automorphism groups of smooth cubic fourfolds. We use essentially the same method as in \cite{CAG}, \cite{Liendo} and \cite{LieFu}. Its starting point is the observation that every finite order endomorphism of a finite-dimensional vector space is semisimple. Note that the same idea, at least in principle, leads to a general algorithm for computing finite order automorphisms (and finite abelian automorphism groups) of {\it arbitrary} (for example, smooth or with some prescribed singularities) hypersurfaces of any degree and dimension possibly satisfying some additional conditions (for example, symplectic as in \cite{LieFu}).\\

Our main result can be formulated as follows.\\

{\bf Theorem.} {\it A finite abelian group acts effectively on a smooth cubic fourfold if and only if it is a subgroup of one of the following groups:
\begin{center}
\begin{multline*}
(\ZZ{3})^{\oplus 5},\;\; (\ZZ{3})^{\oplus 4}\oplus \ZZ{2},\;\; (\ZZ{3})^{\oplus 3}\oplus \ZZ{4},\;\; (\ZZ{3})^{\oplus 3}\oplus (\ZZ{2})^{\oplus 2},\\
(\ZZ{3})^{\oplus 2}\oplus \ZZ{8},\;\; (\ZZ{3})^{\oplus 2}\oplus \ZZ{2} \oplus \ZZ{4},\;\; (\ZZ{3})^{\oplus 2}\oplus (\ZZ{2})^{\oplus 3},\;\; \ZZ{3} \oplus \ZZ{16},\\
\ZZ{3} \oplus \ZZ{2}\oplus \ZZ{8},\;\;\; \ZZ{3} \oplus (\ZZ{4})^{\oplus 2},\;\;\; \ZZ{32},\\
(\ZZ{3})^{\oplus 2}\oplus \ZZ{9},\;\;\; \ZZ{2}\oplus \ZZ{3}\oplus \ZZ{9}, \;\;\; \ZZ{4}\oplus \ZZ{9},\\
(\ZZ{3})^{\oplus 2}\oplus \ZZ{5},\;\; \ZZ{2}\oplus \ZZ{3}\oplus \ZZ{5},\;\; \ZZ{3}\oplus \ZZ{7},\;\;\; \ZZ{3}\oplus \ZZ{11}.
\end{multline*}
\end{center}}

A more precise version of this Theorem (including the description of all possible group actions upto a linear change of variables as well as families of smooth cubic fourfolds admitting those actions) is given in Corollary $1$, Theorem $3$ and Theorem $4$ below.\\

\vspace{2cm}

\vspace{1cm}

\subsection{Notation.}

Let $n\geq 1, d\geq 3$. In this note we say that a group $G$ is an {\it automorphism group} of a degree $d$ and dimension $n$ hypersurface, if there is a group monomorphism $\phi \colon G\hookrightarrow Aut({\mathbb P}^{n+1})$ and a smooth hypersurface $X\subset {{\mathbb P}^{n+1}}$ of degree $d$ such that ${\phi}(G)(X)\subset X$.\\

Note that such an 'automorphism group' $G$ is nothing else but a subgroup of the (finite) group of {\it linear} automorphisms of $X$ denoted by $Lin(X)$ in \cite{Poonen}, \cite{MatsumuraMonsky}. The latter coincides with the usual automorphism group $Aut(X)$ of $X$ considered as an abstract variety unless $(n,d)=(1,3)$ or $(n,d)=(2,4)$ (see \cite{Poonen}, \cite{MatsumuraMonsky} for the proof).\\

From now on we will consider only the case $(n,d)=(4,3)$ (note, however, that our methods should work, at least in principle, for any pair $(n,d)$). Our base field is $\mathbb C$ throughout the note.\\

A primitive $k-$th root of unity will be denoted by $\sqrt[k]{1}$ and we will take $\sqrt[k]{1}=e^{2\pi i /k}$ for the sake of conreteness.\\

We will denote a diagonal matrix $A$ with diagonal entries $a_1, a_2,...,a_s$ by 
$$
A=(a_1,a_2,...,a_s).
$$

Given variables $x,y,...$ we will denote by $f_k(x,y,...)$ an arbitrary $k$-form in $x,y,...$.\\

\section{Diagonalizable abelian automorphism groups.}

We start by computing all {\it diagonalizable} abelian automorphism groups of cubic fourfolds, i.e. finite groups $G$ such that the monomorphism $\phi \colon G\hookrightarrow Aut({\mathbb P}^{5})\cong PGL(6)$ factors through the projection $GL(6)\rightarrow PGL(6)$.\\

This assumption allows us to use our key remark that every element $f\in G$ of a such a group $G$ is conjugate (as an element of $GL(6)$) to a diagonal matrix
$$
diag({\om}^{c_0},{\om}^{c_1},{\om}^{c_2},{\om}^{c_3},{\om}^{c_4},{\om}^{c_5})
$$
(which we will denote by $({\om}^{c_0},{\om}^{c_1},{\om}^{c_2},{\om}^{c_3},{\om}^{c_4},{\om}^{c_5})$ in future), where $\om$ is a $k-$th root of unity (with $k=ord(f)$ being the order of $f$) and $(c_0,c_1,c_2,c_3,c_4,c_5)\in (\mathbb Z / k\mathbb Z)^{\oplus 6}$.\\

The sequence $(c_0,c_1,c_2,c_3,c_4,c_5)$ is defined upto adding an arbitrary sequence of the form $(c,c,c,c,c,c)\in (\mathbb Z / k\mathbb Z)^{\oplus 6}$ (with the same $c \in \mathbb Z / k\mathbb Z$). This sequence was called the {\it signature} of $f\in G$ in \cite{Liendo}.\\

{\bf Theorem 1.} {\it Diagonalizable abelian automorphism groups of smooth cubic fourfolds are exactly the subgroups of the following groups (see Table $1$ below).
\begin{center}
  \includegraphics[width=1.0\textwidth]{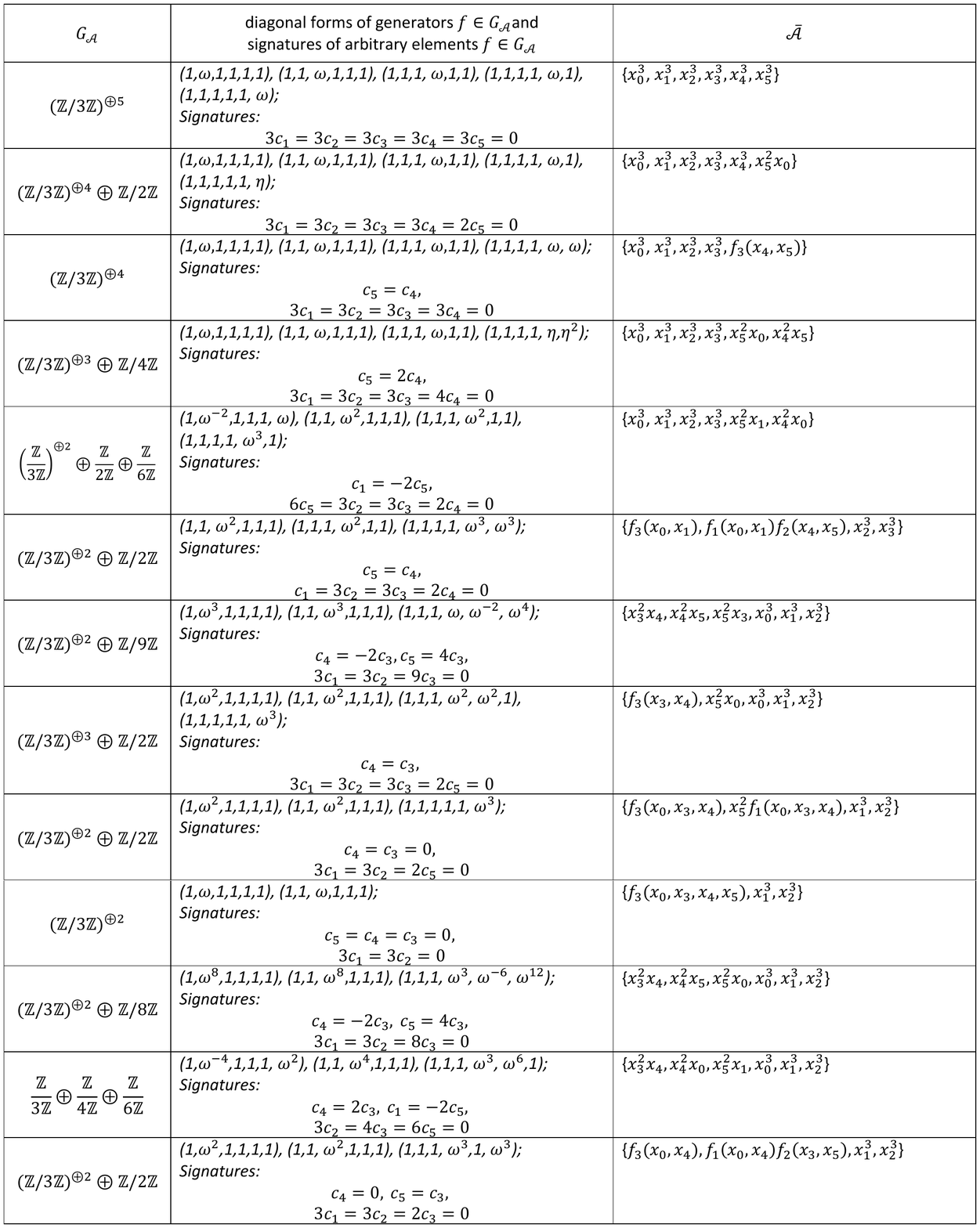}
\end{center}

\begin{center}
  \includegraphics[width=1.0\textwidth]{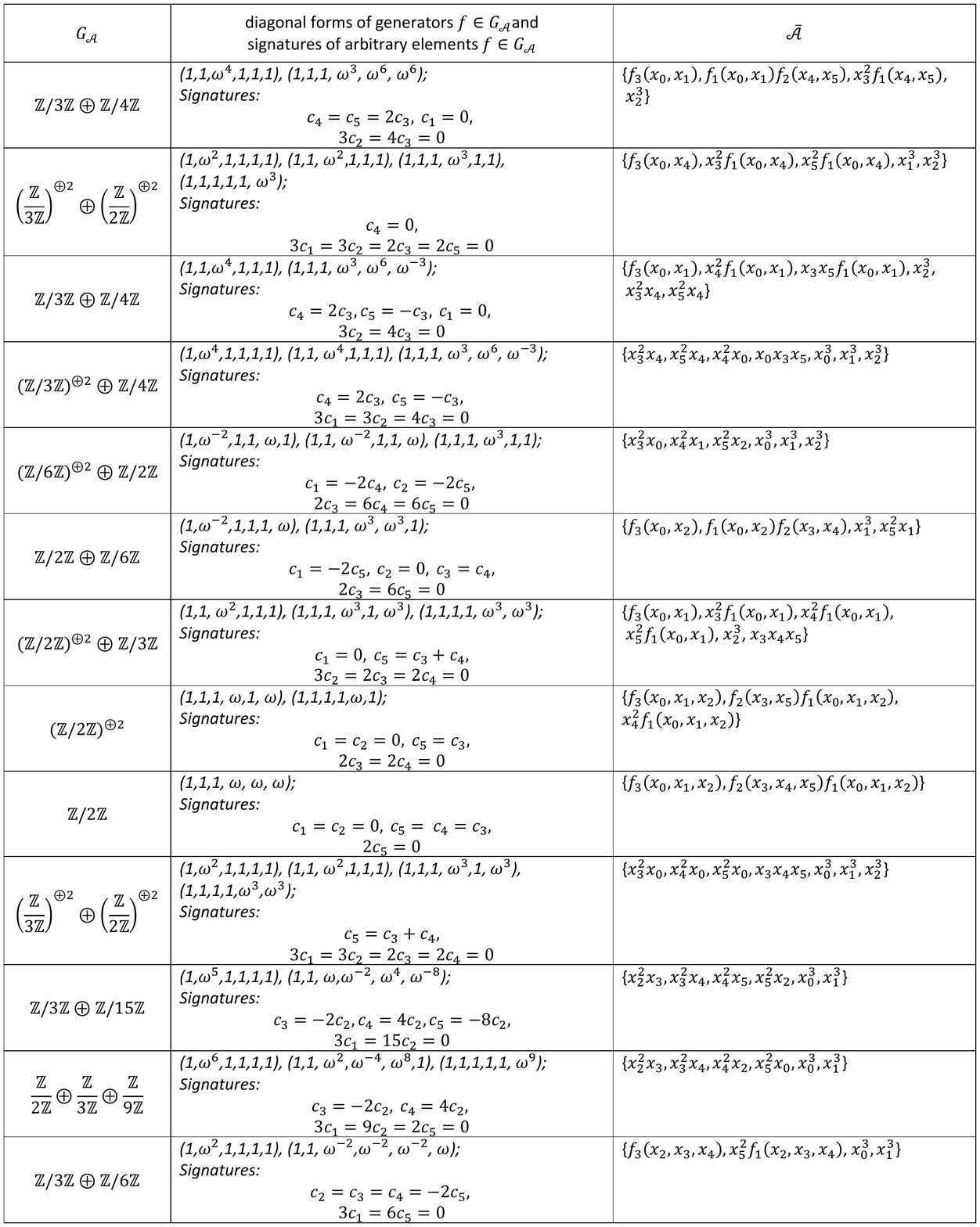}
\end{center}

\begin{center}
  \includegraphics[width=1.0\textwidth]{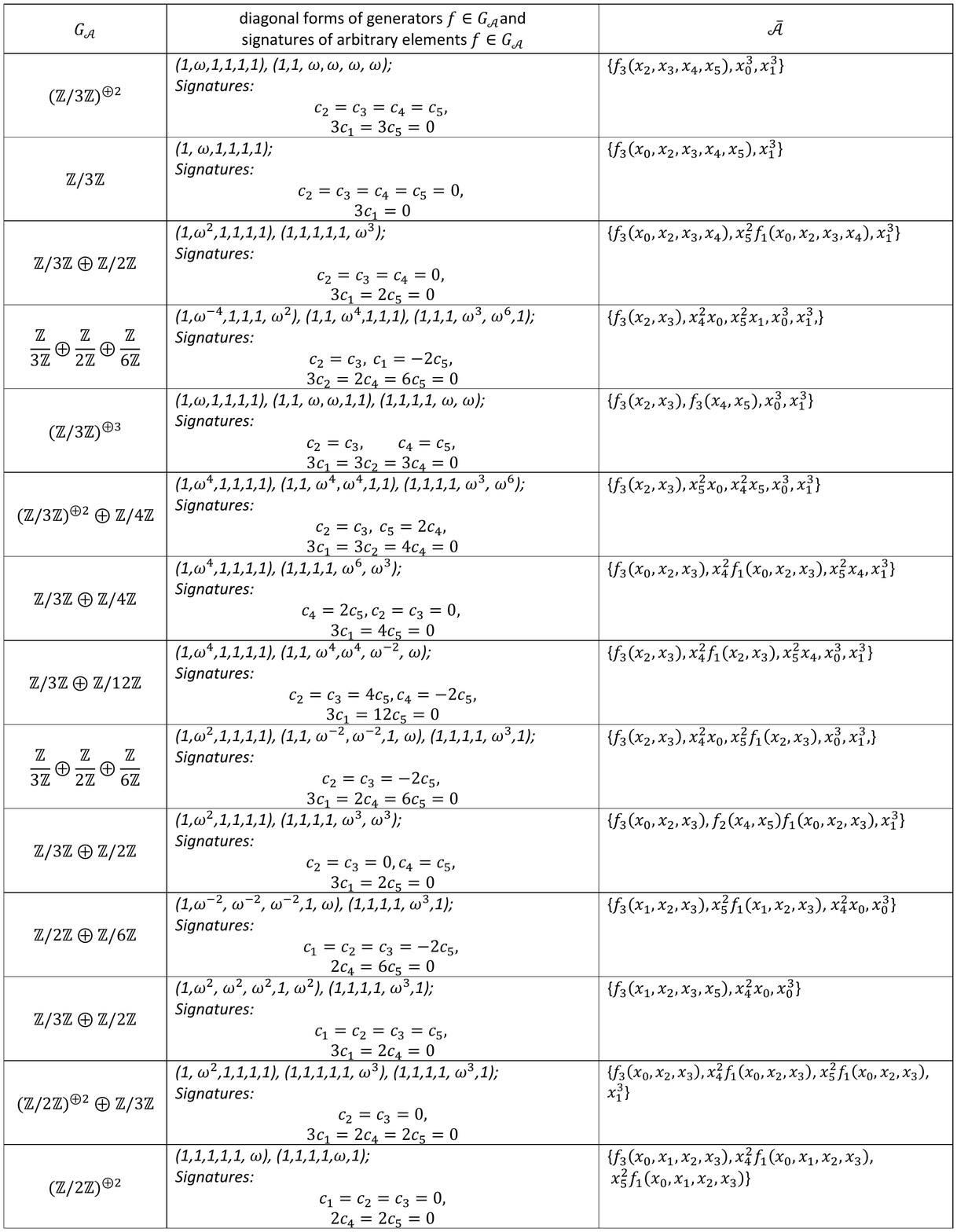}
\end{center}

\begin{center}
  \includegraphics[width=1.0\textwidth]{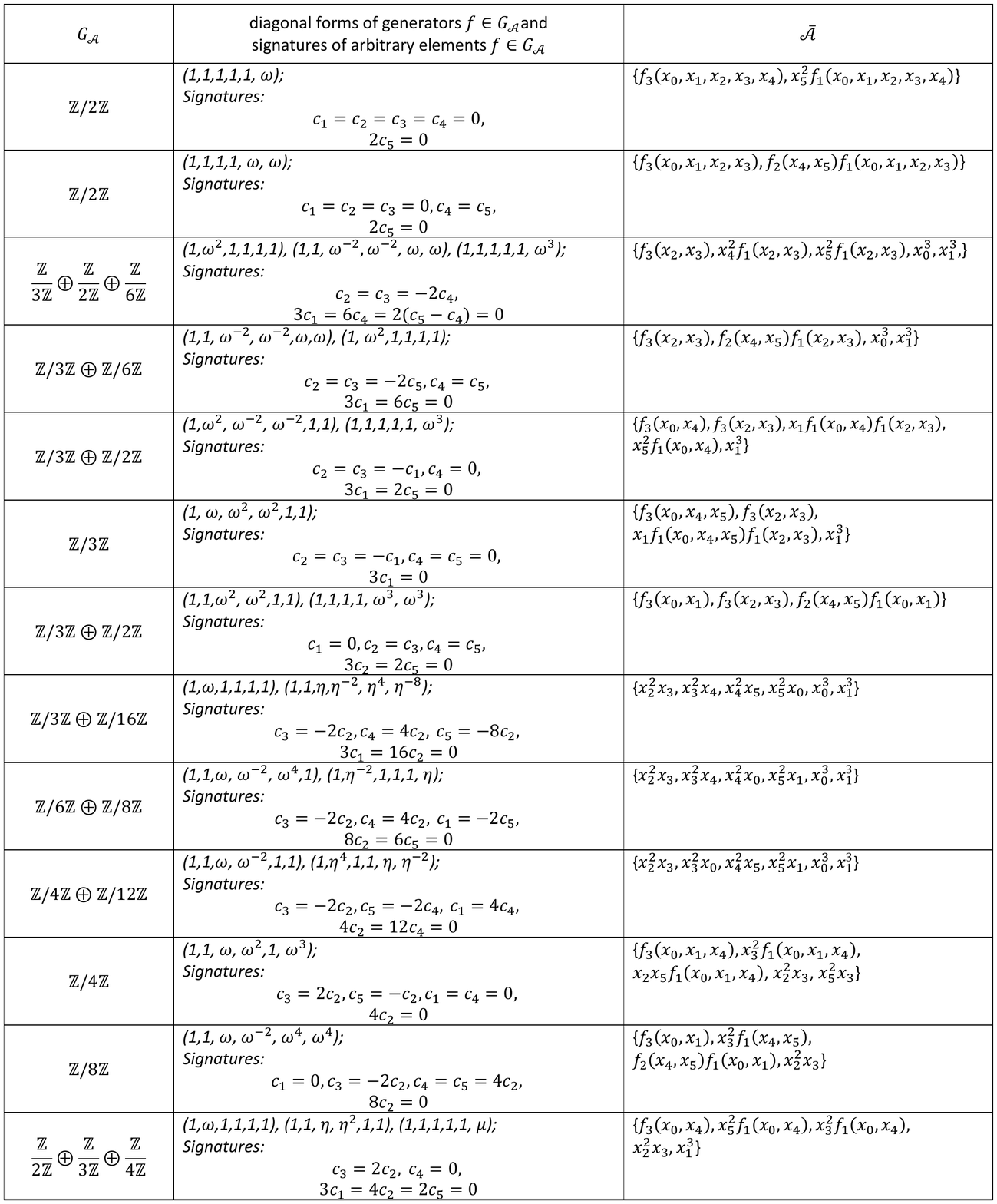}
\end{center}

\begin{center}
  \includegraphics[width=1.0\textwidth]{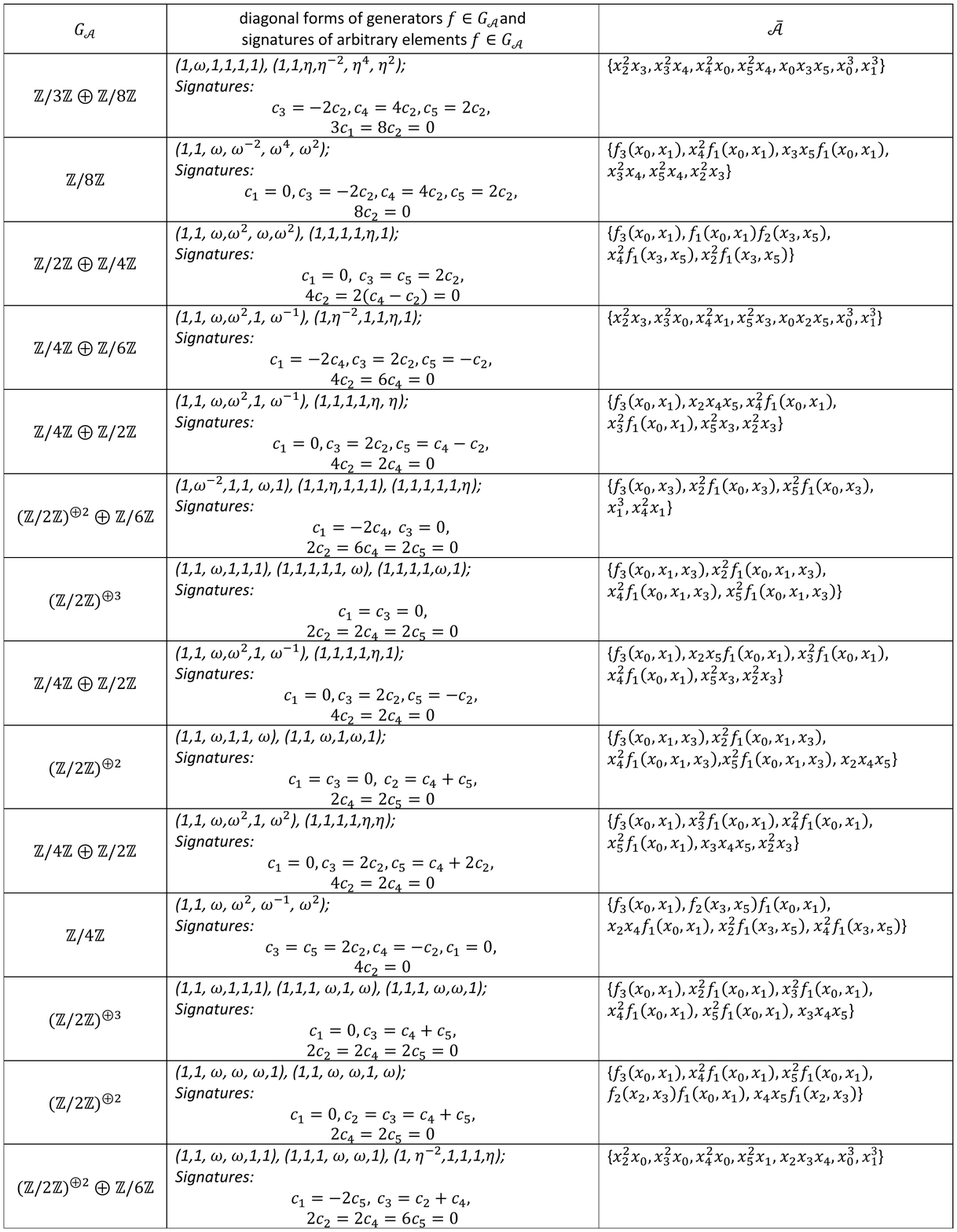}
\end{center}

\begin{center}
  \includegraphics[width=1.0\textwidth]{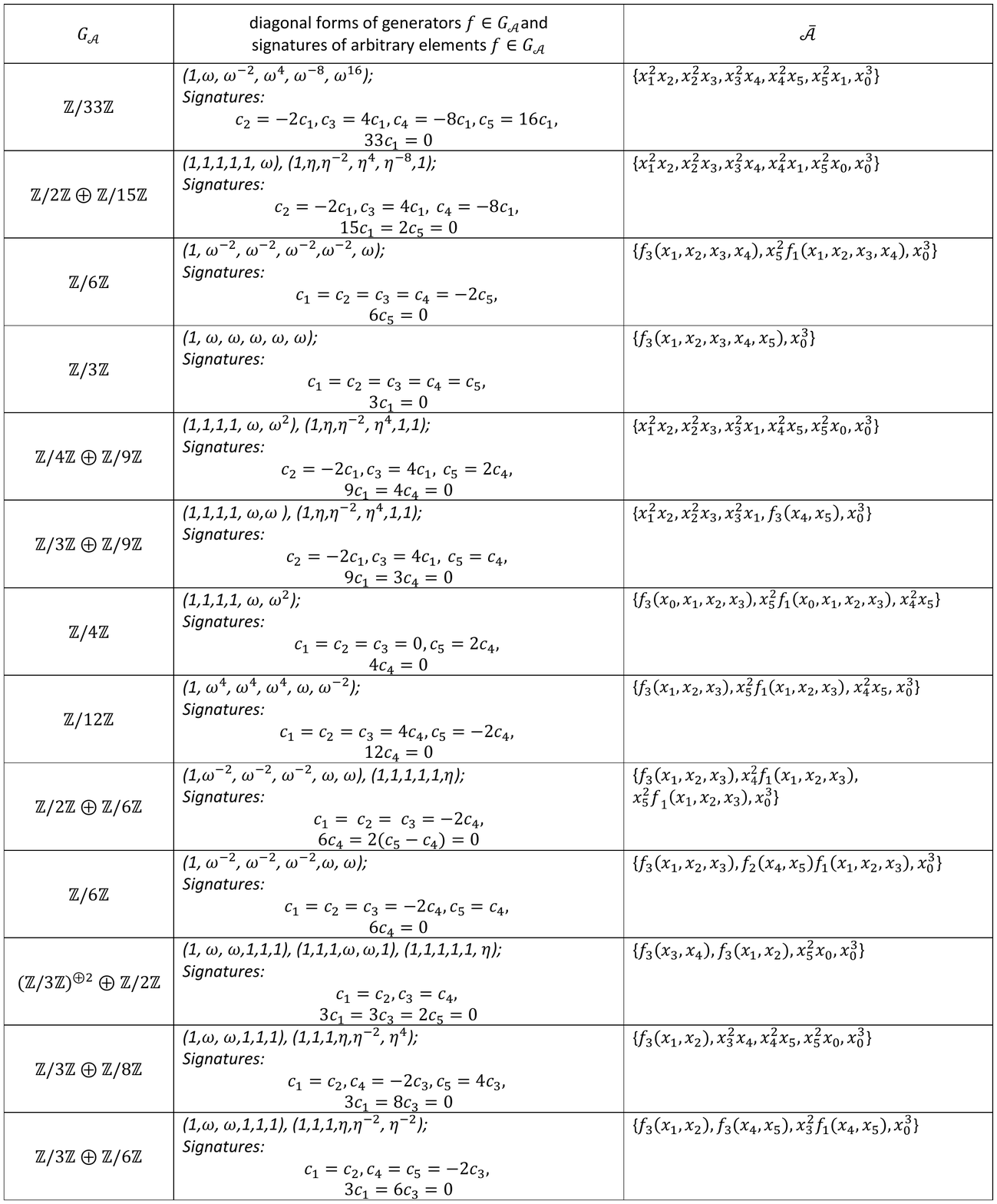}
\end{center}

\begin{center}
  \includegraphics[width=1.0\textwidth]{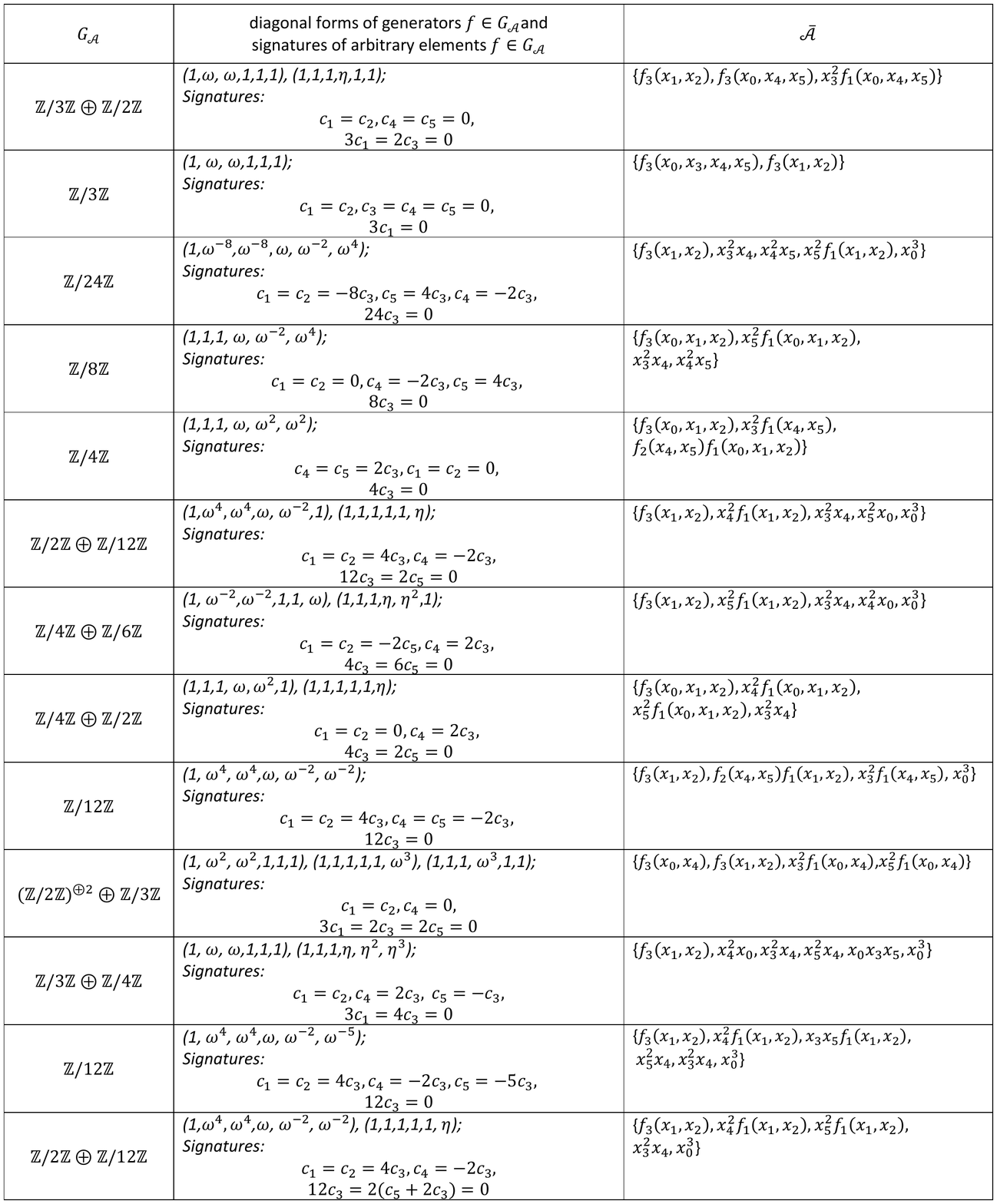}
\end{center}

\begin{center}
  \includegraphics[width=1.0\textwidth]{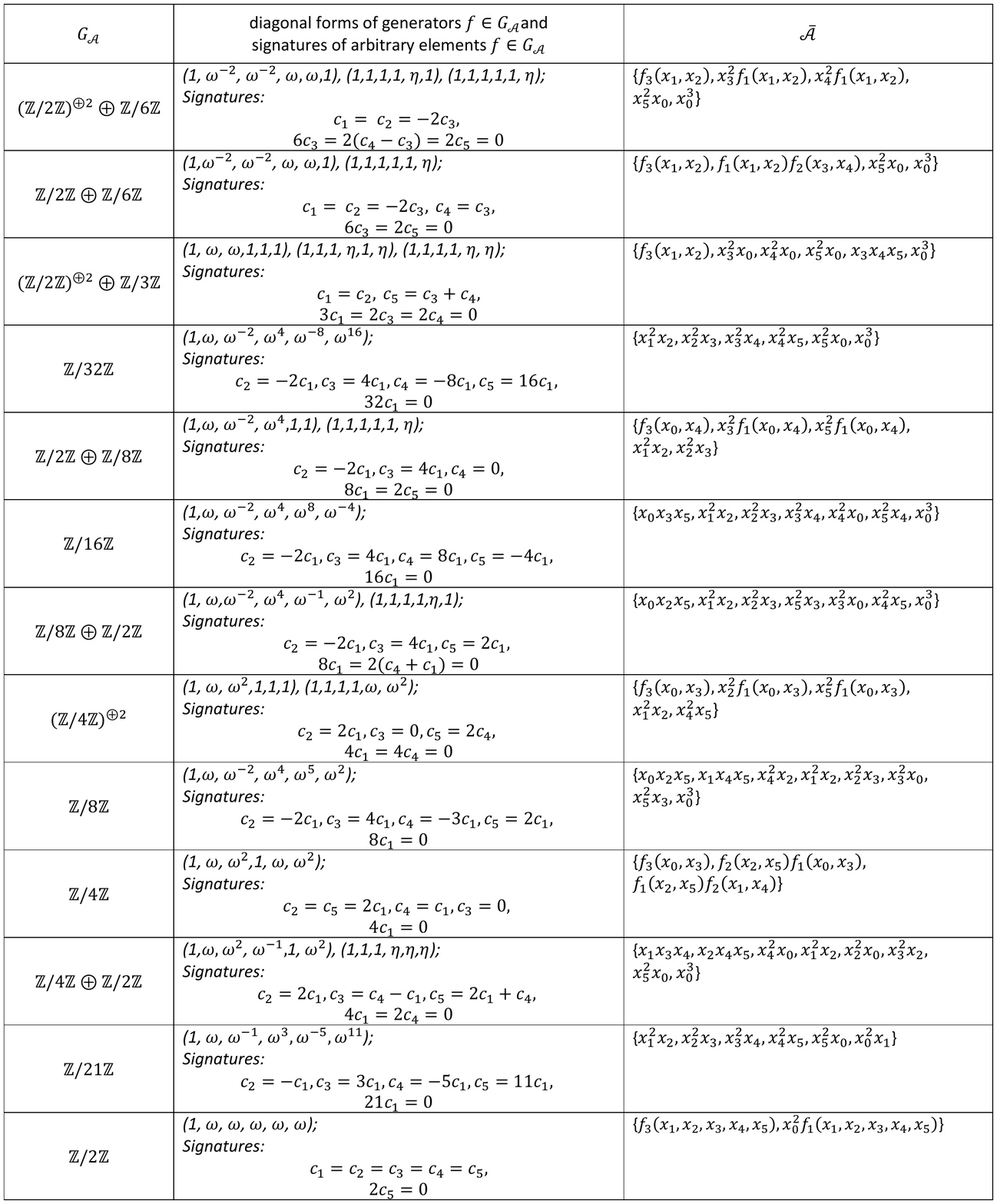}
\end{center}

\begin{center}
  \includegraphics[width=1.0\textwidth]{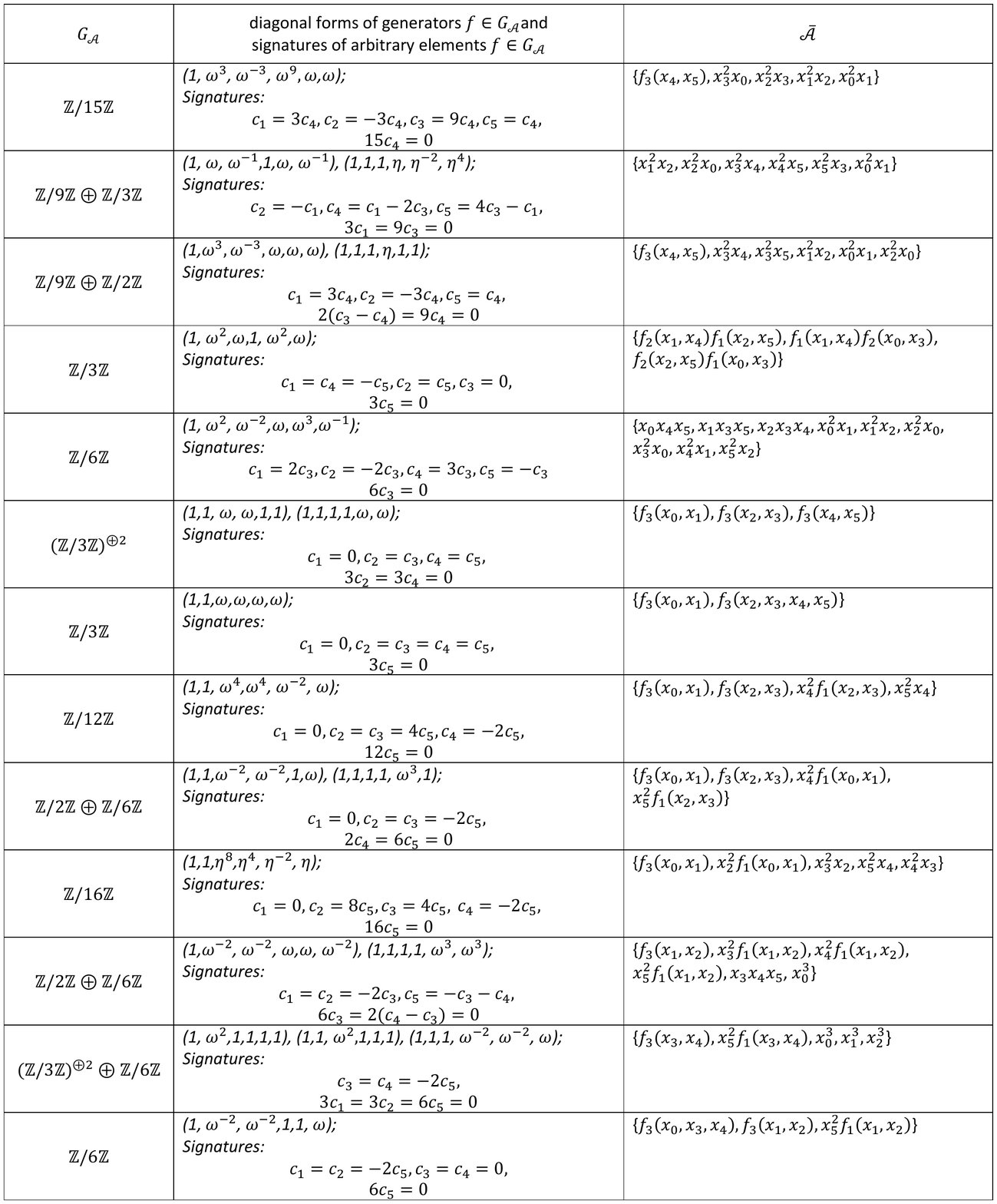}
\end{center}

\begin{center}
  \includegraphics[width=1.0\textwidth]{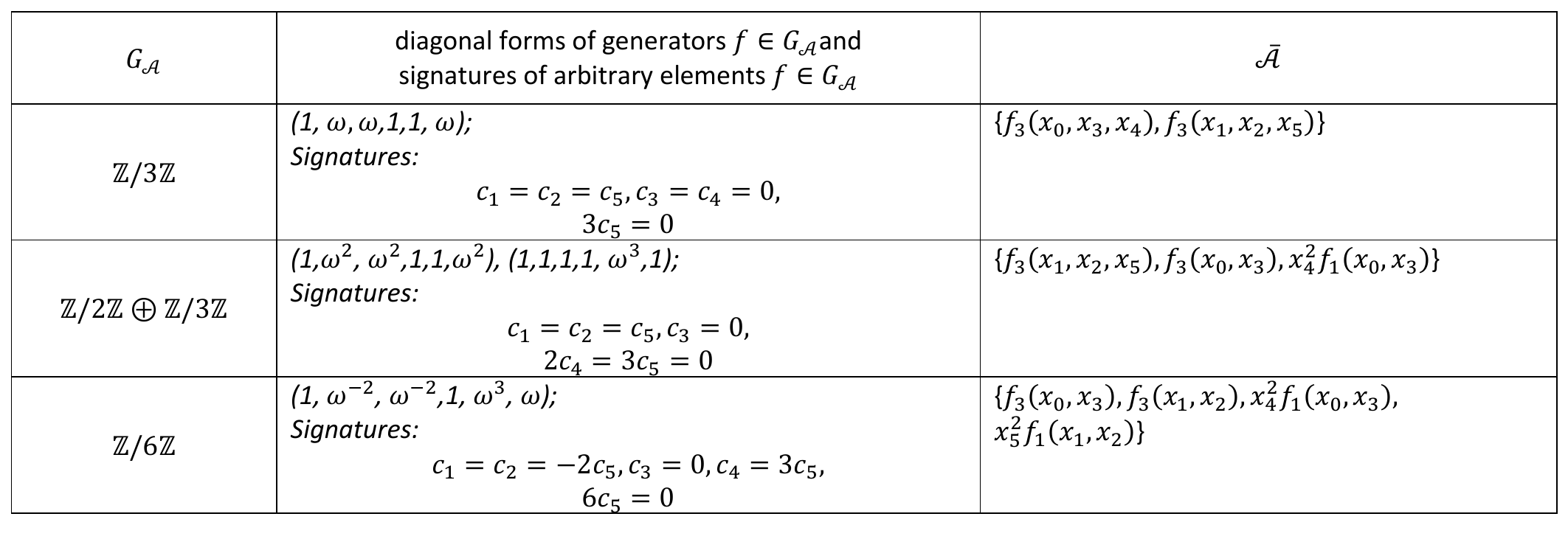}
\end{center}

{\bf Table $1$.} {\sl Diagonalizable abelian automorphism groups of smooth cubic fourfolds.} The first column lists abelian groups $G$ which act effectively on a smooth cubic fourfold. The second column describes generators of $G$ (diagonalized after possibly conjugating by an element of $GL(6)$) as diagonal matrices in $GL(6)$ as well as the conditions on the signatures of arbitrary elements of $G$. The third column gives maximal sets of cubic monomials, which lie in the same eigenspace of each of these generators acting correspondingly on $\text{Sym}^3 ({\mathbb C}^6)$ (the cubic fourfolds with automorphism group $G_{\mathcal A}$ are exactly\footnote{Note that the same group may appear (as a subgroup) in several rows of the table. If it appears as a proper subgroup in a different row, then one has to add to $\AAbar$ in that row all monomials lying in the same eigenspaces of the subgroup elements as the monomials of the original $\AAbar$.\par On the other hand, for each of our group actions one can describe directly the (smooth) sets $\mathcal B\subset \mathcal M$ of monomials lying in each of the eigenspaces of the generators. We will do such computation in the proof of Corollary $1$ for some of the rows of Table $1$ (which will constitute after that Table $2$).\par One may also extend Table $1$ by adding to sets $\AAbar$ various monomials and computing extra restrictions they impose on the elements of the group $\GG$, the resulting subgroups of $\GG$ and the corresponding enlarged sets $\AAbar$.} the ones which in suitable coordinates can be represented as vanishing loci of linear combinations of monomials in $\bar{\mathcal A}$.)}\\

Let us explain our approach to proving Theorem $1$.\\

Let $V={\mathbb C}^6$ with coordinates $x_0,x_1,x_2,x_3,x_4,x_5$. Let $\mathcal M$ be the set of all cubic monomials in $x_i$ (i.e. the set of the corresponding coordinates in $\text{Sym}^3 V$).\\

We will call a subset $\AA \subset \mathcal M$ {\it smooth}, if $\AA$ is the set of monomials of a cubic form $F$ with smooth vanishing locus $X\subset {\mathbb P}(V)$.\\

We will describe all {\it minimal smooth} sets $\AA \subset \mathcal M$. For each of them we will determine the group $G_{\mathcal A}$ of diagonal linear automorphisms $f=({\om}^{c_0},{\om}^{c_1},{\om}^{c_2},{\om}^{c_3},{\om}^{c_4},{\om}^{c_5})$ of $V$, which leave the corresponding cubics $X\subset {\mathbb P}(V)$ invariant. Here $\om$ denotes a primitive $k-$th root of unity and $c_i\in {\mathbb Z}/k{\mathbb Z}$, where $k=ord(f)$ is the order of $f$.\\

Group $G_{\mathcal A}$ will in turn determine a larger set of monomials $\overline{\mathcal A}\subset \mathcal M$ such that $G_{\mathcal A}$ is an automorphism group of a smooth cubic fourfold $X\subset {\mathbb P}(V)$ if and only if\footnote{For the 'only if' part we require that the monomials of a cubic form defining $X$ lie in the same eigenspaces of the elements of $\GG$ as the monomials in $\AA$.} $X$ can be defined in the chosen coordinate system in $V$ by a cubic form all of whose monomials lie in $\overline{\mathcal A}$.\\

More concretely, we define the $f-$weight of a monomial $m=x_0^{i_0}x_1^{i_1}x_2^{i_2}x_3^{i_3}x_4^{i_4}x_5^{i_5}$ to be\footnote{For this definition to make sense, we have to consider $f$ {\it together} with a choice of the signature $(c_0,...,c_5)\in ({\mathbb Z}/k{\mathbb Z})^{\oplus 6}$. However, since we will be working only with one $f$ at a time, the ambiguity of the choice of $c_i$ will not cause any difficulties.} 
$$
wt_f(m)=\sum_{j=0}^{5} i_j\cdot c_j \in \mathbb Z / k \mathbb Z.
$$

Then elements $f$ of $G_{\mathcal A}$ are given by the condition that $f-$weights of all monomials in $\AA$ are the same (as elements of $\mathbb Z / k \mathbb Z$). Note that without loss of generality we may assume that $c_0=0$.\\

This gives us a system of linear equations in $c_0,c_1,c_2,c_3,c_4,c_5$ over  $\mathbb Z / k \mathbb Z$, whose solution set is $G_{\mathcal A}$.\\

The set $\overline{\mathcal A}\subset \mathcal M$ will consist of all monomials  $m\in \mathcal M$ whose $f$-weight is the same as the $f$-weight of monomials in $\AA$ (for any $f\in G_{\mathcal A}$).\\

Now let us describe how we are going to enumerate all minimal smooth sets $\AA \subset \mathcal M$. It is sufficient to do this upto a permutation of coordinates $x_0,x_1,x_2,x_3,x_4,x_5$. We will break down these sets into $7$ groups:
\begin{itemize}
\item sets $\AA$ with exactly $6$ cubes $x_0^3,x_1^3,x_2^3,x_3^3,x_4^3,x_5^3$,
\item sets $\AA$ with exactly $5$ cubes $x_0^3,x_1^3,x_2^3,x_3^3,x_4^3$,
\item sets $\AA$ with exactly $4$ cubes $x_0^3,x_1^3,x_2^3,x_3^3$,
\item sets $\AA$ with exactly $3$ cubes $x_0^3,x_1^3,x_2^3$,
\item sets $\AA$ with exactly $2$ cubes $x_0^3,x_1^3$,
\item sets $\AA$ with exactly $1$ cube $x_0^3$,
\item sets $\AA$ with no cubes .
\end{itemize}

Each of these $7$ groups of sets with exactly $k$ cubes $x_0^3,...,x_{k-1}^3$ will be in turn broken down into subgroups, each of which is specified by the condition that the length of a longest cycle in $\AA$ is exactly $l$.\\

Here we say that $\AA$ {\it has a cycle of length $l$}, if it contains monomials
$$
x_{i_1}^2x_{i_2}, \; x_{i_2}^2x_{i_3},\; ... \; , x_{i_l}^2x_{i_1}
$$
with $i_j\geq k$ and $i_p\neq i_q$ for $p\neq q$.\\

Since we are working upto a permutation of coordinates $x_i$ we may assume that such a longest cycle involves variables $x_k,x_{k+1},...,x_{k+l-1}$.\\

It will be convenient to use graphical notation to denote $\AA$. Each subset $\AA\subset \mathcal M$ can be encoded by an oriented graph with $6$ vertices corresponding to the $6$ coordinates $x_0,x_1,x_2,x_3,x_4,x_5$. If $x_k^3\in \AA$, then the $k$-th vertex will be encircled:

\begin{center}

\end{center}

If $x_k^2x_l \in \AA$, then the graph will contain an arrow from the vertex $k$ to the vertex $l$:
\begin{center}
%
\end{center}

In addition to that, we will use dashed curves in order to connect vertices $k$, $l$ and $m$ whenever $x_kx_lx_m\in \AA$:
\begin{center}
%
\end{center}

We will also use the following shorthands
\begin{center}
%
\end{center}

as well as omit encircled vertices from our pictures if they are not connected by edges (or dashed curves) to other vertices.\\ 

For example, the graph 
\begin{center}
%
\end{center}

represents the set $\AA=\{ x_0^3,x_1^3,x_2^2x_3,x_3^2x_2,x_4^2x_2,x_4^2x_3,x_5^2x_2,x_5^2x_3 \}$. This is a minimal smooth set.\\

In order to keep our pictures comprehensible, we will use the usual symbolism for monomials as well. For example, we will write 
\begin{center}
%
\end{center} 

to say that in order to get a smooth set one has to add to the set described by the graph one of the monomials $a_1,...,a_{p}$, one of the monomials $b_1,...,b_{q}$ and so on.\\

In order to determine when a given subset $\AA\subset \mathcal M$ is smooth, we will use necessary conditions given in the next two lemmas. They follow from the Jacobian criterion.\\

{\bf Lemma 1 (\cite{Liendo}, Lemma 1.3).} {\it In order to be smooth, $\AA$ should contain a monomial $x_i^2x_j$ for any $i$.}\\

{\it Proof:} Let $F=F(x_0,x_1,x_2,x_3,x_4,x_5)$ be a cubic form which contains no monomials of the form $x_0^2x_j$ for any $j$. This means that $F$ can be written as 
$$
F=x_0\cdot P(x_1,x_2,x_3,x_4,x_5)+Q(x_1,x_2,x_3,x_4,x_5)
$$
for some quadratic and cubic forms $P$ and $Q$.\\ 

By the Jacobian criterion its vanishing locus $X=V(F)\subset {\mathbb P}^{5}$ will be singular at the point with homogeneous coordinates $(1:0:0:0:0:0)$. {\it QED}\\

In the next Lemma we will use the following terminology. Assume that $\AA\subset \mathcal M$ contains $k$ cubes $x_0^3,...,x_{k-1}^3$ and no other cubes. We say that $(x_i,x_j)$ for some $k\leq i < j$ is a {\it singular pair}, if
$$
x_i^2x_j, x_j^2x_i \notin \AA
$$
and
$$
\{ x_i^2x_p, x_j^2x_p, x_ix_jx_p \} \cap \AA \neq \emptyset
$$
(i.e. one, two or all three of these three monomials for a given $p$ lie in $\AA$) for at most one value of $p$.\\ 

We say that $(x_i,x_j, x_l)$ for some $k\leq i < j < l$ is a {\it singular triple}, if 
$$
x_ix_jx_l, x_i^2x_j, x_j^2x_i, x_i^2x_l, x_l^2x_i, x_j^2x_l, x_l^2x_j \notin \AA
$$
and
$$
\{  x_i^2x_p, x_j^2x_p, x_l^2x_p, x_ix_jx_p, x_ix_lx_p, x_lx_jx_p \} \cap \AA \neq \emptyset
$$
(i.e. one, two,..., six of these six monomials for a given $p$ lie in $\AA$) for at most two values of $p$.\\ 

{\bf Lemma 2.} {\it In order to be smooth, $\AA$ should contain no singular pairs and no singular triples.}\\

{\it Proof:} Let $F=F(x_0,x_1,x_2,x_3,x_4,x_5)$ be a cubic form which contains a singular pair $(x_0, x_1)$. This means that $F$ can be written as 
$$
F=x_p\cdot P(x_0,x_1)+x_0\cdot Q_0(x_2,x_3,x_4,x_5)+x_1\cdot Q_1(x_2,x_3,x_4,x_5)+ R(x_2,x_3,x_4,x_5)
$$
for some quadratic forms $P, Q_0, Q_1$ and a cubic form $R$ and for some index $p\geq 2$.\\

By the Jacobian criterion its vanishing locus $X=V(F)\subset {\mathbb P}^{5}$ will be singular at a point with homogeneous coordinates $(a:b:0:0:0:0)$, where $P(a,b)=0$.\\ 

Let $F=F(x_0,x_1,x_2,x_3,x_4,x_5)$ be a cubic form which contains a singular triple $(x_0$, $x_1$, $x_2)$. This means that $F$ can be written as 
\begin{multline*}
F=x_p\cdot P_1(x_0,x_1,x_2)+x_q\cdot P_2(x_0,x_1,x_2)+x_0\cdot Q_0(x_3,x_4,x_5)+\\
+x_1\cdot Q_1(x_3,x_4,x_5)+x_2\cdot Q_2(x_3,x_4,x_5)+ R(x_3,x_4,x_5)
\end{multline*}
for some quadratic forms $P_1, P_2, Q_0, Q_1, Q_2$ and a cubic form $R$ and for some indices $p>q\geq 3$.\\

By the Jacobian criterion its vanishing locus $X=V(F)\subset {\mathbb P}^{5}$ will be singular at a point with homogeneous coordinates $(a:b:c:0:0:0)$, where $P_1(a,b,c)=P_2(a,b,c)=0$. {\it QED}\\

In the following subsections we will enumerate all minimal subsets $\AA\subset \mathcal M$ which satisfy the necessary conditions of Lemma 1 and Lemma 2. It will turn out (and is checked in the Appendix) that they all are smooth, i.e. in the case of cubic fourfolds conditions of Lemma~1 and Lemma~2 together are also {\it sufficient } for smoothness of $\AA$.\\

Smoothness of any given set of monomials $\AA\subset \mathcal M$ can be checked (at least in principle) by a discriminant computation. Practically, one can use Macaulay 2 (or a direct computation with the Jacobian criterion).\\

We used the following Macaulay 2 commands in order to verify smoothness of a cubic fourfold $X$ given as the vanishing locus of a cubic form $F$:\\ 

\begin{itemize}
\item $R=QQ[x_0,x_1,x_2,x_3,x_4,x_5]$
\item $F=x_0^3+x_1^3+x_2^2x_3+...$
\item $dim\; singularLocus(ideal(F))$
\end{itemize}

When the result is $0$, one can conclude that $X$ is indeed smooth.\\

Note that in order to prove Theorem $1$ we do not need to verify smoothness of {\it all} our minimal sets $\AA$. It is sufficient to do this only for sets $\AAbar$ corresponding to maximal groups $\GG$. This is done in the proof of Corollary $1$. The complete verification that all minimal sets $\AA$ we obtain are smooth is postponed until the Appendix, where one can find the Macaulay~2 code we used.\\

Finally, note that appearance of several smooth subsets $\AA\subset \mathcal M$ which are {\it not} minimal is harmless to our analysis. For this reason, we will not be concerned with checking minimality of our sets $\AA$. It will be sufficient for us to know that all minimal smooth sets $\AA$ are contained {\it among} those which we list. \\

Now let us do the enumeration.\\

{\it Proof of Theorem $1$:}\\

\subsection{Case of $6$ cubes.}

There is only one minimal smooth set in this group:
\begin{center}
\begin{tikzpicture}

\coordinate (0) at (0,0); \node[below, font=\scriptsize] at (0) {0}; \fill (0) circle (1pt); 
\coordinate (1) at (1,0); \node[below, font=\scriptsize] at (1) {1}; \fill (1) circle (1pt); 
\coordinate (2) at (2,0); \node[below, font=\scriptsize] at (2) {2}; \fill (2) circle (1pt); 
\coordinate (3) at (3,0); \node[below, font=\scriptsize] at (3) {3}; \fill (3) circle (1pt); 
\coordinate (4) at (4,0); \node[below, font=\scriptsize] at (4) {4}; \fill (4) circle (1pt);
\coordinate (5) at (5,0); \node[below, font=\scriptsize] at (5) {5}; \fill (5) circle (1pt);

\draw (0) circle (1mm); \draw (1) circle (1mm); \draw (2) circle (1mm); \draw (3) circle (1mm); \draw (4) circle (1mm); \draw (5) circle (1mm);

\end{tikzpicture}
\end{center}

In other words:
$$
\AA =\{ x_0^3, x_1^3, x_2^3, x_3^3, x_4^3, x_5^3 \}.
$$

Elements $f=(1, {\om}^{c_1}, {\om}^{c_2}, {\om}^{c_3}, {\om}^{c_4}, {\om}^{c_5})\in G_{\AA}$, $\om=\sqrt[d]{1}$, $d=ord(f)$ are given by the conditions:
$$
3c_1=3c_2=3c_3=3c_4=3c_5=0 \;\;(\mbox{in}\; \mathbb Z / d \mathbb Z).
$$

This means that $d=3$ and $\GG\cong (\ZZ{3})^{\oplus 5}$ with generators
$$
f_1=(1,\om,1,1,1,1),\; f_2=(1,1,\om,1,1,1),\; f_3=(1,1,1,\om,1,1),
$$

$$
f_4=(1,1,1,1,\om,1),\;\;\; f_5=(1,1,1,1,1,\om), 
$$
where $\om=\sqr{3}$.\\

Then $\AAbar$ consists of monomials $x_0^{i_0}x_1^{i_1}x_2^{i_2}x_3^{i_3}x_4^{i_4}x_5^{i_5}$ such that 
$$
i_1 \equiv i_2 \equiv i_3 \equiv i_4 \equiv i_5 \equiv 0 \; mod \; 3. 
$$

Since $i_0+i_1 + i_2 + i_3 + i_4 + i_5=3$ and each $i_p\geq 0$ we conclude that $\AAbar=\AA$.\\

\subsection{Case of $5$ cubes.}

There is only one (upto a permutation of $x_i$) such minimal smooth set of monomials:
\begin{center}
\begin{tikzpicture}

\coordinate (5) at (0,0); \node[below, font=\scriptsize] at (5) {5}; \fill (5) circle (1pt); 
\coordinate (0) at (1,0); \node[below, font=\scriptsize] at (0) {0}; \fill (0) circle (1pt); 

\draw (0) circle (1mm);

\draw (5)--(0);

\end{tikzpicture}
\end{center}

In other words:
$$
\AA =\{ x_0^3, x_1^3, x_2^3, x_3^3, x_4^3, x_5^2x_0 \}.
$$

Conditions on $f=(1, {\om}^{c_1}, {\om}^{c_2}, {\om}^{c_3}, {\om}^{c_4}, {\om}^{c_5})\in \GG$, $\om=\sqrt[d]{1}$, $d=ord(f)$ are as follows:
$$
3c_1=3c_2=3c_3=3c_4=0,\; 2c_5=0.
$$

Hence $d=6$ and $f=(1, {\om}^{2c_1}, {\om}^{2c_2}, {\om}^{2c_3}, {\om}^{2c_4}, {\om}^{3c_5})$, $c_1,c_2,c_3,c_4\in \ZZ{3}$, $c_5\in \ZZ{2}$.\\

This means that $\GG\cong (\ZZ{3})^{\oplus 4}\oplus \ZZ{2}$ with generators
$$
f_1=(1,\om,1,1,1,1),\; f_2=(1,1,\om,1,1,1),\; f_3=(1,1,1,\om,1,1),
$$

$$
f_4=(1,1,1,1,\om,1),\;\;\; f_5=(1,1,1,1,1,\eta),
$$
where $\om=\sqr{3}$, $\eta=\sqr{2}$.\\

Then $\AAbar$ consists of monomials $x_0^{i_0}x_1^{i_1}x_2^{i_2}x_3^{i_3}x_4^{i_4}x_5^{i_5}$ such that 
$$
i_1 \equiv i_2 \equiv i_3 \equiv i_4 \equiv 0 \; mod \; 3, \;\;\; i_5\equiv 0 \; mod \; 2. 
$$

This means that $\AAbar=\AA$.\\

\subsection{Case of $4$ cubes.}

Let the cubes be $x_0^3, x_1^3, x_2^3, x_3^3$. The other variables $x_4$ and $x_5$ may form a cycle:
\begin{center}

\end{center}

or (the corresponding vertices) may be connected by a single edge:
\begin{center}
%
\end{center}

or they may be disconnected at all:
\begin{center}
%
\end{center}

If $\AA =\{ x_4^2x_5, x_5^2x_4, x_0^3, x_1^3, x_2^3, x_3^3 \}$, then $f=(1, {\om}^{c_1}, {\om}^{c_2}, {\om}^{c_3}, {\om}^{c_4}, {\om}^{c_5})$, $\om=\sqrt[d]{1}$, $d=ord(f)$, $c_i\in \ZZ{d}$ is in $\GG$ if and only if
$$
3c_1=3c_2=3c_3=2c_4+c_5=c_4+2c_5=0,
$$
i.e. $d=3$ and $f=(1, {\om}^{c_1}, {\om}^{c_2}, {\om}^{c_3}, {\om}^{c_4}, {\om}^{c_4})$, $c_i\in \ZZ{3}$.\\

This means that $\GG\cong (\ZZ{3})^{\oplus 4}$ with generators
$$
f_1=(1,\om,1,1,1,1),\; f_2=(1,1,\om,1,1,1),\; f_3=(1,1,1,\om,1,1),\; f_4=(1,1,1,1,\om,\om),
$$
where $\om=\sqr{3}$.\\

Then $\AAbar$ consists of monomials $x_0^{i_0}x_1^{i_1}x_2^{i_2}x_3^{i_3}x_4^{i_4}x_5^{i_5}$ such that 
$$
i_1 \equiv i_2 \equiv i_3 \equiv i_4 +i_5\equiv  0 \; mod \; 3. 
$$

This means that $\AAbar=\AA \cup \{ f_3(x_4, x_5) \}$.\\

If $\AA =\{ x_4^2x_5, x_5^2x_0, x_0^3, x_1^3, x_2^3, x_3^3 \}$, then $f=(1, {\om}^{c_1}, {\om}^{c_2}, {\om}^{c_3}, {\om}^{c_4}, {\om}^{c_5})$, $\om=\sqrt[d]{1}$, $d=ord(f)$, $c_i\in \ZZ{d}$ is in $\GG$ if and only if
$$
3c_1=3c_2=3c_3=2c_5=2c_4+c_5=0 \;\; \mbox{in}\; \ZZ{d},
$$
i.e. $f=(1, {\om}^{c_1}, {\om}^{c_2}, {\om}^{c_3}, {\eta}^{c_4}, {\eta}^{2c_4})$, where $\om=\sqr{3}$, $\eta=\sqr{4}$, $c_1,c_2,c_3\in \ZZ{3}$, $c_4\in \ZZ{4}$.\\

This means that $\GG\cong (\ZZ{3})^{\oplus 3}\oplus \ZZ{4}$ with generators
$$
f_1=(1,\om,1,1,1,1),\; f_2=(1,1,\om,1,1,1),\; f_3=(1,1,1,\om,1,1),\; f_4=(1,1,1,1,\eta,{\eta}^2).
$$

Then $\AAbar$ consists of monomials $x_0^{i_0}x_1^{i_1}x_2^{i_2}x_3^{i_3}x_4^{i_4}x_5^{i_5}$ such that 
$$
i_1 \equiv i_2 \equiv i_3 \equiv 0 \; mod \; 3, \;\; i_4+2i_5 \equiv 0 \; mod \; 4. 
$$

This means that $\AAbar=\AA$.\\

If $\AA =\{ x_4^2x_0, x_5^2x_1, x_0^3, x_1^3, x_2^3, x_3^3 \}$, then $f=(1, {\om}^{c_1}, {\om}^{c_2}, {\om}^{c_3}, {\om}^{c_4}, {\om}^{c_5})$, $\om=\sqrt[d]{1}$, $d=ord(f)$, $c_i\in \ZZ{d}$ is in $\GG$ if and only if
$$
3c_1=3c_2=3c_3=2c_4=2c_5+c_1=0 \;\; \mbox{in}\; \ZZ{d},
$$
i.e. $d=6$ and $f=(1, {\om}^{-2c_5}, {\om}^{2c_2}, {\om}^{2c_3}, {\om}^{3c_4}, {\om}^{c_5})$, where $\om=\sqr{6}$, $c_2,c_3\in \ZZ{3}$, $c_4\in \ZZ{2}$, $c_5\in \ZZ{6}$.\\

This means that $\GG\cong (\ZZ{3})^{\oplus 2}\oplus \ZZ{2}\oplus \ZZ{6}$ with generators
$$
f_1=(1,{\om}^{-2},1,1,1,\om),\; f_2=(1,1,{\om}^2,1,1,1),\; f_3=(1,1,1,{\om}^2,1,1),\; f_4=(1,1,1,1,{\om}^3,1).
$$

Then $\AAbar$ consists of monomials $x_0^{i_0}x_1^{i_1}x_2^{i_2}x_3^{i_3}x_4^{i_4}x_5^{i_5}$ such that 
$$
i_5\equiv 2i_1 \; mod \; 6, \;\; i_2 \equiv i_3 \equiv 0 \; mod \; 3, \;\; i_4\equiv 0 \; mod \; 2. 
$$

This means that $\AAbar=\AA$.\\

From now on (and until the end of the proof) we will follow exactly the same scheme of analysis, but we will abbreviate it.\\

If $\AA =\{ x_4^2x_0, x_5^2x_0, x_0^3, x_1^3, x_2^3, x_3^3, x_1x_4x_5 \}$, then we get the following conditions on elements $c_i$ of the signature of $f\in \GG$ as above:
$$
3c_1=3c_2=3c_3=c_1+c_4+c_5=2c_4=2c_5=0 \;\; \mbox{in}\; \ZZ{d},
$$
i.e. $d=6$, $c_1=0$ and $f=(1, 1, {\om}^{2c_2}, {\om}^{2c_3}, {\om}^{3c_4}, {\om}^{3c_4})$, where $\om=\sqr{6}$, $c_2, c_3\in \ZZ{3}$, $c_4\in \ZZ{2}$.\\

This means that $\GG\cong (\ZZ{3})^{\oplus 2}\oplus \ZZ{2}$ with generators
$$
f_1=(1,1,{\om}^2,1,1,1),\; f_2=(1,1,1,{\om}^2,1,1),\; f_3=(1,1,1,1,{\om}^3,{\om}^3).
$$

Then we get the following equations describing monomials in $\AAbar$ as above:
$$
2 i_2 \equiv 2i_3 \equiv 3(i_4+i_5)\equiv 0 \; mod \; 6. 
$$

This means that $\AAbar=\AA \cup \{ f_3(x_0,x_1), f_2(x_4,x_5)\cdot f_1(x_0,x_1) \}$.\\

\subsection{Case of $3$ cubes.}

Let the cubes be $x_0^3, x_1^3, x_2^3$. If a longest cycle has length $3$, then $\AA$ is
\begin{center}
\begin{tikzpicture}

\coordinate (0) at (0,0); \node[below, font=\scriptsize] at (0) {0}; \fill (0) circle (1pt); 
\coordinate (1) at (1,0); \node[below, font=\scriptsize] at (1) {1}; \fill (1) circle (1pt);
\coordinate (2) at (2,0); \node[below, font=\scriptsize] at (2) {2}; \fill (2) circle (1pt);
\coordinate (3) at (3,0); \node[below, font=\scriptsize] at (3) {3}; \fill (3) circle (1pt);
\coordinate (4) at (4,1); \node[right, font=\scriptsize] at (4) {4}; \fill (4) circle (1pt);
\coordinate (5) at (4,-1); \node[right, font=\scriptsize] at (5) {5}; \fill (5) circle (1pt);

\arr{3}{4}; \arr{4}{5}; \arr{5}{3}; \draw (0) circle (1mm); \draw (1) circle (1mm);\draw (2) circle (1mm);

\end{tikzpicture}
\end{center}

i.e. $\AA=\{ x_0^3,x_1^3,x_2^3, x_3^2x_4, x_4^2x_5, x_5^2x_3 \}$. Then $3c_1=3c_2=2c_3+c_4=2c_4+c_5=2c_5+c_3=0$, i.e.
$$
c_4=-2c_3,\; c_5=4c_3,\; 9c_3=3c_1=3c_2=0 \;\; \mbox{in}\; \ZZ{d}.
$$

Hence $d=9$ and $f=(1, {\om}^{3c_1}, {\om}^{3c_2}, {\om}^{c_3}, {\om}^{-2c_3}, {\om}^{4c_3})$, where $\om=\sqr{9}$, $c_1,c_2\in \ZZ{3}$, $c_3\in \ZZ{9}$.\\

This means that $\GG\cong (\ZZ{3})^{\oplus 2}\oplus \ZZ{9}$ with generators 
$$
f_1=(1,{\om}^3,1,1,1,1),\; f_2=(1,1,{\om}^3,1,1,1),\; f_3=(1,1,1,\om, {\om}^7,{\om}^4).
$$

Hence 
$$
3i_1\equiv 3i_2\equiv i_3-2i_4+4i_5 \equiv 0 \; mod \; 9.
$$

This means that $\AAbar=\AA$.\\

\subsubsection{Case of $3$ cubes. Length $2$ longest cycle.}

If a longest cycle has length $2$, then $\AA$ is either

\begin{center}

\end{center}

Indeed, let $x_3, x_4$ form a cycle. Then by Lemma 1 (\cite{Liendo}, Lemma 1.3) the $5$-th vertex should be connected either to one of the cubes (say, to the $0$-th vertex) or to the cycle (say, to the $3$-rd vertex).\\

In the former case there are neither singular pairs nor singular triples. In the latter case $(x_4, x_5)$ may be a singular pair. In order to avoid this, either vertices $4$ and $5$ should be connected by an edge (which means that $5$ should be connected to $4$, since we assume that there are no cycles of length $3$), or one of them should be connected to a cube, or there is a dashed curve passing through both of them and a cube.\\

Note that if $5$ is connected to a cube, this will lead to a nonminimal $\AA$. Hence only the $4$-th vertex may be connected to a cube (say, to $0$).\\

If $\AA =\{ x_3^2x_4, x_4^2x_3, x_5^2x_0, x_0^3, x_1^3, x_2^3 \}$, then $3c_1=3c_2=2c_3+c_4=2c_4+c_3=2c_5=0$, i.e.
$$
c_4=c_3,\; 3c_1=3c_2=3c_3=2c_5=0 \;\; \mbox{in}\; \ZZ{d}.
$$
Hence $d=6$ and $f=(1, {\om}^{2c_1}, {\om}^{2c_2}, {\om}^{2c_3}, {\om}^{2c_3}, {\om}^{3c_5})$, where $\om=\sqr{6}$, $c_1,c_2,c_3\in \ZZ{3}$, $c_5\in \ZZ{2}$.\\

This means that $\GG\cong (\ZZ{3})^{\oplus 3}\oplus \ZZ{2}$ with generators
$$
f_1=(1,{\om}^2,1,1,1,1),\; f_2=(1,1,{\om}^2,1,1,1),\; f_3=(1,1,1,{\om}^2,{\om}^2,1),\; f_4=(1,1,1,1,1,{\om}^3).
$$

Then
$$
2 i_1\equiv 2i_2 \equiv 2(i_3+i_4) \equiv 3i_5\equiv 0 \; mod \; 6. 
$$

This means that $\AAbar=\AA\cup \{ f_3(x_3,x_4) \}$.\\

If $\AA =\{ x_3^2x_4, x_4^2x_3, x_5^2x_3, x_4^2x_0, x_0^3, x_1^3, x_2^3 \}$, then $3c_1=3c_2=3c_3=2c_5+c_3=2c_4=0$, $c_3=c_4$
$$
3c_1=3c_2=c_3=c_4=2c_5=0 \;\; \mbox{in}\; \ZZ{d},
$$
i.e. $d=6$ and $f=(1, {\om}^{2c_1}, {\om}^{2c_2}, 1, 1, {\om}^{3c_5})$, where $\om=\sqr{6}$, $c_1,c_2\in \ZZ{3}$, $c_5\in \ZZ{2}$.\\

This means that $\GG\cong (\ZZ{3})^{\oplus 2}\oplus \ZZ{2}$ with generators
$$
f_1=(1,{\om}^2,1,1,1,1),\; f_2=(1,1,{\om}^2,1,1,1),\; f_3=(1,1,1,1,1,{\om}^3).
$$

Hence
$$
2 i_1\equiv 2i_2 \equiv 3i_5\equiv 0 \; mod \; 6. 
$$

This means that $\AAbar=\AA \cup \{  f_3(x_0,x_3,x_4), x_5^2\cdot f_1(x_0,x_3,x_4) \}$.\\

If $\AA =\{ x_3^2x_4, x_4^2x_3, x_5^2x_3, x_0x_4x_5, x_0^3, x_1^3, x_2^3 \}$, then $3c_1=3c_2=3c_3=0$, $c_3=c_4$, $c_3=-2c_5$, $c_4+c_5=0$, i.e.
$$
3c_1=3c_2=0, \; c_3=c_4=c_5=0 \;\; \mbox{in}\; \ZZ{d},
$$
Hence $d=3$ and $f=(1, {\om}^{c_1}, {\om}^{c_2}, 1, 1, 1)$, where $\om=\sqr{3}$, $c_i\in \ZZ{3}$.\\

This means that $\GG\cong (\ZZ{3})^{\oplus 2}$ with generators
$$
f_1=(1,{\om},1,1,1,1),\; f_2=(1,1,{\om},1,1,1).
$$

Hence
$$
i_1\equiv i_2 \equiv 0 \; mod \; 3. 
$$

This means that $\AAbar=\AA \cup \{  f_3(x_0,x_3,x_4, x_5) \}$.\\

If $\AA =\{ x_3^2x_4, x_4^2x_3, x_5^2x_3, x_5^2x_4, x_0^3, x_1^3, x_2^3 \}$, then $3c_1=3c_2=3c_3=0$, $c_3=c_4=-2c_5$, i.e.
$$
3c_1=3c_2=6c_5=0, \; c_3=c_4=-2c_5 \;\; \mbox{in}\; \ZZ{d},
$$
Hence $d=6$ and $f=(1, {\om}^{2c_1}, {\om}^{2c_2}, {\om}^{-2c_5}, {\om}^{-2c_5}, {\om}^{c_5})$, where $\om=\sqr{6}$, $c_1, c_2\in \ZZ{3}$, $c_5\in \ZZ{5}$.\\

This means that $\GG\cong (\ZZ{3})^{\oplus 2}\oplus \ZZ{6}$ with generators
$$
f_1=(1,{\om}^2,1,1,1,1),\; f_2=(1,1,{\om}^2,1,1,1),\; f_3=(1,1,1,{\om}^{-2},{\om}^{-2},\om).
$$

Hence
$$
i_1\equiv i_2 \equiv 0 \; mod \; 3,\; i_5\equiv 2(i_3+i_4) \; mod \; 6. 
$$

This means that $\AAbar=\AA \cup \{  f_3(x_0,x_3,x_4, x_5) \}$.\\

\subsubsection{Case of $3$ cubes. No cycles.}

If $\AA$ has no cycles, then the following graphs are possible: 
\begin{center}

\end{center}

Indeed, consider the subgraph formed by vertices $x_3, x_4, x_5$ which are not cubes. Let us consider the length of a longest path in this subgraph.\\

If this subgraph has a path of length $3$, then we may assume that $x_3$ is connected to $x_4$, $x_4$ is connected to $x_5$. Then by Lemma 1 (\cite{Liendo}, Lemma 1.3) $x_5$  should be connected to a cube (say, to $0$) since there are no cycles.\\ 

If the length of a longest path of this subgraph is $2$, then we may assume that $3$ is connected to $4$, which is then connected to one of the cubes (say, to $0$), since there are neither cycles nor paths of length $3$ by our assumptions. By Lemma 1 (\cite{Liendo}, Lemma 1.3) the remaining vertex $5$ should be connected to either a cube or to the length $2$ path formed by vertices $3$ and $4$.\\

If $5$ is connected to a cube, then this cube may be either the same cube to which the length $2$ path was connected (in which case $(x_4,x_5)$ may form a singular pair) or a different one (in which case there are neither singular pairs nor singular triples). In order to avoid the singular pair $(x_4,x_5)$, we should have that either
\begin{itemize}
\item there is an edge connecting $4$ and $5$, which in fact means that $5$ is connected to $4$ (rather than $4$ to $5$) since there are not length $3$ paths, or
\item one of the vertices $4$, $5$ is connected to either a different cube ($1$ or $2$) or to the vertex $3$ (the latter is not possible since there are neither cycles nor length $3$ paths), or
\item there is a dashed curve connecting $4$, $5$ and a different cube ($1$ or $2$) or the $3$-th vertex.
\end{itemize}

The second option here ($4$ or $5$ is connected to another cube) does not occur, since it would lead to a nonminimal $\AA$.\\

If $5$ is connected to the length $2$ path formed by vertices $3$ and $4$, then in fact it is connected to the vertex $4$ (since there are no paths of length $3$). In this case $(x_3, x_5)$ may form a singular pair. In order to avoid this, we should have that either 
\begin{itemize}
\item there is an edge between $3$ and $5$, which is impossible since there are no paths of length $3$, or
\item one of the vertices $3$, $5$ is connected to a cube, which would lead either to a nonminimal $\AA$ (in case the cube is not $0$) or to one of the already analyzed cases (if the cube is $0$), or
\item there is a dashed curve passing through $3$, $5$ and one of the cubes.
\end{itemize}

Finally, it is possible that the subgraph formed by vertices $3$, $4$, $5$ is totally disconnected. In this case by Lemma 1 (\cite{Liendo}, Lemma 1.3) each of its vertices should be connected to one of the cubes. The cubes may be all different, or only two of them may be different, or all the vertices $3$, $4$, $5$ may be connected to the single cube (say, to $0$).\\

In case $3$, $4$, $5$ are connected to $3$ different cubes (say, $3$ is connected to $0$, $4$ is connected to $1$ and $5$ is connected to $2$), there are neither singular pairs nor singular triples.\\

Suppose vertices $3$ and $4$ are connected to the same cube (say, to $0$) and vertex $5$ is connected to a different cube (say, to $1$). Then $(x_3,x_4)$ may form a singular pair and $(x_3,x_4,x_5)$ may form a singular triple.\\

In order to avoid the singular pair, there should be a dashed curve connecting $3$, $4$ together with either a cube ($1$ or $2$) or the $5$-th vertex. Alternatively, $4$ (or $3$) may be connected to the same cube $1$ to which vertex $5$ is connected. All other options are either proscribed by our assumptions (when $3$, $4$ are connected by an edge to each other or to the vertex $5$) or lead to a nonminimal $\AA$ (when $3$ or $4$ is connected by an edge to the cube $2$).\\

If the dashed curve involves either the cube $2$ or the vertex $5$, this automatically resolves the singular triple $(x_3,x_4,x_5)$ as well.\\

In the case, when the dashed curve passes through the cube $1$, there should be also a dashed curve connecting vertex $5$, one of the vertices $3$, $4$ and the third cube $2$ in order to avoid the singular triple $(x_3,x_4,x_5)$. By the symmetry, we may assume that the dashed curve passes through $2$, $4$ and $5$. All other options would either lead to nonminimal $\AA$ or have already appeared.\\

Indeed, there are no edges between $3$, $4$, $5$ by our assumption. If we introduce a dashed curve through $3$, $4$, $5$, then this will lead to a nonminimal $\AA$, because the dashed curve connecting $3$, $4$ to $5$ has appeared earlier and automatically resolved both the singular pair and the singular triple. Hence in order to resolve the singular triple $(x_3,x_4,x_5)$ we should have either an edge from one of $3$, $4$, $5$ to the third cube $2$ or a dashed curve passing through two of $3$, $4$, $5$ and the cube $2$. An edge from $3$ or $4$ to $2$ will contradict to our assumption that $3$, $4$, $5$ are connected to only two different cubes. If we have an edge from $5$ to $2$ or a dashed curve through $3$, $4$ and $2$, then this will lead us to a nonminimal $\AA$, because such configurations (upto a permutation of $x_i$) have already appeared above. Hence only a dashed curve through $2$, $5$ and one of the vertices $3$, $4$ is possible.\\

This resolves the singular triple $(x_3,x_4,x_5)$, when there is a dashed curve through $3$, $4$ and $1$.\\

Suppose that $4$ is connected to $1$ by an edge. Then again we should have a dashed curve passing though $2$, $3$, $5$ in order to resolve the singular triple $(x_3,x_4,x_5)$.\\

Indeed, there are neither edges between vertices $3$, $4$, $5$ nor the dashed curve connecting them (such a curve would lead to a nonminimal $\AA$ as above). Hence there will be a singular triple $(x_3,x_4,x_5)$, unless either
\begin{itemize}
\item one of the vertices $3$, $4$, $5$ is connected to the cube $2$ (which is impossible by our assumption), or
\item there is a dashed curve through a pair of vertices $3$, $4$, $5$ and the cube $2$.
\end{itemize}

In the second case the dashed curve should pass through $3$, $5$ and $2$, because if it passes through, say, $3$, $4$ and $2$, we will get a nonminimal $\AA$ (such a dashed curve has already appeared above and resolved both the singular pair and the singular triple).\\

The remaining possibility is when $3$, $4$ and $5$ are all connected to one and the same cube (say, to $0$). In this case we may have three singular pairs $(x_3,x_4)$, $(x_4,x_5)$, $(x_5,x_3)$ and one singular triple $(x_3,x_4,x_5)$.\\

A dashed curve through vertices $3$, $4$, $5$ resolves all singular pairs and all singular triples.\\

Suppose that $3$, $4$, $5$ are not connected by such dashed curve. Then by our assumptions each singular pair is resolved by a dashed curve through the corresponding pair of vertices and one of the cubes $1$, $2$ (since no edges are allowed).\\

In fact, each of these cubes $1$, $2$ should lie on a dashed curve passing through a pair of $3$, $4$, $5$, since otherwise $(x_3, x_4, x_5)$ would form a singular triple. Without loss of generality (since we can always permute $x_i$) we may assume that there are dashed curves through $1$, $3$, $4$, through $1$, $4$, $5$ and through $2$, $3$, $5$.\\

This argument gives us the pictures shown above. Now let us do computations.\\

If $\AA =\{ x_3^2x_4, x_4^2x_5, x_5^2x_0, x_0^3, x_1^3, x_2^3 \}$, then $3c_1=3c_2=2c_5=2c_3+c_4=2c_4+c_5=0$, i.e.
$$
c_4=-2c_3,\; c_5=4c_3,\; 3c_1=3c_2=8c_3=0 \;\; \mbox{in}\; \ZZ{d}.
$$

Hence $d=24$ and $f=(1, {\om}^{8c_1}, {\om}^{8c_2}, {\om}^{3c_3}, {\om}^{-6c_3}, {\om}^{12c_3})$, where $\om=\sqr{24}$, $c_1,c_2\in \ZZ{3}$, $c_3\in \ZZ{8}$.\\

This means that $\GG\cong (\ZZ{3})^{\oplus 2}\oplus \ZZ{8}$ with generators
$$
f_1=(1,{\om}^8,1,1,1,1),\; f_2=(1,1,{\om}^8,1,1,1),\; f_3=(1,1,1,{\om}^3,{\om}^{-6},{\om}^{12}).
$$

Hence $i_1\equiv i_2 \equiv 0 \; mod \; 3$ and $i_3-2i_4+4i_5\equiv 0 \; mod \; 8$. This means that $\AAbar=\AA$.\\

If $\AA =\{ x_3^2x_4, x_4^2x_0, x_5^2x_1, x_0^3, x_1^3, x_2^3 \}$, then $3c_1=3c_2=c_1+2c_5=2c_4=2c_3+c_4=0$, i.e.
$$
c_4=2c_3,\; c_1=-2c_5,\; 3c_2=4c_3=6c_5=0 \;\; \mbox{in}\; \ZZ{d}.
$$

Hence $d=12$ and $f=(1, {\om}^{-4c_5}, {\om}^{4c_2}, {\om}^{3c_3}, {\om}^{6c_3}, {\om}^{2c_5})$, where $\om=\sqr{12}$, $c_2\in \ZZ{3}$, $c_3\in \ZZ{4}$, $c_5\in \ZZ{6}$.\\

This means that $\GG\cong \ZZ{3} \oplus \ZZ{4}\oplus \ZZ{6}$ with generators
$$
f_1=(1,{\om}^{-4},1,1,1,{\om}^2),\; f_2=(1,1,{\om}^4,1,1,1),\; f_3=(1,1,1,{\om}^3,{\om}^{6},1).
$$

Hence $i_5\equiv 2i_1 \; mod \; 6$, $i_2 \equiv 0 \; mod \; 3$ and $i_3\equiv 2i_4 \; mod \; 4$. This means that $\AAbar=\AA$.\\

Let $\AAN{0} =\{ x_3^2x_4, x_4^2x_0, x_5^2x_0, x_0^3, x_1^3, x_2^3 \}$. Then $3c_1=3c_2=2c_3+c_4=2c_4=2c_5=0$, i.e.
$$
c_4=2c_3,\; 3c_1=3c_2=4c_3=2c_5=0 \;\; \mbox{in}\; \ZZ{d}.
$$
Hence $d=12$ and $f=(1, {\om}^{4c_1}, {\om}^{4c_2}, {\om}^{3c_3}, {\om}^{6c_3}, {\om}^{6c_5})$, where $\om=\sqr{12}$, $c_1,c_2\in \ZZ{3}$, $c_3\in \ZZ{4}$, $c_5\in \ZZ{2}$.\\

For $\AA=\AAN{0}\cup \{ x_3x_4x_5 \}$ we get an extra condition $c_3+c_4+c_5=0$, i.e. 
$$
c_5=c_3,\; 3c_1=3c_2=2c_3=c_4=0 \;\; \mbox{in}\; \ZZ{d}.
$$
Hence $d=6$ and $f=(1, {\om}^{2c_1}, {\om}^{2c_2}, {\om}^{3c_3}, 1, {\om}^{3c_3})$, where $\om=\sqr{6}$, $c_1,c_2\in \ZZ{3}$, $c_3\in \ZZ{2}$.\\

This means that $\GG\cong (\ZZ{3})^{\oplus 2} \oplus \ZZ{2}$ with generators
$$
f_1=(1,{\om}^{2},1,1,1,1),\; f_2=(1,1,{\om}^2,1,1,1),\; f_3=(1,1,1,{\om}^3,1,{\om}^{3}).
$$

Then $i_1\equiv i_2\equiv 0 \; mod \; 3$ and $i_3+i_5 \equiv 0 \; mod \; 2$. This means that $\AAbar=\AA\cup \{ f_3(x_0,x_4), f_1(x_0,x_4)\cdot f_2(x_3,x_5) \}$.\\

For $\AA=\AAN{0}\cup \{ x_1x_4x_5 \}$ we get an extra condition $c_1+c_4+c_5=0$, i.e.
$$
c_4=c_5=2c_3,\; c_1=3c_2=4c_3=0  \;\; \mbox{in}\; \ZZ{d}.
$$
Hence $d=12$, $f=(1, 1, {\om}^{4c_2}, {\om}^{3c_3}, {\om}^{6c_3}, {\om}^{6c_3})$, where $\om=\sqr{12}$, $c_4\in \ZZ{3}$, $c_3\in \ZZ{4}$.\\

This means that $\GG\cong \ZZ{3} \oplus \ZZ{4}$ with generators
$$
f_1=(1,1,{\om}^{4},1,1,1),\; f_2=(1,1,1,{\om}^3,{\om}^6,{\om}^6).
$$

Then $i_2\equiv 0 \; mod \; 3$ and $i_3\equiv 2(i_4+i_5) \; mod \; 4$. This means that $\AAbar=\AA\cup \{ f_3(x_0,x_1)$, $f_1(x_0,x_1)\cdot f_2(x_4,x_5)$, $x_3^2\cdot f_1(x_4,x_5) \}$.\\

For $\AA=\AAN{0}\cup \{ x_5^2x_4 \}$ we get an extra condition $c_4+2c_5=0$, i.e.
$$
c_4=0,\; 3c_1=3c_2=2c_3=2c_5=0  \;\; \mbox{in}\; \ZZ{d}.
$$
Hence $d=6$, $f=(1, {\om}^{2c_1}, {\om}^{2c_2}, {\om}^{3c_3}, 1, {\om}^{3c_5})$, where $\om=\sqr{6}$, $c_1,c_2\in \ZZ{3}$, $c_3,c_5\in \ZZ{2}$.\\

This means that $\GG\cong (\ZZ{3})^{\oplus 2} \oplus (\ZZ{2})^{\oplus 2}$ with generators
$$
f_1=(1,{\om}^{2},1,1,1,1), \; f_2=(1,1,{\om}^{2},1,1,1),\; f_3=(1,1,1,{\om}^3,1,1),\; f_4=(1,1,1,1,1,{\om}^3).
$$

Then $i_1\equiv i_2\equiv 0 \; mod \; 3$ and $i_3 \equiv i_5\equiv 0 \; mod \; 2$. This means that $\AAbar=\AA\cup \{ x_0^2x_4, x_4^3, x_3^2x_0 \}$.\\

Let $\AAN{0} =\{ x_3^2x_4, x_5^2x_4, x_4^2x_0, x_0^3, x_1^3, x_2^3 \}$. Then $3c_1=3c_2=2c_3+c_4=2c_4=2c_5+c_4=0$, i.e.
$$
c_4=2c_3,\; 3c_1=3c_2=2(c_5-c_3)=4c_3=0 \;\; \mbox{in}\; \ZZ{d}.
$$
Hence $d=12$ and $f=(1, {\om}^{4c_1}, {\om}^{4c_2}, {\om}^{3c_3}, {\om}^{6c_3}, {\om}^{3c_3+6c_5})$, where $\om=\sqr{12}$, $c_1,c_2\in \ZZ{3}$, $c_3\in \ZZ{4}$, $c_5\in \ZZ{2}$.\\

For $\AA=\AAN{0}\cup \{ x_1x_3x_5 \}$ we get an extra condition $c_1+c_3+c_5=0$, i.e.
$$
c_4=2c_3,\; c_5=-c_3,\; c_1=3c_2=4c_3=0\;\; \mbox{in}\; \ZZ{d}.
$$
Then $d=12$, $f=(1, 1, {\om}^{4c_2}, {\om}^{3c_3}, {\om}^{6c_3}, {\om}^{-3c_3})$, where $\om=\sqr{12}$, $c_2\in \ZZ{3}$, $c_3\in \ZZ{4}$.\\

This means that $\GG\cong \ZZ{3} \oplus \ZZ{4}$ with generators
$$
f_1=(1,1,{\om}^{4},1,1,1),\; f_2=(1,1,1,{\om}^3,{\om}^6,{\om}^{-3}).
$$

Hence $i_2\equiv 0 \; mod \; 3$ and $i_3+2i_4-i_5 \equiv 0 \; mod \; 4$. This means that $\AAbar=\AA\cup \{ x_4^2x_1, x_0x_3x_5 \}$.\\

For $\AA=\AAN{0}\cup \{ x_0x_3x_5 \}$ we get an extra condition $c_3+c_5=0$, i.e. 
$$
c_5=-c_3,\; c_4=2c_3,\; 3c_1=3c_2=4c_3=0.
$$

Then $f=(1, {\om}^{4c_1}, {\om}^{4c_2}, {\om}^{3c_3}, {\om}^{6c_3}, {\om}^{-3c_3})$, where $\om=\sqr{12}$, $c_1,c_2\in \ZZ{3}$, $c_3\in \ZZ{4}$.\\

This means that $\GG\cong (\ZZ{3})^{\oplus 2} \oplus \ZZ{4}$ with generators
$$
f_1=(1,{\om}^{4},1,1,1,1),\; f_2=(1,1,{\om}^{4},1,1,1),\; f_3=(1,1,1,{\om}^3,{\om}^6,{\om}^{-3}).
$$

Hence $i_1\equiv i_2\equiv 0 \; mod \; 3$, $i_3+2i_4-i_5 \equiv 0 \; mod \; 4$. This means that $\AAbar=\AA$.\\

If $\AA =\{ x_3^2x_0, x_4^2x_1, x_5^2x_2, x_0^3, x_1^3, x_2^3 \}$, then $3c_1=3c_2=2c_3=2c_4+c_1=2c_5+c_2=0$, i.e.
$$
c_1=-2c_4,\; c_2=-2c_5,\; 6c_4=6c_5=2c_3=0.
$$
Hence $d=6$ and $f=(1, {\om}^{-2c_4}, {\om}^{-2c_5}, {\om}^{3c_3}, {\om}^{c_4}, {\om}^{c_5})$, where $\om=\sqr{6}$, $c_3\in \ZZ{2}$, $c_4,c_5\in \ZZ{6}$.\\

This means that $\GG\cong \ZZ{2}\oplus (\ZZ{6} )^{\oplus 2}$ with generators
$$
f_1=(1,{\om}^{-2},1,1,{\om},1),\; f_2=(1,1,{\om}^{-2},1,1,\om),\; f_3=(1,1,1,{\om}^3,1,1).
$$

Then $i_4\equiv 2i_1 \; mod \; 6$, $i_5 \equiv 2i_2 \; mod \; 6$ and $i_3 \equiv 0 \; mod \; 2$. This means that $\AAbar=\AA$.\\

Let $\AAN{0} =\{ x_3^2x_0, x_4^2x_0, x_5^2x_1, x_0^3, x_1^3, x_2^3 \}$. Then $3c_1=3c_2=2c_3=2c_4=2c_5+c_1=0$, i.e. $c_1=-2c_5$, $2c_3=3c_2=6c_5=2c_4=0$. Hence $f=(1, {\om}^{-2c_5}, {\om}^{2c_2}, {\om}^{3c_3}, {\om}^{3c_4}, {\om}^{c_5})$, where $\om=\sqr{6}$, $c_i\in \ZZ{6}$.\\

For $\AA=\AAN{0}\cup \{ x_2x_3x_4 \}$ we get an extra condition $c_2+c_3+c_4=0$, i.e.
$$
c_2=0, \; c_4=c_3,\; c_1=-2c_5,\; 2c_3=6c_5=0.
$$
Hence $f=(1, {\om}^{-2c_5}, 1, {\om}^{3c_3}, {\om}^{3c_3}, {\om}^{c_5})$, where $\om=\sqr{6}$, $c_3\in \ZZ{2}$, $c_5\in \ZZ{6}$.\\

This means that $\GG\cong \ZZ{6} \oplus \ZZ{2}$ with generators
$$
f_1=(1,{\om}^{-2},1,1,1,\om),\; f_2=(1,1,1,{\om}^3,{\om}^3,1).
$$

Hence $i_3+i_4\equiv 0 \; mod \; 2$, $i_5\equiv 2i_1 \; mod \; 6$. This means that $\AAbar=\AA\cup \{ f_3(x_0,x_2), f_2(x_3,x_4)\cdot f_1(x_0,x_2) \}$.\\

For $\AA=\AAN{0}\cup \{ x_3x_4x_5 \}$ we get an extra condition $c_3+c_4+c_5=0$, i.e.
$$
c_5=c_3+c_4,\; c_1=0,\; 3c_2=2c_3=2c_4=0.
$$
Hence $f=(1, 1, {\om}^{2c_2}, {\om}^{3c_3}, {\om}^{3c_4}, {\om}^{3c_3+3c_4})$, where $\om=\sqr{6}$, $c_2\in \ZZ{3}$, $c_3,c_4\in \ZZ{2}$.\\

This means that $\GG\cong \ZZ{3} \oplus (\ZZ{2})^{\oplus 2}$ with generators
$$
f_1=(1,1,{\om}^{2},1,1,1),\; f_2=(1,1,1,{\om}^3,1,{\om}^3),\; f_3=(1,1,1,1,{\om}^3,{\om}^3).
$$

Then $i_2\equiv 0 \; mod \; 3$, $i_3\equiv i_4\equiv i_5 \; mod \; 2$. This means that $\AAbar=\AA\cup \{ f_3(x_0,x_1), x_3^2\cdot f_1(x_0,x_1), x_4^2\cdot f_1(x_0,x_1), x_5^2\cdot f_1(x_0,x_1) \}$.\\

For $\AA=\AAN{0}\cup \{ x_2x_3x_5, x_4^2x_1 \}$ we get extra conditions $c_1+2c_4=c_2+c_3+c_5=0$, i.e.
$$
c_1=c_2=0, \; c_3=c_5, \; 2c_3=2c_4=0.
$$
Hence $f=(1, 1, 1, {\om}^{c_3}, {\om}^{c_4}, {\om}^{c_3})$, where $\om=\sqr{2}$, $c_i\in \ZZ{2}$.\\

This means that $\GG\cong (\ZZ{2})^{\oplus 2}$ with generators
$$
f_1=(1,1,1,{\om},1,\om),\; f_2=(1,1,1,1,{\om},1).
$$

Then $i_4\equiv i_3+i_5\equiv 0 \; mod \; 2$. This means that $\AAbar=\AA\cup \{ f_3(x_0,x_1,x_2), f_2(x_3,x_5)\cdot f_1(x_0,x_1,x_2), x_4^2\cdot  f_1(x_0,x_1,x_2) \}$.\\

For $\AA=\AAN{0}\cup \{ x_1x_3x_4, x_2x_4x_5 \}$ we get extra conditions $c_1+c_3+c_4=c_2+c_4+c_5=0$, i.e.
$$
c_1=c_2=0,\; c_3=c_4=c_5, \; 2c_5=0.
$$
Hence $f=(1, 1, 1, {\om}^{c_5}, {\om}^{c_5}, {\om}^{c_5})$, where $\om=\sqr{2}$, $c_i\in \ZZ{2}$.\\

This means that $\GG\cong \ZZ{2}$ with a generator
$$
f=(1,1,1,{\om},\om,\om).
$$

Then $i_3+i_4+i_5\equiv 0 \; mod \; 2$. This means that $\AAbar=\AA\cup \{ f_3(x_0,x_1,x_2), f_2(x_3,x_4,x_5)\cdot f_1(x_0,x_1,x_2) \}$.\\

Let $\AAN{0} =\{ x_3^2x_0, x_4^2x_0, x_5^2x_0, x_0^3, x_1^3, x_2^3 \}$. Then $3c_1=3c_2=2c_3=2c_4=2c_5=0$. Hence $f=(1, {\om}^{2c_1}, {\om}^{2c_2}, {\om}^{3c_3}, {\om}^{3c_4}, {\om}^{3c_5})$, where $\om=\sqr{6}$, $c_1,c_2\in \ZZ{3}$, $c_3,c_4,c_5\in \ZZ{2}$.\\

For $\AA=\AAN{0}\cup \{ x_3x_4x_5 \}$ we get an extra condition $c_3+c_4+c_5=0$, i.e.
$$
c_5=c_3+c_4,\; 3c_1=3c_2=2c_3=2c_4=0.
$$
Hence $f=(1, {\om}^{2c_1}, {\om}^{2c_2}, {\om}^{3c_3}, {\om}^{3c_4}, {\om}^{3c_3+3c_4})$, where $\om=\sqr{6}$, $c_1,c_2\in \ZZ{3}$, $c_3,c_4\in \ZZ{2}$.\\

This means that $\GG\cong (\ZZ{2})^{\oplus 2}\oplus (\ZZ{3})^{\oplus 2}$ with generators
$$
f_1=(1,{\om}^{2},1,1,1,1),\; f_2=(1,1,{\om}^2,1,1,1),\; f_3=(1,1,1,{\om}^3,1,{\om}^3),\; f_4=(1,1,1,1,{\om}^3,{\om}^3).
$$

Hence $i_1\equiv i_2\equiv  0 \; mod \; 3$ and $i_3\equiv i_4 \equiv i_5 \; mod \; 2$. This means that $\AAbar=\AA$.\\

For $\AA=\AAN{0}\cup \{ x_1x_3x_4, x_1x_4x_5, x_2x_3x_5 \}$ we get extra conditions $c_1+c_3+c_4=c_1+c_4+c_5=c_2+c_3+c_5=0$, i.e.
$$
c_5=c_3=c_4,\; c_1=c_2=0,\; 2c_3=0.
$$
Hence $f=(1, 1,1, {\om}^{c_3}, {\om}^{c_3}, {\om}^{c_3})$, where $\om=\sqr{2}$, $c_i\in \ZZ{2}$.\\

This case has already appeared.\\

\subsection{Case of $2$ cubes.}

Let the cubes be $x_0^3, x_1^3$. If a longest cycle has length $4$, then $\AA$ is
\begin{center}
\begin{tikzpicture}

\coordinate (0) at (0,.5); \node[below, font=\scriptsize] at (0) {0}; \fill (0) circle (1pt); \draw (0) circle (1mm); 
\coordinate (1) at (1,.5); \node[below, font=\scriptsize] at (1) {1}; \fill (1) circle (1pt); \draw (1) circle (1mm); 
\coordinate (2) at (2,0); \node[below, font=\scriptsize] at (2) {2}; \fill (2) circle (1pt);
\coordinate (4) at (3,1); \node[above, font=\scriptsize] at (4) {4}; \fill (4) circle (1pt);
\coordinate (3) at (2,1); \node[above, font=\scriptsize] at (3) {3}; \fill (3) circle (1pt);
\coordinate (5) at (3,0); \node[below, font=\scriptsize] at (5) {5}; \fill (5) circle (1pt);

\arr{2}{3}; \arr{4}{5}; \arr{3}{4};\arr{5}{2};

\end{tikzpicture}
\end{center}

i.e. $\AA=\{ x_0^3,x_1^3,x_2^2x_3, x_3^2x_4, x_4^2x_5, x_5^2x_2 \}$. Then $3c_1=2c_2+c_3=2c_3+c_4=2c_4+c_5=2c_5+c_2=0 \;\; \mbox{in}\; \ZZ{d}$, i.e.
$$
c_3=-2c_2, c_4=4c_2, c_5=-8c_2, 3c_1=15c_2=0.
$$

Hence $d=15$, $f=(1, {\om}^{5c_1}, {\om}^{c_2}, {\om}^{-2c_2}, {\om}^{4c_2}, {\om}^{-8c_2})$, where $\om=\sqr{15}$.\\

This means that $\GG\cong \ZZ{3}\oplus \ZZ{15}$ with generators 
$$
f_1=(1,{\om}^5,1,1,1,1),\; f_2=(1,1,{\om},{\om}^{-2},{\om}^{4},{\om}^{-8}).
$$

Then $i_1\equiv 0 \; mod \; 3$, $i_2-2i_3+4i_4-8i_5\equiv 0 \; mod \; 15$. This means that $\AAbar=\AA$.\\

\subsubsection{Case of $2$ cubes. Length $3$ longest cycle.}

If a longest cycle has length $3$, then $\AA$ is either

\begin{center}

\end{center}

Suppose $x_2, x_3, x_4$ form a cycle. By Lemma 1 (\cite{Liendo}, Lemma 1.3) the $5$-th vertex should be connected to either a cube (say, to $0$) or to the cycle (say, to vertex $2$). In the former case there are neither singular pairs nor singular triples.\\

In the latter case there may be a singular pair $(x_4, x_5)$. In order to resolve it, we need to have either
\begin{itemize}
\item an edge between vertices $4$ and $5$ (hence an edge from $5$ to $4$ since otherwise we would get a cycle of length $4$), or
\item an edge from vertices $4$, $5$ to the $3$-rd vertex or to a cube, or
\item a dashed curve passing through $4$, $5$ and either $3$ or a cube (say, $0$).
\end{itemize}

Note that in the second case $5$ can not be connected to a cube by an edge, because this would give us a nonminimal $\AA$ (such configuration has appeared earlier).\\

This gives us the pictures shown above. Now let us do computations.\\

If $\AA =\{ x_2^2x_3, x_3^2x_4, x_4^2x_2, x_5^2x_0, x_0^3, x_1^3 \}$, then $3c_1=2c_5=2c_2+c_3=2c_3+c_4=2c_4+c_2=0$ in $\ZZ{d}$, i.e.
$$
c_3=-2c_2, c_4=4c_2, 3c_1=9c_2=2c_5=0.
$$
Hence $d=18$ and $f=(1, {\om}^{6c_1}, {\om}^{2c_2}, {\om}^{-4c_2}, {\om}^{8c_2}, {\om}^{9c_5})$, where $\om=\sqr{18}$.\\

This means that $\GG\cong \ZZ{2}\oplus \ZZ{3} \oplus \ZZ{9}$ with generators
$$
f_1=(1,{\om}^{6},1,1,1,1),\; f_2=(1,1,{\om}^{2},{\om}^{-4},{\om}^{8},1),\; f_3=(1,1,1,1,1,{\om}^9).
$$

Then $i_1\equiv 0 \; mod \; 3$, $i_5 \equiv 0 \; mod \; 2$ and $i_2-2i_3+4i_4 \equiv 0 \; mod \; 9$. This means that $\AAbar=\AA$.\\

Let $\AAN{0} =\{ x_2^2x_3, x_3^2x_4, x_4^2x_2,x_5^2x_2, x_0^3, x_1^3 \}$. Then $3c_1=2c_2+c_3=2c_3+c_4=2c_4+c_2=2c_5+c_2=0$ in $\ZZ{d}$, i.e.
$$
c_2=-2c_5, c_3=4c_5, c_4=-8c_5, 3c_1=18c_5=0.
$$

Hence $d=18$, $f=(1, {\om}^{6c_1}, {\om}^{-2c_5}, {\om}^{4c_5}, {\om}^{-8c_5}, {\om}^{c_5})$, where $\om=\sqr{18}$.\\

This means that $G_{\AAN{0}} \cong \ZZ{18} \oplus \ZZ{3}$ with generators
$$
f_1=(1,{\om}^{6},1,1,1,1),\; f_2=(1,1,{\om}^{-2},{\om}^{4},{\om}^{-8},\om).
$$

For $\AA=\AAN{0}\cup \{ x_5^2x_3 \}$ we get an extra condition $c_3+2c_5=0$, i.e. $6c_5=0$. Hence $d=6$, $f=(1, {\om}^{2c_1}, {\om}^{-2c_5}, {\om}^{-2c_5}, {\om}^{-2c_5}, {\om}^{c_5})$, where $\om=\sqr{6}$.\\

This means that $\GG\cong \ZZ{6} \oplus \ZZ{3}$ with generators
$$
f_1=(1,{\om}^{2},1,1,1,1),\; f_2=(1,1,{\om}^{-2},{\om}^{-2},{\om}^{-2},{\om}).
$$

Then $i_1\equiv 0 \; mod \; 3$ and $i_5\equiv 2(i_2+i_3+i_4)  \; mod \; 6$. This means that $\AAbar=\AA\cup \{ f_3(x_2,x_3,x_4), x_5^2\cdot f_1(x_2,x_3,x_4) \}$.\\

For $\AA=\AAN{0}\cup \{ x_5^2x_4 \}$ and $\AA=\AAN{0}\cup \{ x_4^2x_3 \}$ we get extra conditions $c_4+2c_5=0$ and $c_3+2c_4=0$ respectively, i.e. $6c_5=0$. This leads to the same group as above.\\

For $\AA=\AAN{0}\cup \{ x_3x_4x_5 \}$ we get an extra relation $c_3+c_4+c_5=0$, i.e. $3c_5=0$. Hence $d=6$, $f=(1, {\om}^{c_1}, {\om}^{c_5}, {\om}^{c_5}, {\om}^{c_5}, {\om}^{c_5})$, where $\om=\sqr{3}$.\\

This means that $\GG\cong (\ZZ{3})^{\oplus 2}$ with generators
$$
f_1=(1,{\om},1,1,1,1),\; f_2=(1,1,{\om},{\om},{\om},{\om}).
$$

Then $i_1\equiv i_2+i_3+i_4+i_5 \equiv 0 \; mod \; 3$. This means that $\AAbar=\AA\cup \{ f_3(x_2,x_3,x_4, x_5) \}$.\\

For $\AA=\AAN{0}\cup \{ x_0x_4x_5 \}$ we get an extra relation $c_4+c_5=0$, i.e. $7c_5=0$. Since $18c_5=7c_5=0$ in $\ZZ{d}$ and $(18,7)=1$, we conclude that $c_5=0$. Hence $d=3$, $f=(1, {\om}^{c_1}, 1, 1, 1, 1)$, where $\om=\sqr{3}$.\\

This means that $\GG\cong \ZZ{3}$ with a generator $f=(1,{\om},1,1,1,1)$ and $\AAbar=\AA\cup \{ f_3(x_0,x_2,x_3,x_4, x_5) \}$.\\

For $\AA=\AAN{0}\cup \{ x_4^2x_0 \}$ we get an extra relation $2c_4=0$, i.e. $2c_5=0$. Hence $d=6$, $f=(1, {\om}^{2c_1}, 1, 1, 1, {\om}^{3c_5})$, where $\om=\sqr{6}$.\\

This means that $\GG\cong \ZZ{3} \oplus \ZZ{2}$ with generators
$$
f_1=(1,{\om}^2,1,1,1,1),\; f_2=(1,1,1,1,1,{\om}^3).
$$

Then $i_1\equiv 0 \; mod \; 3$, $i_5\equiv 0 \; mod \; 2$. This means that $\AAbar=\AA\cup \{ f_3(x_0,x_2,x_3,x_4), x_5^2\cdot f_1(x_0,x_2,x_3,x_4) \}$.\\

\subsubsection{Case of $2$ cubes. Length $2$ longest cycle.}

Suppose that the length of a longest cycle is $2$. Let $x_2, x_3$ form such a cycle.\\

Then $\AA$ may be one of the following:

\begin{center}

\end{center}

Indeed, the remaining variables $x_4, x_5$ may also form a cycle (in which case there are neither singular pairs nor singular triples), or they may be connected by an edge (without forming a cycle), or they may be disconnected.\\

If $x_4, x_5$ do not form a cycle, but $5$ is connected to $4$, then by Lemma 1 (\cite{Liendo}, Lemma 1.3) vertex $4$ should be connected either to a cube (say, $0$) or to the cycle formed by $2$ and $3$ (say, $4$ is connected to $2$). In the former case there are neither singular pairs nor singular triples.\\

In the latter case $(x_3, x_4)$ may form a singular pair. They will not form a singular pair, if there is an edge between $3$ and $4$ (i.e. an edge from $4$ to $3$ since there are no cycles of length $3$). Alternatively, there should be either
\begin{itemize}
\item an edge from $3$ or $4$ to a cube or to the $5$-th vertex, or
\item a dashed curve passing through $3$, $4$ and either the $5$-th vertex or one of the cubes (say, $0$).
\end{itemize}

Note that in the first case there can be an edge only from $3$ to a cube (say, to $0$). Otherwise, either $x_4, x_5$ would form a cycle (if $4$ is connected to $5$) or there would be a length $4$ cycle (if $3$ is connected to $5$) or we would get a nonminimal $\AA$ (if $4$ is connected to a cube).\\

It may also happen that vertices $4$, $5$ are disconnected. In this case, by Lemma 1 (\cite{Liendo}, Lemma 1.3) each of them should be connected to either a cube or to the cycle.\\

If  $4$ and $5$ are connected to different cubes, then there are neither singular pairs nor singular triples.\\

Suppose that one of them is connected to a cube (say, $4$ is connected to $0$) and the other one is connected to the cycle (say, $5$ is connected to $2$). In this case $(x_3, x_5)$ may form a singular pair and $(x_3,x_4,x_5)$ may form a singular triple.\\

Both the singular pair and the singular triple will be resolved, if there is an edge between $3$, $5$, or an edge from one of them to  $4$, or a dashed curve through $3$, $4$, $5$. An edge from $3$ to $5$ would lead to a length $3$ cycle (which is not allowed). Neither an edge from $5$ to $4$ is allowed by our assumption. Hence we may have either an edge from $5$ to $3$ or an edge from $3$ to $4$.\\ 

Alternatively, one can resolve both the singular pair and the singular triple by having an edge from $3$ or $5$ to the second cube $1$ (only an edge from $3$ to $1$ is possible) or a dashed curve passing through $3$, $5$ and the second cube $1$.\\

If neither of these alternatives holds, then we need to resolve the singular pair $(x_3, x_5)$ and the singular triple $(x_3,x_4,x_5)$ separately.\\

The singular pair $(x_3, x_5)$ can be resolved by having either an edge from $3$ or $5$ to the cube $0$ or a dashed curve passing through these three vertices.\\

Then the singular triple $(x_3,x_4,x_5)$ can be resolved either by having an edge from $4$ to $3$ or to $5$ (only an edge from $4$ to $3$ is allowed), or by connecting a pair of $3$, $4$, $5$ by a dashed curve or just the vertex $4$ by an edge to the second cube $1$. Note that a dashed curve through vertices $1$, $3$, $5$ is not allowed, because it would lead to a nonminimal $\AA$.\\

Now suppose that $4$ and $5$ are connected to either the cycle or to a cube, but not to both of them.\\

Let us assume first that $4$ is connected to $3$ and $5$ is connected to $2$. Then there may be two singular pairs $(x_2, x_4)$ and $(x_3, x_5)$. We resolve them as above.\\

If $4$ and $5$ are connected to the same cube (say, to $0$), then we may have a singular pair $(x_4, x_5)$. It is analyzed the same way as we did earlier.\\

This gives us the pictures shown above.\\

Alternatively, $\AA$ may be one of the following:

\begin{center}

\end{center}

In these cases both vertices $4$ and $5$ are connected to the same vertex of the cycle (say, to $2$). There may be 
\begin{itemize}
\item three singular pairs: $(x_4,x_5)$, $(x_3,x_5)$, $(x_3,x_4)$, and
\item one singular triple $(x_3,x_4,x_5)$.
\end{itemize}

We may resolve all these singular pairs and the singular triple at once by including a dashed curve through $3$, $4$, $5$.\\

Note that by our assumptions no edges between vertices $3$, $4$, $5$ are allowed. Edges from $4$ or $5$ to cubes are not allowed either. Hence only dashed curves or an edge from $3$ to a cube can be used.\\

After we resolve the singular pairs and the singular triple the same way as we did in the previous examples, we obtain the pictures shown above. Now let us do computations.\\

If $\AA =\{ x_2^2x_3, x_3^2x_2, x_4^2x_0, x_5^2x_1, x_0^3, x_1^3 \}$, then $c_2=c_3$, $3c_1=3c_2=2c_4=2c_5+c_1=0$ in $\ZZ{d}$, i.e.
$$
c_1=-2c_5,\; c_2=c_3,\; 3c_2=2c_4=6c_5=0.
$$
Hence $d=6$ and $f=(1, {\om}^{-2c_5}, {\om}^{2c_2}, {\om}^{2c_2}, {\om}^{3c_4}, {\om}^{c_5})$, where $\om=\sqr{6}$.\\

This means that $\GG\cong \ZZ{2}\oplus \ZZ{3}\oplus \ZZ{6}$ with generators
$$
f_1=(1,{\om}^{-2},1,1,1,{\om}),\; f_2=(1,1,{\om}^{2},{\om}^{2},1,1),\; f_3=(1,1,1,1,{\om}^3,1).
$$

Then $i_5 \equiv 2i_1\; mod \; 6$, $i_2+i_3 \equiv 0 \; mod \; 3$, $i_4 \equiv 0 \; mod \; 2$. This means that $\AAbar=\AA\cup \{ f_3(x_2,x_3) \}$.\\

If $\AA =\{ x_2^2x_3, x_3^2x_2, x_4^2x_5, x_5^2x_4, x_0^3, x_1^3 \}$, then $c_2=c_3$, $c_4=c_5$, $3c_1=3c_2=3c_4=0$ in $\ZZ{d}$, i.e. $d=3$ and $f=(1, {\om}^{c_1}, {\om}^{c_2}, {\om}^{c_2}, {\om}^{c_4}, {\om}^{c_4})$, where $\om=\sqr{3}$.\\

This means that $\GG\cong (\ZZ{3})^{\oplus 3}$ with generators
$$
f_1=(1,{\om},1,1,1,1),\; f_2=(1,1,{\om},{\om},1,1),\; f_3=(1,1,1,1,{\om},\om).
$$

Then $i_1 \equiv i_2+i_3 \equiv i_4+i_5 \equiv 0 \; mod \; 3$. This means that $\AAbar=\AA\cup \{ f_3(x_2,x_3), f_3(x_4,x_5) \}$.\\

If $\AA =\{ x_2^2x_3, x_3^2x_2, x_4^2x_5, x_5^2x_0, x_0^3, x_1^3 \}$, then $c_2=c_3$, $3c_1=3c_2=2c_4+c_5=2c_5=0$ in $\ZZ{d}$, i.e.
$$
c_2=c_3,\; c_5=2c_4,\; 3c_1=3c_2=4c_4=0,
$$
Hence $d=12$ and $f=(1, {\om}^{4c_1}, {\om}^{4c_2}, {\om}^{4c_2}, {\om}^{3c_4}, {\om}^{6c_4})$, where $\om=\sqr{12}$.\\

This means that $\GG\cong (\ZZ{3})^{\oplus 2}\oplus \ZZ{4}$ with generators
$$
f_1=(1,{\om}^4,1,1,1,1),\; f_2=(1,1,{\om}^4,{\om}^4,1,1),\; f_3=(1,1,1,1,{\om}^3,{\om}^6).
$$

Then $i_1 \equiv i_2+i_3 \equiv 0 \; mod \; 3$, $i_4\equiv 2i_5  \; mod \; 4$. This means that $\AAbar=\AA\cup \{ f_3(x_2,x_3) \}$.\\

Let $\AAN{0} =\{ x_2^2x_3, x_3^2x_2, x_4^2x_2, x_5^2x_4, x_0^3, x_1^3 \}$. Then $c_2=c_3$, $3c_1=3c_2=2c_4+c_2=2c_5+c_4=0$ in $\ZZ{d}$, i.e.
$$
c_2=c_3=4c_5,\; c_4=-2c_5,\; 3c_1=12c_5=0.
$$
Hence $d=12$, $f=(1, {\om}^{4c_1}, {\om}^{4c_5}, {\om}^{4c_5}, {\om}^{-2c_5}, {\om}^{c_5})$, where $\om=\sqr{12}$.\\

This means that $G_{\AAN{0}} \cong \ZZ{12} \oplus \ZZ{3}$ with generators
$$
f_1=(1,{\om}^{4},1,1,1,1),\; f_2=(1,1,{\om}^{4},{\om}^{4},{\om}^{-2},\om).
$$

For $\AA=\AAN{0}\cup \{ x_3^2x_0 \}$ we get an extra condition $2c_3=0$, i.e. 
$$
c_2=c_3=0,\; c_4=2c_5,\; 3c_1=4c_5=0 \; \mbox{in} \; \ZZ{d}.
$$
Hence $d=12$, $f=(1, {\om}^{4c_1}, 1, 1, {\om}^{6c_5}, {\om}^{3c_5})$, where $\om=\sqr{12}$.\\

This means that $\GG\cong \ZZ{4} \oplus \ZZ{3}$ with generators
$$
f_1=(1,{\om}^{4},1,1,1,1),\; f_2=(1,1,1,1,{\om}^{6},{\om}^{3}).
$$

Then $i_1\equiv 0 \; mod \; 3$ and $i_5\equiv 2i_4  \; mod \; 4$. This means that $\AAbar=\AA\cup \{ f_3(x_0,x_2,x_3), x_4^2\cdot f_1(x_0,x_2,x_3) \}$.\\

For $\AA=\AAN{0}\cup \{ x_3x_4x_5 \}$ we get an extra condition $c_3+c_4+c_5=0$, i.e. 
$$
c_2=c_3=c_4=c_5,\; 3c_1=3c_5=0 \; \mbox{in} \; \ZZ{d}.
$$
Hence $d=3$, $f=(1, {\om}^{c_1}, {\om}^{c_5}, {\om}^{c_5}, {\om}^{c_5}, {\om}^{c_5})$, where $\om=\sqr{3}$.\\

This case has already appeared.\\

For $\AA=\AAN{0}\cup \{ x_0x_3x_4 \}$ we get an extra condition $c_3+c_4=0$, i.e. 
$$
c_2=c_3=c_4=3c_1=2c_5=0 \; \mbox{in} \; \ZZ{d}.
$$
Hence $d=6$, $f=(1, {\om}^{2c_1}, 1, 1, 1, {\om}^{3c_5})$, where $\om=\sqr{6}$.\\

This case has already appeared.\\

For $\AA=\AAN{0}\cup \{ x_4^2x_3 \}$ the extra condition $c_3+2c_4=0$ is automatically satisfied.\\

Hence $\GG\cong \ZZ{3} \oplus \ZZ{12}$ (as computed above).\\

Then $i_1\equiv 0 \; mod \; 3$ and $i_5-2i_4+4(i_2+i_3)\equiv 0  \; mod \; 12$. This means that  $\AAbar=\AA\cup \{ f_3(x_2,x_3) \}$.\\

Let $\AAN{0} =\{ x_2^2x_3, x_3^2x_2, x_4^2x_0, x_5^2x_2, x_0^3, x_1^3 \}$. Then $c_2=c_3$, $3c_1=3c_2=2c_4=2c_5+c_2=0$ in $\ZZ{d}$, i.e.
$$
c_2=c_3=-2c_5,\; 3c_1=2c_4=6c_5=0.
$$
Hence $d=6$, $f=(1, {\om}^{2c_1}, {\om}^{-2c_5}, {\om}^{-2c_5}, {\om}^{3c_4}, {\om}^{c_5})$, where $\om=\sqr{6}$.\\

This means that $G_{\AAN{0}} \cong \ZZ{6} \oplus \ZZ{3} \oplus \ZZ{2}$ with generators
$$
f_1=(1,{\om}^{2},1,1,1,1),\; f_2=(1,1,{\om}^{-2},{\om}^{-2},1,\om),\; f_3=(1,1,1,1,{\om}^{3},1).
$$

For $\AA=\AAN{0}\cup \{ x_5^2x_3 \}$ the extra condition $c_3+2c_5=0$ is automatically satisfied. Hence $\GG=G_{\AAN{0}} \cong \ZZ{6} \oplus \ZZ{3} \oplus \ZZ{2}$ as above.\\

Then $i_1\equiv 0 \; mod \; 3$, $i_4\equiv 0 \; mod \; 2$ and $i_5\equiv 2(i_2+i_3)  \; mod \; 6$. This means that $\AAbar=\AA\cup \{ f_3(x_2,x_3) \}$.\\

For $\AA=\AAN{0}\cup \{ x_3^2x_4 \}$ we get an extra condition $2c_3+c_4=0$, i.e. 
$$
c_2=c_3=c_4=3c_1=2c_5=0 \; \mbox{in} \; \ZZ{d}.
$$
Hence $d=6$, $f=(1, {\om}^{2c_1}, 1, 1, 1, {\om}^{3c_5})$, where $\om=\sqr{6}$.\\

This case has already appeared.\\

For $\AA=\AAN{0}\cup \{ x_3x_4x_5 \}$ we get an extra condition $c_3+c_4+c_5=0$, i.e. 
$$
c_5=c_4,\; c_2=c_3=3c_1=2c_4=0 \; \mbox{in} \; \ZZ{d}.
$$
Hence $d=6$, $f=(1, {\om}^{2c_1}, 1, 1, {\om}^{3c_5}, {\om}^{3c_5})$, where $\om=\sqr{6}$.\\

This means that $\GG\cong \ZZ{3} \oplus \ZZ{2}$ with generators
$$
f_1=(1,{\om}^2,1,1,1,1),\; f_2=(1,1,1,1,{\om}^3,{\om}^3).
$$
Then $i_1\equiv 0 \; mod \; 3$ and $i_4+i_5\equiv 0  \; mod \; 2$. This means that  $\AAbar=\AA\cup \{ f_3(x_0,x_2,x_3), f_2(x_4,x_5) \cdot f_1(x_0,x_2,x_3) \}$.\\

For $\AA=\AAN{0}\cup \{ x_3^2x_1 \}$ we get an extra condition $c_1+2c_3=0$, i.e. 
$$
c_1=c_2=c_3=-2c_5,\; 2c_4=6c_5=0 \; \mbox{in} \; \ZZ{d}.
$$
Hence $d=6$, $f=(1, {\om}^{-2c_5}, {\om}^{-2c_5}, {\om}^{-2c_5}, {\om}^{3c_4}, {\om}^{c_5})$, where $\om=\sqr{6}$.\\

This means that $\GG\cong \ZZ{6} \oplus \ZZ{2}$ with generators
$$
f_1=(1,{\om}^{-2},{\om}^{-2},{\om}^{-2},1,\om),\; f_2=(1,1,1,1,{\om}^3,1).
$$
Then $i_4\equiv 0 \; mod \; 2$ and $i_5\equiv 2(i_1+i_2+i_3)  \; mod \; 6$. This means that  $\AAbar=\AA\cup \{ f_3(x_1,x_2,x_3), x_5^2 \cdot f_1(x_1,x_2,x_3) \}$.\\

For $\AA=\AAN{0}\cup \{ x_1x_3x_5 \}$ we get an extra condition $c_1+c_3+c_5=0$, i.e. 
$$
c_1=c_2=c_3=c_5,\; 3c_1=2c_4=0 \; \mbox{in} \; \ZZ{d}.
$$
Hence $d=6$, $f=(1, {\om}^{2c_5}, {\om}^{2c_5}, {\om}^{2c_5}, {\om}^{3c_4}, {\om}^{2c_5})$, where $\om=\sqr{6}$.\\

This means that $\GG\cong \ZZ{2} \oplus \ZZ{3}$ with generators
$$
f_1=(1,{\om}^2,{\om}^{2},{\om}^{2},1,{\om}^{2}),\; f_2=(1,1,1,1,{\om}^3,1).
$$
Then $i_4\equiv 0 \; mod \; 2$ and $i_1+i_2+i_3+i_5\equiv 0  \; mod \; 3$. This means that  $\AAbar=\AA\cup \{ f_3(x_1,x_2,x_3,x_5) \}$.\\

For $\AAN{1}=\AAN{0}\cup \{ x_5^2x_0 \}$, $\AAN{1}=\AAN{0}\cup \{ x_3^2x_0 \}$ and for $\AA=\AAN{0}\cup \{ x_3^2x_0,x_4^2x_3 \}$, $\AA=\AAN{0}\cup \{ x_5^2x_0,x_4^2x_3 \}$ we get the same extra condition $2c_5=c_3+2c_4=2c_3=0$, i.e. 
$$
c_2=c_3=0,\; 3c_1=2c_4=2c_5=0 \; \mbox{in} \; \ZZ{d}.
$$
Hence $d=6$, $f=(1, {\om}^{2c_1}, 1, 1, {\om}^{3c_4}, {\om}^{3c_5})$, where $\om=\sqr{6}$.\\

This means that $G_{\AAN{1}}= \GG\cong (\ZZ{2})^{\oplus 2} \oplus \ZZ{3}$ with generators
$$
f_1=(1,{\om}^2,1,1,1,1),\; f_2=(1,1,1,1,{\om}^3,1),\; f_3=(1,1,1,1,1,{\om}^3).
$$
Then $i_4\equiv i_5\equiv 0 \; mod \; 2$ and $i_1\equiv 0  \; mod \; 3$. This means that  $\AAbar=\AA\cup \{ f_3(x_0,x_2,x_3), x_4^2\cdot f_1(x_0,x_2,x_3), x_5^2\cdot f_1(x_0,x_2,x_3) \}$.\\

For $\AA=\AAN{1}\cup \{ x_4^2x_1 \}$ we get an extra condition $c_1=0$, i.e. 
$$
c_1=c_2=c_3=0,\; 2c_4=2c_5=0 \; \mbox{in} \; \ZZ{d}.
$$
Hence $d=2$, $f=(1, 1, 1, 1, {\om}^{c_4}, {\om}^{c_5})$, where $\om=\sqr{2}$.\\

This means that $\GG\cong (\ZZ{2})^{\oplus 2}$ with generators
$$
f_1=(1,1,1,1,{\om},1),\; f_2=(1,1,1,1,1,{\om}).
$$
Then $i_4\equiv i_5\equiv 0 \; mod \; 2$. This means that  $\AAbar=\AA\cup \{ f_3(x_0,x_1,x_2,x_3), x_4^2\cdot f_1(x_0,x_1,x_2,x_3), x_5^2\cdot f_1(x_0,x_1,x_2,x_3) \}$.\\

For $\AA=\AAN{1}\cup \{ x_1x_3x_4 \}$ we get an extra condition $c_1+c_3+c_4=0$, i.e. 
$$
c_1=c_2=c_3=c_4=0,\; 2c_5=0.
$$
Hence $\GG\cong \ZZ{2}$ with a generator
$$
f=(1,1,1,1,1,{\om}),\; \om=\sqr{2},
$$
and  $\AAbar=\AA\cup \{ f_3(x_0,x_1,x_2,x_3,x_4), x_5^2\cdot f_1(x_0,x_1,x_2,x_3,x_4) \}$.\\

For $\AA=\AAN{1}\cup \{ x_1x_4x_5 \}$ we get an extra condition $c_1+c_4+c_5=0$, i.e. 
$$
c_1=c_2=c_3=0,\; c_5=c_4,\; 2c_5=0.
$$
Hence $\GG\cong \ZZ{2}$ with a generator
$$
f=(1,1,1,1,\om,{\om}),\; \om=\sqr{2},
$$
and  $\AAbar=\AA\cup \{ f_3(x_0,x_1,x_2,x_3), f_2(x_4,x_5)\cdot f_1(x_0,x_1,x_2,x_3) \}$.\\

For $\AAN{1}=\AAN{0}\cup \{ x_0x_3x_5 \}$ and for $\AA=\AAN{0}\cup \{ x_0x_3x_5,x_4^2x_3 \}$ we get the same extra condition $c_3+c_5=0$, i.e. 
$$
c_2=c_3=c_5=0,\; 3c_1=2c_4=0.
$$
These conditions already appeared above.\\

For $\AA=\AAN{1}\cup \{ x_4^2x_1 \}$ we get an extra condition $c_1=0$, i.e. 
$$
c_1=c_2=c_3=c_5=0,\; 2c_4=0.
$$
This case has already appeared.\\

For $\AA=\AAN{1}\cup \{ x_1x_3x_4 \}$ and $\AA=\AAN{1}\cup \{ x_1x_4x_5 \}$ we get the same extra condition $c_1+c_4=0$, i.e. 
$$
c_1=c_2=c_3=c_4=c_5=0.
$$
Hence $\GG=Id$ in these cases.\\

Let $\AAN{0} =\{ x_2^2x_3, x_3^2x_2, x_4^2x_3, x_5^2x_2, x_0^3, x_1^3 \}$. Then 
$$
c_2=c_3=-2c_4,\; 3c_1=6c_4=2(c_5-c_4)=0. 
$$

This means that $G_{\AAN{0}} \cong \ZZ{6} \oplus \ZZ{3} \oplus \ZZ{2}$ with generators
$$
f_1=(1,{\om}^{2},1,1,1,1),\; f_2=(1,1,{\om}^{-2},{\om}^{-2},\om,\om),\; f_3=(1,1,1,1,1,{\om}^{3}).
$$

For $\AAN{1}=\AAN{0}\cup \{ x_5^2x_3 \}$ and for $\AA=\AAN{0}\cup \{ x_5^2x_3,x_4^2x_2 \}$ extra conditions $c_3+2c_5=c_2+2c_4=0$ are automatically satisfied.\\

Hence $\GG=G_{\AAN{1}}=G_{\AAN{0}} \cong \ZZ{6} \oplus \ZZ{3} \oplus \ZZ{2}$ as above.\\

Then $i_1\equiv 0 \; mod \; 3$, $i_5\equiv 0 \; mod \; 2$, $i_4+i_5\equiv 2(i_2+i_3) \; mod \; 6$. This means that  $\AAbar=\AA\cup \{ f_3(x_2,x_3), x_4^2\cdot f_1(x_2,x_3), x_5^2\cdot f_1(x_2,x_3) \}$.\\

For $\AA=\AAN{1}\cup \{ x_2x_4x_5 \}$ we get an extra condition $c_2+c_4+c_5=0$, i.e. 
$$
c_2=c_3=-2c_4,\; c_5=c_4,\; 3c_1=6c_4=0.
$$
Hence $\GG\cong \ZZ{6} \oplus \ZZ{3}$ with generators
$$
f_1=(1,{\om}^2,1,1,1,1),\; f_2=(1,1,{\om}^{-2},{\om}^{-2},{\om},\om), \; \om=\sqr{6},
$$
and  $\AAbar=\AA\cup \{ f_3(x_2,x_3), f_2(x_4,x_5)\cdot f_1(x_2,x_3) \}$.\\

For $\AA=\AAN{1}\cup \{ x_2^2x_0 \}$ we get an extra condition $2c_2=0$, i.e. 
$$
c_2=c_3=0, \; 3c_1=2c_4=2c_5=0.
$$
This case was already considered.\\

For $\AA=\AAN{1}\cup \{ x_0x_2x_4 \}$ we get an extra condition $c_2+c_4=0$, i.e. 
$$
c_2=c_3=c_4=0, \; 3c_1=2c_5=0.
$$
This case was already considered.\\

For $\AAN{1}=\AAN{0}\cup \{ x_3^2x_0 \}$ and $\AA=\AAN{0}\cup \{ x_3^2x_0,x_2^2x_0 \}$ we get the same extra condition $2c_4=0$, i.e. 
$$
c_2=c_3=0, \; 3c_1=2c_4=2c_5=0.
$$
These conditions already appeared above. The same happens for $\AA=\AAN{1}\cup \{ x_0x_2x_4 \}$.\\

For $\AA=\AAN{1}\cup \{ x_2x_4x_5 \}$ we get an extra condition $c_4=c_5$, i.e. 
$$
c_2=c_3=0,\; c_4=c_5, \; 3c_1=2c_5=0.
$$
This case was already considered.\\

For $\AA=\AAN{1}\cup \{ x_2^2x_1 \}$ we get an extra condition $c_1=0$, i.e. 
$$
c_1=c_2=c_3=0, \; 2c_4=2c_5=0.
$$
This case was already considered.\\

For $\AA=\AAN{1}\cup \{ x_1x_2x_4 \}$ we get an extra condition $c_1=c_4$, i.e. 
$$
c_1=c_2=c_3=c_4=0, \; 2c_5=0.
$$
This case has already appeared.\\

For $\AAN{1}=\AAN{0}\cup \{ x_3x_4x_5 \}$ and $\AA=\AAN{0}\cup \{ x_3x_4x_5. x_2x_4x_5 \}$ we get the same extra condition $c_5=c_4$, i.e. 
$$
c_2=c_3=-2c_4,\; c_5=c_4,\; 3c_1=6c_4=0.
$$
These conditions already appeared above.\\

For $\AA=\AAN{1}\cup \{ x_0x_2x_4 \}$ we get an extra condition $c_2+c_4=0$, i.e. 
$$
c_2=c_3=c_4=c_5=0,\; 3c_1=0.
$$
This case was already considered.\\

For $\AAN{1}=\AAN{0}\cup \{ x_0x_3x_5 \}$ we get an extra condition $c_3+c_5=0$, i.e. 
$$
c_2=c_3=c_5=0,\; 3c_1=2c_4=0.
$$

For $\AA=\AAN{1}\cup \{ x_2x_4x_5 \}$ and $\AA=\AAN{1}\cup \{ x_0x_2x_4 \}$ we get the same extra condition $c_4=0$, which was already considered above.\\

For $\AA=\AAN{1}\cup \{ x_1x_2x_4 \}$ we get an extra condition $c_1=c_4=0$, which implies that $\GG=Id$.\\

Let $\AAN{0} =\{ x_2^2x_3, x_3^2x_2,x_4^2x_0,x_5^2x_0, x_2x_4x_5, x_0^3, x_1^3 \}$. Then 
$$
c_2=c_3=0,\;c_4=c_5,\;  3c_1=2c_5=0. 
$$

This means that $G_{\AAN{0}} \cong \ZZ{2} \oplus \ZZ{3}$ with generators
$$
f_1=(1,{\om}^{2},1,1,1,1),\; f_2=(1,1,1,1,{\om}^{3},{\om}^{3}),\; \om=\sqr{6}.
$$

For $\AA=\AAN{0}\cup \{ x_3^2x_4 \}$ we get an extra condition $c_4=0$, i.e.
$$
c_2=c_3=c_4=c_5=0,\; 3c_1=0.  
$$

This case has already appeared.\\

For $\AA=\AAN{0}\cup \{ x_3^2x_1 \}$ we get an extra condition $c_1=0$, i.e.
$$
c_1=c_2=c_3=0,\; c_4=c_5,\; 2c_5=0. 
$$

This case has already appeared.\\

For $\AA=\AAN{0}\cup \{ x_3x_4x_5 \}$ the extra condition is automatically satisfied.\\

This case has already appeared.\\

For $\AA=\AAN{0}\cup \{ x_1x_3x_4 \}$ we get an extra condition $c_1+c_3+c_4=0$, i.e.
$$
c_1=c_2=c_3=c_4=c_5=0.
$$
This means that $\GG =Id$.\\

If $\AA =\{ x_2^2x_3, x_3^2x_2, x_4^2x_0,x_5^2x_0,x_1x_4x_5, x_0^3, x_1^3 \}$, then $c_2=c_3$, $3c_1=3c_2=2c_4=2c_5=0$, $c_1+c_4+c_5=0$, i.e.
$$
c_1=0,\; c_2=c_3,\; c_4=c_5,\; 3c_2=2c_4=0. 
$$

This means that $\GG\cong \ZZ{3}\oplus \ZZ{2} $ with generators
$$
f_1=(1,1,{\om}^2,{\om}^2,1,1),\; f_2=(1,1,1,1,{\om}^3,{\om}^3),\; \om=\sqr{6}.
$$
Then $i_2+i_3 \equiv 0 \; mod \; 3$, $i_4+i_5 \equiv 0 \; mod \; 2$. This means that  $\AAbar=\AA\cup \{ f_3(x_2,x_3),f_1(x_0,x_1)\cdot f_2(x_4,x_5), f_3(x_0,x_1) \}$.\\

Let $\AAN{0} =\{ x_2^2x_3, x_3^2x_2,x_4^2x_2,x_5^2x_2, x_0^3, x_1^3 \}$. Then 
$$
c_2=c_3=-2c_5,\; 3c_1=2(c_4-c_5)=6c_5=0. 
$$

This means that $G_{\AAN{0}} \cong \ZZ{2} \oplus \ZZ{3}\oplus \ZZ{6}$ with generators
$$
f_1=(1,{\om}^{2},1,1,1,1),\; f_2=(1,1,1,1,{\om}^{3},1), \; f_3=(1,1,{\om}^{-2},{\om}^{-2},\om,\om),\; \om=\sqr{6}.
$$

If $\AA =\AAN{0}\cup \{ x_3x_4x_5 \}$ we get an extra condition $c_3+c_4+c_5=0$, i.e.
$$
c_2=c_3=-2c_5,\; c_4=c_5,\; 3c_1=6c_5=0
$$

This case was already considered.\\

For $\AA=\AAN{0}\cup \{ x_0x_4x_5, x_0x_3x_5, x_1x_3x_4 \}$, $\AA=\AAN{0}\cup \{ x_0x_4x_5, x_1x_3x_4, x_1x_3x_5 \}$ and $\AA=\AAN{0}\cup \{ x_0x_4x_5, x_1x_4x_5, x_0x_3x_4, x_0x_3x_5 \}$ we get extra conditions which imply that $c_1=c_2=c_3=c_4=c_5=0$. Hence $\GG=Id$ in these cases.\\

For $\AA=\AAN{0}\cup \{ x_3^2x_1, x_0x_4x_5 \}$ we get extra conditions $c_4+c_5=c_1+2c_3=0$, i.e.
$$
c_1=c_2=c_3=0,\; c_4=c_5,\; 2c_5=0.
$$
This case has already appeared above.\\

For $\AAN{1}=\AAN{0}\cup \{ x_3^2x_0, x_1x_3x_4 \}$ we get extra conditions $2c_3=c_1+c_3+c_4=0$, i.e.
$$
c_1=c_2=c_3=c_4=0,\; 2c_5=0.
$$
All such cases were considered above.\\

\subsubsection{Case of $2$ cubes. No cycles.}

If there are no cycles at all, then $\AA$ may be one of the following:

\begin{center}

\end{center}

Let us consider the subgraph formed by the variables $x_2$, $x_3$, $x_4$, $x_5$ which are not cubes. The length of a longest path in this subgraph may be either $4$ or $3$ or $2$ or $1$ (in which case the subgraph is totally disconnected).\\

If $2$ is connected to $3$, $3$ is connected to $4$ and $4$ is connected to $5$, then by Lemma 1 (\cite{Liendo}, Lemma 1.3) the $5$-th vertex should be connected to a cube (since otherwise there would be a cycle). In this case we have neither singular pairs nor singular triples.\\

Suppose that a longest path has length $3$. Let $2$ be connected to $3$ and $3$ be connected to $4$. Then $4$ will have to be connected to a cube (say, to $0$). It can be connected neither to $5$ (since there are no length $4$ paths) nor to $2$ nor to $3$ (since there are no cycles).\\

Then $5$ can be connected either to the other cube $1$, or to the same cube $0$, or to the path (i.e. either to $3$ or to $4$). In the first case there are neither singular pairs nor singular triples.\\

If $5$ is connected to $3$, then there may be one singular pair $(x_2, x_5)$ and one singular triple $(x_2, x_4, x_5)$. They are resolved as usual. Note that the configuration is symmetric with respect to vertices $2$ and $5$.\\

If $5$ is connected to $4$, then there may be one singular pair $(x_3, x_5)$. It is resolved as usual.\\

If $5$ is connected to the same cube $0$ as the length $3$ path, then there may be one singular pair $(x_4, x_5)$ and one singular triple $(x_2, x_4, x_5)$. We repeat the same procedure as above in order to resolve them.\\

If the length of a longest path in the subgraph formed by vertices $2$, $3$, $4$, $5$ is equal to $2$, then we may assume that $2$ is connected to $3$ and $3$ is connected to a cube (say, to $0$). Note that $3$ can not be connected to $4$ or $5$, because this would lead to a path of length $3$.\\

The remaining vertices $4$ and $5$ may also form a length $2$ path. Let us assume that $4$ is connected to $5$. Then by Lemma 1 (\cite{Liendo}, Lemma 1.3) the $5$-th vertex should be connected to a cube. It can not be connected to vertices $2$, $3$, $4$, because there are neither cycles no paths of length $3$ according to our assumptions.\\

If $3$ and $5$ are connected to different cubes, then there are neither singular pairs nor singular triples.\\

This analysis gives us the pictures shown above.\\

Vertices $3$ and $5$ may be also connected to the same cube $0$:

\begin{center}

\end{center}

In this case we may have one singular pair $(x_3, x_5)$, which we resolve as usual.\\

Suppose that the remaining vertices $4$ and $5$ are disconnected.\\

In this case $\AA$ may be one of the following:

\begin{center}
%
\end{center}

By Lemma 1 (\cite{Liendo}, Lemma 1.3) and our assumption the $4$-th and the $5$-th vertices should be connected either to a cube or to the length $2$ path formed by vertices $2$ and $3$.\\

Let us assume first that $5$ is connected to $3$ (it can not be connected to $2$ since there are no length $3$ paths) and $4$ is connected to the second cube $1$. In this case we may have one singular pair $(x_2, x_5)$ and one singular triple $(x_2, x_4,x_5)$. We resolve them as usual. Note that this configuration is symmetric with respect to vertices $2$, $5$.\\

Suppose that $4$ and $5$ are connected to different cubes (say, $5$ is connected to $0$ and $4$ is connected to $1$). Then we may get one singular pair $(x_3, x_5)$ and one singular triple $(x_3, x_4,x_5)$. We resolve them as usual.\\

If $4$ and $5$ are connected to the same cube $1$, then we may get
\begin{itemize} 
\item one singular pair $(x_4, x_5)$ and
\item two singular triples $(x_2, x_4,x_5)$ and $(x_3, x_4,x_5)$.
\end{itemize}

We resolve them as usual. Note that the configuration is symmetric with respect to vertices $4$ and $5$.\\

The next possibility is that one of the vertices $4$, $5$ is connected to the length $2$ path formed by $2$, $3$ and the other one is connected to the same cube $0$ as $3$. We may assume that $4$ is connected to $3$ (it can not be connected to $2$ since there are no length $3$ paths) and $5$ is connected to $0$. In this case there may be 
\begin{itemize} 
\item two singular pairs $(x_2, x_4)$, $(x_3, x_5)$ and
\item one singular triple $(x_2, x_4,x_5)$.
\end{itemize}

We resolve them as usual. Note that the configuration is symmetric with respect to vertices $2$ and $4$.\\

It may also happen that both vertices $4$ and $5$ are connected to the length $2$ path (i.e. to the $3$-rd vertex). In this case there may be 
\begin{itemize} 
\item three singular pairs $(x_2, x_4)$, $(x_2, x_5)$, $(x_4, x_5)$ and
\item one singular triple $(x_2, x_4,x_5)$.
\end{itemize}
We resolve them as usual. Note that the configuration is symmetric with respect to vertices $2$, $4$ and $5$.\\

Finally, both vertices $4$ and $5$ may be connected to the same cube $0$ as vertex $3$. In this case there may be 
\begin{itemize} 
\item three singular pairs $(x_3, x_4)$, $(x_3, x_5)$, $(x_4, x_5)$ and
\item two singular triples $(x_2, x_4,x_5)$, $(x_3, x_4,x_5)$.
\end{itemize}
We resolve them as usual. Note that the configuration is symmetric with respect to vertices $4$ and $5$.\\

This analysis gives us the pictures shown above.\\

The remaining situation we need to consider is when vertices $2$, $3$, $4$, $5$ are totally disconnected. By Lemma 1 (\cite{Liendo}, Lemma 1.3) each of them should be connected to a cube. We get the following pictures:

\begin{center}

\end{center}

Indeed, it may happen that pairs of vertices $2$, $3$, $4$, $5$ are connected to two different cubes (say, $2$ and $3$ are connected to $0$, while $4$ and $5$ are connected to $1$). In this case there may be 
\begin{itemize} 
\item two singular pairs $(x_2, x_3)$, $(x_4, x_5)$ and
\item four singular triples $(x_2, x_4,x_5)$, $(x_3, x_4,x_5)$, $(x_2, x_3,x_4)$, $(x_2, x_3,x_5)$.
\end{itemize}
We resolve them as usual. Note that the configuration is symmetric with respect to vertices $2$ and $3$ as well as $4$ and $5$.\\

It may also happen that three of the vertices $2$, $3$, $4$, $5$ are connected to one cube (say, $2$, $3$ and $4$ are connected to $0$) and the fourth of them is connected to the other cube (i.e. $5$ is connected to $1$). In this case there may be 
\begin{itemize} 
\item three singular pairs $(x_2, x_3)$, $(x_2, x_4)$, $(x_3, x_4)$ and
\item three singular triples $(x_2, x_3, x_5)$, $(x_2, x_4, x_5)$, $(x_3, x_4, x_5)$.
\end{itemize}
We resolve them as usual. Note that the configuration is symmetric with respect to vertices $2$, $3$ and $4$.\\

Finally, all four vertices $2$, $3$, $4$, $5$ may be connected to one and the same cube (say, to $0$). In this case there may be 
\begin{itemize} 
\item six singular pairs $(x_2, x_3)$, $(x_2, x_4)$, $(x_2, x_5)$, $(x_3, x_4)$, $(x_3, x_5)$, $(x_4, x_5)$ and
\item four singular triples $(x_2, x_3, x_4)$, $(x_2, x_3, x_5)$, $(x_2, x_4, x_5)$, $(x_3, x_4, x_5)$.
\end{itemize}
We resolve them as usual. Note that the configuration is symmetric with respect to vertices $2$, $3$, $4$, $5$.\\

This gives us the pictures shown above. Now let us do computations.\\

If $\AA =\{ x_2^2x_3, x_3^2x_4, x_4^2x_5,x_5^2x_0, x_0^3, x_1^3 \}$, then 
$$
c_3=-2c_2,\; c_4=4c_2,\; c_5=-8c_2,\; 3c_1=16c_2=0. 
$$

This means that $\GG\cong \ZZ{3}\oplus \ZZ{16} $ with generators
$$
f_1=(1,{\om},1,1,1,1),\; f_2=(1,1,\eta,{\eta}^{-2},{\eta}^4,{\eta}^{-8}),\; \om=\sqr{3}, \eta=\sqr{16}.
$$
Then $i_1 \equiv 0 \; mod \; 3$, $i_2-2i_3+4i_4-8i_5 \equiv 0 \; mod \; 16$. This means that  $\AAbar=\AA$.\\

If $\AA =\{ x_2^2x_3, x_3^2x_4, x_4^2x_0,x_5^2x_1, x_0^3, x_1^3 \}$, then 
$$
c_3=-2c_2,\; c_4=4c_2,\; c_1=-2c_5,\; 8c_2=6c_5=0. 
$$

This means that $\GG\cong \ZZ{6}\oplus \ZZ{8} $ with generators
$$
f_1=(1,1,{\om},{\om}^{-2},{\om}^4,1),\; f_2=(1,{\eta}^{-2},1,1,1,\eta),\; \om=\sqr{8}, \eta=\sqr{6}.
$$
Then $i_5 \equiv 2i_1 \; mod \; 6$, $i_2-2i_3+4i_4 \equiv 0 \; mod \; 8$. This means that  $\AAbar=\AA$.\\

If $\AA =\{ x_2^2x_3, x_3^2x_0, x_4^2x_5,x_5^2x_1, x_0^3, x_1^3 \}$, then 
$$
c_3=-2c_2,\; c_5=-2c_4,\; c_1=4c_4,\; 4c_2=12c_4=0. 
$$

This means that $\GG\cong \ZZ{4}\oplus \ZZ{12} $ with generators
$$
f_1=(1,1,{\om},{\om}^{-2},1,1),\; f_2=(1,{\eta}^{4},1,1,\eta,{\eta}^{-2}),\; \om=\sqr{4}, \eta=\sqr{12}.
$$
Then $i_2 \equiv 2i_3 \; mod \; 4$, $4i_1+i_4-2i_5 \equiv 0 \; mod \; 12$. This means that  $\AAbar=\AA$.\\

Let $\AAN{0} =\{ x_2^2x_3, x_3^2x_4,x_4^2x_0,x_5^2x_3, x_0^3, x_1^3 \}$. Then 
$$
c_3=-2c_2,\; c_4=4c_2,\; 3c_1=2(c_5-c_2)=8c_2=0. 
$$

For $\AA=\AAN{0}\cup \{ x_0x_2x_5, x_1x_4x_5 \}$ we get extra conditions $c_2+c_5=c_1+c_4+c_5=0$, i.e.
$$
c_1=c_2=c_3=c_4=c_5=0.  
$$
Hence $\GG=Id$.\\

For $\AAN{1}=\AAN{0}\cup \{ x_5^2x_0 \}$ and $\AA=\AAN{0}\cup \{ x_5^2x_4 \}$ we get the same extra condition $2c_5=2c_2+c_4=c_4+2c_5=0$, i.e.
$$
c_3=c_4=0,\; 3c_1=2c_2=2c_5=0. 
$$

This case was already considered.\\

For $\AA=\AAN{1}\cup \{ x_1x_2x_4 \}$ we get an extra condition $c_1+c_2+c_4=0$, i.e.
$$
c_1=c_2=c_3=c_4=0,\; 2c_5=0. 
$$

This case was already considered.\\

For $\AA=\AAN{1}\cup \{ x_1x_4x_5 \}$ we get an extra condition $c_1+c_4+c_5=0$, i.e.
$$
c_1=c_2=c_3=c_4=0,\; 2c_2=0. 
$$

This case was already considered.\\

For $\AA=\AAN{0}\cup \{ x_2x_4x_5 \}$ we get an extra condition $c_2+c_4+c_5=0$, i.e.
$$
c_3=2c_2,\; c_5=-c_2,\; c_4=0,\; 3c_1=4c_2=0. 
$$

This case was already considered.\\

For $\AA=\AAN{0}\cup \{ x_1x_2x_5 \}$ we get an extra condition $c_1+c_2+c_5=0$, i.e.
$$
c_3=2c_2,\; c_5=-c_2,\; c_1=c_4=0,\; 4c_2=0. 
$$

This means that $\GG\cong \ZZ{4}$ with a generator
$$
f=(1,1,{\om},{\om}^{2},1,{\om}^3),\; \om=\sqr{4}.
$$
Then $i_5 \equiv i_2+2i_3 \; mod \; 4$. This means that  $\AAbar=\AA\cup \{  f_3(x_0,x_1,x_4), x_3^2\cdot f_1(x_0,x_1,x_4), x_2x_5\cdot f_1(x_0,x_1,x_4)  \}$.\\

Let $\AAN{0} =\{ x_2^2x_3, x_3^2x_4,x_4^2x_0,x_5^2x_0, x_0^3, x_1^3 \}$. Then 
$$
c_3=-2c_2,\; c_4=4c_2,\; 3c_1=8c_2=2c_5=0. 
$$

For $\AA=\AAN{0}\cup \{ x_2x_4x_5 \}$ we get an extra condition $c_2+c_4+c_5=0$, i.e.
$$
c_5=c_2,\; c_3=c_4=0,\; 3c_1=2c_2=0. 
$$

This case was already considered.\\

For $\AA=\AAN{0}\cup \{ x_1x_4x_5 \}$ we get an extra condition $c_1+c_4+c_5=0$, i.e.
$$
c_1=0,\; c_3=-2c_2,\; c_4=c_5=4c_2,\; 8c_2=0. 
$$

This means that $\GG\cong \ZZ{8}$ with a generator
$$
f_1=(1,1,{\om},{\om}^{-2},{\om}^4,{\om}^4),\; \om=\sqr{8}. 
$$
Then $i_2 \equiv 2i_3+4(i_4+i_5) \; mod \; 8$. This means that  $\AAbar=\AA\cup \{  f_3(x_0,x_1), f_2(x_4,x_5)\cdot f_1(x_0,x_1), x_3^2x_5 \}$.\\

For $\AA=\AAN{0}\cup \{ x_5^2x_4 \}$ we get an extra condition $c_4=0$, i.e.
$$
c_3=2c_2,\; c_4=0,\; 3c_1=4c_2=2c_5=0.
$$

This means that $\GG\cong \ZZ{2}\oplus \ZZ{3} \oplus \ZZ{4}$ with generators
$$
f_1=(1,{\om},1,1,1,1),\; f_2=(1,1,{\eta},{\eta}^2,1,1),\; f_3=(1,1,1,1,1,\mu),\; \om=\sqr{3}, \eta=\sqr{4}, \mu=\sqr{2}. 
$$
Then $i_1\equiv 0 \; mod \; 3$, $i_5\equiv 0 \; mod \; 2$, $i_2 \equiv 2i_3 \; mod \; 4$. This means that  $\AAbar=\AA\cup \{  f_3(x_0,x_4), x_5^2 \cdot f_1(x_0,x_4), x_3^2 \cdot f_1(x_0,x_4) \}$.\\

For $\AAN{1}=\AAN{0}\cup \{ x_3x_4x_5 \}$  we get an extra condition $c_5=-2c_2$, i.e.
$$
c_3=c_5=2c_2,\; c_4=0,\; 3c_1=4c_2=0.
$$
All such cases either have already appeared above or will appear and will be analyzed below.\\

Let $\AAN{0} =\{ x_2^2x_3, x_3^2x_4,x_4^2x_0,x_5^2x_4, x_0^3, x_1^3 \}$. Then 
$$
c_3=-2c_2,\; c_4=4c_2,\; 3c_1=8c_2=2(c_5-2c_2)=0. 
$$

For $\AA=\AAN{0}\cup \{ x_0x_3x_5 \}$ we get an extra condition $c_3+c_5=0$, i.e.
$$
c_3=-2c_2,\; c_5=2c_2,\; c_4=4c_2,\; 3c_1=8c_2=0. 
$$

This means that $\GG\cong \ZZ{3} \oplus \ZZ{8}$ with generators
$$
f_1=(1,{\om},1,1,1,1),\; f_2=(1,1,{\eta},{\eta}^{-2},{\eta}^{4},{\eta}^{2}),\; \om=\sqr{3}, \eta=\sqr{8}. 
$$
Then $i_1\equiv 0 \; mod \; 3$, $i_2-2i_3+2i_5 \equiv 4i_4 \; mod \; 8$. This means that  $\AAbar=\AA$.\\

For $\AA=\AAN{0}\cup \{ x_2x_3x_5 \}$ we get an extra condition $c_2+c_3+c_5=0$, i.e.
$$
c_3=c_4=0,\; c_5=c_2,\; c_4=4c_2,\; 3c_1=2c_2=0. 
$$

This case was already considered.\\

For $\AA=\AAN{0}\cup \{ x_1x_3x_5 \}$ we get an extra condition $c_1+c_3+c_5=0$, i.e.
$$
c_1=0,\; c_3=-2c_2,\; c_4=4c_2,\; c_5=2c_2,\; 8c_2=0. 
$$

This means that $\GG\cong \ZZ{8}$ with a generator
$$
f=(1,1,{\om},{\om}^{-2},{\om}^{4},{\om}^{2}),\; \om=\sqr{8}. 
$$
Then $i_2-2i_3+2i_5 \equiv 4i_4 \; mod \; 8$. This means that  $\AAbar=\AA\cup \{ f_3(x_0,x_1), x_4^2\cdot f_1(x_0,x_1), x_3x_5\cdot f_1(x_0,x_1) \}$.\\

For $\AA=\AAN{0}\cup \{ x_3^2x_0 \}$ we get an extra condition $2c_3=0$, i.e.
$$
c_3=2c_2,\; c_4=0,\; 3c_1=4c_2=2c_5=0. 
$$

This case was already considered.\\

Let $\AAN{0} =\{ x_2^2x_3, x_3^2x_0,x_4^2x_5,x_5^2x_0, x_0^3, x_1^3 \}$. Then 
$$
c_3=2c_2,\; c_5=2c_4,\; 3c_1=4c_2=4c_4=0. 
$$

For $\AA=\AAN{0}\cup \{ x_1x_3x_5 \}$ we get an extra condition $c_1+c_3+c_5=0$, i.e.
$$
c_1=0,\; c_3=c_5=2c_2,\; 4c_2=2(c_4-c_2)=0. 
$$

This means that $\GG\cong \ZZ{2} \oplus \ZZ{4}$ with generators
$$
f_1=(1,1,{\om},{\om}^2,{\om},{\om}^2),\; f_2=(1,1,1,1,\eta ,1),\; \om=\sqr{4}, \eta=\sqr{2}. 
$$
Then $i_4\equiv 0 \; mod \; 2$, $i_2+i_4\equiv 2(i_3+i_5) \; mod \; 4$. This means that  $\AAbar=\AA\cup $ $ \{  f_3(x_0,x_1),$ $ x_2^2x_5,$ $x_4^2x_3, $ $f_1(x_0,x_1)\cdot f_2(x_3,x_5) \} $.\\

For $\AA=\AAN{0}\cup \{ x_3x_4x_5 \}$ we get an extra condition $c_3+c_4+c_5=0$, i.e.
$$
c_3=c_4=2c_2,\; c_5=0,\; 3c_1=4c_2=0. 
$$

This case was already considered.\\

Let $\AAN{0} =\{ x_2^2x_3, x_3^2x_0,x_4^2x_1,x_5^2x_3, x_0^3, x_1^3 \}$. Then 
$$
c_1=-2c_4,\; c_3=2c_2,\; 4c_2=6c_4=2(c_5-c_2)=0. 
$$

For $\AA=\AAN{0}\cup \{ x_0x_2x_5 \}$ we get an extra condition $c_2+c_5=0$, i.e.
$$
c_1=-2c_4,\; c_3=2c_2,\; c_5=-c_2,\; 4c_2=6c_4=0. 
$$

This means that $\GG\cong \ZZ{4} \oplus \ZZ{6}$ with generators
$$
f_1=(1,1,{\om},{\om}^{2},1,{\om}^{-1}),\; f_2=(1,{\eta}^{-2},1,1,\eta ,1),\; \om=\sqr{4}, \eta=\sqr{6}. 
$$
Then $i_4\equiv 2i_1 \; mod \; 6$, $i_2\equiv 2i_3+i_5 \; mod \; 4$. This means that  $\AAbar=\AA$.\\

For $\AA=\AAN{0}\cup \{ x_2x_4x_5 \}$ we get an extra condition $c_2+c_4+c_5=0$, i.e.
$$
c_1=0,\; c_3=2c_2,\; c_5=c_4-c_2,\; 4c_2=2c_4=0.
$$

This means that $\GG\cong \ZZ{4} \oplus \ZZ{2}$ with generators
$$
f_1=(1,1,{\om},{\om}^{2},1,{\om}^{-1}),\; f_2=(1,1,1,1,\eta ,\eta),\; \om=\sqr{4}, \eta=\sqr{2}. 
$$
Then $i_4+i_5 \equiv 0 \; mod \; 2$, $i_2\equiv 2i_3+i_5 \; mod \; 4$. This means that  $\AAbar=\AA\cup \{  f_3(x_0,x_1), x_4^2x_0, x_3^2x_1 \}$.\\

For $\AA=\AAN{0}\cup \{ x_5^2x_0 \}$ we get an extra condition $2c_5=0$, i.e.
$$
c_1=-2c_4,\; c_3=0,\; 2c_2=6c_4=2c_5=0.
$$

This means that $\GG\cong \ZZ{6} \oplus (\ZZ{2})^{\oplus 2}$ with generators
$$
f_1=(1,{\om}^{-2},1,1,{\om},1),\; f_2=(1,1,\eta ,1,1,1),\; f_3=(1,1,1,1,1,\eta ),\; \om=\sqr{6}, \eta=\sqr{2}. 
$$
Then $i_2\equiv i_5 \equiv 0 \; mod \; 2$, $i_4\equiv 2i_1 \; mod \; 6$. This means that  $\AAbar=\AA\cup \{  f_3(x_0,x_3), x_2^2\cdot f_1(x_0,x_3), x_5^2\cdot f_1(x_0,x_3) \}$.\\

For $\AA=\AAN{0}\cup \{ x_4^2x_0, x_2^2x_1 \}$ we get an extra condition $2c_4=c_1+2c_2=0$, i.e.
$$
c_1=c_3=0,\; 2c_2=2c_4=2c_5=0.
$$

This means that $\GG\cong (\ZZ{2})^{\oplus 3}$ with generators
$$
f_1=(1,1,{\om},1,1,1),\; f_2=(1,1,1,1,\om ,1),\; f_3=(1,1,1,1,1,\om ),\; \om=\sqr{2}, 
$$
and  $\AAbar=\AA\cup \{  f_3(x_0,x_1,x_3), x_2^2\cdot f_1(x_0,x_1,x_3), x_4^2\cdot f_1(x_0,x_1,x_3),x_5^2\cdot f_1(x_0,x_1,x_3) \}$.\\

For $\AA=\AAN{0}\cup \{ x_4^2x_0, x_1x_2x_5 \}$ we get an extra condition $2c_4=c_1+c_2+c_5=0$, i.e.
$$
c_1=0,\; c_3=2c_2,\; c_5=-c_2,\; 4c_2=2c_4=0.
$$

This means that $\GG\cong \ZZ{2} \oplus \ZZ{4}$ with generators
$$
f_1=(1,1,{\om},{\om}^2,1,{\om}^{-1}),\; f_2=(1,1,1,1,\eta ,1),\; \om=\sqr{4}, \eta=\sqr{2}, 
$$
Then $i_4\equiv 0 \; mod \; 2$, $i_2\equiv i_5+2i_3 \; mod \; 4$. This means that  $\AAbar=\AA\cup \{  f_3(x_0,x_1), x_3^2\cdot f_1(x_0,x_1), x_0x_2x_5 \}$.\\

For $\AAN{1}=\AAN{0}\cup \{ x_0x_4x_5 \}$ we get an extra conditions $c_4+c_5=0$, i.e.
$$
c_1=c_3=0,\; c_5=c_4,\; 2c_2=2c_4=0.
$$
All such cases have already appeared.\\

Let $\AAN{0} =\{ x_2^2x_3, x_3^2x_0,x_4^2x_1,x_5^2x_0, x_0^3, x_1^3 \}$. Then 
$$
c_1=-2c_4,\; c_3=2c_2,\; 4c_2=6c_4=2c_5=0. 
$$

For $\AA=\AAN{0}\cup \{ x_3^2x_1, x_4^2x_3 \}$ we get extra conditions $c_1=c_3=0$, i.e.
$$
c_1=c_3=0,\; 2c_2=2c_4=2c_5=0. 
$$

This case was already considered.\\

For $\AA=\AAN{0}\cup \{ x_3^2x_1, x_2x_3x_4 \}$ we get extra conditions $c_1=0$, $c_4=c_2$, i.e.
$$
c_1=c_3=0,\; c_4=c_2,\; 2c_2=2c_5=0. 
$$

This case was already considered.\\

For $\AA=\AAN{0}\cup \{ x_3^2x_1, x_2x_4x_5 \}$ we get extra conditions $c_1=c_2+c_4+c_5=0$, i.e.
$$
c_1=c_3=0,\; c_2=c_4+c_5,\; 2c_4=2c_5=0. 
$$

This means that $\GG\cong (\ZZ{2})^{\oplus 2}$ with generators
$$
f_1=(1,1,{\om},1,1,{\om}),\; f_2=(1,1,{\om},1,{\om},1),\; \om=\sqr{2}, 
$$
Then $i_2+i_5\equiv i_2+i_4 \equiv 0 \; mod \; 2$. This means that  $\AAbar=\AA\cup \{  f_3(x_0,x_1,x_3), x_2^2\cdot f_1(x_0,x_1,x_3), x_4^2\cdot f_1(x_0,x_1,x_3), x_5^2\cdot f_1(x_0,x_1,x_3) \}$.\\

Note that adding monomial $x_5^2x_1$ to $\AAN{0}$ leads to the same extra condition $c_1=0$ as adding monomial $x_3^2x_1$.\\

For $\AAN{1}=\AAN{0}\cup \{ x_1x_3x_5 \}$ we get an extra conditions $c_1+c_3+c_5=0$, i.e.
$$
c_1=0,\; c_3=c_5=2c_2,\; 4c_2=2c_4=0. 
$$

For $\AA=\AAN{1}\cup \{ x_4^2x_3 \}$ we get an extra conditions $c_3=0$, i.e.
$$
c_1=c_3=c_5=0,\; 2c_2=2c_4=0. 
$$

This case was already considered.\\

For $\AA=\AAN{1}\cup \{ x_2x_3x_4 \}$ and $\AA=\AAN{1}\cup \{ x_2x_4x_5 \}$ we get the same extra condition $c_4=c_2$, i.e.
$$
c_1=c_3=c_5=0,\; c_4=c_2,\; 2c_2=0. 
$$

This case was already considered.\\

For $\AA=\AAN{0}\cup \{ x_2x_3x_5 \}$ we get an extra condition $c_2+c_3+c_5=0$, i.e.
$$
c_1=-2c_4,\; c_3=0, \; c_5=c_2,\; 2c_2=6c_4=0.
$$

This case was already considered.\\

For $\AA=\AAN{0}\cup \{ x_3x_4x_5 \}$ we get an extra condition $c_3+c_4+c_5=0$, i.e.
$$
c_1=0,\; c_3=2c_2, c_5=c_4+2c_2,\; 4c_2=2c_4=0. 
$$

This means that $\GG\cong \ZZ{2} \oplus \ZZ{4}$ with generators
$$
f_1=(1,1,{\om},{\om}^2,1,{\om}^2),\; f_2=(1,1,1,1,\eta,\eta),\; \om=\sqr{4}, \eta=\sqr{2}, 
$$
Then $i_4+i_5\equiv 0 \; mod \; 2$, $i_2\equiv 2(i_3+i_5) \; mod \; 4$. This means that  $\AAbar=\AA\cup \{  f_3(x_0,x_1), x_3^2\cdot f_1(x_0,x_1), x_4^2\cdot f_1(x_0,x_1), x_5^2\cdot f_1(x_0,x_1) \}$.\\

Let $\AAN{0} =\{ x_2^2x_3, x_3^2x_0,x_4^2x_1,x_5^2x_1, x_0^3, x_1^3 \}$. Then 
$$
c_1=-2c_4,\; c_3=2c_2,\; 4c_2=6c_4=2(c_5-c_4)=0. 
$$

For $\AAN{1}=\AAN{0}\cup \{ x_0x_4x_5 \}$ we get an extra condition $c_4+c_5=0$, i.e.
$$
c_1=0,\; c_3=2c_2,\; c_5=c_4,\; 4c_2=2c_4=0. 
$$

For $\AA=\AAN{1}\cup \{ x_3x_4x_5 \}$ we get an extra condition $c_3=0$, i.e.
$$
c_1=c_3=0,\; c_5=c_4,\; 2c_2=2c_4=0. 
$$
This case was already considered.\\

For $\AA=\AAN{1}\cup \{ x_2x_3x_5 \}$ we get an extra condition $c_5=c_2$, i.e.
$$
c_1=c_3=0,\; c_5=c_4=c_2,\; 2c_2=0. 
$$
This case was already considered.\\

For $\AA=\AAN{0}\cup \{ x_2x_4x_5 \}$ we get an extra condition $c_2+c_4+c_5=0$, i.e.
$$
c_1=c_3=0,\; c_5=c_2+c_4,\; 2c_2=2c_4=0. 
$$
This case was already considered.\\

For $\AAN{1}=\AAN{0}\cup \{ x_3x_4x_5 \}$ we get an extra condition $c_3+c_4+c_5=0$, i.e.
$$
c_1=0,\; c_3=2c_2,\; c_5=c_4+2c_2,\; 4c_2=2c_4=0. 
$$

For $\AA=\AAN{1}\cup \{ x_2^2x_0 \}$ we get an extra condition $2c_2=0$, i.e.
$$
c_1=c_3=0,\; c_5=c_4,\; 2c_2=2c_4=0. 
$$
This case was already considered.\\

For $\AA=\AAN{1}\cup \{ x_0x_2x_4 \}$ we get an extra condition $c_2+c_4=0$, i.e.
$$
c_1=c_3=0,\; c_5=c_4=c_2,\; 2c_2=0. 
$$
This case was already considered.\\

For $\AA=\AAN{1}\cup \{ x_2^2x_4 \}$ we get an extra condition $c_4=2c_2$, i.e.
$$
c_1=c_5=0,\; c_3=c_4=2c_2,\; 4c_2=0.
$$
This case will appear and will be analyzed below.\\

Let $\AAN{0} =\{ x_2^2x_3, x_3^2x_0, x_4^2x_3, x_5^2x_0, x_2x_4x_5, x_0^3, x_1^3 \}$. Then $c_3=2c_2$, $3c_1=4c_2=2(c_4-c_2)=2c_5=0$, $c_2+c_4+c_5=0$, i.e. 
$$
c_3=2c_2,\; c_4=c_5-c_2,\; 3c_1=4c_2=2c_5=0. 
$$

For $\AA=\AAN{0}\cup \{ x_5^2x_3 \}$ we get an extra condition $2c_2=0$, i.e.
$$
c_3=0,\; c_4=c_2+c_5,\; 3c_1=2c_2=2c_5=0.
$$
This case was already considered.\\

For $\AA=\AAN{0}\cup \{ x_1x_3x_5 \}$ we get an extra condition $c_1+c_3+c_5=0$, i.e.
$$
c_1=0,\; c_3=c_5=2c_2,\; c_4=c_2,\; 4c_2=0. 
$$
This case will appear and will be analyzed below.\\

For $\AA=\AAN{0}\cup \{ x_2x_3x_5 \}$ we get an extra condition $c_2+c_3+c_5=0$, i.e.
$$
c_3=c_4=0,\; c_5=c_2,\; 3c_1=2c_2=0. 
$$
This case has already appeared.\\

Let $\AAN{0} =\{ x_2^2x_3, x_3^2x_0, x_4^2x_3, x_5^2x_0, x_1x_2x_4, x_0^3, x_1^3 \}$. Then $c_3=2c_2$, $3c_1=4c_2=2(c_4-c_2)=2c_5=0$, $c_1+c_2+c_4=0$, i.e. 
$$
c_1=0,\; c_3=2c_2,\; c_4=-c_2,\; 4c_2=2c_5=0. 
$$

For $\AA=\AAN{0}\cup \{ x_5^2x_3 \}$ and $\AA=\AAN{0}\cup \{ x_2x_3x_5 \}$ we get extra conditions $2c_2=0$ and $c_2=c_5$ respectively. They both imply that $2c_2=0$, i.e.
$$
c_1=c_3=0,\; c_4=c_2,\; 2c_2=2c_5=0. 
$$
All such cases were already considered.\\

For $\AA=\AAN{0}\cup \{ x_1x_3x_5 \}$ we get an extra condition $c_5=2c_2$, i.e.
$$
c_1=0,\; c_3=c_5=2c_2,\; c_4=-c_2,\; 4c_2=0. 
$$
This means that $\GG\cong \ZZ{4}$ with a generator
$$
f=(1,1,{\om},{\om}^2,{\om}^{-1},{\om}^2),\; \om =\sqr{4}, 
$$
Then $i_2\equiv i_4+ 2(i_3+i_5) \; mod \; 4$. This means that  $\AAbar=\AA \cup $ $\{  f_3(x_0,x_1),$ $ x_5^2x_1,$ $ x_0x_3x_5,$ $ x_4^2x_5,$ $ x_2^2x_5,$ $x_3^2x_1,$ $ x_0x_2x_4 \}$.\\

Let $\AAN{0} =\{ x_2^2x_3, x_3^2x_0, x_4^2x_3, x_5^2x_0, x_1x_4x_5, x_0^3, x_1^3 \}$. Then $c_3=2c_2$, $3c_1=4c_2=2(c_4-c_2)=2c_5=0$, $c_1+c_4+c_5=0$, i.e. 
$$
c_1=c_3=0,\; c_5=c_4,\; 2c_2=2c_4=0. 
$$
All such cases were already considered.\\

Let $\AAN{0} =\{ x_2^2x_3, x_3^2x_0,x_4^2x_3,x_5^2x_3, x_0x_2x_4, x_1x_2x_4, x_0^3, x_1^3 \}$. Then 
$$
c_1=0,\; c_3=2c_2,\; c_4=-c_2,\; 4c_2=2(c_5-c_2)=0. 
$$

For $\AA =\AAN{0}\cup \{ x_1x_2x_5, x_1x_4x_5 \}$ we get extra conditions $c_2=c_4=-c_5$, i.e.
$$
c_1=c_3=0,\; c_4=c_5=c_2,\; 2c_2=0.
$$
This case was already considered.\\

Let $\AAN{0} =\{ x_2^2x_3, x_3^2x_0,x_4^2x_3,x_5^2x_3, x_0x_2x_4, x_1x_4x_5, x_0^3, x_1^3 \}$. Then 
$$
c_1=0,\; c_3=2c_2,\; c_4=-c_2, \; c_5=c_2,\; 4c_2=0. 
$$

For $\AA=\AAN{0}\cup \{ x_0x_2x_5 \}$ and $\AA=\AAN{0}\cup \{ x_1x_2x_5 \}$ we get the same extra condition $c_2+c_5=0$, i.e.
$$
c_1=c_3=0,\; c_4=c_5=c_2,\; 2c_2=0. 
$$
This case was already considered.\\

If $\AA =\{ x_2^2x_3, x_3^2x_0,x_4^2x_3,x_5^2x_3, x_2x_4x_5, x_0^3, x_1^3 \}$. Then 
$$
c_3=0,\; c_5=c_2+c_4,\; 3c_1=2c_2=2c_4=0. 
$$
This case was already considered.\\

Let $\AAN{0} =\{ x_2^2x_3, x_3^2x_0,x_4^2x_0,x_5^2x_0, x_0^3, x_1^3 \}$. Then 
$$
c_3=2c_2,\; 3c_1=4c_2=2c_4=2c_5=0. 
$$

For $\AAN{1}=\AAN{0}\cup \{ x_3x_4x_5 \}$ we get an extra condition $c_3+c_4+c_5=0$, i.e.
$$
c_3=2c_2,\; c_5=2c_2+c_4,\; 3c_1=4c_2=2c_4=0.
$$

For $\AA=\AAN{1}\cup \{ x_2^2x_1 \}$, $\AA=\AAN{1}\cup \{ x_1x_2x_4 \}$, $\AA=\AAN{1}\cup \{ x_1x_2x_5 \}$, $\AA=\AAN{1}\cup \{ x_1x_4x_5 \}$ we get extra conditions $c_1=2c_2$, $c_1+c_2+c_4=0$, $c_1+c_2+c_5=0$, $c_1=2c_2$ respectively. They all imply that $c_1=2c_2=0$, i.e.
$$
c_1=c_3=0,\; c_5=c_4,\; 2c_2=2c_4=0.
$$
This case has already appeared.\\

For $\AA=\AAN{1}\cup \{ x_2x_4x_5 \}$ we get an extra condition $c_2=0$, i.e.
$$
c_2=c_3=0,\; c_5=c_4,\; 3c_1=2c_4=0.
$$
This case has already appeared.\\

For $\AA=\AAN{1}\cup \{ x_2^2x_4 \}$ and $\AA=\AAN{1}\cup \{ x_2^2x_5 \}$ we get extra conditions $c_4=2c_2$ and $c_4=0$ respectively, i.e. upto a permutation of $x_i$ we get
$$
c_3=c_4=2c_2,\; c_5=0,\; 3c_1=4c_2=0.
$$
This case has already appeared.\\

For $\AAN{1}=\AAN{0}\cup \{ x_2x_4x_5, x_1x_3x_5 \}$ we get extra conditions $c_2+c_4+c_5=c_1+c_3+c_5=0$, i.e.
$$
c_1=c_3=c_5=0,\; c_4=c_2,\; 2c_2=0. 
$$
All such cases were already considered.\\

For $\AAN{1}=\AAN{0}\cup \{ x_2x_3x_5, x_1x_4x_5 \}$ we get extra conditions $c_2+c_3+c_5=0$ and $c_1+c_4+c_5=0$, i.e.
$$
c_1=c_3=0,\; c_2=c_4=c_5,\; 2c_2=0. 
$$
All such cases were already considered.\\

Let $\AAN{0} =\{ x_2^2x_0, x_3^2x_0, x_3^2x_1, x_4^2x_1,x_5^2x_1, x_0^3, x_1^3 \}$. Then 
$$
c_1=0,\; 2c_2=2c_3=2c_4=2c_5=0. 
$$

For $\AA=\AAN{0}\cup \{ x_3x_4x_5 \}$ we get an extra condition $c_3+c_4+c_5=0$, i.e.
$$
c_1=0,\; c_3=c_4+c_5,\; 2c_2=2c_4=2c_5=0.
$$
This means that $\GG\cong (\ZZ{2})^{\oplus 3}$ with generators
$$
f_1=(1,1,{\om},1,1,1),\; f_2=(1,1,1,{\om},\om,1),\; f_3=(1,1,1,{\om},1,\om),\; \om =\sqr{2}, 
$$
Then $i_2\equiv i_3+i_4 \equiv i_3+i_5 \equiv 0 \; mod \; 2$. This means that  $\AAbar=\AA \cup \{  f_3(x_0,x_1), x_2^2\cdot f_1(x_0,x_1), x_3^2\cdot f_1(x_0,x_1), x_4^2\cdot f_1(x_0,x_1), x_5^2\cdot f_1(x_0,x_1) \}$.\\

For $\AA=\AAN{0}\cup \{ x_2x_4x_5 \}$ we get an extra condition $c_2+c_4+c_5=0$, i.e.
$$
c_1=0,\; c_2=c_4+c_5,\; 2c_3=2c_4=2c_5=0.
$$
This case was already considered.\\

If $\AA =\{ x_2^2x_0, x_3^2x_0, x_4^2x_1, x_5^2x_1, x_2x_4x_5, x_1x_2x_3, x_0^3, x_1^3 \}$. Then 
$$
c_1=0,\; c_3=c_2=c_4+c_5,\; 2c_4=2c_5=0. 
$$
This means that $\GG\cong (\ZZ{2})^{\oplus 2}$ with generators
$$
f_1=(1,1,{\om},\om,\om,1),\; f_2=(1,1,\om,{\om},1,\om),\; \om =\sqr{2}, 
$$
Then $i_2+i_3+i_4\equiv i_2+i_3+i_5 \equiv 0 \; mod \; 2$. This means that  $\AAbar=\AA \cup \{  f_3(x_0,x_1), x_4^2\cdot f_1(x_0,x_1), x_5^2\cdot f_1(x_0,x_1), f_2(x_2,x_3) \cdot f_1(x_0,x_1), x_3x_4x_5 \}$.\\

If $\AA =\{ x_2^2x_0, x_3^2x_0, x_4^2x_1, x_5^2x_1, x_2x_4x_5, x_2x_3x_4, x_0^3, x_1^3 \}$. Then 
$$
c_1=0,\; c_3=c_5=c_2+c_4,\; 2c_2=2c_4=0. 
$$
This case was already considered.\\

Let $\AAN{0} =\{ x_2^2x_0, x_3^2x_0, x_4^2x_0,x_5^2x_1, x_0^3, x_1^3 \}$. Then 
$$
c_1=-2c_5,\; 2c_2=2c_3=2c_4=6c_5=0. 
$$

For $\AA=\AAN{0}\cup \{ x_2x_3x_4 \}$ we get an extra condition $c_2+c_3+c_4=0$, i.e.
$$
c_1=-2c_5,\; c_3=c_2+c_4,\; 2c_2=2c_4=6c_5=0.
$$
This means that $\GG\cong (\ZZ{2})^{\oplus 2}\oplus \ZZ{6}$ with generators
$$
f_1=(1,1,{\om},\om,1,1),\; f_2=(1,1,1,{\om},\om,1),\; f_3=(1,{\eta}^{-2},1,1,1,\eta),\; \om =\sqr{2}, \eta =\sqr{6}, 
$$
Then $i_2+i_3\equiv i_3+i_4 \equiv 0 \; mod \; 2$, $i_5 \equiv 2i_1 \; mod \; 6$. This means that  $\AAbar=\AA $.\\

For $\AAN{1}=\AAN{0}\cup \{ x_1x_3x_4, x_2x_4x_5 \}$ we get extra conditions $c_1+c_3+c_4=c_2+c_4+c_5=0$, i.e.
$$
c_1=0,\; c_3=c_4=c_2+c_5,\; 2c_2=2c_5=0.
$$
All such cases have already appeared.\\

Let $\AAN{0} =\{ x_2^2x_0, x_3^2x_0, x_4^2x_0,x_5^2x_0, x_2x_3x_4, x_1x_2x_5, x_0^3, x_1^3 \}$. Then 
$$
c_1=0,\; c_5=c_2=c_3+c_4,\; 2c_3=2c_4=0. 
$$

All such cases have already appeared.\\

If $\AA =\{ x_2^2x_0, x_3^2x_0, x_4^2x_0, x_5^2x_0, x_2x_3x_4, x_2x_3x_5, x_2x_4x_5, x_3x_4x_5, x_0^3, x_1^3 \}$. Then 
$$
c_2=c_3=c_4=c_5=0,\; 3c_1=0. 
$$
This case was already considered.\\

\subsection{Case of $1$ cube.}

Let the cube be $x_0^3$. If a longest cycle has length $5$, then $\AA$ is
\begin{center}

\end{center}

i.e. $\AA=\{ x_0^3, x_1^2x_2, x_2^2x_3, x_3^2x_4, x_4^2x_5, x_5^2x_1 \}$. Then $2c_1+c_2=2c_2+c_3=2c_3+c_4=2c_4+c_5=2c_5+c_1=0$, i.e.
$$
c_2=-2c_1,\; c_3=4c_1,\; c_4=-8c_1,\; c_5=16c_1,\; 33c_1=0.
$$

This means that $\GG\cong \ZZ{33}$ with a generator 
$$
f=(1,{\om},{\om}^{-2},{\om}^{4},{\om}^{-8},{\om}^{16}), \; \om=\sqr{33}.
$$

Then $i_1-2i_2+4i_3-8i_4+16i_5\equiv 0 \; mod \; 33$. This means that $\AAbar=\AA$.\\

\subsubsection{Case of $1$ cube. Length $4$ longest cycle.}

If a longest cycle has length $4$, then $\AA$ is either

\begin{center}
%
\end{center}

Indeed, suppose $x_1, x_2, x_3, x_4$ form a length $4$ cycle. By Lemma 1 (\cite{Liendo}, Lemma 1.3) the $5$-th vertex should be connected either to the cube $0$ (in which case there are neither singular pairs nor singular triples) or to the cycle (say, to the $1$-st vertex).\\ 

In the latter case we may have one singular pair $(x_4, x_5)$ and one singular triple $(x_2,x_4,x_5)$. They are resolved as in the earlier examples.\\

This gives us the pictures shown above. Now let us do computations.\\

If $\AA=\{ x_1^2x_2, x_2^2x_3, x_3^2x_4, x_4^2x_1, x_5^2x_0, x_0^3 \}$. Then 
$$
c_2=-2c_1,\; c_3=4c_1,\; c_4=-8c_1,\; 15c_1=2c_5=0.
$$

This means that $\GG\cong \ZZ{2}\oplus \ZZ{15}$ with generators 
$$
f_1=(1,{\om},{\om}^{-2},{\om}^{4},{\om}^{-8},1),\; f_2=(1,1,1,1,1,\eta), \; \om=\sqr{15}, \eta=\sqr{2}.
$$

Then $i_1-2i_2+4i_3-8i_4 \equiv 0 \; mod \; 15$, $i_5 \equiv 0 \; mod \; 2$. This means that $\AAbar=\AA$.\\

Let $\AAN{0} =\{ x_1^2x_2, x_2^2x_3, x_3^2x_4, x_4^2x_1, x_5^2x_1, x_0^3 \}$. Then 
$$
c_1=-2c_5,\; c_2=4c_5,\; c_3=-8c_5,\; c_4=16c_5,\; 30c_5=0.
$$

For $\AA=\AAN{0}\cup \{ x_4^2x_2 \}$, $\AAN{1}=\AAN{0}\cup \{ x_4^2x_3 \}$, $\AA=\AAN{0}\cup \{ x_5^2x_2 \}$ and $\AAN{1}=\AAN{0}\cup \{ x_5^2x_3 \}$ we get the same extra condition $6c_5=0$, i.e.
$$
c_1=c_2=c_3=c_4=-2c_5,\; 6c_5=0.
$$
This means that $\GG\cong \ZZ{6}$ with a generator
$$
f=(1,{\om}^{-2},{\om}^{-2},{\om}^{-2},{\om}^{-2},\om),\; \om =\sqr{6}, 
$$
Then $i_5\equiv 2(i_1+i_2+i_3+i_4) \; mod \; 6$. This means that  $\AAbar=\AA \cup \{  f_3(x_1,x_2,x_3,x_4),x_5^2\cdot f_1(x_1,x_2,x_3,x_4) \}$.\\

For $\AA=\AAN{1}\cup \{ x_2^2x_0 \}$, $\AA=\AAN{1}\cup \{ x_0x_2x_4 \}$, $\AA=\AAN{1}\cup \{ x_0x_2x_5 \}$ we get extra conditions $2c_5=0$, $2c_5=0$, $c_5=0$ respectively. They all imply that $2c_5=0$, i.e.
$$
c_1=c_2=c_3=c_4=0,\; 2c_5=0.
$$
All such cases have already appeared.\\

For $\AA=\AAN{1}\cup \{ x_2^2x_4 \}$ the extra condition is automatically satisfied.\\

This case was already considered.\\

For $\AA=\AAN{1}\cup \{ x_2^2x_5 \}$ we get an extra condition $3c_5=0$, i.e.
$$
c_1=c_2=c_3=c_4=c_5,\; 3c_5=0.
$$
This means that $\GG\cong \ZZ{3}$ with a generator
$$
f=(1,{\om},{\om},{\om},{\om},\om),\; \om =\sqr{3}, 
$$
Then $i_1+i_2+i_3+i_4+i_5\equiv 0 \; mod \; 3$. This means that  $\AAbar=\AA \cup \{  f_3(x_1,x_2,x_3,x_4, x_5) \}$.\\

For $\AA=\AAN{0}\cup \{ x_4^2x_0 \}$ we get an extra condition $2c_5=0$, i.e.
$$
c_1=c_2=c_3=c_4=0,\; 2c_5=0.
$$
This case was already considered.\\

For $\AA=\AAN{0}\cup \{ x_0x_4x_5 \}$ and $\AA=\AAN{0}\cup \{ x_2x_4x_5 \}$ we get the same extra condition $c_5=0$, i.e.
$$
c_1=c_2=c_3=c_4=c_5=0.
$$
Hence $\GG = Id$.\\

For $\AAN{1}=\AAN{0}\cup \{ x_3x_4x_5 \}$ we get an extra condition $c_3+c_4+c_5=0$, i.e.
$$
c_1=c_2=c_3=c_4=c_5,\; 3c_5=0.
$$
All such cases have already appeared.\\

\subsubsection{Case of $1$ cube. Length $3$ longest cycle.}

Suppose that a longest cycle has length $3$. We may assume that $1$ is connected to $2$, $2$ is connected to $3$ and $3$ is connected to $1$.\\

If the remaining vertices $4$ and $5$ are connected, then $\AA$ is one of the following:

\begin{center}

\end{center}

Indeed, let $4$ be connected to $5$. If there is also an edge from $5$ to $4$ (i.e. $4$ and $5$ form a cycle), then there are neither singular pairs nor singular triples.\\ 

Let us assume that $4$ and $5$ do not form a cycle. Then by Lemma 1 (\cite{Liendo}, Lemma 1.3) the $5$-th vertex should be connected to either the cube $0$ or to the length $3$ cycle (say, to $3$).\\

In the former case there are neither singular pairs nor singular triples. In the latter case $(x_2,x_5)$ may form a singular pair. We resolve it as usual.\\

Now let us assume that vertices $4$ and $5$ are disconnected. By Lemma 1 (\cite{Liendo}, Lemma 1.3) each of them should be connected either to the cube $0$ or to the cycle.\\

Let us consider the case when one of the vertices $4$, $5$ is connected to the cube (say, $5$ is connected to $0$) and the other one is connected to the cycle (say, $4$ is connected to $1$).\\

Then $\AA$ will be one of the following:

\begin{center}


\end{center}

\begin{flalign*}
\mbox{where}\; & A=(x_5^2x_3,\; x_5^2x_2,\; x_2x_3x_5,\; x_2x_4x_5).\\
\end{flalign*}

Indeed, in this case we may get
\begin{itemize}
\item one singular pair $(x_3,x_4)$ and 
\item one singular triple $(x_3,x_4,x_5)$.
\end{itemize}

We resolve them as in the earlier examples.\\

Suppose that both vertices $4$ and $5$ are connected to the cycle formed by $1$, $2$, $3$. If they are connected to two different vertices of that cycle (say, $4$ is connected to $2$ and $5$ is connected to $3$), then we may get two singular pairs $(x_1, x_4)$ and $(x_2, x_5)$ and one singular triple $(x_1, x_4, x_5)$.\\

We resolve them as usual and obtain the following pictures:

\begin{center}

\end{center}

\begin{flalign*}
\mbox{where}\; & A=(x_4^2x_1,\; x_1x_4x_5,\; x_1^2x_3,\; x_4^2x_3,\; x_1x_3x_4,\; x_1^2x_0,\; x_0x_1x_4);\\
& B=(x_1x_4x_5,\; x_5^2x_1,\; x_4^2x_1,\; x_1^2x_0,\; x_0x_4x_5,\; x_0x_1x_4,\; x_0x_1x_5).\\
\end{flalign*}

Finally, it may happen that both vertices $4$ and $5$ are connected to the same vertex of the cycle (say, to $1$) or to the cube $0$.\\

In this case $\AA$ is one of the following:

\begin{center}
%
\end{center}

Indeed, if $4$ and $5$ are both connected to $1$, then we may get 
\begin{itemize}
\item three singular pairs $(x_3,x_4)$, $(x_3,x_5)$, $(x_4,x_5)$ and
\item two singular triples $(x_3,x_4, x_5)$, $(x_2, x_4,x_5)$.
\end{itemize}

Note that this configuration is symmetric with respect to vertices $4$ and $5$.\\

If $4$ and $5$ are both connected to the cube $0$, then we may get 
\begin{itemize}
\item one singular pair $(x_4,x_5)$ and
\item three singular triples $(x_1,x_4, x_5)$, $(x_2, x_4,x_5)$, $(x_3, x_4,x_5)$.
\end{itemize}

Note that this configuration is symmetric with respect to vertices $4$ and $5$ as well as with respect to a cyclic permutation of vertices $1$, $2$, $3$.\\

We resolve all these singular pairs and triples as usual and obtain the pictures shown above.\\

Now let us do computations.\\

If $\AA=\{ x_1^2x_2, x_2^2x_3, x_3^2x_1, x_4^2x_5, x_5^2x_0, x_0^3 \}$, then 
$$
c_2=-2c_1,\; c_3=4c_1,\; c_5=2c_4,\; 9c_1=4c_4=0.
$$

This means that $\GG\cong \ZZ{4}\oplus \ZZ{9}$ with generators 
$$
f_1=(1,{\om},{\om}^{-2},{\om}^{4},1,1),\; f_2=(1,1,1,1,\eta,{\eta}^2), \; \om=\sqr{9}, \eta=\sqr{4}.
$$

Then $i_1-2i_2+4i_3 \equiv 0 \; mod \; 9$, $i_4 \equiv 2i_5 \; mod \; 4$. This means that $\AAbar=\AA$.\\

If $\AA=\{ x_1^2x_2, x_2^2x_3, x_3^2x_1, x_4^2x_5, x_5^2x_4, x_0^3 \}$, then 
$$
c_2=-2c_1,\; c_3=4c_1,\; c_5=c_4,\; 9c_1=3c_4=0.
$$

This means that $\GG\cong \ZZ{3}\oplus \ZZ{9}$ with generators 
$$
f_1=(1,{\om},{\om}^{-2},{\om}^{4},1,1),\; f_2=(1,1,1,1,\eta,{\eta}), \; \om=\sqr{9}, \eta=\sqr{3}.
$$

Then $i_1-2i_2+4i_3 \equiv 0 \; mod \; 9$, $i_4+i_5 \equiv 0 \; mod \; 3$. This means that $\AAbar=\AA\cup \{ f_3(x_4,x_5) \}$.\\

Let $\AAN{0} =\{ x_1^2x_2, x_2^2x_3, x_3^2x_1, x_4^2x_5, x_5^2x_3, x_0^3 \}$. Then 
$$
c_2=-2c_1,\; c_3=4c_1,\; c_5=-2c_4,\; 9c_1=4(c_1-c_4)=0.
$$

For $\AA=\AAN{0}\cup \{ x_2^2x_0 \}$ we get an extra condition $c_1=0$, i.e.
$$
c_1=c_2=c_3=0,\; c_5=2c_4,\; 4c_4=0.
$$
This means that $\GG\cong \ZZ{4}$ with a generator 
$$
f=(1,1,1,1,\om,{\om}^2), \; \om=\sqr{4}.
$$

Then $i_4 \equiv 2i_5 \; mod \; 4$. This means that $\AAbar=\AA\cup \{ f_3(x_0,x_1,x_2,x_3), x_5^2\cdot f_1(x_0,x_1,x_2,x_3) \}$.\\

For $\AA=\AAN{0}\cup \{ x_5^2x_2 \}$ and $\AA=\AAN{0}\cup \{ x_2^2x_1 \}$ we get the same extra condition $c_2+2c_5=c_1+2c_2=0$, i.e.
$$
c_1=c_2=c_3=4c_4,\; c_5=-2c_4,\; 12c_4=0.
$$
This means that $\GG\cong \ZZ{12}$ with a generator 
$$
f=(1,{\om}^{4},{\om}^{4},{\om}^{4},\om,{\om}^{-2}), \; \om=\sqr{12}.
$$

Then $i_4+4(i_1+i_2+i_3) \equiv 2i_5 \; mod \; 12$. This means that $\AAbar=\AA\cup \{ f_3(x_1,x_2,x_3), x_5^2\cdot f_1(x_1,x_2,x_3) \}$.\\

For $\AA=\AAN{0}\cup \{ x_0x_2x_5 \}$ we get an extra condition $c_2+c_5=0$, i.e.
$$
c_1=c_2=c_3=0,\; c_5=2c_4,\; 4c_4=0.
$$
This case was already considered.\\

For $\AA=\AAN{0}\cup \{ x_2x_4x_5 \}$ we get an extra condition $c_2+c_4+c_5=0$, i.e.
$$
c_1=c_2=c_3=c_4=c_5,\; 3c_1=0.
$$
This case was already considered.\\

For $\AA=\AAN{0}\cup \{ x_1x_2x_5 \}$ we get an extra condition $c_1+c_2+c_5=0$, i.e.
$$
c_1=c_2=c_3=c_5=-2c_4,\; 6c_4=0.
$$
This case was already considered.\\

Let $\AAN{0} =\{ x_1^2x_2, x_2^2x_3, x_3^2x_1, x_4^2x_1, x_5^2x_0, x_0^3 \}$. Then 
$$
c_1=-2c_4,\; c_2=4c_4,\; c_3=-8c_4,\; 18c_4=2c_5=0.
$$

For $\AA=\AAN{0}\cup \{ x_4^2x_2 \}$ and $\AA=\AAN{0}\cup \{ x_3^2x_2 \}$ we get the same extra condition $6c_4=0$, i.e.
$$
c_1=c_2=c_3=-2c_4,\; 6c_4=2c_5=0.
$$

This case was already considered.\\

For $\AA=\AAN{0}\cup \{ x_2x_3x_4 \}$ we get an extra condition $c_2+c_3+c_4=0$, i.e.
$$
c_1=c_2=c_3=c_4,\; 3c_4=2c_5=0.
$$
This case was already considered.\\

For $\AA=\AAN{0}\cup \{ x_3^2x_5 \}$ we get an extra condition $2c_3+c_5=0$, i.e.
$$
c_1=c_2=c_3=c_5=0,\; 2c_4=0.
$$
This case was already considered.\\

For $\AA=\AAN{0}\cup \{ x_3x_4x_5 \}$ we get an extra condition $c_5=7c_4$, i.e.
$$
c_1=c_2=c_3=0,\; c_5=c_4,\; 2c_4=0.
$$
This case was already considered.\\

For $\AAN{1}=\AAN{0}\cup \{ x_0x_3x_4 \}$, $\AAN{1}=\AAN{0}\cup \{ x_4^2x_0 \}$ and $\AAN{1}=\AAN{0}\cup \{ x_3^2x_0 \}$ we get extra conditions $c_4=0$, $2c_4=0$ and $2c_4=0$ respectively. They all imply that $2c_4=0$, i.e.
$$
c_1=c_2=c_3=0,\; 2c_4=2c_5=0.
$$
All such cases were already considered.\\

Let $\AAN{0} =\{ x_1^2x_2, x_2^2x_3, x_3^2x_1, x_4^2x_2, x_5^2x_2, x_5^2x_3, x_0^3 \}$ or $\AAN{0} =\{ x_1^2x_2,$ $x_2^2x_1,$ $x_2^2x_3,$ $x_3^2x_1,$ $x_4^2x_2,$ $x_5^2x_3,$ $x_0^3 \}$. Then 
$$
c_1=c_2=c_3=-2c_4,\; 6c_4=2(c_5-c_4)=0.
$$

For $\AA=\AAN{0}\cup \{ x_4^2x_1 \}$ the extra condition is automatically satisfied.\\

This means that $\GG\cong \ZZ{2}\oplus \ZZ{6}$ with generators 
$$
f_1=(1,{\om}^{-2},{\om}^{-2},{\om}^{-2},\om,\om),\; f_2=(1,1,1,1,1,{\eta}), \; \om=\sqr{6}, \eta=\sqr{2}.
$$

Then $i_5 \equiv 0 \; mod \; 2$, $i_4+i_5 \equiv 2(i_1+i_2+i_3) \; mod \; 6$. This means that $\AAbar=\AA\cup \{ f_3(x_1,x_2,x_3), x_4^2\cdot f_1(x_1,x_2,x_3), x_5^2\cdot f_1(x_1,x_2,x_3) \}$.\\

For $\AA=\AAN{0}\cup \{ x_1^2x_0 \}$ and $\AA=\AAN{0}\cup \{ x_0x_1x_4 \}$ we get extra conditions $2c_4=0$ and $c_4=0$ respectively. They both immply that $2c_4=0$, i.e.
$$
c_1=c_2=c_3=0,\; 2c_4=2c_5=0.
$$
All such cases were already considered.\\

For $\AA=\AAN{0}\cup \{ x_1x_4x_5 \}$ we get an extra condition $c_5=c_4$, i.e.
$$
c_1=c_2=c_3=-2c_4,\; c_5=c_4,\; 6c_4=0.
$$
This means that $\GG\cong \ZZ{6}$ with a generator 
$$
f=(1,{\om}^{-2},{\om}^{-2},{\om}^{-2},\om,\om), \; \om=\sqr{6}.
$$

Then $i_4+i_5 \equiv 2(i_1+i_2+i_3) \; mod \; 6$. This means that $\AAbar=\AA\cup \{ f_3(x_1,x_2,x_3), f_2(x_4,x_5) \cdot f_1(x_1,x_2,x_3) \}$.\\

Any other extra condition (after adding a monomial $x_0^{i_0}x_1^{i_1}x_2^{i_2}x_3^{i_3}x_4^{i_4}x_5^{i_5}$ to $\AAN{0}$) will have the form $a\cdot c_4+b\cdot c_5=0$, $a,b\in \mathbb Z$, i.e. either $a\cdot c_4=0$ or $c_5=a\cdot c_4$ with some $a\in \{ 0,1,2,3,4,5 \}$. All such cases have already appeared.\\

Let $\AAN{0} =\{ x_1^2x_2, x_2^2x_3, x_3^2x_1, x_4^2x_2, x_5^2x_3, x_2x_4x_5, x_0^3 \}$. Then 
$$
c_1=c_2=c_3=-2c_4,\; c_5=c_4,\; 6c_4=0.
$$
All such cases have already appeared.\\

Let $\AAN{0} =\{ x_1^2x_2, x_2^2x_3, x_3^2x_1, x_4^2x_2, x_2^2x_4, x_5^2x_3, x_0^3 \}$. Then 
$$
c_2=c_4=-2c_1,\; c_3=4c_1,\; 9c_1=2(c_5+2c_1)=0.
$$

For $\AA=\AAN{0}\cup \{ x_4^2x_1 \}$, $\AA=\AAN{0}\cup \{ x_1x_4x_5 \}$, $\AA=\AAN{0}\cup \{ x_1^2x_0 \}$, $\AA=\AAN{0}\cup \{ x_0x_1x_4 \}$, $\AAN{1}=\AAN{0}\cup \{ x_1^2x_3 \}$, $\AAN{1}=\AAN{0}\cup \{ x_1x_3x_4 \}$ we get extra conditions $3c_1=0$, $c_5=c_1$, $2c_1=0$, $c_1=0$, $3c_1=0$, $3c_1=0$ respectively. They all imply that $3c_1=0$, i.e.
$$
c_1=c_2=c_3=c_4=-2c_5,\; 6c_5=0.
$$
All such cases have already appeared.\\

For $\AAN{1}=\AAN{0}\cup \{ x_4^2x_3 \}$ the extra condition is automatically satisfied.\\

For $\AA=\AAN{1}\cup \{ x_5^2x_1 \}$, $\AA=\AAN{1}\cup \{ x_0x_4x_5 \}$, $\AA=\AAN{1}\cup \{ x_0x_1x_5 \}$ we get extra conditions $3c_1=0$, $c_5=2c_1$, $c_5=-c_1$ respectively. They all imply that $3c_1=0$, i.e.
$$
c_1=c_2=c_3=c_4=-2c_5,\; 6c_5=0.
$$
All such cases have already appeared.\\

Let $\AAN{0} =\{ x_1^2x_2, x_2^2x_3, x_3^2x_1, x_4^2x_2, x_5^2x_1, x_5^2x_3, x_0^3 \}$. Then 
$$
c_1=c_2=c_3=-2c_5,\; 2(c_4-c_5)=6c_5=0.
$$
All such cases were already considered.\\

Let $\AAN{0} =\{ x_1^2x_2, x_2^2x_3, x_3^2x_1, x_4^2x_2, x_5^2x_3, x_1x_2x_5, x_0^3 \}$. Then 
$$
c_1=c_2=c_3=c_5=-2c_4,\; 6c_4=0.
$$
All such cases were already considered.\\

Let $\AAN{0} =\{ x_1^2x_2, x_2^2x_3, x_3^2x_1, x_4^2x_2, x_5^2x_3, x_2^2x_0, x_0^3 \}$. Then 
$$
c_1=c_2=c_3=0,\; 2c_4=2c_5=0.
$$
All such cases have already appeared.\\

Let $\AAN{0} =\{ x_1^2x_2, x_2^2x_3, x_3^2x_1, x_4^2x_2, x_5^2x_3, x_0x_2x_5, x_0^3 \}$. Then 
$$
c_1=c_2=c_3=c_5=0,\; 2c_4=0.
$$
All such cases have already appeared.\\

Let $\AAN{0} =\{ x_1^2x_2, x_2^2x_3, x_3^2x_1, x_4^2x_0, x_5^2x_0, x_2x_4x_5, x_0^3 \}$ or $\AAN{0} =\{ x_1^2x_2,$ $x_2^2x_3,$ $x_3^2x_1,$ $x_4^2x_1,$ $x_5^2x_1,$ $x_0x_4x_5,$ $x_0^3 \}$. Then 
$$
c_1=c_2=c_3=0,\; c_4=c_5,\; 2c_5=0.
$$
All such cases were already considered.\\

Let $\AAN{0} =\{ x_1^2x_2, x_2^2x_3, x_3^2x_1, x_4^2x_1, x_5^2x_1, x_3x_4x_5, x_0^3 \}$. Then 
$$
c_1=c_2=c_3=-2c_4,\; c_5=c_4,\; 6c_4=0.
$$
All such cases were already considered.\\

If $\AA =\{ x_1^2x_2, x_2^2x_3, x_3^2x_1, x_4^2x_1, x_5^2x_1, x_3^2x_0, x_2x_4x_5, x_0^3 \}$. Then 
$$
c_1=c_2=c_3=0,\; c_4=c_5,\; 2c_5=0.
$$
This case was already considered.\\

Let $\AAN{0} =\{ x_1^2x_2, x_2^2x_3, x_3^2x_1, x_4^2x_1, x_5^2x_1, x_2x_4x_5, x_0x_3x_5, x_0^3 \}$. Then 
$$
c_1=c_2=c_3=c_4=c_5=0.
$$
Hence $\GG = Id$ for all such cases.\\

\subsubsection{Case of $1$ cube. Length $2$ longest cycle.}

If a longest cycle has length $2$, then $\AA$ may be either

\begin{center}

\end{center}

or one of the following:

\begin{center}
%

\begin{flalign*}
\mbox{where}\; & A=(x_2x_3x_5,\; x_3^2x_2,\; x_5^2x_2,\; x_2^2x_0,\; x_0x_2x_3,\; x_0x_2x_5);\\
& B=(x_2x_3x_5,\; x_0x_2x_5);\\
& C=(x_2x_3x_5,\; x_2^2x_3,\; x_2^2x_5,\; x_2^2x_0,\; x_0x_2x_3,\; x_0x_2x_5).\\
\end{flalign*}

\end{center}

Suppose $x_1,x_2$ form a length $2$ cycle.\\

Consider the subgraph formed by vertices $x_3, x_4, x_5$. This subgraph may also contain a length $2$ cycle (formed, say, by vertices $4$ and $5$). Then by Lemma 1 (\cite{Liendo}, Lemma 1.3) the remaining vertex $3$ should be connected either to the cube $0$ or to one of the two length $2$ cycles (say, to $4$).\\

In the former case, there are neither singular pairs nor singular triples. In the latter case we may get
\begin{itemize}
\item one singular pair $(x_3,x_5)$ and
\item two singular triples $(x_1, x_3,x_5)$, $(x_2, x_3,x_5)$.
\end{itemize}

We resolve them as in the earlier examples.\\

Suppose that the subgraph formed by vertices $3$, $4$, $5$ has no cycles. Then the length of its longest path may be either $3$ or $2$ or $1$ (in which case the subgraph is totally disconnected).\\

If $3$ is connected to $4$ and $4$ is connected to $5$, then by Lemma 1 (\cite{Liendo}, Lemma 1.3) the $5$-th vertex should be connected either to the cube $0$ or to the cycle (say, to the $1$-st vertex).\\

Suppose that it is connected to the cube. Then there are neither singular pairs nor singular triples.\\

This gives us the pictures shown above.\\

If $5$ is connected to $1$, then there may be
\begin{itemize}
\item one singular pair $(x_2,x_5)$ and
\item one singular triple $(x_2, x_3,x_5)$.
\end{itemize}

We resolve them as usual and obtain the following pictures of $\AA$:

\begin{center}

\end{center}

Suppose that the length of a longest path in the subgraph formed by vertices $3$, $4$, $5$ is $2$.\\  

We may assume that $3$ is connected to $4$. By Lemma 1 (\cite{Liendo}, Lemma 1.3) $4$ should be connected either to the cube $0$ or to the cycle formed by vertices $1$, $2$. It can be connected neither to the $5$-th vertex (since there are no length $3$ paths by our assumption) nor to the $3$-rd vertex (since $3$, $4$, $5$ do not form cycles). The remaining vertex $5$ should be also connected either to the cube $0$ or to the cycle formed by $1$, $2$ or to the length $2$ path formed by $3$, $4$.\\

Let one of the vertices $4$, $5$ be connected to the cube and the other one be connected to the cycle (say, to $1$).\\

In these cases $\AA$ will be one of the following:

\begin{center}


\begin{flalign*}
\mbox{where}\; & A=(x_3^2x_2,\; x_3^2x_0,\; x_2x_3x_0,\; x_0x_3x_5);\\
& B=(x_2x_4x_5,\; x_2^2x_4,\; x_5^2x_4,\; x_4^2x_2,\; x_2x_3x_4,\; x_3x_4x_5).\\
\end{flalign*}
\end{center}

Indeed, if $5$ is connected to the cube and $4$ is connected to $1$, then we may get one singular pair $(x_2,x_4)$ and one singular triple $(x_2,x_4,x_5)$.\\

If $4$ is connected to the cube and $5$ is connected to $1$, then we may get
\begin{itemize}
\item one singular pair $(x_2,x_5)$ and 
\item two singular triples $(x_2,x_3,x_5)$, $(x_2,x_4,x_5)$.
\end{itemize}

We resolve these singular pairs and triples as usual and obtain the pictures shown above.\\

Suppose that both $4$ and $5$ are connected to the cycle.\\

Let them be connected to different vertices of the cycle (say, $4$ is connected to $1$ and $5$ is connected to $2$). In this case there may be
\begin{itemize}
\item two singular pairs $(x_1,x_5)$, $(x_2,x_4)$ and 
\item one singular triple $(x_1,x_3,x_5)$.
\end{itemize}

We resolve them as usual and obtain the following pictures:

\begin{center}


\begin{flalign*}
\mbox{where}\; & A=(x_3^2x_5,\; x_3^2x_1,\; x_3^2x_0,\; x_0x_1x_3,\; x_0x_3x_5);\;\;\;\;\; A_1=(x_3^2x_5,\; x_3^2x_1,\; x_0x_1x_3,\; x_0x_3x_5);\\
& B=(x_4^2x_2,\; x_2x_3x_4,\; x_2x_4x_5,\; x_2^2x_0,\; x_0x_2x_4).\\
\end{flalign*}
\end{center}

If $4$ and $5$ are both connected to the cube $0$, then there may be
\begin{itemize}
\item one singular pair $(x_4,x_5)$ and 
\item two singular triples $(x_1,x_4,x_5)$, $(x_2,x_4,x_5)$.
\end{itemize}

We resolve them as usual and obtain the following pictures:

\begin{center}
%
\end{center}

Now suppose that vertices $4$ and $5$ are either both connected to the same vertex of the cycle formed by vertices $1$ and $2$ (say, to $1$) or $5$ is connected to the length $2$ path (i.e. to $4$), while $4$ is connected to the cube $0$.\\

In this case $\AA$ may be one of the following:

\begin{center}
%

\begin{flalign*}
\mbox{where}\; & A=(x_2x_3x_5,\; x_3^2x_2,\; x_3^2x_5,\; x_3^2x_0,\; x_2x_3x_0,\; x_2x_5x_0,\; x_0x_3x_5);\\
& B=(x_0x_2x_4,\; x_0x_2x_5,\; x_0x_4x_5);\\
& C=(x_2x_3x_4,\; x_0x_2x_4);\;\;\;\;\; D=(x_2x_3x_5,\; x_0x_2x_5).\\
\end{flalign*}
\end{center}

If $4$ is connected to $0$ and $5$ is connected to $4$, we may get 
\begin{itemize}
\item one singular pair $(x_3,x_5)$ and 
\item two singular triples $(x_1,x_3,x_5)$, $(x_2,x_3,x_5)$.
\end{itemize}

Note that the configuration is symmetric with respect to vertices $3$ and $5$ as well as with respect to vertices $1$ and $2$.\\

If $4$ and $5$ are both connected to $1$, we may get 
\begin{itemize}
\item three singular pairs $(x_4,x_5)$, $(x_2,x_4)$, $(x_2,x_5)$ and 
\item two singular triples $(x_2,x_4,x_5)$, $(x_2,x_3,x_5)$.
\end{itemize}

After resolving all these singular pairs and triples we obtain the pictures shown above.\\

It may also happen that $5$ is connected to the length $2$ path (i.e. to $4$), while $4$ is connected to the cycle (say, to $1$).\\

In this case $\AA$ will be one of the following:

\begin{center}


\begin{flalign*}
\mbox{where}\; & A=(x_4^2x_2,\; x_2^2x_0,\; x_0x_2x_4,\; x_2x_3x_4);\\
& B=(x_1x_3x_5,\; x_0x_1x_3,\; x_0x_1x_5,\; x_1^2x_0);\\
& C=(x_2^2x_0,\; x_2x_3x_0,\; x_2x_5x_0).\\
\end{flalign*}

\end{center}

We may get 
\begin{itemize}
\item two singular pairs $(x_3,x_5)$, $(x_2,x_4)$ and 
\item two singular triples $(x_1,x_3,x_5)$, $(x_2,x_3,x_5)$.
\end{itemize}

After resolving them we obtain the pictures shown above. Note that this configuration is symmetric with respect to vertices $3$ and $5$.\\

This concludes the analysis of the cases when the length of a longest path in the subgraph formed by vertices $3$, $4$, $5$ is $2$.\\

Suppose now that this subgraph is totally disconnected. By Lemma 1 (\cite{Liendo}, Lemma 1.3) each of its vertices $3$, $4$ and $5$ should be connected to either the cycle formed by $1$ and $2$ or to the cube $0$.\\

Suppose that $4$ is connected to $1$, $3$ is connected to $2$ and $5$ is connected to $0$.\\ 

In this case $\AA$ will be one of the following:

\begin{center}


\begin{flalign*}
\mbox{where}\; & A=(x_5^2x_1,\; x_1x_3x_4,\; x_1x_4x_5,\; x_3x_4x_5);\;\;\;\;\; B=(x_5^2x_2,\; x_2x_3x_4,\; x_2x_3x_5,\; x_3x_4x_5).
\end{flalign*}
\end{center}

Indeed, we may get
\begin{itemize}
\item two singular pairs $(x_1,x_3)$, $(x_2,x_4)$ and 
\item two singular triples $(x_1,x_3,x_5)$, $(x_2,x_4,x_5)$.
\end{itemize}

After resolving these singular pairs and triples as usual we obtain the pictures shown above.\\

Suppose that $5$ is still connected to the cube $0$, but $3$ and $4$ are connected to the same vertex of the cycle (say, to $1$).\\

In this case we may get
\begin{itemize}
\item three singular pairs $(x_2,x_3)$, $(x_2,x_4)$, $(x_3,x_4)$ and 
\item three singular triples $(x_2,x_3, x_5)$, $(x_2,x_4, x_5)$, $(x_3,x_4, x_5)$.
\end{itemize}

We resolve these singular pairs and triples as usual and obtain the following possibilities for $\AA$:

\begin{center}


\begin{flalign*}
\mbox{where}\; & A=(x_2^2x_0,\; x_3^2x_0,\; x_0x_2x_3,\; x_2^2x_5,\; x_2x_3x_5);\\
& B=(x_2^2x_0,\; x_4^2x_0,\; x_0x_2x_4,\; x_2^2x_5,\; x_2x_4x_5);\;\;\;\;\; C=(x_2^2x_0,\; x_3^2x_0,\; x_0x_2x_3).\\
\end{flalign*}
\end{center}

Suppose that $3$ and $4$ are connected to the same vertex of the cycle (say, to $1$) and $5$ is connected to the other vertex (i.e. to $2$).\\

In this case we may get
\begin{itemize}
\item four singular pairs $(x_2,x_3)$, $(x_2,x_4)$, $(x_3,x_4)$, $(x_1,x_5)$ and 
\item two singular triples $(x_2,x_3, x_4)$, $(x_3,x_4, x_5)$.
\end{itemize}

We resolve these singular pairs and triples as usual and obtain the following possibilities for $\AA$:

\begin{center}


\begin{flalign*}
\mbox{where}\; & A=(x_5^2x_1,\; x_1^2x_0,\; x_0x_1x_5,\; x_1x_5x_3,\; x_1x_5x_4);\\
& B=(x_2^2x_0,\; x_0x_2x_3,\; x_2x_3x_5);\;\;\;\;\; C=(x_2^2x_0,\; x_0x_2x_4,\; x_2x_4x_5).\\
\end{flalign*}
\end{center}

Now suppose that two of the vertices $3$, $4$, $5$ (say, $4$ and $5$) are connected to the cube $0$ and the third of them is connected to the cycle (say, $3$ is connected to $1$).\\

In this case we may get
\begin{itemize}
\item two singular pairs $(x_2,x_3)$, $(x_4,x_5)$ and 
\item five singular triples $(x_2,x_3, x_4)$, $(x_2,x_3, x_5)$, $(x_1,x_4, x_5)$, $(x_2,x_4, x_5)$, $(x_3,x_4, x_5)$.
\end{itemize}

We resolve these singular pairs and triples as usual and obtain the following possibilities for $\AA$:

\begin{center}


\begin{flalign*}
\mbox{where}\; & A=(x_2^2x_0,\; x_2^2x_4,\; x_3^2x_0,\; x_0x_2x_3,\; x_2x_3x_4).\\
\end{flalign*}
\end{center}

If all three vertices $3$, $4$, $5$ are connected to the same vertex of the cycle (say, to $1$), then we may get
\begin{itemize}
\item six singular pairs $(x_2,x_3)$, $(x_2,x_4)$, $(x_2,x_5)$, $(x_3,x_4)$, $(x_3,x_5)$, $(x_4,x_5)$ and 
\item four singular triples $(x_3, x_4,x_5)$, $(x_2, x_3,x_4)$, $(x_2, x_3,x_5)$, $(x_2, x_4,x_5)$.
\end{itemize}

We resolve then as in the earlier examples.\\

The resulting possibilities for $\AA$ are as follows:

\begin{center}

\end{center}

Finally, it may happen that all three vertices $3$, $4$, $5$ are connected to the cube $0$.\\

In this case we may get
\begin{itemize}
\item three singular pairs $(x_3,x_4)$, $(x_3,x_5)$, $(x_4,x_5)$ and 
\item seven singular triples $(x_1, x_3,x_4)$, $(x_1, x_3,x_5)$, $(x_1, x_4,x_5)$, $(x_2, x_3,x_4)$, $(x_2, x_3,x_5)$,\\
$(x_2, x_4,x_5)$, $(x_3, x_4,x_5)$.
\end{itemize}

We resolve then as in the earlier examples.\\

The resulting possibilities for $\AA$ are as follows:

\begin{center}
%
\end{center}

If $\AA =\{ x_1^2x_2, x_2^2x_1, x_3^2x_4, x_4^2x_3, x_5^2x_0, x_0^3 \}$. Then 
$$
c_1=c_2,\; c_3=c_4,\; 3c_1=3c_3=2c_5=0.
$$
This means that $\GG\cong \ZZ{2}\oplus (\ZZ{3})^{\oplus 2}$ with generators 
$$
f_1=(1,{\om},{\om},1,1,1), \; f_2=(1,1,1,{\om},{\om},1), \; f_3=(1,1,1,1,1,\eta ), \; \om=\sqr{3}, \eta=\sqr{2}.
$$
Then $i_1+i_2 \equiv i_3+i_4\equiv 0 \; mod \; 3$ and $i_5 \equiv 0 \; mod \; 2$. This means that $\AAbar=\AA\cup \{ f_3(x_1,x_2), f_3(x_3,x_4) \}$.\\

If $\AA =\{ x_1^2x_2, x_2^2x_1, x_3^2x_4, x_4^2x_5, x_5^2x_0, x_0^3 \}$. Then 
$$
c_1=c_2,\; c_4=-2c_3,\; c_5=4c_3,\; 3c_1=8c_3=0.
$$
This means that $\GG\cong \ZZ{8}\oplus \ZZ{3}$ with generators 
$$
f_1=(1,{\om},{\om},1,1,1), \; f_2=(1,1,1,{\eta},{\eta}^{-2},{\eta}^{4}), \; \om=\sqr{3}, \eta=\sqr{8}.
$$
Then $i_1+i_2 \equiv 0 \; mod \; 3$ and $i_3 \equiv 2i_4+4i_5 \; mod \; 8$. This means that $\AAbar=\AA\cup \{ f_3(x_1,x_2) \}$.\\

If $\AA =\{ x_1^2x_2, x_2^2x_1, x_3^2x_4, x_4^2x_5, x_5^2x_4, x_3^2x_5, x_0^3 \}$. Then 
$$
c_1=c_2,\; c_4=c_5=-2c_3,\; 3c_1=6c_3=0.
$$
This means that $\GG\cong \ZZ{6}\oplus \ZZ{3}$ with generators 
$$
f_1=(1,{\om},{\om},1,1,1), \; f_2=(1,1,1,{\eta},{\eta}^{-2},{\eta}^{-2}), \; \om=\sqr{3}, \eta=\sqr{6}.
$$
Then $i_1+i_2 \equiv 0 \; mod \; 3$ and $i_3 \equiv 2(i_4+i_5) \; mod \; 6$. This means that $\AAbar=\AA\cup \{ f_3(x_1,x_2)$, $f_3(x_4,x_5) \}$.\\

Let $\AAN{0} =\{ x_1^2x_2, x_2^2x_1, x_3^2x_1, x_4^2x_5, x_5^2x_4, x_3^2x_4, x_0^3 \}$, $\AAN{0} =\{ x_1^2x_2,$ $x_2^2x_1,$ $x_3^2x_4,$ $x_4^2x_5,$ $x_5^2x_4,$ $x_5^2x_2,$ $x_0^3 \}$ or $\AAN{0} =\{ x_1^2x_2, x_2^2x_1, x_3^2x_4, x_4^2x_5, x_5^2x_1, x_2x_4x_5, x_0^3 \}$. Then 
$$
c_1=c_2=c_4=c_5=-2c_3,\; 6c_3=0.
$$
All such cases were already considered.\\

If $\AA =\{ x_1^2x_2, x_2^2x_1, x_3^2x_4, x_4^2x_5, x_5^2x_4, x_5^2x_0, x_0^3 \}$. Then 
$$
c_1=c_2,\; c_4=c_5=0,\; 3c_1=2c_3=0.
$$
This means that $\GG\cong \ZZ{2}\oplus \ZZ{3}$ with generators 
$$
f_1=(1,{\om},{\om},1,1,1), \; f_2=(1,1,1,{\eta},1,1), \; \om=\sqr{3}, \eta=\sqr{2}.
$$
Then $i_1+i_2 \equiv 0 \; mod \; 3$ and $i_3 \equiv 0 \; mod \; 2$. This means that $\AAbar=\AA\cup \{ f_3(x_0,x_4,x_5), x_3^2\cdot f_1(x_0,x_4,x_5), f_3(x_1,x_2) \}$.\\

If $\AA =\{ x_1^2x_2, x_2^2x_1, x_3^2x_4, x_4^2x_5, x_5^2x_4, x_0x_3x_5, x_0^3 \}$. Then 
$$
c_1=c_2,\; c_3=c_4=c_5=0,\; 3c_1=0.
$$
This means that $\GG\cong \ZZ{3}$ with a generator 
$$
f=(1,{\om},{\om},1,1,1), \; \om=\sqr{3}.
$$
Then $i_1+i_2 \equiv 0 \; mod \; 3$. This means that $\AAbar=\AA\cup \{ f_3(x_0,x_3,x_4,x_5), f_3(x_1,x_2) \}$.\\

Let $\AAN{0} =\{ x_1^2x_2, x_2^2x_1, x_3^2x_4, x_4^2x_5, x_5^2x_4, x_1x_3x_5, x_0^3 \}$. Then 
$$
c_1=c_2=c_3=c_4=c_5,\; 3c_3=0.
$$
All such cases were already considered.\\

If $\AA =\{ x_1^2x_2, x_2^2x_1, x_3^2x_4, x_4^2x_5, x_5^2x_1, x_5^2x_2, x_0^3 \}$. Then 
$$
c_1=c_2=-8c_3,\; c_4=-2c_3,\; c_5=4c_3,\; 24c_3=0.
$$
This means that $\GG\cong \ZZ{24}$ with a generator 
$$
f=(1,{\om}^{-8},{\om}^{-8},\om,{\om}^{-2},{\om}^{4}), \; \om=\sqr{24}.
$$
Then $i_3+4i_5 \equiv 2i_4+8(i_1+i_2) \; mod \; 24$. This means that $\AAbar=\AA\cup \{ f_3(x_1,x_2) \}$.\\

If $\AA =\{ x_1^2x_2, x_2^2x_1, x_3^2x_4, x_4^2x_5, x_5^2x_1, x_2^2x_0, x_0^3 \}$. Then 
$$
c_1=c_2=0,\; c_4=-2c_3,\; c_5=4c_3,\; 8c_3=0.
$$
This means that $\GG\cong \ZZ{8}$ with a generator 
$$
f=(1,1,1,\om, {\om}^{-2},{\om}^{4}), \; \om=\sqr{8}.
$$
Then $i_3 \equiv 2i_4+4i_5 \; mod \; 8$. This means that $\AAbar=\AA\cup \{ f_3(x_0,x_1,x_2), x_5^2\cdot f_1(x_0,x_1,x_2) \}$.\\

If $\AA =\{ x_1^2x_2, x_2^2x_1, x_3^2x_4, x_4^2x_5, x_5^2x_1, x_2x_3x_5, x_0^3 \}$. Then 
$$
c_1=c_2=c_3=c_4=c_5,\; 3c_3=0.
$$
This case was already considered.\\

If $\AA =\{ x_1^2x_2, x_2^2x_1, x_3^2x_4, x_4^2x_5, x_5^2x_1, x_0x_2x_5, x_0^3 \}$. Then 
$$
c_1=c_2=c_5=0,\; c_4=2c_3,\; 4c_3=0.
$$
This case was already considered.\\

If $\AA =\{ x_1^2x_2, x_2^2x_1, x_3^2x_4, x_4^2x_1, x_5^2x_0, x_2x_4x_5, x_0^3 \}$. Then 
$$
c_1=c_2=0,\; c_4=c_5=2c_3,\; 4c_3=0.
$$
This means that $\GG\cong \ZZ{4}$ with a generator 
$$
f=(1,1,1,\om, {\om}^{2},{\om}^{2}), \; \om=\sqr{4}.
$$
Then $i_3 \equiv 2(i_4+i_5) \; mod \; 4$. This means that $\AAbar=\AA\cup \{ f_3(x_0,x_1,x_2), x_3^2\cdot f_1(x_4,x_5), f_2(x_4,x_5)\cdot f_1(x_0,x_1,x_2) \}$.\\

If $\AA =\{ x_1^2x_2, x_2^2x_1, x_3^2x_4, x_4^2x_1, x_5^2x_0, x_2x_3x_4, x_0^3 \}$. Then 
$$
c_1=c_2=c_3=c_4,\; 3c_1=2c_5=0.
$$
This case was already considered.\\

If $\AA =\{ x_1^2x_2, x_2^2x_1, x_3^2x_4, x_4^2x_1, x_5^2x_0, x_2^2x_5, x_0^3 \}$. Then 
$$
c_1=c_2=c_5=0,\; c_4=2c_3,\; 4c_3=0.
$$
This case was already considered.\\

If $\AA =\{ x_1^2x_2, x_2^2x_1, x_3^2x_4, x_4^2x_1, x_4^2x_2, x_5^2x_0, x_0^3 \}$. Then 
$$
c_1=c_2=4c_3,\; c_4=-2c_3,\; 12c_3=2c_5=0.
$$
This means that $\GG\cong \ZZ{12} \oplus \ZZ{2}$ with generators 
$$
f_1=(1,{\om}^{4},{\om}^{4},\om, {\om}^{-2},1), \; f_2=(1,1,1,1,1,\eta), \; \om=\sqr{12}, \eta=\sqr{2}.
$$
Then $i_3+4(i_1+i_2) \equiv 2i_4 \; mod \; 12$ and $i_5 \equiv 0 \; mod \; 2$. This means that $\AAbar=\AA\cup \{ f_3(x_1,x_2) \}$.\\

Let $\AAN{0} =\{ x_1^2x_2, x_2^2x_1, x_3^2x_4, x_4^2x_1, x_5^2x_0, x_2^2x_0, x_0^3 \}$ or $\AAN{0} =\{ x_1^2x_2,$ $x_2^2x_1,$ $x_3^2x_4,$ $x_4^2x_0,$ $x_5^2x_0,$ $x_4^2x_1,$ $x_0^3 \}$. Then 
$$
c_1=c_2=0,\; c_4=2c_3,\; 4c_3=2c_5=0.
$$

For $\AA=\AAN{0}\cup \{ x_5^2x_2 \}$ the extra condition is automatically satisfied.\\

This case will appear and will be analyzed below.\\

For $\AA=\AAN{0}\cup \{ x_5^2x_4 \}$ we get an extra condition $2c_3=0$, i.e.
$$
c_1=c_2=c_4=0,\; 2c_3=2c_5=0.
$$
This case was already considered.\\

For $\AA=\AAN{0}\cup \{ x_2x_3x_5 \}$ and $\AA=\AAN{0}\cup \{ x_3x_4x_5 \}$ we get the same extra condition $c_3=c_5$, i.e.
$$
c_1=c_2=c_4=0,\; c_3=c_5,\; 2c_5=0.
$$
This case was already considered.\\

Let $\AAN{0} =\{ x_1^2x_2, x_2^2x_1, x_3^2x_4, x_4^2x_1, x_5^2x_0, x_0x_2x_4, x_0^3 \}$. Then 
$$
c_1=c_2=c_4=0,\; 2c_3=2c_5=0.
$$
All such cases were already considered.\\

If $\AA =\{ x_1^2x_2, x_2^2x_1, x_3^2x_4, x_4^2x_0, x_5^2x_1, x_5^2x_2, x_0^3 \}$. Then 
$$
c_1=c_2=-2c_5,\; c_4=2c_3,\; 4c_3=6c_5=0.
$$
This means that $\GG\cong \ZZ{6} \oplus \ZZ{4}$ with generators 
$$
f_1=(1,{\om}^{-2},{\om}^{-2},1,1,\om), \; f_2=(1,1,1,\eta, {\eta}^2,1), \; \om=\sqr{6}, \eta=\sqr{4}.
$$
Then $i_5 \equiv 2(i_1+i_2) \; mod \; 6$ and $i_3 \equiv 2i_4 \; mod \; 4$. This means that $\AAbar=\AA\cup \{ f_3(x_1,x_2) \}$.\\

If $\AA =\{ x_1^2x_2, x_2^2x_1, x_3^2x_4, x_4^2x_0, x_5^2x_1, x_2x_3x_5, x_0^3 \}$. Then 
$$
c_1=c_2=c_4=0,\; c_5=c_3,\; 2c_3=0.
$$
This case was already considered.\\

If $\AA =\{ x_1^2x_2, x_2^2x_1, x_3^2x_4, x_4^2x_0, x_5^2x_1, x_2^2x_3, x_0^3 \}$. Then 
$$
c_1=c_2=c_3=c_4=0,\; 2c_5=0.
$$
This case was already considered.\\

Let $\AAN{0} =\{ x_1^2x_2, x_2^2x_1, x_3^2x_4, x_4^2x_0, x_5^2x_1, x_0x_2x_5, x_0^3 \}$. Then 
$$
c_1=c_2=c_5=0,\; c_4=2c_3,\; 4c_3=0.
$$
All such cases were already considered.\\

Let $\AAN{0} =\{ x_1^2x_2, x_2^2x_1, x_3^2x_4, x_4^2x_0, x_5^2x_1, x_2x_4x_5, x_0^3 \}$. Then 
$$
c_1=c_2=0,\; c_4=c_5=2c_3,\; 4c_3=0.
$$
All such cases were already considered.\\

Let $\AAN{0} =\{ x_1^2x_2, x_2^2x_1, x_3^2x_4, x_4^2x_0, x_5^2x_1, x_2^2x_4, x_0^3 \}$. Then 
$$
c_1=c_2=c_4=0,\; 2c_3=2c_5=0.
$$
All such cases were already considered.\\

Let $\AAN{0} =\{ x_1^2x_2, x_2^2x_1, x_3^2x_4, x_4^2x_0, x_5^2x_1, x_2^2x_0, x_0^3 \}$ or $\AAN{0} =\{ x_1^2x_2,$ $x_2^2x_1,$ $x_3^2x_4,$ $x_4^2x_0,$ $x_5^2x_0,$ $x_5^2x_1,$ $x_0^3 \}$. Then 
$$
c_1=c_2=0,\; c_4=2c_3,\; 4c_3=2c_5=0. 
$$

For $\AA=\AAN{0}\cup \{ x_4^2x_2 \}$ the extra condition $c_2+2c_4=0$ is automatically satisfied.\\

This means that $\GG\cong \ZZ{2} \oplus \ZZ{4}$ with generators 
$$
f_1=(1,1,1,\om,{\om}^{2},1), \; f_2=(1,1,1,1,1,\eta), \; \om=\sqr{4}, \eta=\sqr{2}.
$$
Then $i_3 \equiv 2i_4 \; mod \; 4$ and $i_5 \equiv 0 \; mod \; 2$. This means that $\AAbar=\AA\cup \{ f_3(x_0,x_1,x_2), x_4^2\cdot f_1(x_0,x_1,x_2), x_5^2\cdot f_1(x_0,x_1,x_2) \}$.\\

Any further extra condition (after adding to $\AAN{0}$ a monomial $x_0^{i_0}x_1^{i_1}x_2^{i_2}x_3^{i_3}x_4^{i_4}x_5^{i_5}$) has the form $ac_3+bc_5=0$, $a,b\in \mathbb Z$, i.e. either $a\cdot c_3=0$ or $c_5=a\cdot c_3$ for some $a\in \{ 0,1,2,3 \}$. All such cases have already appeared.\\ 

Let $\AAN{0} =\{ x_1^2x_2, x_2^2x_1, x_3^2x_4, x_4^2x_0, x_5^2x_1, x_5^2x_4, x_0^3 \}$. Then 
$$
c_1=c_2=c_4=0,\; 2c_3=2c_5=0. 
$$
All such cases have already appeared.\\

Let $\AAN{0} =\{ x_1^2x_2, x_2^2x_1, x_3^2x_4, x_4^2x_1, x_5^2x_1, x_2x_4x_5, x_0^3 \}$. Then 
$$
c_1=c_2=4c_3,\; c_4=c_5=-2c_3,\; 12c_3=0. 
$$

For $\AA=\AAN{0}\cup \{ x_3^2x_5 \}$ the extra condition $2c_3+c_5=0$ is automatically satisfied.\\

This means that $\GG\cong \ZZ{12}$ with a generator
$$
f=(1,{\om}^{4},{\om}^{4},\om,{\om}^{-2},{\om}^{-2}), \; \om=\sqr{12}.
$$
Then $i_3+4(i_1+i_2) \equiv 2(i_4+i_5) \; mod \; 12$. This means that $\AAbar=\AA\cup \{ f_3(x_1,x_2), f_2(x_4,x_5)\cdot f_1(x_1,x_2) \}$.\\

Any further extra condition (after adding to $\AAN{0}$ a monomial $x_0^{i_0}x_1^{i_1}x_2^{i_2}x_3^{i_3}x_4^{i_4}x_5^{i_5}$) has the form $a\cdot c_3=0$ for some $a\in  \{ 0,1,2,3,4,6 \}$. All such cases have already appeared.\\

If $\AA =\{ x_1^2x_2, x_2^2x_1, x_3^2x_4, x_4^2x_0, x_5^2x_4, x_5^2x_0, x_0^3 \}$. Then 
$$
c_1=c_2,\; c_4=0,\; 3c_1=2c_3=2c_5=0.
$$
This means that $\GG\cong (\ZZ{2})^{\oplus 2} \oplus \ZZ{3}$ with generators 
$$
f_1=(1,1,1,1,1,\om), \; f_2=(1,1,1,\om ,1,1), \; f_3=(1,\eta,\eta, 1,1,1), \; \om=\sqr{2}, \eta=\sqr{3}.
$$
Then $i_1+i_2 \equiv 0 \; mod \; 3$ and $i_3 \equiv i_5 \equiv 0 \; mod \; 2$. This means that $\AAbar=\AA\cup \{ f_3(x_0,x_4)$, $f_3(x_1,x_2)$, $x_3^2x_0 \}$.\\

If $\AA =\{ x_1^2x_2, x_2^2x_1, x_3^2x_4, x_4^2x_0, x_5^2x_0, x_3x_4x_5, x_0^3 \}$. Then 
$$
c_1=c_2,\; c_3=c_5,\; c_4=0,\; 3c_1=2c_5=0.
$$
This case was already considered.\\

If $\AA =\{ x_1^2x_2, x_2^2x_1, x_3^2x_4, x_4^2x_1, x_5^2x_4, x_5^2x_1, x_2^2x_0, x_0^3 \}$. Then 
$$
c_1=c_2=c_4=0,\; 2c_3=2c_5=0.
$$
This case was already considered.\\

Let $\AAN{0} =\{ x_1^2x_2, x_2^2x_1, x_3^2x_4, x_4^2x_0, x_5^2x_0, x_1x_4x_5, x_0^3 \}$. Then 
$$
c_1=c_2=0,\; c_4=c_5=2c_3,\; 4c_3=0. 
$$
All such cases have already appeared.\\ 

If $\AA =\{ x_1^2x_2, x_2^2x_1, x_3^2x_4, x_4^2x_0, x_5^2x_4, x_0x_3x_5, x_0^3 \}$. Then 
$$
c_1=c_2,\; c_4=2c_3,\; c_5=-c_3,\; 3c_1=4c_3=0.
$$
This means that $\GG\cong \ZZ{4} \oplus \ZZ{3}$ with generators 
$$
f_1=(1,\om,\om,1,1,1), \; f_2=(1,1,1,\eta ,{\eta}^2,{\eta}^3), \; \om=\sqr{3}, \eta=\sqr{4}.
$$
Then $i_1+i_2 \equiv 0 \; mod \; 3$ and $i_3 \equiv i_5+2i_4 \; mod \; 4$. This means that $\AAbar=\AA\cup \{ f_3(x_1,x_2) \}$.\\

Let $\AAN{0} =\{ x_1^2x_2, x_2^2x_1, x_3^2x_4, x_4^2x_0, x_5^2x_4, x_1x_3x_5, x_0^3 \}$ or $\AAN{0} =\{ x_1^2x_2,$ $x_2^2x_1,$ $x_3^2x_4,$ $x_4^2x_1,$ $x_5^2x_4,$ $x_0x_3x_5,$ $x_0^3 \}$. Then 
$$
c_1=c_2=0,\; c_4=2c_3,\; c_5=-c_3,\; 4c_3=0. 
$$
All such cases have already appeared.\\

Let $\AAN{0} =\{ x_1^2x_2, x_2^2x_1, x_3^2x_4, x_4^2x_1, x_5^2x_4, x_2x_3x_5, x_0^3 \}$ or $\AAN{0} =\{ x_1^2x_2,$ $x_2^2x_1,$ $x_3^2x_4,$ $x_4^2x_1,$ $x_5^2x_4,$ $x_1x_3x_5,$ $x_0^3 \}$. Then 
$$
c_1=c_2=4c_3,\; c_4=-2c_3,\; c_5=-5c_3,\; 12c_3=0. 
$$

For $\AA=\AAN{0}\cup \{ x_1x_3x_5, x_4^2x_2 \}$ extra conditions $c_1+c_3+c_5=c_2+2c_4=0$ are automatically satisfied.\\

This means that $\GG\cong \ZZ{12}$ with a generator
$$
f=(1,{\om}^{4},{\om}^{4},\om,{\om}^{-2},{\om}^{-5}), \; \om=\sqr{12}.
$$
Then $i_3+4(i_1+i_2) \equiv 2i_4+5i_5 \; mod \; 12$. This means that $\AAbar=\AA\cup \{ f_3(x_1,x_2) \}$.\\

Any further extra condition (after adding to $\AAN{0}$ a monomial $x_0^{i_0}x_1^{i_1}x_2^{i_2}x_3^{i_3}x_4^{i_4}x_5^{i_5}$) has the form $a\cdot c_3=0$ for some $a\in  \{ 0,1,2,3,4,6 \}$. All such cases have already appeared.\\

Let $\AAN{0} =\{ x_1^2x_2, x_2^2x_1, x_3^2x_4, x_4^2x_1, x_5^2x_1, x_5^2x_4, x_0^3 \}$, $\AAN{0} =\{ x_1^2x_2,$ $x_2^2x_1,$ $x_3^2x_4,$ $x_4^2x_1,$ $x_5^2x_2,$ $x_5^2x_4,$ $x_0^3 \}$ or $\AAN{0} =\{ x_1^2x_2, x_2^2x_1, x_3^2x_4, x_4^2x_1, x_5^2x_2, x_1^2x_4, x_0^3 \}$. Then 
$$
c_1=c_2=c_4=-2c_3,\; 6c_3=2(c_5-c_3)=0.
$$

For $\AAN{1}=\AAN{0}\cup \{ x_2x_3x_4 \}$ we get an extra condition $3c_3=0$, i.e.
$$
c_1=c_2=c_3=c_4=-2c_5,\; 6c_5=0. 
$$
All such cases have already appeared.\\

For $\AAN{1}=\AAN{0}\cup \{ x_0x_1x_3 \}$, $\AAN{1}=\AAN{0}\cup \{ x_0x_3x_5 \}$, $\AAN{1}=\AAN{0}\cup \{ x_0x_2x_4 \}$, $\AAN{1}=\AAN{0}\cup \{ x_2^2x_0 \}$ we get extra conditions $c_3=0$, $c_5=-c_3$, $2c_3=0$, $2c_3=0$ respectively. They all imply that $2c_3=0$, i.e.
$$
c_1=c_2=c_4=0,\; 2c_3=2c_5=0.
$$
All such cases have already appeared.\\

For $\AAN{1}=\AAN{0}\cup \{ x_3^2x_5 \}$ and $\AAN{1}=\AAN{0}\cup \{ x_2x_4x_5 \}$ we get the same extra condition $c_5=-2c_3$, i.e.
$$
c_1=c_2=c_4=c_5=-2c_3,\; 6c_3=0. 
$$
All such cases have already appeared.\\

For $\AAN{1}=\AAN{0}\cup \{ x_3^2x_1 \}$ the extra condition $c_1+2c_3=0$ is automatically satisfied.\\

For $\AA=\AAN{1}\cup \{ x_4^2x_2 \}$ the extra condition $c_2+2c_4=0$ is automatically satisfied.\\ 

This case was already considered.\\

Let $\AAN{0} =\{ x_1^2x_2, x_2^2x_1, x_3^2x_4, x_4^2x_1, x_5^2x_1, x_5^2x_2, x_0^3 \}$. Then 
$$
c_1=c_2=4c_3,\; c_4=-2c_3,\; 12c_3=2(c_5+2c_3)=0. 
$$

For $\AA=\AAN{0}\cup \{ x_4^2x_2 \}$ the extra condition $c_2+2c_4=0$ is automatically satisfied.\\

This means that $\GG\cong \ZZ{12}\oplus \ZZ{2}$ with generators
$$
f_1=(1,{\om}^{4},{\om}^{4},\om,{\om}^{-2},{\om}^{-2}), \; f_2=(1,1,1,1,1,\eta),\; \om=\sqr{12}, \eta=\sqr{2}.
$$
Then $i_3+4(i_1+i_2) \equiv 2(i_4+i_5) \; mod \; 12$, $i_5 \equiv 0 \; mod \; 2$. This means that $\AAbar=\AA\cup \{ f_3(x_1,x_2), x_4^2\cdot f_1(x_1,x_2), x_3^2x_5 \}$.\\

For $\AA=\AAN{0}\cup \{ x_2x_4x_5 \}$ we get an extra condition $c_5=-2c_3$, i.e.
$$
c_1=c_2=4c_3,\; c_4=c_5=-2c_3,\; 12c_3=0. 
$$
This case has already appeared.\\

For $\AA=\AAN{0}\cup \{ x_2^2x_0 \}$ we get an extra condition $4c_3=0$, i.e.
$$
c_1=c_2=0,\; c_4=2c_3,\; 4c_3=2c_5=0. 
$$
This case has already appeared.\\ 

For $\AA=\AAN{0}\cup \{ x_0x_2x_4 \}$ we get an extra condition $2c_3=0$, i.e.
$$
c_1=c_2=c_4=0,\; 2c_3=2c_5=0. 
$$
This case has already appeared.\\

For $\AA=\AAN{0}\cup \{ x_2x_3x_4 \}$ we get an extra condition $3c_3=0$, i.e.
$$
c_1=c_2=c_3=c_4=-2c_5,\; 6c_5=0. 
$$
This case has already appeared.\\

Let $\AAN{0} =\{ x_1^2x_2, x_2^2x_1, x_3^2x_4, x_4^2x_1, x_5^2x_2, x_1^2x_0, x_0^3 \}$. Then 
$$
c_1=c_2=0,\; c_4=2c_3,\; 4c_3=2c_5=0. 
$$
All such cases have already appeared.\\

Let $\AAN{0} =\{ x_1^2x_2, x_2^2x_1, x_3^2x_4, x_4^2x_1, x_5^2x_2, x_0x_1x_5, x_0^3 \}$. Then 
$$
c_1=c_2=c_5=0,\; c_4=2c_3,\; 4c_3=0. 
$$
All such cases have already appeared.\\

Let $\AAN{0} =\{ x_1^2x_2, x_2^2x_1, x_3^2x_4, x_4^2x_1, x_5^2x_2, x_1x_4x_5, x_0^3 \}$. Then 
$$
c_1=c_2=4c_3,\; c_4=c_5=-2c_3,\; 12c_3=0. 
$$
All such cases have already appeared.\\

Let $\AAN{0} =\{ x_1^2x_2, x_2^2x_1, x_3^2x_4, x_4^2x_1, x_5^2x_2, x_1x_3x_5, x_0^3 \}$ or $\AAN{0} =\{ x_1^2x_2,$ $x_2^2x_1,$ $x_3^2x_4,$ $x_4^2x_1,$ $x_5^2x_1,$ $x_3x_4x_5,$ $x_0^3 \}$. Then 
$$
c_1=c_2=c_4=-2c_3,\; c_5=c_3,\; 6c_3=0. 
$$
All such cases have already appeared.\\

Let $\AAN{0} =\{ x_1^2x_2, x_2^2x_1, x_3^2x_4, x_4^2x_1, x_5^2x_1, x_0x_4x_5, x_0^3 \}$. Then 
$$
c_1=c_2=0,\; c_4=c_5=2c_3,\; 4c_3=0.
$$
All such cases have already appeared.\\

If $\AA =\{ x_1^2x_2, x_2^2x_1, x_3^2x_1, x_3^2x_2, x_4^2x_1, x_4^2x_2, x_5^2x_0, x_0^3 \}$, then 
$$
c_1=c_2=-2c_3,\; 6c_3=2(c_4-c_3)=2c_5=0.
$$
This means that $\GG\cong \ZZ{6}\oplus (\ZZ{2})^{\oplus 2}$ with generators
$$
f_1=(1,{\om}^{-2},{\om}^{-2},\om,{\om},1), \; f_2=(1,1,1,1,\eta,1),\; f_3=(1,1,1,1,1,\eta),\; \om=\sqr{6}, \eta=\sqr{2}.
$$
Then $i_3+i_4 \equiv 2(i_1+i_2) \; mod \; 6$, $i_4\equiv i_5 \equiv 0 \; mod \; 2$. This means that $\AAbar=\AA\cup \{ f_3(x_1,x_2) \}$.\\

Let $\AAN{0} =\{ x_1^2x_2, x_2^2x_1, x_3^2x_1, x_3^2x_2, x_4^2x_1, x_5^2x_0, x_0^3 \}$. Then 
$$
c_1=c_2=-2c_3,\; 6c_3=2(c_4-c_3)=2c_5=0.
$$

For $\AA=\AAN{0}\cup \{ x_2^2x_5 \}$ we get an extra condition $c_5=-2c_3$, i.e.
$$
c_1=c_2=c_5=0,\; 2c_3=2c_4=0. 
$$
This case has already appeared.\\

For $\AA=\AAN{0}\cup \{ x_2x_3x_4 \}$ we get an extra condition $c_4=c_3$, i.e.
$$
c_1=c_2=-2c_3,\; c_4=c_3,\; 6c_3=2c_5=0.
$$
This means that $\GG\cong \ZZ{6}\oplus \ZZ{2}$ with generators
$$
f_1=(1,{\om}^{-2},{\om}^{-2},\om,{\om},1), \; f_2=(1,1,1,1,1,\eta),\; f_3=(1,1,1,1,1,\eta),\; \om=\sqr{6}, \eta=\sqr{2}.
$$
Then $i_3+i_4 \equiv 2(i_1+i_2) \; mod \; 6$, $i_5 \equiv 0 \; mod \; 2$. This means that $\AAbar=\AA\cup \{ f_3(x_1,x_2), f_2(x_3,x_4)\cdot f_1(x_1,x_2) \}$.\\ 

For $\AA=\AAN{0}\cup \{ x_2x_4x_5 \}$ we get an extra condition $c_5=2c_3-c_4$, i.e.
$$
c_1=c_2=0,\; c_5=c_4,\; 2c_3=2c_4=0.
$$
This case has already appeared.\\

For $\AAN{1}=\AAN{0}\cup \{ x_5^2x_2 \}$ we get an extra condition $c_2=0$, i.e.
$$
c_1=c_2=0,\; 2c_3=2c_4=2c_5=0.
$$
All such cases have already appeared.\\

For $\AAN{1}=\AAN{0}\cup \{ x_2x_3x_5 \}$ we get an extra condition $c_5=c_3$, i.e.
$$
c_1=c_2=0,\; c_5=c_3,\; 2c_3=2c_4=0.
$$
All such cases have already appeared.\\

For $\AAN{1}=\AAN{0}\cup \{ x_3x_4x_5 \}$ we get an extra condition $c_3+c_4+c_5=0$, i.e.
$$
c_1=c_2=0,\; c_5=c_3+c_4,\; 2c_3=2c_4=0.
$$
All such cases have already appeared.\\

Let $\AAN{0} =\{ x_1^2x_2, x_2^2x_1, x_3^2x_2, x_4^2x_1, x_5^2x_0, x_2^2x_0, x_0^3 \}$. Then 
$$
c_1=c_2=0,\; 2c_3=2c_4=2c_5=0.
$$
All such cases have already appeared.\\

Let $\AAN{0} =\{ x_1^2x_2, x_2^2x_1, x_3^2x_2, x_4^2x_1, x_5^2x_0, x_2^2x_5, x_0^3 \}$. Then 
$$
c_1=c_2=c_5=0,\; 2c_3=2c_4=0.
$$
All such cases have already appeared.\\

Let $\AAN{0} =\{ x_1^2x_2, x_2^2x_1, x_3^2x_2, x_4^2x_1, x_5^2x_0, x_2x_3x_4, x_0^3 \}$. Then 
$$
c_1=c_2=-2c_3,\; c_4=c_3,\; 6c_3=2c_5=0.
$$

For $\AA=\AAN{0}\cup \{ x_1x_3x_4 \}$ the extra condition $c_1+c_3+c_4=0$ is automatically satisfied.\\ 

This case has already appeared.\\

Let $\AAN{0} =\{ x_1^2x_2, x_2^2x_1, x_3^2x_2, x_4^2x_1, x_5^2x_0, x_2x_4x_5, x_0^3 \}$. Then 
$$
c_1=c_2=0,\; c_4=c_5,\; 2c_3=2c_5=0.
$$
All such cases have already appeared.\\

Let $\AAN{0} =\{ x_1^2x_2, x_2^2x_1, x_3^2x_1, x_4^2x_1, x_5^2x_0, x_2x_3x_5, x_0^3 \}$. Then 
$$
c_1=c_2=0,\; c_5=c_3,\; 2c_3=2c_4=0.
$$
All such cases have already appeared.\\

Let $\AAN{0} =\{ x_1^2x_2, x_2^2x_1, x_3^2x_1, x_4^2x_1, x_5^2x_0, x_4^2x_0, x_0^3 \}$. Then 
$$
c_1=c_2=0,\; 2c_3=2c_4=2c_5=0.
$$
All such cases have already appeared.\\

Let $\AAN{0} =\{ x_1^2x_2, x_2^2x_1, x_3^2x_1, x_4^2x_1, x_5^2x_0, x_3x_4x_5, x_0^3 \}$. Then 
$$
c_1=c_2=0,\; c_5=c_3+c_4,\; 2c_3=2c_4=0.
$$
All such cases have already appeared.\\

Let $\AAN{0} =\{ x_1^2x_2, x_2^2x_1, x_3^2x_2, x_4^2x_1, x_5^2x_0, x_0x_2x_4, x_0^3 \}$. Then 
$$
c_1=c_2=c_4=0,\; 2c_3=2c_5=0.
$$
All such cases have already appeared.\\

Let $\AAN{0} =\{ x_1^2x_2, x_2^2x_1, x_3^2x_2, x_4^2x_1, x_5^2x_0, x_4^2x_0, x_0^3 \}$. Then 
$$
c_1=c_2=0,\; 2c_3=2c_4=2c_5=0.
$$
All such cases have already appeared.\\

If $\AA =\{ x_1^2x_2, x_2^2x_1, x_3^2x_1, x_4^2x_1, x_5^2x_0, x_2x_3x_4, x_0^3 \}$, then 
$$
c_1=c_2=-2c_3,\; c_4=c_3,\; 6c_3=2c_5=0.
$$
This case has already appeared.\\

Let $\AAN{0} =\{ x_1^2x_2, x_2^2x_1, x_3^2x_1, x_4^2x_1, x_5^2x_2, x_2x_3x_4, x_0^3 \}$. Then 
$$
c_1=c_2=-2c_3,\; c_4=c_3,\; 6c_3=2(c_5-c_3)=0.
$$

For $\AAN{1}=\AAN{0}\cup \{ x_0x_3x_4 \}$, $\AAN{1}=\AAN{0}\cup \{ x_0x_3x_5 \}$, $\AAN{1}=\AAN{0}\cup \{ x_0x_4x_5 \}$ we get extra conditions $2c_3=0$, $x_5=-x_3$, $x_5=-x_3$ respectively. They all imply that $2c_3=0$, i.e.
$$
c_1=c_2=0,\; c_4=c_3,\; 2c_3=2c_5=0.
$$
All such cases have already appeared.\\

For $\AAN{1}=\AAN{0}\cup \{ x_3x_4x_5 \}$ we get an extra condition $c_5=-2c_3$, i.e.
$$
c_1=c_2=c_5=-2c_3,\; c_4=c_3,\; 6c_3=0.
$$
All such cases have already appeared.\\

Let $\AAN{0} =\{ x_1^2x_2, x_2^2x_1, x_3^2x_1, x_3^2x_2, x_4^2x_1, x_5^2x_2, x_0^3 \}$. Then 
$$
c_1=c_2=-2c_3,\; 6c_3=2(c_4-c_3)=2(c_5-c_3)=0.
$$

For $\AAN{1}=\AAN{0}\cup \{ x_3x_4x_5 \}$ we get an extra condition $c_3+c_4+c_5=0$, i.e.
$$
c_1=c_2=-2c_3,\; c_5=-c_3-c_4,\; 6c_3=2(c_4-c_3)=0.
$$

For $\AAN{2}=\AAN{1}\cup \{ x_0x_1x_5 \}$ and $\AAN{2}=\AAN{1}\cup \{ x_1^2x_0 \}$ we get extra conditions $2c_3=0$ and $c_4=3c_3$ respectively. They both imply that $2c_3=0$, i.e.
$$
c_1=c_2=0,\; c_5=c_3+c_4,\; 2c_3=2c_4=0.
$$
All such cases have already appeared.\\

For $\AAN{2}=\AAN{1}\cup \{ x_1x_3x_5 \}$ we get an extra condition $c_4=-2c_3$, i.e.
$$
c_1=c_2=c_4=-2c_3,\; c_5=c_3,\; 6c_3=0.
$$
All such cases have already appeared.\\

For $\AAN{2}=\AAN{1}\cup \{ x_1x_4x_5 \}$ we get an extra condition $3c_3=0$, i.e.
$$
c_1=c_2=c_3=-2c_4,\; c_5=c_4,\; 6c_4=0.
$$
All such cases have already appeared.\\

For $\AAN{2}=\AAN{1}\cup \{ x_0x_2x_4 \}$ and $\AAN{2}=\AAN{1}\cup \{ x_2^2x_0 \}$ we get extra conditions $c_4=2c_3$ and $2c_3=0$ respectively. They both imply that $2c_3=0$, i.e.
$$
c_1=c_2=0,\; c_5=c_3+c_4,\; 2c_3=2c_4=0.
$$
All such cases have already appeared.\\

For $\AAN{2}=\AAN{1}\cup \{ x_2x_4x_5 \}$ we get an extra condition $3c_3=0$, i.e.
$$
c_1=c_2=c_3=-2c_4,\; c_5=c_4,\; 6c_4=0.
$$
All such cases have already appeared.\\

For $\AA =\AAN{1}\cup \{ x_5^2x_1, x_4^2x_2 \}$ extra conditions are automatically satisfied, i.e.
$$
c_1=c_2=-2c_3,\; c_5=-c_3-c_4,\; 6c_3=2(c_4-c_3)=0.
$$
This means that $\GG\cong \ZZ{6}\oplus \ZZ{2}$ with generators
$$
f_1=(1,{\om}^{-2},{\om}^{-2},\om,\om,{\om}^{-2}), \; f_2=(1,1,1,1,\eta,\eta),\; \om=\sqr{6}, \eta=\sqr{2}.
$$
Then $i_3+i_4 \equiv 2(i_1+i_2+i_5) \; mod \; 6$, $i_4+i_5 \equiv 0 \; mod \; 2$. This means that $\AAbar=\AA\cup \{ f_3(x_1,x_2) \}$.\\

For $\AAN{1}=\AAN{0}\cup \{ x_0x_3x_4 \}$, $\AAN{1}=\AAN{0}\cup \{ x_0x_3x_5 \}$, $\AAN{1}=\AAN{0}\cup \{ x_0x_4x_5 \}$ we get extra conditions $c_4=-c_3$, $c_5=-c_3$, $c_5=-c_4$ respectively. Hence upto a permutation of $x_i$ we get
$$
c_1=c_2=0,\; c_4=c_3,\; 2c_3=2c_5=0.
$$
All such cases have already appeared.\\

Let $\AAN{0} =\{ x_1^2x_2, x_2^2x_1, x_3^2x_1, x_4^2x_1, x_5^2x_2, x_0x_3x_4, x_2x_3x_5, x_0^3 \}$. Then 
$$
c_1=c_2=0,\; c_3=c_4=c_5,\; 2c_3=0.
$$
All such cases have already appeared.\\

Let $\AAN{0} =\{ x_1^2x_2, x_2^2x_1, x_3^2x_1, x_4^2x_1, x_5^2x_2, x_3x_4x_5, x_0^3 \}$ or $\AAN{0} =\{ x_1^2x_2,$ $x_2^2x_1,$ $x_3^2x_1,$ $x_4^2x_1,$ $x_5^2x_1,$ $x_3x_4x_5,$ $x_0^3 \}$. Then 
$$
c_1=c_2=-2c_3,\; c_5=-c_3-c_4,\; 6c_3=2(c_4-c_3)=0.
$$

For $\AAN{1}=\AAN{0}\cup \{ x_2^2x_0 \}$, $\AAN{1}=\AAN{0}\cup \{ x_0x_2x_3 \}$ we get extra conditions $2c_3=0$, $c_3=0$ respectively. They both imply that $2c_3=0$, i.e.
$$
c_1=c_2=0,\; c_5=c_3+c_4,\; 2c_3=2c_4=0.
$$
All such cases have already appeared.\\

For $\AAN{1}=\AAN{0}\cup \{ x_2x_3x_5 \}$ we get an extra condition $c_5=c_3$, i.e.
$$
c_1=c_2=c_4=-2c_3,\; c_5=c_3,\; 6c_3=0.
$$
All such cases have already appeared.\\

For $\AAN{1}=\AAN{0}\cup \{ x_2x_3x_4 \}$ we get an extra condition $c_4=c_3$, i.e.
$$
c_1=c_2=c_5=-2c_3,\; c_4=c_3,\; 6c_3=0.
$$
All such cases have already appeared.\\

Let $\AAN{0} =\{ x_1^2x_2, x_2^2x_1, x_3^2x_1, x_4^2x_0, x_5^2x_0, x_3x_4x_5, x_0^3 \}$. Then 
$$
c_1=c_2=0,\; c_5=c_3+c_4,\; 2c_3=2c_4=0.
$$
All such cases have already appeared.\\

Let $\AAN{0} =\{ x_1^2x_2, x_2^2x_1, x_3^2x_1, x_4^2x_0, x_5^2x_0, x_2x_4x_5, x_0^3 \}$. Then 
$$
c_1=c_2=0,\; c_5=c_4,\; 2c_3=2c_4=0.
$$
All such cases have already appeared.\\

Let $\AAN{0} =\{ x_1^2x_2, x_2^2x_1, x_3^2x_1, x_3^2x_2, x_4^2x_0, x_5^2x_0, x_0^3 \}$. Then 
$$
c_1=c_2=-2c_3,\; 6c_3=2c_4=2c_5=0.
$$

For $\AAN{1}=\AAN{0}\cup \{ x_1x_4x_5 \}$ and $\AA=\AAN{0}\cup \{ x_3x_4x_5 \}$ we get extra conditions $c_4+c_5=2c_3$ and $c_4+c_5=-c_3$ respectively. They both imply that $2c_3=0$, i.e.
$$
c_1=c_2=0,\; 2c_3=2c_4=2c_5=0.
$$
All such cases have already appeared.\\

Let $\AAN{0} =\{ x_1^2x_2, x_2^2x_1, x_3^2x_1, x_4^2x_1, x_5^2x_1, x_0x_3x_5, x_0^3 \}$. Then 
$$
c_1=c_2=0,\; c_5=c_3,\; 2c_3=2c_4=0.
$$
All such cases have already appeared.\\

If $\AA =\{ x_1^2x_2, x_2^2x_1, x_3^2x_1, x_4^2x_0, x_5^2x_0, x_1x_4x_5, x_2x_3x_4, x_0^3 \}$ or $\AA =\{ x_1^2x_2,$ $x_2^2x_1,$ $x_3^2x_0,$ $x_4^2x_0,$ $x_5^2x_0,$ $x_1x_3x_5,$ $x_1x_4x_5,$ $x_2x_3x_4,$ $x_0^3 \}$, then 
$$
c_1=c_2=0,\; c_5=c_4=c_3,\; 2c_3=0.
$$
This case has already appeared.\\

If $\AA =\{ x_1^2x_2, x_2^2x_1, x_3^2x_0, x_4^2x_0, x_5^2x_0, x_3x_4x_5, x_0^3 \}$, then 
$$
c_1=c_2,\; c_5=c_3+c_4,\; 3c_1=2c_3=2c_4=0.
$$
This means that $\GG\cong \ZZ{3}\oplus (\ZZ{2})^{\oplus 2}$ with generators
$$
f_1=(1,{\om},{\om},1,1,1), \; f_2=(1,1,1,\eta,1,\eta),\; f_3=(1,1,1,1,\eta,\eta),\; \om=\sqr{3}, \eta=\sqr{2}.
$$
Then $i_1+i_2 \equiv 0 \; mod \; 3$, $i_3+i_5 \equiv i_4+i_5 \equiv 0 \; mod \; 2$. This means that $\AAbar=\AA\cup \{ f_3(x_1,x_2) \}$.\\

\subsubsection{Case of $1$ cube. No cycles.}

If there are no cycles at all, then $\AA$ may be either

\begin{center}

\end{center}

or one of the following:

\begin{center}
%
\end{center}

Consider the subgraph formed by vertices $x_1,x_2,x_3,x_4,x_5$. The maximal length of a path in this subgraph may be either $5$ or $4$ or $3$ or $2$ or $1$ (in which case the subgraph is totally disconnected).\\

If this length is equal to $5$, then we may assume that $1$ is connected to $2$, $2$ is connected to $3$, $3$ is connected to $4$ and $4$ is connected to $5$. By Lemma 1 (\cite{Liendo}, Lemma 1.3) the $5$-th vertex should be connected either to the cube $0$ or to one of the vertices $1$, $2$, $3$, $4$. The latter is not possible since we assume that there are no cycles. In the former case there are neither singular pairs nor singular triples.\\ 

If the maximal length of a path in the subgraph formed by $1$, $2$, $3$, $4$, $5$ is $4$, then we may assume that $1$ is connected to $2$, $2$ is connected to $3$ and $3$ is connected to $4$. By Lemma 1 (\cite{Liendo}, Lemma 1.3) the $4$-th vertex should be connected to the cube $0$. It can be connected neither to $5$ (since there are no length $5$ paths) nor to any of the vertices $1$, $2$, $3$ (since there are no cycles).\\

By Lemma 1 (\cite{Liendo}, Lemma 1.3) the $5$-th vertex should be connected either to the $2$-nd vertex or to the $3$-rd vertex or to the $4$-th vertex or to the cube $0$. It can not be connected to the $1$-st vertex, because there are no length $5$ paths.\\

If the $5$-th vertex is connected to the $2$-nd vertex, we may get 
\begin{itemize}
\item one singular pair $(x_1,x_5)$ and 
\item two singular triples $(x_1,x_5, x_3)$, $(x_1,x_5, x_4)$.
\end{itemize}

Note that this configuration is symmetric with respect to vertices $1$ and $5$.\\

If the $5$-th vertex is connected to the $3$-rd vertex, we may get 
\begin{itemize}
\item one singular pair $(x_2,x_5)$ and 
\item one singular triple $(x_2,x_4, x_5)$.
\end{itemize}

If the $5$-th vertex is connected to the $4$-th vertex, we may get 
\begin{itemize}
\item one singular pair $(x_3,x_5)$ and 
\item one singular triple $(x_1,x_3, x_5)$.
\end{itemize}

If the $5$-th vertex is connected to the cube $0$, we may get 
\begin{itemize}
\item one singular pair $(x_4,x_5)$ and 
\item two singular triples $(x_1,x_4, x_5)$, $(x_2,x_4, x_5)$.
\end{itemize}

All these singular pairs and singular triples can be resolved as in the earlier examples and we obtain the pictures shown above.\\

If the maximal length of a path in the subgraph formed by $1$, $2$, $3$, $4$, $5$ is $3$, then we may assume that $1$ is connected to $2$, $2$ is connected to $3$ and $3$ is connected to the cube $0$. Indeed, $3$ can be connected neither to $4$ or $5$ (since there are no length $4$ paths) nor to $1$ or $2$ (since there are no cycles).\\

The remaining vertices $4$ and $5$ may be connected by an edge or disconnected.\\

If $4$ is connected to $5$, then by Lemma 1 (\cite{Liendo}, Lemma 1.3) the $5$-th vertex should be connected either to $3$ or to the cube $0$. It can be connected neither to $1$ (since there are no length $5$ paths) nor to $2$ (since there are no length $4$ paths) nor to $4$ (since there are no cycles).\\

If the $5$-th vertex is connected to $3$, we may get one singular pair $(x_2,x_5)$.\\

If the $5$-th vertex is connected to the cube $0$, we may get one singular pair $(x_3,x_5)$ and one singular triple $(x_1,x_3, x_5)$.\\

These singular pairs and singular triples can be resolved as usual and we obtain the following pictures of $\AA$:

\begin{center}

\end{center}

If the remaining vertices $4$ and $5$ are disconnected, then by Lemma 1 (\cite{Liendo}, Lemma 1.3) each of them should be connected to either $2$ or $3$ or $0$. They can not be connected to $1$, because there are no length $4$ paths.\\

If they are connected to a different pair of these three vertices $0$, $2$, $3$, then $\AA$ is one of the following:

\begin{center}
%
\end{center}

Indeed, if $4$ is connected to $2$ and $5$ is connected to $3$, then we may get 
\begin{itemize}
\item two singular pairs $(x_1,x_4)$, $(x_2,x_5)$ and 
\item two singular triples $(x_1,x_4, x_3)$, $(x_1,x_4, x_5)$.
\end{itemize}

If $4$ is connected to $2$ and $5$ is connected to $0$, then we may get 
\begin{itemize}
\item two singular pairs $(x_1,x_4)$, $(x_3,x_5)$ and 
\item four singular triples $(x_1,x_4, x_3)$, $(x_1,x_4, x_5)$, $(x_1,x_3, x_5)$, $(x_3,x_4, x_5)$.
\end{itemize}

If $4$ is connected to $3$ and $5$ is connected to $0$, then we may get 
\begin{itemize}
\item two singular pairs $(x_2,x_4)$, $(x_3,x_5)$ and 
\item two singular triples $(x_2,x_4, x_5)$, $(x_1,x_3, x_5)$.
\end{itemize}

We resolve these singular pairs and triples as usual and obtain the pictures shown above.\\

If the remaining vertices $4$ and $5$ are both connected to the cube $0$, then we may get 
\begin{itemize}
\item three singular pairs $(x_3,x_4)$, $(x_3,x_5)$, $(x_4,x_5)$ and 
\item five singular triples $(x_3,x_4, x_5)$, $(x_2,x_4, x_5)$, $(x_1,x_4, x_5)$, $(x_1,x_3, x_4)$, $(x_1,x_3, x_5)$.
\end{itemize}

After resolving them as usual we obtain the following pictures of $\AA$:

\begin{center}

\end{center}

If the remaining vertices $4$ and $5$ are either both connected to $2$ or both connected to $3$, then $\AA$ is one of the following:

\begin{center}
%
\end{center}

Indeed, if $4$ and $5$ are both connected to the $3$-rd vertex, then we may get 
\begin{itemize}
\item three singular pairs $(x_2,x_4)$, $(x_2,x_5)$, $(x_4,x_5)$ and 
\item two singular triples $(x_1,x_4, x_5)$, $(x_2,x_4, x_5)$.
\end{itemize}

Note that this configuration is symmetric with respect to vertices $4$ and $5$.\\

If $4$ and $5$ are both connected to the $2$-nd vertex, then we may get 
\begin{itemize}
\item three singular pairs $(x_1,x_4)$, $(x_1,x_5)$, $(x_4,x_5)$ and 
\item four singular triples $(x_1,x_4, x_5)$, $(x_3,x_4, x_5)$, $(x_1,x_3, x_4)$, $(x_1,x_3, x_5)$.
\end{itemize}

Note that this configuration is symmetric with respect to vertices $1$, $4$ and $5$.\\

After resolving these singular pairs and triples as usual we obtain the pictures shown above.\\

Suppose that the maximal length of a path in the subgraph formed by $1$, $2$, $3$, $4$, $5$ is $2$. We may assume that $1$ is connected to $2$ (which then should be connected to $0$ by our assumptions).\\

The remaining vertices $3$, $4$ and $5$ may also form a length $2$ path. If they do (say, $5$ is connected to $4$), then by Lemma 1 (\cite{Liendo}, Lemma 1.3) the $4$-th vertex can be connected only to the cube $0$ (since otherwise we would get a length $3$ path or a cycle).\\

In this case the $3$-rd vertex can be connected either to the cube $0$ or to $2$ or to $4$.\\ 

If it is connected to $0$, then $\AA$ can be one of the following:

\begin{center}

\end{center}

Indeed, in this case we may get 
\begin{itemize}
\item three singular pairs $(x_2,x_3)$, $(x_2,x_4)$, $(x_3,x_4)$ and 
\item three singular triples $(x_2,x_3, x_4)$, $(x_2,x_3, x_5)$, $(x_1,x_3, x_4)$.
\end{itemize}

These singular pairs and triples can be resolved as usual and we obtain the pictures shown above.\\

If the $3$-rd vertex is connected to $2$ or to $4$, then by symmetry we may assume that it is connected to $2$.\\

In this case we may get 
\begin{itemize}
\item two singular pairs $(x_1,x_3)$, $(x_2,x_4)$ and 
\item two singular triples $(x_1,x_3, x_4)$, $(x_1,x_3, x_5)$.
\end{itemize}

We resolve them as usual and obtain the following pictures:

\begin{center}
\begin{tikzpicture}

\coordinate (0) at (0,0); \node[above, font=\scriptsize] at (0) {0}; \fill (0) circle (1pt); \draw (0) circle (1mm); 
\coordinate (1) at (-2,0); \node[above, font=\scriptsize] at (1) {1}; \fill (1) circle (1pt);
\coordinate (2) at (-1,0); \node[above, font=\scriptsize] at (2) {2}; \fill (2) circle (1pt);
\coordinate (3) at (-1,-1); \node[below, font=\scriptsize] at (3) {3}; \fill (3) circle (1pt);
\coordinate (4) at (1,0); \node[above, font=\scriptsize] at (4) {4}; \fill (4) circle (1pt);
\coordinate (5) at (2,0); \node[above, font=\scriptsize] at (5) {5}; \fill (5) circle (1pt);

\arr{1}{2};\arr{3}{2};\arr{5}{4};\draw(2)--(0);\draw(4)--(0);

\node[right,draw=none,text width=12cm, font=\small, scale=.8] at (2.3,0) { $\cup (x_1^2x_4, x_1x_3x_0, x_1x_3x_4, x_1x_3x_5) \cup (x_1^2x_4, x_3^2x_4, x_4x_3x_5, x_1x_3x_4, x_1x_3x_5$, $x_1x_4x_5)\cup (x_1x_3x_5, x_1x_3x_0, x_0x_3x_5, x_1x_0x_5)\cup (x_2x_4x_5, x_1x_2x_4, x_2x_3x_4)$ };

\end{tikzpicture}
\end{center}

Suppose now that the remaining vertices $3$, $4$ and $5$ do not form length $2$ paths (i.e. they are totally disconnected).\\

By Lemma 1 (\cite{Liendo}, Lemma 1.3) and our assumptions each of them should be connected either to the cube $0$ or to the $2$-nd vertex.\\

If $4$ and $5$ are connected to $0$ and $3$ is connected to $2$, then we may get 
\begin{itemize}
\item four singular pairs $(x_1,x_3)$, $(x_2,x_4)$, $(x_2,x_5)$, $(x_4,x_5)$ and 
\item five singular triples $(x_1,x_3, x_4)$, $(x_1,x_3, x_5)$, $(x_1,x_4, x_5)$, $(x_3,x_4, x_5)$, $(x_2,x_4, x_5)$.
\end{itemize}

Note that this configuration is symmetric with respect to vertices $1$, $3$ as well as with respect to vertices $4$, $5$.\\

If $4$ and $3$ are connected to $2$ and $5$ is connected to $0$, then we may get 
\begin{itemize}
\item four singular pairs $(x_1,x_3)$, $(x_1,x_4)$, $(x_3,x_4)$, $(x_2,x_5)$ and 
\item four singular triples $(x_1,x_3, x_4)$, $(x_1,x_3, x_5)$, $(x_1,x_4, x_5)$, $(x_3,x_4, x_5)$.
\end{itemize}

Note that this configuration is symmetric with respect to vertices $1$, $3$, $4$.\\

We resolve these singular pairs and triples as usual and obtain the following pictures of $\AA$:

\begin{center}



\begin{flalign*}
\mbox{where}\; & B=(x_5^2x_2,\; x_2x_3x_5,\; x_2x_5x_4,\; x_1x_2x_5).\\
\end{flalign*}
\end{center}

If all three vertices $3$, $4$ and $5$ are connected to the cube $0$, then we may get 
\begin{itemize}
\item six singular pairs $(x_2,x_3)$, $(x_2,x_4)$, $(x_2,x_5)$, $(x_3,x_4)$, $(x_3,x_5)$, $(x_4,x_5)$ and 
\item seven singular triples $(x_2,x_3, x_4)$, $(x_2,x_3, x_5)$, $(x_2,x_4, x_5)$, $(x_3,x_4, x_5)$, $(x_1,x_3, x_4)$,\\
$(x_1,x_3, x_5)$, $(x_1,x_4, x_5)$.
\end{itemize}

After we resolve them as usual we obtain the following pictures:

\begin{center}


\begin{flalign*}
\mbox{where}\; & C=(x_1x_5x_2,\; x_3x_5x_2,\; x_4x_5x_2)\cup (x_1x_4x_2,\; x_2x_4x_5,\; x_2x_4x_3).\\
\end{flalign*}
\end{center}

If all three vertices $3$, $4$ and $5$ are connected to the $2$-nd vertex, then we may get 
\begin{itemize}
\item six singular pairs $(x_1,x_3)$, $(x_1,x_4)$, $(x_1,x_5)$, $(x_3,x_4)$, $(x_3,x_5)$, $(x_4,x_5)$ and 
\item four singular triples $(x_1,x_3, x_4)$, $(x_1,x_3, x_5)$, $(x_1,x_4, x_5)$, $(x_3,x_4, x_5)$.
\end{itemize}

After we resolve them as usual we obtain the following pictures:

\begin{center}
%
\end{center}

Finally, it may happen that the subgraph formed by vertices $1$, $2$, $3$, $4$, $5$ is totally disconnected. In this case each of them should be connected to the cube $0$ by Lemma 1 (\cite{Liendo}, Lemma 1.3).\\

We get 
\begin{itemize}
\item ten singular pairs $(x_1,x_2)$, $(x_1,x_3)$, $(x_1,x_4)$, $(x_1,x_5)$, $(x_2,x_3)$, $(x_2,x_4)$, $(x_2,x_5)$, $(x_3,x_4)$, $(x_3,x_5)$, $(x_4,x_5)$ and 
\item ten singular triples $(x_1,x_2, x_3)$, $(x_1,x_2, x_4)$, $(x_1,x_2, x_5)$, $(x_1,x_3, x_4)$, $(x_1,x_3, x_5)$,\\
$(x_1,x_4, x_5)$, $(x_2,x_3, x_4)$, $(x_2,x_3, x_5)$, $(x_2,x_4, x_5)$,  $(x_3,x_4, x_5)$.
\end{itemize}

After we resolve them as usual we obtain the following pictures of $\AA$:

\begin{center}
\begin{tikzpicture}

\coordinate (0) at (0,0); \node[above right, font=\scriptsize] at (0) {0}; \fill (0) circle (1pt); \draw (0) circle (1mm); 
\coordinate (1) at (-1,-1); \node[below left, font=\scriptsize] at (1) {1}; \fill (1) circle (1pt);
\coordinate (2) at (-1,0); \node[above, font=\scriptsize] at (2) {2}; \fill (2) circle (1pt);
\coordinate (3) at (0,1); \node[above, font=\scriptsize] at (3) {3}; \fill (3) circle (1pt);
\coordinate (4) at (1,0); \node[below, font=\scriptsize] at (4) {4}; \fill (4) circle (1pt);
\coordinate (5) at (0,-1); \node[below, font=\scriptsize] at (5) {5}; \fill (5) circle (1pt);

\draw(1)--(0);\draw(2)--(0);\draw(3)--(0);\draw(4)--(0);\draw(5)--(0);

\draw[dashed] plot [smooth, tension=.5] coordinates {(-1,-1.3)(-1,-1)(-1,0)(0,1)(.3,1.3)}; 
\draw[dashed] plot [smooth, tension=.5] coordinates {(-.3,1.3)(0,1)(1,0)(0,-1)(-.3,-1.3)}; 
\draw[dashed] plot [smooth, tension=.5] coordinates {(-1.3,-.7)(-1,-1)(0,-1)(1.5,-.3)(1,0)(.7,.2)}; 

\node[right,draw=none,text width=10cm, font=\small, scale=.8] at (1.7,0) { $\cup (x_1x_2x_5, x_3x_2x_5, x_4x_2x_5) \cup (x_1x_2x_4, x_3x_2x_4, x_4x_2x_5)$  };

\end{tikzpicture}
\end{center}

Now let us computations.\\

If $\AA =\{ x_1^2x_2, x_2^2x_3, x_3^2x_4, x_4^2x_5, x_5^2x_0, x_0^3 \}$, then 
$$
c_2=-2c_1,\; c_3=4c_1,\; c_4=-8c_1,\; c_5=16c_1,\; 32c_1=0.
$$
This means that $\GG\cong \ZZ{32}$ with a generator
$$
f=(1,{\om},{\om}^{-2},{\om}^{4},{\om}^{-8},{\om}^{16}), \; \om=\sqr{32}.
$$
Then $i_1-2i_2+4i_3-8i_4 \equiv 16i_5 \; mod \; 32$. This means that $\AAbar=\AA$.\\ 

Let $\AAN{0} =\{ x_1^2x_2, x_2^2x_3, x_3^2x_4, x_4^2x_0, x_5^2x_2, x_0^3 \}$. Then 
$$
c_2=-2c_1,\; c_3=4c_1,\; c_4=-8c_1,\; 16c_1=2(c_5-c_1)=0.
$$

For $\AA=\AAN{0}\cup \{ x_1x_3x_5 \}$, $\AAN{1}=\AAN{0}\cup \{ x_1x_4x_5 \}$ we get extra conditions $c_5=-5c_1$, $c_5=7c_1$ respectively. They all imply that $4c_1=0$ and $c_5=-c_1$, i.e.
$$
c_2=2c_1,\; c_3=c_4=0,\; c_5=-c_1,\; 4c_1=0.
$$
All such cases have already appeared.\\

For $\AAN{1}=\AAN{0}\cup \{ x_1x_3x_4 \}$, $\AA=\AAN{0}\cup \{ x_5^2x_3 \}$, $\AAN{1}=\AAN{0}\cup \{ x_5^2x_4 \}$, $\AAN{1}=\AAN{0}\cup \{ x_5^2x_0 \}$ we get extra conditions $3c_1=0$, $c_3=-2c_1$, $6c_1=0$, $2c_1=0$ respectively. They all imply that $2c_1=0$, i.e.
$$
c_2=c_3=c_4=0,\; 2c_1=2c_5=0.
$$
All such cases have already appeared.\\

Let $\AAN{0} =\{ x_1^2x_2, x_2^2x_3, x_3^2x_4, x_4^2x_0, x_5^2x_3, x_0^3 \}$. Then 
$$
c_2=-2c_1,\; c_3=4c_1,\; c_4=-8c_1,\; 16c_1=2(c_5+2c_1)=0.
$$

For $\AA=\AAN{0}\cup \{ x_2x_4x_5 \}$ and $\AAN{1}=\AAN{0}\cup \{ x_0x_2x_5 \}$ we get the same extra conditions $c_5=-6c_1$ and $c_5=2c_1$, i.e.
$$
c_2=-2c_1,\; c_3=4c_1,\; c_4=0,\; c_5=2c_1,\; 8c_1=0.
$$
All such cases have already appeared.\\

For $\AA=\AAN{0}\cup \{ x_2^2x_4 \}$, $\AA=\AAN{0}\cup \{ x_1x_2x_5 \}$, $\AAN{1}=\AAN{0}\cup \{ x_2^2x_0 \}$, $\AA=\AAN{0}\cup \{ x_5^2x_4 \}$ we get extra conditions $4c_1=0$, $c_5=c_1$, $4c_1=0$, $4c_1=0$ respectively. They all imply that $4c_1=0$, i.e.
$$
c_2=2c_1,\; c_3=c_4=0,\; 4c_1=2c_5=0.
$$
All such cases have already appeared.\\

Let $\AAN{0} =\{ x_1^2x_2, x_2^2x_3, x_3^2x_4, x_4^2x_0, x_5^2x_0, x_5^2x_3, x_0^3 \}$. Then 
$$
c_2=2c_1,\; c_3=c_4=0,\; 4c_1=2c_5=0.
$$
All such cases have already appeared.\\

Let $\AAN{0} =\{ x_1^2x_2, x_2^2x_3, x_3^2x_4, x_4^2x_0, x_5^2x_4, x_0^3 \}$. Then 
$$
c_2=-2c_1,\; c_3=4c_1,\; c_4=-8c_1,\; 16c_1=2(c_5-4c_1)=0.
$$

For $\AA=\AAN{0}\cup \{ x_5^2x_0 \}$ we get an extra condition $8c_1=0$, i.e.
$$
c_2=-2c_1,\; c_3=4c_1,\; c_4=0,\; 8c_1=2c_5=0.
$$
This means that $\GG\cong \ZZ{8}\oplus \ZZ{2}$ with generators
$$
f_1=(1,{\om},{\om}^{-2},{\om}^{4},1,1), \; f_2=(1,1,1,1,1,\eta), \; \om=\sqr{8}, \eta=\sqr{2}.
$$
Then $i_1-2i_2 \equiv 4i_3 \; mod \; 8$, $i_5 \equiv 0 \; mod \; 2$. This means that $\AAbar=\AA\cup \{ f_3(x_0,x_4), x_3^2\cdot f_1(x_0,x_4) \}$.\\

For $\AA=\AAN{0}\cup \{ x_0x_3x_5 \}$ we get an extra condition $c_5=-4c_1$, i.e.
$$
c_2=-2c_1,\; c_3=4c_1,\; c_4=-8c_1,\; c_5=-4c_1,\; 16c_1=0.
$$
This means that $\GG\cong \ZZ{16}$ with a generator
$$
f=(1,{\om},{\om}^{-2},{\om}^{4},{\om}^{8},{\om}^{-4}), \; \om=\sqr{16}.
$$
Then $i_1-2i_2+4i_3 \equiv 4i_5+8i_4 \; mod \; 16$. This means that $\AAbar=\AA$.\\

For $\AAN{1}=\AAN{0}\cup \{ x_2x_3x_5 \}$, $\AA=\AAN{0}\cup \{ x_1x_3x_5 \}$, $\AA=\AAN{0}\cup \{ x_3^2x_0 \}$ we get extra conditions $c_5=-2c_1$, $c_5=-5c_1$, $8c_1=0$ respectively. They all imply that $8c_1=0$, i.e.
$$
c_2=-2c_1,\; c_3=4c_1,\; c_4=0,\; 8c_1=2c_5=0.
$$
All such cases have already appeared.\\

Let $\AAN{0} =\{ x_1^2x_2, x_2^2x_3, x_3^2x_4, x_4^2x_0, x_5^2x_0, x_0^3 \}$. Then 
$$
c_2=-2c_1,\; c_3=4c_1,\; c_4=8c_1,\; 16c_1=2c_5=0.
$$

For $\AAN{1}=\AAN{0}\cup \{ x_2x_4x_5 \}$, $\AAN{1}=\AAN{0}\cup \{ x_3x_4x_5 \}$, $\AA=\AAN{0}\cup \{ x_1x_4x_5 \}$ we get extra conditions $c_5=6c_1$, $c_5=4c_1$, $c_5=7c_1$. They all imply that $8c_1=0$, i.e.
$$
c_2=-2c_1,\; c_3=4c_1,\; c_4=0,\; 8c_1=2c_5=0.
$$
All such cases have already appeared.\\

Let $\AAN{0} =\{ x_1^2x_2, x_2^2x_3, x_3^2x_0, x_4^2x_5, x_5^2x_3, x_0^3 \}$. Then 
$$
c_2=-2c_1,\; c_3=4c_1,\; c_5=-2c_4,\; 8c_1=4(c_4-c_1)=0.
$$

For $\AA=\AAN{0}\cup \{ x_0x_2x_5 \}$ we get an extra condition $2(c_4+c_1)=0$, i.e.
$$
c_2=-2c_1,\; c_3=4c_1,\; c_5=2c_1,\; 8c_1=2(c_4+c_1)=0.
$$
This means that $\GG\cong \ZZ{8}\oplus \ZZ{2}$ with generators
$$
f_1=(1,{\om},{\om}^{-2},{\om}^{4},{\om}^{-1},{\om}^{2}), \; f_2=(1,1,1,1,\eta,1),\; \om=\sqr{8}, \eta=\sqr{2}.
$$
Then $i_1-2i_2+2i_5 \equiv 4i_3+i_4 \; mod \; 8$, $i_4 \equiv 0 \; mod \; 2$. This means that $\AAbar=\AA$.\\

For $\AA=\AAN{0}\cup \{ x_2x_4x_5 \}$ we get an extra condition $c_4=-2c_1$, i.e.
$$
c_2=c_4=2c_1,\; c_3=c_5=0,\; 4c_1=0.
$$
This case has already appeared.\\

For $\AA=\AAN{0}\cup \{ x_5^2x_0 \}$ we get an extra condition $4c_4=0$, i.e.
$$
c_2=2c_1,\; c_3=0,\; c_5=2c_4,\; 4c_1=4c_4=0.
$$
This means that $\GG\cong (\ZZ{4})^{\oplus 2}$ with generators
$$
f_1=(1,{\om},{\om}^{2},1,1,1), \; f_2=(1,1,1,1,{\om},{\om}^{2}),\; \om=\sqr{4}.
$$
Then $i_1 \equiv 2i_2 \; mod \; 4$, $i_4 \equiv 2i_5 \; mod \; 4$. This means that $\AAbar=\AA\cup \{  f_3(x_0,x_3), x_2^2\cdot f_1(x_0,x_3) \}$.\\

Let $\AAN{0} =\{ x_1^2x_2, x_2^2x_3, x_3^2x_0, x_4^2x_5, x_5^2x_0, x_0^3 \}$. Then 
$$
c_2=-2c_1,\; c_3=4c_1,\; c_5=2c_4,\; 8c_1=4c_4=0.
$$

For $\AA=\AAN{0}\cup \{ x_3x_4x_5 \}$ we get an extra condition $c_4=4c_1$, i.e.
$$
c_2=-2c_1,\; c_3=c_4=4c_1,\; c_5=0,\; 8c_1=0.
$$
This case has already appeared.\\

For $\AA=\AAN{0}\cup \{ x_1x_3x_5 \}$ we get an extra condition $c_5=3c_1$, i.e.
$$
c_2=c_3=0,\; c_1=c_5=2c_4,\; 4c_4=0.
$$
This case has already appeared.\\

For $\AAN{1}=\AAN{0}\cup \{ x_2x_3x_5 \}$ we get an extra condition $c_5=-2c_1$, i.e.
$$
c_2=c_5=2c_1,\; c_3=0,\; 4c_1=2(c_4-c_1)=0.
$$
All such cases either have already appeared or will appear and will be analyzed below.\\

Let $\AAN{0} =\{ x_1^2x_2, x_2^2x_3, x_3^2x_0, x_4^2x_2, x_5^2x_3, x_0^3 \}$. Then 
$$
c_2=-2c_1,\; c_3=4c_1,\; 8c_1=2(c_4-c_1)=2(c_5-2c_1)=0.
$$

For $\AAN{1}=\AAN{0}\cup \{ x_4^2x_3 \}$ and $\AAN{1}=\AAN{0}\cup \{ x_4^2x_0 \}$ we get the same extra condition $c_3+2c_4=2c_4=0$, i.e.
$$
c_2=c_3=0,\; 2c_1=2c_4=2c_5=0.
$$
All such cases have already appeared.\\

For $\AAN{1}=\AAN{0}\cup \{ x_0x_1x_4 \}$ and $\AAN{1}=\AAN{0}\cup \{ x_1x_3x_4 \}$ we get the same extra condition $c_1+c_4=c_1+c_3+c_4=0$, i.e.
$$
c_2=2c_1,\; c_3=0,\; c_4=-c_1\; 4c_1=2c_5=0.
$$
All such cases have already appeared.\\

For $\AAN{1}=\AAN{0}\cup \{ x_1x_4x_5 \}$ we get an extra condition $c_1+c_4+c_5=0$, i.e.
$$
c_2=-2c_1,\; c_3=4c_1,\; c_5=-c_1-c_4,\; 8c_1=2(c_4-c_1)=0.
$$

For $\AA=\AAN{1}\cup \{ x_5^2x_2 \}$, $\AA=\AAN{1}\cup \{ x_1x_2x_5 \}$, $\AA=\AAN{1}\cup \{ x_2x_4x_5 \}$ we get extra conditions $2c_1=0$, $2c_5=c_1$, $3c_1=0$ respectively. They all imply that $2c_1=0$, i.e.
$$
c_2=c_3=0,\; c_5=c_1+c_4,\; 2c_1=2c_4=0.
$$
These cases have already appeared.\\

For $\AA=\AAN{1}\cup \{ x_5^2x_0 \}$ and $\AA=\AAN{1}\cup \{ x_2^2x_0 \}$ we get the same extra condition $4c_1=0$, i.e.
$$
c_2=2c_1,\; c_3=0,\; c_5=-c_1-c_4,\; 4c_1=2(c_4-c_1)=0.
$$
This means that $\GG\cong \ZZ{4} \oplus \ZZ{2}$ with generators
$$
f_1=(1,{\om},{\om}^{2},1,{\om},{\om}^{2}), \; f_2=(1,1,1,1,{\om}^2,{\om}^{2}),\; \om=\sqr{4}.
$$
Then $i_1+i_4 \equiv 2(i_2+i_5) \; mod \; 4$, $i_4+i_5 \equiv 0 \; mod \; 2$. This means that $\AAbar=\AA\cup \{  f_3(x_0,x_3), x_2^2\cdot f_1(x_0,x_3),  x_5^2\cdot f_1(x_0,x_3) \}$.\\

For $\AA=\AAN{1}\cup \{ x_0x_2x_5 \}$ we get the same extra condition $c_4=-3c_1$, i.e.
$$
c_2=-2c_1,\; c_3=4c_1,\; c_4=-3c_1,\; c_5=2c_1,\; 8c_1=0.
$$
This means that $\GG\cong \ZZ{8}$ with a generator
$$
f=(1,{\om},{\om}^{-2},{\om}^{4},{\om}^5,{\om}^{2}), \; \om=\sqr{8}.
$$
Then $i_1+2i_5 \equiv 2i_2+4i_3+3i_4 \; mod \; 8$. This means that $\AAbar=\AA$.\\

Let $\AAN{0} =\{ x_1^2x_2, x_2^2x_3, x_3^2x_0, x_4^2x_2, x_5^2x_0, x_0^3 \}$. Then 
$$
c_2=-2c_1,\; c_3=4c_1,\; 8c_1=2(c_4-c_1)=2c_5=0.
$$

For $\AAN{1}=\AAN{0}\cup \{ x_5^2x_2 \}$ we get an extra condition $2c_1=0$, i.e.
$$
c_2=c_3=0,\; 2c_1=2c_4=2c_5=0.
$$
All such cases have already appeared.\\

For $\AAN{1}=\AAN{0}\cup \{ x_2x_3x_5 \}$ we get an extra condition $c_2+c_3+c_5=0$, i.e.
$$
c_2=c_5=2c_1,\; c_3=0,\; 4c_1=2(c_4-c_1)=0.
$$

For $\AA=\AAN{1}\cup \{ x_1x_4x_5 \}$ we get an extra condition $c_4=c_1$, i.e.
$$
c_2=c_5=2c_1,\; c_3=0,\; c_4=c_1,\; 4c_1=0.
$$
This means that $\GG\cong \ZZ{4}$ with a generator
$$
f=(1,{\om},{\om}^{2},1,{\om},{\om}^{2}), \; \om=\sqr{4}.
$$
Then $i_1+i_4 \equiv 2(i_2+i_5) \; mod \; 4$. This means that $\AAbar=\AA\cup \{ f_3(x_0,x_3), f_2(x_2,x_5)\cdot f_1(x_0,x_3), f_1(x_2,x_5)\cdot f_2(x_1,x_4) \}$.\\

Any other extra condition (after adding to $\AAN{1}$ a monomial $x_0^{i_0}x_1^{i_1}x_2^{i_2}x_3^{i_3}x_4^{i_4}x_5^{i_5}$) has the form $a\cdot c_1+b\cdot c_4=0$, $a,b\in \mathbb Z$, i.e. either $a\cdot c_1=0$ or $c_4=a\cdot c_1$ for some $a\in  \{ 0,1,2,3 \}$. All such cases have already appeared.\\

For $\AAN{1}=\AAN{0}\cup \{ x_1x_3x_5 \}$ we get an extra condition $c_5=3c_1$, i.e.
$$
c_2=c_3=0,\; c_5=c_1,\; 2c_1=2c_4=0.
$$
All such cases have already appeared.\\

For $\AAN{1}=\AAN{0}\cup \{ x_3x_4x_5 \}$ we get an extra condition $c_3+c_4+c_5=0$, i.e.
$$
c_2=c_3=0,\; c_4=c_5,\; 2c_1=2c_5=0.
$$
All such cases have already appeared.\\

Let $\AAN{0} =\{ x_1^2x_2, x_2^2x_3, x_3^2x_0, x_4^2x_3, x_5^2x_0, x_0^3 \}$. Then 
$$
c_2=-2c_1,\; c_3=4c_1,\; 8c_1=2(c_4-2c_1)=2c_5=0.
$$

For $\AAN{1}=\AAN{0}\cup \{ x_5^2x_3 \}$ we get an extra condition $4c_1=0$, i.e.
$$
c_2=2c_1,\; c_3=0,\; 4c_1=2c_4=2c_5=0.
$$

For $\AAN{2}=\AAN{1}\cup \{ x_2x_4x_5 \}$ we get an extra condition $c_5=c_4+2c_1$, i.e.
$$
c_2=2c_1,\; c_3=0,\; c_5=c_4+2c_1,\; 4c_1=2c_4=0.
$$
All such cases have already appeared.\\

For $\AAN{2}=\AAN{1}\cup \{ x_1x_2x_4 \}$, $\AAN{2}=\AAN{1}\cup \{ x_1x_2x_5 \}$, $\AAN{2}=\AAN{1}\cup \{ x_1x_4x_5 \}$ we get extra conditions $c_4=c_1$, $c_5=c_1$, $c_1=c_4+c_5$ respectively. They all imply that $2c_1=0$, i.e.
$$
c_2=c_3=0,\; 2c_1=2c_4=2c_5=0.
$$
All such cases have already appeared.\\

For $\AAN{1}=\AAN{0}\cup \{ x_3x_4x_5 \}$ we get an extra condition $c_5=c_3+c_4$, i.e.
$$
c_2=2c_1,\; c_3=0,\; c_4=c_5,\; 4c_1=2c_5=0.
$$

For $\AAN{2}=\AAN{1}\cup \{ x_1x_2x_4 \}$, $\AAN{2}=\AAN{1}\cup \{ x_1x_2x_5 \}$, $\AAN{2}=\AAN{1}\cup \{ x_1x_4x_5 \}$, $\AA=\AAN{1}\cup \{ x_2x_4x_5 \}$ we get extra conditions $c_4=c_1$, $c_5=c_1$, $c_1=0$, $2c_1=0$ respectively. They all imply that $2c_1=0$, i.e.
$$
c_2=c_3=0,\; c_4=c_5,\; 2c_1=2c_5=0.
$$
All such cases have already appeared.\\

For $\AAN{1}=\AAN{0}\cup \{ x_2x_3x_5 \}$ we get an extra condition $c_5=-2c_1$, i.e.
$$
c_2=c_5=2c_1,\; c_3=0,\; 4c_1=2c_4=0.
$$

For $\AAN{2}=\AAN{1}\cup \{ x_2x_4x_5 \}$ we get an extra condition $c_4=0$, i.e.
$$
c_2=c_5=2c_1,\; c_3=c_4=0,\; 4c_1=0.
$$
All such cases have already appeared.\\

For $\AAN{2}=\AAN{1}\cup \{ x_1x_2x_4 \}$, $\AAN{2}=\AAN{1}\cup \{ x_1x_4x_5 \}$, $\AAN{2}=\AAN{1}\cup \{ x_1x_2x_5 \}$ we get extra conditions $c_4=c_1$, $c_4=c_1$, $c_1=0$. They all imply that $2c_1=0$, i.e.
$$
c_2=c_5=c_3=0,\; 2c_1=2c_4=0.
$$
All such cases have already appeared.\\

For $\AAN{1}=\AAN{0}\cup \{ x_1x_3x_5 \}$ we get an extra condition $c_5=-5c_1$, i.e.
$$
c_2=c_3=0,\; c_5=c_1,\; 2c_1=2c_4=0.
$$
All such cases have already appeared.\\

Let $\AAN{0} =\{ x_1^2x_2, x_2^2x_3, x_3^2x_0, x_4^2x_0, x_5^2x_0, x_0^3 \}$. Then 
$$
c_2=-2c_1,\; c_3=4c_1,\; 8c_1=2c_4=2c_5=0.
$$

For $\AAN{1}=\AAN{0}\cup \{ x_3x_4x_5 \}$ we get an extra condition $c_5=c_4+4c_1$, i.e.
$$
c_2=-2c_1,\; c_3=4c_1,\; c_5=c_4+4c_1,\; 8c_1=2c_4=0.
$$

For $\AA=\AAN{1}\cup \{ x_2^2x_4 \}$ we get an extra condition $c_4=4c_1$, i.e.
$$
c_2=-2c_1,\; c_3=c_4=4c_1,\; c_5=0,\; 8c_1=0.
$$
This case has already appeared.\\

For $\AA=\AAN{1}\cup \{ x_1x_4x_5 \}$, $\AA=\AAN{1}\cup \{ x_1x_2x_4 \}$, $\AA=\AAN{1}\cup \{ x_2x_4x_5 \}$ we get extra conditions $5c_1=0$, $c_4=c_1$, $2c_1=0$ respectively. They all imply that $2c_1=0$, i.e.
$$
c_2=c_3=0,\; c_5=c_4,\; 2c_1=2c_4=0.
$$
All such cases have already appeared.\\

For $\AAN{1}=\AAN{0}\cup \{ x_2x_4x_5 \}$ we get an extra condition $c_5=c_4+2c_1$, i.e.
$$
c_2=2c_1,\; c_3=0,\; c_5=c_4+2c_1,\; 4c_1=2c_4=0.
$$

For $\AAN{2}=\AAN{1}\cup \{ x_2x_3x_4 \}$ we get an extra condition $c_4=2c_1$, i.e.
$$
c_2=c_4=2c_1,\; c_3=c_5=0,\; 4c_1=0.
$$
All such cases have already appeared.\\

For $\AAN{2}=\AAN{1}\cup \{ x_1x_3x_4 \}$ we get an extra condition $c_4=c_1$, i.e.
$$
c_2=c_3=0,\; c_4=c_5=c_1,\; 2c_1=0.
$$
All such cases have already appeared.\\

For $\AAN{1}=\AAN{0}\cup \{ x_2x_3x_4, x_1x_4x_5 \}$ we get extra conditions $c_4=2c_1$, $c_1=c_4+c_5$, i.e.
$$
c_2=c_3=c_4=0,\; c_5=c_1,\; 2c_1=0.
$$
All such cases have already appeared.\\

Let $\AAN{0} =\{ x_1^2x_2, x_2^2x_3, x_3^2x_0, x_4^2x_3, x_5^2x_3, x_0^3 \}$. Then 
$$
c_2=-2c_1,\; c_3=4c_1,\; 8c_1=2(c_4-2c_1)=2(c_5-2c_1)=0.
$$

For $\AAN{1}=\AAN{0}\cup \{ x_2x_4x_5 \}$ we get an extra condition $c_5=2c_1-c_4$, i.e.
$$
c_2=2c_1,\; c_3=0,\; c_5=2c_1+c_4,\; 4c_1=2c_4=0.
$$
All such cases have already appeared.\\

For $\AAN{1}=\AAN{0}\cup \{ x_0x_4x_5 \}$ we get an extra condition $c_5=-c_4$, i.e.
$$
c_2=-2c_1,\; c_3=4c_1,\; c_5=-c_4,\; 8c_1=2(c_4-2c_1)=0.
$$

For $\AAN{2}=\AAN{1}\cup \{ x_1x_4x_5 \}$ and $\AAN{2}=\AAN{1}\cup \{ x_1x_2x_4 \}$ we get extra conditions $c_1=0$ and $c_4=c_1$ respectively. They both imply that $2c_1=0$, i.e.
$$
c_2=c_3=0,\; c_5=c_4,\; 2c_1=2c_4=0.
$$
All such cases have already appeared.\\

For $\AAN{1}=\AAN{0}\cup \{ x_1x_4x_5 \}$ we get an extra condition $c_5=-c_4-c_1$, i.e.
$$
c_2=c_3=0,\; c_5=c_1+c_4,\; 2c_1=2c_4=0.
$$
All such cases have already appeared.\\

Let $\AAN{0} =\{ x_1^2x_2, x_2^2x_3, x_3^2x_0, x_4^2x_2, x_5^2x_2, x_0^3 \}$. Then 
$$
c_2=-2c_1,\; c_3=4c_1,\; 8c_1=2(c_4-c_1)=2(c_5-c_1)=0.
$$

For $\AA=\AAN{0}\cup \{ x_1x_4x_5 \}$ we get an extra condition $c_5=-c_4-c_1$, i.e.
$$
c_2=c_3=0,\; c_5=c_1+c_4,\; 2c_1=2c_4=0.
$$
This case has already appeared.\\

For $\AAN{1}=\AAN{0}\cup \{ x_3x_4x_5 \}$ we get an extra condition $c_5=4c_1-c_4$, i.e.
$$
c_2=2c_1,\; c_3=0,\; c_5=-c_4,\; 4c_1=2(c_4-c_1)=0.
$$

For $\AAN{2}=\AAN{1}\cup \{ x_1x_3x_5 \}$ and $\AAN{2}=\AAN{1}\cup \{ x_0x_1x_5 \}$ we get the same extra condition $c_4=c_1$, i.e.
$$
c_2=2c_1,\; c_3=0,\; c_4=c_1,\; c_5=-c_1,\; 4c_1=0.
$$

For $\AAN{3}=\AAN{2}\cup \{ x_1x_3x_4 \}$ and $\AAN{3}=\AAN{2}\cup \{ x_0x_1x_4 \}$ we get the same extra condition $2c_1=0$, i.e.
$$
c_2=c_3=0,\; c_4=c_5=c_1,\; 2c_1=0.
$$
All such cases have already appeared.\\

For $\AAN{1}=\AAN{0}\cup \{ x_0x_4x_5, x_1x_3x_4 \}$ we get extra conditions $c_5=-c_4$, $c_4=-c_1$, i.e.
$$
c_2=2c_1,\; c_3=0,\; c_4=-c_1,\; c_5=c_1,\; 4c_1=0.
$$

For $\AAN{2}=\AAN{1}\cup \{ x_1x_3x_5 \}$ and $\AAN{2}=\AAN{1}\cup \{ x_0x_1x_5 \}$ we get the same extra condition $2c_1=0$, i.e.
$$
c_2=c_3=0,\; c_4=c_5=c_1,\; 2c_1=0.
$$
All such cases have already appeared.\\

Let $\AAN{0} =\{ x_1^2x_2, x_2^2x_0, x_3^2x_0, x_4^2x_0, x_5^2x_4, x_0^3 \}$. Then 
$$
c_2=2c_1,\; c_4=2c_5,\; 4c_1=2c_3=4c_5=0.
$$

For $\AAN{1}=\AAN{0}\cup \{ x_2x_3x_4 \}$ we get an extra condition $c_3=2c_1+2c_5$, i.e.
$$
c_2=2c_1,\; c_3=2c_1+2c_5,\; c_4=2c_5,\; 4c_1=4c_5=0.
$$

For $\AAN{2}=\AAN{1}\cup \{ x_2x_3x_5 \}$ we get an extra condition $c_5=0$, i.e.
$$
c_2=c_3=2c_1,\; c_4=c_5=0,\; 4c_1=0.
$$
All such cases have already appeared.\\

For $\AAN{2}=\AAN{1}\cup \{ x_3^2x_2 \}$, $\AAN{2}=\AAN{1}\cup \{ x_5^2x_3 \}$ we get the same extra condition $2c_1=0$, i.e.
$$
c_2=0,\; c_3=c_4=2c_5,\; 2c_1=4c_5=0.
$$

For $\AA=\AAN{2}\cup \{ x_1x_3x_4 \}$ we get an extra condition $c_1=0$, i.e.
$$
c_1=c_2=0,\; c_3=c_4=2c_5,\; 4c_5=0.
$$
This case has already appeared.\\

For $\AA=\AAN{2}\cup \{ x_3^2x_4 \}$, $\AA=\AAN{2}\cup \{ x_1^2x_3 \}$, $\AA=\AAN{2}\cup \{ x_1^2x_4 \}$, $\AA=\AAN{2}\cup \{ x_3x_4x_5 \}$, $\AA=\AAN{2}\cup \{ x_1x_3x_5 \}$, $\AA=\AAN{2}\cup \{ x_1x_4x_5 \}$ we get extra conditions $2c_5=0$, $2c_5=0$, $2c_5=0$, $c_5=0$, $c_1=c_5$, $c_1=c_5$ respectively. They all imply that $2c_5=0$, i.e.
$$
c_2=c_3=c_4=0,\; 2c_1=2c_5=0.
$$
All such cases have already appeared.\\

For $\AAN{2}=\AAN{1}\cup \{ x_5^2x_2 \}$ we get an extra condition $2c_1=2c_5$, i.e.
$$
c_2=c_4=2c_1,\; c_3=0,\; 4c_1=2(c_5-c_1)=0.
$$
All such cases have already appeared.\\

For $\AAN{2}=\AAN{1}\cup \{ x_1x_2x_3 \}$ we get an extra condition $c_1=2c_5$, i.e.
$$
c_2=0,\; c_1=c_3=c_4=2c_5,\; 4c_5=0.
$$

For $\AA=\AAN{2}\cup \{ x_1x_3x_4 \}$, $\AA=\AAN{2}\cup \{ x_3^2x_4 \}$, $\AA=\AAN{2}\cup \{ x_1^2x_3 \}$, $\AA=\AAN{2}\cup \{ x_1^2x_4 \}$, $\AA=\AAN{2}\cup \{ x_3x_4x_5 \}$, $\AA=\AAN{2}\cup \{ x_1x_3x_5 \}$, $\AA=\AAN{2}\cup \{ x_1x_4x_5 \}$ we get extra conditions $2c_5=0$, $2c_5=0$, $2c_5=0$, $2c_5=0$, $c_5=0$, $c_5=0$, $c_5=0$ respectively. They all imply that $2c_5=0$, i.e.
$$
c_1=c_2=c_3=c_4=0,\; 2c_5=0.
$$
All such cases have already appeared.\\

For $\AAN{2}=\AAN{1}\cup \{ x_1x_2x_5 \}$ we get an extra condition $c_5=c_1$, i.e.
$$
c_2=c_4=2c_1,\; c_3=0,\; c_5=c_1,\; 4c_1=0.
$$
All such cases have already appeared.\\

For $\AA=\AAN{1}\cup \{ x_1x_3x_5 \}$ we get an extra condition $c_5=-c_1$, i.e.
$$
c_2=c_4=2c_1,\; c_3=0,\; c_5=-c_1,\; 4c_1=0.
$$
This case has already appeared.\\

For $\AAN{1}=\AAN{0}\cup \{ x_3^2x_4 \}$ we get an extra condition $2c_5=0$, i.e.
$$
c_2=2c_1,\; c_4=0,\; 4c_1=2c_3=2c_5=0.
$$

For $\AAN{2}=\AAN{1}\cup \{ x_2x_4x_5 \}$ we get an extra condition $c_5=2c_1$, i.e.
$$
c_2=c_5=2c_1,\; c_4=0,\; 4c_1=2c_3=0.
$$

For $\AA=\AAN{2}\cup \{ x_2x_3x_5 \}$ we get an extra condition $c_3=0$, i.e.
$$
c_2=c_5=2c_1,\; c_3=c_4=0,\; 4c_1=0.
$$
This case has already appeared.\\

For $\AA=\AAN{2}\cup \{ x_3^2x_2 \}$, $\AA=\AAN{2}\cup \{ x_5^2x_2 \}$, $\AA=\AAN{2}\cup \{ x_1x_2x_3 \}$, $\AA=\AAN{2}\cup \{ x_1x_3x_5 \}$, $\AA=\AAN{2}\cup \{ x_1x_2x_5 \}$ we get extra conditions $2c_1=0$, $2c_1=0$, $c_1=c_3$, $c_1=c_3$, $c_1=0$. They all imply that $2c_1=0$, i.e.
$$
c_2=c_4=c_5=0,\; 2c_1=2c_3=0.
$$
All such cases have already appeared.\\

For $\AAN{2}=\AAN{1}\cup \{ x_1x_2x_4 \}$ we get an extra condition $c_1=0$, i.e.
$$
c_1=c_2=c_4=0,\; 2c_3=2c_5=0.
$$
All such cases have already appeared.\\

For $\AAN{1}=\AAN{0}\cup \{ x_2x_3x_5, x_1x_2x_4 \}$ we get extra conditions $c_5=2c_1+c_3$ and $c_1=2c_5$, i.e.
$$
c_1=c_2=c_4=0,\; c_5=c_3,\; 2c_3=0.
$$
All such cases have already appeared.\\

For $\AAN{1}=\AAN{0}\cup \{ x_2x_4x_5, x_1x_2x_3 \}$ we get extra conditions $c_5=2c_1$ and $c_3=c_1$, i.e.
$$
c_2=c_4=c_5=0,\; c_3=c_1,\; 2c_1=0.
$$
All such cases have already appeared.\\

For $\AA=\AAN{0}\cup \{ x_2x_4x_5, x_1x_3x_4, x_2x_3x_5 \}$ we get extra conditions $c_5=2c_1$, $c_3=c_1+2c_5$ and $c_5=2c_1+c_3$, i.e.
$$
c_1=c_2=c_3=c_4=c_5=0.
$$
Hence $\GG = Id$.\\

Let $\AAN{0} =\{ x_1^2x_2, x_2^2x_0, x_3^2x_2, x_4^2x_0, x_5^2x_4, x_0^3 \}$. Then 
$$
c_2=2c_1,\; c_4=2c_5,\; 4c_1=2(c_3-c_1)=4c_5=0.
$$

For $\AAN{1}=\AAN{0}\cup \{ x_2x_4x_5 \}$ we get an extra condition $c_5=2c_1$, i.e.
$$
c_2=c_5=2c_1,\; c_4=0,\; 4c_1=2(c_3-c_1)=0.
$$

For $\AAN{2}=\AAN{1}\cup \{ x_1^2x_4 \}$ we get an extra condition $2c_1=0$, i.e.
$$
c_2=c_4=c_5=0,\; 2c_1=2c_3=0.
$$
All such cases have already appeared.\\

For $\AAN{2}=\AAN{1}\cup \{ x_1x_3x_0 \}$ and $\AAN{2}=\AAN{1}\cup \{ x_1x_3x_4 \}$ we get the same extra condition $c_3=-c_1$, i.e.
$$
c_2=c_5=2c_1,\; c_3=-c_1,\; c_4=0,\; 4c_1=0.
$$
All such cases have already appeared.\\

For $\AAN{2}=\AAN{1}\cup \{ x_1x_3x_5 \}$ we get an extra condition $c_3=c_1$, i.e.
$$
c_2=c_5=2c_1,\; c_3=c_1,\; c_4=0,\; 4c_1=0.
$$
All such cases have already appeared.\\

For $\AAN{1}=\AAN{0}\cup \{ x_1x_2x_4 \}$ we get an extra condition $c_1=2c_5$, i.e.
$$
c_2=0,\; c_1=c_4=2c_5,\; 2c_3=4c_5=0.
$$

For $\AAN{2}=\AAN{1}\cup \{ x_1^2x_4 \}$ we get an extra condition $2c_5=0$, i.e.
$$
c_1=c_2=c_4=0,\; 2c_3=2c_5=0.
$$
All such cases have already appeared.\\

For $\AAN{2}=\AAN{1}\cup \{ x_1x_3x_4 \}$ we get an extra condition $c_3=0$, i.e.
$$
c_2=c_3=0,\; c_1=c_4=2c_5,\; 4c_5=0.
$$
All such cases have already appeared.\\

For $\AAN{2}=\AAN{1}\cup \{ x_1x_3x_5 \}$ we get an extra condition $c_5=c_3$, i.e.
$$
c_1=c_2=c_4=0,\; c_5=c_3,\; 2c_3=0.
$$
All such cases have already appeared.\\

For $\AAN{2}=\AAN{1}\cup \{ x_0x_1x_3 \}$ we get an extra condition $c_3=2c_5$, i.e.
$$
c_2=0,\; c_1=c_3=c_4=2c_5,\; 4c_5=0.
$$

For $\AA=\AAN{2}\cup \{ x_3^2x_4 \}$, $\AA=\AAN{2}\cup \{ x_3x_4x_5 \}$ and $\AA=\AAN{2}\cup \{ x_1x_4x_5 \}$ we get extra conditions $2c_5=0$, $c_5=0$ and $c_5=0$ respectively. They all imply that $2c_5=0$, i.e.
$$
c_1=c_2=c_3=c_4=0,\; 2c_5=0.
$$
All such cases have already appeared.\\

For $\AAN{1}=\AAN{0}\cup \{ x_2x_3x_4 \}$ we get an extra condition $c_3=2c_1+2c_5$, i.e.
$$
c_2=0,\; c_3=c_4=2c_5,\; 2c_1=4c_5=0.
$$

For $\AAN{2}=\AAN{1}\cup \{ x_1^2x_4 \}$ we get an extra condition $2c_5=0$, i.e.
$$
c_2=c_3=c_4=0,\; 2c_1=2c_5=0.
$$
All such cases have already appeared.\\

For $\AAN{2}=\AAN{1}\cup \{ x_1x_3x_4 \}$ we get an extra condition $c_1=0$, i.e.
$$
c_1=c_2=0,\; c_3=c_4=2c_5,\; 4c_5=0.
$$
All such cases have already appeared.\\

For $\AAN{2}=\AAN{1}\cup \{ x_1x_3x_5 \}$ we get an extra condition $c_5=c_1$, i.e.
$$
c_2=c_3=c_4=0,\; c_5=c_1,\; 2c_1=0.
$$
All such cases have already appeared.\\

For $\AAN{2}=\AAN{1}\cup \{ x_0x_1x_3 \}$ we get an extra condition $c_1=2c_5$, i.e.
$$
c_2=0,\; c_1=c_3=c_4=2c_5,\; 4c_5=0.
$$

For $\AA=\AAN{2}\cup \{ x_3^2x_4 \}$, $\AA=\AAN{2}\cup \{ x_3x_4x_5 \}$ and $\AA=\AAN{2}\cup \{ x_1x_4x_5 \}$ we get extra conditions $2c_5=0$, $c_5=0$ and $c_5=0$ respectively. They all imply that $2c_5=0$, i.e.
$$
c_1=c_2=c_3=c_4=0,\; 2c_5=0.
$$
All such cases have already appeared.\\

Let $\AAN{0} =\{ x_1^2x_2, x_2^2x_0, x_3^2x_2, x_4^2x_0, x_5^2x_0, x_2x_4x_5, x_0^3 \}$. Then 
$$
c_2=2c_1,\; c_5=2c_1+c_4,\; 4c_1=2(c_3-c_1)=2c_4=0.
$$

For $\AAN{1}=\AAN{0}\cup \{ x_1^2x_0 \}$ and $\AAN{1}=\AAN{0}\cup \{ x_3^2x_0 \}$ we get the same extra condition $2c_1=0$, i.e.
$$
c_2=0,\; c_5=c_4,\; 2c_1=2c_3=2c_4=0.
$$
This case has already appeared.\\

For $\AA=\AAN{1}\cup \{ x_1x_4x_5 \}$ we get an extra condition $c_1=0$, i.e.
$$
c_1=c_2=0,\; c_5=c_4,\; 2c_3=2c_4=0.
$$
This case has already appeared.\\

For $\AA=\AAN{1}\cup \{ x_1x_3x_4 \}$ we get an extra condition $c_4=c_1+c_3$, i.e.
$$
c_2=0,\; c_5=c_4=c_1+c_3,\; 2c_1=2c_3=0.
$$
This case has already appeared.\\

For $\AA=\AAN{0}\cup \{ x_1x_3x_4 \}$ we get an extra condition $c_3=c_4-c_1$, i.e.
$$
c_2=2c_1,\;c_3=c_4-c_1,\; c_5=2c_1+c_4,\; 4c_1=2c_4=0.
$$
This means that $\GG\cong \ZZ{2} \oplus \ZZ{4}$ with generators
$$
f_1=(1,{\om},{\om}^{2},{\om}^{-1},1,{\om}^{2}), \; f_2=(1,1,1,{\om}^{2},{\om}^{2},{\om}^{2}), \; \om=\sqr{4}.
$$
Then $i_3+i_4+i_5 \equiv 0 \; mod \; 2$ and $i_1 \equiv i_3+2(i_2+i_5) \; mod \; 4$. This means that $\AAbar=\AA$.\\

For $\AAN{1}=\AAN{0}\cup \{ x_0x_1x_3 \}$ we get an extra condition $c_3=-c_1$, i.e.
$$
c_2=2c_1,\; c_3=-c_1,\; c_5=2c_1+c_4,\; 4c_1=2c_4=0.
$$

For $\AA=\AAN{1}\cup \{ x_1x_4x_5 \}$ we get an extra condition $c_1=0$, i.e.
$$
c_1=c_2=c_3=0,\; c_5=c_4,\; 2c_4=0.
$$
This case has already appeared.\\

Let $\AAN{0} =\{ x_1^2x_2, x_2^2x_0, x_3^2x_2, x_4^2x_0, x_4^2x_2, x_5^2x_0, x_0^3 \}$. Then 
$$
c_2=0,\; 2c_1=2c_3=2c_4=2c_5=0.
$$

For $\AAN{1}=\AAN{0}\cup \{ x_1x_2x_5 \}$ we get an extra condition $c_5=c_1$, i.e.
$$
c_2=0,\; c_5=c_1,\; 2c_1=2c_3=2c_4=0.
$$

For $\AAN{2}=\AAN{1}\cup \{ x_1x_4x_5 \}$ and $\AA=\AAN{1}\cup \{ x_1x_3x_5 \}$ we get the same (upto a permutation of $x_i$) extra condition $c_4=0$, i.e.
$$
c_2=c_4=0,\; c_5=c_1,\; 2c_1=2c_3=0.
$$
All such cases have already appeared.\\

For $\AAN{2}=\AAN{1}\cup \{ x_3x_4x_5 \}$ and $\AA=\AAN{1}\cup \{ x_1x_3x_4 \}$ we get the same extra condition $c_1=c_3+c_4$, i.e.
$$
c_2=0,\; c_5=c_1=c_3+c_4,\; 2c_3=2c_4=0.
$$
All such cases have already appeared.\\

For $\AAN{1}=\AAN{0}\cup \{ x_5^2x_2 \}$ the extra condition is automatically satisfied\\

For $\AAN{2}=\AAN{1}\cup \{ x_1x_3x_4 \}$, $\AAN{2}=\AAN{1}\cup \{ x_1x_3x_5 \}$, $\AAN{2}=\AAN{1}\cup \{ x_1x_4x_5 \}$ and $\AAN{2}=\AAN{1}\cup \{ x_3x_4x_5 \}$ we get the same (upto a permutation of $x_i$) extra condition $c_1=c_3+c_4$, i.e.
$$
c_2=0,\; c_1=c_3+c_4,\; 2c_3=2c_4=2c_5=0.
$$
All such cases have already appeared.\\

Let $\AAN{0} =\{ x_1^2x_2, x_2^2x_0, x_3^2x_2, x_4^2x_0, x_5^2x_0, x_1x_2x_5, x_0^3 \}$. Then 
$$
c_2=0,\; c_5=c_1,\; 2c_1=2c_3=2c_4=0.
$$

For $\AAN{1}=\AAN{0}\cup \{ x_1x_4x_5 \}$ we get an extra condition $c_4=0$, i.e.
$$
c_2=c_4=0,\; c_5=c_1,\; 2c_1=2c_3=0.
$$
All such cases have already appeared.\\

For $\AAN{1}=\AAN{0}\cup \{ x_3x_4x_5 \}$ we get an extra condition $c_5=c_3+c_4$, i.e.
$$
c_2=0,\; c_5=c_1=c_3+c_4,\; 2c_3=2c_4=0.
$$
All such cases have already appeared.\\

Let $\AAN{0} =\{ x_1^2x_2, x_2^2x_0, x_3^2x_2, x_4^2x_2, x_5^2x_0, x_1x_3x_4, x_0^3 \}$. Then 
$$
c_2=0,\; c_4=c_1+c_3,\; 2c_1=2c_3=2c_5=0.
$$
All such cases have already appeared.\\

Let $\AAN{0} =\{ x_1^2x_2, x_2^2x_0, x_3^2x_2, x_4^2x_2, x_5^2x_0, x_0x_3x_4, x_0^3 \}$. Then 
$$
c_2=2c_1,\; c_4=-c_3,\; 4c_1=2(c_3-c_1)=2c_5=0.
$$

For $\AAN{1}=\AAN{0}\cup \{ x_0x_1x_3 \}$ we get an extra condition $c_3=-c_1$, i.e.
$$
c_2=2c_1,\; c_3=-c_1,\; c_4=c_1,\; 4c_1=2c_5=0.
$$

For $\AAN{2}=\AAN{1}\cup \{ x_0x_1x_4 \}$ we get an extra condition $2c_1=0$, i.e.
$$
c_2=0,\; c_3=c_4=c_1,\; 2c_1=2c_5=0.
$$

For $\AAN{3}=\AAN{2}\cup \{ x_3x_4x_5 \}$ we get an extra condition $c_5=0$, i.e.
$$
c_2=c_5=0,\; c_3=c_4=c_1,\; 2c_1=0.
$$
All such cases have already appeared.\\

For $\AAN{2}=\AAN{1}\cup \{ x_1x_4x_5 \}$ we get an extra condition $c_5=2c_1$, i.e.
$$
c_2=c_5=2c_1,\; c_3=-c_1,\; c_4=c_1,\; 4c_1=0.
$$

For $\AAN{3}=\AAN{2}\cup \{ x_1x_2x_5 \}$, $\AAN{3}=\AAN{2}\cup \{ x_2x_3x_5 \}$, $\AAN{3}=\AAN{2}\cup \{ x_2x_4x_5 \}$ we get the same extra condition $c_1=0$, i.e.
$$
c_1=c_2=c_3=c_4=c_5=0.
$$
Hence $\GG = Id$ in all such cases.\\

For $\AAN{3}=\AAN{2}\cup \{ x_5^2x_2 \}$ we get an extra condition $2c_1=0$, i.e.
$$
c_2=c_5=0,\; c_3=c_4=c_1,\; 2c_1=0.
$$
All such cases have already appeared.\\

For $\AAN{1}=\AAN{0}\cup \{ x_1x_3x_5 \}$ we get an extra condition $c_3=c_5-c_1$, i.e.
$$
c_2=2c_1,\; c_3=c_5-c_1,\; c_4=c_5+c_1,\; 4c_1=2c_5=0.
$$

For $\AAN{2}=\AAN{1}\cup \{ x_0x_1x_4 \}$ we get an extra condition $c_5=2c_1$, i.e.
$$
c_2=c_5=2c_1,\; c_3=c_1,\; c_4=-c_1,\; 4c_1=0.
$$
All such cases have already appeared.\\

For $\AAN{2}=\AAN{1}\cup \{ x_1x_4x_5 \}$ we get an extra condition $2c_1=0$, i.e.
$$
c_2=0,\; c_3=c_4=c_1+c_5,\; 2c_1=2c_5=0.
$$
All such cases have already appeared.\\

Let $\AAN{0} =\{ x_1^2x_2, x_2^2x_0, x_3^2x_0, x_4^2x_0, x_5^2x_0, x_3x_4x_5, x_0^3 \}$. Then 
$$
c_2=2c_1,\; c_5=c_3+c_4,\; 4c_1=2c_3=2c_4=0.
$$

For $\AAN{1}=\AAN{0}\cup \{ x_1x_2x_3 \}$ we get an extra condition $c_3=c_1$, i.e.
$$
c_2=0,\; c_3=c_1=c_4+c_5,\; 2c_4=2c_5=0.
$$
All such cases have already appeared.\\

For $\AAN{1}=\AAN{0}\cup \{ x_2x_3x_4 \}$ we get an extra condition $c_4=2c_1+c_3$, i.e.
$$
c_2=c_5=2c_1,\; c_4=2c_1+c_3,\; 4c_1=2c_3=0.
$$

For $\AAN{2}=\AAN{1}\cup \{ x_1x_2x_5 \}$ we get an extra condition $c_1=0$, i.e.
$$
c_1=c_2=c_5=0,\; c_4=c_3,\; 2c_3=0.
$$
All such cases have already appeared.\\

For $\AAN{2}=\AAN{1}\cup \{ x_2x_3x_5 \}$ and $\AAN{2}=\AAN{1}\cup \{ x_2x_4x_5 \}$ we get the same (unto a permutation of $x_i$) extra condition $c_3=0$, i.e.
$$
c_2=c_4=c_5=2c_1,\; c_3=0,\; 4c_1=0.
$$

For $\AAN{3}=\AAN{2}\cup \{ x_2x_4x_5 \}$ (for $\AAN{3}=\AAN{2}\cup \{ x_2x_3x_5 \}$ in case of $\AAN{2}=\AAN{1}\cup \{ x_2x_4x_5 \}$ above) we get an extra condition $2c_1=0$, i.e.
$$
c_2=c_3=c_4=c_5=0,\; 2c_1=0.
$$
All such cases have already appeared.\\

For $\AAN{3}=\AAN{2}\cup \{ x_1x_4x_5 \}$, $\AAN{3}=\AAN{2}\cup \{ x_1x_2x_4 \}$, $\AAN{3}=\AAN{2}\cup \{ x_1x_2x_5 \}$ (for $\AAN{3}=\AAN{2}\cup \{ x_1x_3x_5 \}$, $\AAN{3}=\AAN{2}\cup \{ x_1x_2x_3 \}$, $\AAN{3}=\AAN{2}\cup \{ x_1x_2x_5 \}$ in case of $\AAN{2}=\AAN{1}\cup \{ x_2x_4x_5 \}$ above) we get the same extra condition $c_1=0$, i.e.
$$
c_1=c_2=c_3=c_4=c_5=0.
$$
Hence $\GG = Id$ in all such cases.\\

For $\AAN{1}=\AAN{0}\cup \{ x_2x_3x_5 \}$ we get an extra condition $c_4=2c_1$, i.e.
$$
c_2=c_4=2c_1,\; c_5=c_3+2c_1,\; 4c_1=2c_3=0.
$$

For $\AAN{2}=\AAN{1}\cup \{ x_1x_2x_4 \}$ we get an extra condition $c_1=0$, i.e.
$$
c_1=c_2=c_4=0,\; c_5=c_3,\; 2c_3=0.
$$
All such cases have already appeared.\\

For $\AAN{2}=\AAN{1}\cup \{ x_2x_3x_4 \}$ and $\AAN{2}=\AAN{1}\cup \{ x_2x_4x_5 \}$ we get the same (upto a permutation of $x_i$) extra condition $c_3=0$, i.e.
$$
c_2=c_4=c_5=2c_1,\; c_3=0,\; 4c_1=0.
$$

For $\AAN{3}=\AAN{2}\cup \{ x_1x_4x_5 \}$, $\AAN{3}=\AAN{2}\cup \{ x_1x_2x_4 \}$, $\AAN{3}=\AAN{2}\cup \{ x_1x_2x_5 \}$ (for $\AAN{3}=\AAN{2}\cup \{ x_1x_3x_4 \}$, $\AAN{3}=\AAN{2}\cup \{ x_1x_2x_4 \}$, $\AAN{3}=\AAN{2}\cup \{ x_1x_2x_3 \}$ in case of $\AAN{2}=\AAN{1}\cup \{ x_2x_4x_5 \}$ above) we get the same extra condition $c_1=0$, i.e.
$$
c_1=c_2=c_3=c_4=c_5.
$$
Hence $\GG = Id$ in all such cases.\\

For $\AAN{3}=\AAN{2}\cup \{ x_2x_4x_5 \}$ (for $\AAN{3}=\AAN{2}\cup \{ x_2x_3x_4 \}$ in case of $\AAN{2}=\AAN{1}\cup \{ x_2x_4x_5 \}$ above) we get an extra condition $2c_1=0$, i.e.
$$
c_2=c_3=c_4=c_5=0,\; 2c_1=0.
$$
All such cases have already appeared.\\

Let $\AAN{0} =\{ x_1^2x_2, x_2^2x_0, x_3^2x_0, x_4^2x_0, x_5^2x_0, x_1x_3x_5, x_0^3 \}$. Then 
$$
c_2=0,\; c_5=c_1+c_3,\; 2c_1=2c_3=2c_4=0.
$$
All such cases have already appeared.\\

If $\AAN{0} =\{ x_1^2x_0, x_2^2x_0, x_3^2x_0, x_4^2x_0, x_5^2x_0, x_1x_2x_3, x_3x_4x_5, x_1x_4x_5, x_0^3 \}$, then 
$$
c_2=0,\; c_3=c_1=c_4+c_5,\; 2c_4=2c_5=0.
$$
All such cases have already appeared.\\

Let $\AAN{0} =\{ x_1^2x_2, x_2^2x_0, x_3^2x_2, x_4^2x_2, x_5^2x_2, x_1x_3x_4, x_0^3 \}$. Then 
$$
c_2=0,\; c_4=c_1+c_3,\; 2c_1=2c_3=2c_5=0.
$$
All such cases have already appeared.\\

\subsection{Case of no cubes.}

If a longest cycle has length $6$, then $\AA$ is 

\begin{center}
\begin{tikzpicture}

\coordinate (0) at (0,1); \node[above, font=\scriptsize] at (0) {0}; \fill (0) circle (1pt); 
\coordinate (1) at (1,1); \node[above, font=\scriptsize] at (1) {1}; \fill (1) circle (1pt);
\coordinate (2) at (1.5,.5); \node[right, font=\scriptsize] at (2) {2}; \fill (2) circle (1pt);
\coordinate (3) at (1,0); \node[below, font=\scriptsize] at (3) {3}; \fill (3) circle (1pt);
\coordinate (4) at (0,0); \node[below, font=\scriptsize] at (4) {4}; \fill (4) circle (1pt);
\coordinate (5) at (-.5,.5); \node[left, font=\scriptsize] at (5) {5}; \fill (5) circle (1pt);

\arr{0}{1};  \arr{1}{2}; \arr{2}{3}; \arr{3}{4};  \arr{4}{5};  \arr{5}{0}; 

\end{tikzpicture}
\end{center}

i.e. $\AA = \{ x_1^2x_2, x_2^2x_3, x_3^2x_4, x_4^2x_5, x_5^2x_0, x_0^2x_1 \}$. Then $c_1=2c_1+c_2=2c_2+c_3=2c_3+c_4=2c_4+c_5=2c_5$, i.e.
$$
c_2=-c_1,\; c_3=3c_1,\; c_4=-5c_1,\; c_5=11c_1,\; 21c_1=0.
$$
This means that $\GG\cong \ZZ{21}$ with a generator
$$
f=(1,\om, {\om}^{-1}, {\om}^{3}, {\om}^{-5}, {\om}^{11}),\; \om=\sqr{21.}
$$
Then $\AAbar$ contains a monomial $x_0^{i_0}x_1^{i_1}x_2^{i_2}x_3^{i_3}x_4^{i_4}x_5^{i_5}$ is and only if 
$$
i_1-i_2+3i_3-5i_4+11i_5\equiv 1 \; mod \; 21
$$
($1\; mod \; 21$ is the $f$-weight of $x_0^2x_1\in \AA\subset \AAbar$). This means that $\AAbar=\AA$.\\

\subsubsection{Case of no cubes. Length $5$ longest cycle.}

If a longest cycle has length $5$, then $\AA$ is one of the following:

\begin{center}

\end{center}

Indeed, let vertices $1$, $2$, $3$, $4$, $5$ form a longest cycle. By Lemma 1 (\cite{Liendo}, Lemma 1.3) the remaining vertex $0$ should be connected to this cycle. We may assume that $0$ is connected to $1$. Then there may be one singular pair $(x_0,x_5)$ and two singular triples $(x_0,x_2,x_5)$, $(x_0,x_3,x_5)$.\\
 
After resolving them as usual we obtain the pictures shown above.\\

Now let us do computations.\\

Let $\AAN{0} =\{ x_0^2x_1, x_1^2x_2, x_2^2x_3, x_3^2x_4, x_4^2x_5, x_5^2x_1 \}$. Then 
$$
c_2=-c_1,\; c_3=3c_1,\; c_4=-5c_1,\; c_5=11c_1,\; 22c_1=0.
$$

For $\AA=\AAN{0}\cup \{ x_5^2x_2 \}$ and $\AA=\AAN{0}\cup \{ x_0^2x_5 \}$ we get extra conditions $c_2+2c_5=c_1$ and $c_5=c_1$ respectively. They both are equivalent to $2c_1=0$, i.e.
$$
c_1=c_2=c_3=c_4=c_5,\; 2c_1=0.
$$
This means that $\GG\cong \ZZ{2}$ with a generator
$$
f=(1,{\om},{\om},{\om},{\om},{\om}), \; \om=\sqr{2}.
$$
Then $i_1+i_2+i_3+i_4+i_5 \equiv 1 \; mod \; 2$. This means that $\AAbar=\AA\cup \{ f_3(x_1,x_2,x_3,x_4,x_5), x_0^2\cdot f_1(x_1,x_2,x_3,x_4,x_5) \}$.\\

For $\AA=\AAN{0}\cup \{ x_0x_2x_5 \}$, $\AAN{1}=\AAN{0}\cup \{ x_0x_3x_5 \}$, $\AAN{1}=\AAN{0}\cup \{ x_0x_4x_5 \}$ we get extra conditions $c_2+c_5=c_1$, $c_3+c_5=c_1$, $c_4+c_5=c_1$ respectively. They all imply that $c_1=0$, i.e.
$$
c_1=c_2=c_3=c_4=c_5=0.
$$
Hence $\GG = Id$ in all such cases.\\

For $\AAN{1}=\AAN{0}\cup \{ x_0^2x_4 \}$, $\AAN{1}=\AAN{0}\cup \{ x_5^2x_4 \}$, $\AAN{1}=\AAN{0}\cup \{ x_0^2x_3 \}$, $\AAN{1}=\AAN{0}\cup \{ x_5^2x_3 \}$ we get extra conditions $c_4=c_1$, $c_4=c_1$, $c_3=c_1$, $c_3=c_1$ respectively. They all imply that $2c_1=0$, i.e.
$$
c_1=c_2=c_3=c_4=c_5,\; 2c_1=0.
$$
All such cases have already appeared.\\

\subsubsection{Case of no cubes. Length $4$ longest cycle.}

If a longest cycle has length $4$, then $\AA$ may be either 

\begin{center}

\end{center}

or one of the following:

\begin{center}
%

\begin{flalign*}
\mbox{where}\; & A=(x_1^2x_3,\; x_1^2x_4,\; x_1^2x_5,\; x_1x_3x_5,\; x_1x_4x_5).\\
\end{flalign*}
\end{center}

Indeed, let $x_0,x_1,x_2,x_3$ form a length $4$ cycle. The remaining vertices $4$ and $5$ may be either connected or disconnected.\\

If $4$ and $5$ form a (length $2$) cycle, then there are neither singular pairs nor singular triples.\\ 

If $4$ and $5$ do not form a cycle, but $5$ is connected to $4$, then by Lemma 1 (\cite{Liendo}, Lemma 1.3) the $4$-th vertex should be connected to the cycle formed by $0$, $1$, $2$ and $3$.\\

We may assume that $4$ is connected to $0$. Then we may get one singular pair $(x_3,x_4)$ and one singular triple $(x_1, x_3,x_4)$. We resolve them as usual and obtain the pictures shown above.\\

Suppose that vertices $4$ and $5$ are disconnected. By Lemma 1 (\cite{Liendo}, Lemma 1.3) each of them should be connected to the cycle.\\

If they are connected to a pair of opposite vertices of the cycle (say, $4$ is connected to $0$ and $5$ is connected to $2$), then we may get
\begin{itemize}
\item two singular pairs $(x_3,x_4)$, $(x_1,x_5)$ and
\item four singular triples $(x_1, x_3,x_4)$, $(x_1,x_3, x_5)$, $(x_1, x_4,x_5)$, $(x_3, x_4,x_5)$.
\end{itemize}

We resolve them as usual and obtain the following pictures of $\AA$:

\begin{center}


\begin{multline*}
\mbox{where}\; A=(x_1x_3x_5,\; x_1^2x_3,\; x_5^2x_1,\; x_5^2x_3,\; x_1^2x_0,\; x_5^2x_0,\; x_1x_0x_5,\; x_1^2x_4,\; x_1x_4x_5)\cup (x_1x_3x_5,\\
 x_5^2x_1,\; x_5^2x_3,\; x_1^2x_4,\; x_3^2x_5,\; x_1x_3x_4,\; x_1x_4x_5,\; x_3x_4x_5,\; x_1^2x_3)\cup (x_1x_4x_5,\\
 x_5^2x_1,\; x_5^2x_3,\; x_1^2x_4,\; x_4^2x_1,\; x_1^2x_3,\; x_4^2x_3,\; x_1x_3x_4,\; x_1x_3x_5,\; x_3x_4x_5)\cup (x_1x_4x_3,\\
 x_1^2x_3,\; x_4^2x_3,\; x_1^2x_4,\; x_4^2x_1,\; x_3^2x_5,\; x_5x_3x_4,\; x_1x_3x_5,\; x_1x_4x_5)\cup (x_3x_4x_5,\; x_3^2x_5,\\
 x_5^2x_3,\; x_4^2x_3,\; x_4^2x_1,\; x_5^2x_1,\; x_1x_3x_4,\; x_1x_3x_5,\; x_1x_4x_5).
\end{multline*}
\end{center}

Suppose that vertices $4$ and $5$ are connected to a pair of vertices of the cycle which share a side (say, $4$ is connected to $0$ and $5$ is connected to $1$). Then we may get
\begin{itemize}
\item two singular pairs $(x_3,x_4)$, $(x_0,x_5)$ and
\item three singular triples $(x_1, x_3,x_4)$, $(x_3, x_4,x_5)$, $(x_0, x_2,x_5)$.
\end{itemize}

We resolve them as usual and obtain the following pictures of $\AA$:

\begin{center}


\begin{flalign*}
\mbox{where}\; & A=(x_5^2x_0,\; x_0x_4x_5,\; x_0^2x_4,\; x_0^2x_3,\; x_5^2x_3,\; x_0x_3x_5,\; x_0^2x_2,\; x_0x_2x_5)\cup (x_0x_2x_5,\\
& x_0^2x_2,\; x_2^2x_0,\; x_5^2x_0,\; x_2^2x_5,\; x_0^2x_4,\; x_2^2x_4,\; x_0x_2x_4,\; x_2x_4x_5,\; x_0x_4x_5);\\
& B=(x_3x_4x_5,\; x_3^2x_5,\; x_5^2x_3,\; x_2x_4x_5,\; x_2x_3x_4,\; x_2x_3x_5,\; x_3^2x_2,\; x_4^2x_2);\\
& C= (x_3x_4x_5,\; x_3^2x_5,\; x_3^2x_2,\; x_2x_4x_5,\; x_2x_3x_4,\; x_2x_3x_5);\\
& D=(x_1x_3x_4,\; x_3^2x_5,\; x_3^2x_1,\; x_1^2x_5,\; x_4^2x_1,\; x_1^2x_4,\; x_1^2x_3,\; x_1x_4x_5,\; x_3x_4x_5,\; x_1x_3x_5);\\
& E=(x_1x_3x_4,\;x_1^2x_4,\; x_1^2x_3,\; x_1x_4x_5,\; x_3x_4x_5,\; x_1x_3x_5,\; x_3^2x_5,\; x_1^2x_5).\\
\end{flalign*}
\end{center}

Finally, suppose that vertices $4$ and $5$ are connected to the same vertex of the cycle (say, to $0$). Then we may get
\begin{itemize}
\item three singular pairs $(x_3,x_4)$, $(x_3,x_5)$, $(x_4,x_5)$ and
\item five singular triples $(x_1, x_3,x_4)$, $(x_1, x_3,x_5)$, $(x_1, x_4,x_5)$, $(x_2,x_4, x_5)$, $(x_3, x_4,x_5)$.
\end{itemize}

We resolve them as usual and obtain the following pictures of $\AA$:

\begin{center}


\begin{flalign*}
\mbox{where}\; & A=(x_3^2x_1,\; x_3^2x_2,\; x_1x_3x_4,\; x_2x_3x_4)\cup (x_3^2x_1,\; x_3^2x_2,\; x_1x_3x_5,\; x_2x_3x_5);\\
& B=(x_3^2x_1,\; x_1x_4x_5,\; x_1x_3x_4,\; x_1x_3x_5)\cup (x_1x_4x_5,\; x_1^2x_4,\; x_1^2x_5,\; x_1^2x_3,\; x_1x_3x_4,\; x_1x_3x_5).
\end{flalign*}
\end{center}

Now let us do computations.\\

If $\AA =\{ x_0^2x_1, x_1^2x_2, x_2^2x_3, x_3^2x_0, x_4^2x_5, x_5^2x_4 \}$, then 
$$
c_1=3c_4,\; c_2=-3c_4,\; c_3=9c_4,\; c_5=c_4,\; 15c_4=0.
$$
This means that $\GG\cong \ZZ{15}$ with a generator
$$
f=(1,{\om}^{3},{\om}^{-3},{\om}^{9},{\om},{\om}), \; \om=\sqr{15}.
$$
Then $3i_1-3i_2+9i_3+i_4+i_5 \equiv 3 \; mod \; 15$. This means that $\AAbar=\AA\cup \{ f_3(x_4,x_5) \}$.\\

Let $\AAN{0} =\{ x_0^2x_1, x_1^2x_2, x_2^2x_3, x_3^2x_0, x_4^2x_0, x_5^2x_4 \}$. Then 
$$
c_1=4c_5,\; c_2=-4c_5,\; c_3=12c_5,\; c_4=2c_5,\; 20c_5=0.
$$

For $\AA=\AAN{0}\cup \{ x_3^2x_1 \}$, $\AAN{1}=\AAN{0}\cup \{ x_3^2x_2 \}$, $\AA=\AAN{0}\cup \{ x_4^2x_1 \}$, $\AAN{1}=\AAN{0}\cup \{ x_4^2x_2 \}$ we get the same extra condition $c_1=0$, i.e.
$$
c_1=c_2=c_3=0,\; c_4=2c_5,\; 4c_5=0.
$$
All such cases have already appeared.\\

For $\AA=\AAN{0}\cup \{ x_1x_3x_4 \}$, $\AAN{1}=\AAN{0}\cup \{ x_2x_3x_4 \}$, $\AA=\AAN{0}\cup \{ x_3x_4x_5 \}$ we get extra conditions $2c_5=0$, $2c_5=0$, $c_5=0$ respectively. They all imply that $2c_5=0$, i.e.
$$
c_1=c_2=c_3=c_4=0,\; 2c_5=0.
$$
All such cases have already appeared.\\

Let $\AAN{0} =\{ x_0^2x_1, x_1^2x_2, x_2^2x_3, x_3^2x_0, x_4^2x_0, x_5^2x_2 \}$. Then 
$$
c_1=2c_4,\; c_2=-2c_4,\; c_3=6c_4,\; 10c_4=2(c_5-2c_4)=0.
$$

For $\AAN{1}=\AAN{0}\cup \{ x_3^2x_1 \}$, $\AAN{1}=\AAN{0}\cup \{ x_4^2x_1 \}$, $\AAN{1}=\AAN{0}\cup \{ x_4^2x_3 \}$, $\AAN{1}=\AAN{0}\cup \{ x_1x_3x_4 \}$, $\AAN{1}=\AAN{0}\cup \{ x_3^2x_5 \}$, $\AAN{1}=\AAN{0}\cup \{ x_3x_4x_5 \}$  we get extra conditions $2c_4=0$, $2c_4=0$, $2c_4=0$, $c_4=0$, $c_5=0$, $c_5=5c_4$ respectively. They all imply that $2c_4=0$, i.e.
$$
c_1=c_2=c_3=0,\; 2c_4=2c_5=0.
$$
All such cases have already appeared.\\

For $\AAN{1}=\AAN{0}\cup \{ x_4^2x_2 \}$, $\AAN{1}=\AAN{0}\cup \{ x_3^2x_2 \}$, $\AAN{1}=\AAN{0}\cup \{ x_2x_3x_4 \}$ we get extra conditions $2c_4=0$, $2c_4=0$, $c_4=0$ respectively. They all imply that $2c_4=0$, i.e.
$$
c_1=c_2=c_3=0,\; 2c_4=2c_5=0.
$$
All such cases have already appeared.\\

Let $\AAN{0} =\{ x_0^2x_1, x_1^2x_2, x_2^2x_3, x_3^2x_0, x_4^2x_0, x_5^2x_1 \}$. Then 
$$
c_1=2c_4,\; c_2=-2c_4,\; c_3=6c_4,\; 10c_4=2c_5=0.
$$

For $\AAN{1}=\AAN{0}\cup \{ x_3^2x_1 \}$, $\AAN{1}=\AAN{0}\cup \{ x_3^2x_2 \}$, $\AAN{1}=\AAN{0}\cup \{ x_4^2x_1 \}$, $\AAN{1}=\AAN{0}\cup \{ x_4^2x_2 \}$, $\AAN{1}=\AAN{0}\cup \{ x_1x_3x_4 \}$, $\AAN{1}=\AAN{0}\cup \{ x_2x_3x_4 \}$  we get extra conditions $2c_4=0$, $2c_4=0$, $2c_4=0$, $2c_4=0$, $c_4=0$, $c_4=0$ respectively. They all imply that $2c_4=0$, i.e.
$$
c_1=c_2=c_3=0,\; 2c_4=2c_5=0.
$$
All such cases have already appeared.\\

For $\AAN{1}=\AAN{0}\cup \{ x_3x_4x_5 \}$ we get an extra condition $c_5=5c_4$, i.e.
$$
c_1=2c_4,\; c_2=-2c_4,\; c_3=6c_4,\; c_5=5c_4,\; 10c_4=0.
$$

For $\AA=\AAN{1}\cup \{ x_5^2x_0 \}$, $\AAN{2}=\AAN{1}\cup \{ x_0^2x_4 \}$, $\AAN{2}=\AAN{1}\cup \{ x_0x_4x_5 \}$, $\AAN{2}=\AAN{1}\cup \{ x_0^2x_3 \}$, $\AAN{2}=\AAN{1}\cup \{ x_5^2x_3 \}$, $\AAN{2}=\AAN{1}\cup \{ x_0x_3x_5 \}$, $\AAN{2}=\AAN{1}\cup \{ x_0x_2x_5 \}$, $\AAN{2}=\AAN{1}\cup \{ x_0^2x_2 \}$ we get extra conditions $2c_4=0$, $c_4=0$, $2c_4=0$, $2c_4=0$, $2c_4=0$, $c_4=0$, $c_4=0$, $2c_4=0$ respectively. They all imply that $2c_4=0$, i.e.
$$
c_1=c_2=c_3=0,\; c_5=c_4,\; 2c_4=0.
$$
All such cases have already appeared.\\

Let $\AAN{0} =\{ x_0^2x_1, x_1^2x_2, x_2^2x_3, x_3^2x_0, x_4^2x_0, x_5^2x_0 \}$. Then 
$$
c_1=2c_4,\; c_2=-2c_4,\; c_3=6c_4,\; 10c_4=2(c_5-c_4)=0.
$$

For $\AAN{1}=\AAN{0}\cup \{ x_3x_4x_5 \}$, $\AAN{1}=\AAN{0}\cup \{ x_2x_4x_5 \}$, $\AAN{1}=\AAN{0}\cup \{ x_1x_4x_5 \}$ we get extra conditions $c_5=5c_4$, $c_5=3c_4$, $c_5=-c_4$ respectively. They all imply that $2c_4=0$, i.e.
$$
c_1=c_2=c_3=0,\; 2c_4=2c_5=0.
$$
All such cases have already appeared.\\

\subsubsection{Case of no cubes. Length $3$ longest cycle.}

If a longest cycle has length $3$, then $\AA$ may be either

\begin{center}

\end{center}

or one of the following:

\begin{center}
%

\begin{flalign*}
\mbox{where}\; & A=(x_2x_3x_5,\; x_5^2x_2,\; x_2^2x_5,\; x_2^2x_1,\; x_3^2x_1,\; x_5^2x_1,\; x_1x_3x_5,\; x_2x_1x_3,\; x_2x_1x_5).
\end{flalign*}
\end{center}

Indeed, let $x_0,x_1,x_2$ form a length $3$ cycle.\\

Consider the subgraph formed by vertices $x_3,x_4,x_5$. This subgraph may also have a length $3$ cycle. If it does, then we may assume that $3$ is connected to $4$, $4$ is connected to $5$ and $5$ is connected to $3$. Then there are neither singular pairs nor singular triples.\\

Suppose the subgraph formed by vertices $3$, $4$, $5$ does not have length $3$ cycles. Let us consider the length of a longest path of this subgraph. If it is equal to $3$, then we may assume that $3$ is connected to $4$ and $4$ is connected to $5$.\\

By Lemma 1 (\cite{Liendo}, Lemma 1.3) the $5$-th vertex should be connected either to $4$ or to the cycle (say, to $0$). It can not be connected to $3$, because the subgraph formed by $3$, $4$, $5$ has no length $3$ cycles by the assumption. \\

If $5$ is connected to $4$, then we may have one singular pair $(x_5,x_3)$ and three singular triples: $(x_0,x_3,x_5)$, $(x_1,x_3,x_5)$ and $(x_2,x_3,x_5)$.\\

In order to resolve the singular pair, we should have that either
\begin{itemize}
\item $3$ and $5$ are connected by an edge, or
\item $3$ or $5$ is connected to the length $3$ cycle (say, to $0$), or
\item there is a dashed curve passing through $3$, $5$ and the length $3$ cycle (say, through the $2$-nd vertex).
\end{itemize}

In the first case, $3$ should be connected to $5$ (since otherwise the subgraph formed by $3$, $4$, $5$ would have a length $3$ cycle). Then the singular triples are automatically resolved as well.\\

In the other two cases, one resolves the singular triples as usual.\\

This gives us the pictures shown above.\\

Alternatively, the $5$-th vertex may be connected to the length $3$ cycle (say, to $0$). In this case we may get
\begin{itemize}
\item one singular pair $(x_2,x_5)$ and
\item one singular triple $(x_2,x_3,x_5)$.
\end{itemize}

We resolve them as usual and obtain the following pictures of $\AA$:

\begin{center}

\end{center}

Suppose that the length of a longest path of the subgraph formed by vertices $3$, $4$, $5$ is equal to $2$. We may assume that $4$ is connected to $3$. Then by Lemma 1 (\cite{Liendo}, Lemma 1.3) either $3$ is connected to $4$ (i.e. $3$ and $4$ form a length $2$ cycle) or $3$ is connected to the length $3$ cycle (say, to $0$).\\

If $3$ and $4$ form a cycle, then  by Lemma 1 (\cite{Liendo}, Lemma 1.3) the $5$-th vertex should be connected to either of the cycles. It can not be connected to the length $2$ cycle (since otherwise the subgraph formed by vertices $3$, $4$, $5$ would have a length $3$ path). Hence it should be connected to the length $3$ cycle (say, to $0$).\\ 

In this case we may get
\begin{itemize}
\item one singular pair $(x_2,x_5)$ and
\item two singular triples $(x_2,x_3,x_5)$, $(x_2,x_4,x_5)$.
\end{itemize}

We resolve them as usual and obtain the following pictures of $\AA$:

\begin{center}

\end{center}

If $3$ and $4$ do not form a cycle, while  $3$ is connected to $0$, then by Lemma 1 (\cite{Liendo}, Lemma 1.3) the $5$-th vertex should be connected either to $3$ or to the cycle.\\

If the $5$-th vertex is connected to $1$, then we may get 
\begin{itemize}
\item two singular pairs $(x_0,x_5)$, $(x_2,x_3)$ and
\item two singular triples $(x_0,x_4,x_5)$, $(x_2,x_3,x_5)$.
\end{itemize}

If the $5$-th vertex is connected to $2$, then we may get 
\begin{itemize}
\item two singular pairs $(x_1,x_5)$, $(x_2,x_3)$ and
\item two singular triples $(x_1,x_3,x_5)$, $(x_1,x_4,x_5)$.
\end{itemize}

We resolve all these singular pairs and triples as usual and obtain the following pictures of $\AA$:

\begin{center}


\begin{flalign*}
\mbox{where}\; & A=(x_1x_3x_5,\; x_1^2x_3,\; x_5^2x_3,\; x_1x_3x_4,\; x_1x_4x_5,\; x_3x_4x_5)\cup (x_1x_4x_5,\; x_4^2x_1,\; x_4^2x_5,\\
& x_4^2x_0,\; x_1^2x_0,\; x_5^2x_0,\; x_1x_0x_4,\; x_1x_0x_5,\; x_0x_4x_5);\\
& B=(x_5^2x_0,\; x_0^2x_2,\; x_5^2x_2,\; x_0x_2x_5,\; x_0^2x_3,\; x_5^2x_3,\; x_0x_3x_5,\; x_0x_4x_5)\cup (x_0x_4x_5,\\
& x_4^2x_5,\; x_4^2x_0,\; x_5^2x_0,\; x_0^2x_2,\; x_5^2x_2,\; x_4^2x_2,\; x_2x_0x_4,\; x_2x_0x_5,\; x_2x_4x_5);\\
& C=(x_2x_3x_5,\; x_5^2x_3,\; x_3^2x_2,\; x_2x_4x_5,\; x_2x_3x_4,\; x_3x_4x_5);\\
& F=(x_5^2x_0,\; x_5^2x_3,\; x_0x_3x_5,\; x_0x_4x_5,\; x_0^2x_2,\; x_0x_2x_5)\cup (x_0x_4x_5,\; x_4^2x_0,\; x_5^2x_0,\; x_4^2x_5,\\
& x_0^2x_2,\; x_4^2x_2,\; x_2x_0x_4,\; x_2x_0x_5,\; x_2x_4x_5);\\
& G=(x_5^2x_3,\; x_2x_3x_5,\; x_2x_4x_5,\; x_2^2x_1,\; x_1x_2x_5)\cup (x_2x_4x_5,\; x_4^2x_5,\; x_4^2x_2,\; x_2^2x_1,\; x_5^2x_1,\\
& x_4^2x_1,\; x_2x_1x_4,\; x_2x_1x_5,\; x_1x_4x_5);
\end{flalign*}
\end{center}

\begin{center}


\begin{flalign*}
\mbox{where}\; & A=(x_1x_3x_5,\; x_1^2x_3,\; x_5^2x_3,\; x_1x_3x_4,\; x_1x_4x_5,\; x_3x_4x_5)\cup (x_1x_4x_5,\; x_4^2x_1,\; x_4^2x_5,\\
& x_4^2x_0,\; x_1^2x_0,\; x_5^2x_0,\; x_1x_0x_4,\; x_1x_0x_5,\; x_0x_4x_5);\\
& D=(x_2x_3x_5,\; x_5^2x_3,\; x_5^2x_2,\; x_2x_4x_5,\; x_3x_4x_5);\\
& E=(x_5^2x_0,\; x_0^2x_2,\; x_5^2x_2,\; x_0x_2x_5,\; x_0^2x_3,\; x_5^2x_3,\; x_0x_3x_5,\; x_0x_4x_5)\cup (x_4^2x_5,\; x_5^2x_0,\\
& x_4^2x_0,\; x_0^2x_2,\; x_5^2x_2,\; x_4^2x_2,\; x_0x_4x_2,\; x_0x_2x_5,\; x_0x_4x_5,\; x_2x_4x_5);\\
& F=(x_1x_3x_5,\; x_1^2x_0,\; x_5^2x_0,\; x_1^2x_3,\; x_5^2x_3,\; x_0x_1x_5,\; x_1x_4x_5).\\
\end{flalign*}
\end{center}

If the $5$-th vertex is connected to $0$, then we may get 
\begin{itemize}
\item three singular pairs $(x_2,x_3)$, $(x_2,x_5)$, $(x_3,x_5)$ and
\item three singular triples $(x_1,x_3,x_5)$, $(x_2,x_3,x_5)$, $(x_2,x_4,x_5)$.
\end{itemize}

If the $5$-th vertex is connected to $3$, then we may get 
\begin{itemize}
\item two singular pairs $(x_2,x_3)$, $(x_4,x_5)$ and
\item three singular triples $(x_0,x_4,x_5)$, $(x_1,x_4,x_5)$, $(x_2,x_4,x_5)$.
\end{itemize}

We resolve all these singular pairs and triples as usual and obtain the following pictures of $\AA$:

\begin{center}

\end{center}

\begin{center}
\begin{flalign*}
\mbox{where}\; & A=(x_1x_3x_5,\; x_2x_3x_5,\; x_3x_4x_5)\cup (x_1x_2x_5,\; x_2^2x_1,\; x_2x_4x_5,\; x_2x_3x_5);\\
& B=(x_5^2x_0,\; x_4^2x_0,\; x_0x_4x_5,\; x_1x_4x_5,\; x_2x_4x_5)\cup (x_0x_4x_5,\; x_5^2x_0,\; x_4^2x_0,\; x_0^2x_2,\; x_0x_4x_2,\\
& x_0x_2x_5,\; x_2x_4x_5)\cup (x_1x_4x_5,\; x_1^2x_0,\; x_4^2x_0,\; x_5^2x_0,\; x_0x_4x_5,\; x_0x_1x_4,\; x_0x_1x_5);\\
& C=(x_1x_4x_5,\; x_1x_2x_4,\; x_1x_2x_5,\; x_2x_4x_5);\;\;\;\;\; D=(x_2x_4x_5,\; x_3x_4x_5,\; x_2x_3x_5);\\
& E=(x_1x_3x_5,\; x_1^2x_5,\; x_1^2x_3,\; x_1x_3x_4,\; x_1x_4x_5,\; x_3x_4x_5);\\
& F=(x_2x_4x_5,\; x_4^2x_5,\; x_4^2x_2,\; x_4^2x_1,\; x_1x_2x_4,\; x_1x_4x_5,\; x_1x_2x_5).
\end{flalign*}
\end{center}

Finally, suppose that the subgraph formed by vertices $3$, $4$, $5$ is totally disconnected. By Lemma 1 (\cite{Liendo}, Lemma 1.3) each of its vertices should be connected to the length $3$ cycle.\\

Let us assume that the three vertices $3$, $4$, $5$ are connected to three different vertices of the cycle (say, $3$ is connected to $0$, $4$ is connected to $1$ and $5$ is connected to $2$).\\

Then we may get 
\begin{itemize}
\item three singular pairs $(x_0,x_4)$, $(x_1,x_5)$, $(x_2,x_3)$ and
\item three singular triples $(x_0,x_4, x_5)$, $(x_1,x_3,x_5)$, $(x_2,x_3, x_4)$.
\end{itemize}

We resolve them as usual and obtain the following pictures of $\AA$:

\begin{center}


\begin{flalign*}
\mbox{where}\; & A=(x_3^2x_2,\; x_3^2x_1,\; x_2^2x_1,\; x_1x_2x_3,\; x_2^2x_5,\; x_2x_3x_4,\; x_2x_3x_5)\cup (x_1x_4x_5,\; x_5^2x_1,\\
& x_1^2x_4,\; x_1^2x_0,\; x_5^2x_0,\; x_1x_0x_5,\; x_1x_3x_5)\cup (x_2x_3x_4,\; x_3^2x_2,\; x_4^2x_2,\; x_2^2x_5,\\
& x_2x_3x_5,\; x_2x_4x_5,\; x_3x_4x_5)\cup (x_1x_3x_5,\; x_3^2x_1,\; x_5^2x_1,\; x_1^2x_4,\; x_1x_4x_5,\; x_1x_3x_4,\; x_3x_4x_5);\\
& B=(x_0x_4x_5,\; x_0^2x_3,\; x_5^2x_0,\; x_4^2x_0,\; x_0x_3x_5,\; x_0x_3x_4,\; x_3x_4x_5).\\
\end{flalign*}
\end{center}

Suppose that vertices $3$, $4$, $5$ are connected to only two distinct vertices of the cycle (say, $3$ and $4$ are connected to $0$, while $5$ is connected to $1$ or to $2$).\\

If $5$ is connected to $1$, then we may get 
\begin{itemize}
\item four singular pairs $(x_2,x_3)$, $(x_2,x_4)$, $(x_3,x_4)$, $(x_0,x_5)$ and
\item five singular triples $(x_3,x_4, x_5)$, $(x_1,x_3,x_4)$, $(x_2,x_3, x_4)$, $(x_2,x_3, x_5)$, $(x_2,x_4, x_5)$.
\end{itemize}

If $5$ is connected to $2$, then we may get 
\begin{itemize}
\item four singular pairs $(x_2,x_3)$, $(x_2,x_4)$, $(x_3,x_4)$, $(x_1,x_5)$ and
\item five singular triples $(x_3,x_4, x_5)$, $(x_1,x_3,x_4)$, $(x_2,x_3, x_4)$, $(x_1,x_3, x_5)$, $(x_1,x_4, x_5)$.
\end{itemize}

We resolve these singular pairs and triples as usual and obtain the following pictures of $\AA$:

\begin{center}


\begin{flalign*}
\mbox{where}\; & A_1=(x_2^2x_1,\; x_4^2x_1,\; x_2x_4x_5,\; x_1x_2x_4); \;\;\; A_2=(x_2^2x_1,\; x_2^2x_5,\; x_2x_4x_5,\; x_1x_2x_4);\\
& B_1=(x_2^2x_1,\; x_3^2x_1,\; x_2x_3x_5,\; x_1x_2x_3); \;\;\; B_2=(x_2^2x_1,\; x_2^2x_5,\; x_2x_3x_5,\; x_1x_2x_3);\\
& C_1=(x_2x_3x_5,\; x_2x_4x_5,\; x_3x_4x_5); \;\;\;\;\; C_2=(x_2^2x_5,\; x_2x_3x_5,\; x_2x_4x_5,\; x_3x_4x_5);\\
& D=(x_4^2x_2,\; x_3^2x_2,\; x_4^2x_1,\; x_3^2x_1,\; x_2^2x_1,\; x_1x_3x_4,\; x_1x_2x_3,\;x_1x_2x_4);\\
& E_1=(x_5^2x_0,\; x_5^2x_2,\; x_0^2x_2,\; x_0^2x_3,\; x_0^2x_4,\; x_0x_3x_5,\; x_0x_2x_5,\; x_0x_4x_5);\\
& E_2=(x_1^2x_0,\; x_5^2x_0,\; x_1x_3x_5,\; x_1x_4x_5,\; x_0x_1x_5).\\
\end{flalign*}
\end{center}

The remaining possibility is that vertices $3$, $4$, $5$ are connected to one and the same vertex of the length $3$ cycle (say, to $0$).\\

Then we may get 
\begin{itemize}
\item six singular pairs $(x_2,x_3)$, $(x_2,x_4)$, $(x_2,x_5)$, $(x_3,x_4)$, $(x_3,x_5)$, $(x_4,x_5)$ and
\item seven singular triples $(x_3,x_4, x_5)$, $(x_2,x_3,x_4)$, $(x_2,x_3, x_5)$, $(x_2,x_4, x_5)$, $(x_1,x_3, x_4)$,\\
$(x_1,x_3, x_5)$, $(x_1,x_4, x_5)$.
\end{itemize}

After resolving them as in the earlier examples, we obtain the following pictures of $\AA$\footnote{Note that in the last picture we take a usual {\it union} of the set of monomials represented by the graph and a set of three other monomials.}:

\begin{center}

\end{center}

Now let us do computations.\\

If $\AA =\{ x_0^2x_1, x_1^2x_2, x_2^2x_0, x_3^2x_4, x_4^2x_5, x_5^2x_3 \}$. Then 
$$
c_2=-c_1,\; c_4=c_1-2c_3,\; c_5=4c_3-c_1,\; 3c_1=9c_3=0.
$$
This means that $\GG\cong \ZZ{3}\oplus \ZZ{9}$ with generators
$$
f_1=(1,{\om},{\om}^{-1},1,{\om},{\om}^{-1}), \; f_2=(1,1,1,{\eta},{\eta}^{-2},{\eta}^4), \; \om=\sqr{3}, \eta=\sqr{9}.
$$
Then $i_1+i_4 \equiv 1+i_2+i_5 \; mod \; 3$ and $i_3+4i_5 \equiv 2i_4 \; mod \; 9$. This means that $\AAbar=\AA$.\\

At this point the following simplifying remark can be made:
\begin{center}
{\it Whenever $\AAbar$ contains a cube, this case has already appeared earlier.}
\end{center}

Still, for the sake of completeness we will continue our analysis in the same fashion as above.\\

If $\AA =\{ x_0^2x_1, x_1^2x_2, x_2^2x_0, x_3^2x_4, x_4^2x_5, x_5^2x_4, x_3^2x_5 \}$. Then 
$$
c_1=3c_4,\; c_2=-3c_4,\; c_5=c_4,\; 2(c_3-c_4)=9c_4=0.
$$
This means that $\GG\cong \ZZ{2}\oplus \ZZ{9}$ with generators
$$
f_1=(1,{\om}^3,{\om}^{-3},\om,{\om},{\om}), \; f_2=(1,1,1,{\eta},1,1), \; \om=\sqr{9}, \eta=\sqr{2}.
$$
Then $i_3 \equiv 0 \; mod \; 2$ and $3i_1+i_3+i_4+i_5 \equiv 3i_2+3 \; mod \; 9$. This means that $\AAbar=\AA\cup \{ f_3(x_4,x_5) \}$.\\

Let $\AAN{0} =\{ x_0^2x_1, x_1^2x_2, x_2^2x_0, x_3^2x_4, x_4^2x_5, x_5^2x_4, x_2x_3x_5 \}$. Then 
$$
c_1=c_2=c_3=c_4=c_5=0.
$$
Hence $\GG=Id$ in all such cases.\\

Let $\AAN{0} =\{ x_0^2x_1, x_1^2x_2, x_2^2x_0, x_3^2x_0, x_3^2x_4, x_4^2x_5, x_5^2x_4 \}$ or $\AAN{0} =\{ x_0^2x_1,$ $x_1^2x_2,$ $x_2^2x_0,$ $x_3^2x_4,$ $x_4^2x_5,$ $x_5^2x_4,$ $x_5^2x_0 \}$. Then 
$$
c_1=c_2=c_4=c_5=0,\; 2c_3=0.
$$
All such cases have already appeared.\\

Let $\AAN{0} =\{ x_0^2x_1, x_1^2x_2, x_2^2x_0, x_3^2x_4, x_4^2x_5, x_5^2x_0 \}$. Then 
$$
c_1=2c_5,\; c_2=-2c_5,\; c_4=2c_5-2c_3,\; 4c_3+3c_5=6c_5=0.
$$

For $\AA=\AAN{0}\cup \{ x_5^2x_1 \}$, $\AA=\AAN{0}\cup \{ x_2^2x_1 \}$, $\AA=\AAN{0}\cup \{ x_1x_2x_5 \}$, $\AAN{1}=\AAN{0}\cup \{ x_2x_4x_5 \}$ we get extra conditions $2c_5=0$, $2c_5=0$, $c_5=0$, $c_5=-2c_3$ respectively. They all imply that $2c_5=0$, i.e.
$$
c_1=c_2=0,\; c_4=-2c_3,\; c_5=4c_3,\; 8c_3=0.
$$
All such cases have already appeared.\\

For $\AA=\AAN{0}\cup \{ x_2x_3x_5 \}$ we get an extra condition $c_3=3c_5$, i.e.
$$
c_1=c_4=-c_5,\; c_2=c_5,\; c_3=0,\; 3c_5=0.
$$
This means that $\GG\cong \ZZ{3}$ with a generator
$$
f=(1,{\om}^{-1},{\om},1,{\om}^{-1},{\om}), \; \om=\sqr{3}.
$$
Then $i_2+i_5 \equiv 2+i_1+i_4 \; mod \; 3$. This means that $\AAbar=\AA\cup \{ f_2(x_2,x_5)\cdot f_1(x_0,x_3), f_2(x_0,x_3)\cdot f_1(x_1,x_4), f_2(x_1,x_4)\cdot f_1(x_2,x_5) \}$.\\

Let $\AAN{0} =\{ x_0^2x_1, x_1^2x_2, x_2^2x_0, x_3^2x_4, x_4^2x_3, x_5^2x_0 \}$. Then 
$$
c_1=2c_5,\; c_2=-2c_5,\; c_3=c_4,\; 3c_3-2c_5=6c_5=0.
$$

For $\AA=\AAN{0}\cup \{ x_1x_2x_5 \}$, $\AA=\AAN{0}\cup \{ x_2^2x_1 \}$, $\AA=\AAN{0}\cup \{ x_5^2x_1 \}$, $\AAN{1}=\AAN{0}\cup \{ x_2^2x_3 \}$, $\AAN{1}=\AAN{0}\cup \{ x_2x_3x_5 \}$  we get extra conditions $c_5=0$, $2c_5=0$, $2c_5=0$, $c_3=0$, $c_3=3c_5$. They all imply that $2c_5=0$, i.e.
$$
c_1=c_2=0,\; c_3=c_4,\; 3c_3=2c_5=0.
$$
All such cases have already appeared.\\

Let $\AAN{0} =\{ x_0^2x_1, x_1^2x_2, x_2^2x_0, x_3^2x_0, x_4^2x_3, x_5^2x_1 \}$. Then 
$$
c_1=4c_4,\; c_2=-4c_4,\; c_3=2c_4,\; 12c_4=2c_5=0.
$$

For $\AAN{1}=\AAN{0}\cup \{ x_2x_3x_5 \}$ we get an extra condition $c_5=6c_4$, i.e.
$$
c_1=4c_4,\; c_2=-4c_4,\; c_3=2c_4,\; c_5=6c_4,\; 12c_4=0.
$$

For $\AA=\AAN{1}\cup \{ x_5^2x_0 \}$, $\AAN{2}=\AAN{1}\cup \{ x_0^2x_2 \}$, $\AAN{2}=\AAN{1}\cup \{ x_5^2x_2 \}$, $\AAN{2}=\AAN{1}\cup \{ x_5^2x_3 \}$, $\AAN{2}=\AAN{1}\cup \{ x_0^2x_3 \}$, $\AAN{2}=\AAN{1}\cup \{ x_0x_2x_5 \}$, $\AAN{2}=\AAN{1}\cup \{ x_0x_3x_5 \}$ we get extra conditions $4c_4=0$, $4c_4=0$, $4c_4=0$, $2c_4=0$, $2c_4=0$, $2c_4=0$, $4c_4=0$. They all imply that $4c_4=0$, i.e.
$$
c_1=c_2=0,\; c_3=c_5=2c_4,\; 4c_4=0.
$$
All such cases have already appeared.\\

For $\AA=\AAN{1}\cup \{ x_0x_4x_5 \}$ we get an extra condition $3c_4=0$, i.e.
$$
c_1=c_4,\; c_2=c_3=-c_4,\; c_5=0,\; 3c_4=0.
$$
This case has already appeared.\\

For $\AAN{1}=\AAN{0}\cup \{ x_1x_2x_3 \}$, $\AAN{1}=\AAN{0}\cup \{ x_2x_3x_4 \}$, $\AAN{1}=\AAN{0}\cup \{ x_2^2x_1 \}$, $\AAN{1}=\AAN{0}\cup \{ x_3^2x_1 \}$ we get extra conditions $2c_4=0$, $c_4=0$, $4c_4=0$, $4c_4=0$ respectively. They all imply that $4c_4=0$, i.e.
$$
c_1=c_2=0,\; c_3=2c_4,\; 4c_4=2c_5=0.
$$
All such cases have already appeared.\\

Let $\AAN{0} =\{ x_0^2x_1, x_1^2x_2, x_2^2x_0, x_3^2x_0, x_3^2x_1, x_4^2x_3, x_5^2x_2 \}$ or $\AAN{0} =\{ x_0^2x_1,$ $x_1^2x_2,$ $x_2^2x_0,$ $x_3^2x_0,$ $x_3^2x_1,$ $x_4^2x_3,$ $x_5^2x_0 \}$. Then 
$$
c_1=c_2=0,\; c_3=2c_4,\; 4c_4=2c_5=0.
$$
All such cases have already appeared.\\

Let $\AAN{0} =\{ x_0^2x_1, x_1^2x_2, x_2^2x_0, x_3^2x_0, x_4^2x_3, x_5^2x_2 \}$. Then 
$$
c_1=4c_4,\; c_2=-4c_4,\; c_3=2c_4,\; 2(c_5-4c_4)=12c_4=0.
$$

For $\AAN{1}=\AAN{0}\cup \{ x_2^2x_1 \}$, $\AAN{1}=\AAN{0}\cup \{ x_2x_3x_4 \}$, $\AAN{1}=\AAN{0}\cup \{ x_1x_2x_3 \}$, $\AAN{1}=\AAN{0}\cup \{ x_2^2x_5 \}$, $\AAN{1}=\AAN{0}\cup \{ x_2x_3x_5 \}$ we get extra conditions $4c_4=0$, $c_4=0$, $2c_4=0$, $c_5=0$, $c_5=6c_4$ respectively. They all imply that $4c_4=0$, i.e.
$$
c_1=c_2=0,\; c_3=2c_4,\; 4c_4=2c_5=0.
$$
All such cases have already appeared.\\

Let $\AAN{0} =\{ x_0^2x_1, x_1^2x_2, x_2^2x_0, x_3^2x_0, x_4^2x_3, x_5^2x_3 \}$. Then 
$$
c_1=4c_4,\; c_2=-4c_4,\; c_3=2c_4,\; 12c_4=2(c_5-c_4)=0.
$$

For $\AAN{1}=\AAN{0}\cup \{ x_3^2x_1 \}$ and $\AAN{1}=\AAN{0}\cup \{ x_2^2x_1 \}$ we get the same extra condition $4c_4=0$, i.e.
$$
c_1=c_2=0,\; c_3=2c_4,\; 4c_4=2(c_5-c_4)=0.
$$

For $\AAN{2}=\AAN{1}\cup \{ x_1x_4x_5 \}$, $\AAN{2}=\AAN{1}\cup \{ x_0x_4x_5 \}$, $\AAN{2}=\AAN{1}\cup \{ x_2x_4x_5 \}$ we get the same extra condition $c_5=-c_4$, i.e.
$$
c_1=c_2=0,\; c_3=2c_4,\; c_5=-c_4,\; 4c_4=0.
$$
All such cases have already appeared.\\

For $\AAN{2}=\AAN{1}\cup \{ x_5^2x_0 \}$, $\AAN{2}=\AAN{1}\cup \{ x_4^2x_0 \}$ we get the same extra condition $2c_4=0$, i.e.
$$
c_1=c_2=c_3=0,\; 2c_4=2c_5=0.
$$
All such cases have already appeared.\\

For $\AAN{1}=\AAN{0}\cup \{ x_2x_3x_4 \}$ and $\AAN{1}=\AAN{0}\cup \{ x_1x_2x_3 \}$ we get extra conditions $c_4=0$ and $2c_4=0$ respectively. They both imply that $2c_4=0$, i.e.
$$
c_1=c_2=c_3=0,\; 2c_4=2c_5=0.
$$
All such cases have already appeared.\\

Let $\AAN{0} =\{ x_0^2x_1, x_1^2x_2, x_2^2x_0, x_3^2x_0, x_4^2x_3, x_5^2x_0 \}$. Then 
$$
c_1=2c_5,\; c_2=-2c_5,\; c_3=2c_5-2c_4,\; 2(c_5-2c_4)=6c_5=0.
$$

For $\AAN{1}=\AAN{0}\cup \{ x_2x_3x_4 \}$, $\AAN{1}=\AAN{0}\cup \{ x_2^2x_1 \}$, $\AAN{1}=\AAN{0}\cup \{ x_1x_2x_3 \}$, $\AAN{1}=\AAN{0}\cup \{ x_2x_3x_5 \}$ we get extra conditions $c_4=-2c_5$, $4c_5=0$, $2c_4=0$, $c_5=-2c_4$ respectively. They all imply that $4c_4=0$, i.e.
$$
c_1=c_2=0,\; c_3=2c_4,\; 4c_4=2c_5=0.
$$
All such cases have already appeared.\\

Let $\AAN{0} =\{ x_0^2x_1, x_1^2x_2, x_2^2x_0, x_3^2x_0, x_4^2x_1, x_5^2x_2 \}$. Then 
$$
c_1=2c_3,\; c_2=-2c_3,\; 6c_3=2c_4=2(c_5+c_3)=0.
$$

For $\AAN{1}=\AAN{0}\cup \{ x_0x_2x_4 \}$, $\AAN{1}=\AAN{0}\cup \{ x_4^2x_0 \}$, $\AAN{1}=\AAN{0}\cup \{ x_4^2x_2 \}$, $\AAN{1}=\AAN{0}\cup \{ x_0^2x_2 \}$, $\AAN{1}=\AAN{0}\cup \{ x_0x_3x_4 \}$ we get extra conditions $c_4=4c_3$, $2c_3=0$, $4c_3=0$, $4c_3=0$, $c_4=c_3$ respectively. They all imply that $2c_3=0$, i.e.
$$
c_1=c_2=0,\; 2c_3=2c_4=2c_5=0.
$$
All such cases have already appeared.\\

For $\AAN{1}=\AAN{0}\cup \{ x_0x_4x_5 \}$ we get an extra condition $c_5=2c_3+c_4$, i.e.
$$
c_1=2c_3,\; c_2=-2c_3,\; c_5=2c_3+c_4,\; 6c_3=2c_4=0.
$$

For $\AAN{2}=\AAN{1}\cup \{ x_3^2x_1 \}$, $\AAN{2}=\AAN{1}\cup \{ x_5^2x_1 \}$, $\AAN{2}=\AAN{1}\cup \{ x_1^2x_4 \}$, $\AAN{2}=\AAN{1}\cup \{ x_1x_3x_4 \}$, $\AAN{2}=\AAN{1}\cup \{ x_1x_4x_5 \}$, $\AAN{2}=\AAN{1}\cup \{ x_3x_4x_5 \}$ we get extra conditions $2c_3=0$, $2c_3=0$, $c_4=-2c_3$, $c_4=-c_3$, $2c_3=0$, $c_3=0$. They all imply that $2c_3=0$, i.e.
$$
c_1=c_2=0,\; c_5=c_4,\; 2c_3=2c_4=0.
$$
All such cases have already appeared.\\

For $\AAN{2}=\AAN{1}\cup \{ x_1x_3x_5 \}$ we get an extra condition $c_4=3c_3$, i.e.
$$
c_1=2c_3,\; c_2=-2c_3,\; c_4=3c_3,\; c_5=-c_3,\; 6c_3=0.
$$

For $\AAN{3}=\AAN{2}\cup \{ x_3^2x_2 \}$, $\AAN{3}=\AAN{2}\cup \{ x_4^2x_2 \}$, $\AAN{3}=\AAN{2}\cup \{ x_2^2x_5 \}$, $\AAN{3}=\AAN{2}\cup \{ x_2x_3x_5 \}$, $\AAN{3}=\AAN{2}\cup \{ x_2x_4x_5 \}$ we get extra conditions $2c_3=0$, $2c_3=0$, $c_3=0$, $2c_3=0$, $2c_3=0$. They all imply that $2c_3=0$, i.e.
$$
c_1=c_2=0,\; c_4=c_5=c_3,\; 2c_3=0.
$$
All such cases have already appeared.\\

For $\AA=\AAN{2}\cup \{ x_2x_3x_4 \}$ the extra condition is automatically satisfied, i.e.
$$
c_1=2c_3,\; c_2=-2c_3,\; c_4=3c_3,\; c_5=-c_3,\; 6c_3=0.
$$
This means that $\GG\cong \ZZ{6}$ with a generator
$$
f=(1,{\om}^{2},{\om}^{-2},\om,{\om}^{3},{\om}^{-1}), \; \om=\sqr{6}.
$$
Then $2i_1+i_3 \equiv 2+2i_2+i_5+3i_4 \; mod \; 6$. This means that $\AAbar=\AA$.\\

For $\AAN{1}=\AAN{0}\cup \{ x_0^2x_3 \}$ we get an extra condition $c_3=0$, i.e.
$$
c_1=c_2=c_3=0,\; 2c_4=2c_5=0.
$$
All such cases have already appeared.\\

Let $\AAN{0} =\{ x_0^2x_1, x_1^2x_2, x_2^2x_0, x_3^2x_0, x_4^2x_0, x_5^2x_1 \}$. Then 
$$
c_1=2c_3,\; c_2=-2c_3,\; 6c_3=2(c_4-c_3)=2c_5=0.
$$

For $\AAN{1}=\AAN{0}\cup \{ x_2x_3x_4 \}$, $\AAN{1}=\AAN{0}\cup \{ x_3^2x_1 \}$, $\AAN{1}=\AAN{0}\cup \{ x_1x_3x_4 \}$ we get extra conditions $c_4=3c_3$, $2c_3=0$, $c_4=-c_3$ respectively. They all imply that $2c_3=0$, i.e.
$$
c_1=c_2=0,\; 2c_3=2c_4=2c_5=0.
$$
All such cases have already appeared.\\

For $\AAN{1}=\AAN{0}\cup \{ x_3x_4x_5 \}$ we get an extra condition $c_4=c_3+c_5$, i.e.
$$
c_1=2c_3,\; c_2=-2c_3,\; c_4=c_3+c_5,\; 6c_3=2c_5=0.
$$

For $\AAN{2}=\AAN{1}\cup \{ x_1x_3x_4 \}$, $\AAN{2}=\AAN{1}\cup \{ x_1x_2x_3 \}$, $\AAN{2}=\AAN{1}\cup \{ x_1x_2x_4 \}$, $\AAN{2}=\AAN{1}\cup \{ x_2^2x_1 \}$, $\AAN{2}=\AAN{1}\cup \{ x_3^2x_2 \}$, $\AAN{2}=\AAN{1}\cup \{ x_4^2x_2 \}$, $\AAN{2}=\AAN{1}\cup \{ x_3^2x_1 \}$, $\AAN{2}=\AAN{1}\cup \{ x_4^2x_1 \}$ we get extra conditions $c_5=2c_3$, $c_3=0$, $c_5=c_3$, $2c_3=0$, $2c_3=0$, $2c_3=0$, $2c_3=0$, $2c_3=0$ respectively. They all imply that $2c_3=0$, i.e.
$$
c_1=c_2=0,\; c_4=c_3+c_5,\; 2c_3=2c_5=0.
$$
All such cases have already appeared.\\

Let $\AAN{0} =\{ x_0^2x_1, x_1^2x_2, x_2^2x_0, x_3^2x_0, x_4^2x_0, x_5^2x_2 \}$. Then 
$$
c_1=2c_3,\; c_2=-2c_3,\; 6c_3=2(c_4-c_3)=2(c_5+c_3)=0.
$$

For $\AAN{1}=\AAN{0}\cup \{ x_2x_3x_4 \}$, $\AAN{1}=\AAN{0}\cup \{ x_1x_3x_4 \}$, $\AAN{1}=\AAN{0}\cup \{ x_3x_4x_5 \}$ we get extra conditions $c_4=3c_3$, $c_4=-c_3$, $c_3=c_4+c_5$ respectively. They all imply that $2c_3=0$, i.e.
$$
c_1=c_2=0,\; 2c_3=2c_4=2c_5=0.
$$
All such cases have already appeared.\\

Let $\AAN{0} =\{ x_0^2x_1, x_1^2x_2, x_2^2x_0, x_3^2x_0, x_4^2x_0, x_5^2x_0 \}$. Then 
$$
c_1=2c_3,\; c_2=-2c_3,\; 6c_3=2(c_4-c_3)=2(c_5-c_3)=0.
$$

For $\AAN{1}=\AAN{0}\cup \{ x_3x_4x_5 \}$ and $\AAN{1}=\AAN{0}\cup \{ x_1x_3x_4 \}$ we get extra conditions $c_3=c_4+c_5$ and $c_4=-c_3$ respectively. They both imply that $2c_3=0$, i.e.
$$
c_1=c_2=0,\; 2c_3=2c_4=2c_5=0.
$$
All such cases have already appeared.\\

This completes all the cases when a longest cycle has length at least $3$ or $\AAbar$ contains at least one cube.\\ 

At this point we could finish, because of the following observation:
\begin{center}
{\it If $\AA$ is smooth and contains no cycle of length at least $3$, then $\AAbar$ contains at least one cube.}
\end{center}

Indeed, if $x_0, x_1$ form a longest cycle in $\AA$, then $x_0^2x_1, x_1^2x_0\in \AA$. Hence $c_1=0$ and $\AAbar$ consists of monomials with $f$-weights equal to $0$ (which is the $f$-weight of $x_0^2x_1 \in \AA\subset \AAbar$). In particular, $x_0^3\in \AAbar$.\\

Still, for the sake of completeness we will continue our analysis in the same fashion as above.\\

\subsubsection{Case of no cubes. Length $2$ longest cycle.}

If a longest cycle has length $2$, then $\AA$ may be either

\begin{center}

\end{center}

or one of the following:

\begin{center}
%

\begin{flalign*}
\mbox{where}\; & A=(x_3^2x_0,\; x_3^2x_1,\; x_3^2x_4,\; x_5^2x_1,\; x_5^2x_0,\; x_0x_3x_5,\; x_1x_3x_5,\; x_3x_4x_5);\\
& B_0=(x_0x_3x_5,\; x_3^2x_0,\; x_0^2x_3,\; x_0^2x_5,\; x_5^2x_0,\; x_3^2x_4,\; x_0x_3x_4,\; x_0x_4x_5,\; x_3x_4x_5,\; x_0^2x_4);\\
& B_1=(x_1x_3x_5,\; x_3^2x_1,\; x_1^2x_3,\; x_1^2x_5,\; x_5^2x_1,\; x_3^2x_4,\; x_1x_3x_4,\; x_3x_4x_5);\\
& B_4=(x_3^2x_4,\; x_4^2x_3,\; x_3^2x_1,\; x_5^2x_1,\; x_1x_3x_4,\; x_1x_3x_5).\\
\end{flalign*}
\end{center}

Indeed, we may assume that $x_0, x_1$ form a cycle.\\

Consider the subgraph formed by the other vertices $x_2$, $x_3$, $x_4$, $x_5$. If it has cycles, then we may assume that $x_2, x_3$ form a cycle. The remaining vertices $4$, $5$ may be either connected or disconnected.\\

If $5$ is connected to $4$, then by Lemma 1 (\cite{Liendo}, Lemma 1.3) the $4$-th vertex should be connected either to $5$ or to one of the two cycles.\\

If $4$ is connected to $5$ (i.e. $x_4, x_5$ form another cycle), then there are neither singular pairs nor singular triples.\\

Suppose that $4$ and $5$ do not form a cycle and $4$ is connected to $2$. Then we may get
\begin{itemize}
\item one singular pair $(x_3,x_4)$ and 
\item two singular triples $(x_0, x_3,x_4)$, $(x_1, x_3,x_4)$.
\end{itemize}

We can resolve them as in the earlier examples.\\

Suppose that the remaining vertices $4$ and $5$ are disconnected. Then by Lemma 1 (\cite{Liendo}, Lemma 1.3) each of them should be connected to one of the cycles.\\

Suppose they are connected to different cycles (say, $4$ is connected to $0$ and $5$ is connected to $2$). Then we may get
\begin{itemize}
\item two singular pairs $(x_1,x_4)$, $(x_3,x_5)$ and 
\item six singular triples $(x_1, x_2,x_4)$, $(x_1, x_3,x_4)$, $(x_1, x_4,x_5)$, $(x_0, x_3,x_5)$, $(x_1, x_3,x_5)$,\\
$(x_3, x_4,x_5)$.
\end{itemize}

We resolve them as usual and obtain the pictures shown above.\\

Suppose that $4$ and $5$ are connected to the same cycle. Then they may be connected either to different vertices of that cycle or to the same vertex of one cycle (say, to $0$).\\

If $4$ is connected to $1$ and $5$ is connected to $0$, then we may get
\begin{itemize}
\item two singular pairs $(x_0,x_4)$, $(x_1,x_5)$ and 
\item four singular triples $(x_0, x_2,x_4)$, $(x_0, x_3,x_4)$, $(x_1, x_2,x_5)$, $(x_1, x_3,x_5)$.
\end{itemize}
 
If both $4$ and $5$ are connected to $0$, then we may get 
\begin{itemize}
\item three singular pairs $(x_1,x_4)$, $(x_1,x_5)$, $(x_4,x_5)$ and 
\item seven singular triples $(x_1, x_4,x_5)$, $(x_1, x_2,x_4)$, $(x_1, x_3,x_4)$, $(x_1, x_2,x_5)$, $(x_1, x_3,x_5)$,\\
$(x_2, x_4,x_5)$, $(x_3, x_4,x_5)$.
\end{itemize}

They are resolved as usual and we obtain the following pictures of $\AA$:

\begin{center}

\end{center}

Suppose that the subgraph formed by vertices $2$, $3$, $4$, $5$ has no cycles. Then the length of a longest path there may be either $4$ or $3$ or $2$ or $1$ (in which case the subgraph is totally disconnected).\\

If the length of a longest path in the subgraph formed by vertices $2$, $3$, $4$, $5$ is equal to $4$, then we may assume that $5$ is connected to $4$, $4$ is connected to $3$, $3$ is connected to $2$ and $2$ is connected to the length $2$ cycle (say, to $0$), because $2$, $3$, $4$, $5$ do not form cycles.\\ 

In this case we may get 
\begin{itemize}
\item one singular pair $(x_1,x_2)$ and 
\item two singular triples $(x_1, x_2,x_4)$, $(x_1, x_2,x_5)$.
\end{itemize}

We can resolve them as usual.\\

If the length of a longest path in the subgraph formed by vertices $2$, $3$, $4$, $5$ is equal to $3$, then we may assume that $4$ is connected to $3$, $3$ is connected to $2$ and $2$ is connected to the length $2$ cycle (say, to $0$), because $2$, $3$, $4$, $5$ do not form cycles.\\ 

By Lemma 1 (\cite{Liendo}, Lemma 1.3) the $5$-th vertex should be connected either to the cycle or to the length $3$ path.\\

If $5$ is connected to $1$, then we may get 
\begin{itemize}
\item two singular pairs $(x_1,x_2)$, $(x_0,x_5)$ and 
\item three singular triples $(x_1, x_2,x_4)$, $(x_0, x_3,x_5)$, $(x_0, x_4,x_5)$.
\end{itemize}

We can resolve them as usual.\\

This analysis gives us the following possiblities for $\AA$:

\begin{center}


\begin{flalign*}
\mbox{where}\; & A=(x_5^2x_2,\; x_5^2x_3,\; x_0^2x_2,\; x_0x_2x_5,\; x_0x_3x_5,\; x_0x_4x_5)\cup (x_0x_3x_5,\; x_3^2x_5,\; x_5^2x_3,\\
& x_3^2x_0,\; x_0x_4x_5,\; x_0x_3x_4,\; x_3x_4x_5)\cup (x_0x_4x_5,\; x_4^2x_5,\; x_4^2x_2,\; x_4^2x_0,\; x_0^2x_2,\\
& x_5^2x_2,\; x_0x_2x_5,\; x_0x_2x_4,\; x_2x_4x_5);\\
& B=(x_1x_2x_4,\; x_4^2x_5,\; x_4^2x_2,\; x_4^2x_1,\; x_1x_4x_5,\; x_1x_2x_5,\; x_2x_4x_5).\\
\end{flalign*}
\end{center}

If $5$ is connected to the length $3$ path, then it may be connected either to $2$ or to $3$ (since $2$, $3$, $4$, $5$ do not have length $4$ paths).\\

If $5$ is connected to $2$, then we may get 
\begin{itemize}
\item two singular pairs $(x_1,x_2)$, $(x_3,x_5)$ and 
\item three singular triples $(x_0, x_3,x_5)$, $(x_1, x_3,x_5)$, $(x_1, x_2,x_4)$.
\end{itemize}

If $5$ is connected to $3$, then we may get 
\begin{itemize}
\item two singular pairs $(x_1,x_2)$, $(x_4,x_5)$ and 
\item five singular triples $(x_1, x_2,x_4)$, $(x_1, x_2,x_5)$, $(x_0, x_4,x_5)$, $(x_1, x_4,x_5)$, $(x_2, x_4,x_5)$.
\end{itemize}

We resolve them as usual and obtain the following possiblities for $\AA$:

\begin{center}

\end{center}

\begin{center}
\begin{flalign*}
\mbox{where}\; & A=(x_5^2x_0,\; x_3^2x_0,\; x_3^2x_1,\; x_1x_3x_5,\; x_0x_3x_5,\; x_3x_4x_5)\cup (x_0x_3x_5,\; x_5^2x_0,\; x_3^2x_0,\\
& x_0x_4x_5, x_0x_3x_4, x_3x_4x_5)\cup (x_1x_3x_5, x_3^2x_1,\; x_1x_4x_5,\; x_1x_3x_4,\; x_3x_4x_5 );\\
& B=(x_5^2x_0,\; x_4^2x_0,\; x_1x_4x_5,\; x_0x_4x_5,\; x_2x_4x_5)\cup (x_0x_4x_5,\; x_0^2x_2,\; x_5^2x_0,\; x_4^2x_0,\\
& x_0x_2x_5,\; x_0x_2x_4,\; x_2x_4x_5);\\
& B_1= (x_5x_2x_4, x_1x_2x_4, x_1x_4x_5);\\
& C=(x_1x_2x_4, x_4^2x_5,\; x_4^2x_2,\; x_4^2x_1,\; x_1x_4x_5,\; x_1x_2x_5,\; x_2x_4x_5).
\end{flalign*}
\end{center}

If $5$ is connected to $0$, then we may get 
\begin{itemize}
\item three singular pairs $(x_1,x_2)$, $(x_1,x_5)$, $(x_2,x_5)$ and 
\item five singular triples $(x_1, x_2,x_4)$, $(x_1, x_2,x_5)$, $(x_2, x_4,x_5)$, $(x_1, x_3,x_5)$, $(x_1, x_4,x_5)$.
\end{itemize}

We resolve them as usual and obtain the following possiblities for $\AA$:

\begin{center}


\begin{flalign*}
\mbox{where}\; & C=(x_3^2x_1,\; x_3^2x_5,\; x_1x_3x_4,\; x_1x_3x_5,\; x_1x_4x_5,\; x_3x_4x_5).\\
\end{flalign*}
\end{center}

Suppose that the length of a longest path in the subgraph formed by vertices $2$, $3$, $4$, $5$ is equal to $2$. We may assume that $3$ is connected to $2$ and $2$ is connected to the length $2$ cycle (say, to $0$), because $2$, $3$, $4$, $5$ do not form cycles and there are no length $3$ paths.\\ 

The remaining vertices $4$ and $5$ may be either connected or disconnected.\\

If $5$ is connected to $4$, then by Lemma 1 (\cite{Liendo}, Lemma 1.3) the $4$-th vertex should be connected either to the cycle or to the length $2$ path formed by $2$ and $3$ (since $2$, $3$, $4$, $5$ do not form cycles). The latter is not possible, because there are no length $3$ paths. Hence $4$ is connected to the cycle.\\

If $4$ is connected to $1$, then we may get 
\begin{itemize}
\item two singular pairs $(x_1,x_2)$, $(x_0,x_4)$ and 
\item two singular triples $(x_1, x_2,x_5)$, $(x_0, x_3,x_4)$.
\end{itemize}

We can resolve them as usual and obtain the following possiblities for $\AA$:

\begin{center}


\begin{flalign*}
\mbox{where}\; & A=(x_3^2x_0,\; x_3^2x_4,\; x_3^2x_5,\; x_0x_3x_4,\; x_0x_3x_5,\; x_0x_4x_5,\; x_3x_4x_5)\cup (x_0x_2x_4,\; x_0x_3x_4,\\
& x_0x_4x_5);\\
& B=(x_5^2x_1,\; x_5^2x_2,\; x_5^2x_3,\; x_1x_2x_3,\; x_1x_2x_5,\; x_1x_3x_5,\; x_2x_3x_5).\\
\end{flalign*}
\end{center}

$\AA$ may be also one of the following:

\begin{center}
%

\begin{flalign*}
\mbox{where}\; & A=(x_3^2x_0,\; x_0x_3x_4,\; x_0x_3x_5,\; x_0x_4x_5,\; x_3x_4x_5)\cup (x_4^2x_0,\; x_0x_3x_4,\; x_0x_4x_5,\; x_0x_3x_5,\\
& x_3x_4x_5)\cup (x_0^2x_2,\; x_5^2x_2,\; x_0x_3x_5,\; x_0x_2x_5,\; x_0x_4x_5);\\
& B=(x_3^2x_0,\; x_3^2x_1,\; x_4^2x_0,\; x_4^2x_1,\; x_0x_3x_4,\; x_1x_3x_4,\; x_3x_4x_5);\\
& C=(x_3^2x_0,\; x_4^2x_0,\; x_0x_3x_4,\; x_0x_3x_5,\; x_0x_4x_5,\; x_3x_4x_5)\cup (x_3^2x_1,\; x_4^2x_1,\; x_1x_3x_4,\\
& x_1x_3x_5,\; x_1x_4x_5,\; x_3x_4x_5).\\
\end{flalign*}
\end{center}

Indeed, if $4$ is connected to $0$, then we may get 
\begin{itemize}
\item three singular pairs $(x_1,x_2)$, $(x_1,x_4)$, $(x_2,x_4)$ and 
\item three singular triples $(x_1, x_2,x_4)$, $(x_1, x_3,x_4)$, $(x_1, x_2,x_5)$.
\end{itemize}

We can resolve them as usual.\\

If the remaining vertices $4$ and $5$ are disconnected, then by Lemma 1 (\cite{Liendo}, Lemma 1.3) each of them should be connected either to $0$ or to $1$ or to $2$. They can not be connected to $3$, because there are no length $3$ paths.\\

If $4$ is connected to $2$ and $5$ is connected to $1$, then we may get 
\begin{itemize}
\item three singular pairs $(x_1,x_2)$, $(x_3,x_4)$, $(x_0,x_5)$ and 
\item five singular triples $(x_0, x_3,x_5)$, $(x_0, x_4,x_5)$, $(x_0, x_3,x_4)$, $(x_1, x_3,x_4)$, $(x_3, x_4,x_5)$.
\end{itemize}

We can resolve them as usual and obtain the pictures shown above.\\

If $4$ is connected to $0$ and $5$ is connected to $1$, then we may get 
\begin{itemize}
\item four singular pairs $(x_1,x_2)$, $(x_1,x_4)$, $(x_2,x_4)$, $(x_0,x_5)$ and 
\item four singular triples $(x_1, x_2,x_4)$, $(x_1, x_3,x_4)$, $(x_2, x_4,x_5)$, $(x_0, x_3,x_5)$.
\end{itemize}

We can resolve them as usual and obtain the following pictures of $\AA$:

\begin{center}


\begin{flalign*}
\mbox{where}\; & A=(x_5^2x_2,\; x_0^2x_2,\; x_0x_2x_5,\; x_0x_3x_5,\; x_0x_4x_5)\cup (x_0x_3x_5,\; x_3^2x_4,\; x_3^2x_5,\; x_3^2x_0,\\
& x_0x_4x_5,\; x_0x_3x_4,\; x_3x_4x_5); \\
& B=(x_4^2x_1,\; x_1x_2x_4,\; x_1x_3x_4,\; x_1x_4x_5)\cup (x_1x_3x_4,\; x_3^2x_4,\; x_3^2x_5,\; x_3^2x_1,\; x_4^2x_1,\\
& x_1x_4x_5,\; x_1x_3x_5,\; x_3x_4x_5)\cup (x_5^2x_2,\; x_2x_3x_5,\; x_2x_3x_4,\; x_2x_4x_5,\; x_3x_4x_5); \\
& C=(x_4^2x_1,\; x_1x_2x_4,\; x_1x_2x_3,\; x_1x_3x_4,\; x_2x_3x_4);\\
& C_0=(x_4^2x_1,\; x_1x_2x_4,\; x_1x_2x_5,\; x_1x_4x_5,\; x_2x_4x_5);\\
& E=(x_4^2x_1,\; x_5x_2x_4,\; x_1x_2x_4,\; x_2x_3x_4).\\
\end{flalign*}
\end{center}

If $4$ and $5$ are connected to the same vertex of the cycle which is different from the vertex to which $2$ is connected (say, $4$ and $5$ are connected to $0$, while $2$ is connected to $1$), then we may get 
\begin{itemize}
\item four singular pairs $(x_0,x_2)$, $(x_1,x_4)$, $(x_1,x_5)$, $(x_4,x_5)$ and 
\item five singular triples $(x_1, x_3,x_4)$, $(x_1, x_3,x_5)$, $(x_1, x_4,x_5)$, $(x_2, x_4,x_5)$, $(x_3, x_4,x_5)$.
\end{itemize}

We can resolve them as usual and obtain the following pictures of $\AA$:

\begin{center}

\end{center}

\begin{center}
\begin{flalign*}
\mbox{where}\; & A=(x_1x_2x_5,\; x_1x_3x_5)\cup (x_1x_2x_4,\; x_1x_3x_4); \;\;\; B= (x_0x_2x_3,\; x_0x_2x_4,\; x_0x_2x_5);\\
& C=(x_2x_3x_4,\; x_2x_3x_5,\; x_2x_4x_5,\; x_3x_4x_5).
\end{flalign*}
\end{center}

Suppose that $2$ is connected to $0$, while $4$ and $5$ are connected either to $0$ or to the length $2$ path (i.e. to $2$ since there are no length $3$ paths).\\ 

Then $\AA$ is one of the following:

\begin{center}
%

\begin{flalign*}
\mbox{where}\; & A=(x_5^2x_0,\; x_3^2x_0,\; x_0x_3x_5,\; x_3x_4x_5,\; x_1x_3x_5)\cup (x_5^2x_0,\; x_3^2x_0,\; x_0^2x_4,\; x_0x_3x_5,\\
& x_3x_4x_5,\; x_0x_3x_4,\; x_0x_4x_5)\cup (x_1x_3x_5,\; x_3x_4x_5,\; x_1x_3x_4,\; x_1x_4x_5);\\
& A_0=(x_0x_3x_5,\; x_1x_3x_5)\cup (x_0x_4x_5,\; x_1x_4x_5)\cup (x_0x_3x_5,\; x_0x_3x_4,\; x_0x_4x_5);\\
& C=(x_4^2x_2,\; x_1x_3x_4,\; x_1x_4x_5).
\end{flalign*}
\end{center}

Indeed, if $4$ is connected to $0$ and $5$ is connected to $2$, then we may get
\begin{itemize}
\item four singular pairs $(x_1,x_2)$, $(x_1,x_4)$, $(x_2,x_4)$, $(x_3,x_5)$ and 
\item six singular triples $(x_1, x_2,x_4)$, $(x_1, x_3,x_4)$, $(x_1, x_4,x_5)$, $(x_3, x_4,x_5)$, $(x_0, x_3,x_5)$,\\
$(x_1, x_3,x_5)$.
\end{itemize}

If both $4$ and $5$ are connected to $2$, then we may get
\begin{itemize}
\item four singular pairs $(x_1,x_2)$, $(x_3,x_4)$, $(x_3,x_5)$, $(x_4,x_5)$ and 
\item seven singular triples $(x_1, x_3,x_5)$, $(x_0, x_3,x_5)$, $(x_1, x_3,x_4)$, $(x_0, x_3,x_4)$, $(x_3, x_4,x_5)$,\\
$(x_1, x_4,x_5)$, $(x_0, x_4,x_5)$.
\end{itemize}

We resolve these singular pairs and triples as usual and obtain the pictures shown above.\\

If both $4$ and $5$ (as well as $2$) are connected to the same vertex of the cycle (i.e. to $0$), then we may get
\begin{itemize}
\item six singular pairs $(x_1,x_2)$, $(x_1,x_4)$, $(x_1,x_5)$, $(x_2,x_4)$, $(x_2,x_5)$, $(x_4,x_5)$ and 
\item seven singular triples $(x_1, x_2,x_4)$, $(x_1, x_2,x_5)$, $(x_1, x_4,x_5)$, $(x_2, x_4,x_5)$, $(x_1, x_3,x_4)$,\\
$(x_1, x_3,x_5)$, $(x_3, x_4,x_5)$.
\end{itemize}

We resolve them as usual and obtain the following pictures:

\begin{center}


where $B=(x_2x_3x_4, x_1x_2x_4)\cup (x_2x_3x_5, x_1x_2x_5)\cup (x_1x_2x_3, x_1x_2x_4, x_1x_2x_5)$.\\
\end{center}

Suppose that the subgraph formed by vertices $2$, $3$, $4$, $5$ is totally disconnected.\\

By Lemma 1 (\cite{Liendo}, Lemma 1.3) each of its vertices should be connected to the cycle formed by $0$ and $1$.\\

If pairs of the vertices $2$, $3$, $4$, $5$ are connected to two different vertices of the cycle (say, $2$, $3$ are connected to $0$ and $4$, $5$ are connected to $1$), then $\AA$ may be one of the following:

\begin{center}
%
\end{center}

\begin{center}
\begin{flalign*}
\mbox{where}\; & A=(x_2x_3x_4,\; x_2x_3x_5,\; x_2x_4x_5,\; x_3x_4x_5);\\
& C=(x_4^2x_0,\; x_5^2x_0,\; x_2x_4x_5)\cup (x_4^2x_0,\; x_0x_2x_4,\; x_0x_3x_4)\cup (x_5^2x_0,\; x_0x_2x_5,\; x_0x_3x_5)\cup \\
& \cup(x_4^2x_0,\; x_5^2x_0,\; x_3x_4x_5,\; x_0x_3x_4\; x_0x_3x_5).
\end{flalign*}
\end{center}

Indeed, we may get 
\begin{itemize}
\item six singular pairs $(x_1,x_2)$, $(x_1,x_3)$, $(x_2,x_3)$, $(x_0,x_4)$, $(x_0,x_5)$, $(x_4,x_5)$ and
\item six singular triples $(x_1,x_2,x_3)$, $(x_2,x_3,x_4)$, $(x_2,x_3,x_5)$, $(x_0,x_4,x_5)$, $(x_2,x_4,x_5)$,\\
$(x_3,x_4,x_5)$.
\end{itemize}

Note that the configuration is symmetric with respect to vertices $2$ and $3$, with respect to vertices $4$ and $5$ and with respect to triples of vertices $(0,2,3)$ and $(1,4,5)$.\\

We resolve these singular pairs and triples as usual and obtain the pictures shown above.\\

If three of the vertices $2$, $3$, $4$, $5$ are connected to the same vertex of the cycle (say, $3$, $4$ and $5$ are connected to $0$) and the other one is connected to the other vertex of the cycle (i.e. $2$ is connected to $1$), then we may get 
\begin{itemize}
\item seven singular pairs $(x_1,x_3)$, $(x_1,x_4)$, $(x_1,x_5)$, $(x_3,x_4)$, $(x_3,x_5)$, $(x_4,x_5)$, $(x_0,x_2)$ and
\item seven singular triples $(x_1,x_3,x_4)$, $(x_1,x_3,x_5)$, $(x_1,x_4,x_5)$, $(x_3,x_4,x_5)$, $(x_2,x_3,x_4)$,\\
$(x_2,x_3,x_5)$, $(x_2,x_4,x_5)$.
\end{itemize}

They can be resolved as in the earlier examples.\\

Finally, it may happen that all four vertices $2$, $3$, $4$ and $5$ are connected to one and the same vertex of the cycle (say, to $0$). In this case we may get
\begin{itemize}
\item ten singular pairs $(x_1,x_2)$, $(x_1,x_3)$, $(x_1,x_4)$, $(x_1,x_5)$, $(x_2,x_3)$, $(x_2,x_4)$, $(x_2,x_5)$, $(x_3,x_4)$, $(x_3,x_5)$, $(x_4,x_5)$ and
\item ten singular triples $(x_1,x_2,x_3)$, $(x_1,x_2,x_4)$, $(x_1,x_2,x_5)$, $(x_1,x_3,x_4)$, $(x_1,x_3,x_5)$,\\
$(x_1,x_4,x_5)$, $(x_2,x_3,x_4)$, $(x_2,x_3,x_5)$, $(x_2,x_4,x_5)$, $(x_3,x_4,x_5)$.
\end{itemize}

They can be resolved as usual.\\

The resulting pictures are as follows:

\begin{center}

\end{center}

Now let us do computations.\\

If $\AA =\{ x_0^2x_1, x_1^2x_0, x_2^2x_3, x_3^2x_2, x_4^2x_5, x_5^2x_4 \}$, then 
$$
c_1=0,\; c_2=c_3,\; c_4=c_5,\; 3c_2=3c_4=0.
$$
This means that $\GG\cong (\ZZ{3})^{\oplus 2}$ with generators
$$
f_1=(1,1,\om,\om,1,1), \; f_2=(1,1,1,1,\om,\om), \; \om=\sqr{3}.
$$
Then $i_2+i_3 \equiv i_4+i_5 \equiv 0 \; mod \; 3$. This means that $\AAbar=\AA\cup \{ f_3(x_0,x_1), f_3(x_2,x_3), f_3(x_4,x_5) \}$.\\

Let $\AAN{0} =\{ x_0^2x_1, x_1^2x_0, x_2^2x_3, x_3^2x_2, x_4^2x_2, x_5^2x_4 \}$. Then 
$$
c_1=0,\; c_2=c_3=4c_5,\; c_4=-2c_5,\; 12c_5=0.
$$

For $\AA =\AAN{0}\cup \{ x_3x_4x_5 \}$ we get an extra condition $3c_5=0$, i.e.
$$
c_1=0,\; c_2=c_3=c_4=c_5,\; 3c_5=0.
$$
This means that $\GG\cong \ZZ{3}$ with a generator
$$
f=(1,1,\om,\om,\om,\om), \; \om=\sqr{3}.
$$
Then $i_2+i_3+i_4+i_5 \equiv 0 \; mod \; 3$. This means that $\AAbar=\AA\cup \{ f_3(x_0,x_1), f_3(x_2,x_3,x_4,x_5) \}$.\\

For $\AA =\AAN{0}\cup \{ x_4^2x_3 \}$ the extra condition is automatically satisfied.\\

This means that $\GG\cong \ZZ{12}$ with a generator
$$
f=(1,1,{\om}^{4},{\om}^{4},{\om}^{-2},\om), \; \om=\sqr{12}.
$$
Then $4(i_2+i_3)+i_5 \equiv 2i_4 \; mod \; 12$. This means that $\AAbar=\AA\cup \{ f_3(x_0,x_1), f_3(x_2,x_3) \}$.\\

For $\AAN{1} =\AAN{0}\cup \{ x_1x_3x_4 \}$, $\AAN{1} =\AAN{0}\cup \{ x_4^2x_0 \}$, $\AAN{1} =\AAN{0}\cup \{ x_3^2x_1 \}$ we get extra conditions $2c_5=0$, $4c_5=0$, $4c_5=0$ respectively. They all imply that $4c_5=0$, i.e.
$$
c_1=c_2=c_3=0,\; c_4=2c_5,\; 4c_5=0.
$$
All such cases have already appeared.\\

Let $\AAN{0} =\{ x_0^2x_1, x_1^2x_0, x_2^2x_3, x_3^2x_2, x_4^2x_0, x_5^2x_2 \}$. Then 
$$
c_1=0,\; c_2=c_3=-2c_5,\; 2c_4=6c_5=0.
$$

For $\AAN{1} =\AAN{0}\cup \{ x_4^2x_1 \}$ the extra condition is automatically satisfied.\\

For $\AA =\AAN{1}\cup \{ x_5^2x_3 \}$ the extra condition is automatically satisfied.\\

This means that $\GG\cong \ZZ{6}\oplus \ZZ{2}$ with generators
$$
f_1=(1,1,{\om}^{-2},{\om}^{-2},1,\om), \; f_2=(1,1,1,1,\eta,1),\; \om=\sqr{6}, \eta=\sqr{2}.
$$
Then $i_5 \equiv 2(i_2+i_3) \; mod \; 6$ and $i_4\equiv 0 \; mod \; 2$. This means that $\AAbar=\AA\cup \{ f_3(x_0,x_1)$, $f_3(x_2,x_3) \}$.\\

For $\AAN{2} =\AAN{1}\cup \{ x_3^2x_1 \}$, $\AAN{2} =\AAN{1}\cup \{ x_3^2x_4 \}$, $\AAN{2} =\AAN{1}\cup \{ x_3^2x_0 \}$, $\AAN{2} =\AAN{1}\cup \{ x_5^2x_0 \}$, $\AAN{2} =\AAN{1}\cup \{ x_5^2x_1 \}$, $\AAN{2} =\AAN{1}\cup \{ x_0x_3x_5 \}$, $\AAN{2} =\AAN{1}\cup \{ x_1x_3x_5 \}$, $\AAN{2} =\AAN{1}\cup \{ x_3x_4x_5 \}$ we get extra conditions $2c_5=0$, $c_4=4c_5$, $2c_5=0$, $2c_5=0$, $2c_5=0$, $c_5=0$, $c_5=0$, $c_4=c_5$. They all imply that $2c_5=0$, i.e.
$$
c_1=c_2=c_3=0,\; 2c_4=2c_5=0.
$$
All such cases have already appeared.\\

This completes the analysis also for $\AAN{1} =\AAN{0}\cup \{ x_1^2x_3 \}$, $\AAN{1} =\AAN{0}\cup \{ x_1^2x_2 \}$, $\AAN{1} =\AAN{0}\cup \{ x_1x_3x_4 \}$, $\AAN{1} =\AAN{0}\cup \{ x_4^2x_3 \}$, $\AAN{1} =\AAN{0}\cup \{ x_1^2x_5 \}$, $\AAN{1} =\AAN{0}\cup \{ x_4^2x_2 \}$, $\AAN{1} =\AAN{0}\cup \{ x_1x_2x_4 \}$ and $\AAN{1} =\AAN{0}\cup \{ x_1x_4x_5 \}$.\\

Let $\AAN{0} =\{ x_0^2x_1, x_1^2x_0, x_2^2x_3, x_3^2x_2, x_4^2x_0, x_5^2x_0 \}$. Then 
$$
c_1=0,\; c_2=c_3,\; 3c_2=2c_4=2c_5=0.
$$

For $\AA =\AAN{0}\cup \{ x_1x_4x_5 \}$ we get an extra condition $c_5=c_4$, i.e.
$$
c_1=0,\; c_2=c_3,\; c_5=c_4,\; 3c_2=2c_4=0.
$$
This case has already appeared.\\

For $\AAN{1} =\AAN{0}\cup \{ x_2x_4x_5 \}$ we get an extra condition $c_2=c_4+c_5$, i.e.
$$
c_1=c_2=c_3=0,\; c_4=c_5,\; 2c_5=0.
$$
All such cases have already appeared.\\

Let $\AAN{0} =\{ x_0^2x_1, x_1^2x_0, x_2^2x_3, x_3^2x_2, x_4^2x_1, x_5^2x_0 \}$. Then 
$$
c_1=0,\; c_2=c_3,\; 3c_2=2c_4=2c_5=0.
$$

For $\AAN{1} =\AAN{0}\cup \{ x_0x_2x_4 \}$ and $\AAN{1} =\AAN{0}\cup \{ x_0^2x_2 \}$ we get extra conditions $c_2=c_4$ and $c_2=0$ respectively. They both imply that $c_2=0$, i.e.
$$
c_1=c_2=c_3=0,\; 2c_4=2c_5=0.
$$
All such cases have already appeared.\\

For $\AA =\AAN{0}\cup \{ x_0x_4x_5, x_1x_4x_5 \}$ we get an extra condition $c_5=c_4$, i.e.
$$
c_1=0,\; c_2=c_3,\; c_5=c_4,\; 3c_2=2c_4=0.
$$
This case has already appeared.\\

For $\AA =\AAN{0}\cup \{ x_5^2x_1, x_4^2x_0 \}$ extra conditions are automatically satisfied.\\

This case has already appeared.\\

Let $\AAN{0} =\{ x_0^2x_1, x_1^2x_0, x_2^2x_3, x_3^2x_2, x_4^2x_0, x_4^2x_1, x_5^2x_1 \}$. Then 
$$
c_1=0,\; c_2=c_3,\; 3c_2=2c_4=2c_5=0.
$$

For $\AAN{1} =\AAN{0}\cup \{ x_0x_2x_5 \}$ and $\AAN{1} =\AAN{0}\cup \{ x_0^2x_2 \}$ we get extra conditions $c_2=c_5$ and $c_2=0$ respectively. They both imply that $c_2=0$, i.e.
$$
c_1=c_2=c_3=0,\; 2c_4=2c_5=0.
$$
All such cases have already appeared.\\

For $\AA =\AAN{0}\cup \{ x_0x_4x_5 \}$ we get an extra condition $c_5=c_4$, i.e.
$$
c_1=0,\; c_2=c_3,\; c_5=c_4,\; 3c_2=2c_4=0.
$$
This case has already appeared.\\

Let $\AAN{0} =\{ x_0^2x_1, x_1^2x_0, x_2^2x_0, x_3^2x_2, x_4^2x_3, x_5^2x_4 \}$. Then 
$$
c_1=0,\; c_2=-8c_5,\; c_3=4c_5,\; c_4=-2c_5,\; 16c_5=0.
$$

For $\AA =\AAN{0}\cup \{ x_2^2x_1 \}$ the extra condition is automatically satisfied.\\

This means that $\GG\cong \ZZ{16}$ with a generator
$$
f=(1,1,{\om}^{8},{\om}^{4},{\om}^{-2},\om), \; \om=\sqr{16}.
$$
Then $4i_3+i_5 \equiv 2i_4+8i_2 \; mod \; 16$. This means that $\AAbar=\AA\cup \{ f_3(x_0,x_1) \}$.\\

For $\AAN{1} =\AAN{0}\cup \{ x_1x_2x_3 \}$, $\AAN{1} =\AAN{0}\cup \{ x_1x_2x_4 \}$, $\AA=\AAN{0}\cup \{ x_1x_2x_5 \}$ we get extra conditions $4c_5=0$, $2c_5=0$, $c_5=0$ respectively. They all imply that $4c_5=0$, i.e.
$$
c_1=c_2=c_3=0,\; c_4=2c_5,\; 4c_5=0.
$$
All such cases have already appeared.\\

Let $\AAN{0} =\{ x_0^2x_1, x_1^2x_0, x_2^2x_0, x_3^2x_2, x_4^2x_3, x_5^2x_1 \}$. Then 
$$
c_1=0,\; c_2=4c_4,\; c_3=-2c_4,\; 8c_4=2c_5=0.
$$

For $\AAN{1} =\AAN{0}\cup \{ x_5^2x_0 \}$ the extra condition is automatically satisfied.\\

For $\AA =\AAN{1}\cup \{ x_2^2x_1 \}$ the extra condition is automatically satisfied.\\

This case has already appeared.\\

For $\AAN{1} =\AAN{0}\cup \{ x_1x_2x_5 \}$ we get an extra condition $c_5=4c_4$, i.e.
$$
c_1=0,\; c_2=c_5=4c_4,\; c_3=-2c_4,\; 8c_4=0.
$$
All such cases have already appeared.\\

For $\AAN{1} =\AAN{0}\cup \{ x_2^2x_1 \}$ the extra condition is automatically satisfied.\\

For $\AAN{2} =\AAN{1}\cup \{ x_5^2x_2 \}$, $\AAN{2} =\AAN{1}\cup \{ x_5^2x_3 \}$, $\AAN{2} =\AAN{1}\cup \{ x_0^2x_2 \}$, $\AAN{2} =\AAN{1}\cup \{ x_0x_3x_5 \}$, $\AAN{2} =\AAN{1}\cup \{ x_0x_4x_5 \}$ we get extra conditions $4c_4=0$, $2c_4=0$, $4c_4=0$, $c_5=2c_4$, $c_5=-c_4$ respectively. They all imply that $4c_4=0$, i.e.
$$
c_1=c_2=0,\; c_3=2c_4,\; 4c_4=2c_5=0.
$$
All such cases have already appeared.\\

For $\AAN{2} =\AAN{1}\cup \{ x_0x_2x_5 \}$ we get an extra condition $c_5=4c_4$, i.e.
$$
c_1=0,\; c_2=c_5=4c_4,\; c_3=-2c_4,\; 8c_4=0.
$$
All such cases have already appeared.\\

For $\AA =\AAN{0}\cup \{ x_1x_2x_3 \}$ and $\AA =\AAN{0}\cup \{ x_1x_2x_4 \}$ we get extra conditions $2c_4=0$ and $c_4=0$ respectively. They both imply that $2c_4=0$, i.e.
$$
c_1=c_2=c_3=0,\; 2c_4=2c_5=0.
$$
This case has already appeared.\\

Let $\AAN{0} =\{ x_0^2x_1, x_1^2x_0, x_2^2x_0, x_3^2x_2, x_4^2x_3, x_5^2x_3 \}$. Then 
$$
c_1=0,\; c_2=4c_4,\; c_3=-2c_4,\; 8c_4=2(c_5-c_4)=0.
$$

For $\AAN{1} =\AAN{0}\cup \{ x_1x_2x_4 \}$, $\AAN{1} =\AAN{0}\cup \{ x_1x_2x_3 \}$, $\AAN{1} =\AAN{0}\cup \{ x_5^2x_2 \}$, $\AAN{1} =\AAN{0}\cup \{ x_1x_2x_5 \}$ we get extra conditions $c_4=0$, $2c_4=0$, $2c_4=0$, $c_5=4c_4$. They all imply that $2c_4=0$, i.e.
$$
c_1=c_2=c_3=0,\; 2c_4=2c_5=0.
$$
All such cases have already appeared.\\

Let $\AAN{0} =\{ x_0^2x_1, x_1^2x_0, x_2^2x_0, x_3^2x_2, x_4^2x_3, x_5^2x_2 \}$. Then 
$$
c_1=0,\; c_2=4c_4,\; c_3=-2c_4,\; 8c_4=2(c_5-2c_4)=0.
$$

For $\AAN{1} =\AAN{0}\cup \{ x_1x_2x_3 \}$, $\AAN{1} =\AAN{0}\cup \{ x_1x_2x_4 \}$, $\AAN{1} =\AAN{0}\cup \{ x_1x_2x_5 \}$ we get extra conditions $2c_4=0$, $c_4=0$, $c_5=4c_4$. They all imply that $4c_4=0$, i.e.
$$
c_1=c_2=0,\; c_3=2c_4,\; 4c_4=2c_5=0.
$$
All such cases have already appeared.\\

For $\AAN{1} =\AAN{0}\cup \{ x_2^2x_1 \}$ the extra condition is automatically satisfied.\\

For $\AAN{2} =\AAN{1}\cup \{ x_5^2x_0 \}$, $\AAN{2} =\AAN{1}\cup \{ x_3^2x_0 \}$, $\AAN{2} =\AAN{1}\cup \{ x_3^2x_1 \}$, $\AAN{2} =\AAN{1}\cup \{ x_3x_4x_5 \}$ we get extra conditions $4c_4=0$, $4c_4=0$, $4c_4=0$, $c_5=c_4$ respectively. They all imply that $4c_4=0$, i.e.
$$
c_1=c_2=0,\; c_3=2c_4,\; 4c_4=2c_5=0.
$$
All such cases have already appeared.\\

For $\AAN{2} =\AAN{1}\cup \{ x_0x_3x_5 \}$, $\AAN{2} =\AAN{1}\cup \{ x_1x_3x_5 \}$ we get the same extra condition $c_5=2c_4$, i.e.
$$
c_1=0,\; c_2=4c_4,\; c_3=-2c_4,\; c_5=2c_4,\; 8c_4=0.
$$
All such cases have already appeared.\\

Let $\AAN{0} =\{ x_0^2x_1, x_1^2x_0, x_2^2x_0, x_2^2x_1, x_3^2x_2, x_4^2x_3, x_5^2x_3 \}$. Then 
$$
c_1=0,\; c_2=4c_4,\; c_3=-2c_4,\; 8c_4=2(c_5-c_4)=0.
$$

For $\AAN{1} =\AAN{0}\cup \{ x_5^2x_0 \}$ and $\AAN{1} =\AAN{0}\cup \{ x_4^2x_0 \}$ we get the same extra condition $2c_4=0$, i.e.
$$
c_1=c_2=c_3=0,\; 2c_4=2c_5=0.
$$
All such cases have already appeared.\\

For $\AAN{1} =\AAN{0}\cup \{ x_0x_4x_5 \}$, $\AAN{1} =\AAN{0}\cup \{ x_1x_4x_5 \}$ and $\AA =\AAN{0}\cup \{ x_2x_4x_5 \}$ we get the same extra condition $c_5=3c_4=-c_4$, i.e.
$$
c_1=c_2=0,\; c_3=2c_4,\; c_5=-c_4,\; 4c_4=0.
$$
All such cases have already appeared.\\

Let $\AAN{0} =\{ x_0^2x_1, x_1^2x_0, x_2^2x_0, x_3^2x_2, x_4^2x_3, x_5^2x_0 \}$. Then 
$$
c_1=0,\; c_2=4c_4,\; c_3=-2c_4,\; 8c_4=2c_5=0.
$$

For $\AAN{1} =\AAN{0}\cup \{ x_1x_2x_5 \}$ we get an extra condition $c_5=4c_4$, i.e.
$$
c_1=0,\; c_2=c_5=4c_4,\; c_3=-2c_4,\; 8c_4=0.
$$
All such cases have already appeared.\\

For $\AAN{1} =\AAN{0}\cup \{ x_1x_2x_3 \}$ and $\AAN{1} =\AAN{0}\cup \{ x_1x_2x_4 \}$ we get extra conditions $2c_4=0$ and $c_4=0$ respectively. They both imply that $2c_4=0$, i.e.
$$
c_1=c_2=c_3=0,\; 2c_4=2c_5=0.
$$
All such cases have already appeared.\\

Let $\AAN{0} =\{ x_0^2x_1, x_1^2x_0, x_2^2x_0, x_3^2x_2, x_4^2x_1, x_5^2x_4 \}$ or $\AAN{0} =\{ x_0^2x_1, x_1^2x_0, x_2^2x_0, x_3^2x_2, x_4^2x_0, x_5^2x_4 \}$. Then 
$$
c_1=0,\; c_2=2c_3,\; c_4=2c_5,\; 4c_3=4c_5=0.
$$
All such cases have already appeared.\\

Let $\AAN{0} =\{ x_0^2x_1, x_1^2x_0, x_2^2x_0, x_3^2x_2, x_4^2x_2, x_5^2x_1 \}$. Then 
$$
c_1=0,\; c_2=2c_3,\; 4c_3=2(c_4-c_3)=2c_5=0.
$$

For $\AAN{1} =\AAN{0}\cup \{ x_3^2x_0 \}$, $\AAN{1} =\AAN{0}\cup \{ x_4^2x_0 \}$, $\AAN{1} =\AAN{0}\cup \{ x_3^2x_1 \}$, $\AAN{1} =\AAN{0}\cup \{ x_4^2x_1 \}$ we get the same extra condition $2c_3=0$, i.e.
$$
c_1=c_2=0,\; 2c_3=2c_4=2c_5=0.
$$
All such cases have already appeared.\\

For $\AAN{1} =\AAN{0}\cup \{ x_0x_3x_4 \}$ and $\AAN{1} =\AAN{0}\cup \{ x_1x_3x_4 \}$ we get the same extra condition $c_4=-c_3$, i.e.
$$
c_1=0,\; c_2=2c_3,\; c_4=-c_3,\; 4c_3=2c_5=0.
$$
All such cases have already appeared.\\

For $\AAN{1} =\AAN{0}\cup \{ x_3x_4x_5 \}$ we get an extra condition $c_4=c_5-c_3$, i.e.
$$
c_1=0,\; c_2=2c_3,\; c_4=c_5-c_3,\; 4c_3=2c_5=0.
$$
All such cases have already appeared.\\

Let $\AAN{0} =\{ x_0^2x_1, x_1^2x_0, x_2^2x_0, x_3^2x_2, x_4^2x_0, x_5^2x_1 \}$. Then 
$$
c_1=0,\; c_2=2c_3,\; 4c_3=2c_4=2c_5=0.
$$

For $\AAN{1} =\AAN{0}\cup \{ x_5^2x_2 \}$, $\AAN{1} =\AAN{0}\cup \{ x_2x_3x_4 \}$, $\AAN{1} =\AAN{0}\cup \{ x_2x_3x_5 \}$, $\AAN{1} =\AAN{0}\cup \{ x_3x_4x_5 \}$ we get extra conditions $2c_3=0$, $c_4=c_3$, $c_5=c_3$, $c_3=c_4+c_5$. They all imply that $2c_3=0$, i.e.
$$
c_1=c_2=0,\; 2c_3=2c_4=2c_5=0.
$$
All such cases have already appeared.\\

For $\AAN{1} =\AAN{0}\cup \{ x_2x_4x_5 \}$ we get an extra condition $c_5=c_4+2c_3$, i.e.
$$
c_1=0,\; c_2=2c_3,\; c_5=2c_3+c_4,\; 4c_3=2c_4=0.
$$
All such cases have already appeared.\\

Let $\AAN{0} =\{ x_0^2x_1, x_1^2x_0, x_2^2x_1, x_3^2x_2, x_4^2x_0, x_5^2x_0 \}$. Then 
$$
c_1=0,\; c_2=2c_3,\; 4c_3=2c_4=2c_5=0.
$$

For $\AAN{1} =\AAN{0}\cup \{ x_2x_4x_5 \}$ we get an extra condition $c_5=2c_3+c_4$, i.e.
$$
c_1=0,\; c_2=2c_3,\; c_5=2c_3+c_4,\; 4c_3=2c_4=0.
$$

For $\AAN{2} =\AAN{1}\cup \{ x_1x_3x_4 \}$,  $\AAN{2} =\AAN{1}\cup \{ x_1x_3x_5 \}$ and  $\AAN{2} =\AAN{1}\cup \{ x_3x_4x_5 \}$ we get extra conditions $c_4=c_3$, $c_5=c_3$ and $c_3=0$ respectively. They all imply that $2c_3=0$, i.e.
$$
c_1=c_2=0,\; c_5=c_4,\; 2c_3=2c_4=0.
$$
All such cases have already appeared.\\

For $\AAN{1} =\AAN{0}\cup \{ x_3x_4x_5 \}$ we get an extra condition $c_5=c_4-c_3$, i.e.
$$
c_1=c_2=0,\; c_5=c_3+c_4,\; 2c_3=2c_4=0.
$$
All such cases have already appeared.\\

For $\AAN{1} =\AAN{0}\cup \{ x_1x_4x_5 \}$ we get an extra condition $c_5=c_4$, i.e.
$$
c_1=0,\; c_2=2c_3,\; c_5=c_4,\; 4c_3=2c_4=0.
$$

For $\AAN{2} =\AAN{1}\cup \{ x_2x_4x_5 \}$, $\AAN{2} =\AAN{1}\cup \{ x_3x_4x_5 \}$, $\AAN{2} =\AAN{1}\cup \{ x_2x_3x_4 \}$, $\AAN{2} =\AAN{1}\cup \{ x_2x_3x_5 \}$ we get extra conditions $2c_3=0$, $c_3=0$, $c_4=c_3$, $c_5=c_3$ respectively. They all imply that $2c_3=0$, i.e.
$$
c_1=c_2=0,\; c_5=c_4,\; 2c_3=2c_4=0.
$$
All such cases have already appeared.\\

Let $\AAN{0} =\{ x_0^2x_1, x_1^2x_0, x_2^2x_0, x_3^2x_2, x_4^2x_0, x_5^2x_2 \}$. Then 
$$
c_1=0,\; c_2=2c_3,\; 4c_3=2c_4=2(c_5-c_3)=0.
$$

For $\AAN{1} =\AAN{0}\cup \{ x_1x_2x_4 \}$ we get an extra condition $c_4=2c_3$, i.e.
$$
c_1=0,\; c_2=c_4=2c_3,\; 4c_3=2(c_5-c_3)=0.
$$
All such cases have already appeared.\\

For $\AAN{1} =\AAN{0}\cup \{ x_1x_2x_3 \}$ we get an extra condition $c_3=0$, i.e.
$$
c_1=c_2=c_3=0,\; 2c_4=2c_5=0.
$$
All such cases have already appeared.\\

Let $\AAN{0} =\{ x_0^2x_1, x_1^2x_0, x_2^2x_0, x_3^2x_2, x_4^2x_2, x_5^2x_2 \}$. Then 
$$
c_1=0,\; c_2=2c_3,\; 4c_3=2(c_4-c_3)=2(c_5-c_3)=0.
$$

For $\AAN{1} =\AAN{0}\cup \{ x_1x_3x_4 \}$ we get an extra condition $c_4=-c_3$, i.e.
$$
c_1=0,\; c_2=2c_3,\; c_4=-c_3,\; 4c_3=2(c_5-c_3)=0.
$$

For $\AAN{2} =\AAN{1}\cup \{ x_1x_2x_3 \}$ we get an extra condition $c_3=0$, i.e.
$$
c_1=c_2=c_3=c_4=0,\; 2c_5=0.
$$
All such cases have already appeared.\\

For $\AAN{2} =\AAN{1}\cup \{ x_0x_3x_5 \}$ and $\AAN{2} =\AAN{1}\cup \{ x_1x_3x_5 \}$ we get the same extra condition $c_5=-c_3$, i.e.
$$
c_1=0,\; c_2=2c_3,\; c_4=c_5=-c_3,\; 4c_3=0.
$$

For $\AAN{3} =\AAN{2}\cup \{ x_0x_4x_5 \}$ and $\AAN{3} =\AAN{2}\cup \{ x_1x_4x_5 \}$ we get the same extra condition $2c_3=0$, i.e.
$$
c_1=c_2=0,\; c_4=c_5=c_3,\; 2c_3=0.
$$
All such cases have already appeared.\\

For $\AAN{1} =\AAN{0}\cup \{ x_1x_2x_3, x_1x_4x_5 \}$ we get extra conditions $c_3=0$ and $c_5=-c_4$, i.e.
$$
c_1=c_2=c_3=0,\; c_4=c_5,\; 2c_5=0.
$$
All such cases have already appeared.\\

For $\AAN{1} =\AAN{0}\cup \{ x_3x_4x_5 \}$ we get an extra condition $c_5=-c_3-c_4$, i.e.
$$
c_1=c_2=0,\; c_5=c_3+c_4,\; 2c_3=2c_4=0.
$$
All such cases have already appeared.\\

Let $\AAN{0} =\{ x_0^2x_1, x_1^2x_0, x_2^2x_0, x_3^2x_2, x_4^2x_0, x_5^2x_0 \}$. Then 
$$
c_1=0,\; c_2=2c_3,\; 4c_3=2c_4=2c_5=0.
$$

For $\AAN{1} =\AAN{0}\cup \{ x_2x_4x_5 \}$ we get an extra condition $c_5=2c_3+c_4$, i.e.
$$
c_1=0,\; c_2=2c_3,\; c_5=2c_3+c_4,\; 4c_3=2c_4=0.
$$

For $\AAN{2} =\AAN{1}\cup \{ x_1x_4x_5 \}$ and $\AAN{2} =\AAN{1}\cup \{ x_1x_3x_4 \}$ we get extra conditions $2c_3=0$ and $c_4=c_3$ respectively. They both imply that $2c_3=0$, i.e.
$$
c_1=c_2=0,\; c_5=c_4,\; 2c_3=2c_4=0.
$$
All such cases have already appeared.\\

For $\AAN{2} =\AAN{1}\cup \{ x_1x_2x_4 \}$ we get an extra condition $c_4=2c_3$, i.e.
$$
c_1=c_5=0,\; c_2=c_4=2c_3,\; 4c_3=0.
$$
All such cases have already appeared.\\

For $\AAN{1} =\AAN{0}\cup \{ x_3x_4x_5 \}$ we get an extra condition $c_3=c_4+c_5$, i.e.
$$
c_1=c_2=0,\; c_3=c_4+c_5,\; 2c_4=2c_5=0.
$$
All such cases have already appeared.\\

For $\AAN{1} =\AAN{0}\cup \{ x_1x_4x_5, x_2x_3x_5 \}$ we get extra conditions $c_3=c_4=c_5$, i.e.
$$
c_1=c_2=0,\; c_3=c_4=c_5,\; 2c_3=0.
$$
All such cases have already appeared.\\

Let $\AAN{0} =\{ x_0^2x_1, x_1^2x_0, x_2^2x_0, x_3^2x_0, x_4^2x_1, x_5^2x_1 \}$. Then 
$$
c_1=0,\; 2c_2=2c_3=2c_4=2c_5=0.
$$

For $\AAN{1} =\AAN{0}\cup \{ x_2x_3x_4 \}$, $\AAN{1} =\AAN{0}\cup \{ x_2x_3x_5 \}$, $\AAN{1} =\AAN{0}\cup \{ x_2x_4x_5 \}$, $\AAN{1} =\AAN{0}\cup \{ x_3x_4x_5 \}$ we get the same (upto a permutation of $x_i$) extra condition $c_2=c_3+c_4$, i.e.
$$
c_1=0,\; c_2=c_3+c_4,\; 2c_3=2c_4=2c_5=0.
$$
All such cases have already appeared.\\

Let $\AAN{0} =\{ x_0^2x_1, x_1^2x_0, x_2^2x_1, x_3^2x_0, x_4^2x_0, x_5^2x_0 \}$. Then 
$$
c_1=0,\; 2c_2=2c_3=2c_4=2c_5=0.
$$

For $\AAN{1} =\AAN{0}\cup \{ x_3x_4x_5 \}$ we get an extra condition $c_5=c_3+c_4$, i.e.
$$
c_1=0,\; c_5=c_3+c_4,\; 2c_2=2c_3=2c_4=0.
$$
All such cases have already appeared.\\

For $\AAN{1} =\AAN{0}\cup \{ x_1x_3x_4, x_2x_4x_5 \}$ we get extra conditions $c_4=c_3=c_2+c_5$, i.e.
$$
c_1=0,\; c_3=c_4=c_2+c_5,\; 2c_2=2c_5=0.
$$
All such cases have already appeared.\\

For $\AAN{1} =\AAN{0}\cup \{ x_1x_3x_4, x_2x_3x_4 \}$ we get extra conditions $c_3=c_4$ and $c_2=0$, i.e.
$$
c_1=c_2=0,\; c_3=c_4,\; 2c_4=2c_5=0.
$$
All such cases have already appeared.\\

Let $\AAN{0} =\{ x_0^2x_1, x_1^2x_0, x_2^2x_0, x_3^2x_0, x_4^2x_0, x_5^2x_0, x_2x_3x_4 \}$. Then 
$$
c_1=0,\; c_2=c_3+c_4,\; 2c_3=2c_4=2c_5=0.
$$
All such cases have already appeared.\\

As we explained in the beginning of this section, the fact that all sets $\AA$ we obtained are smooth can be checked by a Macaulay 2 computation. We provide a Macaulay 2 code which verifies this in the Appendix, while the proof of Corollary $1$ shows more specifically that all sets $\AAbar$ we found are smooth. This is sufficient for the proof of Theorem $1$. {\it QED}\\

{\bf Corollary 1.} {\it Diagonalizable abelian automorphism groups of smooth cubic fourfolds are exactly the subgroups of the following groups (see Table $2$ below).

\begin{center}
  \includegraphics[width=1.0\textwidth]{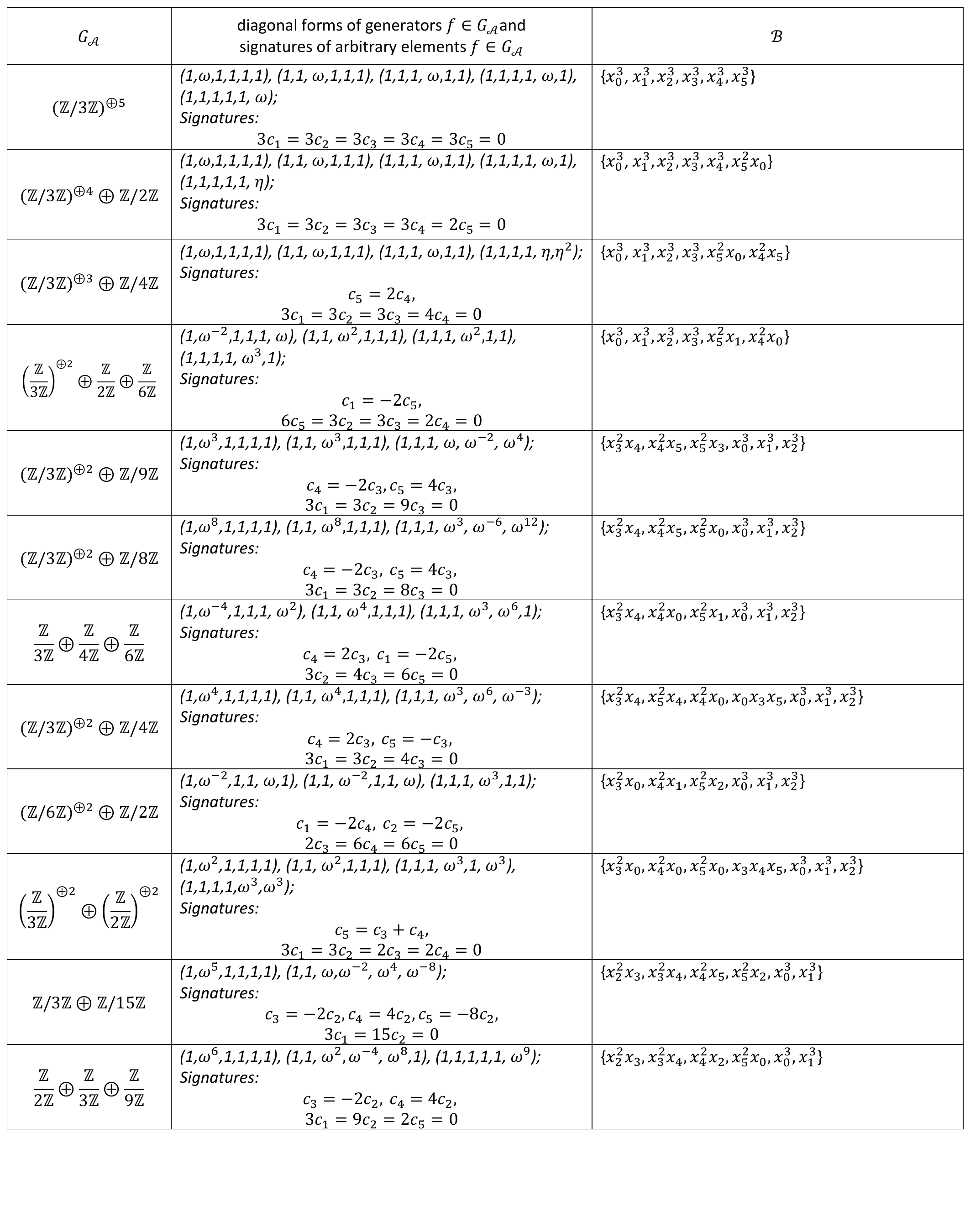}
\end{center}

\begin{center}
  \includegraphics[width=1.0\textwidth]{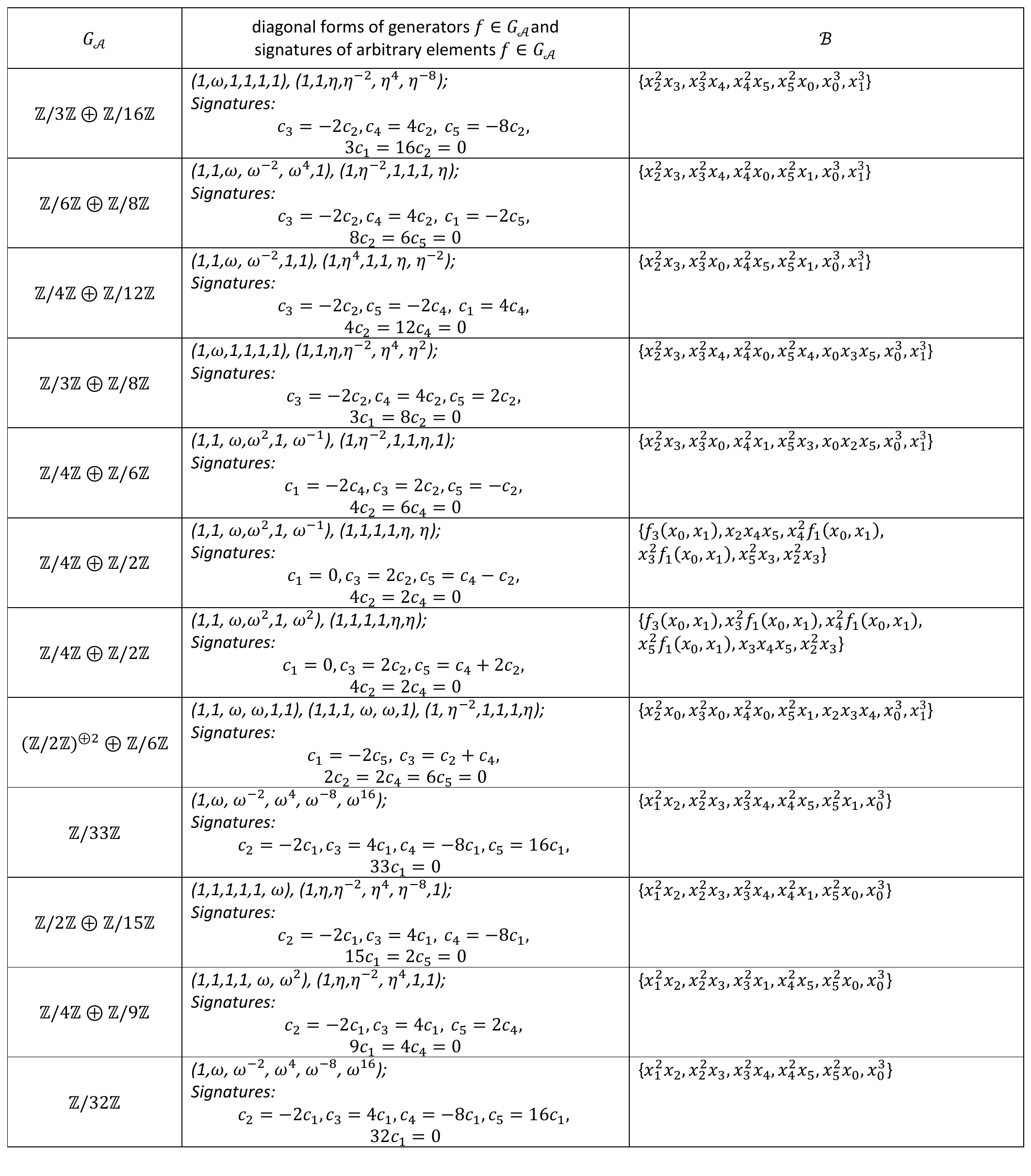}
\end{center}

\begin{center}
  \includegraphics[width=1.0\textwidth]{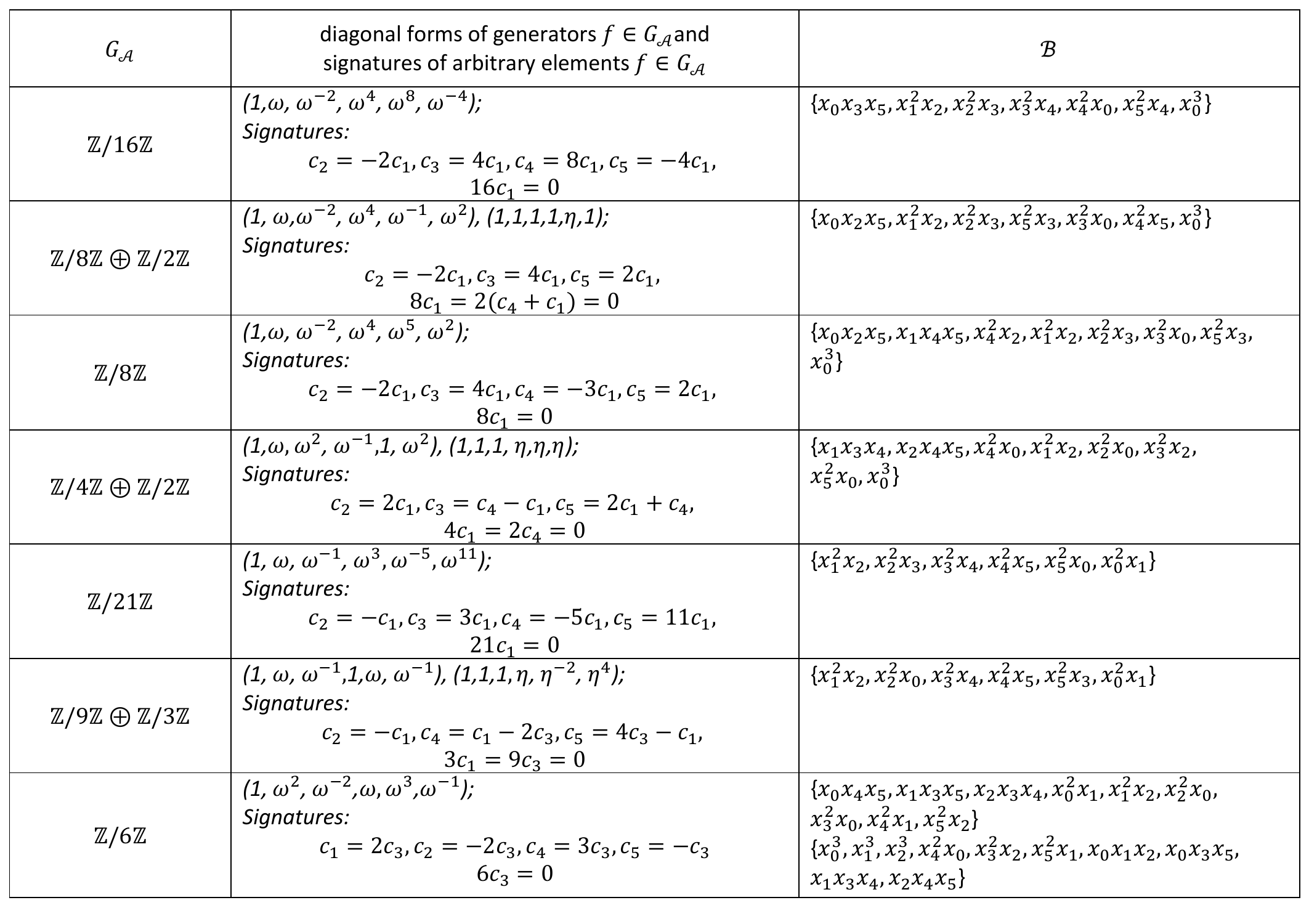}
\end{center}

{\bf Table $2$.} {\sl Diagonalizable abelian automorphism groups of smooth cubic fourfolds.} The first column lists abelian groups $G$ which act effectively on a smooth cubic fourfold. The second column describes generators of $G$ (diagonalized after possibly conjugating by an element of $GL(6)$) as diagonal matrices in $GL(6)$ as well as the conditions on the signatures of arbitrary elements of $G$. The third column gives maximal sets of cubic monomials, which lie in the same eigenspace of each of these generators acting correspondingly on $\text{Sym}^3 ({\mathbb C}^6)$ (the cubic fourfolds with automorphism group $G_{\mathcal A}$ are exactly the ones which in suitable coordinates can be represented as vanishing loci of linear combinations of monomials in one\footnote{In each row of Table $2$ one of the sets $\mathcal B$ is the set $\AAbar$ from Table $1$. There is one more set $\mathcal B$ in the last row. The reason it appears is that the group $\GG$ in the last row is a proper subgroup of the group $\GG$ in the $9$-th row (upto a permutation of $x_i$).} of the sets ${\mathcal B}$.)}\\ 

{\it Proof:} Table $1$ apparently contains some redundancy, because we were not checking strictly minimality of our sets $\AA$. Neither were we checking various possible representations of the same group action. One can reduce this redundancy significantly by noticing that the number of cubes in $\AAbar$ has to decrease along the list in Table $1$. Whenever it increases, this implies that the corresponding group (action) has already appeared earlier in the list (as a subgroup). Let us eliminate such groups $\GG$ from Table $1$.\\

Moreover, let us notice that each of the following groups is a subgroup of another group in Table $1$ (upto a permutation of $x_i$) such that the corresponding sets $\AAbar$ are included one into another:
\begin{itemize}
\item $\GG\cong (\ZZ{3})^{\oplus 2}\oplus \ZZ{2}$ generated by 
$$
(1,1,\om^2,1,1,1),\;(1,1,1,\om^2,1,1),\;(1,1,1,1,\om^3,\om^3), \;\om=\sqr{6}
$$
is a subgroup of the group $\GG\cong (\ZZ{3})^{\oplus 2}\oplus \ZZ{2}\oplus \ZZ{6}$ generated by 
$$
(1,\om^{-2},1,1,1,\om),\;(1,1,\om^2,1,1,1),\;(1,1,1,\om^2,1,1),\;(1,1,1,1,\om^3,1),\; \om=\sqr{6};
$$
\item $\GG\cong \ZZ{3}\oplus \ZZ{4}$ generated by 
$$
(1,1,\om^4,1,1,1),\;(1,1,1,\om^3,\om^6,\om^{6}), \;\om=\sqr{12}
$$
is a subgroup of the group $\GG\cong \ZZ{3}\oplus \ZZ{4}\oplus \ZZ{6}$ generated by 
$$
(1,\om^{-4},1,1,1,\om^2),\;(1,1,\om^4,1,1,1),\;(1,1,1,\om^3,\om^6,1),\;\om=\sqr{12};
$$
\item $\GG\cong \ZZ{3}\oplus \ZZ{4}$ generated by 
$$
(1,1,\om^{4},1,1,1),\;(1,1,1,\om^3,\om^6,\om^{-3}), \; \om=\sqr{12}
$$
is a subgroup of the group $\GG\cong (\ZZ{3})^{\oplus 2}\oplus \ZZ{4}$ generated by
$$
(1,\om^{4},1,1,1,1),\;(1,1,\om^{4},1,1,1),\;(1,1,1,\om^3,\om^6,\om^{-3}),\;\om=\sqr{12};
$$
\item $\GG\cong \ZZ{2}\oplus \ZZ{6}$ generated by 
$$
(1,\om^{-2},1,1,1,\om),\;(1,1,1,\om^3,\om^3,1),\; \om=\sqr{6},
$$
$\GG\cong (\ZZ{2})^{\oplus 2}$ generated by 
$$
(1,1,1,\om^3,1,\om^3),\;(1,1,1,1,\om^3,1),\; \om=\sqr{6},
$$
and $\GG\cong \ZZ{2}$ generated by $(1,1,1,\om^3,\om^3,\om^3)$, $\om=\sqr{6}$ are subgroups (upto a permutation of $x_i$) of the group $\GG\cong (\ZZ{6})^{\oplus 2}\oplus \ZZ{2}$ generated by
$$
(1,\om^{-2},1,1,\om,1),\;(1,1,\om^{-2},1,1,\om),\;(1,1,1,\om^3,1,1),\;\om=\sqr{6};
$$
\item $\GG\cong (\ZZ{2})^{\oplus 2}\oplus \ZZ{3}$ generated by 
$$
(1,1,\om^{2},1,1,1),\;(1,1,1,\om^3,1,\om^3),\;(1,1,1,1,\om^3,\om^3), \; \om=\sqr{6}
$$
is a subgroup of the group $\GG\cong (\ZZ{3})^{\oplus 2}\oplus (\ZZ{2})^{\oplus 2}$ generated by 
$$
(1,\om^{2},1,1,1,1),\;(1,1,\om^{2},1,1,1),\;(1,1,1,\om^3,1,\om^3),\;(1,1,1,1,\om^3,\om^3),\;\om=\sqr{6};
$$
\item $\GG\cong \ZZ{8}$ generated by $(1,1,\om,\om^{-2},\om^4,\om^2)$, $\om=\sqr{8}$ is a subgroup of the group $\GG\cong \ZZ{3}\oplus \ZZ{8}$ generated by 
$$
(1,\eta,1,1,1,1),\;(1,1,\om,\om^{-2},\om^4,\om^2),\; \om=\sqr{8}, \eta=\sqr{3};
$$
\item $\GG\cong \ZZ{2}\oplus \ZZ{4}$ generated by 
$$
(1,1,\om,\om^{2},\om,\om^{2}),\;(1,1,1,1,\eta,1), \;\om=\sqr{4}, \eta=\sqr{2}
$$
and $\GG\cong \ZZ{4}$ generated by $(1,1,\om,\om^{2},\om^{-1},\om^{2})$, $\om=\sqr{4}$ are subgroups (upto a permutation of $x_i$) of the group $\GG\cong \ZZ{8}\oplus \ZZ{2}$ generated by
$$
(1,\om,\om^{-2},\om^4,\om^{-1},\om^2),\;(1,1,1,1,\eta,1),\;\om=\sqr{8}, \eta=\sqr{2};
$$ 
\item $\GG\cong (\ZZ{2})^{\oplus 2}$ generated by 
$$
(1,1,\om,\om,\om,1),\;(1,1,\om,\om,1,\om),\; \om=\sqr{2}
$$
and $\GG\cong (\ZZ{2})^{\oplus 3}$ generated by
$$
(1,1,\om,1,1,1),\;(1,1,1,\om,1,\om),\;(1,1,1,\om,\om,1),\;\om=\sqr{2}
$$
are subgroups (upto a permutation of $x_i$) of the group $\GG\cong (\ZZ{2})^{\oplus 2}\oplus \ZZ{6}$ generated by
$$
(1,1,\om,\om,1,1),\; (1,1,1,\om,\om,1),\; (1,\eta^{-2},1,1,1,\eta), \;\om=\sqr{2}, \eta=\sqr{6};
$$
\item $\GG\cong \ZZ{8}$ generated by $(1,1,\om,\om^{-2},\om^4,\om^4)$, $\om=\sqr{8}$ is a subgroup of the group $\GG\cong \ZZ{6}\oplus \ZZ{8}$ generated by
$$
(1,1,\om,\om^{-2},\om^4,1),\; (1,\eta^{-2},1,1,1,\eta), \;\om=\sqr{8}, \eta=\sqr{6};
$$
\item $\GG\cong \ZZ{4}\oplus \ZZ{2}$ generated by 
$$
(1,1,\om,\om^{2},1,\om^{-1}),\; (1,1,1,1,\eta,1), \;\om=\sqr{4}, \eta=\sqr{2}
$$
is a subgroup of the group $\GG\cong \ZZ{4}\oplus \ZZ{6}$ generated by
$$
(1,1,\om,\om^{2},1,\om^{-1}),\; (1,\eta^{-2},1,1,\eta,1), \;\om=\sqr{4}, \eta=\sqr{6};
$$
\item $\GG\cong \ZZ{3}$ generated by $(1,\om^2,\om,1,\om^2,\om)$, $\om=\sqr{3}$ is a subgroup of the group $\GG\cong \ZZ{9}\oplus \ZZ{3}$ generated by
$$
(1,\om,\om^{-1},1,\om,\om^{-1}),\; (1,1,1,\eta,\eta^{-2},\eta^{4}), \;\om=\sqr{3}, \eta=\sqr{9}.
$$
\end{itemize}

After eliminating these groups from Table $1$ we obtain Table $2$.\\

One can apply a Macaulay 2 computation described in the beginning of this section and check that all sets $\AAbar$ listed\footnote{In the third column of Table $2$ one of the sets $\mathcal B$ is $\AAbar$.} in Table $2$ are smooth (which implies that all sets $\AAbar$ listed in Table $1$ are smooth as well). Indeed, one can take the sum of all monomials with coefficient $1$ and check that the resulting cubic form gives a smooth cubic fourfold in all cases except for the following ones:
\begin{itemize}
\item for the entry number $18$ one can use $F=x_0^3+x_1^3+x_2^2x_3+x_3^2x_0+x_4^2x_1+x_5^2x_3+x_2x_4x_5$,
\item for the entry number $19$ one can use $F=x_0^3+x_1^3+x_2^2x_3+x_3^2x_1+x_4^2x_0-x_4^2x_1+x_5^2x_0+x_3x_4x_5$.
\end{itemize}

Finally, let us compute the third column of Table $2$, i.e. describe all smooth cubic fourfolds which admit each of the group actions listed in Table $2$.\\

Let $\GG$ be one of the groups listed in the table. We have to find all smooth sets $\mathcal B\subset \mathcal M$, whose elements lie in the same eigenspace of each of the generators of $\GG$. By Lemma $1$ we know that any such set $\mathcal B$ for any $p$ contains $x_p^2x_q$ for some $q$. In particular, $x_0^2x_j\in \mathcal B$ for some $j\in \{ 0,1,2,3,4,5 \}$. Let us fix such $j$.\\

We start with the last three rows of Table $2$. In each of these cases $x_0^2x_1\in \AAbar$. Hence if $x_0^2x_1\in \mathcal B$, then $\mathcal B = \AAbar$. So, we may assume that $j\in \{ 0,2,3,4,5 \}$.\\

Let $\GG \cong \ZZ{21}$ be generated by 
$$
f=(1,\om,{\om}^{-1},{\om}^3,{\om}^{-5},{\om}^{11}), \om=\sqr{21}.
$$

Note that 
\begin{itemize}
\item $wt_f(x_0^3)=0\neq wt_f(x_1^2x_k)$ for any $k$,
\item $wt_f(x_0^2x_2)=-1\neq wt_f(x_1^2x_k)$ for any $k$,
\item $wt_f(x_0^2x_3)=3\neq wt_f(x_2^2x_k)$ for any $k$,
\item $wt_f(x_0^2x_4)=-5\neq wt_f(x_1^2x_k)$ for any $k$,
\item $wt_f(x_0^2x_5)=11\neq wt_f(x_1^2x_k)$ for any $k$.
\end{itemize}

This means that either $x_1^2x_k\notin \BB$ for any $k$ or $x_2^2x_k\notin \BB$ for any $k$, i.e. $\BB$ is not smooth by Lemma $1$.\\

Hence $\BB=\AAbar$ is the only possibility for this group $\GG$.\\

Let $\GG \cong \ZZ{9}\oplus \ZZ{3}$ be generated by 
$$
f_1=(1,\om,{\om}^{-1},1,{\om},{\om}^{-1}),\; f_2=(1,1,1,\eta,\eta ^{-2},\eta ^4),\; \om=\sqr{3}, \eta=\sqr{9}.
$$

Note that 
\begin{itemize}
\item $wt_{f_1}(x_0^3)=0\neq wt_{f_1}(x_3^2x_k)$ for any $k\neq 0,3$,
\item $wt_{f_2}(x_0^3)=0\neq wt_{f_2}(x_3^2x_0)=2$,
\item $wt_{f_2}(x_0^3)=0\neq wt_{f_2}(x_3^3)=3$,
\item $wt_{f_1}(x_0^2x_2)=-1\neq wt_{f_1}(x_3^2x_k)$ for any $k\neq 2,5$,
\item $wt_{f_2}(x_0^2x_2)=0\neq wt_{f_2}(x_3^2x_2)=2$,
\item $wt_{f_2}(x_0^2x_2)=0\neq wt_{f_2}(x_3^2x_5)=6$,
\item $wt_{f_1}(x_0^2x_3)=0\neq wt_{f_1}(x_3^2x_k)$ for any $k\neq 0,3$,
\item $wt_{f_2}(x_0^2x_3)=1\neq wt_{f_2}(x_3^2x_0)=2$,
\item $wt_{f_2}(x_0^2x_3)=1\neq wt_{f_2}(x_3^3)=3$,
\item $wt_{f_1}(x_0^2x_4)=1\neq wt_{f_1}(x_3^2x_k)$ for any $k\neq 1,4$,
\item $wt_{f_2}(x_0^2x_4)=-2\neq wt_{f_2}(x_3^2x_1)=2$,
\item $wt_{f_2}(x_0^2x_4)=-2\neq wt_{f_2}(x_3^2x_4)=0$,
\item $wt_{f_1}(x_0^2x_5)=-1\neq wt_{f_1}(x_3^2x_k)$ for any $k\neq 2,5$,
\item $wt_{f_2}(x_0^2x_5)=4\neq wt_{f_2}(x_3^2x_2)=2$,
\item $wt_{f_2}(x_0^2x_5)=4\neq wt_{f_2}(x_3^2x_5)=6$.
\end{itemize}

This means that $x_3^2x_k\notin \BB$ for any $k$, i.e. $\BB$ is not smooth by Lemma $1$.\\

Hence $\BB=\AAbar$ is the only possibility for this group $\GG$.\\

Let $\GG \cong \ZZ{6}$ be generated by 
$$
f=(1,\om ^2,{\om}^{-2},\om,{\om}^3,{\om}^{-1}),\; \om=\sqr{6}.
$$

If $x_0^3\in \BB$, then $\BB$ consists of monomials $m=x_0^{i_0}x_1^{i_1}x_2^{i_2}x_3^{i_3}x_4^{i_4}x_5^{i_5}$ such that $wt_f(m)=wt_f(x_0^3)=0$, i.e.
$$
2i_1-2i_2+i_3+3i_4-i_5\equiv 0 \; mod \; 6.
$$

This means that $\BB=\{ x_0^3$, $x_1^3$, $x_2^3$, $x_4^2x_0$, $x_3^2x_2$, $x_5^2x_1$, $x_0x_1x_2$, $x_0x_3x_5$, $x_1x_3x_4$, $x_2x_4x_5 \}$. This set is smooth.\\

If $x_0^2x_2\in \BB$, then $\BB$ consists of monomials $m=x_0^{i_0}x_1^{i_1}x_2^{i_2}x_3^{i_3}x_4^{i_4}x_5^{i_5}$ such that $wt_f(m)=wt_f(x_0^2x_2)=-2$, i.e.
$$
2i_1-2i_2+i_3+3i_4-i_5\equiv -2 \; mod \; 6.
$$

This means that $\BB=\{ x_4^2x_2, x_0^2x_2, x_5^2x_0, x_1^2x_0, x_2^2x_1, x_3^2x_1, x_1x_4x_5, x_0x_3x_4, x_2x_3x_5 \}$. This set coincides with $\AAbar$ upto a permutation of $x_i$ ($x_1\leftrightarrow x_2$, $x_3\leftrightarrow x_5$) and using $f^{-1}$ as a generator of $\GG$.\\

If $x_0^2x_3\in \BB$, then $\BB$ consists of monomials $m=x_0^{i_0}x_1^{i_1}x_2^{i_2}x_3^{i_3}x_4^{i_4}x_5^{i_5}$ such that $wt_f(m)=wt_f(x_0^2x_3)=1$, i.e.
$$
2i_1-2i_2+i_3+3i_4-i_5\equiv 1 \; mod \; 6.
$$

This means that $\BB=\{ x_4^2x_3, x_0^2x_3, x_5^2x_4, x_1^2x_4, x_2^2x_5, x_3^2x_5, x_1x_2x_3, x_0x_1x_5, x_0x_2x_4 \}$. This set coincides with $\AAbar$ upto a permutation of $x_i$ ($x_0\leftrightarrow x_3$, $x_1\leftrightarrow x_5$, $x_2\leftrightarrow x_4$) and using $({\om}^{-1}f)^{-1}$ as a generator of $\GG$.\\

If $x_0^2x_4\in \BB$, then $\BB$ consists of monomials $m=x_0^{i_0}x_1^{i_1}x_2^{i_2}x_3^{i_3}x_4^{i_4}x_5^{i_5}$ such that $wt_f(m)=wt_f(x_0^2x_4)=3$, i.e.
$$
2i_1-2i_2+i_3+3i_4-i_5\equiv 3 \; mod \; 6.
$$

This means that $\BB=\{ x_3^3, x_4^3, x_5^3, x_0^2x_4, x_2^2x_3, x_1^2x_5, x_0x_1x_3, x_0x_2x_5, x_1x_2x_4, x_3x_4x_5 \}$. This set coincides with an earlier found set $\BB$ upto a permutation of $x_i$ ($x_0\leftrightarrow x_4$, $x_2\leftrightarrow x_3$, $x_1\leftrightarrow x_5$) and using ${\om}^{3}f$ as a generator of $\GG$.\\

If $x_0^2x_5\in \BB$, then $\BB$ consists of monomials $m=x_0^{i_0}x_1^{i_1}x_2^{i_2}x_3^{i_3}x_4^{i_4}x_5^{i_5}$ such that $wt_f(m)=wt_f(x_0^2x_5)=-1$, i.e.
$$
2i_1-2i_2+i_3+3i_4-i_5\equiv -1 \; mod \; 6.
$$

This means that $\BB=\{ x_5^2x_3, x_0^2x_5, x_4^2x_5, x_3^2x_4, x_2^2x_4, x_1^2x_3, x_1x_2x_5, x_0x_1x_4, x_0x_2x_3 \}$. This set coincides with $\AAbar$ upto a permutation of $x_i$ ($x_0\leftrightarrow x_4$, $x_2\leftrightarrow x_3$, $x_1\leftrightarrow x_5$) and using ${\om}^{3}f$ as a generator of $\GG$.\\

We conclude that apart from $\AAbar$ only the following smooth set $\BB$ is possible (upto a permutation of $x_i$) for this group $\GG$:
\begin{itemize}
\item $\BB=\{ x_0^3, x_1^3, x_2^3, x_4^2x_0, x_3^2x_2, x_5^2x_1, x_0x_1x_2, x_0x_3x_5, x_1x_3x_4, x_2x_4x_5 \}$.
\end{itemize}

Now let us consider the other rows of Table $2$. Since $x_0^3\in\AAbar$ in each case, we may assume without loss of generality that $x_0^2x_j\in \BB$ for some $j\in \{ 1,2,3,4,5 \}$.\\

Let $\GG \cong (\ZZ{3})^{\oplus 5}$ be generated by 
\begin{gather*}
f_1=(1,\om,1,1,1,1),\; f_2=(1,1,\om,1,1,1),\; f_3=(1,1,1,\om,1,1),\\
f_4=(1,1,1,1,\om,1),\; f_5=(1,1,1,1,1,\om),\; \om=\sqr{3}.
\end{gather*}

Note that $wt_{f_j}(x_0^2x_j)=1\neq wt_{f_j}(x_j^2x_k)$ for any $k$. This means that $x_j^2x_k\notin \BB$ for any $k$, i.e. $\BB$ is not smooth.\\

Hence $\BB=\AAbar$ is the only possibility for this group $\GG$.\\

Let $\GG \cong (\ZZ{3})^{\oplus 4}\oplus \ZZ{2}$ be generated by 
\begin{gather*}
f_1=(1,\om,1,1,1,1),\; f_2=(1,1,\om,1,1,1),\; f_3=(1,1,1,\om,1,1),\\
f_4=(1,1,1,1,\om,1),\; f_5=(1,1,1,1,1,\eta),\; \om=\sqr{3},\eta=\sqr{2}.
\end{gather*}

Note that 
\begin{itemize}
\item $wt_{f_j}(x_0^2x_j)=1\neq wt_{f_j}(x_j^2x_k)$ for any $j\in \{ 1,2,3,4 \}$ for any $k$,
\item $wt_{f_1}(x_0^2x_5)=0\neq wt_{f_1}(x_1^2x_k)$ for any $k\neq 1$,
\item $wt_{f_5}(x_0^2x_5)=1\neq wt_{f_5}(x_1^3)=0$.
\end{itemize}

This means that either $x_1^2x_k\notin \BB$ for any $k$ or $x_j^2x_k\notin \BB$ for any $k$, i.e. $\BB$ is not smooth by Lemma $1$.\\

Hence $\BB=\AAbar$ is the only possibility for this group $\GG$.\\

Let $\GG \cong (\ZZ{3})^{\oplus 3}\oplus \ZZ{4}$ be generated by 
\begin{gather*}
f_1=(1,\om,1,1,1,1),\; f_2=(1,1,\om,1,1,1),\; f_3=(1,1,1,\om,1,1),\\
f_4=(1,1,1,1,\eta,\eta ^2),\; \om=\sqr{3},\eta=\sqr{4}.
\end{gather*}

Note that 
\begin{itemize}
\item $wt_{f_j}(x_0^2x_j)=1\neq wt_{f_j}(x_j^2x_k)$ for any $j\in \{ 1,2,3, 4 \}$ for any $k$,
\item $wt_{f_4}(x_0^2x_5)=2\neq wt_{f_4}(x_1^2x_k)$ for any $k\neq 5$,
\item $wt_{f_1}(x_0^2x_5)=0\neq wt_{f_1}(x_1^2x_5)=2$.
\end{itemize}

This means that either $x_1^2x_k\notin \BB$ for any $k$ or $x_j^2x_k\notin \BB$ for any $k$, i.e. $\BB$ is not smooth by Lemma $1$.\\

Hence $\BB=\AAbar$ is the only possibility for this group $\GG$.\\

Let $\GG \cong (\ZZ{3})^{\oplus 2}\oplus \ZZ{2}\oplus \ZZ{6}$ be generated by 
\begin{gather*}
f_1=(1,\om^{-2},1,1,1,\om),\; f_2=(1,1,\om^2,1,1,1),\; f_3=(1,1,1,\om^2,1,1),\\
f_4=(1,1,1,1,\om^3,1),\; \om=\sqr{6}.
\end{gather*}

Note that (we compute weights as elements of $\ZZ{6}$ here)
\begin{itemize}
\item $wt_{f_j}(x_0^2x_j)=2\neq wt_{f_j}(x_j^2x_k)$ for any $j\in \{ 2,3 \}$ for any $k$,
\item $wt_{f_2}(x_0^2x_j)=0\neq wt_{f_2}(x_2^2x_k)$ for any $j\in \{ 1,4,5 \}$ for any $k\neq 2$,
\item $wt_{f_1}(x_0^2x_1)=-2\neq wt_{f_1}(x_2^3)=0$,
\item $wt_{f_4}(x_0^2x_4)=3\neq wt_{f_4}(x_2^3)=0$,
\item $wt_{f_1}(x_0^2x_5)=1\neq wt_{f_1}(x_2^3)=0$.
\end{itemize}

This means that either $x_2^2x_k\notin \BB$ for any $k$ or $x_j^2x_k\notin \BB$ for any $k$, i.e. $\BB$ is not smooth by Lemma $1$.\\

Hence $\BB=\AAbar$ is the only possibility for this group $\GG$.\\

Let $\GG \cong (\ZZ{3})^{\oplus 2}\oplus \ZZ{9}$ be generated by 
\begin{gather*}
f_1=(1,\om^{3},1,1,1,1),\; f_2=(1,1,\om^3,1,1,1),\; f_3=(1,1,1,\om ,\om ^{-2},\om ^4),\; \om=\sqr{9}.
\end{gather*}

Note that (we compute weights as elements of $\ZZ{9}$ here)
\begin{itemize}
\item $wt_{f_j}(x_0^2x_j)=3\neq wt_{f_j}(x_j^2x_k)$ for any $j\in \{ 1,2 \}$ for any $k$,
\item $wt_{f_3}(x_0^2x_3)=1\neq wt_{f_3}(x_3^2x_k)$ for any $k$,
\item $wt_{f_3}(x_0^2x_4)=-2\neq wt_{f_3}(x_4^2x_k)$ for any $k$,
\item $wt_{f_3}(x_0^2x_5)=4\neq wt_{f_3}(x_5^2x_k)$ for any $k$.
\end{itemize}

This means that $x_j^2x_k\notin \BB$ for any $k$, i.e. $\BB$ is not smooth by Lemma $1$.\\

Hence $\BB=\AAbar$ is the only possibility for this group $\GG$.\\

Let $\GG \cong (\ZZ{3})^{\oplus 2}\oplus \ZZ{8}$ be generated by 
\begin{gather*}
f_1=(1,\om^{8},1,1,1,1),\; f_2=(1,1,\om^8,1,1,1),\; f_3=(1,1,1,\om ^3,\om ^{-6},\om ^{12}),\; \om=\sqr{24}.
\end{gather*}

Note that (we compute weights as elements of $\ZZ{24}$ here)
\begin{itemize}
\item $wt_{f_j}(x_0^2x_j)=8\neq wt_{f_j}(x_j^2x_k)$ for any $j\in \{ 1,2 \}$ for any $k$,
\item $wt_{f_3}(x_0^2x_3)=3\neq wt_{f_3}(x_3^2x_k)$ for any $k$,
\item $wt_{f_3}(x_0^2x_4)=-6\neq wt_{f_3}(x_4^2x_k)$ for any $k$,
\item $wt_{f_1}(x_0^2x_5)=0\neq wt_{f_1}(x_1^2x_k)$ for any $k\neq 1$,
\item $wt_{f_3}(x_0^2x_5)=12\neq wt_{f_3}(x_1^3)=0$.
\end{itemize}

This means that either $x_j^2x_k\notin \BB$ for any $k$ or $x_1^2x_k\notin \BB$ for any $k$, i.e. $\BB$ is not smooth by Lemma $1$.\\

Hence $\BB=\AAbar$ is the only possibility for this group $\GG$.\\

Let $\GG \cong \ZZ{3} \oplus \ZZ{4}\oplus \ZZ{6}$ be generated by 
\begin{gather*}
f_1=(1,\om^{-4},1,1,1,\om^2),\; f_2=(1,1,\om^4,1,1,1),\; f_3=(1,1,1,\om^3,\om ^6,1),\; \om=\sqr{12}.
\end{gather*}

Note that (we compute weights as elements of $\ZZ{12}$ here)
\begin{itemize}
\item $wt_{f_1}(x_0^2x_1)=-4\neq wt_{f_1}(x_1^2x_k)$ for any $k$,
\item $wt_{f_2}(x_0^2x_2)=4\neq wt_{f_2}(x_2^2x_k)$ for any $k$,
\item $wt_{f_3}(x_0^2x_3)=3\neq wt_{f_3}(x_3^2x_k)$ for any $k$,
\item $wt_{f_1}(x_0^2x_5)=2\neq wt_{f_1}(x_5^2x_k)$ for any $k$,
\item $wt_{f_2}(x_0^2x_4)=0\neq wt_{f_2}(x_2^2x_k)$ for any $k\neq 2$,
\item $wt_{f_3}(x_0^2x_4)=6\neq wt_{f_3}(x_2^3)=0$.
\end{itemize}

This means that either $x_j^2x_k\notin \BB$ for any $k$ or $x_2^2x_k\notin \BB$ for any $k$, i.e. $\BB$ is not smooth by Lemma $1$.\\

Hence $\BB=\AAbar$ is the only possibility for this group $\GG$.\\

Let $\GG \cong (\ZZ{3})^{\oplus 2}\oplus \ZZ{4}$ be generated by 
\begin{gather*}
f_1=(1,\om^{4},1,1,1,1),\; f_2=(1,1,\om^4,1,1,1),\; f_3=(1,1,1,\om ^3,\om ^{6},\om ^{-3}),\; \om=\sqr{12}.
\end{gather*}

Note that (we compute weights as elements of $\ZZ{12}$ here)
\begin{itemize}
\item $wt_{f_1}(x_0^2x_1)=4\neq wt_{f_1}(x_1^2x_k)$ for any $k$,
\item $wt_{f_2}(x_0^2x_2)=4\neq wt_{f_2}(x_2^2x_k)$ for any $k$,
\item $wt_{f_1}(x_0^2x_j)=0\neq wt_{f_1}(x_1^2x_k)$ for any $k \neq 1$ for any $j\in \{ 3,4,5 \}$,
\item $wt_{f_3}(x_0^2x_j)\neq wt_{f_3}(x_1^3)=0$ for any $j\in \{ 3,4,5 \}$.
\end{itemize}

This means that either $x_j^2x_k\notin \BB$ for any $k$ or $x_1^2x_k\notin \BB$ for any $k$, i.e. $\BB$ is not smooth by Lemma $1$.\\

Hence $\BB=\AAbar$ is the only possibility for this group $\GG$.\\

Let $\GG \cong (\ZZ{6})^{\oplus 2}\oplus \ZZ{2}$ be generated by 
\begin{gather*}
f_1=(1,\om^{-2},1,1,\om,1),\; f_2=(1,1,\om^{-2},1,1,\om),\; f_3=(1,1,1,\om ^3,1,1),\; \om=\sqr{6}.
\end{gather*}

Note that (we compute weights as elements of $\ZZ{6}$ here)
\begin{itemize}
\item $wt_{f_1}(x_0^2x_4)=1\neq wt_{f_1}(x_4^2x_k)$ for any $k$,
\item $wt_{f_2}(x_0^2x_5)=1\neq wt_{f_2}(x_5^2x_k)$ for any $k$,
\item $wt_{f_1}(x_0^2x_1)=-2\neq wt_{f_1}(x_1^2x_k)$ for any $k$,
\item $wt_{f_2}(x_0^2x_2)=-2\neq wt_{f_2}(x_2^2x_k)$ for any $k$,
\item $wt_{f_1}(x_0^2x_3)=0\neq wt_{f_1}(x_1^2x_k)$ for any $k\neq 1$,
\item $wt_{f_3}(x_0^2x_3)=3\neq wt_{f_3}(x_1^3)=0$.
\end{itemize}

This means that either $x_j^2x_k\notin \BB$ for any $k$ or $x_1^2x_k\notin \BB$ for any $k$, i.e. $\BB$ is not smooth by Lemma $1$.\\

Hence $\BB=\AAbar$ is the only possibility for this group $\GG$.\\

Let $\GG \cong (\ZZ{3})^{\oplus 2}\oplus (\ZZ{2})^{\oplus 2}$ be generated by 
\begin{gather*}
f_1=(1,\om^{2},1,1,1,1),\; f_2=(1,1,\om^{2},1,1,1),\; f_3=(1,1,1,\om ^3,1,\om^3),\\
f_4=(1,1,1,1,\om ^3,\om^3),\; \om=\sqr{6}.
\end{gather*}

Note that (we compute weights as elements of $\ZZ{6}$ here)
\begin{itemize}
\item $wt_{f_1}(x_0^2x_1)=2\neq wt_{f_1}(x_1^2x_k)$ for any $k$,
\item $wt_{f_2}(x_0^2x_2)=2\neq wt_{f_2}(x_2^2x_k)$ for any $k$,
\item $wt_{f_1}(x_0^2x_j)=0\neq wt_{f_1}(x_1^2x_k)$ for any $k\neq 1$ for any $j\in \{ 3,4 \}$,
\item $wt_{f_j}(x_0^2x_j)=3\neq wt_{f_j}(x_1^3)$ for any $j\in \{ 3,4 \}$,
\item $wt_{f_1}(x_0^2x_5)=0\neq wt_{f_1}(x_1^2x_k)$ for any $k\neq 1$,
\item $wt_{f_3}(x_0^2x_5)=3\neq wt_{f_3}(x_1^3)=0$.
\end{itemize}

This means that either $x_j^2x_k\notin \BB$ for any $k$ or $x_1^2x_k\notin \BB$ for any $k$, i.e. $\BB$ is not smooth by Lemma $1$.\\

Hence $\BB=\AAbar$ is the only possibility for this group $\GG$.\\

Let $\GG \cong \ZZ{3}\oplus \ZZ{15}$ be generated by 
\begin{gather*}
f_1=(1,\om^5,1,1,1,1),\; f_2=(1,1,\om,\om ^{-2},\om ^4,\om ^{-8}),\; \om=\sqr{15}.
\end{gather*}

Note that (we compute weights as elements of $\ZZ{15}$ here)
\begin{itemize}
\item $wt_{f_1}(x_0^2x_1)=5\neq wt_{f_1}(x_1^2x_k)$ for any $k$,
\item $wt_{f_2}(x_0^2x_2)=1\neq wt_{f_2}(x_2^2x_k)$ for any $k$,
\item $wt_{f_2}(x_0^2x_3)=-2\neq wt_{f_2}(x_3^2x_k)$ for any $k$,
\item $wt_{f_2}(x_0^2x_4)=4\neq wt_{f_2}(x_4^2x_k)$ for any $k$,
\item $wt_{f_2}(x_0^2x_5)=-8\neq wt_{f_2}(x_5^2x_k)$ for any $k$.
\end{itemize}

This means that $x_j^2x_k\notin \BB$ for any $k$, i.e. $\BB$ is not smooth by Lemma $1$.\\

Hence $\BB=\AAbar$ is the only possibility for this group $\GG$.\\

Let $\GG \cong \ZZ{2} \oplus \ZZ{3}\oplus \ZZ{9}$ be generated by 
\begin{gather*}
f_1=(1,\om^{6},1,1,1,1),\; f_2=(1,1,\om^2,\om ^{-4},\om ^8,1),\; f_3=(1,1,1,1,1,\om^9),\; \om=\sqr{18}.
\end{gather*}

Note that (we compute weights as elements of $\ZZ{18}$ here)
\begin{itemize}
\item $wt_{f_1}(x_0^2x_1)=6\neq wt_{f_1}(x_1^2x_k)$ for any $k$,
\item $wt_{f_2}(x_0^2x_2)=2\neq wt_{f_2}(x_2^2x_k)$ for any $k$,
\item $wt_{f_2}(x_0^2x_3)=-4\neq wt_{f_2}(x_3^2x_k)$ for any $k$,
\item $wt_{f_2}(x_0^2x_4)=8\neq wt_{f_2}(x_4^2x_k)$ for any $k$,
\item $wt_{f_1}(x_0^2x_5)=0\neq wt_{f_1}(x_1^2x_k)$ for any $k\neq 1$,
\item $wt_{f_3}(x_0^2x_5)=9\neq wt_{f_3}(x_1^3)=0$.
\end{itemize}

This means that either $x_j^2x_k\notin \BB$ for any $k$ or $x_1^2x_k\notin \BB$ for any $k$, i.e. $\BB$ is not smooth by Lemma $1$.\\

Hence $\BB=\AAbar$ is the only possibility for this group $\GG$.\\

Let $\GG \cong \ZZ{3}\oplus \ZZ{16}$ be generated by 
\begin{gather*}
f_1=(1,\om,1,1,1,1),\; f_2=(1,1,\eta,\eta^{-2},\eta ^{4},\eta ^{-8}),\; \om=\sqr{3}, \eta=\sqr{16}.
\end{gather*}

Note that
\begin{itemize}
\item $wt_{f_1}(x_0^2x_1)=1\neq wt_{f_1}(x_1^2x_k)$ for any $k$,
\item $wt_{f_1}(x_0^2x_j)=0\neq wt_{f_1}(x_1^2x_k)$ for any $k\neq 1$ for any $j\in \{ 2,3,4,5 \}$,
\item $wt_{f_2}(x_0^2x_j)\neq wt_{f_2}(x_1^3)=0$ for any $j\in \{ 2,3,4,5 \}$.
\end{itemize}

This means that $x_1^2x_k\notin \BB$ for any $k$, i.e. $\BB$ is not smooth by Lemma $1$.\\

Hence $\BB=\AAbar$ is the only possibility for this group $\GG$.\\

Let $\GG \cong \ZZ{6}\oplus \ZZ{8}$ be generated by 
\begin{gather*}
f_1=(1,1,\om,\om^{-2},\om ^4,1),\; f_2=(1,\eta^{-2},1,1,1,\eta),\; \om=\sqr{8}, \eta=\sqr{6}.
\end{gather*}

Note that
\begin{itemize}
\item $wt_{f_2}(x_0^2x_1)=-2\neq wt_{f_2}(x_1^2x_k)$ for any $k$,
\item $wt_{f_2}(x_0^2x_5)=1\neq wt_{f_2}(x_5^2x_k)$ for any $k$,
\item $wt_{f_1}(x_0^2x_2)=1\neq wt_{f_1}(x_2^2x_k)$ for any $k$,
\item $wt_{f_1}(x_0^2x_3)=-2\neq wt_{f_1}(x_3^2x_k)$ for any $k$,
\item $wt_{f_2}(x_0^2x_4)=0\neq wt_{f_2}(x_1^2x_k)$ for any $k\neq 1$,
\item $wt_{f_1}(x_0^2x_4)=4\neq wt_{f_1}(x_1^3)=0$.
\end{itemize}

This means that either $x_j^2x_k\notin \BB$ for any $k$ or $x_1^2x_k\notin \BB$ for any $k$, i.e. $\BB$ is not smooth by Lemma $1$.\\

Hence $\BB=\AAbar$ is the only possibility for this group $\GG$.\\

Let $\GG \cong \ZZ{4}\oplus \ZZ{12}$ be generated by 
\begin{gather*}
f_1=(1,1,\om,\om^{-2},1,1),\; f_2=(1,\eta^{4},1,1,\eta,\eta^{-2}),\; \om=\sqr{4}, \eta=\sqr{12}.
\end{gather*}

Note that
\begin{itemize}
\item $wt_{f_2}(x_0^2x_1)=4\neq wt_{f_2}(x_1^2x_k)$ for any $k$,
\item $wt_{f_2}(x_0^2x_4)=1\neq wt_{f_2}(x_4^2x_k)$ for any $k$,
\item $wt_{f_2}(x_0^2x_5)=-2\neq wt_{f_2}(x_5^2x_k)$ for any $k$,
\item $wt_{f_1}(x_0^2x_2)=1\neq wt_{f_1}(x_2^2x_k)$ for any $k$,
\item $wt_{f_2}(x_0^2x_3)=0\neq wt_{f_2}(x_1^2x_k)$ for any $k\neq 1$,
\item $wt_{f_1}(x_0^2x_3)=-2\neq wt_{f_1}(x_1^3)=0$.
\end{itemize}

This means that either $x_j^2x_k\notin \BB$ for any $k$ or $x_1^2x_k\notin \BB$ for any $k$, i.e. $\BB$ is not smooth by Lemma $1$.\\

Hence $\BB=\AAbar$ is the only possibility for this group $\GG$.\\

Let $\GG \cong \ZZ{3}\oplus \ZZ{8}$ be generated by 
\begin{gather*}
f_1=(1,\om,1,1,1,1),\; f_2=(1,1,\eta,\eta^{-2},\eta^4,\eta^{2}),\; \om=\sqr{3}, \eta=\sqr{8}.
\end{gather*}

Note that
\begin{itemize}
\item $wt_{f_1}(x_0^2x_1)=1\neq wt_{f_1}(x_1^2x_k)$ for any $k$,
\item $wt_{f_2}(x_0^2x_2)=1\neq wt_{f_2}(x_2^2x_k)$ for any $k$,
\item $wt_{f_1}(x_0^2x_j)=0\neq wt_{f_1}(x_1^2x_k)$ for any $k\neq 1$ for any $j\in \{ 3,4,5 \}$,
\item $wt_{f_2}(x_0^2x_j)\neq wt_{f_2}(x_1^3)=0$ for any $j\in \{ 3,4,5 \}$.
\end{itemize}

This means that either $x_j^2x_k\notin \BB$ for any $k$ or $x_1^2x_k\notin \BB$ for any $k$, i.e. $\BB$ is not smooth by Lemma $1$.\\

Hence $\BB=\AAbar$ is the only possibility for this group $\GG$.\\

Let $\GG \cong \ZZ{4}\oplus \ZZ{6}$ be generated by 
\begin{gather*}
f_1=(1,1,\om,\om ^2,1,\om ^{-1}),\; f_2=(1,\eta^{-2},1,1,\eta,1),\; \om=\sqr{4}, \eta=\sqr{6}.
\end{gather*}

Note that
\begin{itemize}
\item $wt_{f_2}(x_0^2x_1)=-2\neq wt_{f_2}(x_1^2x_k)$ for any $k$,
\item $wt_{f_2}(x_0^2x_4)=1\neq wt_{f_2}(x_4^2x_k)$ for any $k$,
\item $wt_{f_2}(x_0^2x_j)=0\neq wt_{f_2}(x_1^2x_k)$ for any $k\neq 1$ for any $j\in \{ 2,3,5 \}$,
\item $wt_{f_1}(x_0^2x_j)\neq wt_{f_1}(x_1^3)=0$ for any $j\in \{ 2,3,5 \}$.
\end{itemize}

This means that either $x_j^2x_k\notin \BB$ for any $k$ or $x_1^2x_k\notin \BB$ for any $k$, i.e. $\BB$ is not smooth by Lemma $1$.\\

Hence $\BB=\AAbar$ is the only possibility for this group $\GG$.\\

Let $\GG \cong \ZZ{4}\oplus \ZZ{2}$ be generated by 
\begin{gather*}
f_1=(1,1,\om,\om ^2,1,\om ^{-1}),\; f_2=(1,1,1,1,\eta,\eta),\; \om=\sqr{4}, \eta=\sqr{2}.
\end{gather*}

If $j=1$, then $\BB=\AAbar$. Hence we may assume that $j\in \{ 2,3,4,5 \}$.\\

Note that
\begin{itemize}
\item $wt_{f_1}(x_0^2x_2)=1\neq wt_{f_1}(x_2^2x_k)$ for any $k\neq 5$,
\item $wt_{f_2}(x_0^2x_2)=0\neq wt_{f_2}(x_2^2x_5)=1$,
\item $wt_{f_1}(x_0^2x_5)=-1\neq wt_{f_1}(x_5^2x_k)$ for any $k\neq 2$,
\item $wt_{f_2}(x_0^2x_5)=1\neq wt_{f_2}(x_5^2x_2)=0$
\item $wt_{f_2}(x_0^2x_4)=1\neq wt_{f_2}(x_5^2x_k)$ for any $k\neq 4,5 $,
\item $wt_{f_1}(x_0^2x_4)=0\neq wt_{f_1}(x_5^2x_4)=-2$,
\item $wt_{f_1}(x_0^2x_4)=0\neq wt_{f_1}(x_5^3)=-3$.
\end{itemize}

This means that if $j\in \{  2,4,5\}$, then either $x_j^2x_k\notin \BB$ for any $k$ or $x_5^2x_k\notin \BB$ for any $k$, i.e. $\BB$ is not smooth by Lemma $1$.\\

If $j=3$ and $m\in \{ x_0^2x_k, x_1^2x_k, x_0x_1x_k \}$ for some $k\in \{ 0,1,2,3,4,5 \}$, then $m\in \BB$ only if $wt_{f_1}(m)=wt_{f_1}(x_0^2x_3)=2$, i.e. $k=3$. This means that $(x_0,x_1)$ is a singular pair, i.e. $\BB$ is not smooth.\\

Hence $\BB=\AAbar$ is the only possibility for this group $\GG$.\\

Let $\GG \cong \ZZ{4}\oplus \ZZ{2}$ be generated by 
\begin{gather*}
f_1=(1,1,\om,\om ^2,1,\om ^{2}),\; f_2=(1,1,1,1,\eta,\eta),\; \om=\sqr{4}, \eta=\sqr{2}.
\end{gather*}

If $j=1$, then $\BB=\AAbar$. Hence we may assume that $j\in \{ 2,3,4,5 \}$.\\

Note that
\begin{itemize}
\item $wt_{f_1}(x_0^2x_2)=1\neq wt_{f_1}(x_2^2x_k)$ for any $k$.
\end{itemize}

This means that if $j=2$, then $x_2^2x_k\notin \BB$ for any $k$, i.e. $\BB$ is not smooth by Lemma $1$.\\

If $j\in \{ 3,5 \}$ and $m\in \{ x_0^2x_k, x_1^2x_k, x_4^2x_k, x_0x_1x_k, x_0x_4x_k, x_1x_4x_k \}$ for some $k\in \{ 0$, $1$, $2$, $3$, $4$, $5 \}$, then $m\in \BB$ only if $wt_{f_1}(m)=wt_{f_1}(x_0^2x_j)=2$, i.e. $k\in \{ 3,5 \}$. This means that $(x_0,x_1,x_4)$ is a singular triple, i.e. $\BB$ is not smooth.\\

If $j=4$ and $m\in \{ x_0^2x_k, x_1^2x_k, x_2^2x_k, x_0x_1x_k, x_0x_2x_k, x_1x_2x_k \}$ for some $k\in \{ 0,1,2,3,4,5 \}$, then $m\in \BB$ only if $wt_{f_2}(m)=wt_{f_2}(x_0^2x_j)=1$, i.e. $k\in \{ 4,5 \}$. This means that $(x_0,x_1,x_2)$ is a singular triple, i.e. $\BB$ is not smooth.\\

Hence $\BB=\AAbar$ is the only possibility for this group $\GG$.\\

Let $\GG \cong (\ZZ{2})^{\oplus 2}\oplus \ZZ{6}$ be generated by 
\begin{gather*}
f_1=(1,1,\om,\om,1,1),\; f_2=(1,1,1,\om,\om,1),\; f_3=(1,\eta^{-2},1,1,1, \eta),\; \om=\sqr{2},\eta=\sqr{6}.
\end{gather*}

Note that
\begin{itemize}
\item $wt_{f_3}(x_0^2x_1)=-2\neq wt_{f_3}(x_1^2x_k)$ for any $k$,
\item $wt_{f_3}(x_0^2x_5)=1\neq wt_{f_3}(x_5^2x_k)$ for any $k$,
\item $wt_{f_3}(x_0^2x_j)=0\neq wt_{f_3}(x_1^2x_k)$ for any $k\neq 1$ for any $j\in \{ 2,3 \}$,
\item $wt_{f_1}(x_0^2x_j)=1\neq wt_{f_1}(x_1^3)=0$ for any $j\in \{ 2,3 \}$,
\item $wt_{f_3}(x_0^2x_4)=0\neq wt_{f_3}(x_1^2x_k)$ for any $k\neq 1$,
\item $wt_{f_2}(x_0^2x_4)=1\neq wt_{f_2}(x_1^3)=0$.
\end{itemize}

This means that either $x_j^2x_k\notin \BB$ for any $k$ or $x_1^2x_k\notin \BB$ for any $k$, i.e. $\BB$ is not smooth by Lemma $1$.\\

Hence $\BB=\AAbar$ is the only possibility for this group $\GG$.\\

Let $\GG \cong \ZZ{33}$ be generated by 
\begin{gather*}
f=(1,\om ,{\om}^{-2},\om^4,{\om}^{-8},\om ^{16}),\; \om=\sqr{33}.
\end{gather*}

Note that
\begin{itemize}
\item $wt_{f}(x_0^2x_1)=1\neq wt_{f}(x_1^2x_k)$ for any $k$,
\item $wt_{f}(x_0^2x_2)=-2\neq wt_{f}(x_2^2x_k)$ for any $k$,
\item $wt_{f}(x_0^2x_3)=4\neq wt_{f}(x_3^2x_k)$ for any $k$,
\item $wt_{f}(x_0^2x_4)=-8\neq wt_{f}(x_4^2x_k)$ for any $k$,
\item $wt_{f}(x_0^2x_5)=16\neq wt_{f}(x_5^2x_k)$ for any $k$.
\end{itemize}

This means that $x_j^2x_k\notin \BB$ for any $k$, i.e. $\BB$ is not smooth by Lemma $1$.\\

Hence $\BB=\AAbar$ is the only possibility for this group $\GG$.\\

Let $\GG \cong \ZZ{2}\oplus \ZZ{15}$ be generated by 
\begin{gather*}
f_1=(1,1,1,1,1,\om),\; f_2=(1,\eta,\eta^{-2},\eta^4,\eta ^{-8},1),\; \om=\sqr{2}, \eta=\sqr{15}.
\end{gather*}

Note that
\begin{itemize}
\item $wt_{f_2}(x_0^2x_1)=1\neq wt_{f_2}(x_1^2x_k)$ for any $k$,
\item $wt_{f_2}(x_0^2x_2)=-2\neq wt_{f_2}(x_2^2x_k)$ for any $k$,
\item $wt_{f_2}(x_0^2x_3)=4\neq wt_{f_2}(x_3^2x_k)$ for any $k$,
\item $wt_{f_2}(x_0^2x_4)=-8\neq wt_{f_2}(x_4^2x_k)$ for any $k$,
\item $wt_{f_2}(x_0^2x_5)=0\neq wt_{f_2}(x_1^2x_k)$ for any $k\neq 2$,
\item $wt_{f_1}(x_0^2x_5)=1\neq wt_{f_1}(x_1^2x_2)$.
\end{itemize}

This means that either $x_j^2x_k\notin \BB$ for any $k$ or $x_1^2x_k\notin \BB$ for any $k$, i.e. $\BB$ is not smooth by Lemma $1$.\\

Hence $\BB=\AAbar$ is the only possibility for this group $\GG$.\\

Let $\GG \cong \ZZ{4}\oplus \ZZ{9}$ be generated by 
\begin{gather*}
f_1=(1,1,1,1,\om,\om^2),\; f_2=(1,\eta,\eta^{-2},\eta^4,1,1),\; \om=\sqr{4}, \eta=\sqr{9}.
\end{gather*}

Note that
\begin{itemize}
\item $wt_{f_2}(x_0^2x_1)=1\neq wt_{f_2}(x_1^2x_k)$ for any $k$,
\item $wt_{f_2}(x_0^2x_2)=-2\neq wt_{f_2}(x_2^2x_k)$ for any $k$,
\item $wt_{f_2}(x_0^2x_3)=4\neq wt_{f_2}(x_3^2x_k)$ for any $k$,
\item $wt_{f_1}(x_0^2x_4)=1\neq wt_{f_1}(x_4^2x_k)$ for any $k$,
\item $wt_{f_1}(x_0^2x_5)=2\neq wt_{f_1}(x_1^2x_k)$ for any $k\neq 5$,
\item $wt_{f_2}(x_0^2x_5)=0\neq wt_{f_2}(x_1^2x_5)=2$.
\end{itemize}

This means that either $x_j^2x_k\notin \BB$ for any $k$ or $x_1^2x_k\notin \BB$ for any $k$, i.e. $\BB$ is not smooth by Lemma $1$.\\

Hence $\BB=\AAbar$ is the only possibility for this group $\GG$.\\

Let $\GG \cong \ZZ{32}$ be generated by 
\begin{gather*}
f=(1,\om ,{\om}^{-2},\om^4,{\om}^{-8},\om ^{16}),\; \om=\sqr{32}.
\end{gather*}

Note that
\begin{itemize}
\item $wt_{f}(x_0^2x_j)\neq wt_{f}(x_j^2x_k)$ for any $k$ for any $j\in \{ 1,2,3,4 \}$,
\item $wt_{f}(x_0^2x_5)=16\neq wt_{f}(x_1^2x_k)$ for any $k$.
\end{itemize}

This means that either $x_j^2x_k\notin \BB$ for any $k$ or $x_1^2x_k\notin \BB$ for any $k$, i.e. $\BB$ is not smooth by Lemma $1$.\\

Hence $\BB=\AAbar$ is the only possibility for this group $\GG$.\\

Let $\GG \cong \ZZ{16}$ be generated by 
\begin{gather*}
f=(1,\om ,{\om}^{-2},\om^4,{\om}^{8},\om ^{-4}),\; \om=\sqr{16}.
\end{gather*}

Note that
\begin{itemize}
\item $wt_{f}(x_0^2x_j)\neq wt_{f}(x_j^2x_k)$ for any $k$ for any $j\in \{ 1,2 \}$,
\item $wt_{f}(x_0^2x_4)=8\neq wt_{f}(x_1^2x_k)$ for any $k$,
\item $wt_{f}(x_0^2x_3)=4\neq wt_{f}(x_1^2x_k)$ for any $k$,
\item $wt_{f}(x_0^2x_5)=-4\neq wt_{f}(x_1^2x_k)$ for any $k$.
\end{itemize}

This means that either $x_j^2x_k\notin \BB$ for any $k$ or $x_1^2x_k\notin \BB$ for any $k$, i.e. $\BB$ is not smooth by Lemma $1$.\\

Hence $\BB=\AAbar$ is the only possibility for this group $\GG$.\\

Let $\GG \cong \ZZ{8}\oplus \ZZ{2}$ be generated by 
\begin{gather*}
f_1=(1,\om,{\om}^{-2},\om^4,{\om}^{-1},{\om}^{2}),\; f_2=(1,1,1,1,\eta,1),\; \om=\sqr{8}, \eta=\sqr{2}.
\end{gather*}

Note that
\begin{itemize}
\item $wt_{f_2}(x_0^2x_4)=1\neq wt_{f_2}(x_4^2x_k)$ for any $k\neq 4$,
\item $wt_{f_1}(x_0^2x_4)=-1\neq wt_{f_1}(x_4^3)=-3$,
\item $wt_{f_1}(x_0^2x_1)=1\neq wt_{f_1}(x_1^2x_k)$ for any $k\neq 4$,
\item $wt_{f_2}(x_0^2x_1)=0\neq wt_{f_2}(x_1^2x_4)=1$.
\end{itemize}

This means that if $j\in \{ 1,4 \}$, then $x_j^2x_k\notin \BB$ for any $k$, i.e. $\BB$ is not smooth by Lemma~$1$.\\

If $j=2$, then $\BB$ consists of monomials $m=x_0^{i_0}x_1^{i_1}x_2^{i_2}x_3^{i_3}x_4^{i_4}x_5^{i_5}$ such that $wt_{f_1}(m)=wt_{f_1}(x_0^2x_2)=-2$ and $wt_{f_2}(m)=wt_{f_2}(x_0^2x_2)=0$, i.e.
$$
i_1-2i_2+4i_3-i_4+2i_5\equiv -2 \; mod \; 8,\;\;\;\;\; i_4\equiv 0 \; mod \; 2.
$$

This means that $\BB=\{ x_5^3, x_4^2x_0, x_1^2x_3, x_3^2x_2, x_0^2x_2, x_2^2x_5, x_0x_3x_5 \}$. This set coincides with $\AAbar$ upto a permutation of $x_i$ ($x_0\leftrightarrow x_5$, $x_2\leftrightarrow x_3$) and using $({\om}^{-2}f_1)^{-1}\cdot f_2$ and $f_2$ as generators of $\GG$.\\

If $j=3$, then $\BB$ consists of monomials $m=x_0^{i_0}x_1^{i_1}x_2^{i_2}x_3^{i_3}x_4^{i_4}x_5^{i_5}$ such that $wt_{f_1}(m)=wt_{f_1}(x_0^2x_3)=4$ and $wt_{f_2}(m)=wt_{f_2}(x_0^2x_3)=0$, i.e.
$$
i_1-2i_2+4i_3-i_4+2i_5\equiv 4 \; mod \; 8,\;\;\;\;\; i_4\equiv 0 \; mod \; 2.
$$

This means that $\BB=\{ x_3^3, x_4^2x_2, x_1^2x_5, x_2^2x_0, x_5^2x_0, x_0^2x_3, x_2x_3x_5 \}$. This set coincides with $\AAbar$ upto a permutation of $x_i$ ($x_0\leftrightarrow x_3$, $x_2\leftrightarrow x_5$) and using $({\om}^{4}f_1)^{5}$ and $f_2$ as generators of $\GG$.\\

If $j=5$, then $\BB$ consists of monomials $m=x_0^{i_0}x_1^{i_1}x_2^{i_2}x_3^{i_3}x_4^{i_4}x_5^{i_5}$ such that $wt_{f_1}(m)=wt_{f_1}(x_0^2x_5)=2$ and $wt_{f_2}(m)=wt_{f_2}(x_0^2x_5)=0$, i.e.
$$
i_1-2i_2+4i_3-i_4+2i_5\equiv 2 \; mod \; 8,\;\;\;\;\; i_4\equiv 0 \; mod \; 2.
$$

This means that $\BB=\{ x_2^3, x_4^2x_3, x_1^2x_0, x_3^2x_5, x_5^2x_2, x_0^2x_5, x_0x_2x_3 \}$. This set coincides with $\AAbar$ upto a permutation of $x_i$ ($x_0\leftrightarrow x_2$, $x_3\leftrightarrow x_5$) and using $({\om}^{2}f_1)^{3}\cdot f_2$ and $f_2$ as generators of $\GG$.\\

Hence $\BB=\AAbar$ is the only possibility for this group $\GG$ (upto a permutation of $x_i$).\\

Let $\GG \cong \ZZ{8}$ be generated by 
\begin{gather*}
f=(1,\om ,{\om}^{-2},\om^4,{\om}^{5},\om ^{2}),\; \om=\sqr{8}.
\end{gather*}

Note that
\begin{itemize}
\item $wt_{f}(x_0^2x_1)=1\neq wt_{f}(x_1^2x_k)$ for any $k$,
\item $wt_{f}(x_0^2x_4)=5\neq wt_{f}(x_4^2x_k)$ for any $k$.
\end{itemize}

This means that if $j\in \{ 1,4 \}$, then $x_j^2x_k\notin \BB$ for any $k$, i.e. $\BB$ is not smooth by Lemma~$1$.\\

If $j=2$, then $\BB$ consists of monomials $m=x_0^{i_0}x_1^{i_1}x_2^{i_2}x_3^{i_3}x_4^{i_4}x_5^{i_5}$ such that $wt_{f}(m)=wt_{f}(x_0^2x_2)=-2$, i.e.
$$
i_1-2i_2+4i_3-3i_4+2i_5\equiv -2 \; mod \; 8.
$$

This means that $\BB=\{ x_5^3, x_2^2x_5, x_4^2x_3, x_3^2x_2, x_1^2x_3, x_0^2x_2, x_0x_3x_5, x_0x_1x_4 \}$. This set coincides with $\AAbar$ upto a permutation of $x_i$ ($x_0\leftrightarrow x_5$, $x_2\leftrightarrow x_3$) and using $({\om}^{-2}f)^{-1}$ as a generator of $\GG$.\\

If $j=3$, then $\BB$ consists of monomials $m=x_0^{i_0}x_1^{i_1}x_2^{i_2}x_3^{i_3}x_4^{i_4}x_5^{i_5}$ such that $wt_{f}(m)=wt_{f}(x_0^2x_3)=4$, i.e.
$$
i_1-2i_2+4i_3-3i_4+2i_5\equiv 4 \; mod \; 8.
$$

This means that $\BB=\{ x_3^3, x_5^2x_0, x_1^2x_5, x_4^2x_5, x_0^2x_3, x_2^2x_0, x_2x_3x_5, x_1x_2x_4 \}$. This set coincides with $\AAbar$ upto a permutation of $x_i$ ($x_0\leftrightarrow x_3$, $x_2\leftrightarrow x_5$) and using $({\om}^{4}f)^{5}$ as a generator of $\GG$.\\

If $j=5$, then $\BB$ consists of monomials $m=x_0^{i_0}x_1^{i_1}x_2^{i_2}x_3^{i_3}x_4^{i_4}x_5^{i_5}$ such that $wt_{f}(m)=wt_{f}(x_0^2x_5)=2$, i.e.
$$
i_1-2i_2+4i_3-3i_4+2i_5\equiv 2 \; mod \; 8.
$$

This means that $\BB=\{ x_2^3, x_5^2x_2, x_3^2x_5, x_0^2x_5, x_1^2x_0, x_4^2x_0, x_0x_2x_3, x_1x_3x_4 \}$. This set coincides with $\AAbar$ upto a permutation of $x_i$ ($x_0\leftrightarrow x_2$, $x_3\leftrightarrow x_5$) and using $({\om}^{2}f)^{3}$ as a generator of $\GG$.\\

Hence $\BB=\AAbar$ is the only possibility for this group $\GG$ (upto a permutation of $x_i$).\\

Let $\GG \cong \ZZ{4}\oplus \ZZ{2}$ be generated by 
\begin{gather*}
f_1=(1,\om,{\om}^{2},\om^{-1},1,{\om}^{2}),\; f_2=(1,1,1,\eta,\eta,\eta),\; \om=\sqr{4}, \eta=\sqr{2}.
\end{gather*}

Note that
\begin{itemize}
\item $wt_{f_1}(x_0^2x_1)=1\neq wt_{f_1}(x_1^2x_k)$ for any $k\neq 3$,
\item $wt_{f_2}(x_0^2x_1)=0\neq wt_{f_2}(x_1^2x_3)=1$,
\item $wt_{f_1}(x_0^2x_3)=-1\neq wt_{f_1}(x_3^2x_k)$ for any $k\neq 1$,
\item $wt_{f_2}(x_0^2x_3)=1\neq wt_{f_2}(x_3^2x_1)=0$.
\end{itemize}

This means that if $j\in \{ 1,3 \}$, then $x_j^2x_k\notin \BB$ for any $k$, i.e. $\BB$ is not smooth by Lemma~$1$.\\

If $j=2$, then $\BB$ consists of monomials $m=x_0^{i_0}x_1^{i_1}x_2^{i_2}x_3^{i_3}x_4^{i_4}x_5^{i_5}$ such that $wt_{f_1}(m)=wt_{f_1}(x_0^2x_2)=2$ and $wt_{f_2}(m)=wt_{f_2}(x_0^2x_2)=0$, i.e.
$$
i_1+2i_2-i_3+2i_5\equiv 2 \; mod \; 4,\;\;\;\;\; i_3+i_4+i_5\equiv 0 \; mod \; 2.
$$

This means that $\BB=\{ x_2^3, x_3^2x_0, x_1^2x_0, x_5^2x_2, x_4^2x_2, x_0^2x_2, x_1x_3x_5, x_0x_4x_5 \}$. This set coincides with $\AAbar$ upto a permutation of $x_i$ ($x_0\leftrightarrow x_2$, $x_4\leftrightarrow x_5$) and using $({\om}^{2}f_1)^{-1}$ and $f_2$ as generators of $\GG$.\\

If $j=4$, then $\BB$ consists of monomials $m=x_0^{i_0}x_1^{i_1}x_2^{i_2}x_3^{i_3}x_4^{i_4}x_5^{i_5}$ such that $wt_{f_1}(m)=wt_{f_1}(x_0^2x_4)=0$ and $wt_{f_2}(m)=wt_{f_2}(x_0^2x_4)=1$, i.e.
$$
i_1+2i_2-i_3+2i_5\equiv 0 \; mod \; 4,\;\;\;\;\; i_3+i_4+i_5\equiv 1 \; mod \; 2.
$$

This means that $\BB=\{ x_4^3, x_3^2x_5, x_1^2x_5, x_2^2x_4, x_5^2x_4, x_0^2x_4, x_0x_2x_5, x_0x_1x_3 \}$. This set coincides with $\AAbar$ upto a permutation of $x_i$ ($x_0\leftrightarrow x_4$, $x_2\leftrightarrow x_5$) and using $f_1$ and $\eta f_2\cdot f_1^2$ as generators of $\GG$.\\

If $j=5$, then $\BB$ consists of monomials $m=x_0^{i_0}x_1^{i_1}x_2^{i_2}x_3^{i_3}x_4^{i_4}x_5^{i_5}$ such that $wt_{f_1}(m)=wt_{f_1}(x_0^2x_5)=2$ and $wt_{f_2}(m)=wt_{f_2}(x_0^2x_5)=1$, i.e.
$$
i_1+2i_2-i_3+2i_5\equiv 2 \; mod \; 4,\;\;\;\;\; i_3+i_4+i_5\equiv 1 \; mod \; 2.
$$

This means that $\BB=\{ x_5^3, x_0^2x_5, x_4^2x_5, x_2^2x_5, x_3^2x_4, x_1^2x_4, x_0x_2x_4, x_1x_2x_3 \}$. This set coincides with $\AAbar$ upto a permutation of $x_i$ ($x_0\leftrightarrow x_5$, $x_2\leftrightarrow x_4$) and using $({\om}^{2}f_1)^{-1}$ and $\eta f_2\cdot f_1^2$ as generators of $\GG$.\\

Hence $\BB=\AAbar$ is the only possibility for this group $\GG$ (upto a permutation of $x_i$).\\

{\it QED}\\

\newpage

\section{Arbitrary abelian automorphism groups.}

We will use a classification of finite abelian subgroups of $PGL(6)=Aut({\mathbb P}^5)$. Let $n\geq 2$ be an integer. The notion of an Ad-subgroup of $GL(n)$ is defined in \cite{LieGradings2}.\\

{\bf Definition 1. (\cite{LieGradings2})} A subgroup of semisimple matrices $K\subset GL(n)$ is an {\it Ad-subgroup}, if
\begin{enumerate}
\item[(i)] the commutator of $K$ lies in the center of $GL(n)$: $(K,K)\subset {\mathbb C}^{*}\cdot Id\subset GL(n)$,
\item[(ii)] $K$ is maximal with this property.
\end{enumerate}

Every finite abelian subgroup $G\subset PGL(n)$ is contained in the image of an Ad-subgroup of $GL(n)$ under the projection $GL(n)\rightarrow PGL(n)$. Indeed, the preimage of $G$ will consist of semisimple elements of $GL(n)$ and will satisfy condition (i) of Definition 1.\\

Ad-subgroups of $GL(n)$ are classified in \cite{LieGradings2}. Let us recall this classification. We will modify the notation slightly to make it more convenient for us. See also \cite{LieGradings1} and \cite{GangHan}.\\

{\bf Definition 2. (\cite{LieGradings2}, section 2.7)}\footnote{In \cite{LieGradings2} the authors define Pauli groups to be subgroups of $GL(n)$. The Pauli subgroup as we define it here is the image of the Pauli subgroup defined in \cite{LieGradings2} under the projection $GL(n)\rightarrow PGL(n)$.} The {\it Pauli subgroup} ${\mathcal P}_n\subset PGL(n)$ is the subgroup of $PGL(n)$ generated by the following two $n\times n$ matrices:
$$
P_n=\begin{pmatrix}
0 & 1 & 0 & \ldots & 0\\
0 & 0 & 1 & \ldots & 0\\
 & \ldots &  & \ldots & \\
0 & 0 & 0 & \ldots & 1\\
1 & 0 & 0 & \ldots & 0\\
\end{pmatrix}\; \mbox{and}\; 
W_n=\begin{pmatrix}
1 & 0 & 0 & \ldots & 0\\
0 & {\om}_n & 0 & \ldots & 0\\
0 & 0 & {\om}_n^2 & \ldots & 0\\
\vdots & \vdots & \vdots & \ddots & \vdots \\
0 & 0 & 0 & \ldots & {\om}_n^{n-1}\\
\end{pmatrix},
$$
where ${\om}_n=\sqr{n}$.\\

Note that ${\mathcal P}_n \cong \ZZ{n}\oplus \ZZ{n}$ with generators $P_n$ and $W_n$ (\cite{GangHan}, section 5).\\

Denote by $D_n\subset PGL(n)$ the subgroup generated by diagonal $n\times n$ matrices. We also set ${\mathcal P}_1=D_1=Id$.\\

The classification we need is given by the following theorem of Havl{\' i}cek, Patera, Pelantova (\cite{LieGradings1}, \cite{LieGradings2}).\\

{\bf Theorem 2. (Havl{\' i}cek-Patera-Pelantova, \cite{LieGradings2}, Theorem 3.2)}\footnote{Theorem 3.2 of \cite{LieGradings2} is more precise than what we state here.} {\it Any Ad-subgroup of $GL(n)$ is conjugate to the preimage under the projection $GL(n)\rightarrow PGL(n)$ of a subgroup 
$$
D_{k_0}\otimes {\mathcal P}_{k_1}\otimes \cdots \otimes {\mathcal P}_{k_s}\subset PGL(n),
$$
where $n=k_0\cdot k_1\cdot \cdots \cdot k_s$.}\\

It implies immediately a classification of finite abelian subgroups of $PGL(n)$.\\

{\bf Corollary 2.} {\it Every finite abelian subgroup of $PGL(n)$ is (upto a conjugation) contained in a subgroup
$$
D_{k_0}\otimes {\mathcal P}_{k_1}\otimes \cdots \otimes {\mathcal P}_{k_s}\subset PGL(n)
$$
for some decomposition $n=k_0\cdot k_1\cdot \cdots \cdot k_s$.}\\

In the case $n=6$ we obtain the following possibilities for a finite abelian subgroup $G\subset PGL(6)=Aut({\mathbb P}^5)$ (upto a conjugation):
\begin{enumerate}
\item[(1)] $G\subset D_6$,
\item[(2)] $G\subset D_3\otimes {\mathcal P}_{2}$,
\item[(3)] $G\subset D_2\otimes {\mathcal P}_{3}$,
\item[(4)] $G\subset {\mathcal P}_{6}$.
\end{enumerate}

The first possibility was considered in Theorem $1$. Let us consider the remaining three.\\

\subsection{Case of $G\subset D_3\otimes {\mathcal P}_{2}$.}

In this case
$$
{\mathcal P}_{2}=\{ {\sigma}_0= \begin{pmatrix}
1 & 0 \\
0 & 1\\
\end{pmatrix},\;
{\sigma}_1= \begin{pmatrix}
0 & 1 \\
1 & 0\\
\end{pmatrix}, \; 
{\sigma}_2= \begin{pmatrix}
1 & 0 \\
0 & -1\\
\end{pmatrix}, \; 
{\sigma}_3= \begin{pmatrix}
0 & -1 \\
1 & 0\\
\end{pmatrix} \}\subset PGL(2).
$$

We will also denote by ${\sigma}_i$ its natural preimage in $GL(2)$.\\

Recall that ${\mathcal P}_{2} \cong \ZZ{2}\oplus \ZZ{2}$ (\cite{GangHan}).\\

Let us denote by $x_0,x_1,x_2$ the coordinates corresponding to the factor $D_3$ and by $y_1,y_2$ the coordinates corresponding to the factor ${\mathcal P}_{2}$ of the product $D_3\otimes {\mathcal P}_{2}\subset PGL(6)=Aut({\mathbb P}^5)$. Then we can take 
$$
z_0=x_0y_1,\; z_1=x_1y_1,\; z_2=x_2y_1,\; z_3=x_0y_2,\; z_4=x_1y_2,\; z_5=x_2y_2
$$
as the homogeneous coordinates on ${\mathbb P}^5$.\\

Let $\pi\colon G\subset D_3\otimes {\mathcal P}_{2} \rightarrow {\mathcal P}_{2}$ be the projection onto the second factor\footnote{Notice that the homomorphism $GL(3)\times GL(2)\rightarrow PGL(6), (A,B)\mapsto A\otimes B$ factors through the embedding $PGL(3)\times PGL(2)\subset PGL(6)$.}.\\

{\bf Lemma 3.} {\it If $\pi$ is not surjective, then $G$ is conjugate to a subgroup of $D_6\subset PGL(6)$.}\\

{\it Proof:} It is enough to show that ${\pi}(G)\subset {\mathcal P}_{2}\subset PGL(2)$ is conjugate to a subgroup of $D_2\subset PGL(2)$.\\

Since ${\sigma}_1\cdot {\sigma}_2={\sigma}_3$ and ${\pi}(G)=Id$ or ${\pi}(G)\cong \ZZ{2}$ by the assumption, ${\pi}(G)$ contains at most one of the elements ${\sigma}_1, {\sigma}_2, {\sigma}_3\in {\mathcal P}_{2}$.\\

This means that ${\pi}(G)$ is generated by a single semisimple matrix ${\sigma}_i$ for some $i\in \{ 0,1,2,3 \}$. Hence ${\pi}(G)$ is indeed conjugate to a subgroup of $D_2\subset PGL(2)$. {\it QED}\\

{\bf Lemma 4.} {\it $\pi$ is never surjective in the case of $G\subset D_3\otimes {\mathcal P}_{2}$.}\\

{\it Proof:}

Let $X\subset {\mathbb P}^5$ be a smooth cubic fourfold given by a cubic form $F=F(z_0,z_1,z_2,z_3,z_4,z_5)$.\\

Let us assume that ${\pi}(G)={\mathcal P}_{2}$. This implies that $G$ contains elements $f_1=(1,{\om}^{c_1},{\om}^{c_2})\times {\sigma}_1$ and $f_2=(1,{\om}^{d_1},{\om}^{d_2})\times {\sigma}_2$, where $\om=\sqr{d}$, $c_i,d_i\in \ZZ{d}$ for some $d\geq 2$.\\

Suppose that $F$ contains a monomial $z_0^3=(x_0y_1)^3$ or $z_0^2z_3=(x_0y_1)^2(x_0y_2)$.\\

Since $f_1(F)\in {\mathbb C}^{*}\cdot F$ we conclude that $F$ should also contain a monomial $z_3^3=(x_0y_2)^3$ or $z_3^2z_0=(x_0y_2)^2(x_0y_1)$ respectively. In other words, we can write
$$
F=a_1 (x_0y_1)^3+a_2 (x_0y_2)^3+F'\; \mbox{or}
$$

$$
F=a_1 (x_0y_1)^2(x_0y_2)+a_2 (x_0y_2)^2(x_0y_1)+F'
$$
for some $a_1,a_2\in {\mathbb C}^{*}$ and for some cubic form $F'$ which does not contain the singled out monomials.

Since $f_2(F)=\lambda F$ for some $\lambda \in {\mathbb C}^{*}$ we conclude that in particular
$$
f_2(a_1 (x_0y_1)^3+a_2 (x_0y_2)^3)=\lambda (a_1 (x_0y_1)^3+a_2 (x_0y_2)^3),\; \mbox{or}
$$

$$
f_2(a_1 (x_0y_1)^2(x_0y_2)+a_2 (x_0y_2)^2(x_0y_1))=\lambda (a_1 (x_0y_1)^2(x_0y_2)+a_2 (x_0y_2)^2(x_0y_1))
$$
respectively.\\

Since ${\sigma}_2(y_1)=y_1$ and ${\sigma}_2(y_2)=-y_2$ by the definition, one arrives at a contradiction.\\

Since $X\subset {\mathbb P}^5$ is smooth, Lemma 1 from the beginning of this note (Lemma 1.3 in \cite{Liendo}) implies that (upto a permutation of $x_i$) $F$ should contain either a monomial $(x_0y_1)^2(x_1y_1)$ or a monomial $(x_0y_1)^2(x_1y_2)$. This allows us to decompose $F$ in the same way as above as follows:
$$
F=a_1 (x_0y_1)^2(x_1y_1)+a_2 (x_0y_2)^2(x_1y_2)+F'\; \mbox{or}
$$

$$
F=a_1 (x_0y_1)^2(x_1y_2)+a_2 (x_0y_2)^2(x_1y_1)+F'
$$
for some $a_1,a_2\in {\mathbb C}^{*}$ and for some cubic form $F'$ which does not contain the singled out monomials.\\

Then again we should have that either 
$$
f_2(a_1 (x_0y_1)^2(x_1y_1)+a_2 (x_0y_2)^2(x_1y_2))=\lambda (a_1 (x_0y_1)^2(x_1y_1)+a_2 (x_0y_2)^2(x_1y_2)),\; \mbox{or}
$$

$$
f_2(a_1 (x_0y_1)^2(x_1y_2)+a_2 (x_0y_2)^2(x_1y_1))=\lambda (a_1 (x_0y_1)^2(x_1y_2)+a_2 (x_0y_2)^2(x_1y_1))
$$
respectively  for some $\lambda\in {\mathbb C}^{*}$.\\

Since
$$
f_2((x_0y_1)^2(x_1y_1))={\om}^{d_1}\cdot (x_0y_1)^2(x_1y_1),\;\; f_2((x_0y_2)^2(x_1y_2))=-{\om}^{d_1}\cdot (x_0y_2)^2(x_1y_2),
$$

$$
f_2((x_0y_1)^2(x_1y_2))=-{\om}^{d_1}\cdot (x_0y_1)^2(x_1y_2),\;\; f_2((x_0y_2)^2(x_1y_1))={\om}^{d_1}\cdot (x_0y_2)^2(x_1y_1),
$$

we arrive at a contradiction again. {\it QED}\\

Hence the second possibility (the case of $G\subset D_3\otimes {\mathcal P}_{2}$) does not give us new automorphism groups.\\

\subsection{Case of $G\subset D_2\otimes {\mathcal P}_{3}$.}

In this case ${\mathcal P}_{3}=\{ P_3^i\cdot W_3^j \; \mid \; i,j=0,1,2 \}\subset PGL(3)$, where
$$
P_3= \begin{pmatrix}
0 & 1 & 0 \\
0 & 0 & 1\\
1 & 0 & 0\\
\end{pmatrix},\;
W_3= \begin{pmatrix}
1 & 0 & 0\\
0 & {\om}_3 & 0\\
0 & 0 & {\om}_3^2\\
\end{pmatrix}
$$
and ${\om}_3=\sqr{3}$.\\

Recall that ${\mathcal P}_{3} \cong \ZZ{3}\oplus \ZZ{3}$ (\cite{GangHan}).\\

Let us denote by $x_0,x_1$ the coordinates corresponding to the factor $D_2$ and by $y_1,y_2,y_3$ the coordinates corresponding to the factor ${\mathcal P}_{3}$ of the product $D_2\otimes {\mathcal P}_{3}\subset PGL(6)=Aut({\mathbb P}^5)$. Then we can take 
$$
z_0=x_0y_1,\; z_1=x_1y_1,\; z_2=x_0y_2,\; z_3=x_1y_2,\; z_4=x_0y_3,\; z_5=x_1y_3
$$
as the homogeneous coordinates on ${\mathbb P}^5$.\\

Let $\pi\colon G\subset D_2\otimes {\mathcal P}_{3} \rightarrow {\mathcal P}_{3}$ be the projection onto the second factor. Let $G_0=ker(\pi)$.\\

{\bf Lemma 5.} {\it If $\pi$ is not surjective, then $G$ is conjugate to a subgroup of $D_6\subset PGL(6)$.}\\

{\it Proof:} Without loss of generality we may assume that ${\pi}(G)\subset {\mathcal P}_{3}\subset PGL(3)$ is generated by $P_3^iW_3^j$ for some $i,j$. Hence ${\pi}(G)$ lifts to a subgroup of $GL(3)$ generated by a single semisimple matrix. Hence it is conjugate to a subgroup of $D_3\subset PGL(3)$. {\it QED}\\

Hence we can restrict our attention to the situations, where $\pi$ is surjective.\\

Let $X\subset {\mathbb P}^5$ be a smooth cubic fourfold given by a cubic form $F=F(z_0,z_1,z_2,z_3,z_4,z_5)$.\\

Let $F_x$ be the set of cubic monomials $x_{i}x_{j}x_{k}$ (in variables $x_0, x_1$) such that a monomial $(x_{i}y_{i'})(x_{j}y_{j'})(x_{k}y_{k'})$ appears in $F$ with a nonzero coefficient for some $i',j',k'$.\\

By Lemma 1 (\cite{Liendo}, Lemma 1.3) the cubic form $F$ should contain for any $i$ either $(x_0y_i)^2\cdot (x_0y_j)$ or $(x_0y_i)^2\cdot (x_1y_j)$ for some $j$.\\

We have two cases (upto a permutation of $x_i$):
\begin{itemize}
\item[(A)] $F$ contains $(x_0y_i)^2(x_0y_j)$ for some $i,j$,
\item[(B)] $F$ does not contain $(x_0y_i)^2(x_0y_j)$, $(x_1y_i)^2(x_1y_j)$ for any $i,j$.
\end{itemize}

In the first case $F_x$ contains $x_0^3$. In the second case $F$ for any $i$ contains $(x_0y_i)^2(x_1y_j)$ and $(x_1y_i)^2(x_0y_{j'})$ for some $j,j'$. In particular, in the second case $F_x$ contains both $x_0^2x_1$ and $x_1^2x_0$.\\

\subsubsection{Case A. $F$ contains $(x_0y_i)^2(x_0y_j)$ for some $i,j$. }

Let $F_0=F(z_0,0,z_2,0,z_4,0)$. In other words, we remove from $F$ all monomials containing $x_1$. The resulting cubic form $F_0$ is nonzero by the assumption.\\

Let $F_y$ be the set of cubic monomials $y_{i}y_{j}y_{k}$ (in variables $y_1$, $y_2$, $y_3$) such that a monomial $(x_{i'}y_{i})(x_{j'}y_{j})(x_{k'}y_{k})$ appears in $F_0$ with a nonzero coefficient for some $i',j',k'$.\\

{\bf Lemma 6A.} {\it Suppose that $G$ acts effectively on $X$, $F_x$ contains $x_0^3$ and $\pi$ is surjective. Then $F_y$ contains monomials only from exactly one of the following sets:
\begin{itemize}
\item $y_1^3, y_2^3, y_3^3, y_1y_2y_3$,
\item $y_1^2y_2, y_2^2y_3, y_3^2y_1$,
\item $y_1^2y_3, y_3^2y_2, y_2^2y_1$.
\end{itemize}}
\par
{\it Proof:} Let $y_1^{i_1}y_2^{i_2}y_3^{i_3}$ be a monomial in $F_y$.\\ 

Let $f=(1,{\om}^{c_1})\times W_3\in G$. Then the condition $f(F)\in {\mathbb C}^{*}\cdot F$ implies that the residue class $i_2+2i_3 \; mod \; 3$ is the same for all monomials in $F_y$.\\

If this residue class is $0\; mod \; 3$, then $i_2\equiv i_3\; mod \; 3$, which means exactly that $(i_1,i_2,i_3)\in \{ (3,0,0), (0,3,0), (0,0,3), (1,1,1) \}$.\\

If this residue class is $1\; mod \; 3$, then $i_2\equiv 1+i_3\; mod \; 3$, which means exactly that $(i_1,i_2,i_3)\in \{ (1,0,2), (0,2,1), (2,1,0) \}$.\\

If this residue class is $-1\; mod \; 3$, then $i_3\equiv 1+i_2\; mod \; 3$, which means exactly that $(i_1,i_2,i_3)\in \{ (1,2,0), (0,1,2), (2,0,1) \}$. {\it QED}\\

This Lemma implies that we need to consider three subcases in Case A:
\begin{itemize}
\item[(A1)] $F_y$ is contained in the set of monomials $y_1^3, y_2^3, y_3^3, y_1y_2y_3$,
\item[(A2)] $F_y$ consists of monomials $y_1^2y_2, y_2^2y_3, y_3^2y_1$,
\item[(A3)] $F_y$ consists of monomials $y_1^2y_3, y_3^2y_2, y_2^2y_1$.
\end{itemize}

\paragraph{Case A1. $F$ contains $(x_0y_i)^2(x_0y_j)$ for some $i,j$, $F_y$ is contained in the set of monomials $y_1^3, y_2^3, y_3^3$, $y_1y_2y_3$. }

Let us assume that either 
\begin{multline*}
F=a_1\cdot (x_0y_1)^3 + a_2\cdot (x_0y_2)^3 + a_3\cdot (x_0y_3)^3 + \\
+b_1\cdot (x_1y_1)^3 + b_2\cdot (x_1y_2)^3 + b_3\cdot (x_1y_3)^3 + F',
\end{multline*}
or
\begin{multline*}
F=a_1\cdot (x_0y_1)^3 + a_2\cdot (x_0y_2)^3 + a_3\cdot (x_0y_3)^3 + \\
+b_1\cdot (x_1y_1)^2(x_0y_1) + b_2\cdot (x_1y_2)^2(x_0y_2) + b_3\cdot (x_1y_3)^2(x_0y_3) + F',
\end{multline*}
or
\begin{multline*}
F=a_1\cdot (x_0y_1)^3 + a_2\cdot (x_0y_2)^3 + a_3\cdot (x_0y_3)^3 + \\
+b_1\cdot (x_1y_1)^2(x_0y_2) + b_2\cdot (x_1y_2)^2(x_0y_3) + b_3\cdot (x_1y_3)^2(x_0y_1) + F',
\end{multline*}
or
\begin{multline*}
F=a_1\cdot (x_0y_1)^3 + a_2\cdot (x_0y_2)^3 + a_3\cdot (x_0y_3)^3 + \\
+b_1\cdot (x_1y_1)^2(x_0y_3) + b_2\cdot (x_1y_2)^2(x_0y_1) + b_3\cdot (x_1y_3)^2(x_0y_2) + F',
\end{multline*}
where $F'$ does not contain the singled out monomials and $a_i,b_i\in {\mathbb C}^{*}$.\\

Let us find $G_0=ker(\pi)$. If $f_0=(1,{\om}_0^{c_1})\times 1\in G_0$, then $3c_1=0$ in the first case and $2c_1=0$ in the other three cases.\\

Hence in the first case $G_0\cong \ZZ{3}$ with a generator $f_0=(1,{\om}_0)\times 1$, where ${\om}_0=\sqr{3}$.\\

In the other three cases $G_0\cong \ZZ{2}$ with a generator $f_0=(1,{\om}_0)\times 1$, where ${\om}_0=\sqr{2}$.\\

Let $f_1=(1,{\om}^{c_1})\times P_3$ and $f_2=(1,{\om}^{d_1})\times W_3$ be elements of $G$.\\

Then 
\begin{multline*}
f_1(F)=a_1\cdot (x_0y_2)^3 + a_2\cdot (x_0y_3)^3 + a_3\cdot (x_0y_1)^3 + \\
+{\om}^{3c_1}\cdot (b_1\cdot (x_1y_2)^3 + b_2\cdot (x_1y_3)^3 + b_3\cdot (x_1y_1)^3) + f_1(F'),
\end{multline*}
or
\begin{multline*}
f_1(F)=a_1\cdot (x_0y_2)^3 + a_2\cdot (x_0y_3)^3 + a_3\cdot (x_0y_1)^3 + \\
+{\om}^{2c_1}\cdot  ( b_1\cdot (x_1y_2)^2(x_0y_2) + b_2\cdot (x_1y_3)^2(x_0y_3) + b_3\cdot (x_1y_1)^2(x_0y_1)) + f_1(F'),
\end{multline*}
or
\begin{multline*}
f_1(F)=a_1\cdot (x_0y_2)^3 + a_2\cdot (x_0y_3)^3 + a_3\cdot (x_0y_1)^3 + \\
+{\om}^{2c_1}\cdot  ( b_1\cdot (x_1y_2)^2(x_0y_3) + b_2\cdot (x_1y_3)^2(x_0y_1) + b_3\cdot (x_1y_1)^2(x_0y_2)) + f_1(F'),
\end{multline*}
or
\begin{multline*}
f_1(F)=a_1\cdot (x_0y_2)^3 + a_2\cdot (x_0y_3)^3 + a_3\cdot (x_0y_1)^3 + \\
+{\om}^{2c_1}\cdot  ( b_1\cdot (x_1y_2)^2(x_0y_1) + b_2\cdot (x_1y_3)^2(x_0y_2) + b_3\cdot (x_1y_1)^2(x_0y_3) )+ f_1(F'),
\end{multline*}
respectively and 

\begin{multline*}
f_2(F)=a_1\cdot (x_0y_1)^3 + a_2\cdot (x_0y_2)^3 + a_3\cdot (x_0y_3)^3 + \\
+{\om}^{3d_1}\cdot (b_1\cdot (x_1y_1)^3 + b_2\cdot (x_1y_2)^3 + b_3\cdot (x_1y_3)^3) + f_2(F'),
\end{multline*}
or
\begin{multline*}
f_2(F)=a_1\cdot (x_0y_1)^3 + a_2\cdot (x_0y_2)^3 + a_3\cdot (x_0y_3)^3 + \\
+{\om}^{2d_1}\cdot (b_1\cdot (x_1y_1)^2(x_0y_1) + b_2\cdot (x_1y_2)^2(x_0y_2) + b_3\cdot (x_1y_3)^2(x_0y_3)) + f_2(F'),
\end{multline*}
or
\begin{multline*}
f_2(F)=a_1\cdot (x_0y_1)^3 + a_2\cdot (x_0y_2)^3 + a_3\cdot (x_0y_3)^3 + \\
+{\om}^{2d_1}\cdot {\om}_3\cdot (b_1\cdot (x_1y_1)^2(x_0y_2) + b_2\cdot (x_1y_2)^2(x_0y_3) + b_3\cdot (x_1y_3)^2(x_0y_1)) + f_2(F'),
\end{multline*}
or
\begin{multline*}
f_2(F)=a_1\cdot (x_0y_1)^3 + a_2\cdot (x_0y_2)^3 + a_3\cdot (x_0y_3)^3 + \\
+{\om}^{2d_1}\cdot {\om}_3^{2}\cdot (b_1\cdot (x_1y_1)^2(x_0y_3) + b_2\cdot (x_1y_2)^2(x_0y_1) + b_3\cdot (x_1y_3)^2(x_0y_2)) + f_2(F'),
\end{multline*}
respectively.\\ 

Since $f_1(F)\in {\mathbb C}^{*}\cdot F$, this implies that 
\begin{itemize}
\item $a_3=\lambda a_1$, $a_2=(\lambda )^2 a_1$, $b_3=\lambda \cdot {\om}^{-3c_1} \cdot b_1$, $b_2=(\lambda )^2 \cdot {\om}^{-6c_1} \cdot b_1$, ${\om}^{9c_1}=1$ for some $\lambda\in \mathbb C$ such that ${\lambda }^3=1$ in the first case, and
\item $a_3=\lambda a_1$, $a_2=(\lambda )^2 a_1$, $b_3=\lambda \cdot {\om}^{-2c_1} \cdot b_1$, $b_2=(\lambda )^2 \cdot {\om}^{-4c_1} \cdot b_1$, ${\om}^{6c_1}=1$ for some $\lambda\in \mathbb C$ such that ${\lambda }^3=1$ in the other three cases.
\end{itemize}

Since $f_2(F)\in {\mathbb C}^{*}\cdot F$, this implies that $f_2(F)=F$ and so
\begin{itemize}
\item ${\om}^{3d_1}=1$,
\item ${\om}^{2d_1}=1$,
\item ${\om}^{2d_1}={\om}_3^2$,
\item ${\om}^{2d_1}={\om}_3$,
\end{itemize}
respectively in our four cases.\\

Since we are interested in $f_1, f_2\in G$ only modulo $G_0$ we can take
\begin{itemize}
\item $f_1=(1,{\om}^i)\times P_3$, $f_2=(1,1)\times W_3$ with some $i\in \{ 0,1,2 \}$ and $\om = \sqr{9}$,
\item $f_1=(1,{\om}^i)\times P_3$, $f_2=(1,1)\times W_3$ with some $i\in \{ 0,1,2 \}$ and $\om = \sqr{6}$,
\item $f_1=(1,{\om}^i)\times P_3$, $f_2=(1,{\om}^j)\times W_3$ with some $i,j\in \{ 0,1,2 \}$, $j\equiv 2\; mod \; 3$ and $\om = \sqr{6}$,
\item $f_1=(1,{\om}^i)\times P_3$, $f_2=(1,{\om}^j)\times W_3$ with some $i,j\in \{ 0,1,2 \}$, $j\equiv 1\; mod \; 3$ and $\om = \sqr{6}$
\end{itemize}
respectively in our four cases.\\

Note that in the fourth case $f_2^3=f_0$. If $i=1$, then $f_1^3=f_0$ in all four cases. If $i=2$, then $f_1^6=f_0$ in the first case.\\

We conclude that if $i=0$, then 
\begin{itemize}
\item $G\cong (\ZZ{3})^{\oplus 3}$ with generators $f_0=(1,{\om}_0)\times 1$, $f_1=(1,1)\times P_3$, $f_2=(1,1)\times W_3$, where ${\om}_0=\sqr{3}, \om = \sqr{9}$,
\item $G\cong \ZZ{2}\oplus (\ZZ{3})^{\oplus 2}$ with generators $f_0=(1,{\om}_0)\times 1$, $f_1=(1,1)\times P_3$, $f_2=(1,1)\times W_3$, where ${\om}_0=\sqr{2}, \om = \sqr{6}$,
\item $G\cong \ZZ{2}\oplus (\ZZ{3})^{\oplus 2}$ with generators $f_0=(1,{\om}_0)\times 1$, $f_1=(1,1)\times P_3$, $f_2=(1,{\om}^2)\times W_3$, where ${\om}_0=\sqr{2}, \om = \sqr{6}$,
\item $G\cong \ZZ{3}\oplus \ZZ{6}$ with generators $f_1=(1,1)\times P_3$, $f_2=(1,{\om})\times W_3$, where ${\om}_0=\sqr{2}, \om = \sqr{6}$
\end{itemize}
respectively in our four cases.\\

If $i=1$, then 
\begin{itemize}
\item $G\cong \ZZ{9}\oplus \ZZ{3}$ with generators $f_1=(1,{\om})\times P_3$, $f_2=(1,1)\times W_3$, where $\om = \sqr{9}$,
\item $G\cong \ZZ{6}\oplus \ZZ{3}$ with generators $f_1=(1,{\om})\times P_3$, $f_2=(1,1)\times W_3$, where $\om = \sqr{6}$,
\item $G\cong \ZZ{6}\oplus \ZZ{3}$ with generators $f_1=(1,{\om})\times P_3$, $f_2=(1,{\om}^2)\times W_3$, where $\om = \sqr{6}$,
\item $G\cong \ZZ{6}\oplus \ZZ{3}$ with generators $f_1=(1,{\om})\times P_3$, $f_2f_1=(1,{\om}^2)\times W_3P_3$, where $\om = \sqr{6}$
\end{itemize}
respectively in our four cases.\\

If $i=2$, then 
\begin{itemize}
\item $G\cong \ZZ{9}\oplus \ZZ{3}$ with generators $f_1=(1,{\om}^2)\times P_3$, $f_2=(1,1)\times W_3$, where ${\om}_0=\sqr{3}, \om = \sqr{9}$,
\item $G\cong \ZZ{2}\oplus (\ZZ{3})^{\oplus 2}$ with generators $f_0=(1,{\om}_0)\times 1$, $f_1=(1,{\om}^2)\times P_3$, $f_2=(1,1)\times W_3$, where ${\om}_0=\sqr{2}, \om = \sqr{6}$,
\item $G\cong \ZZ{2}\oplus (\ZZ{3})^{\oplus 2}$ with generators $f_0=(1,{\om}_0)\times 1$, $f_1=(1,{\om}^2)\times P_3$, $f_2=(1,{\om}^2)\times W_3$, where ${\om}_0=\sqr{2}, \om = \sqr{6}$,
\item $G\cong \ZZ{6}\oplus \ZZ{3}$ with generators $f_1=(1,{\om}^2)\times P_3$, $f_2=(1,{\om})\times W_3$, where ${\om}_0=\sqr{2}, \om = \sqr{6}$
\end{itemize}
respectively in our four cases.\\

In order to describe all smooth cubic fourfolds which admit these group actions, we need to determine which cubic forms $F'$ can appear.\\

Since $F_x$ contains the monomial $x_0^3$ by our assumption, the invariance under $f_0$ requires that $F_x$ contains only monomials $x_0^3, x_1^3$ in the first case and only monomials $x_0^3, x_1^2x_0$ in the other three cases.\\

The invariance under $f_2$ requires that in the first two cases only monomials $y_1^3$, $y_2^3$, $y_3^3$, $y_1y_2y_3$ appear in $F$.\\

Let $F'_0=F'_0(z_0,z_1,z_2,z_3,z_4,z_5)$ be a cubic form which is obtained from $F'$ by removing all cubic monomials in $z_0, z_2, z_4$.\\

Let us denote by $F'_y$ the set of cubic monomials $y_{i}y_{j}y_{k}$ (in variables $y_1$, $y_2$, $y_3$) such that a monomial $(x_{i'}y_{i})(x_{j'}y_{j})(x_{k'}y_{k})$ appears in $F'_0$ with a nonzero coefficient for some $i',j',k'$.\\

The invariance under $f_2$ requires that $y_1^{i_1}y_2^{i_2}y_3^{i_3}$ is contained in $F'_y$ only if $i_2\equiv 1+i_3\; mod \; 3$ (in the third case) or $i_3\equiv 1+i_2\; mod \; 3$ (in the fourth case).\\

Hence $F'$ may contain only monomials
\begin{itemize}
\item $(x_0y_1)(x_0y_2)(x_0y_3)$, $(x_1y_1)(x_1y_2)(x_1y_3)$ in the first case,
\item $(x_0y_1)(x_0y_2)(x_0y_3)$, $(x_1y_1)(x_1y_2)(x_0y_3)$, $(x_1y_1)(x_0y_2)(x_1y_3)$, $(x_0y_1)(x_1y_2)(x_1y_3)$ in the second case,
\item $(x_0y_1)(x_0y_2)(x_0y_3)$, $(x_1y_1)(x_1y_2)(x_0y_1)$, $(x_1y_2)(x_1y_3)(x_0y_2)$, $(x_1y_3)(x_1y_1)(x_0y_3)$ in the third case,
\item $(x_0y_1)(x_0y_2)(x_0y_3)$, $(x_1y_1)(x_1y_3)(x_0y_1)$, $(x_1y_2)(x_1y_1)(x_0y_2)$, $(x_1y_3)(x_1y_2)(x_0y_3)$ in the fourth case.
\end{itemize}

We conclude that either
\begin{multline*}
F=a_1\cdot (x_0y_1)^3 + a_2\cdot (x_0y_2)^3 + a_3\cdot (x_0y_3)^3 + \\
+b_1\cdot (x_1y_1)^3 + b_2\cdot (x_1y_2)^3 + b_3\cdot (x_1y_3)^3 + \\
+c_0\cdot (x_0y_1)(x_0y_2)(x_0y_3)+c_1\cdot (x_1y_1)(x_1y_2)(x_1y_3),
\end{multline*}
or
\begin{multline*}
F=a_1\cdot (x_0y_1)^3 + a_2\cdot (x_0y_2)^3 + a_3\cdot (x_0y_3)^3 + \\
+b_1\cdot (x_1y_1)^2(x_0y_1) + b_2\cdot (x_1y_2)^2(x_0y_2) + b_3\cdot (x_1y_3)^2(x_0y_3) +\\
+c_0\cdot (x_0y_1)(x_0y_2)(x_0y_3) +\\
+c_1\cdot (x_1y_1)(x_1y_2)(x_0y_3)+c_2\cdot (x_1y_1)(x_0y_2)(x_1y_3) +c_3\cdot (x_0y_1)(x_1y_2)(x_1y_3),
\end{multline*}
or
\begin{multline*}
F=a_1\cdot (x_0y_1)^3 + a_2\cdot (x_0y_2)^3 + a_3\cdot (x_0y_3)^3 + \\
+b_1\cdot (x_1y_1)^2(x_0y_2) + b_2\cdot (x_1y_2)^2(x_0y_3) + b_3\cdot (x_1y_3)^2(x_0y_1) + \\
+c_0\cdot (x_0y_1)(x_0y_2)(x_0y_3) +\\
+c_1\cdot (x_1y_1)(x_1y_2)(x_0y_1)+c_2\cdot (x_1y_2)(x_1y_3)(x_0y_2) +c_3\cdot (x_1y_3)(x_1y_1)(x_0y_3),
\end{multline*}
or
\begin{multline*}
F=a_1\cdot (x_0y_1)^3 + a_2\cdot (x_0y_2)^3 + a_3\cdot (x_0y_3)^3 + \\
+b_1\cdot (x_1y_1)^2(x_0y_3) + b_2\cdot (x_1y_2)^2(x_0y_1) + b_3\cdot (x_1y_3)^2(x_0y_2) + \\
+c_0\cdot (x_0y_1)(x_0y_2)(x_0y_3) +\\
+c_1\cdot (x_1y_1)(x_1y_3)(x_0y_1)+c_2\cdot (x_1y_2)(x_1y_1)(x_0y_2) +c_3\cdot (x_1y_3)(x_1y_2)(x_0y_3),
\end{multline*}

respectively in our four cases, where $a_i, b_i\in {\mathbb C}^{*}$, $(a_3/a_1)^3=1$, $a_2=a_3^2/a_1$ and

\begin{itemize}
\item $b_3= {\om}^{-3i} \cdot a_3b_1/a_1$, $b_2={\om}^{3i} \cdot a_1b_1/a_3$ in the first case,
\item $b_3= {\om}^{-2i} \cdot a_3b_1/a_1$, $b_2={\om}^{2i} \cdot a_1b_1/a_3$, $c_3={\om}^{2i} \cdot a_1c_1/a_3$, $c_2={\om}^{-2i} \cdot a_3c_1/a_1$ in the second case,
\item $b_3= {\om}^{-2i} \cdot a_3b_1/a_1$, $b_2={\om}^{2i} \cdot a_1b_1/a_3$, $c_2={\om}^{2i} \cdot a_1c_1/a_3$, $c_3={\om}^{-2i} \cdot a_3c_1/a_1$ in the third and the fourth cases.
\end{itemize}

Moreover, in any of the cases if $c_0\neq 0$, then $a_1=a_3$. Also, in the first case $a_3={\om}^{3i}\cdot a_1$ whenever $c_1\neq 0$.\\

It is immediate that cubic fourfolds given by the generic equations above are smooth (for all four cases).\\

The remaining possibilities in Case A1 are the following two:
\begin{multline*}
F=a_1\cdot (x_0y_1)^3 + a_2\cdot (x_0y_2)^3 + a_3\cdot (x_0y_3)^3 + \\
+b_1\cdot (x_1y_1)^2(x_1y_2) + b_2\cdot (x_1y_2)^2(x_1y_3) + b_3\cdot (x_1y_3)^2(x_1y_1) + F',
\end{multline*}
or
\begin{multline*}
F=a_1\cdot (x_0y_1)^3 + a_2\cdot (x_0y_2)^3 + a_3\cdot (x_0y_3)^3 + \\
+b_1\cdot (x_1y_1)^2(x_1y_3) + b_2\cdot (x_1y_2)^2(x_1y_1) + b_3\cdot (x_1y_3)^2(x_1y_2) + F',
\end{multline*}
where $F'$ does not contain the singled out monomials and $a_i,b_i\in {\mathbb C}^{*}$.\\

Let us find $G_0=ker(\pi)$. If $f_0=(1,{\om}_0^{c_1})\times 1\in G_0$, then $3c_1=0$. This means that $G_0\cong \ZZ{3}$ with a generator $f_0=(1,{\om}_0)\times 1$, where ${\om}_0=\sqr{3}$.\\

Let $f_1=(1,{\om}^{c_1})\times P_3$ and $f_2=(1,{\om}^{d_1})\times W_3$ be elements of $G$.\\

Then 
\begin{multline*}
f_1(F)=a_1\cdot (x_0y_2)^3 + a_2\cdot (x_0y_3)^3 + a_3\cdot (x_0y_1)^3 + \\
+{\om}^{3c_1}\cdot (b_1\cdot (x_1y_2)^2(x_1y_3) + b_2\cdot (x_1y_3)^2(x_1y_1) + b_3\cdot (x_1y_1)^2(x_1y_2)) + f_1(F'),
\end{multline*}
or
\begin{multline*}
f_1(F)=a_1\cdot (x_0y_2)^3 + a_2\cdot (x_0y_3)^3 + a_3\cdot (x_0y_1)^3 + \\
+{\om}^{3c_1}\cdot  ( b_1\cdot (x_1y_2)^2(x_1y_1) + b_2\cdot (x_1y_3)^2(x_1y_2) + b_3\cdot (x_1y_1)^2(x_1y_3)) + f_1(F')
\end{multline*}

respectively and 

\begin{multline*}
f_2(F)=a_1\cdot (x_0y_1)^3 + a_2\cdot (x_0y_2)^3 + a_3\cdot (x_0y_3)^3 + \\
+{\om}^{3d_1}\cdot {\om}_3\cdot (b_1\cdot (x_1y_1)^2(x_1y_2) + b_2\cdot (x_1y_2)^2(x_1y_3) + b_3\cdot (x_1y_3)^2(x_1y_1)) + f_2(F'),
\end{multline*}
or
\begin{multline*}
f_2(F)=a_1\cdot (x_0y_1)^3 + a_2\cdot (x_0y_2)^3 + a_3\cdot (x_0y_3)^3 + \\
+{\om}^{3d_1}\cdot {\om}_3^2\cdot (b_1\cdot (x_1y_1)^2(x_1y_3) + b_2\cdot (x_1y_2)^2(x_1y_1) + b_3\cdot (x_1y_3)^2(x_1y_2)) + f_2(F')
\end{multline*}

respectively.\\ 

Since $f_1(F)\in {\mathbb C}^{*}\cdot F$, this implies that in both cases $a_3=\lambda a_1$, $a_2=(\lambda )^2 a_1$, $b_3=\lambda \cdot {\om}^{-3c_1} \cdot b_1$, $b_2=(\lambda )^2 \cdot {\om}^{-6c_1} \cdot b_1$, ${\om}^{9c_1}=1$ for some $\lambda\in \mathbb C$ such that ${\lambda }^3=1$.\\

Since $f_2(F)\in {\mathbb C}^{*}\cdot F$, this implies that $f_2(F)=F$ and so
\begin{itemize}
\item ${\om}^{3d_1}={\om}_3^2$,
\item ${\om}^{3d_1}={\om}_3$
\end{itemize}
respectively in our two cases.\\

Since we are interested in $f_1, f_2\in G$ only modulo $G_0$ we can take
\begin{itemize}
\item $f_1=(1,{\om}^i)\times P_3$, $f_2=(1,{\om}^j)\times W_3$ with some $i,j\in \{ 0,1,2 \}$, $j\equiv 2\; mod \; 3$ and $\om = \sqr{9}$,
\item $f_1=(1,{\om}^i)\times P_3$, $f_2=(1,{\om}^j)\times W_3$ with some $i,j\in \{ 0,1,2 \}$, $j\equiv 1\; mod \; 3$ and $\om = \sqr{9}$
\end{itemize}
respectively in our two cases.\\

Note that in the first case $f_2^6=f_0$ and in the second case $f_2^3=f_0$.\\

We conclude that if $i=0$, then 
\begin{itemize}
\item $G\cong \ZZ{3}\oplus \ZZ{9}$ with generators $f_1=(1,1)\times P_3$, $f_2=(1,{\om}^2)\times W_3$, where $\om = \sqr{9}$,
\item $G\cong \ZZ{3}\oplus \ZZ{9}$ with generators $f_1=(1,1)\times P_3$, $f_2=(1,{\om})\times W_3$, where $\om = \sqr{9}$
\end{itemize}
respectively in our two cases.\\

If $i=1$, then 
\begin{itemize}
\item $G\cong \ZZ{3}\oplus \ZZ{9}$ with generators $f_1=(1,{\om})\times P_3$, $f_1f_2=(1,{\om}^3)\times P_3W_3$, where $\om = \sqr{9}$,
\item $G\cong \ZZ{3}\oplus \ZZ{9}$ with generators $f_1=(1,{\om})\times P_3$, $f_1^2f_2=(1,{\om}^3)\times (P_3)^2(W_3)$, where $\om = \sqr{9}$
\end{itemize}
respectively in our two cases.\\

If $i=2$, then 
\begin{itemize}
\item $G\cong \ZZ{3}\oplus \ZZ{9}$ with generators $f_1=(1,{\om}^2)\times P_3$, $f_1^2f_2=(1,{\om}^6)\times (P_3)^2(W_3)$, where $\om = \sqr{9}$,
\item $G\cong \ZZ{3}\oplus \ZZ{9}$ with generators $f_1=(1,{\om}^2)\times P_3$, $f_1f_2=(1,{\om}^3)\times P_3W_3$, where $\om = \sqr{9}$
\end{itemize}
respectively in our two cases.\\

In order to describe all smooth cubic fourfolds which admit these group actions, we need to determine which cubic forms $F'$ can appear.\\

Since $F_x$ contains the monomial $x_0^3$ by our assumption, the invariance under $f_0$ requires that $F_x$ contains only monomials $x_0^3, x_1^3$.\\

Let $F'_0=F'_0(z_0,z_1,z_2,z_3,z_4,z_5)$ be a cubic form which is obtained from $F'$ by removing all cubic monomials in $z_0, z_2, z_4$.\\

Let us denote by $F'_y$ the set of cubic monomials $y_{i}y_{j}y_{k}$ (in variables $y_1$, $y_2$, $y_3$) such that a monomial $(x_{i'}y_{i})(x_{j'}y_{j})(x_{k'}y_{k})$ appears in $F'_0$ with a nonzero coefficient for some $i',j',k'$.\\

The invariance under $f_2$ requires that $y_1^{i_1}y_2^{i_2}y_3^{i_3}$ is contained in $F'_y$ only if $i_2\equiv 1+i_3\; mod \; 3$ (in the first case) or $i_3\equiv 1+i_2\; mod \; 3$ (in the second case).\\

Hence $F'$ may contain only the monomial $(x_0y_1)(x_0y_2)(x_0y_3)$.\\

We conclude that either
\begin{multline*}
F=a_1\cdot (x_0y_1)^3 + a_2\cdot (x_0y_2)^3 + a_3\cdot (x_0y_3)^3 + \\
+b_1\cdot (x_1y_1)^2(x_1y_2) + b_2\cdot (x_1y_2)^2(x_1y_3) + b_3\cdot (x_1y_3)^2(x_1y_1) + \\
+c_0\cdot (x_0y_1)(x_0y_2)(x_0y_3),
\end{multline*}
or
\begin{multline*}
F=a_1\cdot (x_0y_1)^3 + a_2\cdot (x_0y_2)^3 + a_3\cdot (x_0y_3)^3 + \\
+b_1\cdot (x_1y_1)^2(x_1y_3) + b_2\cdot (x_1y_2)^2(x_1y_1) + b_3\cdot (x_1y_3)^2(x_1y_2) + \\
+c_0\cdot (x_0y_1)(x_0y_2)(x_0y_3)
\end{multline*}

respectively in our two cases, where $a_i, b_i\in {\mathbb C}^{*}$, $(a_3/a_1)^3=1$, $a_2=a_3^2/a_1$, $b_3={\om}^{-3i} \cdot a_3b_1/a_1$, $b_2={\om}^{3i} \cdot a_1b_1/a_3$.\\

Moreover, in any of the cases if $c_0\neq 0$, then $a_1=a_3$.\\

It is immediate that cubic fourfolds given by the generic equations above are smooth (for both cases).\\

\paragraph{Case A2. $F$ contains $(x_0y_i)^2(x_0y_j)$ for some $i,j$, $F_y$ consists of monomials $y_1^2y_2, y_2^2y_3$, $y_3^2y_1$. }

Our assumptions imply that either
\begin{multline*}
F=a_1\cdot (x_0y_1)^2(x_0y_2) + a_2\cdot (x_0y_2)^2(x_0y_3) + a_3\cdot (x_0y_3)^2(x_0y_1) + \\
+b_1\cdot (x_1y_1)^2(x_1y_2) + b_2\cdot (x_1y_2)^2(x_1y_3) + b_3\cdot (x_1y_3)^2(x_1y_1) + F',
\end{multline*}
or
\begin{multline*}
F=a_1\cdot (x_0y_1)^2(x_0y_2) + a_2\cdot (x_0y_2)^2(x_0y_3) + a_3\cdot (x_0y_3)^2(x_0y_1) + \\
+b_1\cdot (x_1y_1)^2(x_1y_3) + b_2\cdot (x_1y_2)^2(x_1y_1) + b_3\cdot (x_1y_3)^2(x_1y_2) + F',
\end{multline*}
or
\begin{multline*}
F=a_1\cdot (x_0y_1)^2(x_0y_2) + a_2\cdot (x_0y_2)^2(x_0y_3) + a_3\cdot (x_0y_3)^2(x_0y_1) + \\
+b_1\cdot (x_1y_1)^2(x_0y_1) + b_2\cdot (x_1y_2)^2(x_0y_2) + b_3\cdot (x_1y_3)^2(x_0y_3) + F',
\end{multline*}
or
\begin{multline*}
F=a_1\cdot (x_0y_1)^2(x_0y_2) + a_2\cdot (x_0y_2)^2(x_0y_3) + a_3\cdot (x_0y_3)^2(x_0y_1) + \\
+b_1\cdot (x_1y_1)^2(x_0y_2) + b_2\cdot (x_1y_2)^2(x_0y_3) + b_3\cdot (x_1y_3)^2(x_0y_1) + F',
\end{multline*}
or
\begin{multline*}
F=a_1\cdot (x_0y_1)^2(x_0y_2) + a_2\cdot (x_0y_2)^2(x_0y_3) + a_3\cdot (x_0y_3)^2(x_0y_1) + \\
+b_1\cdot (x_1y_1)^2(x_0y_3) + b_2\cdot (x_1y_2)^2(x_0y_1) + b_3\cdot (x_1y_3)^2(x_0y_2) + F',
\end{multline*}
where $F'$ does not contain the singled out monomials and $a_i,b_i\in {\mathbb C}^{*}$.\\

Let us find $G_0=ker(\pi)$. If $f_0=(1,{\om}_0^{c_1})\times 1\in G_0$, then $3c_1=0$ in the first and the second cases and $2c_1=0$ in the other three cases.\\

Hence in the first and the second cases $G_0\cong \ZZ{3}$ with a generator $f_0=(1,{\om}_0)\times 1$, where ${\om}_0=\sqr{3}$.\\

In the other three cases $G_0\cong \ZZ{2}$ with a generator $f_0=(1,{\om}_0)\times 1$, where ${\om}_0=\sqr{2}$.\\

Let $f_1=(1,{\om}^{c_1})\times P_3$ and $f_2=(1,{\om}^{d_1})\times W_3$ be elements of $G$.\\

Then 
\begin{multline*}
f_1(F)=a_1\cdot (x_0y_2)^2(x_0y_3) + a_2\cdot (x_0y_3)^2(x_0y_1) + a_3\cdot (x_0y_1)^2(x_0y_2) + \\
+{\om}^{3c_1}\cdot (b_1\cdot (x_1y_2)^2(x_1y_3) + b_2\cdot (x_1y_3)^2(x_1y_1) + b_3\cdot (x_1y_1)^2(x_1y_2)) + f_1(F'),
\end{multline*}
or
\begin{multline*}
f_1(F)=a_1\cdot (x_0y_2)^2(x_0y_3) + a_2\cdot (x_0y_3)^2(x_0y_1) + a_3\cdot (x_0y_1)^2(x_0y_2) + \\
+{\om}^{3c_1}\cdot (b_1\cdot (x_1y_2)^2(x_1y_1) + b_2\cdot (x_1y_3)^2(x_1y_2) + b_3\cdot (x_1y_1)^2(x_1y_3)) + f_1(F'),
\end{multline*}
or
\begin{multline*}
f_1(F)=a_1\cdot (x_0y_2)^2(x_0y_3) + a_2\cdot (x_0y_3)^2(x_0y_1) + a_3\cdot (x_0y_1)^2(x_0y_2) + \\
+{\om}^{2c_1}\cdot  ( b_1\cdot (x_1y_2)^2(x_0y_2) + b_2\cdot (x_1y_3)^2(x_0y_3) + b_3\cdot (x_1y_1)^2(x_0y_1) ) + f_1(F'),
\end{multline*}
or
\begin{multline*}
f_1(F)=a_1\cdot (x_0y_2)^2(x_0y_3) + a_2\cdot (x_0y_3)^2(x_0y_1) + a_3\cdot (x_0y_1)^2(x_0y_2) + \\
+{\om}^{2c_1}\cdot  ( b_1\cdot (x_1y_2)^2(x_0y_3) + b_2\cdot (x_1y_3)^2(x_0y_1) + b_3\cdot (x_1y_1)^2(x_0y_2) ) + f_1(F'),
\end{multline*}
or
\begin{multline*}
f_1(F)=a_1\cdot (x_0y_2)^2(x_0y_3) + a_2\cdot (x_0y_3)^2(x_0y_1) + a_3\cdot (x_0y_1)^2(x_0y_2) + \\
+{\om}^{2c_1}\cdot  ( b_1\cdot (x_1y_2)^2(x_0y_1) + b_2\cdot (x_1y_3)^2(x_0y_2) + b_3\cdot (x_1y_1)^2(x_0y_3) )+ f_1(F')
\end{multline*}

respectively and

\begin{multline*}
f_2(F)={\om}_3\cdot ( a_1\cdot (x_0y_1)^2(x_0y_2) + a_2\cdot (x_0y_2)^2(x_0y_3) + a_3\cdot (x_0y_3)^2(x_0y_1) ) + \\
+{\om}_3\cdot {\om}^{3d_1}\cdot ( b_1\cdot (x_1y_1)^2(x_1y_2) + b_2\cdot (x_1y_2)^2(x_1y_3) + b_3\cdot (x_1y_3)^2(x_1y_1) ) + f_2(F'),
\end{multline*}
or
\begin{multline*}
f_2(F)={\om}_3\cdot ( a_1\cdot (x_0y_1)^2(x_0y_2) + a_2\cdot (x_0y_2)^2(x_0y_3) + a_3\cdot (x_0y_3)^2(x_0y_1) ) + \\
+{\om}_3^2 \cdot {\om}^{3d_1}\cdot ( b_1\cdot (x_1y_1)^2(x_1y_3) + b_2\cdot (x_1y_2)^2(x_1y_1) + b_3\cdot (x_1y_3)^2(x_1y_2) ) + f_2(F'),
\end{multline*}
or
\begin{multline*}
f_2(F)={\om}_3\cdot ( a_1\cdot (x_0y_1)^2(x_0y_2) + a_2\cdot (x_0y_2)^2(x_0y_3) + a_3\cdot (x_0y_3)^2(x_0y_1) ) + \\
+{\om}^{2d_1}\cdot ( b_1\cdot (x_1y_1)^2(x_0y_1) + b_2\cdot (x_1y_2)^2(x_0y_2) + b_3\cdot (x_1y_3)^2(x_0y_3) ) + f_2(F'),
\end{multline*}
or
\begin{multline*}
f_2(F)={\om}_3\cdot ( a_1\cdot (x_0y_1)^2(x_0y_2) + a_2\cdot (x_0y_2)^2(x_0y_3) + a_3\cdot (x_0y_3)^2(x_0y_1) ) + \\
+{\om}_3\cdot {\om}^{2d_1}\cdot ( b_1\cdot (x_1y_1)^2(x_0y_2) + b_2\cdot (x_1y_2)^2(x_0y_3) + b_3\cdot (x_1y_3)^2(x_0y_1) ) + f_2(F'),
\end{multline*}
or
\begin{multline*}
f_2(F)={\om}_3\cdot ( a_1\cdot (x_0y_1)^2(x_0y_2) + a_2\cdot (x_0y_2)^2(x_0y_3) + a_3\cdot (x_0y_3)^2(x_0y_1) ) + \\
+{\om}_3^{2}\cdot {\om}^{2d_1}\cdot ( b_1\cdot (x_1y_1)^2(x_0y_3) + b_2\cdot (x_1y_2)^2(x_0y_1) + b_3\cdot (x_1y_3)^2(x_0y_2) ) + f_2(F')
\end{multline*}
respectively.\\

Since $f_1(F)\in {\mathbb C}^{*}\cdot F$, this implies that $a_3=\lambda a_1$, $a_2=(\lambda )^2 a_1$ and
\begin{itemize}
\item $b_3=\lambda \cdot {\om}^{-3c_1} \cdot b_1$, $b_2=(\lambda )^2 \cdot {\om}^{-6c_1} \cdot b_1$, ${\om}^{9c_1}=1$ in the first and the second cases,
\item $b_3=\lambda \cdot {\om}^{-2c_1} \cdot b_1$, $b_2=(\lambda )^2 \cdot {\om}^{-4c_1} \cdot b_1$, ${\om}^{6c_1}=1$ in the other three cases
\end{itemize}

for some $\lambda\in \mathbb C$ such that ${\lambda }^3=1$.\\

Since $f_2(F)\in {\mathbb C}^{*}\cdot F$, this implies that $f_2(F)={\om}_3\cdot F$ and so
\begin{itemize}
\item ${\om}^{3d_1}=1$,
\item ${\om}^{3d_1}={\om}_3^2$,
\item ${\om}^{2d_1}={\om}_3$,
\item ${\om}^{2d_1}=1$,
\item ${\om}^{2d_1}={\om}_3^2$
\end{itemize}
respectively in our five cases.\\

Since we are interested in $f_1, f_2\in G$ only modulo $G_0$ we can take
\begin{itemize}
\item $f_1=(1,{\om}^i)\times P_3$, $f_2=(1,1)\times W_3$ with some $i\in \{ 0,1,2 \}$ and $\om = \sqr{9}$,
\item $f_1=(1,{\om}^i)\times P_3$, $f_2=(1,{\om}^j)\times W_3$ with some $i,j\in \{ 0,1,2 \}$, $j\equiv 2\; mod \; 3$ and $\om = \sqr{9}$,
\item $f_1=(1,{\om}^i)\times P_3$, $f_2=(1,{\om}^j)\times W_3$ with some $i,j\in \{ 0,1,2 \}$, $j\equiv 1\; mod \; 3$ and $\om = \sqr{6}$,
\item $f_1=(1,{\om}^i)\times P_3$, $f_2=(1,1)\times W_3$ with some $i\in \{ 0,1,2 \}$ and $\om = \sqr{6}$,
\item $f_1=(1,{\om}^i)\times P_3$, $f_2=(1,{\om}^j)\times W_3$ with some $i,j\in \{ 0,1,2 \}$, $j\equiv 2\; mod \; 3$ and $\om = \sqr{6}$
\end{itemize}
respectively in our five cases.\\

Note that in the second case $f_2^6=f_0$, in the third case $f_2^3=f_0$. If $i=1$, then in all cases $f_1^3=f_0$. If $i=2$, then in the first and the second cases $f_1^6=f_0$.\\ 

We conclude that if $i=0$, then 
\begin{itemize}
\item $G\cong (\ZZ{3})^{\oplus 3}$ with generators $f_0=(1,{\om}_0)\times 1$, $f_1=(1,1)\times P_3$, $f_2=(1,1)\times W_3$, where ${\om}_0=\sqr{3}, \om = \sqr{9}$,
\item $G\cong \ZZ{9}\oplus \ZZ{3}$ with generators $f_1=(1,1)\times P_3$, $f_2=(1,{\om}^2)\times W_3$, where ${\om}_0=\sqr{3}, \om = \sqr{9}$,
\item $G\cong \ZZ{6}\oplus \ZZ{3}$ with generators $f_1=(1,1)\times P_3$, $f_2=(1,{\om})\times W_3$, where ${\om}_0=\sqr{2}, \om = \sqr{6}$,
\item $G\cong \ZZ{2}\oplus (\ZZ{3})^{\oplus 2}$ with generators $f_0=(1,{\om}_0)\times 1$, $f_1=(1,1)\times P_3$, $f_2=(1,1)\times W_3$, where ${\om}_0=\sqr{2}, \om = \sqr{6}$,
\item $G\cong \ZZ{2}\oplus (\ZZ{3})^{\oplus 2}$ with generators $f_0=(1,{\om}_0)\times 1$, $f_1=(1,1)\times P_3$, $f_2=(1,{\om}^2)\times W_3$, where ${\om}_0=\sqr{2}, \om = \sqr{6}$
\end{itemize}
respectively in our five cases.\\

If $i=1$, then 
\begin{itemize}
\item $G\cong \ZZ{9}\oplus \ZZ{3}$ with generators $f_1=(1,{\om})\times P_3$, $f_2=(1,1)\times W_3$, where $\om = \sqr{9}$,
\item $G\cong \ZZ{9}\oplus \ZZ{3}$ with generators $f_1=(1,{\om})\times P_3$, $f_1f_2=(1,{\om}^3)\times P_3W_3$, where $\om = \sqr{9}$,
\item $G\cong \ZZ{6}\oplus \ZZ{3}$ with generators $f_1=(1,{\om})\times P_3$, $f_1f_2=(1,{\om}^2)\times P_3W_3$, where $\om = \sqr{6}$,
\item $G\cong \ZZ{6}\oplus \ZZ{3}$ with generators $f_1=(1,{\om})\times P_3$, $f_2=(1,1)\times W_3$, where $\om = \sqr{6}$,
\item $G\cong \ZZ{6}\oplus \ZZ{3}$ with generators $f_1=(1,{\om})\times P_3$, $f_2=(1,{\om}^2)\times W_3$, where $\om = \sqr{6}$
\end{itemize}
respectively in our five cases.\\

If $i=2$, then 
\begin{itemize}
\item $G\cong \ZZ{9}\oplus \ZZ{3}$ with generators $f_1=(1,{\om}^2)\times P_3$, $f_2=(1,1)\times W_3$, where ${\om}_0=\sqr{3}, \om = \sqr{9}$,
\item $G\cong \ZZ{9}\oplus \ZZ{3}$ with generators $f_1=(1,{\om}^2)\times P_3$, $f_1^{-1}f_2=(1,1)\times (P_3)^{2}W_3$, where ${\om}_0=\sqr{3}, \om = \sqr{9}$,
\item $G\cong \ZZ{6}\oplus \ZZ{3}$ with generators $f_1=(1,{\om}^2)\times P_3$, $f_2=(1,{\om})\times W_3$, where ${\om}_0=\sqr{2}, \om = \sqr{6}$,
\item $G\cong \ZZ{2}\oplus (\ZZ{3})^{\oplus 2}$ with generators $f_0=(1,{\om}_0)\times 1$, $f_1=(1,{\om}^2)\times P_3$, $f_2=(1,1)\times W_3$, where ${\om}_0=\sqr{2}, \om = \sqr{6}$,
\item $G\cong \ZZ{2}\oplus (\ZZ{3})^{\oplus 2}$ with generators $f_0=(1,{\om}_0)\times 1$, $f_1=(1,{\om}^2)\times P_3$, $f_2=(1,{\om}^2)\times W_3$, where ${\om}_0=\sqr{2}, \om = \sqr{6}$
\end{itemize}
respectively in our five cases.\\

In order to describe all smooth cubic fourfolds which admit these group actions, we need to determine which cubic forms $F'$ can appear.\\

Since $F_x$ contains the monomial $x_0^3$ by our assumption, the invariance under $f_0$ requires that $F_x$ contains only monomials $x_0^3, x_1^3$ in the first and the second cases and only monomials $x_0^3, x_1^2x_0$ in the other three cases.\\

Let $F'_0=F'_0(z_0,z_1,z_2,z_3,z_4,z_5)$ be a cubic form which is obtained from $F'$ by removing all cubic monomials in $z_0, z_2, z_4$.\\

Let us denote by $F'_y$ the set of cubic monomials $y_{i}y_{j}y_{k}$ (in variables $y_1$, $y_2$, $y_3$) such that a monomial $(x_{i'}y_{i})(x_{j'}y_{j})(x_{k'}y_{k})$ appears in $F'_0$ with a nonzero coefficient for some $i',j',k'$.\\

The invariance under $f_2$ requires that a monomial $y_1^{i_1}y_2^{i_2}y_3^{i_3}$ is contained in $F'_y$ only if 
\begin{itemize}
\item $i_2\equiv 1+i_3\; mod \; 3$ in the first and the fourth cases,
\item $i_3\equiv 1+i_2\; mod \; 3$ in the second and the fifth cases,
\item $i_3\equiv i_2\; mod \; 3$ in the third case.
\end{itemize}

Hence $F'$ may contain only monomials
\begin{itemize}
\item $(x_1y_1)(x_1y_2)(x_0y_3)$, $(x_1y_1)(x_0y_2)(x_1y_3)$, $(x_0y_1)(x_1y_2)(x_1y_3)$ in the third case,\\
\item $(x_1y_1)(x_1y_2)(x_0y_1)$, $(x_1y_2)(x_1y_3)(x_0y_2)$, $(x_1y_3)(x_1y_1)(x_0y_3)$ in the fourth case,\\
\item $(x_1y_1)(x_1y_3)(x_0y_1)$, $(x_1y_2)(x_1y_1)(x_0y_2)$, $(x_1y_3)(x_1y_2)(x_0y_3)$ in the fifth case.\\
\end{itemize}

In the first and the second cases $F'=0$.\\

We conclude that either
\begin{multline*}
F=a_1\cdot (x_0y_1)^2(x_0y_2) + a_2\cdot (x_0y_2)^2(x_0y_3) + a_3\cdot (x_0y_3)^2(x_0y_1) + \\
+b_1\cdot (x_1y_1)^2(x_1y_2) + b_2\cdot (x_1y_2)^2(x_1y_3) + b_3\cdot (x_1y_3)^2(x_1y_1),
\end{multline*}
or
\begin{multline*}
F=a_1\cdot (x_0y_1)^2(x_0y_2) + a_2\cdot (x_0y_2)^2(x_0y_3) + a_3\cdot (x_0y_3)^2(x_0y_1) + \\
+b_1\cdot (x_1y_1)^2(x_1y_3) + b_2\cdot (x_1y_2)^2(x_1y_1) + b_3\cdot (x_1y_3)^2(x_1y_2),
\end{multline*}
or
\begin{multline*}
F=a_1\cdot (x_0y_1)^2(x_0y_2) + a_2\cdot (x_0y_2)^2(x_0y_3) + a_3\cdot (x_0y_3)^2(x_0y_1) + \\
+b_1\cdot (x_1y_1)^2(x_0y_1) + b_2\cdot (x_1y_2)^2(x_0y_2) + b_3\cdot (x_1y_3)^2(x_0y_3) + \\
+c_1\cdot (x_1y_1)(x_1y_2)(x_0y_3) +c_2\cdot  (x_1y_1)(x_0y_2)(x_1y_3)  +c_3\cdot (x_0y_1)(x_1y_2)(x_1y_3),
\end{multline*}
or
\begin{multline*}
F=a_1\cdot (x_0y_1)^2(x_0y_2) + a_2\cdot (x_0y_2)^2(x_0y_3) + a_3\cdot (x_0y_3)^2(x_0y_1) + \\
+b_1\cdot (x_1y_1)^2(x_0y_2) + b_2\cdot (x_1y_2)^2(x_0y_3) + b_3\cdot (x_1y_3)^2(x_0y_1) + \\
+c_1\cdot (x_1y_1)(x_1y_2)(x_0y_1) +c_2\cdot  (x_1y_2)(x_1y_3)(x_0y_2)  +c_3\cdot (x_1y_3)(x_1y_1)(x_0y_3),
\end{multline*}
or   
\begin{multline*}
F=a_1\cdot (x_0y_1)^2(x_0y_2) + a_2\cdot (x_0y_2)^2(x_0y_3) + a_3\cdot (x_0y_3)^2(x_0y_1) + \\
+b_1\cdot (x_1y_1)^2(x_0y_3) + b_2\cdot (x_1y_2)^2(x_0y_1) + b_3\cdot (x_1y_3)^2(x_0y_2) + \\
+c_1\cdot (x_1y_1)(x_1y_3)(x_0y_1) +c_2\cdot  (x_1y_2)(x_1y_1)(x_0y_2)  +c_3\cdot (x_1y_3)(x_1y_2)(x_0y_3)
\end{multline*}

respectively in our five cases, where $a_i, b_i\in {\mathbb C}^{*}$, $(a_3/a_1)^3=1$, $a_2=a_3^2/a_1$ and

\begin{itemize}
\item $b_3={\om}^{-3i} \cdot a_3b_1/a_1$, $b_2={\om}^{3i} \cdot a_1b_1/a_3$ in the first and the second cases,
\item $b_3= {\om}^{-2i} \cdot a_3b_1/a_1$, $b_2={\om}^{2i} \cdot a_1b_1/a_3$, $c_3={\om}^{2i} \cdot a_1c_1/a_3$, $c_2={\om}^{-2i} \cdot a_3c_1/a_1$ in the third case,
\item $b_3= {\om}^{-2i} \cdot a_3b_1/a_1$, $b_2={\om}^{2i} \cdot a_1b_1/a_3$, $c_2={\om}^{2i} \cdot a_1c_1/a_3$, $c_3={\om}^{-2i} \cdot a_3c_1/a_1$ in the fourth and the fifth cases.
\end{itemize}

It is immediate in the first and in the second cases that cubic fourfolds given by the generic equations above are smooth. In order to check smoothness in the other cases let us take $a_1=1$. Then $(a_3)^3=1$, $a_2=(a_3)^2$ and

\begin{itemize}
\item $b_3= {\om}^{-2i} \cdot a_3b_1$, $b_2={\om}^{2i} \cdot a_3^2b_1$, $c_3={\om}^{2i} \cdot a_3^2c_1$, $c_2={\om}^{-2i} \cdot a_3c_1$ in the third case,
\item $b_3= {\om}^{-2i} \cdot a_3b_1$, $b_2={\om}^{2i} \cdot a_3^2b_1$, $c_2={\om}^{2i} \cdot a_3^2c_1$, $c_3={\om}^{-2i} \cdot a_3c_1$ in the fourth and the fifth cases.
\end{itemize}

This means that in the third case
\begin{multline*}
F= z_0^2z_2 + a_3^2\cdot z_2^2z_4 + a_3\cdot z_4^2z_0 + \\
+b_1\cdot ( z_1^2z_0 + {\om}_3^{i} \cdot a_3^2 \cdot z_3^2z_2 + {\om}_3^{-i} \cdot a_3 \cdot z_5^2z_4  ) + \\
+c_1\cdot ( z_1z_3z_4 +{\om}_3^{-i} \cdot a_3 \cdot  z_1z_2z_5  +{\om}_3^{i} \cdot a_3^2 \cdot z_0z_3z_5 ),
\end{multline*}

in the fourth case

\begin{multline*}
F=z_0^2z_2 + a_3^2\cdot z_2^2z_4 + a_3\cdot z_4^2z_0 + \\
+b_1\cdot ( z_1^2z_2 + {\om}_3^{i} \cdot a_3^2 \cdot z_3^2z_4 + {\om}_3^{-i} \cdot a_3 \cdot z_5^2z_0 ) + \\
+c_1\cdot ( z_0z_1z_3 +{\om}_3^{i} \cdot a_3^2 \cdot  z_2z_3z_5  +{\om}_3^{-i} \cdot a_3 \cdot z_1z_4z_5),
\end{multline*}

in the fifth case
  
\begin{multline*}
F=z_0^2z_2 + a_3^2\cdot z_2^2z_4 + a_3\cdot z_4^2z_0 + \\
+b_1\cdot ( z_1^2z_4 + {\om}_3^{i} \cdot a_3^2 \cdot z_3^2z_0 + {\om}_3^{-i} \cdot a_3 \cdot z_5^2z_2 ) + \\
+c_1\cdot ( z_0z_1z_5 +{\om}_3^{i} \cdot a_3^2 \cdot  z_1z_2z_3  +{\om}_3^{-i} \cdot a_3 \cdot z_3z_4z_5).
\end{multline*}

Let us notice that these three cubic forms are the same upto a permutation of $z_1, z_3, z_5$ (and rescaling of $b_1, c_1$). Hence it is enough to consider only the third case.\\

By rescaling variables $z_1$, $z_3$, $z_4$ and $z_5$ and parameters $b_1, c_1$ we can rewrite the cubic form from the third case as follows:
\begin{multline*}
F= z_0^2z_2 + z_2^2z_4 + z_4^2z_0 + \\
+b_1\cdot ( z_1^2z_0 + z_3^2z_2 + z_5^2z_4  ) + c_1\cdot ( z_1z_3z_4 + z_1z_2z_5  +(-1)^{i} \cdot z_0z_3z_5 ),
\end{multline*}

Now one can check using Macaulay 2 that the corresponding cubic fourfold is smooth, if one takes $b_1=1$, $c_1=2$ for both values of $(-1)^{i}$.\\

We conclude that cubic fourfolds given by the generic equations above are smooth (in all five cases).\\

\paragraph{Case A3. $F$ contains $(x_0y_i)^2(x_0y_j)$ for some $i,j$, $F_y$ consists of monomials $y_1^2y_3, y_3^2y_2$, $y_2^2y_1$. }

Our assumptions imply that either
\begin{multline*}
F=a_1\cdot (x_0y_1)^2(x_0y_3) + a_2\cdot (x_0y_2)^2(x_0y_1) + a_3\cdot (x_0y_3)^2(x_0y_2) + \\
+b_1\cdot (x_1y_1)^2(x_1y_2) + b_2\cdot (x_1y_2)^2(x_1y_3) + b_3\cdot (x_1y_3)^2(x_1y_1) + F',
\end{multline*}
or
\begin{multline*}
F=a_1\cdot (x_0y_1)^2(x_0y_3) + a_2\cdot (x_0y_2)^2(x_0y_1) + a_3\cdot (x_0y_3)^2(x_0y_2) + \\
+b_1\cdot (x_1y_1)^2(x_1y_3) + b_2\cdot (x_1y_2)^2(x_1y_1) + b_3\cdot (x_1y_3)^2(x_1y_2) + F',
\end{multline*}
or
\begin{multline*}
F=a_1\cdot (x_0y_1)^2(x_0y_3) + a_2\cdot (x_0y_2)^2(x_0y_1) + a_3\cdot (x_0y_3)^2(x_0y_2) + \\
+b_1\cdot (x_1y_1)^2(x_0y_1) + b_2\cdot (x_1y_2)^2(x_0y_2) + b_3\cdot (x_1y_3)^2(x_0y_3) + F',
\end{multline*}
or
\begin{multline*}
F=a_1\cdot (x_0y_1)^2(x_0y_3) + a_2\cdot (x_0y_2)^2(x_0y_1) + a_3\cdot (x_0y_3)^2(x_0y_2) + \\
+b_1\cdot (x_1y_1)^2(x_0y_2) + b_2\cdot (x_1y_2)^2(x_0y_3) + b_3\cdot (x_1y_3)^2(x_0y_1) + F',
\end{multline*}
or
\begin{multline*}
F=a_1\cdot (x_0y_1)^2(x_0y_3) + a_2\cdot (x_0y_2)^2(x_0y_1) + a_3\cdot (x_0y_3)^2(x_0y_2) + \\
+b_1\cdot (x_1y_1)^2(x_0y_3) + b_2\cdot (x_1y_2)^2(x_0y_1) + b_3\cdot (x_1y_3)^2(x_0y_2) + F',
\end{multline*}
where $F'$ does not contain the singled out monomials and $a_i,b_i\in {\mathbb C}^{*}$.\\

Let us find $G_0=ker(\pi)$. If $f_0=(1,{\om}_0^{c_1})\times 1\in G_0$, then $3c_1=0$ in the first and the second cases and $2c_1=0$ in the other three cases.\\

Hence in the first and the second cases $G_0\cong \ZZ{3}$ with a generator $f_0=(1,{\om}_0)\times 1$, where ${\om}_0=\sqr{3}$.\\

In the other three cases $G_0\cong \ZZ{2}$ with a generator $f_0=(1,{\om}_0)\times 1$, where ${\om}_0=\sqr{2}$.\\

Let $f_1=(1,{\om}^{c_1})\times P_3$ and $f_2=(1,{\om}^{d_1})\times W_3$ be elements of $G$.\\

Then 
\begin{multline*}
f_1(F)=a_1\cdot (x_0y_2)^2(x_0y_1) + a_2\cdot (x_0y_3)^2(x_0y_2) + a_3\cdot (x_0y_1)^2(x_0y_3) + \\
+{\om}^{3c_1}\cdot (b_1\cdot (x_1y_2)^2(x_1y_3) + b_2\cdot (x_1y_3)^2(x_1y_1) + b_3\cdot (x_1y_1)^2(x_1y_2)) + f_1(F'),
\end{multline*}
or
\begin{multline*}
f_1(F)=a_1\cdot (x_0y_2)^2(x_0y_1) + a_2\cdot (x_0y_3)^2(x_0y_2) + a_3\cdot (x_0y_1)^2(x_0y_3) + \\
+{\om}^{3c_1}\cdot (b_1\cdot (x_1y_2)^2(x_1y_1) + b_2\cdot (x_1y_3)^2(x_1y_2) + b_3\cdot (x_1y_1)^2(x_1y_3)) + f_1(F'),
\end{multline*}
or
\begin{multline*}
f_1(F)=a_1\cdot (x_0y_2)^2(x_0y_1) + a_2\cdot (x_0y_3)^2(x_0y_2) + a_3\cdot (x_0y_1)^2(x_0y_3) + \\
+{\om}^{2c_1}\cdot  ( b_1\cdot (x_1y_2)^2(x_0y_2) + b_2\cdot (x_1y_3)^2(x_0y_3) + b_3\cdot (x_1y_1)^2(x_0y_1) ) + f_1(F'),
\end{multline*}
or
\begin{multline*}
f_1(F)=a_1\cdot (x_0y_2)^2(x_0y_1) + a_2\cdot (x_0y_3)^2(x_0y_2) + a_3\cdot (x_0y_1)^2(x_0y_3) + \\
+{\om}^{2c_1}\cdot  ( b_1\cdot (x_1y_2)^2(x_0y_3) + b_2\cdot (x_1y_3)^2(x_0y_1) + b_3\cdot (x_1y_1)^2(x_0y_2) ) + f_1(F'),
\end{multline*}
or
\begin{multline*}
f_1(F)=a_1\cdot (x_0y_2)^2(x_0y_1) + a_2\cdot (x_0y_3)^2(x_0y_2) + a_3\cdot (x_0y_1)^2(x_0y_3) + \\
+{\om}^{2c_1}\cdot  ( b_1\cdot (x_1y_2)^2(x_0y_1) + b_2\cdot (x_1y_3)^2(x_0y_2) + b_3\cdot (x_1y_1)^2(x_0y_3) )+ f_1(F')
\end{multline*}

respectively and

\begin{multline*}
f_2(F)={\om}_3^2\cdot ( a_1\cdot (x_0y_1)^2(x_0y_3) + a_2\cdot (x_0y_2)^2(x_0y_1) + a_3\cdot (x_0y_3)^2(x_0y_2) ) + \\
+{\om}_3\cdot {\om}^{3d_1}\cdot ( b_1\cdot (x_1y_1)^2(x_1y_2) + b_2\cdot (x_1y_2)^2(x_1y_3) + b_3\cdot (x_1y_3)^2(x_1y_1) ) + f_2(F'),
\end{multline*}
or
\begin{multline*}
f_2(F)={\om}_3^2\cdot ( a_1\cdot (x_0y_1)^2(x_0y_3) + a_2\cdot (x_0y_2)^2(x_0y_1) + a_3\cdot (x_0y_3)^2(x_0y_2) ) + \\
+{\om}_3^2 \cdot {\om}^{3d_1}\cdot ( b_1\cdot (x_1y_1)^2(x_1y_3) + b_2\cdot (x_1y_2)^2(x_1y_1) + b_3\cdot (x_1y_3)^2(x_1y_2) ) + f_2(F'),
\end{multline*}
or
\begin{multline*}
f_2(F)={\om}_3^2\cdot ( a_1\cdot (x_0y_1)^2(x_0y_3) + a_2\cdot (x_0y_2)^2(x_0y_1) + a_3\cdot (x_0y_3)^2(x_0y_2) ) + \\
+{\om}^{2d_1}\cdot ( b_1\cdot (x_1y_1)^2(x_0y_1) + b_2\cdot (x_1y_2)^2(x_0y_2) + b_3\cdot (x_1y_3)^2(x_0y_3) ) + f_2(F'),
\end{multline*}
or
\begin{multline*}
f_2(F)={\om}_3^2\cdot ( a_1\cdot (x_0y_1)^2(x_0y_3) + a_2\cdot (x_0y_2)^2(x_0y_1) + a_3\cdot (x_0y_3)^2(x_0y_2) ) + \\
+{\om}_3\cdot {\om}^{2d_1}\cdot ( b_1\cdot (x_1y_1)^2(x_0y_2) + b_2\cdot (x_1y_2)^2(x_0y_3) + b_3\cdot (x_1y_3)^2(x_0y_1) ) + f_2(F'),
\end{multline*}
or
\begin{multline*}
f_2(F)={\om}_3^2\cdot ( a_1\cdot (x_0y_1)^2(x_0y_3) + a_2\cdot (x_0y_2)^2(x_0y_1) + a_3\cdot (x_0y_3)^2(x_0y_2) ) + \\
+{\om}_3^{2}\cdot {\om}^{2d_1}\cdot ( b_1\cdot (x_1y_1)^2(x_0y_3) + b_2\cdot (x_1y_2)^2(x_0y_1) + b_3\cdot (x_1y_3)^2(x_0y_2) ) + f_2(F')
\end{multline*}
respectively.\\

Since $f_1(F)\in {\mathbb C}^{*}\cdot F$, this implies that $a_3=\lambda a_1$, $a_2=(\lambda )^2 a_1$ and
\begin{itemize}
\item $b_3=\lambda \cdot {\om}^{-3c_1} \cdot b_1$, $b_2=(\lambda )^2 \cdot {\om}^{-6c_1} \cdot b_1$, ${\om}^{9c_1}=1$ in the first and the second cases,
\item $b_3=\lambda \cdot {\om}^{-2c_1} \cdot b_1$, $b_2=(\lambda )^2 \cdot {\om}^{-4c_1} \cdot b_1$, ${\om}^{6c_1}=1$ in the other three cases
\end{itemize}

for some $\lambda\in \mathbb C$ such that ${\lambda }^3=1$.\\

Since $f_2(F)\in {\mathbb C}^{*}\cdot F$, this implies that $f_2(F)={\om}^2_3\cdot F$ and so
\begin{itemize}
\item ${\om}^{3d_1}={\om}_3$,
\item ${\om}^{3d_1}=1$,
\item ${\om}^{2d_1}={\om}_3^2$,
\item ${\om}^{2d_1}={\om}_3$,
\item ${\om}^{2d_1}=1$
\end{itemize}
respectively in our five cases.\\

Since we are interested in $f_1, f_2\in G$ only modulo $G_0$ we can take
\begin{itemize}
\item $f_1=(1,{\om}^i)\times P_3$, $f_2=(1,{\om}^j)\times W_3$ with some $i,j\in \{ 0,1,2 \}$, $j\equiv 1\; mod \; 3$ and $\om = \sqr{9}$,
\item $f_1=(1,{\om}^i)\times P_3$, $f_2=(1,1)\times W_3$ with some $i\in \{ 0,1,2 \}$ and $\om = \sqr{9}$,
\item $f_1=(1,{\om}^i)\times P_3$, $f_2=(1,{\om}^j)\times W_3$ with some $i,j\in \{ 0,1,2 \}$, $j\equiv 2\; mod \; 3$ and $\om = \sqr{6}$,
\item $f_1=(1,{\om}^i)\times P_3$, $f_2=(1,{\om}^j)\times W_3$ with some $i,j\in \{ 0,1,2 \}$, $j\equiv 1\; mod \; 3$ and $\om = \sqr{6}$,
\item $f_1=(1,{\om}^i)\times P_3$, $f_2=(1,1)\times W_3$ with some $i\in \{ 0,1,2 \}$ and $\om = \sqr{6}$
\end{itemize}
respectively in our five cases.\\

Note that in the first and the fourth cases $f_2^3=f_0$. If $i=1$, then in all cases $f_1^3=f_0$. If $i=2$, then in the first and the second cases $f_1^6=f_0$.\\ 

We conclude that if $i=0$, then
\begin{itemize}
\item $G\cong \ZZ{9}\oplus \ZZ{3}$ with generators $f_1=(1,1)\times P_3$, $f_2=(1,{\om})\times W_3$, where ${\om}_0=\sqr{3}, \om = \sqr{9}$,
\item $G\cong (\ZZ{3})^{\oplus 3}$ with generators $f_0=(1,{\om}_0)\times 1$, $f_1=(1,1)\times P_3$, $f_2=(1,1)\times W_3$, where ${\om}_0=\sqr{3}, \om = \sqr{9}$,
\item $G\cong \ZZ{2}\oplus (\ZZ{3})^{\oplus 2}$ with generators $f_0=(1,{\om}_0)\times 1$, $f_1=(1,1)\times P_3$, $f_2=(1,{\om}^2)\times W_3$, where ${\om}_0=\sqr{2}, \om = \sqr{6}$,
\item $G\cong \ZZ{6}\oplus \ZZ{3}$ with generators $f_1=(1,1)\times P_3$, $f_2=(1,{\om})\times W_3$, where ${\om}_0=\sqr{2}, \om = \sqr{6}$,
\item $G\cong \ZZ{2}\oplus (\ZZ{3})^{\oplus 2}$ with generators $f_0=(1,{\om}_0)\times 1$, $f_1=(1,1)\times P_3$, $f_2=(1,1)\times W_3$, where ${\om}_0=\sqr{2}, \om = \sqr{6}$
\end{itemize}
respectively in our five cases.\\

If $i=1$, then 
\begin{itemize}
\item $G\cong \ZZ{9}\oplus \ZZ{3}$ with generators $f_1=(1,{\om})\times P_3$, $f_1^2f_2=(1,{\om}^3)\times (P_3)^2W_3$, where $\om = \sqr{9}$,
\item $G\cong \ZZ{9}\oplus \ZZ{3}$ with generators $f_1=(1,{\om})\times P_3$, $f_2=(1,1)\times W_3$, where $\om = \sqr{9}$,
\item $G\cong \ZZ{6}\oplus \ZZ{3}$ with generators $f_1=(1,{\om})\times P_3$, $f_2=(1,{\om}^2)\times W_3$, where $\om = \sqr{6}$,
\item $G\cong \ZZ{6}\oplus \ZZ{3}$ with generators $f_1=(1,{\om})\times P_3$, $f_1f_2=(1,{\om}^2)\times P_3W_3$, where $\om = \sqr{6}$,
\item $G\cong \ZZ{6}\oplus \ZZ{3}$ with generators $f_1=(1,{\om})\times P_3$, $f_2=(1,1)\times W_3$, where $\om = \sqr{6}$
\end{itemize}
respectively in our five cases.\\

If $i=2$, then 
\begin{itemize}
\item $G\cong \ZZ{9}\oplus \ZZ{3}$ with generators $f_1=(1,{\om}^2)\times P_3$, $f_1f_2=(1,{\om}^3)\times P_3W_3$, where ${\om}_0=\sqr{3}, \om = \sqr{9}$,
\item $G\cong \ZZ{9}\oplus \ZZ{3}$ with generators $f_1=(1,{\om}^2)\times P_3$, $f_2=(1,1)\times W_3$, where ${\om}_0=\sqr{3}, \om = \sqr{9}$,
\item $G\cong \ZZ{2}\oplus (\ZZ{3})^{\oplus 2}$ with generators $f_0=(1,{\om}_0)\times 1$, $f_1=(1,{\om}^2)\times P_3$, $f_2=(1,{\om}^2)\times W_3$, where ${\om}_0=\sqr{2}, \om = \sqr{6}$,
\item $G\cong \ZZ{6}\oplus \ZZ{3}$ with generators $f_1=(1,{\om}^2)\times P_3$, $f_2=(1,{\om})\times W_3$, where ${\om}_0=\sqr{2}, \om = \sqr{6}$,
\item $G\cong \ZZ{2}\oplus (\ZZ{3})^{\oplus 2}$ with generators $f_0=(1,{\om}_0)\times 1$, $f_1=(1,{\om}^2)\times P_3$, $f_2=(1,1)\times W_3$, where ${\om}_0=\sqr{2}, \om = \sqr{6}$
\end{itemize}
respectively in our five cases.\\

In order to describe all smooth cubic fourfolds which admit these group actions, we need to determine which cubic forms $F'$ can appear.\\

Since $F_x$ contains the monomial $x_0^3$ by our assumption, the invariance under $f_0$ requires that $F_x$ contains only monomials $x_0^3, x_1^3$ in the first and the second cases and only monomials $x_0^3, x_1^2x_0$ in the other three cases.\\

Let $F'_0=F'_0(z_0,z_1,z_2,z_3,z_4,z_5)$ be a cubic form which is obtained from $F'$ by removing all cubic monomials in $z_0, z_2, z_4$.\\

Let us denote by $F'_y$ the set of cubic monomials $y_{i}y_{j}y_{k}$ (in variables $y_1$, $y_2$, $y_3$) such that a monomial $(x_{i'}y_{i})(x_{j'}y_{j})(x_{k'}y_{k})$ appears in $F'_0$ with a nonzero coefficient for some $i',j',k'$.\\

The invariance under $f_2$ requires that a monomial $y_1^{i_1}y_2^{i_2}y_3^{i_3}$ is contained in $F'_y$ only if $j+i_2+2i_3\equiv 2\; mod \; 3$, i.e.
\begin{itemize}
\item $i_2\equiv 1+i_3\; mod \; 3$ in the first and the fourth cases,
\item $i_3\equiv 1+i_2\; mod \; 3$ in the second and the fifth cases ($j=0$),
\item $i_3\equiv i_2\; mod \; 3$ in the third case.
\end{itemize}

Hence $F'$ may contain only monomials
\begin{itemize}
\item $(x_1y_1)(x_1y_2)(x_0y_3)$, $(x_1y_1)(x_0y_2)(x_1y_3)$, $(x_0y_1)(x_1y_2)(x_1y_3)$ in the third case,\\
\item $(x_1y_1)(x_1y_2)(x_0y_1)$, $(x_1y_2)(x_1y_3)(x_0y_2)$, $(x_1y_3)(x_1y_1)(x_0y_3)$ in the fourth case,\\
\item $(x_1y_1)(x_1y_3)(x_0y_1)$, $(x_1y_2)(x_1y_1)(x_0y_2)$, $(x_1y_3)(x_1y_2)(x_0y_3)$ in the fifth case.\\
\end{itemize}

In the first and the second cases $F'=0$.\\

We conclude that either
\begin{multline*}
F=a_1\cdot (x_0y_1)^2(x_0y_3) + a_2\cdot (x_0y_2)^2(x_0y_1) + a_3\cdot (x_0y_3)^2(x_0y_2) + \\
+b_1\cdot (x_1y_1)^2(x_1y_2) + b_2\cdot (x_1y_2)^2(x_1y_3) + b_3\cdot (x_1y_3)^2(x_1y_1),
\end{multline*}
or
\begin{multline*}
F=a_1\cdot (x_0y_1)^2(x_0y_3) + a_2\cdot (x_0y_2)^2(x_0y_1) + a_3\cdot (x_0y_3)^2(x_0y_2) + \\
+b_1\cdot (x_1y_1)^2(x_1y_3) + b_2\cdot (x_1y_2)^2(x_1y_1) + b_3\cdot (x_1y_3)^2(x_1y_2),
\end{multline*}
or
\begin{multline*}
F=a_1\cdot (x_0y_1)^2(x_0y_3) + a_2\cdot (x_0y_2)^2(x_0y_1) + a_3\cdot (x_0y_3)^2(x_0y_2) + \\
+b_1\cdot (x_1y_1)^2(x_0y_1) + b_2\cdot (x_1y_2)^2(x_0y_2) + b_3\cdot (x_1y_3)^2(x_0y_3) + \\
+c_1\cdot (x_1y_1)(x_1y_2)(x_0y_3) +c_2\cdot  (x_1y_1)(x_0y_2)(x_1y_3)  +c_3\cdot (x_0y_1)(x_1y_2)(x_1y_3),
\end{multline*}
or
\begin{multline*}
F=a_1\cdot (x_0y_1)^2(x_0y_3) + a_2\cdot (x_0y_2)^2(x_0y_1) + a_3\cdot (x_0y_3)^2(x_0y_2) + \\
+b_1\cdot (x_1y_1)^2(x_0y_2) + b_2\cdot (x_1y_2)^2(x_0y_3) + b_3\cdot (x_1y_3)^2(x_0y_1) + \\
+c_1\cdot (x_1y_1)(x_1y_2)(x_0y_1) +c_2\cdot  (x_1y_2)(x_1y_3)(x_0y_2)  +c_3\cdot (x_1y_3)(x_1y_1)(x_0y_3),
\end{multline*}
or   
\begin{multline*}
F=a_1\cdot (x_0y_1)^2(x_0y_3) + a_2\cdot (x_0y_2)^2(x_0y_1) + a_3\cdot (x_0y_3)^2(x_0y_2) + \\
+b_1\cdot (x_1y_1)^2(x_0y_3) + b_2\cdot (x_1y_2)^2(x_0y_1) + b_3\cdot (x_1y_3)^2(x_0y_2) + \\
+c_1\cdot (x_1y_1)(x_1y_3)(x_0y_1) +c_2\cdot  (x_1y_2)(x_1y_1)(x_0y_2)  +c_3\cdot (x_1y_3)(x_1y_2)(x_0y_3),
\end{multline*}

respectively in our five cases, where $a_i, b_i\in {\mathbb C}^{*}$, $(a_3/a_1)^3=1$, $a_2=a_3^2/a_1$ and

\begin{itemize}
\item $b_3={\om}^{-3i} \cdot a_3b_1/a_1$, $b_2={\om}^{3i} \cdot a_1b_1/a_3$ in the first and the second cases,
\item $b_3= {\om}^{-2i} \cdot a_3b_1/a_1$, $b_2={\om}^{2i} \cdot a_1b_1/a_3$, $c_3={\om}^{2i} \cdot a_1c_1/a_3$, $c_2={\om}^{-2i} \cdot a_3c_1/a_1$ in the third case,
\item $b_3= {\om}^{-2i} \cdot a_3b_1/a_1$, $b_2={\om}^{2i} \cdot a_1b_1/a_3$, $c_2={\om}^{2i} \cdot a_1c_1/a_3$, $c_3={\om}^{-2i} \cdot a_3c_1/a_1$ in the fourth and the fifth cases.
\end{itemize}

It is immediate in the first and in the second cases that cubic fourfolds given by the generic equations above are smooth. In order to check smoothness in the other cases let us take $a_1=1$. Then $(a_3)^3=1$, $a_2=(a_3)^2$ and

\begin{itemize}
\item $b_3= {\om}^{-2i} \cdot a_3b_1$, $b_2={\om}^{2i} \cdot a_3^2b_1$, $c_3={\om}^{2i} \cdot a_3^2c_1$, $c_2={\om}^{-2i} \cdot a_3c_1$ in the third case,
\item $b_3= {\om}^{-2i} \cdot a_3b_1$, $b_2={\om}^{2i} \cdot a_3^2b_1$, $c_2={\om}^{2i} \cdot a_3^2c_1$, $c_3={\om}^{-2i} \cdot a_3c_1$ in the fourth and the fifth cases.
\end{itemize}

This means that in the third case
\begin{multline*}
F= z_0^2z_4 + a_3^2\cdot z_2^2z_0 + a_3\cdot z_4^2z_2 + \\
+b_1\cdot ( z_1^2z_0 + {\om}_3^{i} \cdot a_3^2 \cdot z_3^2z_2 + {\om}_3^{-i} \cdot a_3 \cdot z_5^2z_4  ) + \\
+c_1\cdot ( z_1z_3z_4 +{\om}_3^{-i} \cdot a_3 \cdot  z_1z_2z_5  +{\om}_3^{i} \cdot a_3^2 \cdot z_0z_3z_5 ),
\end{multline*}

in the fourth case

\begin{multline*}
F= z_0^2z_4 + a_3^2\cdot z_2^2z_0 + a_3\cdot z_4^2z_2 + \\
+b_1\cdot ( z_1^2z_2 + {\om}_3^{i} \cdot a_3^2 \cdot z_3^2z_4 + {\om}_3^{-i} \cdot a_3 \cdot z_5^2z_0 ) + \\
+c_1\cdot ( z_0z_1z_3 +{\om}_3^{i} \cdot a_3^2 \cdot  z_2z_3z_5  +{\om}_3^{-i} \cdot a_3 \cdot z_1z_4z_5),
\end{multline*}

in the fifth case
  
\begin{multline*}
F= z_0^2z_4 + a_3^2\cdot z_2^2z_0 + a_3\cdot z_4^2z_2 + \\
+b_1\cdot ( z_1^2z_4 + {\om}_3^{i} \cdot a_3^2 \cdot z_3^2z_0 + {\om}_3^{-i} \cdot a_3 \cdot z_5^2z_2 ) + \\
+c_1\cdot ( z_0z_1z_5 +{\om}_3^{i} \cdot a_3^2 \cdot  z_1z_2z_3  +{\om}_3^{-i} \cdot a_3 \cdot z_3z_4z_5).
\end{multline*}

Let us notice that these three cubic forms are the same upto a permutation of $z_1, z_3, z_5$ (and rescaling of $b_1, c_1$). Hence it is enough to consider only the third case.\\

By rescaling variables $z_1$, $z_3$, $z_4$ and $z_5$ and parameters $b_1, c_1$ we can rewrite the cubic form from the third case as follows:
\begin{multline*}
F= z_0^2z_4 + z_2^2z_0 + z_4^2z_2 + \\
+b_1\cdot ( z_1^2z_0 + z_3^2z_2 + z_5^2z_4  ) + c_1\cdot ( z_1z_3z_4 + z_1z_2z_5  +(-1)^{i} \cdot z_0z_3z_5 ),
\end{multline*}

Now one can check using Macaulay 2 that the corresponding cubic fourfold is smooth, if one takes $b_1=1$, $c_1=2$ for both values of $(-1)^{i}$. In fact, essentially the same computation was already done in Case 2A.\\

We conclude that cubic fourfolds given by the generic equations above are smooth (in all five cases).\\

\subsubsection{Case B. $F$ does not contain $(x_0y_i)^2(x_0y_j)$, $(x_1y_i)^2(x_1y_j)$ for any $i,j$. }

Let $F_0$ be a cubic form obtained from $F$ by omitting all monomials containing $x_1^2x_0$, $x_0^3$ and $x_1^3$. It is nonzero by our assumption.\\

Let $F_y$ be the set of cubic monomials $y_{i}y_{j}y_{k}$ (in variables $y_1$, $y_2$, $y_3$) such that a monomial $(x_{i'}y_{i})(x_{j'}y_{j})(x_{k'}y_{k})$ appears in $F_0$ with a nonzero coefficient for some $i',j',k'$.\\

{\bf Lemma 6B.} {\it Suppose that $G$ acts effectively on $X$, $F_x$ contains $x_0^2x_1$ and $\pi$ is surjective. Then $F_y$ contains monomials only from exactly one of the following sets:
\begin{itemize}
\item $y_1^3, y_2^3, y_3^3, y_1y_2y_3$,
\item $y_1^2y_2, y_2^2y_3, y_3^2y_1$,
\item $y_1^2y_3, y_3^2y_2, y_2^2y_1$.
\end{itemize}}
\par
{\it Proof:} Let $y_1^{i_1}y_2^{i_2}y_3^{i_3}$ be a monomial in $F_y$.\\ 

Let $f=(1,{\om}^{c_1})\times W_3\in G$. Then the condition $f(F)\in {\mathbb C}^{*}\cdot F$ implies that the residue class $i_2+2i_3 \; mod \; 3$ is the same for all monomials in $F_y$.\\

If this residue class is $0\; mod \; 3$, then $i_2\equiv i_3\; mod \; 3$, which means exactly that $(i_1,i_2,i_3)\in \{ (3,0,0), (0,3,0), (0,0,3), (1,1,1) \}$.\\

If this residue class is $1\; mod \; 3$, then $i_2\equiv 1+i_3\; mod \; 3$, which means exactly that $(i_1,i_2,i_3)\in \{ (1,0,2), (0,2,1), (2,1,0) \}$.\\

If this residue class is $-1\; mod \; 3$, then $i_3\equiv 1+i_2\; mod \; 3$, which means exactly that $(i_1,i_2,i_3)\in \{ (1,2,0), (0,1,2), (2,0,1) \}$. {\it QED}\\

This Lemma implies that we need to consider three cases:
\begin{itemize}
\item[(B1)] $F_y$ is contained in the set of monomials $y_1^3, y_2^3, y_3^3, y_1y_2y_3$,
\item[(B2)] $F_y$ consists of monomials $y_1^2y_2, y_2^2y_3, y_3^2y_1$,
\item[(B3)] $F_y$ consists of monomials $y_1^2y_3, y_3^2y_2, y_2^2y_1$.
\end{itemize}

\paragraph{Case B1. $F$ does not contain $(x_0y_i)^2(x_0y_j)$, $(x_1y_i)^2(x_1y_j)$ for any $i,j$, $F_y$ is contained in the set of monomials $y_1^3, y_2^3, y_3^3, y_1y_2y_3$. }

Our assumptions imply that either
\begin{multline*}
F=a_1\cdot (x_0y_1)^2(x_1y_1) + a_2\cdot (x_0y_2)^2(x_1y_2) + a_3\cdot (x_0y_3)^2(x_1y_3) + \\
+b_1\cdot (x_1y_1)^2(x_0y_1) + b_2\cdot (x_1y_2)^2(x_0y_2) + b_3\cdot (x_1y_3)^2(x_0y_3) + F',
\end{multline*}
or
\begin{multline*}
F=a_1\cdot (x_0y_1)^2(x_1y_1) + a_2\cdot (x_0y_2)^2(x_1y_2) + a_3\cdot (x_0y_3)^2(x_1y_3) + \\
+b_1\cdot (x_1y_1)^2(x_0y_2) + b_2\cdot (x_1y_2)^2(x_0y_3) + b_3\cdot (x_1y_3)^2(x_0y_1) + F',
\end{multline*}
or
\begin{multline*}
F=a_1\cdot (x_0y_1)^2(x_1y_1) + a_2\cdot (x_0y_2)^2(x_1y_2) + a_3\cdot (x_0y_3)^2(x_1y_3) + \\
+b_1\cdot (x_1y_1)^2(x_0y_3) + b_2\cdot (x_1y_2)^2(x_0y_1) + b_3\cdot (x_1y_3)^2(x_0y_2) + F',
\end{multline*}
where $F'$ does not contain the singled out monomials and $a_i,b_i\in {\mathbb C}^{*}$.\\

Let us find $G_0=ker(\pi)$. If $f_0=(1,{\om}_0^{c_1})\times 1\in G_0$, then $c_1=0$. This means that $G_0=0$. Hence $\pi\colon G\rightarrow {\mathcal P}_3\cong \ZZ{3}\oplus \ZZ{3}$ is isomorphism (since we assume that it is surjective).\\

Let $f_1=(1,{\om}^{c_1})\times P_3$ and $f_2=(1,{\om}^{d_1})\times W_3$ be generators of $G$.\\

Then 
\begin{multline*}
f_1(F)={\om}^{c_1}\cdot ( a_1\cdot (x_0y_2)^2(x_1y_2) + a_2\cdot (x_0y_3)^2(x_1y_3) + a_3\cdot (x_0y_1)^2(x_1y_1) ) + \\
+{\om}^{2c_1}\cdot ( b_1\cdot (x_1y_2)^2(x_0y_2) + b_2\cdot (x_1y_3)^2(x_0y_3) + b_3\cdot (x_1y_1)^2(x_0y_1) ) + f_1(F'),
\end{multline*}
or
\begin{multline*}
f_1(F)={\om}^{c_1}\cdot ( a_1\cdot (x_0y_2)^2(x_1y_2) + a_2\cdot (x_0y_3)^2(x_1y_3) + a_3\cdot (x_0y_1)^2(x_1y_1) ) + \\
+{\om}^{2c_1}\cdot  ( b_1\cdot (x_1y_2)^2(x_0y_3) + b_2\cdot (x_1y_3)^2(x_0y_1) + b_3\cdot (x_1y_1)^2(x_0y_2) ) + f_1(F'),
\end{multline*}
or
\begin{multline*}
f_1(F)={\om}^{c_1}\cdot ( a_1\cdot (x_0y_2)^2(x_1y_2) + a_2\cdot (x_0y_3)^2(x_1y_3) + a_3\cdot (x_0y_1)^2(x_1y_1) ) + \\
+{\om}^{2c_1}\cdot  ( b_1\cdot (x_1y_2)^2(x_0y_1) + b_2\cdot (x_1y_3)^2(x_0y_2) + b_3\cdot (x_1y_1)^2(x_0y_3)) + f_1(F')
\end{multline*}

respectively and 

\begin{multline*}
f_2(F)={\om}^{d_1}\cdot  ( a_1\cdot (x_0y_1)^2(x_1y_1) + a_2\cdot (x_0y_2)^2(x_1y_2) + a_3\cdot (x_0y_3)^2(x_1y_3) ) + \\
+{\om}^{2d_1}\cdot ( b_1\cdot (x_1y_1)^2(x_0y_1) + b_2\cdot (x_1y_2)^2(x_0y_2) + b_3\cdot (x_1y_3)^2(x_0y_3) ) + f_2(F'),
\end{multline*}
or
\begin{multline*}
f_2(F)={\om}^{d_1}\cdot  ( a_1\cdot (x_0y_1)^2(x_1y_1) + a_2\cdot (x_0y_2)^2(x_1y_2) + a_3\cdot (x_0y_3)^2(x_1y_3) ) + \\
+{\om}^{2d_1}\cdot {\om}_3\cdot (b_1\cdot (x_1y_1)^2(x_0y_2) + b_2\cdot (x_1y_2)^2(x_0y_3) + b_3\cdot (x_1y_3)^2(x_0y_1)) + f_2(F'),
\end{multline*}
or
\begin{multline*}
f_2(F)={\om}^{d_1}\cdot  ( a_1\cdot (x_0y_1)^2(x_1y_1) + a_2\cdot (x_0y_2)^2(x_1y_2) + a_3\cdot (x_0y_3)^2(x_1y_3) ) + \\
+{\om}^{2d_1}\cdot {\om}_3^2\cdot (b_1\cdot (x_1y_1)^2(x_0y_3) + b_2\cdot (x_1y_2)^2(x_0y_1) + b_3\cdot (x_1y_3)^2(x_0y_2)) + f_2(F')
\end{multline*}

respectively.\\ 

Since $f_1(F)\in {\mathbb C}^{*}\cdot F$, this implies that 
$$
a_3={\om}^{-c_1}\cdot \lambda a_1,\; a_2={\om}^{-2c_1}\cdot (\lambda )^2 a_1,\; b_3={\om}^{-2c_1}\cdot \lambda b_1,\; b_2={\om}^{-4c_1}\cdot (\lambda )^2 b_1, \; {\om}^{3c_1}=1 
$$
for some $\lambda\in \mathbb C$ such that ${\lambda }^3=1$.\\

Since $f_2(F)\in {\mathbb C}^{*}\cdot F$, this implies that
\begin{itemize}
\item ${\om}^{d_1}=1$,
\item ${\om}^{d_1}={\om}_3^2$,
\item ${\om}^{d_1}={\om}_3$
\end{itemize}
respectively in our three cases.\\

This means that we can take $\om={\om}_3$ and 
\begin{itemize}
\item $f_1=(1,{\om}^i)\times P_3$, $f_2=(1,1)\times W_3$ with some $i\in \{ 0,1,2 \}$ and $\om = \sqr{3}$,
\item $f_1=(1,{\om}^i)\times P_3$, $f_2=(1,{\om}^2)\times W_3$ with some $i\in \{ 0,1,2 \}$ and $\om = \sqr{3}$,
\item $f_1=(1,{\om}^i)\times P_3$, $f_2=(1,{\om})\times W_3$ with some $i \in \{ 0,1,2 \}$ and $\om = \sqr{3}$
\end{itemize}
respectively in our three cases.\\

We conclude that $G\cong (\ZZ{3})^{\oplus 2}$ with generators $f_1=(1,{\om}^i)\times P_3$, $f_2=(1,{\om}^j)\times W_3$ with some $i,j\in \{ 0,1,2 \}$, where $\om = \sqr{3}$. Here $j=0$ corresponds to the first case, $j=2$ corresponds to the second case and $j=1$ corresponds to the third case.\\

In order to describe all smooth cubic fourfolds which admit these group actions, we need to determine which cubic forms $F'$ can appear.\\

By our assumption $x_0^3$ and $x_1^3$ may appear in $F$ only in monomials $(x_0y_1)(x_0y_2)(x_0y_3)$, $(x_1y_1)(x_1y_2)(x_1y_3)$. The invariance under $f_2$ requires that they are possible only in the first case. Moreover, if any of them appears in $F'$ with a nonzero coefficient, then $f_1(F)=F$.\\

Let $F'_0=F'_0(z_0,z_1,z_2,z_3,z_4,z_5)$ be a cubic form which is obtained from $F'$ by removing all monomials containing $x_0^3$, $x_0^2x_1$, $x_1^3$.\\

Let us denote by $F'_y$ the set of cubic monomials $y_{i}y_{j}y_{k}$ (in variables $y_1$, $y_2$, $y_3$) such that a monomial $(x_{i'}y_{i})(x_{j'}y_{j})(x_{k'}y_{k})$ appears in $F'_0$ with a nonzero coefficient for some $i',j',k'$.\\

The invariance under $f_2$ requires that $y_1^{i_1}y_2^{i_2}y_3^{i_3}$ is contained in $F'_y$ only if $i_2\equiv i_3\; mod \; 3$ (in the first case) or  $i_2\equiv 1+i_3\; mod \; 3$ (in the second case) or $i_3\equiv 1+i_2\; mod \; 3$ (in the third case).\\

Hence $F'$ may contain only monomials
\begin{itemize}
\item $(x_1y_1)(x_0y_2)(x_0y_3)$, $(x_0y_1)(x_1y_2)(x_0y_3)$, $(x_0y_1)(x_0y_2)(x_1y_3)$, $(x_0y_1)(x_1y_2)(x_1y_3)$, \\
$(x_1y_1)(x_0y_2)(x_1y_3)$, $(x_1y_1)(x_1y_2)(x_0y_3)$, $(x_0y_1)(x_0y_2)(x_0y_3)$, $(x_1y_1)(x_1y_2)(x_1y_3)$ in the first case,
\item $(x_1y_1)(x_0y_2)(x_0y_3)$, $(x_0y_1)(x_1y_2)(x_0y_3)$, $(x_0y_1)(x_0y_2)(x_1y_3)$, $(x_1y_1)(x_1y_2)(x_0y_1)$, \\
$(x_1y_2)(x_1y_3)(x_0y_2)$, $(x_1y_3)(x_1y_1)(x_0y_3)$ in the second case,
\item $(x_1y_1)(x_0y_2)(x_0y_3)$, $(x_0y_1)(x_1y_2)(x_0y_3)$, $(x_0y_1)(x_0y_2)(x_1y_3)$, $(x_1y_1)(x_1y_3)(x_0y_1)$, \\
$(x_1y_2)(x_1y_1)(x_0y_2)$, $(x_1y_3)(x_1y_2)(x_0y_3)$  in the third case.
\end{itemize}

We conclude that either
\begin{multline*}
F=a_1\cdot (x_0y_1)^2(x_1y_1) + a_2\cdot (x_0y_2)^2(x_1y_2) + a_3\cdot (x_0y_3)^2(x_1y_3) + \\
+b_1\cdot (x_1y_1)^2(x_0y_1) + b_2\cdot (x_1y_2)^2(x_0y_2) + b_3\cdot (x_1y_3)^2(x_0y_3) + \\
+c_1\cdot (x_1y_1)(x_0y_2)(x_0y_3)+c_2\cdot (x_0y_1)(x_1y_2)(x_0y_3)+c_3\cdot (x_0y_1)(x_0y_2)(x_1y_3)+\\
+d_1\cdot (x_0y_1)(x_1y_2)(x_1y_3)  +d_2\cdot  (x_1y_1)(x_0y_2)(x_1y_3) +d_3\cdot (x_1y_1)(x_1y_2)(x_0y_3)+\\
e_0\cdot (x_0y_1)(x_0y_2)(x_0y_3) +e_1\cdot (x_1y_1)(x_1y_2)(x_1y_3),
\end{multline*}
or
\begin{multline*}
F=a_1\cdot (x_0y_1)^2(x_1y_1) + a_2\cdot (x_0y_2)^2(x_1y_2) + a_3\cdot (x_0y_3)^2(x_1y_3) + \\
+b_1\cdot (x_1y_1)^2(x_0y_2) + b_2\cdot (x_1y_2)^2(x_0y_3) + b_3\cdot (x_1y_3)^2(x_0y_1) + \\
+c_1\cdot (x_1y_1)(x_0y_2)(x_0y_3)+c_2\cdot (x_0y_1)(x_1y_2)(x_0y_3)+c_3\cdot (x_0y_1)(x_0y_2)(x_1y_3)+\\
+d_1\cdot  (x_1y_1)(x_1y_2)(x_0y_1) +d_2\cdot  (x_1y_2)(x_1y_3)(x_0y_2) +d_3\cdot (x_1y_3)(x_1y_1)(x_0y_3) ,
\end{multline*}
or
\begin{multline*}
F=a_1\cdot (x_0y_1)^2(x_1y_1) + a_2\cdot (x_0y_2)^2(x_1y_2) + a_3\cdot (x_0y_3)^2(x_1y_3) + \\
+b_1\cdot (x_1y_1)^2(x_0y_3) + b_2\cdot (x_1y_2)^2(x_0y_1) + b_3\cdot (x_1y_3)^2(x_0y_2) + \\
+c_1\cdot (x_1y_1)(x_0y_2)(x_0y_3)+c_2\cdot (x_0y_1)(x_1y_2)(x_0y_3)+c_3\cdot (x_0y_1)(x_0y_2)(x_1y_3)+\\
+d_1\cdot  (x_1y_1)(x_1y_3)(x_0y_1) +d_2\cdot  (x_1y_2)(x_1y_1)(x_0y_2) +d_3\cdot (x_1y_3)(x_1y_2)(x_0y_3)
\end{multline*}

respectively in our three cases, where $a_i, b_i\in {\mathbb C}^{*}$, $(a_3/a_1)^3=1$, $a_2=a_3^2/a_1$ and
$$
b_3= {\om}^{-i} \cdot a_3b_1/a_1,\;\; b_2={\om}^{i} \cdot a_1b_1/a_3,\;\; c_3= a_3c_1/a_1,\;\; c_2=a_1c_1/a_3,
$$

$$
d_3= {\om}^{-i} \cdot a_3d_1/a_1,\; d_2={\om}^{i} \cdot a_1d_1/a_3.
$$

Moreover, if $e_0\neq 0$ or $e_1\neq 0$, then $a_1={\om}^i\cdot a_3$ in the first case.\\

It is immediate that cubic fourfolds given by the generic equations above are smooth (for all three cases).\\

In the first case this follows from the observation that the cubic fourfold given by the form
$$
F(z_0,z_1,z_2,z_3,z_4,z_5)=e_0z_0^2z_1+e_1z_1^2z_0+e_2z_2^2z_3+e_3z_3^2z_2+e_4z_4^2z_5+e_5z_5^2z_4
$$

is smooth for any $e_i\in {\mathbb C}^{*}$.\\

In the second and the third cases this follows from the observation that the cubic fourfold given by the form
$$
F(z_0,z_1,z_2,z_3,z_4,z_5)=e_0z_0^2z_1+e_1z_1^2z_2+e_2z_2^2z_3+e_3z_3^2z_4+e_4z_4^2z_5+e_5z_5^2z_0
$$

is smooth for any $e_i\in {\mathbb C}^{*}$.\\

\paragraph{Case B2. $F$ does not contain $(x_0y_i)^2(x_0y_j)$, $(x_1y_i)^2(x_1y_j)$ for any $i,j$, $F_y$ consists of monomials $y_1^2y_2, y_2^2y_3, y_3^2y_1$. }

Our assumptions imply that either
\begin{multline*}
F=a_1\cdot (x_0y_1)^2(x_1y_2) + a_2\cdot (x_0y_2)^2(x_1y_3) + a_3\cdot (x_0y_3)^2(x_1y_1) + \\
+b_1\cdot (x_1y_1)^2(x_0y_2) + b_2\cdot (x_1y_2)^2(x_0y_3) + b_3\cdot (x_1y_3)^2(x_0y_1) + F',
\end{multline*}
or
\begin{multline*}
F=a_1\cdot (x_0y_1)^2(x_1y_2) + a_2\cdot (x_0y_2)^2(x_1y_3) + a_3\cdot (x_0y_3)^2(x_1y_1) + \\
+b_1\cdot (x_1y_1)^2(x_0y_3) + b_2\cdot (x_1y_2)^2(x_0y_1) + b_3\cdot (x_1y_3)^2(x_0y_2) + F',
\end{multline*}
where $F'$ does not contain the singled out monomials and $a_i,b_i\in {\mathbb C}^{*}$.\\

Let us find $G_0=ker(\pi)$. If $f_0=(1,{\om}_0^{c_1})\times 1\in G_0$, then $c_1=0$. This means that $G_0=0$. Hence $\pi\colon G\rightarrow {\mathcal P}_3\cong \ZZ{3}\oplus \ZZ{3}$ is isomorphism (since we assume that it is surjective).\\

Let $f_1=(1,{\om}^{c_1})\times P_3$ and $f_2=(1,{\om}^{d_1})\times W_3$ be generators of $G$.\\

Then 
\begin{multline*}
f_1(F)={\om}^{c_1}\cdot ( a_1\cdot (x_0y_2)^2(x_1y_3) + a_2\cdot (x_0y_3)^2(x_1y_1) + a_3\cdot (x_0y_1)^2(x_1y_2) ) + \\
+{\om}^{2c_1}\cdot ( b_1\cdot (x_1y_2)^2(x_0y_3) + b_2\cdot (x_1y_3)^2(x_0y_1) + b_3\cdot (x_1y_1)^2(x_0y_2) ) + f_1(F'),
\end{multline*}
or
\begin{multline*}
f_1(F)={\om}^{c_1}\cdot ( a_1\cdot (x_0y_2)^2(x_1y_3) + a_2\cdot (x_0y_3)^2(x_1y_1) + a_3\cdot (x_0y_1)^2(x_1y_2) ) + \\
+{\om}^{2c_1}\cdot  ( b_1\cdot (x_1y_2)^2(x_0y_1) + b_2\cdot (x_1y_3)^2(x_0y_2) + b_3\cdot (x_1y_1)^2(x_0y_3) ) + f_1(F')
\end{multline*}

respectively and 

\begin{multline*}
f_2(F)={\om}^{d_1}\cdot {\om}_3\cdot  ( a_1\cdot (x_0y_1)^2(x_1y_2) + a_2\cdot (x_0y_2)^2(x_1y_3) + a_3\cdot (x_0y_3)^2(x_1y_1) ) + \\
+{\om}^{2d_1}\cdot {\om}_3\cdot ( b_1\cdot (x_1y_1)^2(x_0y_2) + b_2\cdot (x_1y_2)^2(x_0y_3) + b_3\cdot (x_1y_3)^2(x_0y_1) ) + f_2(F'),
\end{multline*}
or
\begin{multline*}
f_2(F)={\om}^{d_1}\cdot {\om}_3\cdot ( a_1\cdot (x_0y_1)^2(x_1y_2) + a_2\cdot (x_0y_2)^2(x_1y_3) + a_3\cdot (x_0y_3)^2(x_1y_1) ) + \\
+{\om}^{2d_1}\cdot {\om}_3^2 \cdot ( b_1\cdot (x_1y_1)^2(x_0y_3) + b_2\cdot (x_1y_2)^2(x_0y_1) + b_3\cdot (x_1y_3)^2(x_0y_2) ) + f_2(F')
\end{multline*}

respectively.\\ 

Since $f_1(F)\in {\mathbb C}^{*}\cdot F$, this implies that 
$$
a_3={\om}^{-c_1}\cdot \lambda a_1,\; a_2={\om}^{-2c_1}\cdot (\lambda )^2 a_1,\; b_3={\om}^{-2c_1}\cdot \lambda b_1,\; b_2={\om}^{-4c_1}\cdot (\lambda )^2 b_1, \; {\om}^{3c_1}=1 
$$
for some $\lambda\in \mathbb C$ such that ${\lambda }^3=1$.\\

Since $f_2(F)\in {\mathbb C}^{*}\cdot F$, this implies that
\begin{itemize}
\item ${\om}^{d_1}=1$,
\item ${\om}^{d_1}={\om}_3^2$
\end{itemize}
respectively in our two cases.\\

This means that we can take $\om={\om}_3$ and 
\begin{itemize}
\item $f_1=(1,{\om}^i)\times P_3$, $f_2=(1,1)\times W_3$ with some $i\in \{ 0,1,2 \}$ and $\om = \sqr{3}$,
\item $f_1=(1,{\om}^i)\times P_3$, $f_2=(1,{\om}^2)\times W_3$ with some $i \in \{ 0,1,2 \}$ and $\om = \sqr{3}$
\end{itemize}
respectively in our two cases.\\

We conclude that $G\cong (\ZZ{3})^{\oplus 2}$ with generators $f_1=(1,{\om}^i)\times P_3$, $f_2=(1,{\om}^j)\times W_3$ with some $i \in \{ 0,1,2 \}$ and some $j\in \{ 0,2 \}$, where $\om = \sqr{3}$. Here $j=0$ corresponds to the first case and $j=2$ corresponds to the second case.\\

In order to describe all smooth cubic fourfolds which admit these group actions, we need to determine which cubic forms $F'$ can appear.\\

By our assumption $x_0^3$ and $x_1^3$ may appear in $F$ only in monomials $(x_0y_1)(x_0y_2)(x_0y_3)$, $(x_1y_1)(x_1y_2)(x_1y_3)$. The invariance under $f_2$ requires that they are possible only in the second case. Moreover, if any of them appears in $F'$ with a nonzero coefficient, then $f_1(F)=F$.\\

Let $F'_0=F'_0(z_0,z_1,z_2,z_3,z_4,z_5)$ be a cubic form which is obtained from $F'$ by removing all monomials containing $x_0^3$, $x_0^2x_1$, $x_1^3$.\\

Let us denote by $F'_y$ the set of cubic monomials $y_{i}y_{j}y_{k}$ (in variables $y_1$, $y_2$, $y_3$) such that a monomial $(x_{i'}y_{i})(x_{j'}y_{j})(x_{k'}y_{k})$ appears in $F'_0$ with a nonzero coefficient for some $i',j',k'$.\\

The invariance under $f_2$ requires that $y_1^{i_1}y_2^{i_2}y_3^{i_3}$ is contained in $F'_y$ only if $i_2\equiv 1+i_3\; mod \; 3$ (in the first case) or $i_3\equiv 1+i_2\; mod \; 3$ (in the second case).\\

Hence $F'$ may contain only monomials
\begin{itemize}
\item $(x_0y_1)(x_0y_2)(x_1y_1)$, $(x_0y_2)(x_0y_3)(x_1y_2)$, $(x_0y_3)(x_0y_1)(x_1y_3)$, $(x_0y_1)(x_1y_2)(x_1y_1)$, \\
$(x_0y_2)(x_1y_3)(x_1y_2)$, $(x_0y_3)(x_1y_1)(x_1y_3)$ in the first case,
\item $(x_0y_1)(x_0y_2)(x_1y_1)$, $(x_0y_2)(x_0y_3)(x_1y_2)$, $(x_0y_3)(x_0y_1)(x_1y_3)$, $(x_0y_1)(x_1y_3)(x_1y_1)$, \\
$(x_0y_2)(x_1y_1)(x_1y_2)$, $(x_0y_3)(x_1y_2)(x_1y_3)$, $(x_0y_1)(x_0y_2)(x_0y_3)$, $(x_1y_1)(x_1y_2)(x_1y_3)$ in the second case.
\end{itemize}

We conclude that either
\begin{multline*}
F=a_1\cdot (x_0y_1)^2(x_1y_2) + a_2\cdot (x_0y_2)^2(x_1y_3) + a_3\cdot (x_0y_3)^2(x_1y_1) + \\
+b_1\cdot (x_1y_1)^2(x_0y_2) + b_2\cdot (x_1y_2)^2(x_0y_3) + b_3\cdot (x_1y_3)^2(x_0y_1) + \\
+c_1\cdot (x_0y_1)(x_0y_2)(x_1y_1) +c_2\cdot (x_0y_2)(x_0y_3)(x_1y_2)+c_3\cdot (x_0y_3)(x_0y_1)(x_1y_3) + \\
+d_1\cdot (x_0y_1)(x_1y_2)(x_1y_1) +d_2\cdot (x_0y_2)(x_1y_3)(x_1y_2)+d_3\cdot (x_0y_3)(x_1y_1)(x_1y_3),
\end{multline*}
or
\begin{multline*}
F=a_1\cdot (x_0y_1)^2(x_1y_2) + a_2\cdot (x_0y_2)^2(x_1y_3) + a_3\cdot (x_0y_3)^2(x_1y_1) + \\
+b_1\cdot (x_1y_1)^2(x_0y_3) + b_2\cdot (x_1y_2)^2(x_0y_1) + b_3\cdot (x_1y_3)^2(x_0y_2) + \\
+c_1\cdot (x_0y_1)(x_0y_2)(x_1y_1) +c_2\cdot (x_0y_2)(x_0y_3)(x_1y_2)+c_3\cdot (x_0y_3)(x_0y_1)(x_1y_3) + \\
+d_1\cdot (x_0y_1)(x_1y_3)(x_1y_1) +d_2\cdot (x_0y_2)(x_1y_1)(x_1y_2)+d_3\cdot (x_0y_3)(x_1y_2)(x_1y_3)+\\
+e_0\cdot (x_0y_1)(x_0y_2)(x_0y_3)+ e_1\cdot (x_1y_1)(x_1y_2)(x_1y_3)
\end{multline*}

respectively in our two cases, where $a_i, b_i\in {\mathbb C}^{*}$, $(a_3/a_1)^3=1$, $a_2=a_3^2/a_1$ and
$$
b_3= {\om}^{-i} \cdot a_3b_1/a_1,\;\; b_2={\om}^{i} \cdot a_1b_1/a_3,\;\; c_3= a_3c_1/a_1,\;\; c_2=a_1c_1/a_3,
$$

$$
d_3= {\om}^{-i} \cdot a_3d_1/a_1,\; d_2={\om}^{i} \cdot a_1d_1/a_3.
$$

Moreover, if $e_0\neq 0$ or $e_1\neq 0$, then $a_1={\om}^i \cdot a_3$ in the second case.\\

As it was in Case B1, it is immediate that cubic fourfolds given by the generic equations above are smooth (for both cases).\\

\paragraph{Case B3. $F$ does not contain $(x_0y_i)^2(x_0y_j)$, $(x_1y_i)^2(x_1y_j)$ for any $i,j$, $F_y$ consists of monomials $y_1^2y_3, y_3^2y_2, y_2^2y_1$. }

Our assumptions imply that either
\begin{multline*}
F=a_1\cdot (x_0y_1)^2(x_1y_3) + a_2\cdot (x_0y_3)^2(x_1y_2) + a_3\cdot (x_0y_2)^2(x_1y_1) + \\
+b_1\cdot (x_1y_1)^2(x_0y_2) + b_2\cdot (x_1y_2)^2(x_0y_3) + b_3\cdot (x_1y_3)^2(x_0y_1) + F',
\end{multline*}
or
\begin{multline*}
F=a_1\cdot (x_0y_1)^2(x_1y_3) + a_2\cdot (x_0y_3)^2(x_1y_2) + a_3\cdot (x_0y_2)^2(x_1y_1) + \\
+b_1\cdot (x_1y_1)^2(x_0y_3) + b_2\cdot (x_1y_2)^2(x_0y_1) + b_3\cdot (x_1y_3)^2(x_0y_2) + F',
\end{multline*}
where $F'$ does not contain the singled out monomials and $a_i,b_i\in {\mathbb C}^{*}$.\\

Let us find $G_0=ker(\pi)$. If $f_0=(1,{\om}_0^{c_1})\times 1\in G_0$, then $c_1=0$. This means that $G_0=0$. Hence $\pi\colon G\rightarrow {\mathcal P}_3\cong \ZZ{3}\oplus \ZZ{3}$ is isomorphism (since we assume that it is surjective).\\

Let $f_1=(1,{\om}^{c_1})\times P_3$ and $f_2=(1,{\om}^{d_1})\times W_3$ be generators of $G$.\\

Then 
\begin{multline*}
f_1(F)={\om}^{c_1}\cdot ( a_1\cdot (x_0y_2)^2(x_1y_1) + a_2\cdot (x_0y_1)^2(x_1y_3) + a_3\cdot (x_0y_3)^2(x_1y_2) ) + \\
+{\om}^{2c_1}\cdot ( b_1\cdot (x_1y_2)^2(x_0y_3) + b_2\cdot (x_1y_3)^2(x_0y_1) + b_3\cdot (x_1y_1)^2(x_0y_2) ) + f_1(F'),
\end{multline*}
or
\begin{multline*}
f_1(F)={\om}^{c_1}\cdot ( a_1\cdot (x_0y_2)^2(x_1y_1) + a_2\cdot (x_0y_1)^2(x_1y_3) + a_3\cdot (x_0y_3)^2(x_1y_2) ) + \\
+{\om}^{2c_1}\cdot  ( b_1\cdot (x_1y_2)^2(x_0y_1) + b_2\cdot (x_1y_3)^2(x_0y_2) + b_3\cdot (x_1y_1)^2(x_0y_3) ) + f_1(F')
\end{multline*}

respectively and 

\begin{multline*}
f_2(F)={\om}^{d_1}\cdot {\om}_3^2\cdot  ( a_1\cdot (x_0y_1)^2(x_1y_3) + a_2\cdot (x_0y_3)^2(x_1y_2) + a_3\cdot (x_0y_2)^2(x_1y_1) ) + \\
+{\om}^{2d_1}\cdot {\om}_3\cdot ( b_1\cdot (x_1y_1)^2(x_0y_2) + b_2\cdot (x_1y_2)^2(x_0y_3) + b_3\cdot (x_1y_3)^2(x_0y_1) ) + f_2(F'),
\end{multline*}
or
\begin{multline*}
f_2(F)={\om}^{d_1}\cdot {\om}_3^2\cdot ( a_1\cdot (x_0y_1)^2(x_1y_3) + a_2\cdot (x_0y_3)^2(x_1y_2) + a_3\cdot (x_0y_2)^2(x_1y_1) ) + \\
+{\om}^{2d_1}\cdot {\om}_3^2 \cdot ( b_1\cdot (x_1y_1)^2(x_0y_3) + b_2\cdot (x_1y_2)^2(x_0y_1) + b_3\cdot (x_1y_3)^2(x_0y_2) ) + f_2(F')
\end{multline*}

respectively.\\ 

Since $f_1(F)\in {\mathbb C}^{*}\cdot F$, this implies that 
$$
a_2={\om}^{-c_1}\cdot \lambda a_1,\; a_3={\om}^{-2c_1}\cdot (\lambda )^2 a_1,\; b_3={\om}^{-2c_1}\cdot \lambda b_1,\; b_2={\om}^{-4c_1}\cdot (\lambda )^2 b_1, \; {\om}^{3c_1}=1 
$$
for some $\lambda\in \mathbb C$ such that ${\lambda }^3=1$.\\

Since $f_2(F)\in {\mathbb C}^{*}\cdot F$, this implies that
\begin{itemize}
\item ${\om}^{d_1}={\om}_3$,
\item ${\om}^{d_1}=1$
\end{itemize}
respectively in our two cases.\\

This means that we can take $\om={\om}_3$ and 
\begin{itemize}
\item $f_1=(1,{\om}^i)\times P_3$, $f_2=(1,{\om})\times W_3$ with some $i\in \{ 0,1,2 \}$ and $\om = \sqr{3}$,
\item $f_1=(1,{\om}^i)\times P_3$, $f_2=(1,1)\times W_3$ with some $i \in \{ 0,1,2 \}$ and $\om = \sqr{3}$
\end{itemize}
respectively in our two cases.\\

We conclude that $G\cong (\ZZ{3})^{\oplus 2}$ with generators $f_1=(1,{\om}^i)\times P_3$, $f_2=(1,{\om}^j)\times W_3$ with some $i \in \{ 0,1,2 \}$ and some $j\in \{ 0,1 \}$, where $\om = \sqr{3}$. Here $j=1$ corresponds to the first case and $j=0$ corresponds to the second case.\\

In order to describe all smooth cubic fourfolds which admit these group actions, we need to determine which cubic forms $F'$ can appear.\\

By our assumption $x_0^3$ and $x_1^3$ may appear in $F$ only in monomials $(x_0y_1)(x_0y_2)(x_0y_3)$, $(x_1y_1)(x_1y_2)(x_1y_3)$. The invariance under $f_2$ requires that they are possible only in the first case. Moreover, if any of them appears in $F'$ with a nonzero coefficient, then $f_1(F)=F$.\\

Let $F'_0=F'_0(z_0,z_1,z_2,z_3,z_4,z_5)$ be a cubic form which is obtained from $F'$ by removing all monomials containing $x_0^3$, $x_0^2x_1$, $x_1^3$.\\

Let us denote by $F'_y$ the set of cubic monomials $y_{i}y_{j}y_{k}$ (in variables $y_1$, $y_2$, $y_3$) such that a monomial $(x_{i'}y_{i})(x_{j'}y_{j})(x_{k'}y_{k})$ appears in $F'_0$ with a nonzero coefficient for some $i',j',k'$.\\

The invariance under $f_2$ requires that $y_1^{i_1}y_2^{i_2}y_3^{i_3}$ is contained in $F'_y$ only if $2j+i_2+2i_3\equiv j+2 \; mod \; 3$, i.e. $i_2\equiv 1+i_3\; mod \; 3$ (in the first case) or $i_3\equiv 1+i_2\; mod \; 3$ (in the second case).\\

Hence $F'$ may contain only monomials
\begin{itemize}
\item $(x_0y_1)(x_0y_3)(x_1y_1)$, $(x_0y_2)(x_0y_1)(x_1y_2)$, $(x_0y_3)(x_0y_2)(x_1y_3)$, $(x_0y_1)(x_1y_2)(x_1y_1)$, \\
$(x_0y_2)(x_1y_3)(x_1y_2)$, $(x_0y_3)(x_1y_1)(x_1y_3)$, $(x_0y_1)(x_0y_2)(x_0y_3)$, $(x_1y_1)(x_1y_2)(x_1y_3)$ in the first case,
\item $(x_0y_1)(x_0y_3)(x_1y_1)$, $(x_0y_2)(x_0y_1)(x_1y_2)$, $(x_0y_3)(x_0y_2)(x_1y_3)$, $(x_0y_1)(x_1y_3)(x_1y_1)$, \\
$(x_0y_2)(x_1y_1)(x_1y_2)$, $(x_0y_3)(x_1y_2)(x_1y_3)$ in the second case.
\end{itemize}

We conclude that either
\begin{multline*}
F=a_1\cdot (x_0y_1)^2(x_1y_3) + a_2\cdot (x_0y_3)^2(x_1y_2) + a_3\cdot (x_0y_2)^2(x_1y_1) + \\
+b_1\cdot (x_1y_1)^2(x_0y_2) + b_2\cdot (x_1y_2)^2(x_0y_3) + b_3\cdot (x_1y_3)^2(x_0y_1) + \\
+c_1\cdot (x_0y_1)(x_0y_3)(x_1y_1) +c_2\cdot (x_0y_2)(x_0y_1)(x_1y_2)+c_3\cdot (x_0y_3)(x_0y_2)(x_1y_3) + \\
+d_1\cdot (x_0y_1)(x_1y_2)(x_1y_1) +d_2\cdot (x_0y_2)(x_1y_3)(x_1y_2)+d_3\cdot (x_0y_3)(x_1y_1)(x_1y_3)+\\
+e_0\cdot (x_0y_1)(x_0y_2)(x_0y_3)+e_1\cdot (x_1y_1)(x_1y_2)(x_1y_3),
\end{multline*}
or
\begin{multline*}
F=a_1\cdot (x_0y_1)^2(x_1y_3) + a_2\cdot (x_0y_3)^2(x_1y_2) + a_3\cdot (x_0y_2)^2(x_1y_1) + \\
+b_1\cdot (x_1y_1)^2(x_0y_3) + b_2\cdot (x_1y_2)^2(x_0y_1) + b_3\cdot (x_1y_3)^2(x_0y_2) + \\
+c_1\cdot (x_0y_1)(x_0y_3)(x_1y_1) +c_2\cdot (x_0y_2)(x_0y_1)(x_1y_2)+c_3\cdot (x_0y_3)(x_0y_2)(x_1y_3) + \\
+d_1\cdot (x_0y_1)(x_1y_3)(x_1y_1) +d_2\cdot (x_0y_2)(x_1y_1)(x_1y_2)+d_3\cdot (x_0y_3)(x_1y_2)(x_1y_3)
\end{multline*}

respectively in our two cases, where $a_i, b_i\in {\mathbb C}^{*}$, $(a_2/a_1)^3=1$, $a_3=a_2^2/a_1$ and
$$
b_3= {\om}^{-i} \cdot a_2b_1/a_1,\;\; b_2={\om}^{i} \cdot a_1b_1/a_2,\;\; c_3= a_2c_1/a_1,\;\; c_2=a_1c_1/a_2,
$$

$$
d_3= {\om}^{-i} \cdot a_2d_1/a_1,\; d_2={\om}^{i} \cdot a_1d_1/a_2.
$$

Moreover, if $e_0\neq 0$ or $e_1\neq 0$, then $a_1={\om}^i\cdot a_2$ in the first case.\\

As it was in Case B1, it is immediate that cubic fourfolds given by the generic equations above are smooth (for both cases).\\

\vspace{2cm}

We can summarize our findings in this subsection in the following theorem.\\

{\bf Theorem 3.} {\it Abelian automorphism groups of smooth cubic fourfolds which are conjugate to subgroups of $D_2\otimes {\mathcal P}_{3}\subset PGL(6)=Aut({\mathbb P}^5)$ are either diagonalizable or the subgroups of the following groups (see Table $3$ below).

\begin{center}
  \includegraphics[width=1.0\textwidth]{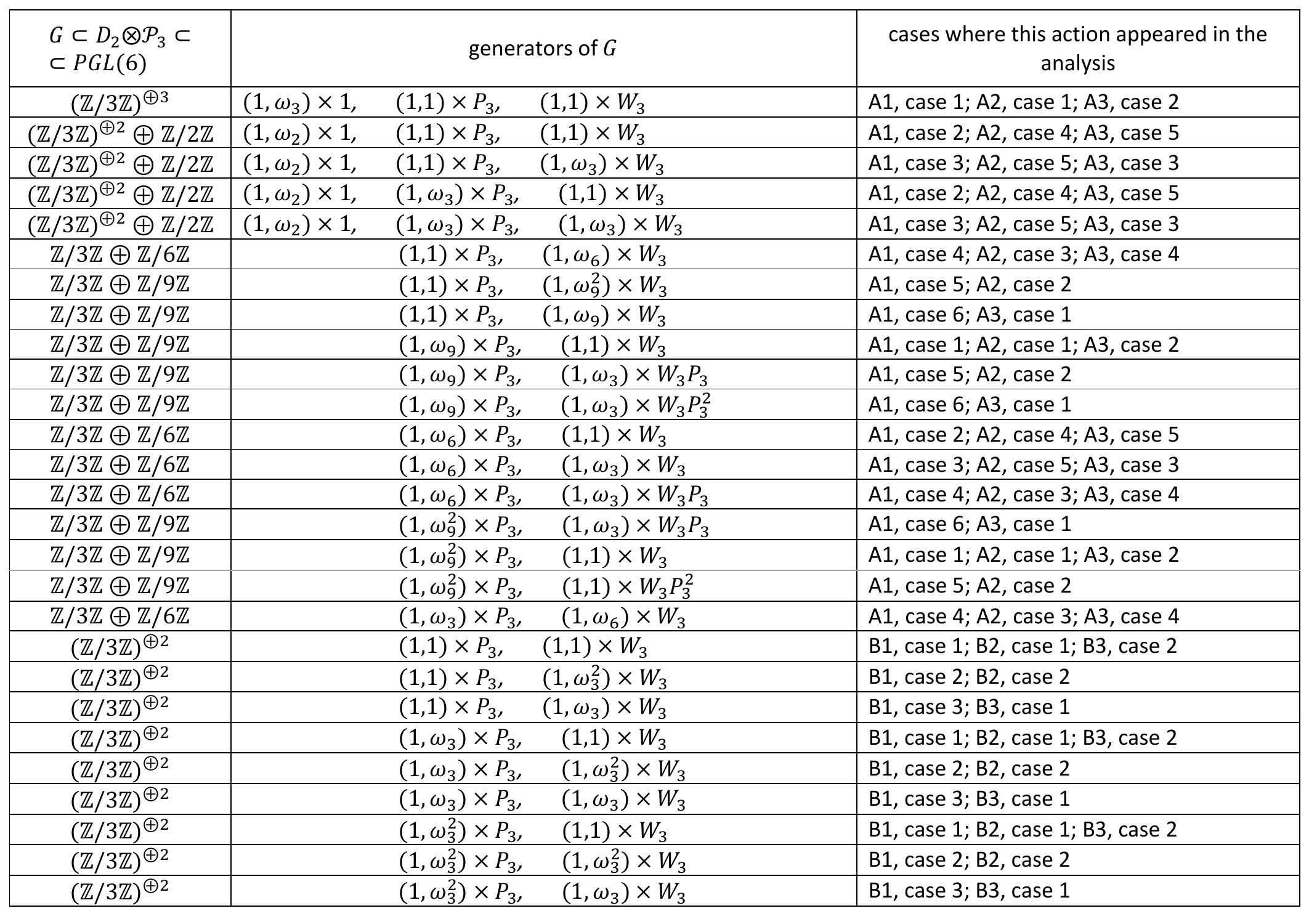}
\end{center}

{\bf Table $3$.} {\sl Abelian automorphism groups of smooth cubic fourfolds which are conjugate to subgroups of $D_2\otimes {\mathcal P}_{3}\subset PGL(6)$.} The first column lists abelian groups $G$ which act effectively on a smooth cubic fourfold. The second column describes generators of $G$ (using the notation of the present subsection). Explicit descriptions of cubic fourfolds which admit these group actions are given in the main part of this subsection (see the analysis of cases A1, A2, A3, B1, B2 and B3 above and directions in the third column). Note that we had $6$ subcases total in Case A1. ${\om}_p$ denotes $\sqr{p}$.}\\

\subsection{Case of $G\subset {\mathcal P}_{6}$.}

In this case ${\mathcal P}_{6}=\{ P_6^i\cdot W_6^j \; \mid \; i,j=0,1,2,3,4,5 \}\subset PGL(6)$, where
$$
P_6= \begin{pmatrix}
0 & 1 & 0 & 0 & 0 & 0 \\
0 & 0 & 1 & 0 & 0 & 0 \\
0 & 0 & 0 & 1 & 0 & 0 \\
0 & 0 & 0 & 0 & 1 & 0 \\
0 & 0 & 0 & 0 & 0 & 1 \\
1 & 0 & 0 & 0 & 0 & 0\\
\end{pmatrix},\;
W_6= \begin{pmatrix}
1 & 0 & 0 & 0 & 0 & 0\\
0 & {\om}_6 & 0 & 0 & 0 & 0\\
0 & 0 & {\om}_6^2 & 0 & 0 & 0\\
0 & 0 & 0 & {\om}_6^3 & 0 & 0\\
0 & 0 & 0 & 0 & {\om}_6^4 & 0\\
0 & 0 & 0 & 0 & 0 & {\om}_6^5\\
\end{pmatrix}
$$
and ${\om}_6=\sqr{6}$.\\

Recall that ${\mathcal P}_{6} \cong \ZZ{6}\oplus \ZZ{6}$ (\cite{GangHan}).\\

In this subsection we will denote ${\om}_p=\sqr{p}$.\\

Let $G\subset {\mathcal P}_{6} \subset PGL(6)$.\\

{\bf Lemma 7.} {\it G is conjugate to a finite subgroup of $D_6$, $D_2\otimes {\mathcal P}_{3}$ or $D_3\otimes {\mathcal P}_{2}$ unless $P_6\cdot W^i_6\in G$ for some $i$.}\\

{\it Proof:} Without loss of generality we may assume that $G$ is one of the following subgroups of ${\mathcal P}_{6}$:
\begin{itemize}
\item $G_0=\{ W^i_6 \; \mid \; i=0,1,2,3,4,5 \} \subset {\mathcal P}_{6}$,
\item $G_2=\{ W^i_6,\; P_6^2\cdot W^i_6,\; P_6^4\cdot W^i_6 \; \mid \; i=0,1,2,3,4,5 \} \subset {\mathcal P}_{6}$,
\item $G_3=\{ W^i_6,\; P_6^3\cdot W^i_6 \; \mid \; i=0,1,2,3,4,5 \} \subset {\mathcal P}_{6}$.
\end{itemize}

Let us denote by $z_0,z_1,z_2,z_3,z_4,z_5$ the homogeneous coordinates in ${\mathbb P}^5$ and identify $PGL(6)=Aut({\mathbb P}^5)$.\\

It is immediate that $G_0\subset D_6\subset PGL(6)$.\\

Let us consider $G_2$.\\

Since $P_6(z_i)=z_{i+1}$, $P_6(z_5)=z_0$, $i=0,1,2,3,4$, we have:
$$
P^2_6(z_0)=z_2,\; \; P^2_6(z_2)=z_4,\; \; P^2_6(z_4)=z_0, 
$$

$$
P^2_6(z_1)=z_3,\; \; P^2_6(z_3)=z_5,\; \; P^2_6(z_5)=z_1. 
$$

This means that we can write
$$
z_0=x_0y_1,\; z_1=x_1y_1,\; z_2=x_0y_2,\; z_3=x_1y_2,\; z_4=x_0y_3,\; z_5=x_1y_3, 
$$

and then 
$$
P^2_6=(1,1)\times P_3,\;\;\;\; W_6=(1,{\om}_6)\times W_3,
$$

where on the right hand sides we have elements of $D_2\otimes {\mathcal P}_{3}\subset PGL(6)$ with respect to the decomposition $z_k=x_iy_j$ above.\\

In other words, $G_2$ is indeed conjugate to a finite subgroup of $D_2\otimes {\mathcal P}_{3}\subset PGL(6)$.\\

Let us consider $G_3$.\\

Note that
$$
P^3_6(z_0)=z_3,\; \; P^3_6(z_1)=z_4,\; \; P^3_6(z_2)=z_5, 
$$

$$
P^3_6(z_3)=z_1,\; \; P^3_6(z_4)=z_1,\; \; P^3_6(z_5)=z_2. 
$$

This means that we can write
$$
z_0=x_0y_1,\; z_1=x_1y_1,\; z_2=x_2y_1,\; z_3=x_0y_2,\; z_4=x_1y_2,\; z_5=x_2y_2, 
$$

and then 
$$
P^3_6=(1,1,1)\times P_2,\;\;\;\; W_6=(1,{\om}_6,{\om}^2_6)\times W_2,
$$

where on the right hand sides we have elements of $D_3\otimes {\mathcal P}_{2}\subset PGL(6)$ with respect to the decomposition $z_k=x_iy_j$ above.\\

In other words, $G_3$ is indeed conjugate to a finite subgroup of $D_3\otimes {\mathcal P}_{2}\subset PGL(6)$. {\it QED}\\

Hence without loss of generality we may assume that $P_6\cdot W^i_6\in G$ for some $i$.\\

Let $X\subset {\mathbb P}^5$ be a smooth cubic fourfold given by a cubic form $F=F(z_0,z_1,z_2,z_3,z_4,z_5)$, on which a finite abelian group $G$ acts effectively.\\

We will consider two cases:
\begin{itemize}
\item[(C)] $F$ contains $z^3_i$ for some $i$,
\item[(D)] $F$ does not contain $z^3_i$ for any $i$.
\end{itemize}

\subsubsection{Case C. $F$ contains $z^3_i$ for some $i$.}

In this case we can write
$$
F=a_0\cdot z^3_0+ a_1\cdot z^3_1+ a_2\cdot z^3_2+ a_3\cdot z^3_3+ a_4\cdot z^3_4+ a_5\cdot z^3_5+F',
$$
where $a_i\in {\mathbb C}^{*}$ and $F'$ does not contain the singled out monomials.\\

Then 
$$
P_6\cdot W^i_6(F)=  a_0\cdot z^3_1+ {\om}^{3i}_6 \cdot a_1\cdot z^3_2+ {\om}^{6i}_6 \cdot a_2\cdot z^3_3+ {\om}^{9i}_6 \cdot a_3\cdot z^3_4+ {\om}^{12i}_6 \cdot a_4\cdot z^3_5+ {\om}^{15i}_6 \cdot a_5\cdot z^3_0 +P_6\cdot W^i_6(F'),
$$
which should be proportional to $F$. Hence 
$$
a_5=(-1)^i\cdot \lambda\cdot a_0,\; \; a_4=(-1)^i\cdot {\lambda}^2 \cdot a_0,\; \; a_3={\lambda}^3 \cdot a_0,
$$

$$
a_2={\lambda}^4 \cdot a_0,\; \; a_1=(-1)^i\cdot {\lambda}^5 \cdot a_0
$$
for some $\lambda\in \mathbb C$ such that ${\lambda}^6=(-1)^i$.\\

In other words, we need to take $a_5$ such that $(a_5/a_0)^6=(-1)^i$ and then
$$
a_4=(-1)^i\cdot a_5^2/a_0,\;\; a_3=(-1)^i\cdot a_5^3/a_0^2,\;\; a_2=(-1)^i\cdot a_0^3/a_5^2,\;\; a_1=(-1)^i\cdot a_0^2/a_5.
$$

On the other hand, 
$$
W^j_6(F)= a_0\cdot z^3_0+ (-1)^{j} \cdot a_1\cdot z^3_1+  a_2\cdot z^3_2+ (-1)^{j} \cdot a_3\cdot z^3_3+  a_4\cdot z^3_4+ (-1)^{j} \cdot a_5\cdot z^3_5 + W^j_6(F'),
$$
which means that $W^j_6(F)\in {\mathbb C}^{*}\cdot F$ only if $j$ is even.\\

We conclude that $G\subset {\mathcal P}_{6}$ is generated by $P_6\cdot W^i_6$ and $W^2_6$ for some $i\in \{ 0,1 \}$.\\

In particular, if $i=0$, then $G\cong \ZZ{6}\oplus \ZZ{3}$ is generated by $P_6$ and $W^2_6$.\\

If $i=1$, then $G\cong \ZZ{6}\oplus \ZZ{3}$ is generated by $P_6\cdot W_6$ and $W^2_6$.\\

In order to describe all cubic fourfolds which admit this group action, we need to find all possible cubic forms $F'$.\\ 

If a monomial $z_{p}z_{q}z_{r}$ is present in $F'$, then the invariance with respect to $W^2_6$ requires that $p+q+r\equiv 0\;  mod \; 3$.\\

Since the set of such monomials is invariant with respect to $P_6$, we conclude that
$$
F'=\sum_{0\leq p\leq q\leq r \leq 5,\;\; p+q+r\equiv 0 \; mod \; 3}^{\prime} a_{pqr} \cdot z_{p}z_{q}z_{r}
$$
with suitable coefficients $a_{pqr}$. The prime signifies that the sum is over $p,q,r$, which are {\it not all the same} (in order to exclude monomials $z_p^3$ which were written in $F$ separately).\\

Coefficients $a_{pqr}$ can be determined as follows. We have to require that $P_6\cdot W^i_6(F')={\lambda}\cdot F'$ with the same $\lambda = (-1)^i\cdot a_5/a_0$ as above.\\

Let $S_1,S_2,...$ be the orbits of the natural action of $P_6$ on the set of monomials $z_{p}z_{q}z_{r}$ such that $p+q+r\equiv 0 \; mod \; 3$ and $p,q,r$ are not all the same.\\

Then for every such orbit $S_{\alpha}$ the condition $P_6\cdot W^i_6(F')={\lambda}\cdot F'$ determines the coefficients $a_{pqr}$ of monomials in $S_{\alpha}$ upto a common factor.\\

Namely, if $m_1,m_2,...,m_k$ are elements of $S_{\alpha}$ such that $P_6(m_p)=m_{p+1}$, $p=1,2,...,k-1$, $P_6(m_k)=m_{1}$, $W_6(m_p)=(-1)^{c_p}\cdot m_p$ (where $c_p\in \mathbb Z$) and $b_p$ is the coefficient of $m_p$ in the sum for $F'$ above, then
$$
b_k=\lambda \cdot (-1)^{c_k\cdot i}\cdot b_{1},\; \; b_p=\lambda \cdot (-1)^{c_p\cdot i}\cdot b_{p+1},\; \; p=1,2,...,k-1.
$$

This gives by induction:
$$
b_p={\lambda }^{k+1-p} \cdot (-1)^{(c_p+...+c_k)\cdot i}\cdot b_{1},\; \; p=1,...,k.
$$

Note that these formulas determine $b_2,...,b_k$ in terms of $b_1$, which has to be zero unless $\lambda = (-1)^i\cdot a_5/a_0$ satisfies the following condition:
$$
{\lambda}^k=(-1)^{(c_1+...+c_k)\cdot i},
$$

i.e.
$$
(a_5/a_0)^k=(-1)^{(c_1+...+c_k-k)\cdot i}.
$$

If this condition is satisfied, then the formulas above determine $b_1,b_2,...,b_k$ upto a common factor. If it is not, then the formulas above say that $b_1=b_2=...=b_k=0$, i.e. monomials $m_1,m_2,...,m_k$ from $S_{\alpha}$ do not appear in $F'$.\\

Note that since $P^6_6=Id$, only the following values of $k$ are possible: $k=1$, $k=2$, $k=3$, $k=6$.\\

In the case $k=1$ the orbit $S_{\alpha}$ consists of just one monomial $m_1$ such that $P_6(m_1)=m_1$ and $W_6(m_1)=(-1)^{c_1}\cdot m_1$ for some $c_1\in \mathbb Z$. In this case the coefficient $b_1$ of this monomial in $F'$ is zero unless $\lambda =(-1)^{c_1\cdot i}$, i.e. $a_5=(-1)^{(c_1-1)\cdot i}\cdot a_0$.\\

Our earlier condition $(a_5/a_0)^6=(-1)^i$ implies then that $i=0$ and
$$
a_0=a_1=a_2=a_3=a_4=a_5
$$
in this case.\\ 

If this is satisfied, then $b_1$ may be arbitrary in this case.\\

This procedure determines coefficients $a_{pqr}$ in $F'$.\\

It is immediate that generic such cubic forms $F$ give smooth cubic fourfolds.\\

\subsubsection{Case D. $F$ does not contain $z^3_i$ for any $i$.}

By Lemma 1 (\cite{Liendo}, Lemma 1.3) $F$ should contain $z_0^2z_i$ for some $i\in \{ 1,2,3,4,5 \}$.\\

Since $P_6\cdot W^i_6(F)\in {\mathbb C}^{*}\cdot F$ and $P_6\cdot W^i_6(z_{p}z_{q}z_{r})\in {\mathbb C}^{*}\cdot z_{p+1}z_{q+1}z_{r+1}$ (adding indices modulo $6$), we should have that either
\begin{itemize}
\item $F=a_0\cdot z_0^2z_1+ a_1\cdot z_1^2z_2 + a_2\cdot z_2^2z_3 + a_3\cdot z_3^2z_4 + a_4\cdot z_4^2z_5 + a_5\cdot z_5^2z_0 +F'$, or
\item $F=a_0\cdot z_0^2z_2+ a_1\cdot z_1^2z_3 + a_2\cdot z_2^2z_4 + a_3\cdot z_3^2z_5 + a_4\cdot z_4^2z_0 + a_5\cdot z_5^2z_1 +F'$, or
\item $F=a_0\cdot z_0^2z_3+ a_1\cdot z_1^2z_4 + a_2\cdot z_2^2z_5 + a_3\cdot z_3^2z_0 + a_4\cdot z_4^2z_1 + a_5\cdot z_5^2z_2 +F'$, or
\item $F=a_0\cdot z_0^2z_4+ a_1\cdot z_1^2z_5 + a_2\cdot z_2^2z_0 + a_3\cdot z_3^2z_1 + a_4\cdot z_4^2z_2 + a_5\cdot z_5^2z_3 +F'$, or
\item $F=a_0\cdot z_0^2z_5+ a_1\cdot z_1^2z_0 + a_2\cdot z_2^2z_1 + a_3\cdot z_3^2z_2 + a_4\cdot z_4^2z_3 + a_5\cdot z_5^2z_4 +F'$,
\end{itemize}
where $a_i\in {\mathbb C}^{*}$ and $F'$ does not contain the singled out monomials.\\

Note that 
\begin{itemize}
\item $W^j_6(F)={\om}^{j}_6 \cdot a_0\cdot z_0^2z_1+ {\om}^{4j}_6 \cdot a_1\cdot z_1^2z_2 + {\om}^{7j}_6 \cdot a_2\cdot z_2^2z_3 + {\om}^{10j}_6 \cdot a_3\cdot z_3^2z_4 + {\om}^{13j}_6 \cdot a_4\cdot z_4^2z_5 + {\om}^{10j}_6 \cdot a_5\cdot z_5^2z_0 +W^j_6(F')$, or
\item $W^j_6(F)={\om}^{2j}_6 \cdot a_0\cdot z_0^2z_2+ {\om}^{5j}_6 \cdot a_1\cdot z_1^2z_3 + {\om}^{8j}_6 \cdot a_2\cdot z_2^2z_4 + {\om}^{11j}_6 \cdot a_3\cdot z_3^2z_5 + {\om}^{8j}_6 \cdot a_4\cdot z_4^2z_0 + {\om}^{11j}_6 \cdot a_5\cdot z_5^2z_1 +W^j_6(F')$, or
\item $W^j_6(F)={\om}^{3j}_6 \cdot a_0\cdot z_0^2z_3+ {\om}^{6j}_6 \cdot a_1\cdot z_1^2z_4 + {\om}^{9j}_6 \cdot a_2\cdot z_2^2z_5 + {\om}^{6j}_6 \cdot a_3\cdot z_3^2z_0 + {\om}^{9j}_6 \cdot a_4\cdot z_4^2z_1 + {\om}^{12j}_6 \cdot a_5\cdot z_5^2z_2 +W^j_6(F')$, or
\item $W^j_6(F)={\om}^{4j}_6 \cdot a_0\cdot z_0^2z_4+ {\om}^{7j}_6 \cdot a_1\cdot z_1^2z_5 + {\om}^{4j}_6 \cdot a_2\cdot z_2^2z_0 + {\om}^{7j}_6 \cdot a_3\cdot z_3^2z_1 + {\om}^{10j}_6 \cdot a_4\cdot z_4^2z_2 + {\om}^{13j}_6 \cdot a_5\cdot z_5^2z_3 +W^j_6(F')$, or
\item $W^j_6(F)={\om}^{5j}_6 \cdot a_0\cdot z_0^2z_5+ {\om}^{2j}_6 \cdot a_1\cdot z_1^2z_0 + {\om}^{5j}_6 \cdot a_2\cdot z_2^2z_1 + {\om}^{8j}_6 \cdot a_3\cdot z_3^2z_2 + {\om}^{11j}_6 \cdot a_4\cdot z_4^2z_3 + {\om}^{14j}_6 \cdot a_5\cdot z_5^2z_4 +W^j_6(F')$
\end{itemize}
respectively in our five cases.\\

Hence $W^j_6(F)\in {\mathbb C}^{*}\cdot F$ only if $3j\equiv 0 \; mod \; 6$ (for all five cases).\\

This means that $W^2_6\in G$, but $W_6\notin G$ (for all five cases). Moreover, $W^2_6(F)=\mu\cdot F$, where 
\begin{itemize}
\item $\mu={\om}_3$ in the first and the fourth cases,
\item $\mu={\om}^2_3$ in the second and the fifth cases,
\item $\mu=1$ in the third case.
\end{itemize}

As we noted above, we may assume that $P_6\cdot W^i_6\in G$ for some $i$. Since $W^2_6\in G$, we may take $i\in \{ 0,1 \}$.\\

Note that
\begin{itemize}
\item $P_6W^i_6(F)={\om}^{i}_6 \cdot a_0\cdot z_1^2z_2+ {\om}^{4i}_6 \cdot a_1\cdot z_2^2z_3 + {\om}^{7i}_6 \cdot a_2\cdot z_3^2z_4 + {\om}^{10i}_6 \cdot a_3\cdot z_4^2z_5 + {\om}^{13i}_6 \cdot a_4\cdot z_5^2z_0 + {\om}^{10i}_6 \cdot a_5\cdot z_0^2z_1 +P_6W^i_6(F')$, or
\item $P_6W^i_6(F)={\om}^{2i}_6 \cdot a_0\cdot z_1^2z_3+ {\om}^{5i}_6 \cdot a_1\cdot z_2^2z_4 + {\om}^{8i}_6 \cdot a_2\cdot z_3^2z_5 + {\om}^{11i}_6 \cdot a_3\cdot z_4^2z_0 + {\om}^{8i}_6 \cdot a_4\cdot z_5^2z_1 + {\om}^{11i}_6 \cdot a_5\cdot z_0^2z_2 +P_6W^i_6(F')$, or
\item $P_6W^i_6(F)={\om}^{3i}_6 \cdot a_0\cdot z_1^2z_4+ {\om}^{6i}_6 \cdot a_1\cdot z_2^2z_5 + {\om}^{9i}_6 \cdot a_2\cdot z_3^2z_0 + {\om}^{6i}_6 \cdot a_3\cdot z_4^2z_1 + {\om}^{9i}_6 \cdot a_4\cdot z_5^2z_2 + {\om}^{12i}_6 \cdot a_5\cdot z_0^2z_3 +P_6W^i_6(F')$, or
\item $P_6W^i_6(F)={\om}^{4i}_6 \cdot a_0\cdot z_1^2z_5+ {\om}^{7i}_6 \cdot a_1\cdot z_2^2z_0 + {\om}^{4i}_6 \cdot a_2\cdot z_3^2z_1 + {\om}^{7i}_6 \cdot a_3\cdot z_4^2z_2 + {\om}^{10i}_6 \cdot a_4\cdot z_5^2z_3 + {\om}^{13i}_6 \cdot a_5\cdot z_0^2z_4 +P_6W^i_6(F')$, or
\item $P_6W^i_6(F)={\om}^{5i}_6 \cdot a_0\cdot z_1^2z_0+ {\om}^{2i}_6 \cdot a_1\cdot z_2^2z_1 + {\om}^{5i}_6 \cdot a_2\cdot z_3^2z_2 + {\om}^{8i}_6 \cdot a_3\cdot z_4^2z_3 + {\om}^{11i}_6 \cdot a_4\cdot z_5^2z_4 + {\om}^{14i}_6 \cdot a_5\cdot z_0^2z_5 +P_6W^i_6(F')$
\end{itemize}
respectively in our five cases.\\

In other words,
\begin{itemize}
\item $P_6W^i_6(F)={\om}^{i}_6 \cdot (a_0\cdot z_1^2z_2 - a_1\cdot z_2^2z_3 + a_2\cdot z_3^2z_4 - a_3\cdot z_4^2z_5 + a_4\cdot z_5^2z_0 - a_5\cdot z_0^2z_1)+P_6W^i_6(F')$, or
\item $P_6W^i_6(F)={\om}^{2i}_6 \cdot (a_0\cdot z_1^2z_3 - a_1\cdot z_2^2z_4 +  a_2\cdot z_3^2z_5 - a_3\cdot z_4^2z_0 + a_4\cdot z_5^2z_1 - a_5\cdot z_0^2z_2)+P_6W^i_6(F')$, or
\item $P_6W^i_6(F)={\om}^{3i}_6 \cdot (a_0\cdot z_1^2z_4 - a_1\cdot z_2^2z_5 + a_2\cdot z_3^2z_0 - a_3\cdot z_4^2z_1 + a_4\cdot z_5^2z_2 - a_5\cdot z_0^2z_3)+P_6W^i_6(F')$, or
\item $P_6W^i_6(F)={\om}^{4i}_6 \cdot (a_0\cdot z_1^2z_5 - a_1\cdot z_2^2z_0 + a_2\cdot z_3^2z_1 - a_3\cdot z_4^2z_2 + a_4\cdot z_5^2z_3 - a_5\cdot z_0^2z_4)+P_6W^i_6(F')$, or
\item $P_6W^i_6(F)={\om}^{5i}_6 \cdot (a_0\cdot z_1^2z_0 - a_1\cdot z_2^2z_1 + a_2\cdot z_3^2z_2 - a_3\cdot z_4^2z_3 + a_4\cdot z_5^2z_4 - a_5\cdot z_0^2z_5) +P_6W^i_6(F')$
\end{itemize}
respectively in our five cases.\\

Hence $P_6W^i_6(F)\in {\mathbb C}^{*}\cdot F$ if and only if 
$$
a_1=-{\lambda}^5\cdot a_0,\;\; a_2={\lambda}^4\cdot a_0,\;\; a_3={\lambda}^3\cdot a_0,\;\; a_4=-{\lambda}^2\cdot a_0,\;\; a_5=-{\lambda}\cdot a_0 
$$
for some $\lambda\in {\mathbb C}^{*}$ such that ${\lambda}^6=-1$ in each of our five cases.\\

In other words, we need to take $a_5$ such that $(a_5/a_0)^6=-1$ and then
$$
a_1=-a_0^2/a_5,\;\; a_2=-a_0^3/a_5^2,\;\; a_3=-a_5^3/a_0^2,\;\; a_4=-a_5^2/a_0.
$$

Then $P_6W^i_6(F)={\om}^{s\cdot i}_6\cdot \lambda \cdot  F$ in the $s$-th case with $\lambda = -a_5/a_0$.\\

This means that $G=\{  W^{2j}_6,\; P^k_6\cdot W^{i+2j}_6 \;\mid \; j=0,1,2,\; k=1,2,3,4,5   \}\cong \ZZ{6}\oplus \ZZ{3}$ with generators $W^2_6$ and $P_6\cdot W^i_6$ for some $i\in \{ 0,1 \}$ (in all five cases).\\

Now let us describe all smooth cubic fourfolds which admit this group action. In order to do this, we need to determine which forms $F'$ can appear in the equations above.\\

If a monomial $z_{p}z_{q}z_{r}$ is present in $F'$, then the invariance with respect to $W^2_6$ requires that $W^2_6(z_{p}z_{q}z_{r})={\om}^{2(p+q+r)}_6\cdot z_{p}z_{q}z_{r}=\mu \cdot z_{p}z_{q}z_{r}$, where the coefficient $\mu$ was specified above.\\

We conclude that $p+q+r\equiv \tau\;  mod \; 3$ has to be satisfied, where 
\begin{itemize}
\item $\tau\equiv 1\;  mod \; 3$ in the first and the fourth cases,
\item $\tau\equiv 2\;  mod \; 3$ in the second and the fifth cases,
\item $\tau\equiv 0\;  mod \; 3$ in the third case.
\end{itemize}

Since the set of such monomials is invariant with respect to $P_6$, we conclude that
$$
F'=\sum_{0\leq p\leq q\leq r\leq 5,\;\; p+q+r\equiv \tau\; mod \; 3}^{\prime} a_{pqr} \cdot z_{p}z_{q}z_{r}
$$
with suitable coefficients $a_{pqr}$. The prime signifies that the sum is over triples $p\leq q\leq r$ excluding those for which $p=q$ and $r\equiv p+s\; mod \; 6$ or $q=r$ and $p\equiv q+s\; mod \; 6$ in the $s$-th case (in order to exclude the monomials which were already written separately in $F$).\\

Coefficients $a_{pqr}$ can be determined as follows. We have to require that $P_6\cdot W^i_6(F')={\om}^{s\cdot i}_6 \cdot {\lambda}\cdot F'$ with the same $\lambda = -a_5/a_0$ as above (in the $s$-th case).\\

Let $S_1,S_2,...$ be the orbits of the natural action of $P_6$ on the set of monomials $z_{p}z_{q}z_{r}$ such that $p+q+r\equiv s \; mod \; 3$ (in the $s$-th case) and $p\leq q\leq r$ satisfy the exclusion condition above.\\

Then for every such orbit $S_{\alpha}$ the condition $P_6\cdot W^i_6(F')={\om}^{s\cdot i}_6 \cdot {\lambda}\cdot F'$ determines the coefficients $a_{pqr}$ of monomials in $S_{\alpha}$ upto a common factor.\\

Namely, if $m_1,m_2,...,m_k$ are elements of $S_{\alpha}$ such that $P_6(m_p)=m_{p+1}$, $p=1,2,...,k-1$, $P_6(m_k)=m_{1}$, $W_6(m_p)={\om}^{c_p}_6\cdot m_p$ (where $c_p\in \mathbb Z$) and $b_p$ is the coefficient of $m_p$ in the sum for $F'$ above, then
$$
b_k={\om}^{(s-c_k)\cdot i}_6 \cdot \lambda \cdot b_{1},\; \; b_p={\om}^{(s-c_p)\cdot i}_6 \cdot\lambda \cdot b_{p+1},\; \; p=1,2,...,k-1
$$
(in the $s$-th case).\\

This gives by induction:
$$
b_p=({\om}^{s\cdot i}_6 \cdot \lambda )^{k+1-p} \cdot {\om}^{-(c_p+...+c_k)\cdot i}_6\cdot b_{1},\; \; p=1,...,k
$$
(in the $s$-th case).\\

Note that these formulas determine $b_2,...,b_k$ in terms of $b_1$, which has to be zero unless $\lambda = -a_5/a_0$ satisfies the following condition:
$$
{\lambda}^k={\om}^{(c_1+...+c_k-s\cdot k)\cdot i}_6,
$$

i.e.
$$
(a_5/a_0)^k=(-1)^k\cdot {\om}^{(c_1+...+c_k-s\cdot k)\cdot i}_6
$$
(in the $s$-th case).\\

If this condition is satisfied, then the formulas above determine $b_1,b_2,...,b_k$ upto a common factor. If it is not, then the formulas above say that $b_1=b_2=...=b_k=0$, i.e. monomials $m_1,m_2,...,m_k$ from $S_{\alpha}$ do not appear in $F'$.\\

Note that since $P^6_6=Id$, only the following values of $k$ are possible: $k=1$, $k=2$, $k=3$, $k=6$.\\

In the case $k=1$ the orbit $S_{\alpha}$ consists of just one monomial $m_1$ such that $P_6(m_1)=m_1$ and $W_6(m_1)={\om}^{c_1}_6\cdot m_1$ for some $c_1\in \mathbb Z$. In this case the coefficient $b_1$ of this monomial in $F'$ is zero unless $\lambda ={\om}^{(c_1-s)\cdot i}_6$, i.e. $a_5=-{\om}^{(c_1-s)\cdot i}_6 \cdot a_0$ (in the $s$-th case). However, our earlier condition $(a_5/a_0)^6=-1$ is inconsistent with this requirement. Hence $b_1=0$, i.e. such monomials $m_1$ do not appear in $F'$.\\

This procedure determines coefficients $a_{pqr}$ in $F'$.\\

It is immediate that generic such cubic forms $F$ give smooth cubic fourfolds.\\

Note that the cubic fourfold given by the form
$$
F(z_0,z_1,z_2,z_3,z_4,z_5)=e_0z_0^2z_1+e_1z_1^2z_2+e_2z_2^2z_0+e_3z_3^2z_4+e_4z_4^2z_5+e_5z_5^2z_3
$$

is smooth for any $e_i\in {\mathbb C}^{*}$.\\

\vspace{2cm}

We can summarize our findings in this subsection in the following theorem.\\

{\bf Theorem 4.} {\it Abelian automorphism groups of smooth cubic fourfolds which are conjugate to subgroups of ${\mathcal P}_{6}\subset PGL(6)=Aut({\mathbb P}^5)$ are either diagonalizable or conjugate to subgroups of $D_3\otimes {\mathcal P}_{2}$ or conjugate to subgroups of $D_2\otimes {\mathcal P}_{3}$ or conjugate to the subgroup of ${\mathcal P}_{6}\subset PGL(6)$, which is abstractly isomorphic to $\ZZ{3}\oplus \ZZ{6}$ generated by matrices $W^2_6$ and $P_6\cdot W^i_6$ for some $i\in \{ 0,1 \}$.}\\

Smooth cubic fourfolds which admit the last group action are described earlier in the present subsection.\\

In particular, Theorem $4$ implies that the full Pauli subgroup ${\mathcal P}_{6}\subset PGL(6)$ does not act effectively on a smooth cubic fourfold.\\

\newpage

\section{Appendix. Smoothness.}

We used Macaulay 2, version 1.2.\\

We ran the following Macaulay 2 code:\\
 
{\scriptsize 

TotalSing=0;\\
R=QQ[x0,x1,x2,x3,x4,x5];\\
print "This is Case of Six Cubes";\\
F = x5\^{ }3 + x4\^{ }3 + x3\^{ }3 + x2\^{ }3 + x1\^{ }3 + x0\^{ }3; TotalSing=TotalSing + dim singularLocus(ideal(F));\\
print "This is Case of Five Cubes";\\
F = x5\^{ }2*x0 + x4\^{ }3 + x3\^{ }3 + x2\^{ }3 + x1\^{ }3 + x0\^{ }3; TotalSing=TotalSing + dim singularLocus(ideal(F));\\ 
print "This is Case of Four Cubes";\\
F = x5\^{ }2*x4 + x4\^{ }2*x5 + x3\^{ }3 + x2\^{ }3 + x1\^{ }3 + x0\^{ }3; TotalSing=TotalSing + dim singularLocus(ideal(F));\\
F = x5\^{ }2*x0 + x4\^{ }2*x5 + x3\^{ }3 + x2\^{ }3 + x1\^{ }3 + x0\^{ }3; TotalSing=TotalSing + dim singularLocus(ideal(F));\\ 
F = x5\^{ }2*x1 + x4\^{ }2*x0 + x3\^{ }3 + x2\^{ }3 + x1\^{ }3 + x0\^{ }3; TotalSing=TotalSing + dim singularLocus(ideal(F));\\ 
F = x5\^{ }2*x0 + x4\^{ }2*x0 + x3\^{ }3 + x2\^{ }3 + x1\^{ }3 + x0\^{ }3 + x1*x4*x5; TotalSing=TotalSing + dim singularLocus(ideal(F));\\ 
print "This is Case of Three Cubes";\\
F = x2\^{ }3 + x1\^{ }3 + x0\^{ }3 + x3\^{ }2*x4 + x4\^{ }2*x5 + x5\^{ }2*x3; TotalSing=TotalSing + dim singularLocus(ideal(F));\\ 
print "This is Case of Three Cubes, Length 2 longest cycle";\\
F = x2\^{ }3 + x1\^{ }3 + x0\^{ }3 + x3\^{ }2*x4 + x4\^{ }2*x3 + x5\^{ }2*x0; TotalSing=TotalSing + dim singularLocus(ideal(F));\\ 
F0 = x2\^{ }3 + 2*x1\^{ }3 + 3*x0\^{ }3 + 5*x3\^{ }2*x4 + 7*x4\^{ }2*x3 + 11*x5\^{ }2*x3;\\ 
F=F0 + 13*x0*x4*x4; TotalSing=TotalSing + dim singularLocus(ideal(F));\\ 
F=F0 + 13*x5*x5*x4; TotalSing=TotalSing + dim singularLocus(ideal(F));\\ 
F=F0 + 13*x0*x4*x5; TotalSing=TotalSing + dim singularLocus(ideal(F));\\ 
print "This is Case of Three Cubes, No cycles";\\
F = x2\^{ }3 + x1\^{ }3 + x0\^{ }3 + x3\^{ }2*x4 + x4\^{ }2*x5 + x5\^{ }2*x0; TotalSing=TotalSing + dim singularLocus(ideal(F));\\ 
F = x2\^{ }3 + x1\^{ }3 + x0\^{ }3 + x3\^{ }2*x4 + x4\^{ }2*x0 + x5\^{ }2*x1; TotalSing=TotalSing + dim singularLocus(ideal(F));\\ 
F0 = x2\^{ }3 + 2*x1\^{ }3 + 3*x0\^{ }3 + 5*x3\^{ }2*x4 + 7*x4\^{ }2*x0 + 11*x5\^{ }2*x0;\\ 
F=F0 + 13*x3*x4*x5; TotalSing=TotalSing + dim singularLocus(ideal(F));\\ 
F=F0 + 13*x1*x4*x5; TotalSing=TotalSing + dim singularLocus(ideal(F));\\ 
F=F0 + 13*x4*x4*x5; TotalSing=TotalSing + dim singularLocus(ideal(F));\\ 
F0 = x2\^{ }3 + 2*x1\^{ }3 + 3*x0\^{ }3 + 5*x3\^{ }2*x4 + 7*x4\^{ }2*x0 + 11*x5\^{ }2*x4; \\
F=F0 + 13*x1*x3*x5; TotalSing=TotalSing + dim singularLocus(ideal(F));\\ 
F=F0 + 13*x0*x3*x5; TotalSing=TotalSing + dim singularLocus(ideal(F));\\ 
F = x2\^{ }3 + x1\^{ }3 + x0\^{ }3 + x3\^{ }2*x0 + x4\^{ }2*x1 + x5\^{ }2*x2; TotalSing=TotalSing + dim singularLocus(ideal(F));\\ 
F0 = x2\^{ }3 + 2*x1\^{ }3 + 3*x0\^{ }3 + 5*x3\^{ }2*x0 + 7*x4\^{ }2*x0 + 11*x5\^{ }2*x1; \\
F=F0 + 13*x3*x4*x5; TotalSing=TotalSing + dim singularLocus(ideal(F));\\ 
F=F0 + 13*x2*x3*x5 + 17*x1*x4*x4; TotalSing=TotalSing + dim singularLocus(ideal(F));\\ 
F=F0 + 13*x2*x4*x5 + 17*x1*x3*x4; TotalSing=TotalSing + dim singularLocus(ideal(F));\\ 
F=F0 + 13*x2*x3*x4; TotalSing=TotalSing + dim singularLocus(ideal(F));\\ 
F0 = x2\^{ }3 + 2*x1\^{ }3 + 3*x0\^{ }3 + 5*x3\^{ }2*x0 + 7*x4\^{ }2*x0 + 11*x5\^{ }2*x0; \\
F=F0 + 13*x3*x4*x5; TotalSing=TotalSing + dim singularLocus(ideal(F));\\ 
F=F0 + 13*x1*x4*x5 + 17*x1*x3*x4 + 23*x2*x3*x5; TotalSing=TotalSing + dim singularLocus(ideal(F));\\ 
print "This is Case of Two Cubes";\\
F = x1\^{ }3 + x0\^{ }3 + x2\^{ }2*x3 + x3\^{ }2*x4 + x4\^{ }2*x5 + x5\^{ }2*x2; TotalSing=TotalSing + dim singularLocus(ideal(F));\\ 
print "This is Case of Two Cubes, Length 3 longest cycle";\\
F = x1\^{ }3 + x0\^{ }3 + x2\^{ }2*x3 + x3\^{ }2*x4 + x4\^{ }2*x2 + x5\^{ }2*x0; TotalSing=TotalSing + dim singularLocus(ideal(F));\\ 
F0 = x1\^{ }3 + x0\^{ }3 + x2\^{ }2*x3 + x3\^{ }2*x4 + x4\^{ }2*x2 + x5\^{ }2*x2; \\
F=F0 + x5*x3*x4; TotalSing=TotalSing + dim singularLocus(ideal(F));\\ 
F=F0 + x5*x0*x4; TotalSing=TotalSing + dim singularLocus(ideal(F));\\ 
F=F0 + x4*x4*x0; TotalSing=TotalSing + dim singularLocus(ideal(F));\\ 
F=F0 + x4*x4*x3; TotalSing=TotalSing + dim singularLocus(ideal(F));\\ 
F=F0 + x4*x5*x5; TotalSing=TotalSing + dim singularLocus(ideal(F));\\ 
F=F0 + x3*x5*x5; TotalSing=TotalSing + dim singularLocus(ideal(F));\\ 
print "This is Case of Two Cubes, Length 2 longest cycle";\\
F1=[ x3\^{ }2*x0, x0*x3*x4, x3*x4*x5 ];\\
F2=[ x3\^{ }2*x1, x3\^{ }2*x4, x1*x3*x5, x3*x4*x5 ];\\
F3=[ x3\^{ }2*x0, x5\^{ }2*x0, x0*x3*x5 ];\\
F4=[ x4\^{ }2*x1, x4\^{ }2*x3, x1*x3*x4, x1*x4*x5 ];\\
F5=[ x4\^{ }2*x2, x2\^{ }2*x0, x0*x2*x4, x2*x4*x5 ];\\
F6=[ x2\^{ }2*x1, x2\^{ }2*x0, x0*x2*x4, x1*x2*x4, x2*x4*x5 ];\\
F7=[ x0*x2*x4, x2*x4*x5 ];\\
F8=[ x0*x2*x4, x1*x2*x4, x2*x4*x5 ];\\
F9=[ x3\^{ }2*x1, x3\^{ }2*x4, x1*x3*x4, x3*x4*x5 ];\\
F = x1\^{ }3 + x0\^{ }3 + x2\^{ }2*x3 + x3\^{ }2*x2 + x4\^{ }2*x0 + x5\^{ }2*x1; TotalSing=TotalSing + dim singularLocus(ideal(F));\\ 
F = x1\^{ }3 + x0\^{ }3 + x2\^{ }2*x3 + x3\^{ }2*x2 + x4\^{ }2*x5 + x5\^{ }2*x4; TotalSing=TotalSing + dim singularLocus(ideal(F));\\ 
F = x1\^{ }3 + x0\^{ }3 + x2\^{ }2*x3 + x3\^{ }2*x2 + x4\^{ }2*x0 + x5\^{ }2*x4; TotalSing=TotalSing + dim singularLocus(ideal(F));\\ 
F = x1\^{ }3 + 2*x0\^{ }3 + 3*x2\^{ }2*x3 + 5*x3\^{ }2*x2 + 7*x4\^{ }2*x2  + 11*x4\^{ }2*x3 + 13*x5\^{ }2*x4;\\ TotalSing=TotalSing + dim singularLocus(ideal(F));\\ 
F = x1\^{ }3 + 2*x0\^{ }3 + 3*x2\^{ }2*x3 + 5*x3\^{ }2*x2 + 7*x5\^{ }2*x2  + 11*x5\^{ }2*x3 + 13*x4\^{ }2*x0;\\ TotalSing=TotalSing + dim singularLocus(ideal(F));\\ 
F0 = x1\^{ }3 + x0\^{ }3 + 2*x2\^{ }2*x3 + 3*x3\^{ }2*x2 + 7*x4\^{ }2*x2 + 11*x5\^{ }2*x4; \\
for iF1 from 0 to (\#F1-1) list (\\
F=F0 + 19*F1\#iF1; TotalSing=TotalSing + dim singularLocus(ideal(F));\\ 
);\\
F0 = x1\^{ }3 + x0\^{ }3 + 2*x2\^{ }2*x3 + 3*x3\^{ }2*x2 + 7*x4\^{ }2*x0 + 11*x5\^{ }2*x2; \\
for iF2 from 0 to (\#F2-1) list (\\
F=F0 + 19*F2\#iF2; TotalSing=TotalSing + dim singularLocus(ideal(F));\\ 
);\\
for iF3 from 0 to (\#F3-1) list (\\
for iF4 from 0 to (\#F4-1) list (\\
F=F0 + 19*F3\#iF3 + 23*F4\#iF4; TotalSing=TotalSing + dim singularLocus(ideal(F));\\ 
));\\
F00 = x1\^{ }3 + x0\^{ }3 + 2*x2\^{ }2*x3 + 3*x3\^{ }2*x2 + 7*x4\^{ }2*x3 + 11*x5\^{ }2*x2; \\
F0=F00 + 23*x5\^{ }2*x3; \\
for iF5 from 0 to (\#F5-1) list ( \\
F=F0 + 19*F5\#iF5; TotalSing=TotalSing + dim singularLocus(ideal(F));\\ 
);\\
F0=F00 + 23*x3\^{ }2*x0; \\
for iF6 from 0 to (\#F6-1) list ( \\
F=F0 + 19*F6\#iF6; TotalSing=TotalSing + dim singularLocus(ideal(F));\\ 
);\\
F0=F00 + 23*x3*x4*x5; \\
for iF7 from 0 to (\#F7-1) list ( \\
F=F0 + 19*F7\#iF7; TotalSing=TotalSing + dim singularLocus(ideal(F));\\ 
);\\
F0=F00 + 23*x3*x0*x5; \\
for iF8 from 0 to (\#F8-1) list ( \\
F=F0 + 19*F8\#iF8; TotalSing=TotalSing + dim singularLocus(ideal(F));\\ 
);\\
F00 = x1\^{ }3 + x0\^{ }3 + 2*x2\^{ }2*x3 + 3*x3\^{ }2*x2 + 7*x4\^{ }2*x0 + 11*x5\^{ }2*x0; \\
F0=F00 + 23*x2*x4*x5; \\
for iF9 from 0 to (\#F9-1) list ( \\
F=F0 + 19*F9\#iF9; TotalSing=TotalSing + dim singularLocus(ideal(F));\\ 
);\\
F=F00 + 19*x1*x4*x5; TotalSing=TotalSing + dim singularLocus(ideal(F));\\ 
F0 = x1\^{ }3 + x0\^{ }3 + 2*x2\^{ }2*x3 + 3*x3\^{ }2*x2 + 7*x4\^{ }2*x2 + 11*x5\^{ }2*x2;\\
F=F0 + 23*x3*x4*x5; TotalSing=TotalSing + dim singularLocus(ideal(F));\\ 
F=F0 + 23*x0*x4*x5 + 19*x1*x3*x4 + 29*x3*x3*x0; TotalSing=TotalSing + dim singularLocus(ideal(F));\\ 
F=F0 + 23*x0*x4*x5 + 19*x3*x3*x1; TotalSing=TotalSing + dim singularLocus(ideal(F));\\ 
F=F0 + 23*x0*x4*x5 + 19*x4*x3*x1 + 29*x5*x3*x1; TotalSing=TotalSing + dim singularLocus(ideal(F));\\ 
F=F0 + 23*x0*x4*x5 + 19*x4*x3*x1 + 29*x5*x3*x0; TotalSing=TotalSing + dim singularLocus(ideal(F));\\ 
F=F0 + 23*x1*x4*x5 + 19*x4*x3*x0 + 29*x5*x3*x0 + 31*x0*x4*x5; TotalSing=TotalSing + dim singularLocus(ideal(F));\\ 
print "This is Case of Two Cubes, No cycles";\\
F1=[ x1*x2*x4, x1*x4*x5 ];\\
F2=[ x2\^{ }2*x1, x2\^{ }2*x4, x2\^{ }2*x5, x1*x2*x4, x1*x2*x5 ];\\
F3=[ x2\^{ }2*x1, x1*x2*x5 ];\\
F4=[ x2\^{ }2*x1, x5\^{ }2*x1, x1*x2*x5 ];\\
F5=[ x3\^{ }2*x1, x5\^{ }2*x1, x1*x3*x5 ];\\
F6=[ x4\^{ }2*x3, x2*x3*x4, x2*x4*x5 ];\\
F7=[ x2*x3*x5, x3*x4*x5 ];\\
F8=[ x2\^{ }2*x0, x2\^{ }2*x4, x0*x2*x4 ];\\
F9=[ x2\^{ }2*x0, x4\^{ }2*x0, x0*x2*x4 ];\\
F10=[ x5\^{ }2*x3, x1*x3*x5, x2*x3*x5, x3*x4*x5 ];\\
A=[ x5\^{ }2*x3, x1*x3*x5, x2*x3*x5 ];\\
F11=[ x1*x2*x5, x0*x2*x5 ];\\
F12=[ x1*x3*x4, x2*x3*x4 ];\\
F13=[ x2\^{ }2*x1, x2\^{ }2*x4, x2\^{ }2*x5, x1*x2*x4, x1*x2*x5, x1*x4*x5, x2*x4*x5 ];\\
F14=[ x1*x3*x4, x2*x3*x4 ];\\
F15=[ x1*x2*x3, x2*x3*x5 ];\\
F = x1\^{ }3 + x0\^{ }3 + x2\^{ }2*x3 + x3\^{ }2*x4 + x4\^{ }2*x5 + x5\^{ }2*x0; TotalSing=TotalSing + dim singularLocus(ideal(F));\\ 
F = x1\^{ }3 + x0\^{ }3 + x2\^{ }2*x3 + x3\^{ }2*x4 + x4\^{ }2*x0 + x5\^{ }2*x1; TotalSing=TotalSing + dim singularLocus(ideal(F));\\ 
F = x1\^{ }3 + x0\^{ }3 + x2\^{ }2*x3 + x3\^{ }2*x0 + x4\^{ }2*x5 + x5\^{ }2*x1; TotalSing=TotalSing + dim singularLocus(ideal(F));\\ 
F00 = x1\^{ }3 + x0\^{ }3 + 2*x2\^{ }2*x3 + 3*x3\^{ }2*x4 + 7*x4\^{ }2*x0 + 11*x5\^{ }2*x3; \\
F=F00 + 19*x5*x4*x1 + 17*x5*x2*x0; TotalSing=TotalSing + dim singularLocus(ideal(F));\\ 
F=F00 + 19*x5*x5*x4; TotalSing=TotalSing + dim singularLocus(ideal(F));\\ 
F=F00 + 19*x2*x5*x4; \\
TotalSing=TotalSing + dim singularLocus(ideal(F));\\ 
F=F00 + 19*x2*x5*x1; TotalSing=TotalSing + dim singularLocus(ideal(F));\\ 
F0=F00 + 23*x0*x5*x5;\\
for iF1 from 0 to (\#F1-1) list ( \\
F=F0 + 19*F1\#iF1; TotalSing=TotalSing + dim singularLocus(ideal(F));\\ 
);\\
F00 = x1\^{ }3 + x0\^{ }3 + 2*x2\^{ }2*x3 + 3*x3\^{ }2*x4 + 7*x4\^{ }2*x0 + 11*x5\^{ }2*x4; \\
F=F00 + 19*x2*x3*x5; TotalSing=TotalSing + dim singularLocus(ideal(F));\\ 
F=F00 + 19*x0*x5*x5; TotalSing=TotalSing + dim singularLocus(ideal(F));\\ 
F=F00 + 19*x0*x3*x3; TotalSing=TotalSing + dim singularLocus(ideal(F));\\ 
F=F00 + 19*x0*x3*x5; TotalSing=TotalSing + dim singularLocus(ideal(F));\\ 
F=F00 + 19*x1*x3*x5; TotalSing=TotalSing + dim singularLocus(ideal(F));\\ 
F00 = x1\^{ }3 + x0\^{ }3 + 2*x2\^{ }2*x3 + 3*x3\^{ }2*x4 + 7*x4\^{ }2*x0 + 11*x5\^{ }2*x0; \\
F0=F00 + 23*x3*x4*x5;\\
for iF2 from 0 to (\#F2-1) list ( \\
F=F0 + 19*F2\#iF2; TotalSing=TotalSing + dim singularLocus(ideal(F));\\ 
);\\
F=F00 + 19*x2*x4*x5; TotalSing=TotalSing + dim singularLocus(ideal(F));\\ 
F=F00 + 19*x1*x4*x5; TotalSing=TotalSing + dim singularLocus(ideal(F));\\ 
F00 = x1\^{ }3 + x0\^{ }3 + 2*x2\^{ }2*x3 + 3*x3\^{ }2*x0 + 7*x4\^{ }2*x5 + 11*x5\^{ }2*x0;\\
F=F00 + 19*x1*x3*x5; TotalSing=TotalSing + dim singularLocus(ideal(F));\\ 
F=F00 + 19*x3*x4*x5; \\
TotalSing=TotalSing + dim singularLocus(ideal(F));\\ 
F00 = x1\^{ }3 + x0\^{ }3 + 2*x2\^{ }2*x3 + 3*x3\^{ }2*x0 + 7*x4\^{ }2*x1 + 11*x5\^{ }2*x3;\\
F=F00 + 19*x0*x2*x5; TotalSing=TotalSing + dim singularLocus(ideal(F));\\ 
F=F00 + 19*x4*x2*x5; TotalSing=TotalSing + dim singularLocus(ideal(F));\\ 
F=F00 + 19*x0*x5*x5; TotalSing=TotalSing + dim singularLocus(ideal(F));\\ 
F0=F00 + 23*x0*x4*x4;\\
for iF3 from 0 to (\#F3-1) list ( \\
F=F0 + 19*F3\#iF3; TotalSing=TotalSing + dim singularLocus(ideal(F));\\ 
);\\
F0=F00 + 23*x0*x4*x5;\\
for iF4 from 0 to (\#F4-1) list ( \\
F=F0 + 19*F4\#iF4; TotalSing=TotalSing + dim singularLocus(ideal(F));\\ 
);\\
F00 = x1\^{ }3 + x0\^{ }3 + 2*x2\^{ }2*x3 + 3*x3\^{ }2*x0 + 7*x4\^{ }2*x1 + 11*x5\^{ }2*x0;\\
F=F00 + 19*x3*x4*x5; TotalSing=TotalSing + dim singularLocus(ideal(F));\\ 
F=F00 + 19*x2*x3*x5; TotalSing=TotalSing + dim singularLocus(ideal(F));\\ 
F0=F00;\\
for iF5 from 0 to (\#F5-1) list ( \\
for iF6 from 0 to (\#F6-1) list ( \\
F=F0 + 19*F5\#iF5 + 23*F6\#iF6; TotalSing=TotalSing + dim singularLocus(ideal(F));\\ 
));\\
F00 = x1\^{ }3 + x0\^{ }3 + 2*x2\^{ }2*x3 + 3*x3\^{ }2*x0 + 7*x4\^{ }2*x1 + 11*x5\^{ }2*x1;\\
F=F00 + 19*x2*x4*x5; TotalSing=TotalSing + dim singularLocus(ideal(F));\\ 
F0=F00 + 23*x0*x4*x5;\\
for iF7 from 0 to (\#F7-1) list ( \\
F=F0 + 19*F7\#iF7; TotalSing=TotalSing + dim singularLocus(ideal(F));\\ 
);\\
F0=F00 + 23*x3*x4*x5;\\
for iF8 from 0 to (\#F8-1) list ( \\
F=F0 + 19*F8\#iF8; TotalSing=TotalSing + dim singularLocus(ideal(F));\\ 
);\\
F00 = x1\^{ }3 + x0\^{ }3 + 2*x2\^{ }2*x3 + 3*x3\^{ }2*x0 + 7*x4\^{ }2*x3 + 11*x5\^{ }2*x0;\\
F0=F00 + 23*x1*x4*x5;\\
for iF9 from 0 to (\#F9-1) list ( \\
for iF10 from 0 to (\#F10-1) list ( \\
F=F0 + 19*F9\#iF9 + 23*F10\#iF10; TotalSing=TotalSing + dim singularLocus(ideal(F));\\ 
));\\
F0=F00 + 23*x2*x4*x5;\\
for iA from 0 to (\#A-1) list ( \\
F=F0 + 19*A\#iA; TotalSing=TotalSing + dim singularLocus(ideal(F));\\ 
);\\
F0=F00 + 23*x2*x4*x1;\\
for iA from 0 to (\#A-1) list ( \\
F=F0 + 19*A\#iA; TotalSing=TotalSing + dim singularLocus(ideal(F));\\ 
);\\
F00 = x1\^{ }3 + x0\^{ }3 + 2*x2\^{ }2*x3 + 3*x3\^{ }2*x0 + 7*x4\^{ }2*x3 + 11*x5\^{ }2*x3;\\
F0=F00 + 23*x0*x2*x4 + 31*x1*x4*x5;\\
for iF11 from 0 to (\#F11-1) list ( \\
F=F0 + 19*F11\#iF11; TotalSing=TotalSing + dim singularLocus(ideal(F));\\ 
);\\
F=F00 + 23*x2*x4*x5; TotalSing=TotalSing + dim singularLocus(ideal(F));\\ 
F=F00 + 23*x2*x4*x0 + 23*x2*x4*x1 + 23*x2*x5*x1 + 23*x1*x4*x5; TotalSing=TotalSing + dim singularLocus(ideal(F));\\ 
F00 = x1\^{ }3 + x0\^{ }3 + 2*x2\^{ }2*x3 + 3*x3\^{ }2*x0 + 7*x4\^{ }2*x0 + 11*x5\^{ }2*x0;\\
F0=F00 + 23*x5*x2*x4 + 31*x1*x3*x5;\\
for iF12 from 0 to (\#F12-1) list ( \\
F=F0 + 19*F12\#iF12; TotalSing=TotalSing + dim singularLocus(ideal(F));\\ 
);\\
F0=F00 + 23*x5*x2*x3 + 31*x1*x4*x5;\\
for iF14 from 0 to (\#F14-1) list ( \\
F=F0 + 19*F14\#iF14; TotalSing=TotalSing + dim singularLocus(ideal(F));\\ 
);\\
F0=F00 + 23*x3*x4*x5;\\
for iF13 from 0 to (\#F13-1) list ( \\
F=F0 + 19*F13\#iF13; TotalSing=TotalSing + dim singularLocus(ideal(F));\\ 
);\\
F00 = x1\^{ }3 + x0\^{ }3 + 2*x2\^{ }2*x0 + 3*x3\^{ }2*x0 + 7*x4\^{ }2*x1 + 11*x5\^{ }2*x1;\\
F=F00 + 19*x3*x4*x5 + 23*x3*x3*x1; TotalSing=TotalSing + dim singularLocus(ideal(F));\\ 
F=F00 + 19*x2*x4*x5 + 23*x3*x3*x1; TotalSing=TotalSing + dim singularLocus(ideal(F));\\ 
F=F00 + 19*x2*x4*x5 + 23*x2*x3*x4; TotalSing=TotalSing + dim singularLocus(ideal(F));\\ 
F=F00 + 19*x2*x4*x5 + 23*x1*x2*x3; TotalSing=TotalSing + dim singularLocus(ideal(F));\\ 
F00 = x1\^{ }3 + x0\^{ }3 + 2*x2\^{ }2*x0 + 3*x3\^{ }2*x0 + 7*x4\^{ }2*x0 + 11*x5\^{ }2*x1;\\
F=F00 + 19*x2*x3*x4; TotalSing=TotalSing + dim singularLocus(ideal(F));\\ 
F0=F00 + 19*x2*x4*x5 + 23*x4*x3*x1; \\
for iF15 from 0 to (\#F15-1) list ( \\
F=F0 + 19*F15\#iF15; TotalSing=TotalSing + dim singularLocus(ideal(F));\\ 
);\\
F00 = x1\^{ }3 + x0\^{ }3 + 2*x2\^{ }2*x0 + 3*x3\^{ }2*x0 + 7*x4\^{ }2*x0 + 11*x5\^{ }2*x0;\\
F=F00 + 19*x2*x3*x4 + 17*x1*x3*x5 + 23*x1*x2*x5 + 29*x1*x4*x5; TotalSing=TotalSing + dim singularLocus(ideal(F));\\ 
F=F00 + 19*x2*x3*x4 + 23*x2*x4*x5 + 17*x2*x3*x5 + 29*x3*x4*x5; TotalSing=TotalSing + dim singularLocus(ideal(F));\\ 
F=F00 + 19*x2*x3*x4 + 23*x1*x2*x5 + 17*x3*x4*x5; TotalSing=TotalSing + dim singularLocus(ideal(F));\\ 
print "This is Case of One Cube";\\
F = x0\^{ }3 + x1\^{ }2*x2 + x2\^{ }2*x3 + x3\^{ }2*x4 + x4\^{ }2*x5 + x5\^{ }2*x1; TotalSing=TotalSing + dim singularLocus(ideal(F));\\ 
F = x0\^{ }3 + x1\^{ }2*x2 + x2\^{ }2*x3 + x3\^{ }2*x4 + x4\^{ }2*x1 + x5\^{ }2*x0; TotalSing=TotalSing + dim singularLocus(ideal(F));\\ 
print "This is Case of One Cube, Length 4 longest cycle";\\
A=[ x2\^{ }2*x0, x2\^{ }2*x4, x2\^{ }2*x5, x0*x2*x4, x0*x2*x5 ];\\
F1=[ x0*x4*x5, x2*x4*x5 ];\\
F00 = x0\^{ }3 + 2*x1\^{ }2*x2 + 3*x2\^{ }2*x3 + 7*x3\^{ }2*x4 + 11*x4\^{ }2*x1 + 13*x5\^{ }2*x1; \\
F0=F00;\\
for iF1 from 0 to (\#F1-1) list ( \\
F=F0 + 19*F1\#iF1; TotalSing=TotalSing + dim singularLocus(ideal(F));\\ 
); \\
F0=F00 + 17*x3*x4*x5;\\
for iA from 0 to (\#A-1) list ( \\
F=F0 + 19*A\#iA; TotalSing=TotalSing + dim singularLocus(ideal(F));\\ 
);\\
F0=F00 + 17*x3*x4*x4;\\
for iA from 0 to (\#A-1) list ( \\
F=F0 + 19*A\#iA; TotalSing=TotalSing + dim singularLocus(ideal(F));\\ 
);\\
F0=F00 + 17*x3*x5*x5;\\
for iA from 0 to (\#A-1) list ( \\
F=F0 + 19*A\#iA; TotalSing=TotalSing + dim singularLocus(ideal(F));\\ 
);\\
F=F00 + 17*x4*x4*x0; TotalSing=TotalSing + dim singularLocus(ideal(F));\\
F=F00 + 17*x4*x5*x5; TotalSing=TotalSing + dim singularLocus(ideal(F));\\ 
F=F00 + 17*x4*x4*x2; TotalSing=TotalSing + dim singularLocus(ideal(F));\\ 
print "This is Case of One Cube, Length 3 longest cycle";\\
A0=[ x5\^{ }2*x2, x5\^{ }2*x3, x2*x3*x5, x2*x4*x5 ];\\
A=[ x4\^{ }2*x1, x4\^{ }2*x3, x1\^{ }2*x0, x1\^{ }2*x3, x1*x3*x4, x1*x4*x5, x1*x4*x0 ];\\
B=[ x4\^{ }2*x1, x5\^{ }2*x1, x1\^{ }2*x0, x0*x1*x4, x0*x1*x5, x0*x4*x5, x1*x4*x5 ];\\
F1=[ x2\^{ }2*x4, x2\^{ }2*x0, x0*x2*x4, x0*x4*x5, x2*x4*x5 ];\\
F2=[ x3\^{ }2*x2, x3\^{ }2*x0, x0*x3*x4, x2*x3*x4 ];\\
F3=[ x3\^{ }2*x2, x3\^{ }2*x0, x0*x3*x5, x2*x3*x5 ];\\
F4=[ x3\^{ }2*x2, x2*x3*x4, x2*x4*x5 ];\\
F5=[ x1\^{ }2*x4, x1\^{ }2*x3, x1*x3*x4, x1*x4*x5, x3*x4*x5 ];\\
F6=[ x3\^{ }2*x2, x2*x3*x4, x0*x3*x4 ];\\
F = x0\^{ }3 + 2*x1\^{ }2*x2 + 3*x2\^{ }2*x3 + 7*x3\^{ }2*x1 + 11*x4\^{ }2*x5 + 13*x5\^{ }2*x0;\\ TotalSing=TotalSing + dim singularLocus(ideal(F));\\ 
F = x0\^{ }3 + 2*x1\^{ }2*x2 + 3*x2\^{ }2*x3 + 7*x3\^{ }2*x1 + 11*x4\^{ }2*x5 + 13*x5\^{ }2*x4;\\ TotalSing=TotalSing + dim singularLocus(ideal(F));\\ 
F0 = x0\^{ }3 + 2*x1\^{ }2*x2 + 3*x2\^{ }2*x3 + 7*x3\^{ }2*x1 + 11*x4\^{ }2*x5 + 13*x5\^{ }2*x3; \\
F=F0 + 17*x1*x2*x5; TotalSing=TotalSing + dim singularLocus(ideal(F));\\ 
F=F0 + 17*x0*x2*x5; TotalSing=TotalSing + dim singularLocus(ideal(F));\\ 
F=F0 + 17*x4*x2*x5; TotalSing=TotalSing + dim singularLocus(ideal(F));\\ 
F=F0 + 17*x5*x5*x2; TotalSing=TotalSing + dim singularLocus(ideal(F));\\ 
F=F0 + 17*x0*x2*x2; TotalSing=TotalSing + dim singularLocus(ideal(F));\\ 
F=F0 + 17*x1*x2*x2; TotalSing=TotalSing + dim singularLocus(ideal(F));\\ 
F00 = x0\^{ }3 + 2*x1\^{ }2*x2 + 3*x2\^{ }2*x3 + 7*x3\^{ }2*x1 + 11*x4\^{ }2*x1 + 13*x5\^{ }2*x0; \\
F0=F00 + 17*x3*x3*x0;\\
for iA0 from 0 to (\#A0-1) list ( \\
F=F0 + 19*A0\#iA0; TotalSing=TotalSing + dim singularLocus(ideal(F));\\ 
); \\
F0=F00 + 17*x4*x3*x0;\\
for iA0 from 0 to (\#A0-1) list ( \\
F=F0 + 19*A0\#iA0; TotalSing=TotalSing + dim singularLocus(ideal(F));\\ 
); \\
F0=F00 + 17*x4*x4*x0;\\
for iA0 from 0 to (\#A0-1) list ( \\
F=F0 + 19*A0\#iA0; TotalSing=TotalSing + dim singularLocus(ideal(F));\\ 
); \\
F=F00 + 17*x3*x4*x5; TotalSing=TotalSing + dim singularLocus(ideal(F));\\ 
F=F00 + 17*x3*x3*x5; TotalSing=TotalSing + dim singularLocus(ideal(F));\\ 
F=F00 + 17*x2*x3*x4; TotalSing=TotalSing + dim singularLocus(ideal(F));\\ 
F=F00 + 17*x2*x3*x3; TotalSing=TotalSing + dim singularLocus(ideal(F));\\ 
F=F00 + 17*x2*x4*x4; TotalSing=TotalSing + dim singularLocus(ideal(F));\\ 
F00 = x0\^{ }3 + 2*x1\^{ }2*x2 + 3*x2\^{ }2*x3 + 7*x3\^{ }2*x1 + 11*x4\^{ }2*x2 + 13*x5\^{ }2*x3; \\
F0=F00 + 17*x5*x5*x2;\\
for iA from 0 to (\#A-1) list ( \\
for iB from 0 to (\#B-1) list ( \\
F=F0 + 19*A\#iA + 23*B\#iB; TotalSing=TotalSing + dim singularLocus(ideal(F));\\ 
));\\
F0=F00 + 17*x5*x4*x2;\\
for iA from 0 to (\#A-1) list ( \\
for iB from 0 to (\#B-1) list ( \\
F=F0 + 19*A\#iA + 23*B\#iB; TotalSing=TotalSing + dim singularLocus(ideal(F));\\ 
));\\
F0=F00 + 17*x4*x2*x2;\\
for iA from 0 to (\#A-1) list ( \\
for iB from 0 to (\#B-1) list ( \\
F=F0 + 19*A\#iA + 23*B\#iB; TotalSing=TotalSing + dim singularLocus(ideal(F));\\ 
));\\
F0=F00 + 17*x1*x2*x2;\\
for iA from 0 to (\#A-1) list ( \\
for iB from 0 to (\#B-1) list ( \\
F=F0 + 19*A\#iA + 23*B\#iB; TotalSing=TotalSing + dim singularLocus(ideal(F));\\ 
));\\
F0=F00 + 17*x1*x2*x5;\\
for iA from 0 to (\#A-1) list ( \\
for iB from 0 to (\#B-1) list ( \\
F=F0 + 19*A\#iA + 23*B\#iB; TotalSing=TotalSing + dim singularLocus(ideal(F));\\ 
));\\
F0=F00 + 17*x2*x2*x0;\\
for iA from 0 to (\#A-1) list ( \\
for iB from 0 to (\#B-1) list ( \\
F=F0 + 19*A\#iA + 23*B\#iB; TotalSing=TotalSing + dim singularLocus(ideal(F));\\ 
));\\
F0=F00 + 17*x5*x2*x0;\\
for iA from 0 to (\#A-1) list ( \\
for iB from 0 to (\#B-1) list ( \\
F=F0 + 19*A\#iA + 23*B\#iB; TotalSing=TotalSing + dim singularLocus(ideal(F));\\ 
));\\
F0=F00 + 17*x5*x5*x1;\\
for iA from 0 to (\#A-1) list ( \\
F=F0 + 19*A\#iA; TotalSing=TotalSing + dim singularLocus(ideal(F));\\ 
); \\
F00 = x0\^{ }3 + 2*x1\^{ }2*x2 + 3*x2\^{ }2*x3 + 7*x3\^{ }2*x1 + 11*x4\^{ }2*x1 + 13*x5\^{ }2*x1; \\
F0=F00 + 17*x5*x4*x0;\\
for iF2 from 0 to (\#F2-1) list ( \\
for iF3 from 0 to (\#F3-1) list ( \\
for iF4 from 0 to (\#F4-1) list ( \\
F=F0 + 19*F2\#iF2 + 23*F3\#iF3 + 29*F4\#iF4; TotalSing=TotalSing + dim singularLocus(ideal(F));\\ 
)));\\
F0=F00 + 17*x5*x3*x0 + 31*x2*x4*x5;\\
for iF6 from 0 to (\#F6-1) list ( \\
F=F0 + 19*F6\#iF6; TotalSing=TotalSing + dim singularLocus(ideal(F));\\ 
); \\
F0=F00 + 17*x3*x4*x5;\\
for iF1 from 0 to (\#F1-1) list ( \\
F=F0 + 19*F1\#iF1; TotalSing=TotalSing + dim singularLocus(ideal(F));\\ 
); \\
F=F00 + 17*x2*x4*x5 + 19*x3*x3*x0; TotalSing=TotalSing + dim singularLocus(ideal(F));\\ 
F00 = x0\^{ }3 + 2*x1\^{ }2*x2 + 3*x2\^{ }2*x3 + 7*x3\^{ }2*x1 + 11*x4\^{ }2*x0 + 13*x5\^{ }2*x0; \\
F0=F00 + 17*x2*x4*x5;\\
for iF5 from 0 to (\#F5-1) list ( \\
F=F0 + 19*F5\#iF5; TotalSing=TotalSing + dim singularLocus(ideal(F));\\ 
); \\
print "This is Case of One Cube, Length 2 longest cycle";\\
A=[ x3\^{ }2*x2, x5\^{ }2*x2, x2\^{ }2*x0, x0*x2*x3, x0*x2*x5, x2*x3*x5 ];\\
B=[ x0*x2*x5, x2*x3*x5 ];\\
C=[ x2\^{ }2*x3, x2\^{ }2*x5, x2\^{ }2*x0, x0*x2*x3, x0*x2*x5, x2*x3*x5 ];\\
F1=[ x3\^{ }2*x2, x3\^{ }2*x5, x3\^{ }2*x0, x0*x2*x3, x0*x3*x5 ];\\
AA=[ x3\^{ }2*x2, x3\^{ }2*x0, x0*x2*x3, x0*x3*x5 ];\\
BB=[ x2\^{ }2*x4, x5\^{ }2*x4, x4\^{ }2*x2, x2*x3*x4, x2*x4*x5, x3*x4*x5 ];\\
AAA=[ x5\^{ }2*x2, x5\^{ }2*x4, x2*x3*x5, x3*x4*x5 ];\\
F = x0\^{ }3 + 2*x1\^{ }2*x2 + 3*x2\^{ }2*x1 + 7*x3\^{ }2*x4 + 11*x4\^{ }2*x3 + 13*x5\^{ }2*x0;\\ TotalSing=TotalSing + dim singularLocus(ideal(F));\\ 
F = x0\^{ }3 + 2*x1\^{ }2*x2 + 3*x2\^{ }2*x1 + 7*x3\^{ }2*x4 + 11*x4\^{ }2*x5 + 13*x5\^{ }2*x0;\\ TotalSing=TotalSing + dim singularLocus(ideal(F));\\ 
F00 = x0\^{ }3 + 2*x1\^{ }2*x2 + 3*x2\^{ }2*x1 + 7*x3\^{ }2*x4 + 11*x4\^{ }2*x5 + 13*x5\^{ }2*x4; \\
F=F00 + 17*x3*x3*x5; TotalSing=TotalSing + dim singularLocus(ideal(F));\\ 
F=F00 + 17*x0*x5*x5; TotalSing=TotalSing + dim singularLocus(ideal(F));\\ 
F=F00 + 17*x0*x3*x5; TotalSing=TotalSing + dim singularLocus(ideal(F));\\ 
F0=F00 + 17*x1*x5*x5;\\
for iA from 0 to (\#A-1) list ( \\
F=F0 + 19*A\#iA; TotalSing=TotalSing + dim singularLocus(ideal(F));\\ 
);\\
F0=F00 + 17*x1*x3*x3;\\
for iB from 0 to (\#B-1) list ( \\
F=F0 + 19*B\#iB; TotalSing=TotalSing + dim singularLocus(ideal(F));\\ 
);\\
F0=F00 + 17*x1*x3*x5;\\
for iC from 0 to (\#C-1) list ( \\
F=F0 + 19*C\#iC; TotalSing=TotalSing + dim singularLocus(ideal(F));\\ 
);\\
F00 = x0\^{ }3 + 2*x1\^{ }2*x2 + 3*x2\^{ }2*x1 + 7*x3\^{ }2*x4 + 11*x4\^{ }2*x5 + 13*x5\^{ }2*x1; \\
F=F00 + 17*x2*x5*x5; TotalSing=TotalSing + dim singularLocus(ideal(F));\\ 
F=F00 + 17*x2*x2*x0; TotalSing=TotalSing + dim singularLocus(ideal(F));\\ 
F=F00 + 17*x2*x3*x5; TotalSing=TotalSing + dim singularLocus(ideal(F));\\ 
F=F00 + 17*x0*x2*x5; TotalSing=TotalSing + dim singularLocus(ideal(F));\\ 
F0=F00 + 17*x2*x4*x5;\\
for iF1 from 0 to (\#F1-1) list ( \\
F=F0 + 19*F1\#iF1; TotalSing=TotalSing + dim singularLocus(ideal(F));\\ 
);\\
F00 = x0\^{ }3 + 2*x1\^{ }2*x2 + 3*x2\^{ }2*x1 + 7*x3\^{ }2*x4 + 11*x4\^{ }2*x1 + 13*x5\^{ }2*x0; \\
F=F00 + 17*x2*x2*x5; TotalSing=TotalSing + dim singularLocus(ideal(F));\\ 
F=F00 + 17*x2*x4*x4; TotalSing=TotalSing + dim singularLocus(ideal(F));\\ 
F=F00 + 17*x2*x3*x4; TotalSing=TotalSing + dim singularLocus(ideal(F));\\ 
F=F00 + 17*x2*x4*x5; TotalSing=TotalSing + dim singularLocus(ideal(F));\\ 
F0=F00 + 17*x0*x2*x2;\\
for iAAA from 0 to (\#AAA-1) list ( \\
F=F0 + 19*AAA\#iAAA; TotalSing=TotalSing + dim singularLocus(ideal(F));\\ 
);\\
F0=F00 + 17*x0*x2*x4;\\
for iAAA from 0 to (\#AAA-1) list ( \\
F=F0 + 19*AAA\#iAAA; TotalSing=TotalSing + dim singularLocus(ideal(F));\\ 
);\\
F0=F00 + 17*x0*x4*x4;\\
for iAAA from 0 to (\#AAA-1) list ( \\
F=F0 + 19*AAA\#iAAA; TotalSing=TotalSing + dim singularLocus(ideal(F));\\ 
);\\
F00 = x0\^{ }3 + 2*x1\^{ }2*x2 + 3*x2\^{ }2*x1 + 7*x3\^{ }2*x4 + 11*x4\^{ }2*x0 + 13*x5\^{ }2*x1; \\
F=F00 + 17*x2*x3*x5; TotalSing=TotalSing + dim singularLocus(ideal(F));\\ 
F=F00 + 17*x2*x5*x5; TotalSing=TotalSing + dim singularLocus(ideal(F));\\ 
F=F00 + 17*x2*x2*x3; TotalSing=TotalSing + dim singularLocus(ideal(F));\\ 
F0=F00 + 17*x4*x2*x2;\\
for iAA from 0 to (\#AA-1) list ( \\
F=F0 + 19*AA\#iAA; TotalSing=TotalSing + dim singularLocus(ideal(F));\\ 
);\\
F0=F00 + 17*x2*x4*x5;\\
for iAA from 0 to (\#AA-1) list ( \\
F=F0 + 19*AA\#iAA; TotalSing=TotalSing + dim singularLocus(ideal(F));\\
);\\
F0=F00 + 17*x4*x5*x5;\\
for iAA from 0 to (\#AA-1) list ( \\
F=F0 + 19*AA\#iAA; TotalSing=TotalSing + dim singularLocus(ideal(F));\\ 
);\\
F0=F00 + 17*x0*x5*x5;\\
for iBB from 0 to (\#BB-1) list ( \\
F=F0 + 19*BB\#iBB; TotalSing=TotalSing + dim singularLocus(ideal(F));\\ 
);\\
F0=F00 + 17*x0*x2*x2;\\
for iBB from 0 to (\#BB-1) list ( \\
F=F0 + 19*BB\#iBB; TotalSing=TotalSing + dim singularLocus(ideal(F));\\ 
);\\
F0=F00 + 17*x0*x2*x5;\\
for iBB from 0 to (\#BB-1) list ( \\
F=F0 + 19*BB\#iBB; TotalSing=TotalSing + dim singularLocus(ideal(F));\\ 
);\\
A=[ x3\^{ }2*x0, x3\^{ }2*x1, x3\^{ }2*x5, x0*x1*x3, x0*x3*x5 ];\\
A1=[ x3\^{ }2*x1, x3\^{ }2*x5, x0*x1*x3, x0*x3*x5 ];\\
B=[ x4\^{ }2*x2, x2\^{ }2*x0, x0*x2*x4, x2*x3*x4, x2*x4*x5 ];\\
F1=[ x2\^{ }2*x3, x2\^{ }2*x5, x2\^{ }2*x4, x2*x3*x4, x2*x3*x5, x2*x4*x5 ];\\
F2=[ x2\^{ }2*x3, x2\^{ }2*x0, x2*x3*x0, x2*x3*x5 ];\\
AA=[ x3\^{ }2*x0, x3\^{ }2*x2, x3\^{ }2*x5, x0*x2*x3, x0*x2*x5, x0*x3*x5, x2*x3*x5 ];\\
BB=[ x0*x2*x4, x0*x2*x5, x0*x4*x5 ];\\
C2=[ x0*x2*x4, x2*x3*x4 ];\\
DD=[ x0*x2*x5, x2*x3*x5 ];\\
F3=[ x2*x3*x4, x2*x3*x5, x3*x4*x5 ];\\
AAA=[ x2\^{ }2*x0, x4\^{ }2*x2, x0*x2*x4, x2*x3*x4 ];\\
BBB=[ x1\^{ }2*x0, x0*x1*x3, x0*x1*x5, x1*x3*x5 ];\\
C3=[ x2\^{ }2*x0, x0*x2*x3, x0*x2*x5 ];\\
F00 = x0\^{ }3 + 2*x1\^{ }2*x2 + 3*x2\^{ }2*x1 + 7*x3\^{ }2*x4 + 11*x4\^{ }2*x1 + 13*x5\^{ }2*x2; \\
F0=F00 + 17*x1*x4*x5;\\
for iA from 0 to (\#A-1) list ( \\
for iB from 0 to (\#B-1) list ( \\
F=F0 + 19*A\#iA + 23*B\#iB; TotalSing=TotalSing + dim singularLocus(ideal(F));\\ 
));\\
F0=F00 + 17*x1*x1*x4;\\
for iA1 from 0 to (\#A1-1) list ( \\
for iB from 0 to (\#B-1) list ( \\
F=F0 + 19*A1\#iA1 + 23*B\#iB; TotalSing=TotalSing + dim singularLocus(ideal(F));\\ 
));\\
F0=F00 + 17*x5*x5*x4;\\
for iA1 from 0 to (\#A1-1) list ( \\
for iB from 0 to (\#B-1) list ( \\
F=F0 + 19*A1\#iA1 + 23*B\#iB; TotalSing=TotalSing + dim singularLocus(ideal(F));\\ 
));\\
F0=F00 + 17*x0*x1*x5;\\
for iB from 0 to (\#B-1) list ( \\
F=F0 + 19*B\#iB; TotalSing=TotalSing + dim singularLocus(ideal(F));\\ 
);\\
F0=F00 + 17*x0*x1*x1;\\
for iB from 0 to (\#B-1) list ( \\
F=F0 + 19*B\#iB; TotalSing=TotalSing + dim singularLocus(ideal(F));\\ 
);\\
F0=F00 + 17*x1*x3*x5;\\
for iB from 0 to (\#B-1) list ( \\
F=F0 + 19*B\#iB; TotalSing=TotalSing + dim singularLocus(ideal(F));\\ 
);\\
F0=F00 + 17*x1*x5*x5;\\
for iB from 0 to (\#B-1) list ( \\
F=F0 + 19*B\#iB; TotalSing=TotalSing + dim singularLocus(ideal(F));\\ 
);\\
F00 = x0\^{ }3 + 2*x1\^{ }2*x2 + 3*x2\^{ }2*x1 + 7*x3\^{ }2*x4 + 11*x4\^{ }2*x0 + 13*x5\^{ }2*x0; \\
F=F00 + 17*x4*x5*x5; TotalSing=TotalSing + dim singularLocus(ideal(F));\\ 
F=F00 + 17*x3*x4*x5; TotalSing=TotalSing + dim singularLocus(ideal(F));\\ 
F0=F00 + 17*x1*x4*x5;\\
for iF1 from 0 to (\#F1-1) list ( \\
F=F0 + 19*F1\#iF1; TotalSing=TotalSing + dim singularLocus(ideal(F));\\ 
);\\
F00 = x0\^{ }3 + 2*x1\^{ }2*x2 + 3*x2\^{ }2*x1 + 7*x3\^{ }2*x4 + 11*x4\^{ }2*x0 + 13*x5\^{ }2*x4; \\
F=F00 + 17*x0*x3*x5; TotalSing=TotalSing + dim singularLocus(ideal(F));\\ 
F0=F00 + 17*x1*x3*x5;\\
for iF2 from 0 to (\#F2-1) list ( \\
F=F0 + 19*F2\#iF2; TotalSing=TotalSing + dim singularLocus(ideal(F));\\ 
);\\
F00 = x0\^{ }3 + 2*x1\^{ }2*x2 + 3*x2\^{ }2*x1 + 7*x3\^{ }2*x4 + 11*x4\^{ }2*x1 + 13*x5\^{ }2*x1; \\
F=F00 + 17*x2*x4*x5 + 19*x2*x2*x0; TotalSing=TotalSing + dim singularLocus(ideal(F));\\
F=F00 + 17*x3*x4*x5 + 19*x2*x2*x0; TotalSing=TotalSing + dim singularLocus(ideal(F));\\ 
F=F00 + 17*x4*x5*x5 + 19*x2*x2*x0; TotalSing=TotalSing + dim singularLocus(ideal(F));\\ 
F0=F00 + 17*x0*x5*x4;\\
for iC2 from 0 to (\#C2-1) list ( \\
for iDD from 0 to (\#DD-1) list ( \\
for iF3 from 0 to (\#F3-1) list ( \\
F=F0 + 19*C2\#iC2 + 23*DD\#iDD + 29*F3\#iF3; TotalSing=TotalSing + dim singularLocus(ideal(F));\\ 
)));\\
F0=F00 + 17*x3*x5*x4;\\
for iC2 from 0 to (\#C2-1) list ( \\
for iDD from 0 to (\#DD-1) list ( \\
for iBB from 0 to (\#BB-1) list ( \\
F=F0 + 19*C2\#iC2 + 23*DD\#iDD + 29*BB\#iBB; TotalSing=TotalSing + dim singularLocus(ideal(F));\\ 
)));\\
F0=F00 + 17*x5*x5*x4;\\
for iC2 from 0 to (\#C2-1) list ( \\
for iAA from 0 to (\#AA-1) list ( \\
F=F0 + 19*C2\#iC2 + 23*AA\#iAA; TotalSing=TotalSing + dim singularLocus(ideal(F));\\ 
));\\
F0=F00 + 17*x2*x4*x5;\\
for iAA from 0 to (\#AA-1) list ( \\
F=F0 + 19*AA\#iAA; TotalSing=TotalSing + dim singularLocus(ideal(F));\\ 
); \\
F0=F00 + 17*x0*x4*x5 + 23*x0*x2*x2;\\
for iF3 from 0 to (\#F3-1) list ( \\
F=F0 + 19*F3\#iF3; TotalSing=TotalSing + dim singularLocus(ideal(F));\\ 
); \\
F00 = x0\^{ }3 + 2*x1\^{ }2*x2 + 3*x2\^{ }2*x1 + 7*x3\^{ }2*x4 + 11*x4\^{ }2*x1 + 13*x5\^{ }2*x4; \\
F0=F00 + 17*x2*x3*x5;\\
for iBBB from 0 to (\#BBB-1) list ( \\
for iAAA from 0 to (\#AAA-1) list ( \\
F=F0 + 19*BBB\#iBBB + 23*AAA\#iAAA; TotalSing=TotalSing + dim singularLocus(ideal(F));\\ 
));\\
F0=F00 + 17*x1*x3*x5;\\
for iC3 from 0 to (\#C3-1) list ( \\
for iAAA from 0 to (\#AAA-1) list ( \\
F=F0 + 19*C3\#iC3 + 23*AAA\#iAAA; TotalSing=TotalSing + dim singularLocus(ideal(F));\\ 
));\\
F0=F00 + 17*x0*x3*x5;\\
for iAAA from 0 to (\#AAA-1) list ( \\
F=F0 + 19*AAA\#iAAA; TotalSing=TotalSing + dim singularLocus(ideal(F));\\ 
); \\
A=[ x5\^{ }2*x1, x1*x3*x4, x1*x4*x5, x3*x4*x5 ];\\
B=[ x5\^{ }2*x2, x2*x3*x4, x2*x3*x5, x3*x4*x5 ];\\
F1=[ x1\^{ }2*x0, x1\^{ }2*x5, x3\^{ }2*x0, x3\^{ }2*x1, x1*x3*x4, x1*x3*x5, x0*x1*x3 ];\\
F2=[ x1\^{ }2*x5, x5\^{ }2*x1, x3\^{ }2*x1, x1*x3*x4, x1*x3*x5, x1*x4*x5, x3*x4*x5 ];\\
F3=[ x1\^{ }2*x0, x1\^{ }2*x5, x3\^{ }2*x0, x1*x3*x4, x1*x3*x5, x0*x1*x3 ];\\
F4=[ x1\^{ }2*x5, x5\^{ }2*x1, x1*x3*x4, x1*x3*x5, x1*x4*x5, x3*x4*x5 ];\\
F5=[ x1\^{ }2*x0, x3\^{ }2*x0, x1*x3*x4, x1*x3*x5, x0*x1*x3 ];\\
F6=[ x5\^{ }2*x1, x1*x3*x4, x1*x3*x5, x1*x4*x5, x3*x4*x5 ];\\
F7=[ x1\^{ }2*x0, x3\^{ }2*x0, x1*x3*x4, x0*x1*x3 ];\\
F8=[ x1\^{ }2*x0, x3\^{ }2*x0, x1*x3*x4 ];\\
F9=[ x3\^{ }2*x0, x1*x3*x4 ];\\
F10=[ x2*x3*x5, x2*x4*x5 ];\\
A1=[ x2\^{ }2*x0, x3\^{ }2*x0, x2\^{ }2*x5, x0*x2*x3, x2*x3*x5 ];\\
B1=[ x2\^{ }2*x0, x4\^{ }2*x0, x2\^{ }2*x5, x0*x2*x4, x2*x4*x5 ];\\
C1=[ x2\^{ }2*x0, x3\^{ }2*x0, x0*x2*x3 ];\\
AA=[ x2\^{ }2*x0, x3\^{ }2*x0, x2\^{ }2*x4, x0*x2*x3, x2*x3*x4 ];\\
A2=[ x1\^{ }2*x0, x5\^{ }2*x1, x0*x1*x5, x1*x3*x5, x1*x4*x5 ];\\
B2=[ x2\^{ }2*x0, x0*x2*x3, x2*x3*x5 ];\\
C2=[ x2\^{ }2*x0, x0*x2*x4, x2*x4*x5 ];\\
G1=[ x0*x3*x4, x0*x3*x5, x0*x4*x5, x3*x4*x5 ];\\
G2=[ x0*x3*x4, x0*x4*x5, x3*x4*x5 ];\\
G3=[ x2\^{ }2*x0, x0*x3*x4, x0*x2*x3 ];\\
G4=[ x2\^{ }2*x0, x2\^{ }2*x4, x3\^{ }2*x0, x0*x2*x3 ];\\
G5=[ x1\^{ }2*x4, x1\^{ }2*x5, x1*x3*x4, x1*x3*x5 ];\\
G6=[ x2\^{ }2*x4, x2*x3*x4, x2*x4*x5 ];\\
G7=[ x2*x3*x5, x0*x2*x5 ];\\
G8=[ x0*x2*x4, x0*x2*x5, x0*x4*x5, x2*x4*x5 ];\\
G9=[ x2*x3*x5, x0*x2*x3 ];\\
F00 = x0\^{ }3 + 2*x1\^{ }2*x2 + 3*x2\^{ }2*x1 + 7*x3\^{ }2*x2 + 11*x4\^{ }2*x1 + 13*x5\^{ }2*x0; \\
F0=F00 + 17*x2*x4*x4;\\
for iF1 from 0 to (\#F1-1) list ( \\
for iF2 from 0 to (\#F2-1) list ( \\
F=F0 + 19*F1\#iF1 + 23*F2\#iF2; TotalSing=TotalSing + dim singularLocus(ideal(F));\\ 
));\\
F0=F00 + 17*x2*x2*x5;\\
for iF3 from 0 to (\#F3-1) list ( \\
for iF4 from 0 to (\#F4-1) list ( \\
F=F0 + 19*F3\#iF3 + 23*F4\#iF4; TotalSing=TotalSing + dim singularLocus(ideal(F));\\ 
));\\
F0=F00 + 17*x2*x4*x5;\\
for iF5 from 0 to (\#F5-1) list ( \\
for iF6 from 0 to (\#F6-1) list ( \\
F=F0 + 19*F5\#iF5 + 23*F6\#iF6; TotalSing=TotalSing + dim singularLocus(ideal(F));\\ 
));\\
F0=F00 + 17*x2*x4*x0;\\
for iF7 from 0 to (\#F7-1) list ( \\
for iA from 0 to (\#A-1) list ( \\
for iB from 0 to (\#B-1) list ( \\
F=F0 + 19*F7\#iF7 + 23*A\#iA + 29*B\#iB; TotalSing=TotalSing + dim singularLocus(ideal(F));\\ 
)));\\
F0=F00 + 17*x2*x2*x0;\\
for iF8 from 0 to (\#F8-1) list ( \\
for iA from 0 to (\#A-1) list ( \\
for iB from 0 to (\#B-1) list ( \\
F=F0 + 19*F8\#iF8 + 23*A\#iA + 29*B\#iB; TotalSing=TotalSing + dim singularLocus(ideal(F));\\ 
)));\\
F0=F00 + 17*x4*x4*x0;\\
for iF9 from 0 to (\#F9-1) list ( \\
for iA from 0 to (\#A-1) list ( \\
for iB from 0 to (\#B-1) list ( \\
F=F0 + 19*F9\#iF9 + 23*A\#iA + 29*B\#iB; TotalSing=TotalSing + dim singularLocus(ideal(F));\\ 
)));\\
F=F00 + 17*x3*x4*x1 + 19*x3*x4*x2; TotalSing=TotalSing + dim singularLocus(ideal(F));\\ 
F00 = x0\^{ }3 + 2*x1\^{ }2*x2 + 3*x2\^{ }2*x1 + 7*x3\^{ }2*x1 + 11*x4\^{ }2*x1 + 13*x5\^{ }2*x0; \\
F0=F00 + 17*x0*x4*x4;\\
for iF10 from 0 to (\#F10-1) list ( \\
for iA1 from 0 to (\#A1-1) list ( \\
F=F0 + 19*F10\#iF10 + 23*A1\#iA1; TotalSing=TotalSing + dim singularLocus(ideal(F));\\ 
));\\
F0=F00 + 17*x3*x4*x5;\\
for iB1 from 0 to (\#B1-1) list ( \\
for iC1 from 0 to (\#C1-1) list ( \\
F=F0 + 19*B1\#iB1 + 23*C1\#iC1; TotalSing=TotalSing + dim singularLocus(ideal(F));\\ 
));\\
F0=F00 + 17*x3*x4*x0 + 31*x3*x5*x2;\\
for iB1 from 0 to (\#B1-1) list ( \\
F=F0 + 19*B1\#iB1; TotalSing=TotalSing + dim singularLocus(ideal(F));\\ 
); \\
F=F00 + 17*x3*x4*x2; TotalSing=TotalSing + dim singularLocus(ideal(F));\\ 
F=F00 + 17*x3*x4*x0 + 19*x3*x4*x5 + 21*x2*x2*x5; TotalSing=TotalSing + dim singularLocus(ideal(F));\\ 
F00 = x0\^{ }3 + 2*x1\^{ }2*x2 + 3*x2\^{ }2*x1 + 7*x3\^{ }2*x1 + 11*x4\^{ }2*x1 + 13*x5\^{ }2*x2; \\
F0=F00 + 17*x2*x3*x4;\\
for iG1 from 0 to (\#G1-1) list ( \\
for iA2 from 0 to (\#A2-1) list ( \\
F=F0 + 19*G1\#iG1 + 23*A2\#iA2; TotalSing=TotalSing + dim singularLocus(ideal(F));\\ 
));\\
F0=F00 + 17*x2*x3*x3 + 29*x2*x4*x4;\\
for iG2 from 0 to (\#G2-1) list ( \\
for iA2 from 0 to (\#A2-1) list ( \\
F=F0 + 19*G2\#iG2 + 23*A2\#iA2; TotalSing=TotalSing + dim singularLocus(ideal(F));\\ 
));\\
F0=F00 + 17*x2*x3*x3;\\
for iG1 from 0 to (\#G1-1) list ( \\
for iA2 from 0 to (\#A2-1) list ( \\
for iC2 from 0 to (\#C2-1) list ( \\
F=F0 + 19*G1\#iG1 + 23*A2\#iA2 + 29*C2\#iC2; TotalSing=TotalSing + dim singularLocus(ideal(F));\\ 
)));\\
F0=F00 + 17*x3*x4*x5;\\
for iG3 from 0 to (\#G3-1) list ( \\
for iA2 from 0 to (\#A2-1) list ( \\
for iB2 from 0 to (\#B2-1) list ( \\
for iC2 from 0 to (\#C2-1) list ( \\
F=F0 + 19*G3\#iG3 + 23*A2\#iA2 + 29*C2\#iC2 + 31*B2\#iB2; TotalSing=TotalSing + dim singularLocus(ideal(F));\\ 
))));\\
F0=F00 + 17*x2*x3*x5 + 31*x0*x3*x4;\\
for iA2 from 0 to (\#A2-1) list ( \\
for iC2 from 0 to (\#C2-1) list ( \\
F=F0 + 23*A2\#iA2 + 29*C2\#iC2; TotalSing=TotalSing + dim singularLocus(ideal(F));\\ 
));\\
F00 = x0\^{ }3 + 2*x1\^{ }2*x2 + 3*x2\^{ }2*x1 + 7*x3\^{ }2*x1 + 11*x4\^{ }2*x0 + 13*x5\^{ }2*x0; \\
F0=F00 + 17*x5*x2*x4;\\
for iG5 from 0 to (\#G5-1) list ( \\
for iAA from 0 to (\#AA-1) list ( \\
F=F0 + 19*G5\#iG5 + 23*AA\#iAA; TotalSing=TotalSing + dim singularLocus(ideal(F));\\ 
));\\
F0=F00 + 17*x5*x3*x4;\\
for iAA from 0 to (\#AA-1) list ( \\
F=F0 + 23*AA\#iAA; TotalSing=TotalSing + dim singularLocus(ideal(F));\\ 
); \\
F0=F00 + 17*x5*x2*x4 + 31*x5*x1*x4;\\
for iG4 from 0 to (\#G4-1) list ( \\
F=F0 + 19*G4\#iG4; TotalSing=TotalSing + dim singularLocus(ideal(F));\\ 
); \\
F0=F00 + 17*x5*x1*x4 + 23*x3*x3*x2;\\
for iG6 from 0 to (\#G6-1) list ( \\
F=F0 + 19*G6\#iG6; TotalSing=TotalSing + dim singularLocus(ideal(F));\\ 
); \\
F=F00 + 17*x5*x1*x4 + 23*x4*x3*x2; TotalSing=TotalSing + dim singularLocus(ideal(F));\\ 
F=F00 + 17*x5*x3*x4 + 23*x3*x3*x2; TotalSing=TotalSing + dim singularLocus(ideal(F));\\ 
F00 = x0\^{ }3 + 2*x1\^{ }2*x2 + 3*x2\^{ }2*x1 + 7*x3\^{ }2*x1 + 11*x4\^{ }2*x1 + 13*x5\^{ }2*x1; \\
F0=F00 + 17*x3*x2*x4 + 29*x3*x5*x4;\\
for iG7 from 0 to (\#G7-1) list ( \\
for iG8 from 0 to (\#G8-1) list ( \\
F=F0 + 19*G7\#iG7 + 23*G8\#iG8; TotalSing=TotalSing + dim singularLocus(ideal(F));\\ 
));\\
F0=F00 + 17*x3*x0*x4 + 29*x3*x5*x0 + 31*x2*x4*x5;\\
for iG9 from 0 to (\#G9-1) list ( \\
F=F0 + 19*G9\#iG9; TotalSing=TotalSing + dim singularLocus(ideal(F));\\ 
); \\
F=F00 + 17*x3*x2*x4 + 29*x3*x5*x0 + 31*x2*x4*x5; TotalSing=TotalSing + dim singularLocus(ideal(F));\\ 
F=F00 + 17*x2*x2*x0 + 31*x3*x4*x5; TotalSing=TotalSing + dim singularLocus(ideal(F));\\ 
F=F00 + 17*x5*x2*x0 + 31*x3*x4*x5 + 19*x2*x3*x0 + 23*x2*x4*x0; TotalSing=TotalSing + dim singularLocus(ideal(F));\\ 
F=F00 + 17*x5*x3*x0 + 31*x3*x2*x5 + 19*x4*x5*x0 + 23*x3*x4*x0 + 29*x2*x4*x0;\\
TotalSing=TotalSing + dim singularLocus(ideal(F));\\ 
F00 = x0\^{ }3 + 2*x1\^{ }2*x2 + 3*x2\^{ }2*x1 + 7*x3\^{ }2*x0 + 11*x4\^{ }2*x0 + 13*x5\^{ }2*x0; \\
F=F00 + 17*x5*x3*x1 + 31*x4*x1*x5 + 19*x2*x3*x4; TotalSing=TotalSing + dim singularLocus(ideal(F));\\ 
F=F00 + 17*x5*x3*x4; TotalSing=TotalSing + dim singularLocus(ideal(F));\\ 
print "This is Case of One Cube, No cycles";\\
A=[ x1*x2*x4, x1*x2*x5, x1*x4*x5 ];\\
F1=[ x3\^{ }2*x0, x0*x1*x5, x0*x3*x5 ];\\
F2=[ x1\^{ }2*x0, x3\^{ }2*x0, x5\^{ }2*x0, x0*x1*x3, x0*x1*x5, x0*x3*x5 ];\\
F3=[ x1*x3*x4, x1*x4*x5, x3*x4*x5 ];\\
F4=[ x1\^{ }2*x0, x1\^{ }2*x3, x1\^{ }2*x5, x0*x1*x3, x0*x1*x5, x0*x3*x5 ];\\
F5=[ x2\^{ }2*x4, x2\^{ }2*x5, x1*x2*x4, x1*x2*x5 ];\\
F6=[ x1\^{ }2*x3, x1\^{ }2*x4, x1\^{ }2*x5, x1*x3*x4, x1*x3*x5, x3*x4*x5 ];\\
F7=[ x1\^{ }2*x3, x1\^{ }2*x4, x1\^{ }2*x5, x1*x3*x4, x1*x4*x5 ];\\
F8=[ x2\^{ }2*x4, x1*x2*x4, x1*x4*x5, x2*x4*x5 ];\\
F9=[ x1*x3*x4, x2*x3*x4 ];\\
F10=[ x1*x3*x5, x2*x3*x5 ];\\
F11=[ x1*x3*x4, x1*x4*x5 ];\\
F12=[ x1*x3*x5, x2*x3*x5 ];\\
G1=[ x1\^{ }2*x0, x1\^{ }2*x4, x1*x0*x4, x1*x4*x5, x0*x4*x5 ];\\
G2=[ x2\^{ }2*x0, x1*x2*x4, x0*x2*x4 ];\\
G3=[ x2\^{ }2*x0, x1*x2*x5, x0*x2*x5 ];\\
G4=[ x1*x2*x4, x1*x4*x5 ];\\
G5=[ x0*x2*x5, x1*x2*x5 ];\\
G6=[ x0*x1*x5, x1*x3*x5 ];\\
G7=[ x0*x1*x4, x1*x3*x4 ];\\
G8=[ x0*x1*x5, x1*x3*x5 ];\\
G9=[ x0*x1*x4, x0*x4*x5 ];\\
F = x0\^{ }3 + x1\^{ }2*x2 + x2\^{ }2*x3 + x3\^{ }2*x4 + x4\^{ }2*x5 + x5\^{ }2*x0; TotalSing=TotalSing + dim singularLocus(ideal(F));\\ 
F00 = x0\^{ }3 + 2*x1\^{ }2*x2 + 3*x2\^{ }2*x3 + 7*x3\^{ }2*x4 + 11*x4\^{ }2*x0 + 13*x5\^{ }2*x2; \\
F=F00 + 19*x5*x3*x1; TotalSing=TotalSing + dim singularLocus(ideal(F));\\ 
F=F00 + 19*x0*x1*x5 + 17*x1*x3*x4; TotalSing=TotalSing + dim singularLocus(ideal(F));\\ 
F=F00 + 19*x5*x5*x3; TotalSing=TotalSing + dim singularLocus(ideal(F));\\ 
F0=F00 + 23*x1*x4*x5;\\
for iF1 from 0 to (\#F1-1) list ( \\
F=F0 + 19*F1\#iF1; TotalSing=TotalSing + dim singularLocus(ideal(F));\\ 
);\\
F0=F00 + 23*x4*x5*x5;\\
for iF2 from 0 to (\#F2-1) list ( \\
F=F0 + 19*F2\#iF2; TotalSing=TotalSing + dim singularLocus(ideal(F));\\ 
);\\
F0=F00 + 23*x0*x5*x5;\\
for iF3 from 0 to (\#F3-1) list ( \\
F=F0 + 19*F3\#iF3; TotalSing=TotalSing + dim singularLocus(ideal(F));\\ 
);\\
F00 = x0\^{ }3 + 2*x1\^{ }2*x2 + 3*x2\^{ }2*x3 + 7*x3\^{ }2*x4 + 11*x4\^{ }2*x0 + 13*x5\^{ }2*x3; \\
F=F00 + 19*x2*x2*x4; TotalSing=TotalSing + dim singularLocus(ideal(F));\\ 
F=F00 + 19*x2*x4*x5; TotalSing=TotalSing + dim singularLocus(ideal(F));\\ 
F=F00 + 19*x1*x2*x5; TotalSing=TotalSing + dim singularLocus(ideal(F));\\ 
F=F00 + 19*x4*x5*x5; TotalSing=TotalSing + dim singularLocus(ideal(F));\\ 
F0=F00 + 23*x0*x2*x2;\\
for iA from 0 to (\#A-1) list ( \\
F=F0 + 19*A\#iA; TotalSing=TotalSing + dim singularLocus(ideal(F));\\ 
);\\
F0=F00 + 23*x0*x2*x5;\\
for iA from 0 to (\#A-1) list ( \\
F=F0 + 19*A\#iA; TotalSing=TotalSing + dim singularLocus(ideal(F));\\ 
);\\
F0=F00 + 23*x0*x5*x5;\\
for iA from 0 to (\#A-1) list ( \\
F=F0 + 19*A\#iA; TotalSing=TotalSing + dim singularLocus(ideal(F));\\ 
);\\
F00 = x0\^{ }3 + 2*x1\^{ }2*x2 + 3*x2\^{ }2*x3 + 7*x3\^{ }2*x4 + 11*x4\^{ }2*x0 + 13*x5\^{ }2*x4; \\
F=F00 + 19*x0*x5*x5; TotalSing=TotalSing + dim singularLocus(ideal(F));\\ 
F=F00 + 19*x0*x3*x5; TotalSing=TotalSing + dim singularLocus(ideal(F));\\ 
F=F00 + 19*x1*x3*x5; TotalSing=TotalSing + dim singularLocus(ideal(F));\\ 
F00 = x0\^{ }3 + 2*x1\^{ }2*x2 + 3*x2\^{ }2*x3 + 7*x3\^{ }2*x4 + 11*x4\^{ }2*x0 + 13*x5\^{ }2*x4; \\
F=F00 + 19*x0*x3*x3; TotalSing=TotalSing + dim singularLocus(ideal(F));\\ 
F0=F00 + 23*x2*x3*x5;\\
for iF4 from 0 to (\#F4-1) list ( \\
F=F0 + 19*F4\#iF4; TotalSing=TotalSing + dim singularLocus(ideal(F));\\ 
);\\
F00 = x0\^{ }3 + 2*x1\^{ }2*x2 + 3*x2\^{ }2*x3 + 7*x3\^{ }2*x4 + 11*x4\^{ }2*x0 + 13*x5\^{ }2*x0; \\
F=F00 + 19*x1*x4*x5; TotalSing=TotalSing + dim singularLocus(ideal(F));\\ 
F0=F00 + 23*x3*x4*x5;\\
for iF5 from 0 to (\#F5-1) list ( \\
F=F0 + 19*F5\#iF5; TotalSing=TotalSing + dim singularLocus(ideal(F));\\ 
);\\
F0=F00 + 23*x2*x4*x5;\\
for iF6 from 0 to (\#F6-1) list ( \\
F=F0 + 19*F6\#iF6; TotalSing=TotalSing + dim singularLocus(ideal(F));\\ 
);\\
F00 = x0\^{ }3 + 2*x1\^{ }2*x2 + 3*x2\^{ }2*x3 + 7*x3\^{ }2*x0 + 11*x4\^{ }2*x5 + 13*x5\^{ }2*x3; \\
F=F00 + 19*x0*x2*x5; TotalSing=TotalSing + dim singularLocus(ideal(F));\\ 
F=F00 + 19*x0*x5*x5; TotalSing=TotalSing + dim singularLocus(ideal(F));\\ 
F=F00 + 19*x2*x4*x5; TotalSing=TotalSing + dim singularLocus(ideal(F));\\ 
F00 = x0\^{ }3 + 2*x1\^{ }2*x2 + 3*x2\^{ }2*x3 + 7*x3\^{ }2*x0 + 11*x4\^{ }2*x5 + 13*x5\^{ }2*x0; \\
F=F00 + 19*x3*x4*x5; TotalSing=TotalSing + dim singularLocus(ideal(F));\\ 
F=F00 + 19*x1*x3*x5; TotalSing=TotalSing + dim singularLocus(ideal(F));\\ 
F0=F00 + 23*x2*x3*x5;\\
for iF7 from 0 to (\#F7-1) list ( \\
F=F0 + 19*F7\#iF7; TotalSing=TotalSing + dim singularLocus(ideal(F));\\ 
);\\
F00 = x0\^{ }3 + 2*x1\^{ }2*x2 + 3*x2\^{ }2*x3 + 7*x3\^{ }2*x0 + 11*x4\^{ }2*x0 + 13*x5\^{ }2*x0; \\
F0=F00 + 23*x3*x4*x5;\\
for iF8 from 0 to (\#F8-1) list ( \\
F=F0 + 19*F8\#iF8; TotalSing=TotalSing + dim singularLocus(ideal(F));\\ 
);\\
F0=F00 + 23*x2*x4*x5;\\
for iF9 from 0 to (\#F9-1) list ( \\
for iF10 from 0 to (\#F10-1) list ( \\
for iF11 from 0 to (\#F11-1) list ( \\
F=F0 + 19*F9\#iF9 + 29*F10\#iF10 + 31*F11\#iF11; TotalSing=TotalSing + dim singularLocus(ideal(F));\\ 
))); \\
F0=F00 + 23*x1*x4*x5 + 41*x2*x3*x4;\\
for iF12 from 0 to (\#F12-1) list ( \\
F=F0 + 19*F12\#iF12; TotalSing=TotalSing + dim singularLocus(ideal(F));\\ 
);\\
F00 = x0\^{ }3 + 2*x1\^{ }2*x2 + 3*x2\^{ }2*x3 + 7*x3\^{ }2*x0 + 11*x4\^{ }2*x3 + 13*x5\^{ }2*x3; \\
F0=F00 + 23*x2*x4*x5;\\
for iG1 from 0 to (\#G1-1) list ( \\
F=F0 + 19*G1\#iG1; TotalSing=TotalSing + dim singularLocus(ideal(F));\\ 
);\\
F0=F00 + 23*x0*x4*x5;\\
for iG2 from 0 to (\#G2-1) list ( \\
for iG3 from 0 to (\#G3-1) list ( \\
for iG4 from 0 to (\#G4-1) list ( \\
F=F0 + 19*G2\#iG2 + 19*G3\#iG3 + 19*G4\#iG4; TotalSing=TotalSing + dim singularLocus(ideal(F));\\ 
)));\\
F0=F00 + 23*x1*x4*x5 + 41*x0*x2*x4;\\
for iG5 from 0 to (\#G5-1) list ( \\
F=F0 + 19*G5\#iG5; TotalSing=TotalSing + dim singularLocus(ideal(F));\\ 
);\\
F=F00 + 23*x1*x4*x5 + 41*x0*x2*x2; TotalSing=TotalSing + dim singularLocus(ideal(F));\\ 
F00 = x0\^{ }3 + 2*x1\^{ }2*x2 + 3*x2\^{ }2*x3 + 7*x3\^{ }2*x0 + 11*x4\^{ }2*x2 + 13*x5\^{ }2*x2;\\
F=F00 + 23*x1*x4*x5; TotalSing=TotalSing + dim singularLocus(ideal(F));\\ 
F0=F00 + 23*x0*x4*x5 + 41*x1*x3*x4;\\
for iG6 from 0 to (\#G6-1) list ( \\
F=F0 + 19*G6\#iG6; TotalSing=TotalSing + dim singularLocus(ideal(F));\\ 
);\\
F0=F00 + 23*x3*x4*x5;\\
for iG7 from 0 to (\#G7-1) list ( \\
for iG8 from 0 to (\#G8-1) list ( \\
for iG9 from 0 to (\#G9-1) list ( \\
F=F0 + 19*G7\#iG7 + 19*G8\#iG8 + 19*G9\#iG9; TotalSing=TotalSing + dim singularLocus(ideal(F));\\ 
)));\\
F1=[ x4\^{ }2*x0, x4\^{ }2*x3, x0*x1*x4, x1*x3*x4, x1*x4*x5 ];\\
F2=[ x5\^{ }2*x0, x5\^{ }2*x2, x2\^{ }2*x0, x0*x2*x5, x1*x2*x5, x2*x4*x5 ];\\
F3=[ x1\^{ }2*x3, x4\^{ }2*x3, x1*x3*x4, x1*x3*x5, x1*x4*x5, x3*x4*x5 ];\\
F4=[ x1\^{ }2*x0, x4\^{ }2*x0, x5\^{ }2*x0, x0*x1*x4, x0*x1*x5, x0*x4*x5, x1*x4*x5 ];\\
F5=[ x1\^{ }2*x0, x1\^{ }2*x3, x0*x1*x4, x1*x3*x4, x1*x4*x5 ];\\
F6=[ x5\^{ }2*x2, x1*x3*x5, x2*x3*x5, x3*x4*x5 ];\\
F7=[ x1\^{ }2*x3, x1*x3*x4, x1*x3*x5, x1*x4*x5, x3*x4*x5 ];\\
F8=[ x4\^{ }2*x3, x1*x3*x4, x1*x3*x5, x1*x4*x5, x3*x4*x5 ];\\
F9=[ x5\^{ }2*x3, x1*x3*x5, x2*x3*x5, x3*x4*x5 ];\\
F10=[ x2\^{ }2*x0, x4\^{ }2*x0, x0*x2*x4, x1*x2*x4, x2*x4*x5 ];\\
F11=[ x1*x2*x4, x1*x2*x5, x1*x4*x5, x2*x4*x5 ];\\
F12=[ x1\^{ }2*x3, x1\^{ }2*x4, x1\^{ }2*x5, x5\^{ }2*x3, x1*x3*x4, x1*x3*x5, x1*x4*x5, x3*x4*x5 ];\\
G1=[ x5\^{ }2*x2, x5\^{ }2*x3, x3\^{ }2*x2, x1*x2*x3, x1*x2*x5, x1*x3*x5, x2*x3*x5 ];\\
G2=[ x1\^{ }2*x3, x1\^{ }2*x4, x3\^{ }2*x4, x1*x3*x4, x1*x3*x5, x1*x4*x5, x3*x4*x5 ];\\
G3=[ x5\^{ }2*x2, x3\^{ }2*x2, x1*x2*x3, x1*x2*x5, x1*x3*x5, x2*x3*x5 ];\\
G4=[ x1*x2*x4, x2*x4*x5 ];\\
G5=[ x1*x3*x4, x3*x4*x5 ];\\
G6=[ x1\^{ }2*x4, x0*x1*x3, x1*x3*x4, x1*x3*x5 ];\\
G7=[ x1\^{ }2*x4, x3\^{ }2*x4, x1*x3*x5, x1*x3*x4, x1*x4*x5, x3*x4*x5 ];\\
G8=[ x1*x3*x5, x1*x3*x0, x1*x0*x5, x3*x0*x5 ];\\
G9=[ x1*x2*x4, x2*x3*x4, x2*x4*x5 ];\\
H1=[ x1*x3*x4, x1*x4*x5 ];\\
H2=[ x1\^{ }2*x0, x3\^{ }2*x0, x0*x1*x3, x1*x3*x4 ];\\
H3=[ x5\^{ }2*x2, x1*x2*x5 ];\\
H4=[ x2*x3*x4, x2*x3*x5, x3*x4*x5 ];\\
H5=[ x1*x4*x5, x3*x4*x5 ];\\
H6=[ x1*x2*x4, x2*x3*x4 ];\\
A1=[ x1*x3*x4, x1*x3*x5, x1*x4*x5, x3*x4*x5 ];\\
A2=[ x1\^{ }2*x0, x3\^{ }2*x0, x0*x1*x3, x1*x3*x4, x1*x3*x5 ];\\
H7=[ x1*x3*x0, x1*x3*x5 ];\\
H8=[ x1*x4*x0, x1*x4*x5 ];\\
H9=[ x3*x4*x5, x1*x3*x5 ];\\
B=[ x5\^{ }2*x2, x1*x2*x5, x2*x3*x5, x2*x4*x5 ];\\
F00 = x0\^{ }3 + 2*x1\^{ }2*x2 + 3*x2\^{ }2*x3 + 7*x3\^{ }2*x0; \\
F0=F00 + 23*x2*x4*x4 + 11*x3*x5*x5;\\
for iF1 from 0 to (\#F1-1) list ( \\
for iF2 from 0 to (\#F2-1) list ( \\
for iF3 from 0 to (\#F3-1) list ( \\
for iF4 from 0 to (\#F4-1) list ( \\
F=F0 + 19*F1\#iF1 + 29*F2\#iF2 + 31*F3\#iF3 + 37*F4\#iF4; TotalSing=TotalSing + dim singularLocus(ideal(F));\\ 
)))); \\
F0=F00 + 23*x2*x4*x4 + 11*x0*x5*x5;\\
for iF5 from 0 to (\#F5-1) list ( \\
for iF6 from 0 to (\#F6-1) list ( \\
for iF7 from 0 to (\#F7-1) list ( \\
for iF8 from 0 to (\#F8-1) list ( \\
F=F0 + 19*F5\#iF5 + 29*F6\#iF6 + 31*F7\#iF7 + 37*F8\#iF8; TotalSing=TotalSing + dim singularLocus(ideal(F));\\ 
))));\\
F0=F00 + 23*x3*x4*x4 + 11*x0*x5*x5;\\
for iF9 from 0 to (\#F9-1) list ( \\
for iF10 from 0 to (\#F10-1) list ( \\
for iF11 from 0 to (\#F11-1) list ( \\
for iF12 from 0 to (\#F12-1) list ( \\
F=F0 + 19*F9\#iF9 + 29*F10\#iF10 + 31*F11\#iF11 + 37*F12\#iF12; TotalSing=TotalSing + dim singularLocus(ideal(F));\\ 
))));\\
F00 = x0\^{ }3 + 2*x1\^{ }2*x2 + 3*x2\^{ }2*x0 + 7*x3\^{ }2*x0 + 11*x4\^{ }2*x0 + 13*x5\^{ }2*x4; \\
F0=F00 + 23*x2*x3*x4 ;\\
for iG1 from 0 to (\#G1-1) list ( \\
for iG2 from 0 to (\#G2-1) list ( \\
F=F0 + 19*G1\#iG1 + 29*G2\#iG2; TotalSing=TotalSing + dim singularLocus(ideal(F));\\ 
));\\
F0=F00 + 23*x3*x3*x4 ;\\
for iG3 from 0 to (\#G3-1) list ( \\
for iG4 from 0 to (\#G4-1) list ( \\
F=F0 + 19*G3\#iG3 + 29*G4\#iG4; TotalSing=TotalSing + dim singularLocus(ideal(F));\\ 
));\\
F0=F00 + 23*x1*x2*x3  + 27*x2*x4*x5;\\
for iG5 from 0 to (\#G5-1) list ( \\
F=F0 + 19*G5\#iG5; TotalSing=TotalSing + dim singularLocus(ideal(F));\\ 
);\\
F0=F00 + 23*x1*x2*x4  + 27*x2*x3*x5;\\
for iG5 from 0 to (\#G5-1) list ( \\
F=F0 + 19*G5\#iG5; TotalSing=TotalSing + dim singularLocus(ideal(F));\\ 
);\\
F=F00 + 23*x1*x3*x4  + 27*x2*x4*x5  + 29*x2*x3*x5; TotalSing=TotalSing + dim singularLocus(ideal(F));\\ 
F0 = x0\^{ }3 + 2*x1\^{ }2*x2 + 3*x2\^{ }2*x0 + 7*x3\^{ }2*x2 + 11*x4\^{ }2*x0 + 13*x5\^{ }2*x4; \\
for iG6 from 0 to (\#G6-1) list ( \\
for iG7 from 0 to (\#G7-1) list ( \\
for iG8 from 0 to (\#G8-1) list ( \\
for iG9 from 0 to (\#G9-1) list ( \\
F=F0 + 19*G6\#iG6 + 29*G7\#iG7 + 31*G8\#iG8 + 37*G9\#iG9; TotalSing=TotalSing + dim singularLocus(ideal(F));\\ 
))));\\
F00 = x0\^{ }3 + 2*x1\^{ }2*x2 + 3*x2\^{ }2*x0 + 7*x3\^{ }2*x2 + 11*x4\^{ }2*x2 + 13*x5\^{ }2*x0; \\
F0=F00 + 17*x0*x3*x4;\\
for iH7 from 0 to (\#H7-1) list ( \\
for iH8 from 0 to (\#H8-1) list ( \\
for iH9 from 0 to (\#H9-1) list ( \\
for iB from 0 to (\#B-1) list ( \\
F=F0 + 19*B\#iB + 29*H7\#iH7 + 31*H8\#iH8 + 37*H9\#iH9; TotalSing=TotalSing + dim singularLocus(ideal(F));\\ 
))));\\
F0=F00 + 17*x1*x3*x4;\\
for iB from 0 to (\#B-1) list ( \\
F=F0 + 19*B\#iB; TotalSing=TotalSing + dim singularLocus(ideal(F));\\ 
);\\
F0=F00 + 17*x0*x4*x4;\\
for iH3 from 0 to (\#H3-1) list ( \\
for iA1 from 0 to (\#A1-1) list ( \\
for iA2 from 0 to (\#A2-1) list ( \\
F=F0 + 19*A1\#iA1 + 29*A2\#iA2 + 31*H3\#iH3; TotalSing=TotalSing + dim singularLocus(ideal(F));\\ 
))); \\
F00 = x0\^{ }3 + 2*x1\^{ }2*x2 + 3*x2\^{ }2*x0 + 7*x3\^{ }2*x2 + 11*x4\^{ }2*x0 + 13*x5\^{ }2*x0; \\
F0=F00 + 17*x2*x4*x5;\\
for iH1 from 0 to (\#H1-1) list ( \\
for iH2 from 0 to (\#H2-1) list ( \\
F=F0 + 19*H1\#iH1 + 31*H2\#iH2; TotalSing=TotalSing + dim singularLocus(ideal(F));\\ 
));\\
F0=F00 + 17*x1*x2*x5;\\
for iH4 from 0 to (\#H4-1) list ( \\
for iH5 from 0 to (\#H5-1) list ( \\
for iH6 from 0 to (\#H6-1) list ( \\
for iA1 from 0 to (\#A1-1) list ( \\
for iA2 from 0 to (\#A2-1) list ( \\
F=F0 + 19*H4\#iH4 + 31*H5\#iH5 + 23*H6\#iH6 + 29*A1\#iA1 + 37*A2\#iA2;\\ TotalSing=TotalSing + dim singularLocus(ideal(F));\\ 
)))));\\
C1=[ x1*x2*x5, x2*x3*x5, x2*x4*x5 ];\\
C2=[ x1*x2*x4, x2*x3*x4, x2*x4*x5 ];\\
F1=[ x2*x3*x4, x2*x3*x5 ];\\
F2=[ x2*x3*x4, x1*x3*x4 ];\\
F3=[ x1*x4*x5, x2*x4*x5 ];\\
F4=[ x1*x2*x3, x1*x2*x5, x1*x3*x5, x2*x3*x5 ];\\
F5=[ x1*x2*x3, x1*x2*x4, x1*x3*x4, x2*x3*x4 ];\\
F6=[ x1*x2*x5, x1*x2*x4, x1*x4*x5, x2*x4*x5 ];\\
F7=[ x1*x2*x3, x2*x3*x4, x2*x3*x5 ];\\
F8=[ x0*x1*x5, x1*x3*x5 ];\\
F9=[ x0*x3*x5, x1*x3*x5, x3*x4*x5 ];\\
F10=[ x0*x4*x5, x1*x4*x5, x3*x4*x5 ];\\
F11=[ x0*x1*x4, x0*x1*x5, x0*x4*x5, x1*x4*x5 ];\\
F12=[ x0*x3*x4, x0*x3*x5, x0*x4*x5, x3*x4*x5 ];\\
F13=[ x1*x2*x5, x2*x3*x5, x2*x4*x5 ];\\
F14=[ x1*x2*x4, x2*x3*x4, x2*x4*x5 ];\\
F00 = x0\^{ }3 + 2*x1\^{ }2*x2 + 3*x2\^{ }2*x0 + 7*x3\^{ }2*x0 + 11*x4\^{ }2*x0 + 13*x5\^{ }2*x0; \\
F0=F00 + 17*x1*x3*x5;\\
for iF1 from 0 to (\#F1-1) list ( \\
for iF2 from 0 to (\#F2-1) list ( \\
for iF3 from 0 to (\#F3-1) list ( \\
for iC1 from 0 to (\#C1-1) list ( \\
for iC2 from 0 to (\#C2-1) list ( \\
F=F0 + 19*F1\#iF1 + 29*F2\#iF2 + 31*F3\#iF3 + 37*C1\#iC1 + 41*C2\#iC2;\\ TotalSing=TotalSing + dim singularLocus(ideal(F));\\ 
))))); \\
F0=F00 + 17*x3*x4*x5;\\
for iF4 from 0 to (\#F4-1) list ( \\
for iF5 from 0 to (\#F5-1) list ( \\
for iF6 from 0 to (\#F6-1) list ( \\
for iF7 from 0 to (\#F7-1) list ( \\
for iC1 from 0 to (\#C1-1) list ( \\
for iC2 from 0 to (\#C2-1) list ( \\
F=F0 + 19*F4\#iF4 + 29*F5\#iF5 + 31*F6\#iF6 + 37*C1\#iC1 + 41*C2\#iC2 + 47*F7\#iF7;\\ TotalSing=TotalSing + dim singularLocus(ideal(F));\\ 
)))))); \\
F00 = x0\^{ }3 + 2*x1\^{ }2*x2 + 3*x2\^{ }2*x0 + 7*x3\^{ }2*x2 + 11*x4\^{ }2*x2 + 13*x5\^{ }2*x2; \\
F0=F00 + 17*x1*x3*x4;\\
for iF8 from 0 to (\#F8-1) list ( \\
for iF9 from 0 to (\#F9-1) list ( \\
for iF10 from 0 to (\#F10-1) list ( \\
for iF11 from 0 to (\#F11-1) list ( \\
for iF12 from 0 to (\#F12-1) list ( \\
F=F0 + 19*F8\#iF8 + 29*F9\#iF9 + 31*F10\#iF10 + 37*F11\#iF11 + 41*F12\#iF12;\\ TotalSing=TotalSing + dim singularLocus(ideal(F));\\ 
))))); \\
F00 = x0\^{ }3 + 2*x1\^{ }2*x0 + 3*x2\^{ }2*x0 + 7*x3\^{ }2*x0 + 11*x4\^{ }2*x0 + 13*x5\^{ }2*x0; \\
F0=F00 + 17*x1*x2*x3 + 43*x1*x4*x5 + 23*x3*x4*x5;\\
for iF13 from 0 to (\#F13-1) list ( \\
for iF14 from 0 to (\#F14-1) list ( \\
F=F0 + 31*F13\#iF13 + 37*F14\#iF14; TotalSing=TotalSing + dim singularLocus(ideal(F));\\ 
));\\
print "This is Case of No Cubes";\\
F = x0\^{ }2*x1 + x1\^{ }2*x2 + x2\^{ }2*x3 + x3\^{ }2*x4 + x4\^{ }2*x5 + x5\^{ }2*x0; TotalSing=TotalSing + dim singularLocus(ideal(F));\\ 
print "This is Case of No Cubes, Length 5 longest cycle";\\
F1=[ x2\^{ }2*x0, x2\^{ }2*x4, x2\^{ }2*x5, x0*x2*x4, x0*x4*x5, x2*x4*x5 ];\\
F2=[ x3\^{ }2*x0, x3\^{ }2*x2, x3\^{ }2*x5, x0*x2*x3, x0*x2*x5, x2*x3*x5 ];\\
F3=[ x3\^{ }2*x0, x3\^{ }2*x2, x3\^{ }2*x5, x5\^{ }2*x3, x0*x2*x3, x0*x2*x5, x0*x3*x5, x2*x3*x5 ];\\
F4=[ x2\^{ }2*x0, x2\^{ }2*x4, x2\^{ }2*x5, x0\^{ }2*x4, x5\^{ }2*x4, x0*x2*x4, x0*x4*x5, x2*x4*x5 ];\\
F00 = x0\^{ }2*x1 + 2*x1\^{ }2*x2 + 3*x2\^{ }2*x3 + 7*x3\^{ }2*x4 + 11*x4\^{ }2*x5 + 13*x5\^{ }2*x1; \\
F=F00 + 19*x2*x5*x5; TotalSing=TotalSing + dim singularLocus(ideal(F));\\ 
F=F00 + 19*x0*x0*x5; TotalSing=TotalSing + dim singularLocus(ideal(F));\\ 
F=F00 + 19*x0*x2*x5; TotalSing=TotalSing + dim singularLocus(ideal(F));\\ 
F0=F00 + 19*x0*x3*x5;\\
for iF1 from 0 to (\#F1-1) list ( \\
F=F0 + 23*F1\#iF1; TotalSing=TotalSing + dim singularLocus(ideal(F));\\ 
); \\
F0=F00 + 19*x5*x5*x3;\\
for iF1 from 0 to (\#F1-1) list ( \\
F=F0 + 23*F1\#iF1; TotalSing=TotalSing + dim singularLocus(ideal(F));\\ 
);\\
F0=F00 + 19*x0*x4*x5;\\
for iF2 from 0 to (\#F2-1) list ( \\
F=F0 + 23*F2\#iF2; TotalSing=TotalSing + dim singularLocus(ideal(F));\\ 
);\\
F0=F00 + 19*x5*x5*x4;\\
for iF3 from 0 to (\#F3-1) list ( \\
F=F0 + 23*F3\#iF3; TotalSing=TotalSing + dim singularLocus(ideal(F));\\ 
);\\
F0=F00 + 19*x0*x0*x4;\\
for iF3 from 0 to (\#F3-1) list ( \\
F=F0 + 23*F3\#iF3; TotalSing=TotalSing + dim singularLocus(ideal(F));\\ 
);\\
F0=F00 + 19*x0*x0*x3;\\
for iF4 from 0 to (\#F4-1) list ( \\
F=F0 + 23*F4\#iF4; TotalSing=TotalSing + dim singularLocus(ideal(F));\\ 
);\\
print "This is Case of No Cubes, Length 4 longest cycle";\\
A=[ x1\^{ }2*x3, x1\^{ }2*x4, x1\^{ }2*x5, x1*x3*x5, x1*x4*x5 ];\\
AA=[ x1\^{ }2*x3, x1\^{ }2*x4, x1\^{ }2*x0, x5\^{ }2*x3, x5\^{ }2*x0, x5\^{ }2*x1, x1*x0*x5, x1*x3*x5, x1*x4*x5 ];\\
F1=[ x1\^{ }2*x4, x1*x4*x5, x3*x4*x5 ];\\
B=[ x3\^{ }2*x5, x1\^{ }2*x4, x5\^{ }2*x3, x5\^{ }2*x1, x1*x3*x4, x1*x3*x5, x1*x4*x5, x3*x4*x5 ];\\
C=[ x1\^{ }2*x3, x1\^{ }2*x4, x4\^{ }2*x3, x4\^{ }2*x1, x5\^{ }2*x3, x5\^{ }2*x1, x1*x3*x4, x1*x3*x5, x1*x4*x5, x3*x4*x5 ];\\
A1=[ x1\^{ }2*x3, x1\^{ }2*x4, x1\^{ }2*x0, x5\^{ }2*x3, x5\^{ }2*x0, x5\^{ }2*x1, x1*x0*x5, x1*x3*x5, x1*x4*x5 ];\\
A23=[ x1\^{ }2*x3, x1\^{ }2*x4, x5\^{ }2*x3, x5\^{ }2*x1, x1*x3*x4, x1*x3*x5, x1*x4*x5, x3*x4*x5 ];\\
A2=[ x3\^{ }2*x5 ];\\
A3=[ x4\^{ }2*x3, x4\^{ }2*x1 ];\\
A45=[ x3\^{ }2*x5, x4\^{ }2*x3, x4\^{ }2*x1, x1*x3*x4, x1*x3*x5, x1*x4*x5, x3*x4*x5 ];\\
A4=[ x1\^{ }2*x3, x1\^{ }2*x4 ];\\
A5=[ x5\^{ }2*x3, x5\^{ }2*x1 ];\\
F = x0\^{ }2*x1 + x1\^{ }2*x2 + x2\^{ }2*x3 + x3\^{ }2*x0 + x4\^{ }2*x5 + x5\^{ }2*x4; TotalSing=TotalSing + dim singularLocus(ideal(F));\\ 
F00 = x0\^{ }2*x1 + 2*x1\^{ }2*x2 + 3*x2\^{ }2*x3 + 7*x3\^{ }2*x0 + 11*x4\^{ }2*x0 + 13*x5\^{ }2*x4; \\
F=F00 + 19*x3*x3*x1; TotalSing=TotalSing + dim singularLocus(ideal(F));\\ 
F=F00 + 19*x4*x4*x1; TotalSing=TotalSing + dim singularLocus(ideal(F));\\ 
F=F00 + 19*x4*x3*x1; TotalSing=TotalSing + dim singularLocus(ideal(F));\\ 
F=F00 + 19*x4*x3*x5; TotalSing=TotalSing + dim singularLocus(ideal(F));\\ 
F0=F00 + 23*x2*x4*x4;\\
for iA from 0 to (\#A-1) list ( \\
F=F0 + 19*A\#iA; TotalSing=TotalSing + dim singularLocus(ideal(F));\\ 
);\\
F0=F00 + 23*x2*x3*x4;\\
for iA from 0 to (\#A-1) list ( \\
F=F0 + 19*A\#iA; TotalSing=TotalSing + dim singularLocus(ideal(F));\\ 
);\\
F0=F00 + 23*x2*x3*x3;\\
for iA from 0 to (\#A-1) list ( \\
F=F0 + 19*A\#iA; TotalSing=TotalSing + dim singularLocus(ideal(F));\\ 
);\\
F00 = x0\^{ }2*x1 + 2*x1\^{ }2*x2 + 3*x2\^{ }2*x3 + 7*x3\^{ }2*x0 + 11*x4\^{ }2*x0 + 13*x5\^{ }2*x2; \\
F0 = F00 + 17*x1*x3*x3;\\
for iAA from 0 to (\#AA-1) list ( \\
for iC from 0 to (\#C-1) list ( \\
F=F0 + 19*AA\#iAA + 23*C\#iC; TotalSing=TotalSing + dim singularLocus(ideal(F));\\ 
)); \\
F0 = F00 + 17*x1*x4*x4;\\
for iAA from 0 to (\#AA-1) list ( \\
for iB from 0 to (\#B-1) list ( \\
F=F0 + 19*AA\#iAA + 23*B\#iB; TotalSing=TotalSing + dim singularLocus(ideal(F));\\ 
)); \\
F0 = F00 + 17*x3*x4*x4;\\
for iAA from 0 to (\#AA-1) list ( \\
for iB from 0 to (\#B-1) list ( \\
F=F0 + 19*AA\#iAA + 23*B\#iB; TotalSing=TotalSing + dim singularLocus(ideal(F));\\ 
)); \\
F0 = F00 + 17*x3*x3*x5;\\
for iAA from 0 to (\#AA-1) list ( \\
for iF1 from 0 to (\#F1-1) list ( \\
F=F0 + 19*AA\#iAA + 23*F1\#iF1; TotalSing=TotalSing + dim singularLocus(ideal(F));\\ 
)); \\
F0 = F00 + 17*x3*x4*x5;\\
for iAA from 0 to (\#AA-1) list ( \\
F=F0 + 19*AA\#iAA; TotalSing=TotalSing + dim singularLocus(ideal(F));\\ 
);\\
F0 = F00 + 17*x3*x4*x1;\\
for iAA from 0 to (\#AA-1) list ( \\
F=F0 + 19*AA\#iAA; TotalSing=TotalSing + dim singularLocus(ideal(F));\\ 
);\\
F00 = x0\^{ }2*x1 + 2*x1\^{ }2*x2 + 3*x2\^{ }2*x3 + 7*x3\^{ }2*x0 + 11*x4\^{ }2*x0 + 13*x5\^{ }2*x2; \\
F0 = F00 + 17*x2*x4*x4;\\
for iA1 from 0 to (\#A1-1) list ( \\
for iA2 from 0 to (\#A2-1) list ( \\
for iA3 from 0 to (\#A3-1) list ( \\
for iA4 from 0 to (\#A4-1) list ( \\
for iA5 from 0 to (\#A5-1) list ( \\
F=F0 + 19*A1\#iA1 + 23*A2\#iA2 + 29*A3\#iA3 + 31*A4\#iA4 + 37*A5\#iA5;\\ TotalSing=TotalSing + dim singularLocus(ideal(F));\\ 
))))); \\
for iA1 from 0 to (\#A1-1) list ( \\
for iA23 from 0 to (\#A23-1) list ( \\
for iA4 from 0 to (\#A4-1) list ( \\
for iA5 from 0 to (\#A5-1) list ( \\
F=F0 + 19*A1\#iA1 + 23*A23\#iA23 + 31*A4\#iA4 + 37*A5\#iA5; TotalSing=TotalSing + dim singularLocus(ideal(F));\\ 
))));\\
for iA1 from 0 to (\#A1-1) list ( \\
for iA2 from 0 to (\#A2-1) list ( \\
for iA3 from 0 to (\#A3-1) list ( \\
for iA45 from 0 to (\#A45-1) list ( \\
F=F0 + 19*A1\#iA1 + 23*A2\#iA2 + 29*A3\#iA3 + 31*A45\#iA45; TotalSing=TotalSing + dim singularLocus(ideal(F));\\ 
))));\\
for iA1 from 0 to (\#A1-1) list ( \\
for iA23 from 0 to (\#A23-1) list ( \\
for iA45 from 0 to (\#A45-1) list ( \\
F=F0 + 19*A1\#iA1 + 23*A23\#iA23 + 31*A45\#iA45; TotalSing=TotalSing + dim singularLocus(ideal(F));\\ 
)));\\
F0 = F00 + 17*x2*x3*x3;\\
for iA1 from 0 to (\#A1-1) list ( \\
for iA2 from 0 to (\#A2-1) list ( \\
for iA3 from 0 to (\#A3-1) list ( \\
for iA4 from 0 to (\#A4-1) list ( \\
for iA5 from 0 to (\#A5-1) list ( \\
F=F0 + 19*A1\#iA1 + 23*A2\#iA2 + 29*A3\#iA3 + 31*A4\#iA4 + 37*A5\#iA5;\\ TotalSing=TotalSing + dim singularLocus(ideal(F));\\ 
)))));\\ 
for iA1 from 0 to (\#A1-1) list ( \\
for iA23 from 0 to (\#A23-1) list ( \\
for iA4 from 0 to (\#A4-1) list ( \\
for iA5 from 0 to (\#A5-1) list ( \\
F=F0 + 19*A1\#iA1 + 23*A23\#iA23 + 31*A4\#iA4 + 37*A5\#iA5; TotalSing=TotalSing + dim singularLocus(ideal(F));\\ 
))));\\
for iA1 from 0 to (\#A1-1) list ( \\
for iA2 from 0 to (\#A2-1) list ( \\
for iA3 from 0 to (\#A3-1) list ( \\
for iA45 from 0 to (\#A45-1) list ( \\
F=F0 + 19*A1\#iA1 + 23*A2\#iA2 + 29*A3\#iA3 + 31*A45\#iA45; TotalSing=TotalSing + dim singularLocus(ideal(F));\\ 
))));\\
for iA1 from 0 to (\#A1-1) list ( \\
for iA23 from 0 to (\#A23-1) list ( \\
for iA45 from 0 to (\#A45-1) list ( \\
F=F0 + 19*A1\#iA1 + 23*A23\#iA23 + 31*A45\#iA45; TotalSing=TotalSing + dim singularLocus(ideal(F));\\ 
)));\\
F0 = F00 + 17*x2*x3*x4;\\
for iA1 from 0 to (\#A1-1) list ( \\
for iA2 from 0 to (\#A2-1) list ( \\
for iA3 from 0 to (\#A3-1) list ( \\
for iA4 from 0 to (\#A4-1) list ( \\
for iA5 from 0 to (\#A5-1) list ( \\
F=F0 + 19*A1\#iA1 + 23*A2\#iA2 + 29*A3\#iA3 + 31*A4\#iA4 + 37*A5\#iA5;\\ TotalSing=TotalSing + dim singularLocus(ideal(F));\\ 
)))));\\ 
for iA1 from 0 to (\#A1-1) list ( \\
for iA23 from 0 to (\#A23-1) list ( \\
for iA4 from 0 to (\#A4-1) list ( \\
for iA5 from 0 to (\#A5-1) list ( \\
F=F0 + 19*A1\#iA1 + 23*A23\#iA23 + 31*A4\#iA4 + 37*A5\#iA5; TotalSing=TotalSing + dim singularLocus(ideal(F));\\ 
))));\\
for iA1 from 0 to (\#A1-1) list ( \\
for iA2 from 0 to (\#A2-1) list ( \\
for iA3 from 0 to (\#A3-1) list ( \\
for iA45 from 0 to (\#A45-1) list ( \\
F=F0 + 19*A1\#iA1 + 23*A2\#iA2 + 29*A3\#iA3 + 31*A45\#iA45; TotalSing=TotalSing + dim singularLocus(ideal(F));\\ 
))));\\
for iA1 from 0 to (\#A1-1) list ( \\
for iA23 from 0 to (\#A23-1) list ( \\
for iA45 from 0 to (\#A45-1) list ( \\
F=F0 + 19*A1\#iA1 + 23*A23\#iA23 + 31*A45\#iA45; TotalSing=TotalSing + dim singularLocus(ideal(F));\\ 
)));\\
A1=[ x5\^{ }2*x0, x5\^{ }2*x3, x0\^{ }2*x2, x0\^{ }2*x3, x0\^{ }2*x4, x0*x2*x5, x0*x3*x5, x0*x4*x5 ];\\
A2=[ x5\^{ }2*x0, x2\^{ }2*x0, x0\^{ }2*x2, x2\^{ }2*x5, x2\^{ }2*x4, x0\^{ }2*x4, x0*x2*x4, x0*x2*x5, x0*x4*x5, x2*x4*x5 ];\\
B=[ x3\^{ }2*x5, x3\^{ }2*x2, x5\^{ }2*x3, x4\^{ }2*x2, x2*x3*x4, x2*x3*x5, x2*x4*x5, x3*x4*x5 ];\\
C=[ x3\^{ }2*x2, x3\^{ }2*x5, x2*x3*x4, x2*x3*x5, x2*x4*x5, x3*x4*x5 ];\\
D=[ x4\^{ }2*x1, x3\^{ }2*x1, x3\^{ }2*x5, x1\^{ }2*x5, x1\^{ }2*x4, x1\^{ }2*x3, x1*x3*x4, x1*x3*x5, x1*x4*x5, x3*x4*x5 ];\\
E=[ x3\^{ }2*x5, x1\^{ }2*x5, x1\^{ }2*x4, x1\^{ }2*x3, x1*x3*x4, x1*x3*x5, x1*x4*x5, x3*x4*x5 ];\\
F00 = x0\^{ }2*x1 + 2*x1\^{ }2*x2 + 3*x2\^{ }2*x3 + 7*x3\^{ }2*x0 + 11*x4\^{ }2*x0 + 13*x5\^{ }2*x1; \\
F0=F00 + 19*x3*x3*x1;\\
for iA1 from 0 to (\#A1-1) list ( \\
for iA2 from 0 to (\#A2-1) list ( \\
for iB from 0 to (\#B-1) list ( \\
F=F0 + 19*A1\#iA1 + 23*A2\#iA2 + 31*B\#iB; TotalSing=TotalSing + dim singularLocus(ideal(F));\\ 
)));\\ 
F0=F00 + 19*x4*x3*x1;\\
for iA1 from 0 to (\#A1-1) list ( \\
for iA2 from 0 to (\#A2-1) list ( \\
for iB from 0 to (\#B-1) list ( \\
F=F0 + 19*A1\#iA1 + 23*A2\#iA2 + 31*B\#iB; TotalSing=TotalSing + dim singularLocus(ideal(F));\\ 
)));\\ 
F0=F00 + 19*x4*x4*x1;\\
for iA1 from 0 to (\#A1-1) list ( \\
for iA2 from 0 to (\#A2-1) list ( \\
for iC from 0 to (\#C-1) list ( \\
F=F0 + 19*A1\#iA1 + 23*A2\#iA2 + 31*C\#iC; TotalSing=TotalSing + dim singularLocus(ideal(F));\\ 
)));\\ 
F0=F00 + 19*x4*x4*x2;\\
for iA1 from 0 to (\#A1-1) list ( \\
for iA2 from 0 to (\#A2-1) list ( \\
for iD from 0 to (\#D-1) list ( \\
F=F0 + 19*A1\#iA1 + 23*A2\#iA2 + 31*D\#iD; TotalSing=TotalSing + dim singularLocus(ideal(F));\\ 
)));\\ 
F0=F00 + 19*x3*x3*x2;\\
for iA1 from 0 to (\#A1-1) list ( \\
for iA2 from 0 to (\#A2-1) list ( \\
for iE from 0 to (\#E-1) list ( \\
F=F0 + 19*A1\#iA1 + 23*A2\#iA2 + 31*E\#iE; TotalSing=TotalSing + dim singularLocus(ideal(F));\\ 
)));\\ 
F0=F00 + 19*x4*x3*x2;\\
for iA1 from 0 to (\#A1-1) list ( \\
for iA2 from 0 to (\#A2-1) list ( \\
for iE from 0 to (\#E-1) list ( \\
F=F0 + 19*A1\#iA1 + 23*A2\#iA2 + 31*E\#iE; TotalSing=TotalSing + dim singularLocus(ideal(F));\\ 
)));\\ 
F0=F00 + 19*x4*x3*x5;\\
for iA1 from 0 to (\#A1-1) list ( \\
for iA2 from 0 to (\#A2-1) list ( \\
F=F0 + 19*A1\#iA1 + 23*A2\#iA2; TotalSing=TotalSing + dim singularLocus(ideal(F));\\ 
));\\
A1=[ x3\^{ }2*x2, x3\^{ }2*x1, x1*x3*x4, x2*x3*x4 ];\\
A2=[ x3\^{ }2*x2, x3\^{ }2*x1, x1*x3*x5, x2*x3*x5 ];\\
B1=[ x3\^{ }2*x1, x1*x3*x4, x1*x3*x5, x1*x4*x5 ];
B2=[ x1\^{ }2*x3, x1\^{ }2*x4, x1\^{ }2*x5, x1*x3*x4, x1*x3*x5, x1*x4*x5 ];\\
F1=[ x2\^{ }2*x1, x2\^{ }2*x4, x2\^{ }2*x5, x1*x2*x4, x1*x2*x5, x1*x4*x5, x2*x4*x5 ];\\
F2=[ x3\^{ }2*x2, x2*x3*x4, x2*x3*x5 ];\\
F00 = x0\^{ }2*x1 + 2*x1\^{ }2*x2 + 3*x2\^{ }2*x3 + 7*x3\^{ }2*x0 + 11*x4\^{ }2*x0 + 13*x5\^{ }2*x0; \\
F0=F00 + 19*x3*x4*x5;\\
for iF1 from 0 to (\#F1-1) list ( \\
F=F0 + 19*F1\#iF1; TotalSing=TotalSing + dim singularLocus(ideal(F));\\ 
); \\
F0=F00 + 19*x2*x4*x5;\\
for iA1 from 0 to (\#A1-1) list ( \\
for iA2 from 0 to (\#A2-1) list ( \\
for iB1 from 0 to (\#B1-1) list ( \\
for iB2 from 0 to (\#B2-1) list ( \\
F=F0 + 19*A1\#iA1 + 23*A2\#iA2 + 31*B1\#iB1 + 37*B2\#iB2; TotalSing=TotalSing + dim singularLocus(ideal(F));\\ 
))));\\ 
F0=F00 + 19*x1*x4*x5;\\
for iA1 from 0 to (\#A1-1) list ( \\
for iA2 from 0 to (\#A2-1) list ( \\
for iF2 from 0 to (\#F2-1) list ( \\
F=F0 + 19*A1\#iA1 + 23*A2\#iA2 + 31*F2\#iF2; TotalSing=TotalSing + dim singularLocus(ideal(F));\\ 
)));\\
print "This is Case of No Cubes, Length 3 longest cycle";\\
A=[ x5\^{ }2*x2, x2\^{ }2*x5, x2\^{ }2*x1, x3\^{ }2*x1, x5\^{ }2*x1, x1*x2*x3, x1*x2*x5, x1*x3*x5, x2*x3*x5 ];\\
F1=[ x1\^{ }2*x0, x1\^{ }2*x5, x1\^{ }2*x3, x1*x0*x3, x1*x0*x5, x1*x3*x5, x0*x3*x5 ];\\
F2=[ x5\^{ }2*x2, x2\^{ }2*x1, x1*x2*x3, x1*x2*x5, x1*x3*x5, x2*x3*x5 ];\\
F3=[ x3\^{ }2*x1, x3\^{ }2*x2, x3\^{ }2*x5, x1*x2*x3, x1*x3*x5 ];\\
F4=[ x4\^{ }2*x1, x4\^{ }2*x2, x2\^{ }2*x4, x1*x2*x4, x1*x4*x5, x2*x4*x5 ];\\
F5=[ x4\^{ }2*x1, x4\^{ }2*x2, x4\^{ }2*x5, x1*x2*x4, x1*x4*x5, x2*x4*x5 ];\\
F = x0\^{ }2*x1 + x1\^{ }2*x2 + x2\^{ }2*x0 + x3\^{ }2*x4 + x4\^{ }2*x5 + x5\^{ }2*x3;\\ TotalSing=TotalSing + dim singularLocus(ideal(F));\\ 
F00 = x0\^{ }2*x1 + 2*x1\^{ }2*x2 + 3*x2\^{ }2*x0 + 7*x3\^{ }2*x4 + 11*x4\^{ }2*x5 + 13*x5\^{ }2*x4; \\
F=F00 + 19*x3*x3*x5; TotalSing=TotalSing + dim singularLocus(ideal(F));\\ 
F0=F00 + 23*x0*x3*x3;\\
for iA from 0 to (\#A-1) list ( \\
F=F0 + 19*A\#iA; TotalSing=TotalSing + dim singularLocus(ideal(F));\\ 
);\\
F0=F00 + 23*x2*x3*x5;\\
for iF1 from 0 to (\#F1-1) list ( \\
F=F0 + 19*F1\#iF1; TotalSing=TotalSing + dim singularLocus(ideal(F));\\ 
);\\
F0=F00 + 23*x0*x5*x5;\\
for iF2 from 0 to (\#F2-1) list ( \\
F=F0 + 19*F2\#iF2; TotalSing=TotalSing + dim singularLocus(ideal(F));\\ 
);\\
F00 = x0\^{ }2*x1 + 2*x1\^{ }2*x2 + 3*x2\^{ }2*x0 + 7*x3\^{ }2*x4 + 11*x4\^{ }2*x5 + 13*x5\^{ }2*x0; \\
F=F00 + 23*x1*x2*x2; TotalSing=TotalSing + dim singularLocus(ideal(F));\\ 
F=F00 + 23*x1*x5*x5; TotalSing=TotalSing + dim singularLocus(ideal(F));\\ 
F=F00 + 23*x1*x2*x5; TotalSing=TotalSing + dim singularLocus(ideal(F));\\ 
F=F00 + 23*x2*x3*x5; TotalSing=TotalSing + dim singularLocus(ideal(F));\\ 
F0=F00 + 23*x2*x4*x5;\\
for iF3 from 0 to (\#F3-1) list ( \\
F=F0 + 19*F3\#iF3; TotalSing=TotalSing + dim singularLocus(ideal(F));\\ 
);\\
F00 = x0\^{ }2*x1 + 2*x1\^{ }2*x2 + 3*x2\^{ }2*x0 + 7*x3\^{ }2*x4 + 11*x4\^{ }2*x3 + 13*x5\^{ }2*x0; \\
F=F00 + 23*x1*x2*x5; TotalSing=TotalSing + dim singularLocus(ideal(F));\\ 
F=F00 + 23*x1*x5*x5; TotalSing=TotalSing + dim singularLocus(ideal(F));\\ 
F=F00 + 23*x1*x2*x2; TotalSing=TotalSing + dim singularLocus(ideal(F));\\ 
F0=F00 + 23*x2*x2*x3;\\
for iF4 from 0 to (\#F4-1) list ( \\
F=F0 + 19*F4\#iF4; TotalSing=TotalSing + dim singularLocus(ideal(F));\\ 
);\\
F0=F00 + 23*x2*x3*x5;\\
for iF5 from 0 to (\#F5-1) list ( \\
F=F0 + 19*F5\#iF5; TotalSing=TotalSing + dim singularLocus(ideal(F));\\ 
);\\
A1=[ x1\^{ }2*x3, x5\^{ }2*x3, x1*x3*x4, x1*x3*x5, x1*x4*x5, x3*x4*x5 ];\\
A2=[ x4\^{ }2*x0, x4\^{ }2*x1, x4\^{ }2*x5, x1\^{ }2*x0, x5\^{ }2*x0, x1*x0*x4, x1*x0*x5, x0*x4*x5, x1*x4*x5 ];\\
B1=[ x5\^{ }2*x0, x5\^{ }2*x2, x5\^{ }2*x3, x0\^{ }2*x2, x0\^{ }2*x3, x0*x2*x5, x0*x3*x5, x0*x4*x5 ];\\
B2=[ x4\^{ }2*x0, x4\^{ }2*x2, x4\^{ }2*x5, x5\^{ }2*x0, x5\^{ }2*x2, x0\^{ }2*x2, x0*x2*x4, x0*x2*x5, x0*x4*x5, x2*x4*x5 ];\\
C=[ x5\^{ }2*x3, x3\^{ }2*x2, x2*x3*x4, x2*x3*x5, x2*x4*x5, x3*x4*x5 ];\\
D=[ x5\^{ }2*x3, x5\^{ }2*x2, x2*x3*x5, x2*x4*x5, x3*x4*x5 ];\\
E1=[ x5\^{ }2*x0, x5\^{ }2*x3, x5\^{ }2*x2, x0\^{ }2*x3, x0\^{ }2*x2, x0*x2*x5, x0*x3*x5, x0*x4*x5 ];\\
E2=[ x5\^{ }2*x0, x5\^{ }2*x2, x4\^{ }2*x0, x4\^{ }2*x2, x0\^{ }2*x2, x4\^{ }2*x5, x0*x2*x5, x0*x4*x5, x0*x2*x4, x2*x4*x5 ];\\
FF=[ x5\^{ }2*x3, x5\^{ }2*x0, x1\^{ }2*x3, x1\^{ }2*x0, x1*x3*x5, x0*x1*x5, x1*x4*x5 ];\\
F1=[ x5\^{ }2*x0, x5\^{ }2*x3, x0\^{ }2*x2, x0*x2*x5, x0*x3*x5, x0*x4*x5 ];\\
F2=[ x5\^{ }2*x0, x0\^{ }2*x2, x4\^{ }2*x0, x4\^{ }2*x2, x4\^{ }2*x5, x0*x2*x4, x0*x2*x5, x0*x4*x5, x2*x4*x5 ];\\
G1=[ x5\^{ }2*x3, x2\^{ }2*x1, x1*x2*x5, x2*x3*x5, x2*x4*x5 ];\\
G2=[ x5\^{ }2*x1, x2\^{ }2*x1, x4\^{ }2*x1, x4\^{ }2*x2, x4\^{ }2*x5, x1*x2*x4, x1*x2*x5, x1*x4*x5, x2*x4*x5 ];\\
H1=[ x1\^{ }2*x0, x1\^{ }2*x3, x5\^{ }2*x0, x5\^{ }2*x1, x5\^{ }2*x3, x0*x1*x5, x1*x3*x5, x1*x4*x5 ];\\
H2=[ x1\^{ }2*x0, x5\^{ }2*x0, x5\^{ }2*x1, x4\^{ }2*x0, x4\^{ }2*x1, x4\^{ }2*x5, x0*x1*x4, x0*x1*x5, x0*x4*x5, x1*x4*x5 ];\\
H3=[ x1\^{ }2*x0, x1\^{ }2*x3, x5\^{ }2*x0, x5\^{ }2*x3, x0*x1*x5, x1*x3*x5, x1*x4*x5 ];\\
F00 = x0\^{ }2*x1 + 2*x1\^{ }2*x2 + 3*x2\^{ }2*x0 + 7*x3\^{ }2*x0 + 11*x4\^{ }2*x3 + 13*x5\^{ }2*x1; \\
F0=F00 + 17*x2*x3*x4;\\
for iB1 from 0 to (\#B1-1) list ( \\
for iB2 from 0 to (\#B2-1) list ( \\
F=F0 + 19*B1\#iB1 + 23*B2\#iB2; TotalSing=TotalSing + dim singularLocus(ideal(F));\\ 
));\\
F0=F00 + 17*x1*x3*x3;\\
for iF1 from 0 to (\#F1-1) list ( \\
for iF2 from 0 to (\#F2-1) list ( \\
for iC from 0 to (\#C-1) list ( \\
F=F0 + 19*F1\#iF1 + 23*F2\#iF2 + 31*C\#iC; TotalSing=TotalSing + dim singularLocus(ideal(F));\\ 
)));\\
F0=F00 + 17*x1*x2*x2;\\
for iE1 from 0 to (\#E1-1) list ( \\
for iE2 from 0 to (\#E2-1) list ( \\
for iD from 0 to (\#D-1) list ( \\
F=F0 + 19*E1\#iE1 + 23*E2\#iE2 + 31*D\#iD; TotalSing=TotalSing + dim singularLocus(ideal(F));\\ 
)));\\
F0=F00 + 17*x1*x2*x3;\\
for iE1 from 0 to (\#E1-1) list ( \\
for iE2 from 0 to (\#E2-1) list ( \\
for iD from 0 to (\#D-1) list ( \\
F=F0 + 19*E1\#iE1 + 23*E2\#iE2 + 31*D\#iD; TotalSing=TotalSing + dim singularLocus(ideal(F));\\ 
)));\\
F0=F00 + 17*x5*x2*x3;\\
for iE1 from 0 to (\#E1-1) list ( \\
for iE2 from 0 to (\#E2-1) list ( \\
F=F0 + 19*E1\#iE1 + 23*E2\#iE2; TotalSing=TotalSing + dim singularLocus(ideal(F));\\ 
));\\
F0 = x0\^{ }2*x1 + 2*x1\^{ }2*x2 + 3*x2\^{ }2*x0 + 7*x3\^{ }2*x0 + 11*x4\^{ }2*x3 + 13*x5\^{ }2*x0 + 17*x3\^{ }2*x1; \\
for iG1 from 0 to (\#G1-1) list ( \\
for iG2 from 0 to (\#G2-1) list ( \\
for iC from 0 to (\#C-1) list ( \\
F=F0 + 19*G1\#iG1 + 23*G2\#iG2 + 31*C\#iC; TotalSing=TotalSing + dim singularLocus(ideal(F));\\ 
)));\\
F00 = x0\^{ }2*x1 + 2*x1\^{ }2*x2 + 3*x2\^{ }2*x0 + 7*x3\^{ }2*x0 + 11*x4\^{ }2*x3 + 13*x5\^{ }2*x2; \\
F0=F00 + 17*x1*x3*x3;\\
for iH1 from 0 to (\#H1-1) list ( \\
for iH2 from 0 to (\#H2-1) list ( \\
F=F0 + 19*H1\#iH1 + 23*H2\#iH2; TotalSing=TotalSing + dim singularLocus(ideal(F));\\ 
));\\
F0=F00 + 17*x2*x3*x4;\\
for iA1 from 0 to (\#A1-1) list ( \\
for iA2 from 0 to (\#A2-1) list ( \\
for iH3 from 0 to (\#H3-1) list ( \\
F=F0 + 19*A1\#iA1 + 23*A2\#iA2 + 31*H3\#iH3; TotalSing=TotalSing + dim singularLocus(ideal(F));\\ 
)));\\
F0=F00 + 17*x2*x2*x1;\\
for iA1 from 0 to (\#A1-1) list ( \\
for iA2 from 0 to (\#A2-1) list ( \\
for iFF from 0 to (\#FF-1) list ( \\
F=F0 + 19*A1\#iA1 + 23*A2\#iA2 + 31*FF\#iFF; TotalSing=TotalSing + dim singularLocus(ideal(F));\\ 
)));\\
F0=F00 + 17*x1*x2*x3;\\
for iA1 from 0 to (\#A1-1) list ( \\
for iA2 from 0 to (\#A2-1) list ( \\
for iFF from 0 to (\#FF-1) list ( \\
F=F0 + 19*A1\#iA1 + 23*A2\#iA2 + 31*FF\#iFF; TotalSing=TotalSing + dim singularLocus(ideal(F));\\ 
)));\\
F0=F00 + 17*x5*x2*x3;\\
for iA1 from 0 to (\#A1-1) list ( \\
for iA2 from 0 to (\#A2-1) list ( \\
for iFF from 0 to (\#FF-1) list ( \\
F=F0 + 19*A1\#iA1 + 23*A2\#iA2 + 31*FF\#iFF; TotalSing=TotalSing + dim singularLocus(ideal(F));\\ 
)));\\
F0=F00 + 17*x5*x2*x2;\\
for iA1 from 0 to (\#A1-1) list ( \\
for iA2 from 0 to (\#A2-1) list ( \\
for iFF from 0 to (\#FF-1) list ( \\
F=F0 + 19*A1\#iA1 + 23*A2\#iA2 + 31*FF\#iFF; TotalSing=TotalSing + dim singularLocus(ideal(F));\\ 
)));\\
A1=[ x1*x3*x5, x2*x3*x5, x3*x4*x5 ];\\
A2=[ x2\^{ }2*x1, x1*x2*x5, x2*x3*x5, x2*x4*x5 ];\\
B1=[ x5\^{ }2*x0, x4\^{ }2*x0, x0*x4*x5, x1*x4*x5, x2*x4*x5 ];\\
B2=[ x5\^{ }2*x0, x4\^{ }2*x0, x0\^{ }2*x2, x0*x2*x4, x0*x2*x5, x0*x4*x5, x2*x4*x5 ];\\
B3=[ x5\^{ }2*x0, x4\^{ }2*x0, x1\^{ }2*x0, x0*x1*x4, x0*x1*x5, x0*x4*x5, x1*x4*x5 ];\\
C=[ x1*x2*x4, x1*x2*x5, x2*x4*x5, x1*x4*x5 ];\\
D=[ x2*x3*x5, x2*x4*x5, x3*x4*x5 ];\\
E=[ x1\^{ }2*x5, x1\^{ }2*x3, x1*x3*x4, x1*x3*x5, x1*x4*x5, x3*x4*x5 ];\\
FF=[ x4\^{ }2*x1, x4\^{ }2*x2, x4\^{ }2*x5, x1*x2*x4, x1*x2*x5, x1*x4*x5, x2*x4*x5 ];\\
H1=[ x0\^{ }2*x2, x0*x2*x4, x0*x2*x5, x0*x4*x5, x2*x4*x5 ];\\
H2=[ x1\^{ }2*x0, x0*x1*x4, x0*x1*x5, x0*x4*x5, x1*x4*x5 ];\\
H3=[ x2\^{ }2*x1, x1*x2*x4, x1*x2*x5, x2*x4*x5, x1*x4*x5 ];\\
H4=[ x1*x4*x5, x2*x4*x5, x0*x4*x5 ];\\
H5=[ x2\^{ }2*x1, x1*x2*x3, x1*x2*x5, x2*x3*x5, x1*x3*x5 ];\\
H6=[ x1\^{ }2*x3, x1\^{ }2*x5, x1*x3*x4, x1*x3*x5, x1*x4*x5, x3*x4*x5 ];\\
H7=[ x2\^{ }2*x1, x4\^{ }2*x1, x4\^{ }2*x2, x4\^{ }2*x5, x1*x2*x4, x1*x2*x5, x1*x4*x5, x2*x4*x5 ];\\
F00 = x0\^{ }2*x1 + 2*x1\^{ }2*x2 + 3*x2\^{ }2*x0 + 7*x3\^{ }2*x0 + 11*x4\^{ }2*x3 + 13*x5\^{ }2*x3; \\
F0=F00 + 17*x1*x3*x3;\\
for iH1 from 0 to (\#H1-1) list ( \\
for iH2 from 0 to (\#H2-1) list ( \\
for iH3 from 0 to (\#H3-1) list ( \\
for iH4 from 0 to (\#H4-1) list ( \\
F=F0 + 19*H1\#iH1 + 23*H2\#iH2 + 29*H3\#iH3 + 31*H4\#iH4; TotalSing=TotalSing + dim singularLocus(ideal(F));\\ 
))));\\
F0=F00 + 17*x2*x3*x4;\\
for iB1 from 0 to (\#B1-1) list ( \\
for iB2 from 0 to (\#B2-1) list ( \\
for iB3 from 0 to (\#B3-1) list ( \\
for iC from 0 to (\#C-1) list ( \\
F=F0 + 19*B1\#iB1 + 23*B2\#iB2 + 29*B3\#iB3 + 31*C\#iC; TotalSing=TotalSing + dim singularLocus(ideal(F));\\ 
))));\\
F0=F00 + 17*x2*x2*x1;\\
for iB1 from 0 to (\#B1-1) list ( \\
for iB2 from 0 to (\#B2-1) list ( \\
for iB3 from 0 to (\#B3-1) list ( \\
F=F0 + 19*B1\#iB1 + 23*B2\#iB2 + 29*B3\#iB3; TotalSing=TotalSing + dim singularLocus(ideal(F));\\ 
)));\\
F0=F00 + 17*x1*x2*x3;\\
for iB1 from 0 to (\#B1-1) list ( \\
for iB2 from 0 to (\#B2-1) list ( \\
for iB3 from 0 to (\#B3-1) list ( \\
for iC from 0 to (\#C-1) list ( \\
F=F0 + 19*B1\#iB1 + 23*B2\#iB2 + 29*B3\#iB3 + 31*C\#iC; TotalSing=TotalSing + dim singularLocus(ideal(F));\\ 
))));\\
F00 = x0\^{ }2*x1 + 2*x1\^{ }2*x2 + 3*x2\^{ }2*x0 + 7*x3\^{ }2*x0 + 11*x4\^{ }2*x3 + 13*x5\^{ }2*x0; \\
F0=F00 + 17*x2*x3*x4;\\
for iA1 from 0 to (\#A1-1) list ( \\
for iA2 from 0 to (\#A2-1) list ( \\
for iH5 from 0 to (\#H5-1) list ( \\
for iH6 from 0 to (\#H6-1) list ( \\
for iH7 from 0 to (\#H7-1) list ( \\
F=F0 + 19*A1\#iA1 + 23*A2\#iA2 + 29*H5\#iH5 + 31*H6\#iH6 + 37*H7\#iH7;\\ TotalSing=TotalSing + dim singularLocus(ideal(F));\\ 
)))));\\
F0=F00 + 17*x2*x2*x1;\\
for iA1 from 0 to (\#A1-1) list ( \\
for iA2 from 0 to (\#A2-1) list ( \\
for iD from 0 to (\#D-1) list ( \\
for iE from 0 to (\#E-1) list ( \\
F=F0 + 19*A1\#iA1 + 23*A2\#iA2 + 29*D\#iD + 31*E\#iE; TotalSing=TotalSing + dim singularLocus(ideal(F));\\ 
))));\\
F0=F00 + 17*x1*x2*x3;\\
for iA1 from 0 to (\#A1-1) list ( \\
for iA2 from 0 to (\#A2-1) list ( \\
for iD from 0 to (\#D-1) list ( \\
for iE from 0 to (\#E-1) list ( \\
for iFF from 0 to (\#FF-1) list ( \\
F=F0 + 19*A1\#iA1 + 23*A2\#iA2 + 29*D\#iD + 31*E\#iE + 37*FF\#iFF;\\ TotalSing=TotalSing + dim singularLocus(ideal(F));\\ 
)))));\\
F0=F00 + 17*x5*x2*x3;\\
for iA1 from 0 to (\#A1-1) list ( \\
for iA2 from 0 to (\#A2-1) list ( \\
for iE from 0 to (\#E-1) list ( \\
for iFF from 0 to (\#FF-1) list ( \\
F=F0 + 19*A1\#iA1 + 23*A2\#iA2 + 31*E\#iE + 37*FF\#iFF; TotalSing=TotalSing + dim singularLocus(ideal(F));\\ 
))));\\
AA13=[ x3\^{ }2*x2, x2\^{ }2*x5, x2*x3*x4, x2*x3*x5 ];\\
AA24=[ x5\^{ }2*x1, x1\^{ }2*x4, x1*x3*x5, x1*x4*x5 ];\\
A1=[ x3\^{ }2*x1, x2\^{ }2*x1, x1*x2*x3 ];\\
A2=[ x5\^{ }2*x0, x1\^{ }2*x0, x0*x1*x5 ];\\
A3=[ x4\^{ }2*x2, x2*x4*x5, x3*x4*x5 ];\\
A4=[ x3\^{ }2*x1, x1*x3*x4, x3*x4*x5 ];\\
B=[ x0\^{ }2*x3, x5\^{ }2*x0, x4\^{ }2*x0, x0*x3*x4, x0*x3*x5, x0*x4*x5, x3*x4*x5 ];\\
F00 = x0\^{ }2*x1 + 2*x1\^{ }2*x2 + 3*x2\^{ }2*x0 + 7*x3\^{ }2*x0 + 11*x4\^{ }2*x1 + 13*x5\^{ }2*x2; \\
F0=F00 + 17*x0*x2*x4; \\
for iA1 from 0 to (\#A1-1) list ( \\
for iA2 from 0 to (\#A2-1) list ( \\
for iA3 from 0 to (\#A3-1) list ( \\
for iA4 from 0 to (\#A4-1) list ( \\
for iB from 0 to (\#B-1) list ( \\
F=F0 + 19*A1\#iA1 + 23*A2\#iA2 + 29*A3\#iA3 + 31*A4\#iA4 + 37*B\#iB;\\ TotalSing=TotalSing + dim singularLocus(ideal(F));\\ 
)))));\\
for iAA13 from 0 to (\#AA13-1) list ( \\
for iA2 from 0 to (\#A2-1) list ( \\
for iA4 from 0 to (\#A4-1) list ( \\
for iB from 0 to (\#B-1) list ( \\
F=F0 + 19*AA13\#iAA13 + 23*A2\#iA2 + 31*A4\#iA4 + 37*B\#iB; TotalSing=TotalSing + dim singularLocus(ideal(F));\\ 
))));\\
for iA1 from 0 to (\#A1-1) list ( \\
for iAA24 from 0 to (\#AA24-1) list ( \\
for iA3 from 0 to (\#A3-1) list ( \\
for iB from 0 to (\#B-1) list ( \\
F=F0 + 19*A1\#iA1 + 23*AA24\#iAA24 + 29*A3\#iA3 + 37*B\#iB; TotalSing=TotalSing + dim singularLocus(ideal(F));\\ 
))));\\
for iAA13 from 0 to (\#AA13-1) list ( \\
for iAA24 from 0 to (\#AA24-1) list ( \\
for iB from 0 to (\#B-1) list ( \\
F=F0 + 19*AA13\#iAA13 + 23*AA24\#iAA24 + 37*B\#iB; TotalSing=TotalSing + dim singularLocus(ideal(F));\\ 
)));\\ 
F0=F00 + 17*x4*x4*x2; \\
for iA1 from 0 to (\#A1-1) list ( \\
for iA2 from 0 to (\#A2-1) list ( \\
for iA3 from 0 to (\#A3-1) list ( \\
for iA4 from 0 to (\#A4-1) list ( \\
for iB from 0 to (\#B-1) list ( \\
F=F0 + 19*A1\#iA1 + 23*A2\#iA2 + 29*A3\#iA3 + 31*A4\#iA4 + 37*B\#iB;\\ TotalSing=TotalSing + dim singularLocus(ideal(F));\\ 
)))));\\
for iAA13 from 0 to (\#AA13-1) list ( \\
for iA2 from 0 to (\#A2-1) list ( \\
for iA4 from 0 to (\#A4-1) list ( \\
for iB from 0 to (\#B-1) list ( \\
F=F0 + 19*AA13\#iAA13 + 23*A2\#iA2 + 31*A4\#iA4 + 37*B\#iB; TotalSing=TotalSing + dim singularLocus(ideal(F));\\ 
))));\\
for iA1 from 0 to (\#A1-1) list ( \\
for iAA24 from 0 to (\#AA24-1) list ( \\
for iA3 from 0 to (\#A3-1) list ( \\
for iB from 0 to (\#B-1) list ( \\
F=F0 + 19*A1\#iA1 + 23*AA24\#iAA24 + 29*A3\#iA3 + 37*B\#iB; TotalSing=TotalSing + dim singularLocus(ideal(F));\\ 
))));\\
for iAA13 from 0 to (\#AA13-1) list ( \\
for iAA24 from 0 to (\#AA24-1) list ( \\
for iB from 0 to (\#B-1) list ( \\
F=F0 + 19*AA13\#iAA13 + 23*AA24\#iAA24 + 37*B\#iB; TotalSing=TotalSing + dim singularLocus(ideal(F));\\ 
)));\\ 
F0=F00 + 17*x0*x0*x2; \\
for iA1 from 0 to (\#A1-1) list ( \\
for iA2 from 0 to (\#A2-1) list ( \\
for iA3 from 0 to (\#A3-1) list ( \\
for iA4 from 0 to (\#A4-1) list ( \\
for iB from 0 to (\#B-1) list ( \\
F=F0 + 19*A1\#iA1 + 23*A2\#iA2 + 29*A3\#iA3 + 31*A4\#iA4 + 37*B\#iB;\\ TotalSing=TotalSing + dim singularLocus(ideal(F));\\ 
)))));\\
for iAA13 from 0 to (\#AA13-1) list ( \\
for iA2 from 0 to (\#A2-1) list ( \\
for iA4 from 0 to (\#A4-1) list ( \\
for iB from 0 to (\#B-1) list ( \\
F=F0 + 19*AA13\#iAA13 + 23*A2\#iA2 + 31*A4\#iA4 + 37*B\#iB; TotalSing=TotalSing + dim singularLocus(ideal(F));\\ 
))));\\
for iA1 from 0 to (\#A1-1) list ( \\
for iAA24 from 0 to (\#AA24-1) list ( \\
for iA3 from 0 to (\#A3-1) list ( \\
for iB from 0 to (\#B-1) list ( \\
F=F0 + 19*A1\#iA1 + 23*AA24\#iAA24 + 29*A3\#iA3 + 37*B\#iB; TotalSing=TotalSing + dim singularLocus(ideal(F));\\ 
))));\\
for iAA13 from 0 to (\#AA13-1) list ( \\
for iAA24 from 0 to (\#AA24-1) list ( \\
for iB from 0 to (\#B-1) list ( \\
F=F0 + 19*AA13\#iAA13 + 23*AA24\#iAA24 + 37*B\#iB; TotalSing=TotalSing + dim singularLocus(ideal(F));\\ 
)));\\
F0=F00 + 17*x0*x4*x4; \\
for iA1 from 0 to (\#A1-1) list ( \\
for iA2 from 0 to (\#A2-1) list ( \\
for iA3 from 0 to (\#A3-1) list ( \\
for iA4 from 0 to (\#A4-1) list ( \\
F=F0 + 19*A1\#iA1 + 23*A2\#iA2 + 29*A3\#iA3 + 31*A4\#iA4; TotalSing=TotalSing + dim singularLocus(ideal(F));\\ 
))));\\
for iAA13 from 0 to (\#AA13-1) list ( \\
for iA2 from 0 to (\#A2-1) list ( \\
for iA4 from 0 to (\#A4-1) list ( \\
F=F0 + 19*AA13\#iAA13 + 23*A2\#iA2 + 31*A4\#iA4; TotalSing=TotalSing + dim singularLocus(ideal(F));\\ 
)));\\
for iA1 from 0 to (\#A1-1) list ( \\
for iAA24 from 0 to (\#AA24-1) list ( \\
for iA3 from 0 to (\#A3-1) list ( \\
F=F0 + 19*A1\#iA1 + 23*AA24\#iAA24 + 29*A3\#iA3; TotalSing=TotalSing + dim singularLocus(ideal(F));\\ 
)));\\ 
for iAA13 from 0 to (\#AA13-1) list ( \\
for iAA24 from 0 to (\#AA24-1) list ( \\
F=F0 + 19*AA13\#iAA13 + 23*AA24\#iAA24; TotalSing=TotalSing + dim singularLocus(ideal(F));\\ 
));\\
F0=F00 + 17*x0*x4*x5; \\
for iA1 from 0 to (\#A1-1) list ( \\
for iA2 from 0 to (\#A2-1) list ( \\
for iA3 from 0 to (\#A3-1) list ( \\
for iA4 from 0 to (\#A4-1) list ( \\
F=F0 + 19*A1\#iA1 + 23*A2\#iA2 + 29*A3\#iA3 + 31*A4\#iA4; TotalSing=TotalSing + dim singularLocus(ideal(F));\\ 
))));\\
for iAA13 from 0 to (\#AA13-1) list ( \\
for iA2 from 0 to (\#A2-1) list ( \\
for iA4 from 0 to (\#A4-1) list ( \\
F=F0 + 19*AA13\#iAA13 + 23*A2\#iA2 + 31*A4\#iA4; TotalSing=TotalSing + dim singularLocus(ideal(F));\\ 
)));\\
for iA1 from 0 to (\#A1-1) list ( \\
for iAA24 from 0 to (\#AA24-1) list ( \\
for iA3 from 0 to (\#A3-1) list ( \\
F=F0 + 19*A1\#iA1 + 23*AA24\#iAA24 + 29*A3\#iA3; TotalSing=TotalSing + dim singularLocus(ideal(F));\\ 
)));\\ 
for iAA13 from 0 to (\#AA13-1) list ( \\
for iAA24 from 0 to (\#AA24-1) list ( \\
F=F0 + 19*AA13\#iAA13 + 23*AA24\#iAA24; TotalSing=TotalSing + dim singularLocus(ideal(F));\\ 
));\\
F0=F00 + 17*x0*x3*x4; \\
for iA1 from 0 to (\#A1-1) list ( \\
for iA2 from 0 to (\#A2-1) list ( \\
for iA3 from 0 to (\#A3-1) list ( \\
for iA4 from 0 to (\#A4-1) list ( \\
F=F0 + 19*A1\#iA1 + 23*A2\#iA2 + 29*A3\#iA3 + 31*A4\#iA4; TotalSing=TotalSing + dim singularLocus(ideal(F));\\ 
))));\\
for iAA13 from 0 to (\#AA13-1) list ( \\
for iA2 from 0 to (\#A2-1) list ( \\
for iA4 from 0 to (\#A4-1) list ( \\
F=F0 + 19*AA13\#iAA13 + 23*A2\#iA2 + 31*A4\#iA4; TotalSing=TotalSing + dim singularLocus(ideal(F));\\ 
)));\\
for iA1 from 0 to (\#A1-1) list ( \\
for iAA24 from 0 to (\#AA24-1) list ( \\
for iA3 from 0 to (\#A3-1) list ( \\
F=F0 + 19*A1\#iA1 + 23*AA24\#iAA24 + 29*A3\#iA3; TotalSing=TotalSing + dim singularLocus(ideal(F));\\ 
)));\\ 
for iAA13 from 0 to (\#AA13-1) list ( \\
for iAA24 from 0 to (\#AA24-1) list ( \\
F=F0 + 19*AA13\#iAA13 + 23*AA24\#iAA24; TotalSing=TotalSing + dim singularLocus(ideal(F));\\ 
));\\
F0=F00 + 17*x0*x0*x3; \\
for iA1 from 0 to (\#A1-1) list ( \\
for iA2 from 0 to (\#A2-1) list ( \\
for iA3 from 0 to (\#A3-1) list ( \\
for iA4 from 0 to (\#A4-1) list ( \\
F=F0 + 19*A1\#iA1 + 23*A2\#iA2 + 29*A3\#iA3 + 31*A4\#iA4; TotalSing=TotalSing + dim singularLocus(ideal(F));\\ 
))));\\
for iAA13 from 0 to (\#AA13-1) list ( \\
for iA2 from 0 to (\#A2-1) list ( \\
for iA4 from 0 to (\#A4-1) list ( \\
F=F0 + 19*AA13\#iAA13 + 23*A2\#iA2 + 31*A4\#iA4; TotalSing=TotalSing + dim singularLocus(ideal(F));\\ 
)));\\
for iA1 from 0 to (\#A1-1) list ( \\
for iAA24 from 0 to (\#AA24-1) list ( \\
for iA3 from 0 to (\#A3-1) list ( \\
F=F0 + 19*A1\#iA1 + 23*AA24\#iAA24 + 29*A3\#iA3; TotalSing=TotalSing + dim singularLocus(ideal(F));\\ 
)));\\ 
for iAA13 from 0 to (\#AA13-1) list ( \\
for iAA24 from 0 to (\#AA24-1) list ( \\
F=F0 + 19*AA13\#iAA13 + 23*AA24\#iAA24; TotalSing=TotalSing + dim singularLocus(ideal(F));\\ 
));\\
G=[ x1\^{ }2*x3, x3\^{ }2*x1, x4\^{ }2*x1, x1\^{ }2*x4, x1\^{ }2*x5, x1*x3*x4, x1*x3*x5, x1*x4*x5, x3*x4*x5 ];\\
F1=[ x1*x3*x4, x1*x3*x5, x1*x4*x5, x3*x4*x5 ];\\
A1=[ x4\^{ }2*x1, x2\^{ }2*x1, x1*x2*x4, x2*x4*x5 ];\\
A2=[ x2\^{ }2*x5, x2\^{ }2*x1, x1*x2*x4, x2*x4*x5 ];\\
B1=[ x3\^{ }2*x1, x2\^{ }2*x1, x1*x2*x3, x2*x3*x5 ];\\
B2=[ x2\^{ }2*x5, x2\^{ }2*x1, x1*x2*x3, x2*x3*x5 ];\\
C1=[ x2*x3*x5, x2*x4*x5, x3*x4*x5 ];\\
C2=[ x2\^{ }2*x5, x2*x3*x5, x2*x4*x5, x3*x4*x5 ];\\
D=[ x3\^{ }2*x2, x3\^{ }2*x1, x4\^{ }2*x1, x4\^{ }2*x2, x2\^{ }2*x1, x1*x3*x4, x1*x2*x3, x1*x2*x4 ];\\
E1=[ x5\^{ }2*x2, x5\^{ }2*x0, x0\^{ }2*x2, x0\^{ }2*x3, x0\^{ }2*x4, x0*x2*x5, x0*x3*x5, x0*x4*x5 ];\\
E2=[ x5\^{ }2*x0, x1\^{ }2*x0, x0*x1*x5, x1*x3*x5, x1*x4*x5 ];\\
F00 = x0\^{ }2*x1 + 2*x1\^{ }2*x2 + 3*x2\^{ }2*x0 + 7*x3\^{ }2*x0 + 11*x4\^{ }2*x0 + 13*x5\^{ }2*x1; \\
F0=F00 + 17*x2*x3*x4; \\
for iE1 from 0 to (\#E1-1) list ( \\
for iG from 0 to (\#G-1) list ( \\
F=F0 + 19*E1\#iE1 + 23*G\#iG; TotalSing=TotalSing + dim singularLocus(ideal(F));\\ 
));\\
F0=F00 + 17*x1*x3*x4; \\
for iE1 from 0 to (\#E1-1) list ( \\
for iA1 from 0 to (\#A1-1) list ( \\
for iB1 from 0 to (\#B1-1) list ( \\
for iC1 from 0 to (\#C1-1) list ( \\
F=F0 + 19*E1\#iE1 + 23*A1\#iA1 + 29*B1\#iB1 + 31*C1\#iC1; TotalSing=TotalSing + dim singularLocus(ideal(F));\\ 
))));\\
F0=F00 + 17*x1*x3*x3; \\
for iE1 from 0 to (\#E1-1) list ( \\
for iA1 from 0 to (\#A1-1) list ( \\
for iC1 from 0 to (\#C1-1) list ( \\
F=F0 + 19*E1\#iE1 + 23*A1\#iA1 + 29*C1\#iC1; TotalSing=TotalSing + dim singularLocus(ideal(F));\\ 
)));\\
F0=F00 + 17*x5*x3*x4; \\
for iE1 from 0 to (\#E1-1) list ( \\
for iA1 from 0 to (\#A1-1) list ( \\
for iB1 from 0 to (\#B1-1) list ( \\
for iD from 0 to (\#D-1) list ( \\
F=F0 + 19*E1\#iE1 + 23*A1\#iA1 + 29*B1\#iB1 + 31*D\#iD; TotalSing=TotalSing + dim singularLocus(ideal(F));\\ 
))));\\
F00 = x0\^{ }2*x1 + 2*x1\^{ }2*x2 + 3*x2\^{ }2*x0 + 7*x3\^{ }2*x0 + 11*x4\^{ }2*x0 + 13*x5\^{ }2*x2; \\
F0=F00 + 17*x2*x3*x4; \\
for iE2 from 0 to (\#E2-1) list ( \\
for iF1 from 0 to (\#F1-1) list ( \\
F=F0 + 19*E2\#iE2 + 23*F1\#iF1; TotalSing=TotalSing + dim singularLocus(ideal(F));\\ 
));\\
F0=F00 + 17*x1*x3*x4; \\
for iE2 from 0 to (\#E2-1) list ( \\
for iA2 from 0 to (\#A2-1) list ( \\
for iB2 from 0 to (\#B2-1) list ( \\
for iC2 from 0 to (\#C2-1) list ( \\
F=F0 + 19*E2\#iE2 + 23*A2\#iA2 + 29*B2\#iB2 + 31*C2\#iC2; TotalSing=TotalSing + dim singularLocus(ideal(F));\\ 
))));\\
F0=F00 + 17*x5*x3*x4; \\
for iE2 from 0 to (\#E2-1) list ( \\
for iA2 from 0 to (\#A2-1) list ( \\
for iB2 from 0 to (\#B2-1) list ( \\
for iD from 0 to (\#D-1) list ( \\
F=F0 + 19*E2\#iE2 + 23*A2\#iA2 + 29*B2\#iB2 + 31*D\#iD; TotalSing=TotalSing + dim singularLocus(ideal(F));\\ 
))));\\
F1=[ x1*x2*x4, x1*x2*x5, x1*x4*x5, x2*x4*x5 ];\\
F2=[ x1*x2*x3, x1*x2*x4, x1*x3*x4, x2*x3*x4 ];\\
F3=[ x1*x2*x5, x2*x3*x5, x2*x4*x5 ];\\
F4=[ x1*x2*x4, x2*x3*x4, x2*x4*x5 ];\\
F5=[ x1*x2*x3, x2*x3*x5 ];\\
F6=[ x2*x4*x5, x2*x3*x4, x2*x3*x5 ];\\
F7=[ x2*x4*x5, x1*x4*x5 ];\\
F8=[ x2*x3*x5, x1*x3*x5 ];\\
F9=[ x1*x2*x3, x2*x3*x4, x2*x3*x5 ];\\
F10=[ x1*x3*x5, x2*x3*x5 ];\\
F00 = x0\^{ }2*x1 + 2*x1\^{ }2*x2 + 3*x2\^{ }2*x0 + 7*x3\^{ }2*x0 + 11*x4\^{ }2*x0 + 13*x5\^{ }2*x0; \\
F0=F00 + 17*x3*x4*x5; \\
for iF1 from 0 to (\#F1-1) list ( \\
for iF2 from 0 to (\#F2-1) list ( \\
for iF3 from 0 to (\#F3-1) list ( \\
for iF4 from 0 to (\#F4-1) list ( \\
for iF5 from 0 to (\#F5-1) list ( \\
F=F0 + 19*F1\#iF1 + 23*F2\#iF2 + 29*F3\#iF3 + 31*F4\#iF4 + 37*F5\#iF5;\\ TotalSing=TotalSing + dim singularLocus(ideal(F));\\ 
)))));\\
F0=F00 + 17*x3*x4*x1 + 41*x2*x2*x1; \\
for iF6 from 0 to (\#F6-1) list ( \\
for iF7 from 0 to (\#F7-1) list ( \\
for iF8 from 0 to (\#F8-1) list ( \\
F=F0 + 19*F6\#iF6 + 23*F7\#iF7 + 29*F8\#iF8; TotalSing=TotalSing + dim singularLocus(ideal(F));\\ 
)));\\
F0=F00 + 17*x2*x4*x5 + 41*x1*x3*x4; \\
for iF9 from 0 to (\#F9-1) list ( \\
for iF10 from 0 to (\#F10-1) list ( \\
F=F0 + 19*F9\#iF9 + 23*F10\#iF10; TotalSing=TotalSing + dim singularLocus(ideal(F));\\ 
));\\ 
F=F00 + 17*x1*x2*x5 + 19*x1*x3*x4+ 19*x1*x4*x5+23*x2*x3*x4+29*x1*x3*x5;\\ TotalSing=TotalSing + dim singularLocus(ideal(F));\\ 
F=F00 + 17*x3*x4*x5 + 41*x1*x2*x2; TotalSing=TotalSing + dim singularLocus(ideal(F));\\ 
print "This is Case of No Cubes, Length 2 longest cycle";\\
A=[ x3\^{ }2*x0, x3\^{ }2*x1, x3\^{ }2*x4, x5\^{ }2*x0, x5\^{ }2*x1, x0*x3*x5, x1*x3*x5, x3*x4*x5 ];\\
B0=[ x3\^{ }2*x0, x0\^{ }2*x3, x0\^{ }2*x4, x3\^{ }2*x4, x5\^{ }2*x0, x0\^{ }2*x5, x0*x3*x4, x0*x3*x5, x0*x4*x5, x3*x4*x5 ];\\
B1=[ x3\^{ }2*x1, x1\^{ }2*x3, x1\^{ }2*x5, x3\^{ }2*x4, x5\^{ }2*x1, x1*x3*x4, x1*x3*x5, x3*x4*x5 ];\\
B4=[ x3\^{ }2*x1, x4\^{ }2*x3, x3\^{ }2*x4, x5\^{ }2*x1, x1*x3*x4, x1*x3*x5 ];\\
F1=[ x3\^{ }2*x1, x1\^{ }2*x3, x1*x3*x4, x1*x3*x5, x1*x4*x5 ]; \\
F2=[ x0\^{ }2*x3, x0\^{ }2*x4, x0\^{ }2*x5, x0*x3*x4, x0*x3*x5, x0*x4*x5 ]; \\
F3=[ x4\^{ }2*x3, x3\^{ }2*x4, x1*x3*x4, x1*x3*x5 ]; \\
F4=[ x4\^{ }2*x3, x5\^{ }2*x1, x1*x3*x4, x1*x3*x5 ]; \\
F5=[ x4\^{ }2*x2, x2\^{ }2*x1, x1*x2*x4, x1*x2*x5, x2*x4*x5 ]; \\
F6=[ x3\^{ }2*x0, x0*x3*x4, x0*x3*x5, x0*x4*x5 ]; \\
F7=[ x1\^{ }2*x2, x4\^{ }2*x2, x1*x2*x4, x1*x2*x5, x2*x4*x5 ]; \\
F8=[ x2\^{ }2*x1, x2\^{ }2*x4, x2\^{ }2*x5, x4\^{ }2*x2, x1*x2*x4, x1*x2*x5, x2*x4*x5 ]; \\
F9=[ x1*x3*x4, x1*x3*x5 ];\\
F = x0\^{ }2*x1 + x1\^{ }2*x0 + x2\^{ }2*x3 + x3\^{ }2*x2 + x4\^{ }2*x5 + x5\^{ }2*x4; TotalSing=TotalSing + dim singularLocus(ideal(F));\\ 
F00 = x0\^{ }2*x1 + 2*x1\^{ }2*x0 + 3*x2\^{ }2*x3 + 7*x3\^{ }2*x2 + 11*x4\^{ }2*x2 + 13*x5\^{ }2*x4; \\
F=F00 + 61*x3*x4*x5; TotalSing=TotalSing + dim singularLocus(ideal(F));\\ 
F=F00 + 61*x3*x4*x4; TotalSing=TotalSing + dim singularLocus(ideal(F));\\ 
F0=F00 + 61*x0*x4*x4;\\
for iF1 from 0 to (\#F1-1) list ( \\
F=F0 + 19*F1\#iF1; TotalSing=TotalSing + dim singularLocus(ideal(F));\\ 
);\\
F0=F00 + 61*x1*x3*x4;\\
for iF2 from 0 to (\#F2-1) list ( \\
F=F0 + 19*F2\#iF2; TotalSing=TotalSing + dim singularLocus(ideal(F));\\ 
);\\
F0=F00 + 61*x1*x3*x3;\\
for iF6 from 0 to (\#F6-1) list ( \\
F=F0 + 19*F6\#iF6; TotalSing=TotalSing + dim singularLocus(ideal(F));\\ 
);\\
F00 = x0\^{ }2*x1 + 2*x1\^{ }2*x0 + 3*x2\^{ }2*x3 + 7*x3\^{ }2*x2 + 11*x4\^{ }2*x0 + 13*x5\^{ }2*x2; \\
F=F00 + 61*x1*x4*x4 + 71*x3*x5*x5; TotalSing=TotalSing + dim singularLocus(ideal(F));\\ 
F0=F00 + 61*x1*x1*x2;\\
for iF3 from 0 to (\#F3-1) list ( \\
for iF4 from 0 to (\#F4-1) list ( \\
for iA from 0 to (\#A-1) list ( \\
for iB0 from 0 to (\#B0-1) list ( \\
for iB1 from 0 to (\#B1-1) list ( \\
for iB4 from 0 to (\#B4-1) list ( \\
F=F0 + 19*F3\#iF3 + 23*F4\#iF4 + 29*A\#iA + 73*B0\#iB0 + 31*B1\#iB1 + 37*B4\#iB4;\\ TotalSing=TotalSing + dim singularLocus(ideal(F));\\ 
))))));\\
F0=F00 + 61*x3*x4*x4;\\
for iF5 from 0 to (\#F5-1) list ( \\
for iA from 0 to (\#A-1) list ( \\
for iB0 from 0 to (\#B0-1) list ( \\
for iB1 from 0 to (\#B1-1) list ( \\
F=F0 + 23*F5\#iF5 + 29*A\#iA + 73*B0\#iB0 + 31*B1\#iB1; TotalSing=TotalSing + dim singularLocus(ideal(F));\\ 
))));\\ 
F0=F00 + 61*x3*x1*x1;\\
for iF7 from 0 to (\#F7-1) list ( \\
for iA from 0 to (\#A-1) list ( \\
for iB0 from 0 to (\#B0-1) list ( \\
for iB4 from 0 to (\#B4-1) list ( \\
F=F0 + 23*F7\#iF7 + 29*A\#iA + 73*B0\#iB0 + 31*B4\#iB4; TotalSing=TotalSing + dim singularLocus(ideal(F));\\ 
))));\\ 
F0=F00 + 61*x3*x4*x1;\\
for iF8 from 0 to (\#F8-1) list ( \\
for iA from 0 to (\#A-1) list ( \\
for iB0 from 0 to (\#B0-1) list ( \\
F=F0 + 23*F8\#iF8 + 29*A\#iA + 73*B0\#iB0; TotalSing=TotalSing + dim singularLocus(ideal(F));\\ 
)));\\ 
F0=F00 + 61*x5*x1*x1;\\
for iA from 0 to (\#A-1) list ( \\
for iB0 from 0 to (\#B0-1) list ( \\
for iB4 from 0 to (\#B4-1) list ( \\
F=F0 + 23*B4\#iB4 + 29*A\#iA + 73*B0\#iB0; TotalSing=TotalSing + dim singularLocus(ideal(F));\\ 
)));\\ 
F0=F00 + 61*x4*x4*x1;\\
for iA from 0 to (\#A-1) list ( \\
for iB0 from 0 to (\#B0-1) list ( \\
for iB1 from 0 to (\#B1-1) list ( \\
F=F0 + 23*B1\#iB1 + 29*A\#iA + 73*B0\#iB0; TotalSing=TotalSing + dim singularLocus(ideal(F));\\ 
)));\\ 
F0=F00 + 61*x2*x4*x1;\\
for iA from 0 to (\#A-1) list ( \\
for iB0 from 0 to (\#B0-1) list ( \\
for iF9 from 0 to (\#F9-1) list ( \\
F=F0 + 23*F9\#iF9 + 29*A\#iA + 73*B0\#iB0; TotalSing=TotalSing + dim singularLocus(ideal(F));\\ 
)));\\ 
F0=F00 + 61*x2*x4*x4;\\
for iA from 0 to (\#A-1) list ( \\
for iB0 from 0 to (\#B0-1) list ( \\
for iF9 from 0 to (\#F9-1) list ( \\
F=F0 + 23*F9\#iF9 + 29*A\#iA + 73*B0\#iB0; TotalSing=TotalSing + dim singularLocus(ideal(F));\\ 
)));\\ 
F0=F00 + 61*x1*x4*x5;\\
for iA from 0 to (\#A-1) list ( \\
for iB0 from 0 to (\#B0-1) list ( \\
F=F0 +  29*A\#iA + 73*B0\#iB0; TotalSing=TotalSing + dim singularLocus(ideal(F));\\ 
));\\ 
F1=[ x0\^{ }2*x2, x0*x2*x5, x0*x4*x5 ]; \\
F2=[ x0\^{ }2*x3, x3\^{ }2*x5, x3\^{ }2*x4, x3\^{ }2*x0, x0*x3*x4, x0*x3*x5, x0*x4*x5, x3*x4*x5 ]; \\
F3=[ x1\^{ }2*x2, x1\^{ }2*x3, x1*x2*x5, x1*x3*x5, x1*x4*x5 ]; \\
F4=[ x1\^{ }2*x2, x1*x2*x4, x1*x2*x5, x1*x4*x5, x2*x4*x5 ]; \\
F5=[ x1\^{ }2*x3, x1*x3*x4, x1*x3*x5, x1*x4*x5, x3*x4*x5 ]; \\
F6=[ x0\^{ }2*x3, x0*x3*x4, x0*x3*x5, x0*x4*x5, x3*x4*x5 ];\\
F7=[ x1*x2*x5, x1*x3*x5, x1*x4*x5 ];\\
F8=[ x2\^{ }2*x1, x2\^{ }2*x4, x2\^{ }2*x5, x1*x2*x4, x1*x2*x5, x1*x4*x5, x2*x4*x5 ];\\
F9=[ x3\^{ }2*x1, x3\^{ }2*x4, x3\^{ }2*x5, x1*x3*x4, x1*x3*x5, x1*x4*x5, x3*x4*x5 ];\\
F10=[ x3\^{ }2*x0, x3\^{ }2*x4, x3\^{ }2*x5, x0*x3*x4, x0*x3*x5, x0*x4*x5, x3*x4*x5 ];\\
G1=[ x1\^{ }2*x3, x1\^{ }2*x2, x1*x2*x4, x1*x3*x4 ];\\
G2=[ x1\^{ }2*x3, x1\^{ }2*x2, x1*x2*x5, x1*x3*x5 ];\\
G3=[ x1\^{ }2*x3, x1*x3*x4, x5*x3*x4 ];\\
G4=[ x3\^{ }2*x1, x3\^{ }2*x4, x3\^{ }2*x5, x1*x3*x4, x1*x3*x5, x3*x4*x5 ];\\
F0 = x0\^{ }2*x1 + 2*x1\^{ }2*x0 + 3*x2\^{ }2*x3 + 7*x3\^{ }2*x2 + 11*x4\^{ }2*x0 + 12*x4\^{ }2*x1 + 13*x5\^{ }2*x1; \\
for iF1 from 0 to (\#F1-1) list ( \\
for iF2 from 0 to (\#F2-1) list ( \\
F=F0 + 19*F1\#iF1 + 23*F2\#iF2; TotalSing=TotalSing + dim singularLocus(ideal(F));\\ 
));\\
F0 = x0\^{ }2*x1 + 2*x1\^{ }2*x0 + 3*x2\^{ }2*x3 + 7*x3\^{ }2*x2 + 11*x4\^{ }2*x1 + 12*x0\^{ }2*x2 + 13*x5\^{ }2*x0; \\
for iF3 from 0 to (\#F3-1) list ( \\
for iF4 from 0 to (\#F4-1) list ( \\
for iF5 from 0 to (\#F5-1) list ( \\
for iF6 from 0 to (\#F6-1) list ( \\
F=F0 + 19*F3\#iF3 + 23*F4\#iF4 + 29*F5\#iF5 + 73*F6\#iF6; TotalSing=TotalSing + dim singularLocus(ideal(F));\\ 
))));\\
F0 = x0\^{ }2*x1 + 2*x1\^{ }2*x0 + 3*x2\^{ }2*x3 + 7*x3\^{ }2*x2 + 11*x4\^{ }2*x1 + 12*x0*x2*x4 + 13*x5\^{ }2*x0; \\
for iF7 from 0 to (\#F7-1) list ( \\
for iF8 from 0 to (\#F8-1) list ( \\
for iF9 from 0 to (\#F9-1) list ( \\
for iF10 from 0 to (\#F10-1) list ( \\
F=F0 + 19*F7\#iF7 + 23*F8\#iF8 + 29*F9\#iF9 + 73*F10\#iF10; TotalSing=TotalSing + dim singularLocus(ideal(F));\\ 
))));\\
F0 = x0\^{ }2*x1 + 2*x1\^{ }2*x0 + 3*x2\^{ }2*x3 + 7*x3\^{ }2*x2 + 11*x4\^{ }2*x0 + 12*x5*x2*x4 + 13*x5\^{ }2*x0; \\
for iG1 from 0 to (\#G1-1) list ( \\
for iG2 from 0 to (\#G2-1) list ( \\
for iG3 from 0 to (\#G3-1) list ( \\
for iG4 from 0 to (\#G4-1) list ( \\
F=F0 + 19*G1\#iG1 + 23*G2\#iG2 + 29*G3\#iG3 + 73*G4\#iG4; TotalSing=TotalSing + dim singularLocus(ideal(F));\\ 
))));\\
F = x0\^{ }2*x1 + 2*x1\^{ }2*x0 + 3*x2\^{ }2*x3 + 7*x3\^{ }2*x2 + 11*x4\^{ }2*x0 + 31*x4\^{ }2*x1 + 13*x5\^{ }2*x0 + 23*x5\^{ }2*x1;\\ TotalSing=TotalSing + dim singularLocus(ideal(F));\\ 
F = x0\^{ }2*x1 + 2*x1\^{ }2*x0 + 3*x2\^{ }2*x3 + 7*x3\^{ }2*x2 + 11*x4\^{ }2*x1 + 31*x4*x5*x0 + 13*x5\^{ }2*x0 + 23*x4*x5*x1;\\ TotalSing=TotalSing + dim singularLocus(ideal(F));\\ 
F = x0\^{ }2*x1 + 2*x1\^{ }2*x0 + 3*x2\^{ }2*x3 + 7*x3\^{ }2*x2 + 11*x4\^{ }2*x0 + 31*x4*x5*x1 + 13*x5\^{ }2*x0;\\ TotalSing=TotalSing + dim singularLocus(ideal(F));\\
A1=[ x5\^{ }2*x2, x5\^{ }2*x3, x0\^{ }2*x2, x0*x2*x5, x0*x3*x5, x0*x4*x5 ]; \\
A2=[ x5\^{ }2*x3, x3\^{ }2*x5, x3\^{ }2*x0, x0*x3*x4, x0*x3*x5, x0*x4*x5, x3*x4*x5 ]; \\
A3=[ x4\^{ }2*x0, x4\^{ }2*x2, x4\^{ }2*x5, x0\^{ }2*x2, x5\^{ }2*x2, x0*x2*x4, x0*x2*x5, x0*x4*x5, x2*x4*x5 ]; \\
B=[ x4\^{ }2*x1, x4\^{ }2*x2, x4\^{ }2*x5, x1*x2*x4, x1*x2*x5, x1*x4*x5, x2*x4*x5 ];\\
F1=[ x4\^{ }2*x1, x4\^{ }2*x2, x1*x4*x5, x2*x4*x5 ];\\
F2=[ x5\^{ }2*x1, x5\^{ }2*x2, x5\^{ }2*x3, x1*x2*x3, x2*x3*x5, x1*x3*x5 ];\\
F00 = x0\^{ }2*x1 + 2*x1\^{ }2*x0 + 3*x2\^{ }2*x0 + 7*x3\^{ }2*x2 + 11*x4\^{ }2*x3 + 13*x5\^{ }2*x1; \\
F0 = F00 + 17*x1*x2*x5;\\
for iA1 from 0 to (\#A1-1) list ( \\
for iA2 from 0 to (\#A2-1) list ( \\
for iA3 from 0 to (\#A3-1) list ( \\
F=F0 + 19*A1\#iA1 + 23*A2\#iA2 + 31*A3\#iA3; TotalSing=TotalSing + dim singularLocus(ideal(F));\\ 
)));\\
F0 = F00 + 17*x1*x2*x2;\\
for iA1 from 0 to (\#A1-1) list ( \\
for iA2 from 0 to (\#A2-1) list ( \\
for iA3 from 0 to (\#A3-1) list ( \\
F=F0 + 19*A1\#iA1 + 23*A2\#iA2 + 31*A3\#iA3; TotalSing=TotalSing + dim singularLocus(ideal(F));\\ 
)));\\
F0 = F00 + 17*x1*x2*x4;\\
for iA1 from 0 to (\#A1-1) list ( \\
for iA2 from 0 to (\#A2-1) list ( \\
for iA3 from 0 to (\#A3-1) list ( \\
F=F0 + 19*A1\#iA1 + 23*A2\#iA2 + 31*A3\#iA3; TotalSing=TotalSing + dim singularLocus(ideal(F));\\ 
)));\\
F0 = F00 + 17*x1*x2*x3;\\
for iA1 from 0 to (\#A1-1) list ( \\
for iA2 from 0 to (\#A2-1) list ( \\
for iA3 from 0 to (\#A3-1) list ( \\
for iB from 0 to (\#B-1) list ( \\
F=F0 + 19*A1\#iA1 + 23*A2\#iA2 + 31*A3\#iA3 + 37*B\#iB; TotalSing=TotalSing + dim singularLocus(ideal(F));\\ 
))));\\
F0 = F00 + 17*x1*x2*x3 + 23*x0*x5*x5;\\
for iB from 0 to (\#B-1) list ( \\
F=F0 + 19*B\#iB; TotalSing=TotalSing + dim singularLocus(ideal(F));\\ 
);\\
F = F00 + 17*x1*x2*x2 + 23*x0*x5*x5; TotalSing=TotalSing + dim singularLocus(ideal(F));\\ 
F = F00 + 17*x1*x2*x4 + 23*x0*x5*x5; TotalSing=TotalSing + dim singularLocus(ideal(F));\\ 
F = F00 + 17*x1*x2*x5 + 23*x0*x5*x5; TotalSing=TotalSing + dim singularLocus(ideal(F));\\ 
F00 = x0\^{ }2*x1 + 2*x1\^{ }2*x0 + 3*x2\^{ }2*x0 + 7*x3\^{ }2*x2 + 11*x4\^{ }2*x3 + 13*x5\^{ }2*x4; \\
F0 = F00 + 17*x1*x2*x3;\\
for iF1 from 0 to (\#F1-1) list ( \\
F=F0 + 19*F1\#iF1; TotalSing=TotalSing + dim singularLocus(ideal(F));\\ 
);\\
F0 = F00 + 17*x1*x2*x4;\\
for iF2 from 0 to (\#F2-1) list ( \\
F=F0 + 19*F2\#iF2; TotalSing=TotalSing + dim singularLocus(ideal(F));\\ 
);\\
F = F00 + 17*x1*x2*x2; TotalSing=TotalSing + dim singularLocus(ideal(F));\\ 
F = F00 + 17*x1*x2*x5; TotalSing=TotalSing + dim singularLocus(ideal(F));\\
C=[ x3\^{ }2*x1, x3\^{ }2*x5, x1*x3*x4, x1*x3*x5, x1*x4*x5, x3*x4*x5 ]; \\
F0 = x0\^{ }2*x1 + 2*x1\^{ }2*x0 + 3*x2\^{ }2*x0 + 7*x3\^{ }2*x2 + 11*x4\^{ }2*x3 + 13*x5\^{ }2*x0 + 23*x1*x2*x5; \\
for iC from 0 to (\#C-1) list ( \\
F=F0 + 19*C\#iC; TotalSing=TotalSing + dim singularLocus(ideal(F));\\ 
);\\
F1=[ x1*x3*x5, x1*x4*x5 ]; \\
F2=[ x2*x3*x5, x2*x4*x5 ];\\
F3=[ x1*x3*x5, x2*x3*x5 ];\\
F0 = x0\^{ }2*x1 + 2*x1\^{ }2*x0 + 3*x2\^{ }2*x0 + 7*x3\^{ }2*x2 + 11*x4\^{ }2*x3 + 13*x5\^{ }2*x0 + 23*x1*x2*x4; \\
for iF1 from 0 to (\#F1-1) list ( \\
for iF2 from 0 to (\#F2-1) list ( \\
for iF3 from 0 to (\#F3-1) list ( \\
F=F0 + 19*F1\#iF1 + 29*F2\#iF2 + 31*F3\#iF3; TotalSing=TotalSing + dim singularLocus(ideal(F));\\ 
)));\\ 
F1=[ x1*x3*x5, x1*x4*x5 ]; \\
F2=[ x2*x3*x5, x2*x4*x5 ];\\
F3=[ x1*x2*x4, x1*x4*x5, x2*x4*x5 ]; \\
F4=[ x4\^{ }2*x1, x4\^{ }2*x2, x4\^{ }2*x5, x1*x2*x4, x1*x2*x5, x1*x4*x5, x2*x4*x5 ]; \\
F0 = x0\^{ }2*x1 + 2*x1\^{ }2*x0 + 3*x2\^{ }2*x0 + 7*x3\^{ }2*x2 + 11*x4\^{ }2*x3 + 13*x5\^{ }2*x0 + 23*x1*x2*x3; \\
for iF1 from 0 to (\#F1-1) list ( \\
for iF2 from 0 to (\#F2-1) list ( \\
for iF3 from 0 to (\#F3-1) list ( \\
for iF4 from 0 to (\#F4-1) list ( \\
F=F0 + 19*F1\#iF1 + 29*F2\#iF2 + 31*F3\#iF3 + 37*F4\#iF4; TotalSing=TotalSing + dim singularLocus(ideal(F));\\ 
))));\\
F1=[ x4\^{ }2*x2, x1*x4*x5, x2*x4*x5 ];\\
F0 = x0\^{ }2*x1 + 2*x1\^{ }2*x0 + 3*x2\^{ }2*x0 + 7*x3\^{ }2*x2 + 11*x4\^{ }2*x3 + 13*x5\^{ }2*x3 + 61*x5\^{ }2*x2; \\
for iF1 from 0 to (\#F1-1) list ( \\
F=F0 +31*x1*x2*x3 + 19*F1\#iF1; TotalSing=TotalSing + dim singularLocus(ideal(F));\\ 
);\\
F=F0 +31*x1*x2*x2; TotalSing=TotalSing + dim singularLocus(ideal(F));\\ 
F=F0 +31*x1*x2*x5; TotalSing=TotalSing + dim singularLocus(ideal(F));\\ 
F=F0 +31*x1*x2*x4; TotalSing=TotalSing + dim singularLocus(ideal(F));\\ 
AA1=[ x3\^{ }2*x0, x5\^{ }2*x0, x0*x3*x5, x3*x4*x5 ]; \\
AA2=[ x3\^{ }2*x1, x1*x3*x5, x1*x3*x4, x1*x4*x5, x3*x4*x5 ];\\
A1=[ x3\^{ }2*x1, x1*x3*x5 ];\\
A2=[ x0*x4*x5, x3*x4*x0 ];\\
BB=[ x5\^{ }2*x0, x4\^{ }2*x0, x0*x4*x5, x2*x4*x5 ]; \\
B1=[ x0\^{ }2*x2, x0*x2*x5, x0*x2*x4 ];\\
B2=[ x1*x4*x5 ];\\
B11=[ x1*x4*x5, x2*x4*x5, x1*x2*x4 ];\\
C=[ x4\^{ }2*x1, x4\^{ }2*x2, x4\^{ }2*x5, x1*x2*x4, x1*x4*x5, x1*x2*x5, x2*x4*x5 ];\\
F00 = x0\^{ }2*x1 + 2*x1\^{ }2*x0 + 3*x2\^{ }2*x0 + 7*x3\^{ }2*x2 + 11*x4\^{ }2*x3 + 13*x5\^{ }2*x3; \\
F0 = F00 + 17*x1*x2*x2;\\
for iB1 from 0 to (\#B1-1) list ( \\
for iB2 from 0 to (\#B2-1) list ( \\
for iB11 from 0 to (\#B11-1) list ( \\
F=F0 + 19*B1\#iB1 + 23*B2\#iB2 + 31*B11\#iB11; TotalSing=TotalSing + dim singularLocus(ideal(F));\\ 
)));\\
for iBB from 0 to (\#BB-1) list ( \\
for iB11 from 0 to (\#B11-1) list ( \\
F=F0 + 19*BB\#iBB  + 31*B11\#iB11; TotalSing=TotalSing + dim singularLocus(ideal(F));\\ 
));\\ 
F0 = F00 + 17*x1*x2*x3;\\
for iB1 from 0 to (\#B1-1) list ( \\
for iB2 from 0 to (\#B2-1) list ( \\
for iB11 from 0 to (\#B11-1) list ( \\
F=F0 + 19*B1\#iB1 + 23*B2\#iB2 + 31*B11\#iB11; TotalSing=TotalSing + dim singularLocus(ideal(F));\\ 
)));\\
for iBB from 0 to (\#BB-1) list ( \\
for iB11 from 0 to (\#B11-1) list ( \\
F=F0 + 19*BB\#iBB  + 31*B11\#iB11; TotalSing=TotalSing + dim singularLocus(ideal(F));\\ 
));\\ 
F0 = F00 + 17*x1*x2*x4;\\
for iB1 from 0 to (\#B1-1) list ( \\
for iB2 from 0 to (\#B2-1) list ( \\
F=F0 + 19*B1\#iB1 + 23*B2\#iB2; TotalSing=TotalSing + dim singularLocus(ideal(F));\\ 
));\\
for iBB from 0 to (\#BB-1) list ( \\
F=F0 + 19*BB\#iBB; TotalSing=TotalSing + dim singularLocus(ideal(F));\\ 
);\\
F00 = x0\^{ }2*x1 + 2*x1\^{ }2*x0 + 3*x2\^{ }2*x0 + 7*x3\^{ }2*x2 + 11*x4\^{ }2*x3 + 13*x5\^{ }2*x2; \\
F0 = F00 + 17*x1*x2*x2;\\
for iA1 from 0 to (\#A1-1) list ( \\
for iA2 from 0 to (\#A2-1) list ( \\
F=F0 + 19*A1\#iA1 + 23*A2\#iA2; TotalSing=TotalSing + dim singularLocus(ideal(F));\\ 
));\\
for iAA1 from 0 to (\#AA1-1) list ( \\
for iAA2 from 0 to (\#AA2-1) list ( \\
F=F0 + 19*AA1\#iAA1 + 23*AA2\#iAA2; TotalSing=TotalSing + dim singularLocus(ideal(F));\\ 
));\\
F0 = F00 + 17*x1*x2*x5;\\
for iA1 from 0 to (\#A1-1) list ( \\
for iA2 from 0 to (\#A2-1) list ( \\
F=F0 + 19*A1\#iA1 + 23*A2\#iA2; TotalSing=TotalSing + dim singularLocus(ideal(F));\\ 
));\\
for iAA1 from 0 to (\#AA1-1) list ( \\
for iAA2 from 0 to (\#AA2-1) list ( \\
F=F0 + 19*AA1\#iAA1 + 23*AA2\#iAA2; TotalSing=TotalSing + dim singularLocus(ideal(F));\\ 
));\\
F0 = F00 + 17*x1*x2*x4;\\
for iA1 from 0 to (\#A1-1) list ( \\
for iA2 from 0 to (\#A2-1) list ( \\
F=F0 + 19*A1\#iA1 + 23*A2\#iA2; TotalSing=TotalSing + dim singularLocus(ideal(F));\\ 
));\\
for iAA1 from 0 to (\#AA1-1) list ( \\
for iAA2 from 0 to (\#AA2-1) list ( \\
F=F0 + 19*AA1\#iAA1 + 23*AA2\#iAA2; TotalSing=TotalSing + dim singularLocus(ideal(F));\\ 
));\\
F0 = F00 + 17*x1*x2*x3;\\
for iA1 from 0 to (\#A1-1) list ( \\
for iA2 from 0 to (\#A2-1) list ( \\
for iC from 0 to (\#C-1) list ( \\
F=F0 + 19*A1\#iA1 + 23*A2\#iA2 + 31*C\#iC; TotalSing=TotalSing + dim singularLocus(ideal(F));\\ 
)));\\
for iAA1 from 0 to (\#AA1-1) list ( \\
for iAA2 from 0 to (\#AA2-1) list ( \\
for iC from 0 to (\#C-1) list ( \\
F=F0 + 19*AA1\#iAA1 + 23*AA2\#iAA2 + 31*C\#iC; TotalSing=TotalSing + dim singularLocus(ideal(F));\\ 
)));\\
AA=[ x0*x4*x4, x0*x3*x4, x0*x4*x5 ]; \\
A1=[ x0*x2*x4 ];\\
A2=[ x3\^{ }2*x0, x3\^{ }2*x4, x3\^{ }2*x5, x0*x3*x5, x3*x4*x5 ];\\
B=[ x5\^{ }2*x1, x5\^{ }2*x2, x5\^{ }2*x3, x1*x2*x3, x1*x2*x5, x1*x3*x5, x2*x3*x5 ]; \\
F0=x0\^{ }2*x1 + 2*x1\^{ }2*x0 + 3*x2\^{ }2*x0 + 7*x3\^{ }2*x2 + 11*x4\^{ }2*x1 + 13*x5\^{ }2*x4 + 17*x1*x2*x3;\\
for iA1 from 0 to (\#A1-1) list ( \\
for iA2 from 0 to (\#A2-1) list ( \\
F=F0 + 19*A1\#iA1 + 20*A2\#iA2; TotalSing=TotalSing + dim singularLocus(ideal(F));\\ 
));\\
for iAA from 0 to (\#AA-1) list ( \\
F=F0 + 19*AA\#iAA; TotalSing=TotalSing + dim singularLocus(ideal(F));\\ 
);\\
F0=x0\^{ }2*x1 + 2*x1\^{ }2*x0 + 3*x2\^{ }2*x0 + 7*x3\^{ }2*x2 + 11*x4\^{ }2*x1 + 13*x5\^{ }2*x4 + 17*x1*x2*x4;\\
for iA1 from 0 to (\#A1-1) list ( \\
for iA2 from 0 to (\#A2-1) list ( \\
for iB from 0 to (\#B-1) list ( \\
F=F0 + 19*A1\#iA1 + 20*A2\#iA2 + 23*B\#iB; TotalSing=TotalSing + dim singularLocus(ideal(F));\\ 
)));\\
for iAA from 0 to (\#AA-1) list ( \\
for iB from 0 to (\#B-1) list ( \\
F=F0 + 19*AA\#iAA + 23*B\#iB; TotalSing=TotalSing + dim singularLocus(ideal(F));\\ 
));\\
F0=x0\^{ }2*x1 + 2*x1\^{ }2*x0 + 3*x2\^{ }2*x0 + 7*x3\^{ }2*x2 + 11*x4\^{ }2*x1 + 13*x5\^{ }2*x4 + 17*x1*x2*x5;\\
for iA1 from 0 to (\#A1-1) list ( \\
for iA2 from 0 to (\#A2-1) list ( \\
F=F0 + 19*A1\#iA1 + 20*A2\#iA2; TotalSing=TotalSing + dim singularLocus(ideal(F));\\ 
));\\
for iAA from 0 to (\#AA-1) list ( \\
F=F0 + 19*AA\#iAA; TotalSing=TotalSing + dim singularLocus(ideal(F));\\ 
);\\
F=x0\^{ }2*x1 + 2*x1\^{ }2*x0 + 3*x2\^{ }2*x0 + 73*x2\^{ }2*x1 + 7*x3\^{ }2*x2 + 11*x4\^{ }2*x1 + 117*x4\^{ }2*x0 + 13*x5\^{ }2*x4;\\ TotalSing=TotalSing + dim singularLocus(ideal(F));\\ 
F1=[ x1*x3*x4, x1*x4*x5 ]; \\
F2=[ x1*x2*x3, x1*x2*x5 ];\\
F3=[ x1*x4*x5, x2*x4*x5, x1*x2*x5 ];\\
F0=x0\^{ }2*x1 + 2*x1\^{ }2*x0 + 3*x2\^{ }2*x0 + 7*x3\^{ }2*x2 + 11*x4\^{ }2*x0 + 13*x5\^{ }2*x4 + 17*x2*x3*x4;\\
for iF1 from 0 to (\#F1-1) list ( \\
for iF2 from 0 to (\#F2-1) list ( \\
for iF3 from 0 to (\#F3-1) list ( \\
F=F0 + 19*F1\#iF1 + 20*F2\#iF2 + 23*F3\#iF3; TotalSing=TotalSing + dim singularLocus(ideal(F));\\ 
)));\\
F1=[ x3\^{ }2*x1, x3\^{ }2*x4, x3\^{ }2*x5, x1*x3*x4, x1*x4*x5, x1*x3*x5, x3*x4*x5 ]; \\
F2=[ x5\^{ }2*x1, x5\^{ }2*x2, x5\^{ }2*x3, x1*x2*x3, x1*x2*x5, x1*x3*x5, x2*x3*x5 ]; \\
F0=x0\^{ }2*x1 + 2*x1\^{ }2*x0 + 3*x2\^{ }2*x0 + 7*x3\^{ }2*x2 + 11*x4\^{ }2*x0 + 13*x5\^{ }2*x4 + 17*x1*x2*x4;\\
for iF1 from 0 to (\#F1-1) list ( \\
for iF2 from 0 to (\#F2-1) list ( \\
F=F0 + 19*F1\#iF1 + 23*F2\#iF2; TotalSing=TotalSing + dim singularLocus(ideal(F));\\ 
));\\  
AA=[ x0*x5*x5, x0*x3*x5, x0*x4*x5 ]; \\
A1=[ x3\^{ }2*x0, x0*x3*x4, x3*x4*x5 ];\\
A2=[ x4\^{ }2*x0, x0*x3*x4, x3*x4*x5 ];\\
A3=[ x0\^{ }2*x2, x5\^{ }2*x2, x0*x2*x5 ];\\
B=[ x3\^{ }2*x1, x4\^{ }2*x1, x1*x3*x4 ]; \\
C1=[ x0*x3*x5, x0*x4*x5 ]; \\
C2=[ x3\^{ }2*x1, x4\^{ }2*x1, x1*x3*x4, x1*x3*x5, x1*x4*x5, x3*x4*x5 ];\\
BC1=[ x3\^{ }2*x0, x4\^{ }2*x0, x0*x3*x4, x3*x4*x5 ];\\
F0=x0\^{ }2*x1 + 2*x1\^{ }2*x0 + 3*x2\^{ }2*x0 + 7*x3\^{ }2*x2 + 11*x4\^{ }2*x2 + 13*x5\^{ }2*x1 + 17*x1*x2*x5;\\
for iA1 from 0 to (\#A1-1) list ( \\
for iA2 from 0 to (\#A2-1) list ( \\
for iA3 from 0 to (\#A3-1) list ( \\
for iB from 0 to (\#B-1) list ( \\
for iC1 from 0 to (\#C1-1) list ( \\
for iC2 from 0 to (\#C2-1) list ( \\
F=F0 + 19*A1\#iA1 + 20*A2\#iA2 + 21*A3\#iA3 + 23*B\#iB + 31*C1\#iC1 + 37*C2\#iC2;\\ TotalSing=TotalSing + dim singularLocus(ideal(F));\\ 
))))));\\
for iAA from 0 to (\#AA-1) list ( \\
for iB from 0 to (\#B-1) list ( \\
for iC1 from 0 to (\#C1-1) list ( \\
for iC2 from 0 to (\#C2-1) list ( \\
F=F0 + 19*AA\#iAA + 23*B\#iB + 31*C1\#iC1 + 37*C2\#iC2; TotalSing=TotalSing + dim singularLocus(ideal(F));\\ 
))));\\
for iA1 from 0 to (\#A1-1) list ( \\
for iA2 from 0 to (\#A2-1) list ( \\
for iA3 from 0 to (\#A3-1) list ( \\
for iBC1 from 0 to (\#BC1-1) list ( \\
for iC2 from 0 to (\#C2-1) list ( \\
F=F0 + 19*A1\#iA1 + 20*A2\#iA2 + 21*A3\#iA3 + 31*BC1\#iBC1 + 37*C2\#iC2;\\ TotalSing=TotalSing + dim singularLocus(ideal(F));\\ 
)))));\\
for iAA from 0 to (\#AA-1) list ( \\
for iBC1 from 0 to (\#BC1-1) list ( \\
for iC2 from 0 to (\#C2-1) list ( \\
F=F0 + 19*AA\#iAA + 31*BC1\#iBC1 + 37*C2\#iC2; TotalSing=TotalSing + dim singularLocus(ideal(F));\\ 
)));\\
F0=x0\^{ }2*x1 + 2*x1\^{ }2*x0 + 3*x2\^{ }2*x0 + 7*x3\^{ }2*x2 + 11*x4\^{ }2*x2 + 13*x5\^{ }2*x1 + 17*x1*x2*x3;\\
for iA1 from 0 to (\#A1-1) list ( \\
for iA2 from 0 to (\#A2-1) list ( \\
for iA3 from 0 to (\#A3-1) list ( \\
for iB from 0 to (\#B-1) list ( \\
for iC1 from 0 to (\#C1-1) list ( \\
for iC2 from 0 to (\#C2-1) list ( \\
F=F0 + 19*A1\#iA1 + 20*A2\#iA2 + 21*A3\#iA3 + 23*B\#iB + 31*C1\#iC1 + 37*C2\#iC2;\\ TotalSing=TotalSing + dim singularLocus(ideal(F));\\ 
))))));\\
for iAA from 0 to (\#AA-1) list ( \\
for iB from 0 to (\#B-1) list ( \\
for iC1 from 0 to (\#C1-1) list ( \\
for iC2 from 0 to (\#C2-1) list ( \\
F=F0 + 19*AA\#iAA + 23*B\#iB + 31*C1\#iC1 + 37*C2\#iC2; TotalSing=TotalSing + dim singularLocus(ideal(F));\\ 
))));\\
for iA1 from 0 to (\#A1-1) list ( \\
for iA2 from 0 to (\#A2-1) list ( \\
for iA3 from 0 to (\#A3-1) list ( \\
for iBC1 from 0 to (\#BC1-1) list ( \\
for iC2 from 0 to (\#C2-1) list ( \\
F=F0 + 19*A1\#iA1 + 20*A2\#iA2 + 21*A3\#iA3 + 31*BC1\#iBC1 + 37*C2\#iC2;\\ TotalSing=TotalSing + dim singularLocus(ideal(F));\\ 
)))));\\
for iAA from 0 to (\#AA-1) list ( \\
for iBC1 from 0 to (\#BC1-1) list ( \\
for iC2 from 0 to (\#C2-1) list ( \\
F=F0 + 19*AA\#iAA + 31*BC1\#iBC1 + 37*C2\#iC2; TotalSing=TotalSing + dim singularLocus(ideal(F));\\ 
)));\\
F0=x0\^{ }2*x1 + 2*x1\^{ }2*x0 + 3*x2\^{ }2*x0 + 7*x3\^{ }2*x2 + 11*x4\^{ }2*x2 + 13*x5\^{ }2*x1 + 17*x1*x2*x2;\\
for iA1 from 0 to (\#A1-1) list ( \\
for iA2 from 0 to (\#A2-1) list ( \\
for iA3 from 0 to (\#A3-1) list ( \\
for iB from 0 to (\#B-1) list ( \\
for iC1 from 0 to (\#C1-1) list ( \\
for iC2 from 0 to (\#C2-1) list ( \\
F=F0 + 19*A1\#iA1 + 20*A2\#iA2 + 21*A3\#iA3 + 23*B\#iB + 31*C1\#iC1 + 37*C2\#iC2;\\ TotalSing=TotalSing + dim singularLocus(ideal(F));\\ 
))))));\\
for iAA from 0 to (\#AA-1) list ( \\
for iB from 0 to (\#B-1) list ( \\
for iC1 from 0 to (\#C1-1) list ( \\
for iC2 from 0 to (\#C2-1) list ( \\
F=F0 + 19*AA\#iAA + 23*B\#iB + 31*C1\#iC1 + 37*C2\#iC2;\\ TotalSing=TotalSing + dim singularLocus(ideal(F));\\ 
))));\\
for iA1 from 0 to (\#A1-1) list ( \\
for iA2 from 0 to (\#A2-1) list ( \\
for iA3 from 0 to (\#A3-1) list ( \\
for iBC1 from 0 to (\#BC1-1) list ( \\
for iC2 from 0 to (\#C2-1) list ( \\
F=F0 + 19*A1\#iA1 + 20*A2\#iA2 + 21*A3\#iA3 + 31*BC1\#iBC1 + 37*C2\#iC2;\\ TotalSing=TotalSing + dim singularLocus(ideal(F));\\ 
)))));\\
for iAA from 0 to (\#AA-1) list ( \\
for iBC1 from 0 to (\#BC1-1) list ( \\
for iC2 from 0 to (\#C2-1) list ( \\
F=F0 + 19*AA\#iAA + 31*BC1\#iBC1 + 37*C2\#iC2; TotalSing=TotalSing + dim singularLocus(ideal(F));\\ 
)));\\
AA=[ x5\^{ }2*x0, x0*x3*x5, x0*x4*x5 ]; \\
A1=[ x5\^{ }2*x2, x0\^{ }2*x2, x0*x2*x5 ];\\
A2=[ x3\^{ }2*x0, x3\^{ }2*x4, x3\^{ }2*x5, x0*x3*x4, x3*x4*x5 ];\\
BB=[ x4\^{ }2*x1, x1*x3*x4, x1*x4*x5 ]; \\
B1=[ x1*x2*x4 ];\\
B2=[ x3\^{ }2*x1, x3\^{ }2*x4, x3\^{ }2*x5, x1*x3*x5, x3*x4*x5 ];\\
B3=[ x5\^{ }2*x2, x2*x3*x5, x2*x3*x4, x2*x4*x5, x3*x4*x5 ];\\
C=[ x1*x2*x3, x1*x3*x4, x2*x3*x4 ];\\
C0=[ x1*x2*x5, x1*x4*x5, x2*x4*x5 ];\\
E=[ x2*x3*x4, x2*x4*x5 ];\\
EC=[ x4\^{ }2*x1, x1*x2*x4 ]; \\
F0=x0\^{ }2*x1 + 2*x1\^{ }2*x0 + 3*x2\^{ }2*x0 + 7*x3\^{ }2*x2 + 11*x4\^{ }2*x0 + 13*x5\^{ }2*x1 + 17*x1*x2*x3;\\
for iAA from 0 to (\#AA-1) list ( \\
for iBB from 0 to (\#BB-1) list ( \\
for iB3 from 0 to (\#B3-1) list ( \\
for iEC from 0 to (\#EC-1) list ( \\
F=F0 + 19*AA\#iAA + 23*BB\#iBB + 31*B3\#iB3 + 37*EC\#iEC; TotalSing=TotalSing + dim singularLocus(ideal(F));\\ 
))));\\
for iA1 from 0 to (\#A1-1) list ( \\
for iA2 from 0 to (\#A2-1) list ( \\
for iBB from 0 to (\#BB-1) list ( \\
for iB3 from 0 to (\#B3-1) list ( \\
for iEC from 0 to (\#EC-1) list ( \\
F=F0 + 19*A1\#iA1 + 20*A2\#iA2 + 23*BB\#iBB + 31*B3\#iB3 + 37*EC\#iEC;\\ TotalSing=TotalSing + dim singularLocus(ideal(F));\\ 
)))));\\
for iAA from 0 to (\#AA-1) list ( \\
for iB1 from 0 to (\#B1-1) list ( \\
for iB2 from 0 to (\#B2-1) list ( \\
for iB3 from 0 to (\#B3-1) list ( \\
for iEC from 0 to (\#EC-1) list ( \\
F=F0 + 19*AA\#iAA + 23*B1\#iB1 + 29*B2\#iB2 + 31*B3\#iB3 + 37*EC\#iEC;\\ TotalSing=TotalSing + dim singularLocus(ideal(F));\\ 
)))));\\
for iA1 from 0 to (\#A1-1) list ( \\
for iA2 from 0 to (\#A2-1) list ( \\
for iB1 from 0 to (\#B1-1) list ( \\
for iB2 from 0 to (\#B2-1) list ( \\
for iB3 from 0 to (\#B3-1) list ( \\
for iEC from 0 to (\#EC-1) list ( \\
F=F0 + 19*A1\#iA1 + 20*A2\#iA2 + 23*B1\#iB1 + 29*B2\#iB2 + 31*B3\#iB3 + 37*EC\#iEC;\\ TotalSing=TotalSing + dim singularLocus(ideal(F));\\ 
))))));\\
for iAA from 0 to (\#AA-1) list ( \\
for iBB from 0 to (\#BB-1) list ( \\
for iB3 from 0 to (\#B3-1) list ( \\
for iC0 from 0 to (\#C0-1) list ( \\
for iE from 0 to (\#E-1) list ( \\
F=F0 + 19*AA\#iAA + 23*BB\#iBB + 31*B3\#iB3 + 37*C0\#iC0 + 41*E\#iE;\\ TotalSing=TotalSing + dim singularLocus(ideal(F));\\ 
)))));\\
for iA1 from 0 to (\#A1-1) list ( \\
for iA2 from 0 to (\#A2-1) list ( \\
for iBB from 0 to (\#BB-1) list ( \\
for iB3 from 0 to (\#B3-1) list ( \\
for iC0 from 0 to (\#C0-1) list ( \\
for iE from 0 to (\#E-1) list ( \\
F=F0 + 19*A1\#iA1 + 20*A2\#iA2 + 23*BB\#iBB + 31*B3\#iB3 + 37*C0\#iC0 + 41*E\#iE;\\ TotalSing=TotalSing + dim singularLocus(ideal(F));\\ 
))))));\\
for iAA from 0 to (\#AA-1) list ( \\
for iB1 from 0 to (\#B1-1) list ( \\
for iB2 from 0 to (\#B2-1) list ( \\
for iB3 from 0 to (\#B3-1) list ( \\
for iC0 from 0 to (\#C0-1) list ( \\
for iE from 0 to (\#E-1) list ( \\
F=F0 + 19*AA\#iAA + 23*B1\#iB1 + 29*B2\#iB2 + 31*B3\#iB3 + 37*C0\#iC0 + 41*E\#iE;\\ TotalSing=TotalSing + dim singularLocus(ideal(F));\\ 
))))));\\
for iA1 from 0 to (\#A1-1) list ( \\
for iA2 from 0 to (\#A2-1) list ( \\
for iB1 from 0 to (\#B1-1) list ( \\
for iB2 from 0 to (\#B2-1) list ( \\
for iB3 from 0 to (\#B3-1) list ( \\
for iC0 from 0 to (\#C0-1) list ( \\
for iE from 0 to (\#E-1) list ( \\
F=F0 + 19*A1\#iA1 + 20*A2\#iA2 + 23*B1\#iB1 + 29*B2\#iB2 + 31*B3\#iB3 + 37*C0\#iC0 + 41*E\#iE;\\ TotalSing=TotalSing + dim singularLocus(ideal(F));\\ 
)))))));\\
F0=x0\^{ }2*x1 + 2*x1\^{ }2*x0 + 3*x2\^{ }2*x0 + 7*x3\^{ }2*x2 + 11*x4\^{ }2*x0 + 13*x5\^{ }2*x1 + 17*x1*x2*x5;\\
for iAA from 0 to (\#AA-1) list ( \\
for iBB from 0 to (\#BB-1) list ( \\
for iB3 from 0 to (\#B3-1) list ( \\
for iEC from 0 to (\#EC-1) list ( \\
F=F0 + 19*AA\#iAA + 23*BB\#iBB + 31*B3\#iB3 + 37*EC\#iEC; TotalSing=TotalSing + dim singularLocus(ideal(F));\\ 
))));\\
for iA1 from 0 to (\#A1-1) list ( \\
for iA2 from 0 to (\#A2-1) list ( \\
for iBB from 0 to (\#BB-1) list ( \\
for iB3 from 0 to (\#B3-1) list ( \\
for iEC from 0 to (\#EC-1) list ( \\
F=F0 + 19*A1\#iA1 + 20*A2\#iA2 + 23*BB\#iBB + 31*B3\#iB3 + 37*EC\#iEC;\\ TotalSing=TotalSing + dim singularLocus(ideal(F));\\ 
)))));\\
for iAA from 0 to (\#AA-1) list ( \\
for iB1 from 0 to (\#B1-1) list ( \\
for iB2 from 0 to (\#B2-1) list ( \\
for iB3 from 0 to (\#B3-1) list ( \\
for iEC from 0 to (\#EC-1) list ( \\
F=F0 + 19*AA\#iAA + 23*B1\#iB1 + 29*B2\#iB2 + 31*B3\#iB3 + 37*EC\#iEC;\\ TotalSing=TotalSing + dim singularLocus(ideal(F));\\ 
)))));\\
for iA1 from 0 to (\#A1-1) list ( \\
for iA2 from 0 to (\#A2-1) list ( \\
for iB1 from 0 to (\#B1-1) list ( \\
for iB2 from 0 to (\#B2-1) list ( \\
for iB3 from 0 to (\#B3-1) list ( \\
for iEC from 0 to (\#EC-1) list ( \\
F=F0 + 19*A1\#iA1 + 20*A2\#iA2 + 23*B1\#iB1 + 29*B2\#iB2 + 31*B3\#iB3 + 37*EC\#iEC;\\ TotalSing=TotalSing + dim singularLocus(ideal(F));\\ 
))))));\\
for iAA from 0 to (\#AA-1) list ( \\
for iBB from 0 to (\#BB-1) list ( \\
for iB3 from 0 to (\#B3-1) list ( \\
for iC from 0 to (\#C-1) list ( \\
for iE from 0 to (\#E-1) list ( \\
F=F0 + 19*AA\#iAA + 23*BB\#iBB + 31*B3\#iB3 + 37*C\#iC + 41*E\#iE;\\ TotalSing=TotalSing + dim singularLocus(ideal(F));\\ 
)))));\\
for iA1 from 0 to (\#A1-1) list ( \\
for iA2 from 0 to (\#A2-1) list ( \\
for iBB from 0 to (\#BB-1) list ( \\
for iB3 from 0 to (\#B3-1) list ( \\
for iC from 0 to (\#C-1) list ( \\
for iE from 0 to (\#E-1) list ( \\
F=F0 + 19*A1\#iA1 + 20*A2\#iA2 + 23*BB\#iBB + 31*B3\#iB3 + 37*C\#iC + 41*E\#iE;\\ TotalSing=TotalSing + dim singularLocus(ideal(F));\\ 
))))));\\
for iAA from 0 to (\#AA-1) list ( \\
for iB1 from 0 to (\#B1-1) list ( \\
for iB2 from 0 to (\#B2-1) list ( \\
for iB3 from 0 to (\#B3-1) list ( \\
for iC from 0 to (\#C-1) list ( \\
for iE from 0 to (\#E-1) list ( \\
F=F0 + 19*AA\#iAA + 23*B1\#iB1 + 29*B2\#iB2 + 31*B3\#iB3 + 37*C\#iC + 41*E\#iE;\\ TotalSing=TotalSing + dim singularLocus(ideal(F));\\ 
))))));\\
for iA1 from 0 to (\#A1-1) list ( \\
for iA2 from 0 to (\#A2-1) list ( \\
for iB1 from 0 to (\#B1-1) list ( \\
for iB2 from 0 to (\#B2-1) list ( \\
for iB3 from 0 to (\#B3-1) list ( \\
for iC from 0 to (\#C-1) list ( \\
for iE from 0 to (\#E-1) list ( \\
F=F0 + 19*A1\#iA1 + 20*A2\#iA2 + 23*B1\#iB1 + 29*B2\#iB2 + 31*B3\#iB3 + 37*C\#iC + 41*E\#iE;\\ TotalSing=TotalSing + dim singularLocus(ideal(F));\\ 
)))))));\\
F0=x0\^{ }2*x1 + 2*x1\^{ }2*x0 + 3*x2\^{ }2*x0 + 7*x3\^{ }2*x2 + 11*x4\^{ }2*x0 + 13*x5\^{ }2*x1 + 17*x1*x2*x4;\\
for iAA from 0 to (\#AA-1) list ( \\
for iBB from 0 to (\#BB-1) list ( \\
for iB3 from 0 to (\#B3-1) list ( \\
F=F0 + 19*AA\#iAA + 23*BB\#iBB + 31*B3\#iB3; TotalSing=TotalSing + dim singularLocus(ideal(F));\\ 
)));\\
for iA1 from 0 to (\#A1-1) list ( \\
for iA2 from 0 to (\#A2-1) list ( \\
for iBB from 0 to (\#BB-1) list ( \\
for iB3 from 0 to (\#B3-1) list ( \\
F=F0 + 19*A1\#iA1 + 20*A2\#iA2 + 23*BB\#iBB + 31*B3\#iB3; TotalSing=TotalSing + dim singularLocus(ideal(F));\\ 
))));\\
for iAA from 0 to (\#AA-1) list ( \\
for iB1 from 0 to (\#B1-1) list ( \\
for iB2 from 0 to (\#B2-1) list ( \\
for iB3 from 0 to (\#B3-1) list ( \\
F=F0 + 19*AA\#iAA + 23*B1\#iB1 + 29*B2\#iB2 + 31*B3\#iB3; TotalSing=TotalSing + dim singularLocus(ideal(F));\\ 
))));\\
for iA1 from 0 to (\#A1-1) list ( \\
for iA2 from 0 to (\#A2-1) list ( \\
for iB1 from 0 to (\#B1-1) list ( \\
for iB2 from 0 to (\#B2-1) list ( \\
for iB3 from 0 to (\#B3-1) list ( \\
F=F0 + 19*A1\#iA1 + 20*A2\#iA2 + 23*B1\#iB1 + 29*B2\#iB2 + 31*B3\#iB3;\\ TotalSing=TotalSing + dim singularLocus(ideal(F));\\ 
)))));\\
F0=x0\^{ }2*x1 + 2*x1\^{ }2*x0 + 3*x2\^{ }2*x0 + 7*x3\^{ }2*x2 + 11*x4\^{ }2*x0 + 13*x5\^{ }2*x1 + 17*x1*x2*x2;\\
for iAA from 0 to (\#AA-1) list ( \\
for iBB from 0 to (\#BB-1) list ( \\
for iB3 from 0 to (\#B3-1) list ( \\
F=F0 + 19*AA\#iAA + 23*BB\#iBB + 31*B3\#iB3; TotalSing=TotalSing + dim singularLocus(ideal(F));\\ 
)));\\
for iA1 from 0 to (\#A1-1) list ( \\
for iA2 from 0 to (\#A2-1) list ( \\
for iBB from 0 to (\#BB-1) list ( \\
for iB3 from 0 to (\#B3-1) list ( \\
F=F0 + 19*A1\#iA1 + 20*A2\#iA2 + 23*BB\#iBB + 31*B3\#iB3; TotalSing=TotalSing + dim singularLocus(ideal(F));\\ 
))));\\
for iAA from 0 to (\#AA-1) list ( \\
for iB1 from 0 to (\#B1-1) list ( \\
for iB2 from 0 to (\#B2-1) list ( \\
for iB3 from 0 to (\#B3-1) list ( \\
F=F0 + 19*AA\#iAA + 23*B1\#iB1 + 29*B2\#iB2 + 31*B3\#iB3; TotalSing=TotalSing + dim singularLocus(ideal(F));\\ 
))));\\
for iA1 from 0 to (\#A1-1) list ( \\
for iA2 from 0 to (\#A2-1) list ( \\
for iB1 from 0 to (\#B1-1) list ( \\
for iB2 from 0 to (\#B2-1) list ( \\
for iB3 from 0 to (\#B3-1) list ( \\
F=F0 + 19*A1\#iA1 + 20*A2\#iA2 + 23*B1\#iB1 + 29*B2\#iB2 + 31*B3\#iB3;\\ TotalSing=TotalSing + dim singularLocus(ideal(F));\\ 
)))));\\
F1=[ x1*x3*x4, x1*x3*x5, x3*x4*x5 ]; \\
F2=[ x1*x3*x4, x1*x3*x5, x3*x4*x5 ]; \\
F3=[ x1*x2*x5, x1*x3*x5 ]; \\
F4=[ x1*x2*x5, x1*x3*x5 ]; \\
A1=[ x1*x2*x5, x1*x3*x5 ]; \\
A2=[ x1*x2*x4, x1*x3*x4 ]; \\
B=[ x0*x2*x3, x0*x2*x4, x0*x2*x5 ];\\
C=[ x2*x3*x4, x2*x4*x5, x2*x3*x5, x3*x4*x5 ]; \\
F00 = x0\^{ }2*x1 + 2*x1\^{ }2*x0 + 3*x2\^{ }2*x1 + 7*x3\^{ }2*x2 + 11*x4\^{ }2*x0 + 13*x5\^{ }2*x0; \\
F0=F00 + 17*x2*x4*x5;\\
for iA1 from 0 to (\#A1-1) list ( \\
for iA2 from 0 to (\#A2-1) list ( \\
for iB from 0 to (\#B-1) list ( \\
for iF1 from 0 to (\#F1-1) list ( \\
F=F0 + 19*A1\#iA1 + 23*A2\#iA2 + 31*F1\#iF1 + 29*B\#iB; TotalSing=TotalSing + dim singularLocus(ideal(F));\\ 
))));\\
F0=F00 + 17*x2*x4*x5 + 37*x2*x2*x0;\\
for iA1 from 0 to (\#A1-1) list ( \\
for iA2 from 0 to (\#A2-1) list ( \\
for iF2 from 0 to (\#F2-1) list ( \\
F=F0 + 19*A1\#iA1 + 23*A2\#iA2 + 29*F2\#iF2; TotalSing=TotalSing + dim singularLocus(ideal(F));\\ 
)));\\
F0=F00 + 17*x1*x4*x5;\\
for iB from 0 to (\#B-1) list ( \\
for iC from 0 to (\#C-1) list ( \\
F=F0 + 19*B\#iB + 23*C\#iC; TotalSing=TotalSing + dim singularLocus(ideal(F));\\ 
));\\ 
F0=F00 + 17*x3*x4*x5+ 41*x1*x2*x4;\\
for iB from 0 to (\#B-1) list ( \\
for iF4 from 0 to (\#F4-1) list ( \\
F=F0 + 19*B\#iB + 23*F4\#iF4; TotalSing=TotalSing + dim singularLocus(ideal(F));\\ 
));\\ 
F0=F00 + 17*x3*x4*x5+ 41*x1*x2*x4+23*x0*x2*x2;\\
for iF3 from 0 to (\#F3-1) list ( \\
F=F0 + 19*F3\#iF3; TotalSing=TotalSing + dim singularLocus(ideal(F));\\ 
);\\
F0=F00 + 17*x1*x4*x5+ 23*x0*x2*x2;\\
for iC from 0 to (\#C-1) list ( \\
F=F0 + 19*C\#iC; TotalSing=TotalSing + dim singularLocus(ideal(F));\\ 
);\\
F1=[ x4\^{ }2*x2, x1*x2*x5, x1*x4*x5, x2*x4*x5 ]; \\
F2=[ x4\^{ }2*x2, x2*x3*x4, x2*x4*x5 ]; \\
F3=[ x0*x3*x4, x0*x3*x5, x0*x4*x5 ];\\
F4=[ x0*x3*x4, x1*x3*x4 ];\\
F5=[ x0*x3*x5, x1*x3*x5 ];\\
A01=[ x0*x3*x4, x0*x3*x5, x0*x4*x5 ];\\
A02=[ x0*x4*x5, x1*x4*x5 ];\\
A03=[ x0*x3*x5, x1*x3*x5 ];\\
A1=[ x3\^{ }2*x0, x5\^{ }2*x0, x0*x3*x5, x1*x3*x5, x3*x4*x5 ];\\
A2=[ x3\^{ }2*x0, x5\^{ }2*x0, x0\^{ }2*x4, x0*x3*x4, x0*x3*x5, x0*x4*x5, x3*x4*x5 ];\\
A3=[ x1*x3*x4, x1*x3*x5, x1*x4*x5, x3*x4*x5 ];\\
C=[ x4\^{ }2*x2, x1*x3*x4, x1*x4*x5 ]; \\
F00 = x0\^{ }2*x1 + 2*x1\^{ }2*x0 + 3*x2\^{ }2*x0 + 7*x3\^{ }2*x2 + 11*x4\^{ }2*x2 + 13*x5\^{ }2*x2; \\
F0=F00 + 17*x1*x2*x3 + 61*x1*x3*x4;\\
for iA01 from 0 to (\#A01-1) list ( \\
for iA02 from 0 to (\#A02-1) list ( \\
for iA03 from 0 to (\#A03-1) list ( \\
F=F0 + 19*A01\#iA01 + 23*A02\#iA02 + 29*A03\#iA03; TotalSing=TotalSing + dim singularLocus(ideal(F));\\ 
)));\\
F0=F00 + 17*x1*x2*x2 + 61*x1*x3*x4;\\
for iA01 from 0 to (\#A01-1) list ( \\
for iA02 from 0 to (\#A02-1) list ( \\
for iA03 from 0 to (\#A03-1) list ( \\
F=F0 + 19*A01\#iA01 + 23*A02\#iA02 + 29*A03\#iA03; TotalSing=TotalSing + dim singularLocus(ideal(F));\\ 
)));\\
F0=F00 + 17*x1*x2*x3 + 61*x1*x4*x5;\\
for iF3 from 0 to (\#F3-1) list ( \\
for iF4 from 0 to (\#F4-1) list ( \\
for iF5 from 0 to (\#F5-1) list ( \\
F=F0 + 19*F3\#iF3 + 23*F4\#iF4 + 29*F5\#iF5; TotalSing=TotalSing + dim singularLocus(ideal(F));\\ 
)));\\
F=F00 + 17*x1*x2*x3 + 61*x3*x4*x5; TotalSing=TotalSing + dim singularLocus(ideal(F));\\ 
F=F00 + 17*x1*x2*x2 + 61*x3*x4*x5; TotalSing=TotalSing + dim singularLocus(ideal(F));\\ 
F00 = x0\^{ }2*x1 + 2*x1\^{ }2*x0 + 3*x2\^{ }2*x0 + 7*x3\^{ }2*x2 + 11*x4\^{ }2*x0 + 13*x5\^{ }2*x2; \\
F0=F00 + 17*x1*x2*x4;\\
for iA1 from 0 to (\#A1-1) list ( \\
for iA2 from 0 to (\#A2-1) list ( \\
for iA3 from 0 to (\#A3-1) list ( \\
F=F0 + 19*A1\#iA1 + 23*A2\#iA2 + 29*A3\#iA3; TotalSing=TotalSing + dim singularLocus(ideal(F));\\ 
)));\\
F0=F00 + 17*x1*x2*x3;\\
for iA1 from 0 to (\#A1-1) list ( \\
for iA2 from 0 to (\#A2-1) list ( \\
for iA3 from 0 to (\#A3-1) list ( \\
for iF1 from 0 to (\#F1-1) list ( \\
for iF2 from 0 to (\#F2-1) list ( \\
for iC from 0 to (\#C-1) list ( \\
F=F0 + 19*A1\#iA1 + 23*A2\#iA2 + 29*A3\#iA3 + 31*F1\#iF1 + 37*F2\#iF2 + 41*C\#iC;\\ TotalSing=TotalSing + dim singularLocus(ideal(F));\\ 
))))));\\
F1=[ x2\^{ }2*x1, x3\^{ }2*x1, x1*x2*x5, x1*x3*x5, x2*x3*x5 ]; \\
F2=[ x3\^{ }2*x1, x1*x3*x4, x1*x3*x5 ]; \\
F3=[ x2\^{ }2*x1, x1*x2*x4, x1*x2*x5 ]; \\
F4=[ x2\^{ }2*x1, x3\^{ }2*x1, x2*x3*x4 ]; \\
F5=[ x1*x2*x3, x1*x2*x5, x1*x3*x5, x2*x3*x5 ];\\
F6=[ x1*x3*x4, x1*x3*x5, x1*x4*x5, x3*x4*x5 ];\\
F7=[ x1*x2*x5, x1*x3*x5, x1*x4*x5 ];\\
F8=[ x1*x2*x4, x1*x3*x4, x1*x4*x5 ];\\
F9=[ x1*x2*x3, x1*x2*x4 ];\\
F10=[ x1*x2*x5, x1*x3*x5, x1*x4*x5 ];\\
F11=[ x1*x4*x5, x1*x2*x4 ];\\
A=[ x2*x3*x4, x2*x3*x5, x2*x4*x5, x3*x4*x5 ];\\
B1=[ x2*x3*x4, x1*x2*x4 ];\\
B2=[ x2*x3*x5, x1*x2*x5 ];\\
B3=[ x1*x2*x3, x1*x2*x4, x1*x2*x5 ];\\
C1=[ x4\^{ }2*x0, x5\^{ }2*x0, x2*x4*x5 ]; \\
C2=[ x4\^{ }2*x0, x2*x4*x0, x3*x4*x0 ]; \\
C3=[ x5\^{ }2*x0, x2*x5*x0, x3*x5*x0 ];\\
C4=[ x4\^{ }2*x0, x3*x4*x5, x3*x4*x0, x5\^{ }2*x0, x3*x5*x0 ];\\
F00 = x0\^{ }2*x1 + 2*x1\^{ }2*x0 + 3*x2\^{ }2*x0 + 7*x3\^{ }2*x0 + 11*x4\^{ }2*x1 + 13*x5\^{ }2*x1; \\
F0=F00;\\
for iF1 from 0 to (\#F1-1) list ( \\
for iF2 from 0 to (\#F2-1) list ( \\
for iF3 from 0 to (\#F3-1) list ( \\
for iF4 from 0 to (\#F4-1) list ( \\
for iA from 0 to (\#A-1) list ( \\
for iC1 from 0 to (\#C1-1) list ( \\
for iC2 from 0 to (\#C2-1) list ( \\
for iC3 from 0 to (\#C3-1) list ( \\
for iC4 from 0 to (\#C4-1) list ( \\
F=F0 + 19*F1\#iF1 + 23*F2\#iF2 + 29*F3\#iF3 + 31*F4\#iF4 + 37*A\#iA + 41*C1\#iC1 + 43*C2\#iC2 + 47*C3\#iC3 + 51*C4\#iC4;\\ TotalSing=TotalSing + dim singularLocus(ideal(F));\\
)))))))));\\
F0=F00+17*x1*x2*x3;\\
for iA from 0 to (\#A-1) list ( \\
for iC1 from 0 to (\#C1-1) list ( \\
for iC2 from 0 to (\#C2-1) list ( \\
for iC3 from 0 to (\#C3-1) list ( \\
for iC4 from 0 to (\#C4-1) list ( \\
F=F0 + 19*C1\#iC1 + 23*C2\#iC2 + 29*C3\#iC3 + 31*C4\#iC4 + 37*A\#iA;\\ TotalSing=TotalSing + dim singularLocus(ideal(F));\\ 
)))));\\
F0=F00+17*x1*x2*x3+41*x0*x4*x5;\\
for iA from 0 to (\#A-1) list ( \\
F=F0 + 19*A\#iA; TotalSing=TotalSing + dim singularLocus(ideal(F));\\ 
);\\
F00 = x0\^{ }2*x1 + 2*x1\^{ }2*x0 + 3*x2\^{ }2*x0 + 7*x3\^{ }2*x2 + 11*x4\^{ }2*x0 + 13*x5\^{ }2*x0; \\
F0=F00 + 17*x2*x4*x5;\\
for iF5 from 0 to (\#F5-1) list ( \\
for iF6 from 0 to (\#F6-1) list ( \\
for iF7 from 0 to (\#F7-1) list ( \\
for iF8 from 0 to (\#F8-1) list ( \\
for iF9 from 0 to (\#F9-1) list ( \\
F=F0 + 19*F5\#iF5 + 23*F6\#iF6 + 29*F7\#iF7 + 31*F8\#iF8 + 37*F9\#iF9;\\ TotalSing=TotalSing + dim singularLocus(ideal(F));\\ 
)))));\\
F0=F00 + 17*x3*x4*x5;\\
for iF10 from 0 to (\#F10-1) list ( \\
for iF11 from 0 to (\#F11-1) list ( \\
for iB1 from 0 to (\#B1-1) list ( \\
for iB2 from 0 to (\#B2-1) list ( \\
for iB3 from 0 to (\#B3-1) list ( \\
F=F0 + 19*F10\#iF10 + 23*F11\#iF11 + 29*B1\#iB1 + 31*B2\#iB2 + 37*B3\#iB3;\\ TotalSing=TotalSing + dim singularLocus(ideal(F));\\ 
)))));\\
F0=F00 + 17*x1*x4*x5 + 41*x2*x3*x5;\\
for iB1 from 0 to (\#B1-1) list ( \\
for iB2 from 0 to (\#B2-1) list ( \\
for iB3 from 0 to (\#B3-1) list ( \\
F=F0 + 19*B1\#iB1 + 31*B2\#iB2 + 23*B3\#iB3; TotalSing=TotalSing + dim singularLocus(ideal(F));\\ 
)));\\
F1=[ x2\^{ }2*x0, x0*x2*x3, x0*x2*x4, x0*x2*x5 ]; \\
F2=[ x1*x2*x3, x1*x2*x4, x1*x3*x4, x2*x3*x4 ]; \\
F3=[ x1*x2*x3, x1*x2*x5, x1*x3*x5, x2*x3*x5 ]; \\
F4=[ x1*x2*x4, x1*x2*x5, x1*x4*x5, x2*x4*x5 ];\\
F5=[ x1*x2*x3, x1*x3*x4, x1*x3*x5 ];\\
F6=[ x1*x2*x4, x1*x3*x4, x1*x4*x5 ];\\
F7=[ x1*x2*x5, x1*x3*x5, x1*x4*x5 ];\\
F8=[ x2\^{ }2*x0, x0*x2*x4, x0*x2*x5, x0*x2*x3 ]; \\
F9=[ x1*x2*x5, x1*x3*x5, x1*x4*x5 ];\\
F10=[ x1*x3*x5, x2*x3*x5 ];\\
F11=[ x2\^{ }2*x0, x0*x2*x4, x0*x2*x5, x0*x2*x3 ]; \\
G1=[ x1*x2*x5, x2*x3*x5 ]; \\
G2=[ x1*x3*x5, x2*x3*x5, x3*x4*x5 ]; \\
G3=[ x1*x4*x5, x2*x4*x5, x3*x4*x5 ]; \\
G4=[ x1*x4*x5, x1*x2*x5, x1*x3*x5 ]; \\
G5=[ x1*x4*x5, x1*x2*x4, x1*x3*x4 ]; \\
G6=[ x1*x3*x5, x1*x2*x3, x1*x3*x4 ]; \\
G7=[ x1*x2*x3, x1*x2*x5, x1*x2*x4 ]; \\
G8=[ x1*x4*x5, x2*x4*x5, x1*x2*x5, x1*x2*x4 ]; \\
G9=[ x1*x4*x5, x3*x4*x5, x1*x3*x5, x1*x3*x4 ];\\
F00 = x0\^{ }2*x1 + 2*x1\^{ }2*x0 + 3*x2\^{ }2*x1 + 7*x3\^{ }2*x0 + 11*x4\^{ }2*x0 + 13*x5\^{ }2*x0; \\
F0=F00 + 61*x3*x4*x5;\\
for iF1 from 0 to (\#F1-1) list ( \\
for iF2 from 0 to (\#F2-1) list ( \\
for iF3 from 0 to (\#F3-1) list ( \\
for iF4 from 0 to (\#F4-1) list ( \\
for iF5 from 0 to (\#F5-1) list ( \\
for iF6 from 0 to (\#F6-1) list ( \\
for iF7 from 0 to (\#F7-1) list ( \\
F=F0 + 19*F1\#iF1 + 23*F2\#iF2 + 29*F3\#iF3 + 31*F4\#iF4 + 37*F5\#iF5 + 41*F6\#iF6 + 43*F7\#iF7;\\ TotalSing=TotalSing + dim singularLocus(ideal(F));\\ 
)))))));\\
F0=F00 + 61*x3*x4*x1+17*x2*x4*x5;\\
for iF8 from 0 to (\#F8-1) list ( \\
for iF9 from 0 to (\#F9-1) list ( \\
for iF10 from 0 to (\#F10-1) list ( \\
F=F0 + 19*F8\#iF8 + 23*F9\#iF9 + 29*F10\#iF10; TotalSing=TotalSing + dim singularLocus(ideal(F));\\ 
)));\\
F0=F00 + 61*x3*x4*x1+17*x2*x3*x4 + 31*x3*x5*x1+37*x1*x4*x5;\\
for iF11 from 0 to (\#F11-1) list ( \\
F=F0 + 19*F11\#iF11; TotalSing=TotalSing + dim singularLocus(ideal(F));\\ 
); \\
F00 = x0\^{ }2*x1 + 2*x1\^{ }2*x0 + 3*x2\^{ }2*x0 + 7*x3\^{ }2*x0 + 11*x4\^{ }2*x0 + 13*x5\^{ }2*x0; \\
F0=F00 + 61*x3*x4*x2;\\
for iG1 from 0 to (\#G1-1) list ( \\
for iG2 from 0 to (\#G2-1) list ( \\
for iG3 from 0 to (\#G3-1) list ( \\
for iG4 from 0 to (\#G4-1) list ( \\
for iG5 from 0 to (\#G5-1) list ( \\
for iG6 from 0 to (\#G6-1) list ( \\
for iG7 from 0 to (\#G7-1) list ( \\
for iG8 from 0 to (\#G8-1) list ( \\
for iG9 from 0 to (\#G9-1) list ( \\
F=F0 + 19*G1\#iG1 + 23*G2\#iG2 + 29*G3\#iG3 + 31*G4\#iG4 + 37*G5\#iG5 + 41*G6\#iG6 + 43*G7\#iG7 + 47*G8\#iG8 + 51*G9\#iG9;\\ TotalSing=TotalSing + dim singularLocus(ideal(F));\\ 
)))))))));\\
print concatenate("TotalSing=", toString TotalSing);\\
if TotalSing==0 then print "EVERYTHING WAS SMOOTH!";\\
if TotalSing$>0$ then print "ATTENTION! SOMETHING WAS NOT SMOOTH!";\\

} 

\vspace{2cm}

We observed output $0$, from which we concluded that all sets $\AA$ we found during the proof of Theorem $1$ are indeed smooth.\\

\section{Acknowledgement}

We learned about the problem of describing automorphism groups of cubic fourfolds from Yuri Zarhin and a computation of the prime order automorphisms of cubic fourfolds from paper \cite{Liendo}, which together with paper \cite{LieFu} contains the key ideas (in particular, the remark about the semisimplicity of a finite order linear endomorphism) used in this note.\\

\bibliographystyle{ams-plain}

\bibliography{AbelianAutomorphismsOfCubic4folds}

\end{document}